# Diagram of Representations of Universal Algebras

Aleks Kleyn

Aleks_Kleyn@MailAPS.org
http://AleksKleyn.dyndns-home.com:4080/
http://sites.google.com/site/AleksKleyn/
http://arxiv.org/a/kleyn_a_1
http://AleksKleyn.blogspot.com/

ABSTRACT. Theory of representations of universal algebra is a natural development of the theory of universal algebra. In the book, I considered representation of universal algebra, diagram of representations and examples of representation. Morphism of the representation is the map that conserve the structure of the representation. Exploring of morphisms of the representation leads to the concepts of generating set and basis of representation.

# Contents

CHAPTER 1

# Preface

## 1.1. Representation Theory

In my papers, I often explore problems relating to the representation of universal algebra. Initially it was small sketches which I repeatedly corrected and rewrote. However gradually there were new observations. As a result, auxiliary tool became a consistent theory.

I realized this when I was writing book [10], and I decided to dedicate a separate book to the questions related with representation of universal algebra. Exploring of the theory of representations of universal algebra shows that this theory has a lot of common with theory of universal algebra.

The definition of vector space as representation of field in the Abelian group was the main impetus of deeper exploring of representations of universal algebra. I put attention that this definition changes role of linear map. It was found that linear map is the map that preserves the structure of the representation. It is easy to generalize this structure for an arbitrary representation of universal algebra. Thus I came to the notion of morphism of representation.

The set of regular automorphisms of vector space forms a group. This group is single transitive on the set of basises of vector space. This statement is the foundation of the theory of invariants of vector space.

The natural question arises. Can we generalize this structure to arbitrary representation? The basis is not the only set that forms the vector space. If we add an arbitrary vector to the set of vectors of basis, then a new set also generates the same vector space, however this set is not basis. This statement is initial point where I started exploring of generating set of representation. Generating set of representation is one more interesting parallel between theory of representations and theory of universal algebra.

The set of automorphisms of representations is loop. Nonassociativity of the product is the source of numerous questions which require additional research. All these questions lead to the need to understand the theory of invariants of a given representation.

If we consider the theory of representations of universal algebra as an extension of the theory of universal algebra, then why not consider the representation of one representation in another representation. Thus the concept of the tower representations appeared. The most amazing fact is the statement that all maps in the tower of representations are coordinated.

## 1.2. On the Edge of Theory

Over the years, I believed that representation theory is the main tool to study covariance principle. However, in the process of writing this book, I suddenly found myself on the edge of representation theory. It was extremely important event.

More precisely, it was two different discoveries, interconnected by topic of non commutative addition. At the beginning, I discovered that I can study affine geometry on affine manifold. (This is not new discovery. I think people have known about this since Descartes and Gauss). The most important for me here was the statement that sum is not defined for every pairs of vectors. I met similar problem when I was studying basis manifold of Minkowski space ([11]). If connection on affine manifold has nonzero torsion, then sum of vectors becomes non commutative.

Later, I decided to study representation of ring in non-abelian group. Although algebra is closed relative operation, I see opportunity for further development of representation theory. We can use the definition of basis from this book; however some important details will be hiden. I am interested in the version that elements of basis may have a given order; but right now I do not have a clear idea of what may follow from this assumption.

# Preliminary Definitions

This chapter contains definitions and theorems which are necessary for an understanding of the text of this book. So the reader may read the statements from this chapter in process of reading the main text of the book.

## 2.1. Equivalence

DEFINITION 2.1.1. *Correspondence* $\Phi \in A \times A$ *is called* **equivalence**, *if*[2.1]

2.1.1.1: *correspondence* $\Phi$ *is* **reflexive**

$$(a, a) \in \Phi$$

2.1.1.2: *correspondence* $\Phi$ *is* **symmetric**

$$(a, b) \in \Phi \Rightarrow (b, a) \in \Phi$$

2.1.1.3: *correspondence* $\Phi$ *is* **transitive**

$$(a, b), (b, c) \in \Phi \Rightarrow (a, c) \in \Phi$$

$\square$

THEOREM 2.1.2. *For the map*

$$f : A \to B$$

*the set*

$$(2.1.1) \qquad \ker f = \{(a, b) : a, b \in A, f(a) = f(b)\}$$

*is equivalence and is called* **kernel of map**.[2.2]

PROOF.

LEMMA 2.1.3. *Correspondence* $\ker f$ *is reflexive.*

PROOF. From the equality

$$f(a) = f(a)$$

and from the definition (2.1.1), it follows that

$$(2.1.2) \qquad\qquad\qquad (a, a) \in \ker f$$

The lemma follows from the statement (2.1.2) and from the definition 2.1.1.1. $\odot$

LEMMA 2.1.4. *Correspondence* $\ker f$ *is symmetric.*

---

PROOF. The equality

(2.1.3)                          $$f(a) = f(b)$$

follows from the statement

$$(a, b) \in \ker f$$

and from the definition (2.1.1). The equality

(2.1.4)                          $$f(b) = f(a)$$

follows from the equality (2.1.3). The statement

$$(b, a) \in \ker f$$

follows from the equality (2.1.4) and from the definition (2.1.1). Therefore, we proved the statement

(2.1.5)                  $$(a, b) \in \ker f \Rightarrow (b, a) \in \ker f$$

The lemma follows from the statement (2.1.5) and from the definition 2.1.1.2.   ⊙

LEMMA 2.1.5. *Correspondence* $\ker f$ *is transitive.*

PROOF. The equality

(2.1.6)                          $$f(a) = f(b)$$

follows from the statement

$$(a, b) \in \ker f$$

and from the definition (2.1.1). The equality

(2.1.7)                          $$f(b) = f(c)$$

follows from the statement

$$(b, c) \in \ker f$$

and from the definition (2.1.1). The equality

(2.1.8)                          $$f(a) = f(c)$$

follows from equalities (2.1.6), (2.1.7). The statement

$$(a, c) \in \ker f$$

follows from the equality (2.1.8) and from the definition (2.1.1). Therefore, we proved the statement

(2.1.9)              $$(a, b), (b, c) \in \ker f \Rightarrow (a, c) \in \ker f$$

The lemma follows from the statement (2.1.9) and from the definition 2.1.1.2.   ⊙

The theorem follows from lemmas 2.1.3, 2.1.4, 2.1.5 and from the definition 2.1.1.                                                                          □

THEOREM 2.1.6. *Let $N$ be equivalence on the set $A$. Consider category $\mathcal{A}$ whose objects are maps* [2.3]

$$f_1 : A \to S_1 \quad \ker f_1 \supseteq N$$

$$f_2 : A \to S_2 \quad \ker f_2 \supseteq N$$

---

[2.3] The statement of lemma is similar to the statement on p. [2]-119.

*We define morphism $f_1 \to f_2$ to be map $h : S_1 \to S_2$ making following diagram commutative*

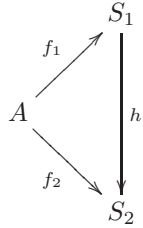

*The map*

$$\operatorname{nat} N : A \to A/N$$

*is universally repelling in the category $\mathcal{A}$.* [2.4]

PROOF. Consider diagram

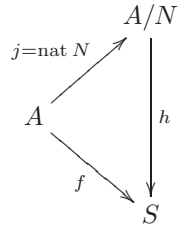

(2.1.10)                              $\ker f \supseteq N$

From the statement (2.1.10) and the equality

$$j(a_1) = j(a_2)$$

it follows that

$$f(a_1) = f(a_2)$$

Therefore, we can uniquely define the map $h$ using the equality

$$h(j(b)) = f(b)$$

$\square$

## 2.2. Universal Algebra

DEFINITION 2.2.1. *For any sets* [2.5] *$A$, $B$, **Cartesian power** $B^A$ is the set of maps*

$$f : A \to B$$

$\square$

DEFINITION 2.2.2. *For any $n \geq 0$, a map* [2.6]

$$\omega : A^n \to A$$

*is called $n$-**ary operation on set** $A$ or just **operation on set** $A$. For any $a_1$, ..., $a_n \in A$, we use either notation $\omega(a_1, ..., a_n)$, $a_1...a_n\omega$ to denote image of map $\omega$.* $\square$

---

REMARK 2.2.3. *According to definitions 2.2.1, 2.2.2, n-ari operation $\omega \in A^{A^n}$.*

$\square$

DEFINITION 2.2.4. *An* **operator domain** *is the set of operators*[2.7] $\Omega$ *with a map*

$$a : \Omega \to N$$

*If $\omega \in \Omega$, then $a(\omega)$ is called the* **arity** *of operator $\omega$. If $a(\omega) = n$, then operator $\omega$ is called n-ary. We use notation*

$$\Omega(n) = \{\omega \in \Omega : a(\omega) = n\}$$

*for the set of n-ary operators.*

$\square$

DEFINITION 2.2.5. *Let $A$ be a set. Let $\Omega$ be an operator domain.*[2.8] *The family of maps*

$$\Omega(n) \to A^{A^n} \quad n \in N$$

*is called $\Omega$-algebra structure on $A$. The set $A$ with $\Omega$-algebra structure is called $\Omega$-**algebra** $A_\Omega$ or* **universal algebra**. *The set $A$ is called* **carrier of** $\Omega$-algebra.

$\square$

The operator domain $\Omega$ describes a set of $\Omega$-algebras. An element of the set $\Omega$ is called operator, because an operation assumes certain set. According to the remark 2.2.3 and the definition 2.2.5, for each operator $\omega \in \Omega(n)$, we match $n$-ary operation $\omega$ on $A$.

THEOREM 2.2.6. *Let the set $B$ be $\Omega$-algebra. Then the set $B^A$ of maps*

$$f : A \to B$$

*also is $\Omega$-algebra.*

PROOF. Let $\omega \in \Omega(n)$. For maps $f_1, ..., f_n \in B^A$, we define the operation $\omega$ by the equality

$$(f_1...f_n\omega)(x) = f_1(x)...f_n(x)\omega$$

$\square$

DEFINITION 2.2.7. *Let $B \subseteq A$. Since, for any $b_1, ..., b_n \in B$, $b_1...b_n\omega \in B$, then we say that $B$* **is closed with respect to** $\omega$ *or that $B$* **admits operation** $\omega$.

$\square$

DEFINITION 2.2.8. *$\Omega$-algebra $B_\Omega$ is* **subalgebra** *of $\Omega$-algebra $A_\Omega$ if following statements are true*[2.9]

2.2.8.1: $B \subseteq A$.

2.2.8.2: *if operator $\omega \in \Omega$ defines operations $\omega_A$ on $A$ and $\omega_B$ on $B$, then*

$$\omega_A|B = \omega_B$$

$\square$

---

[2.7] I follow the definition (1), page [14]-48.

[2.8] I follow the definition (2), page [14]-48.

[2.9] I follow the definition on page [14]-48.

Definition 2.2.9. *Let $A$, $B$ be $\Omega$-algebras and $\omega \in \Omega(n)$. The map* [2.10]

$$f : A \to B$$

**is compatible with operation** $\omega$, *if, for all $a_1, ..., a_n \in A$,*

(2.2.1) $$f(a_1)...f(a_n)\omega = f(a_1...a_n\omega)$$

*The map $f$ is called* **homomorphism** *from $\Omega$-algebra $A$ to $\Omega$-algebra $B$, if $f$ is compatible with each $\omega \in \Omega$. We use notation* $\mathrm{Hom}(\Omega; A \to B)$ *for the set of homomorphisms from $\Omega$-algebra $A$ to $\Omega$-algebra $B$.*                    □

Theorem 2.2.10. *Since operator domain is empty, then a homomorphism from $\Omega$-algebra $A$ to $\Omega$-algebra $B$ is a map*

$$f : A \to B$$

*Therefore, $\mathrm{Hom}(\emptyset; A \to B) = B^A$.*

Proof. The theorem follows from definitions 2.2.1, 2.2.9.                    □

Definition 2.2.11. *Homomorphism $f$ is called* [2.11] **isomorphism** *between $A$ and $B$, if correspondence $f^{-1}$ is homomorphism. If there is an isomorphism between $A$ and $B$, then we say that $A$ and $B$ are isomorphic and write* $A \cong B$. *An injective homomorphis is called* **monomorphism**. *A surjective homomorphis is called* **epimorphism**.                    □

Definition 2.2.12. *A homomorphism in which source and target are the same algebra is called* **endomorphism**. *We use notation* $\mathrm{End}(\Omega; A)$ *for the set of endomorphisms of $\Omega$-algebra $A$. An endomorphism which is also an isomorphism is called* **automorphism**.                    □

Theorem 2.2.13. $\mathrm{End}(\Omega; A) = \mathrm{Hom}(\Omega; A \to A)$

Proof. The theorem follows from the definitions 2.2.9, 2.2.12.                    □

Theorem 2.2.14. *Since operator domain is empty, then an endomorphism of the set $A$ is a map*

$$t : A \to A$$

*Therefore, $\mathrm{End}(\emptyset; A) = A^A$.*

Proof. The theorem follows from the theorems 2.2.10, 2.2.13.                    □

Definition 2.2.15. *If there is a monomorphism from $\Omega$-algebra $A$ to $\Omega$-algebra $B$, then we say that $A$* **can be embeded in** $B$.                    □

Definition 2.2.16. *If there is an epimorphism from $A$ to $B$, then $B$ is called* **homomorphic image** *of algebra $A$.*                    □

---

[2.10] I follow the definition on page [14]-49.

[2.11] I follow the definition on page [14]-49.

### 2.3. Cartesian Product of Universal Algebras

DEFINITION 2.3.1. *Let $\mathcal{A}$ be a category. Let $\{B_i, i \in I\}$ be the set of objects of $\mathcal{A}$. Object*

$$P = \prod_{i \in I} B_i$$

*and set of morphisms*

$$\{f_i : P \to B_i, i \in I\}$$

*is called a* **product of set of objects** $\{B_i, i \in I\}$ **in category** $\mathcal{A}$ [2.12] *if for any object $R$ and set of morphisms*

$$\{g_i : R \to B_i, i \in I\}$$

*there exists a unique morphism*

$$h : R \to P$$

*such that diagram*

$$
\begin{array}{ccc}
P & \xrightarrow{f_i} & B_i \\
\uparrow{\scriptstyle h} & \nearrow{\scriptstyle g_i} & \\
R & &
\end{array}
\qquad f_i \circ h = g_i
$$

*is commutative for all $i \in I$.*

*If $|I| = n$, then we also will use notation*

$$P = \prod_{i=1}^{n} B_i = B_1 \times ... \times B_n$$

*for product of set of objects $\{B_i, i \in I\}$ in $\mathcal{A}$.*  □

EXAMPLE 2.3.2. *Let $\mathcal{S}$ be the category of sets.* [2.13] *According to the definition 2.3.1, Cartesian product*

$$A = \prod_{i \in I} A_i$$

*of family of sets $(A_i, i \in I)$ and family of projections on the $i$-th factor*

$$p_i : A \to A_i$$

*are product in the category $\mathcal{S}$.*  □

THEOREM 2.3.3. *The product exists in the category $\mathcal{A}$ of $\Omega$-algebras. Let $\Omega$-algebra $A$ and family of morphisms*

$$p_i : A \to A_i \quad i \in I$$

*be product in the category $\mathcal{A}$. Then*

2.3.3.1: *The set $A$ is Cartesian product of family of sets $(A_i, i \in I)$*
2.3.3.2: *The homomorphism of $\Omega$-algebra*

$$p_i : A \to A_i$$

*is projection on $i$-th factor.*
2.3.3.3: *We can represent any $A$-number $a$ as tuple $(p_i(a), i \in I)$ of $A_i$-numbers.*

---

[2.12] I made definition according to [2], page 58.
[2.13] See also the example in [2], page 59.

2.3.3.4: *Let $\omega \in \Omega$ be n-ary operation. Then operation $\omega$ is defined component-wise*

(2.3.1) $$a_1...a_n\omega = (a_{1i}...a_{ni}\omega, i \in I)$$

*where $a_1 = (a_{1i}, i \in I), ..., a_n = (a_{ni}, i \in I)$ .*

Proof. Let

$$A = \prod_{i \in I} A_i$$

be Cartesian product of family of sets $(A_i, i \in I)$ and, for each $i \in I$, the map

$$p_i : A \to A_i$$

be projection on the $i$-th factor. Consider the diagram of morphisms in category of sets $\mathcal{S}$

(2.3.2) $$A \xrightarrow{p_i} A_i \qquad p_i \circ \omega = g_i$$

where the map $g_i$ is defined by the equality

$$g_i(a_1, ..., a_n) = p_i(a_1)...p_i(a_n)\omega$$

According to the definition 2.3.1, the map $\omega$ is defined uniquely from the set of diagrams (2.3.2)

(2.3.3) $$a_1...a_n\omega = (p_i(a_1)...p_i(a_n)\omega, i \in I)$$

The equality (2.3.1) follows from the equality (2.3.3). ☐

Definition 2.3.4. *If $\Omega$-algebra $A$ and family of morphisms*

$$p_i : A \to A_i \quad i \in I$$

*is product in the category $\mathcal{A}$, then $\Omega$-algebra $A$ is called* **direct** *or* **Cartesian product of $\Omega$-algebras** $(A_i, i \in I)$ . ☐

Theorem 2.3.5. *Let set $A$ be Cartesian product of sets $(A_i, i \in I)$ and set $B$ be Cartesian product of sets $(B_i, i \in I)$ . For each $i \in I$, let*

$$f_i : A_i \to B_i$$

*be the map from the set $A_i$ into the set $B_i$. For each $i \in I$, consider commutative diagram*

(2.3.4) $$B \xrightarrow{p_i'} B_i$$

*where maps $p_i$, $p_i'$ are projection on the $i$-th factor. The set of commutative diagrams (2.3.4) uniquely defines map*

$$f : A \to B$$

$$f(a_i, i \in I) = (f_i(a_i), i \in I)$$

Proof. For each $i \in I$, consider commutative diagram

(2.3.5)

Let $a \in A$. According to the statement 2.3.3.3, we can represent $A$-number $a$ as tuple of $A_i$-numbers

(2.3.6) $$a = (a_i, i \in I) \quad a_i = p_i(a) \in A_i$$

Let

(2.3.7) $$b = f(a) \in B$$

According to the statement 2.3.3.3, we can represent $B$-number $b$ as tuple of $B_i$-numbers

(2.3.8) $$b = (b_i, i \in I) \quad b_i = p'_i(b) \in B_i$$

From commutativity of diagram (1) and from equalities (2.3.7), (2.3.8), it follows that

(2.3.9) $$b_i = g_i(b)$$

From commutativity of diagram (2) and from the equality (2.3.6), it follows that

$$b_i = f_i(a_i)$$

$\square$

Theorem 2.3.6. *Let $\Omega$-algebra $A$ be Cartesian product of $\Omega$-algebras $(A_i, i \in I)$ and $\Omega$-algebra $B$ be Cartesian product of $\Omega$-algebras $(B_i, i \in I)$. For each $i \in I$, let the map*

$$f_i : A_i \to B_i$$

*be homomorphism of $\Omega$-algebra. Then the map*

$$f : A \to B$$

*defined by the equality*

(2.3.10) $$f(a_i, i \in I) = (f_i(a_i), i \in I)$$

*is homomorphism of $\Omega$-algebra.*

Proof. Let $\omega \in \Omega$ be n-ary operation. Let $a_1 = (a_{1i}, i \in I)$, ..., $a_n = (a_{ni}, i \in I)$, $b_1 = (b_{1i}, i \in I)$, ..., $b_n = (b_{ni}, i \in I)$. From equalities (2.3.1), (2.3.10), it follows that

$$f(a_1...a_n\omega) = f(a_{1i}...a_{ni}\omega, i \in I)$$
$$= (f_i(a_{1i}...a_{ni}\omega), i \in I)$$
$$= ((f_i(a_{1i}))...(f_i(a_{ni})), i \in I)$$
$$= (b_{1i}...b_{ni}\omega, i \in I)$$
$$f(a_1)...f(a_n)\omega = b_1...b_n\omega = (b_{1i}...b_{ni}\omega, i \in I)$$

□

Definition 2.3.7. *Equivalence on $\Omega$-algebra $A$, which is subalgebra of $\Omega$-algebra $A^2$, is called* **congruence** [2.14] *on $A$.*        □

Theorem 2.3.8 (first isomorphism theorem). *Let*

$$f : A \to B$$

*be homomorphism of $\Omega$-algebras with kernel $s$. Then there is decomposition*

$$A/\ker f \xrightarrow{\ q\ } f(A) \qquad f = p \circ q \circ r$$

$$
\begin{array}{ccc}
A/\ker f & \xrightarrow{\ q\ } & f(A) \\
{\scriptstyle p}\big\uparrow & & \big\downarrow{\scriptstyle r} \\
A & \xrightarrow{\ f\ } & B
\end{array}
$$

2.3.8.1: *The* **kernel of homomorphism** $\ker f = f \circ f^{-1}$ *is a congruence on $\Omega$-algebra $A$.*

2.3.8.2: *The set $A/\ker f$ is $\Omega$-algebra.*

2.3.8.3: *The map*

$$p : a \in A \to a^{\ker f} \in A/\ker f$$

*is epimorphism and is called* **natural homomorphism**.

2.3.8.4: *The map*

$$q : p(a) \in A/\ker f \to f(a) \in f(A)$$

*is the isomorphism*

2.3.8.5: *The map*

$$r : f(a) \in f(A) \to f(a) \in B$$

*is the monomorphism*

Proof. The statement 2.3.8.1 follows from the proposition II.3.4 ([14], page 58). Statements 2.3.8.2, 2.3.8.3 follow from the theorem II.3.5 ([14], page 58) and from the following definition. Statements 2.3.8.4, 2.3.8.5 follow from the theorem II.3.7 ([14], page 60).        □

## 2.4. Semigroup

Usually the operation $\omega \in \Omega(2)$ is called product

$$ab\omega = ab$$

or sum

$$ab\omega = a + b$$

Definition 2.4.1. *Let $A$ be $\Omega$-algebra and $\omega \in \Omega(2)$. A-number $e$ is called* **neutral element of operation** $\omega$, *when for any A-number $a$ following equations are true*

$$(2.4.1) \qquad\qquad e a \omega = a$$

$$(2.4.2) \qquad\qquad a e \omega = a$$

□

---

[2.14] I follow the definition on page [14]-57.

Definition 2.4.2. *Let $A$ be $\Omega$-algebra. The operation $\omega \in \Omega(2)$ is called* **associative** *if the following equality is true*

$$a(bc\omega)\omega = (ab\omega)c\omega$$

<div align="right">□</div>

Definition 2.4.3. *Let $A$ be $\Omega$-algebra. The operation $\omega \in \Omega(2)$ is called* **commutative** *if the following equality is true*

$$ab\omega = ba\omega$$

<div align="right">□</div>

Definition 2.4.4. *Let $\Omega = \{\omega\}$. If the operation $\omega \in \Omega(2)$ is associative, then $\Omega$-algebra is called* **semigroup***. If the operation in the semigroup is commutative, then the semigroup is called* **Abelian semigroup***.* □

# Representation of Universal Algebra

## 3.1. Representation of Universal Algebra

DEFINITION 3.1.1. *Let the set $A_2$ be $\Omega_2$-algebra. Let the set of transformations* $\text{End}(\Omega_2, A_2)$ *be $\Omega_1$-algebra. The homomorphism*

$$f : A_1 \to \text{End}(\Omega_2; A_2)$$

*of $\Omega_1$-algebra $A_1$ into $\Omega_1$-algebra* $\text{End}(\Omega_2, A_2)$ *is called* **representation of $\Omega_1$-algebra $A_1$ or $A_1$-representation** *in $\Omega_2$-algebra $A_2$.* □

Diagram

$$
\begin{array}{ccc}
A_2 & \xrightarrow{\ f(a)\ } & A_2 \\
& \Big\Uparrow{\scriptstyle f} & \\
& A_1 &
\end{array}
$$

means that we consider the representation of $\Omega_1$-algebra $A_1$. The map $f(a)$ is image of $a \in A_1$. We also use record

$$f : A_1 \relbar\joinrel\twoheadrightarrow A_2$$

to denote the representation of $\Omega_1$-algebra $A_1$ in $\Omega_2$-algebra $A_2$.

There are several ways to describe the representation. We can define the map $f$ keeping in mind that the domain is $\Omega_1$-algebra $A_1$ and range is $\Omega_1$-algebra $\text{End}(\Omega_2, A_2)$. Either we can specify $\Omega_1$-algebra $A_1$ and $\Omega_2$-algebra $A_2$ keeping in mind that we know the structure of the map $f$.[3.1]

DEFINITION 3.1.2. *Let the map*

$$f : A_1 \to \text{End}(\Omega_2; A_2)$$

*be an isomorphism of the $\Omega_1$-algebra $A_1$ into* $\text{End}(\Omega_2, A_2)$. *Then the representation*

$$f : A_1 \relbar\joinrel\twoheadrightarrow A_2$$

*of the $\Omega_1$-algebra $A_1$ is called* **effective**.[3.2] □

THEOREM 3.1.3. *The representation*

$$f : A_1 \relbar\joinrel\twoheadrightarrow A_2$$

---

[3.1] For instance, we consider vector space $V$ over field $D$ (section 9.3).

[3.2] See similar definition of effective representation of group in [18], page 16, [19], page 111, [15], page 51 (Cohn calls such representation faithful). See also the theorem 5.4.2.

*is effective iff the statement $a_1 \neq b_1$, $a_1$, $b_1 \in A_1$, implies that there exists $a_2 \in A_2$ such that* [3.3]

$$f(a_1)(a_2) \neq f(b_1)(a_2)$$

PROOF. Let the representation $f$ be effective and $a_1 \neq b_1$. If for any $a_2 \in A_2$, the equality

$$f(a_1)(a_2) = f(b_1)(a_2)$$

is true, then

$$f(a_1) = f(b_1)$$

This contradicts to the statement that the representation $f$ is effective.

Let the statement $a_1 \neq b_1$, $a_1$, $b_1 \in A_1$, imply that there exists $a_2 \in A_2$ such that

$$f(a_1)(a_2) \neq f(b_1)(a_2)$$

Therefore, the statement $a_1 \neq b_1$, $a_1$, $b_1 \in A_1$, implies that

$$f(a) \neq f(b)$$

According to the definition 3.1.2, the representation $f$ is effective.            □

DEFINITION 3.1.4. *The representation*

$$f : A_1 \longrightarrow\!\!\!* A_2$$

*of the $\Omega_1$-algebra $A_1$ is called* **free**, [3.4] *if the statement*

$$f(a_1)(a_2) = f(b_1)(a_2)$$

*for any $a_2 \in A_2$ implies that $a_1 = b_1$.*            □

THEOREM 3.1.5. *The representation*

$$f : A_1 \longrightarrow\!\!\!* A_2$$

*of the $\Omega_1$-algebra $A_1$ is called* **free**, *if the statement $f(a_1) = f(b_1)$ implies that $a_1 = b_1$.*

PROOF. The statement $f(a_1) = f(b_1)$ is true iff when

$$f(a_1)(a_2) = f(b_1)(a_2)$$

for any $a_2 \in A_2$.            □

THEOREM 3.1.6. *Free representation is effective.*

---

[3.3] In case of group, the theorem 3.1.3 has the following form. *The representation*

$$f : A_1 \longrightarrow\!\!\!* A_2$$

*is effective iff, for any $A_1$-number $a_1 \neq e$, there exists $a_2 \in A_2$ such that*

$$f(a_1)(a_2) \neq a_2$$

[3.4] See similar definition of free representation of group in [18], page 16. See also the theorem 5.5.2.

Proof. Let the map

$$f : A_1 \longrightarrow_* A_2$$

be free representation. Let $a$, $b \in A_1$. According to the definition 3.1.4, the statement

$$f(a_1)(a_2) = f(b_1)(a_2)$$

for any $a_2 \in A_2$ implies that $a_1 = b_1$. Therefore, if $a_1 \neq b_1$, then there exists $a_2 \in A_2$ such that

$$f(a_1)(a_2) \neq f(b_1)(a_2)$$

According to the theorem 3.1.3, the representation $f$ is effective. $\square$

Remark 3.1.7. *Representation of rotation group in affine space is effective. However this representation is not free, since origin is fixed point of every transformation.* $\square$

Definition 3.1.8. *The representation*

$$f : A_1 \longrightarrow_* A_2$$

*of $\Omega_1$-algebra is called* **transitive** [3.5] *if for any $a$, $b \in A_2$, exists such $g$ that*

$$a = f(g)(b)$$

*The representation of $\Omega_1$-algebra is called* **single transitive** *if it is transitive and free.* $\square$

Theorem 3.1.9. *Representation is single transitive iff for any $a, b \in A_2$ exists one and only one $g \in A_1$ such that $a = f(g)(b)$*

Proof. The theorem follows from definitions 3.1.4 and 3.1.8. $\square$

Theorem 3.1.10. *Let*

$$f : A_1 \longrightarrow_* A_2$$

*be a single transitive representation of $\Omega_1$-algebra $A_1$ in $\Omega_2$-algebra $A_2$. There is the structure of $\Omega_1$-algebra on the set $A_2$.*

Proof. Let $b \in A_2$, $\omega \in \Omega_1(n)$. For any $A_2$-numbers $b_1, ..., b_n$, there exist $A_1$-numbers $a_1, ..., a_n$ such that

$$b_1 = f(a_1)(b) \quad ... \quad b_n = f(a_n)(b)$$

We introduce the operation $\omega$ on the set $A_2$ by the equality

(3.1.1) $$b_1...b_n\omega = f(a_1...a_n\omega)(b)$$

We also require that choice of $A_2$-number $b$ does not depend on operation $\omega$. $\square$

Theorem 3.1.11. *Let*

$$f : A_1 \longrightarrow_* A_2$$

*be an effective representation of $\Omega_1$-algebra $A_1$ in $\Omega_2$-algebra $A_2$. Let $\omega \in \Omega_1(n) \cap \Omega_2(n)$. Then*

(3.1.2) $$f(a_1...a_n\omega)(b) = f(a_1)(b)...f(a_n)(b)\omega$$

---

[3.5] See similar definition of transitive representation of group in [19], page 110, [15], page 51.

## 3.2. Morphism of Representations of Universal Algebra

THEOREM 3.2.1. *Let $A_1$ and $B_1$ be $\Omega_1$-algebras. Representation of $\Omega_1$-algebra $B_1$*

$$g : B_1 \xrightarrow{\;*\;} A_2$$

*and homomorphism of $\Omega_1$-algebra*

$$h : A_1 \to B_1$$

*define representation $f$ of $\Omega_1$-algebra $A_1$*

(3.2.1)

$$
\begin{array}{ccc}
A_1 & \xrightarrow{\;\;f\;\;} & \mathrm{End}(\Omega_2; A_2) \\
 & \searrow{\scriptstyle h} \quad \nearrow{\scriptstyle g} & \\
 & B_1 &
\end{array}
$$

PROOF. Since map $g$ is homomorphism of $\Omega_1$-algebra $B_1$ into $\Omega_1$-algebra $\mathrm{End}(\Omega_2, A_2)$, the map $f$ is homomorphism of $\Omega_1$-algebra $A_1$ into $\Omega_1$-algebra $\mathrm{End}(\Omega_2, A_2)$. ☐

We also use diagram

$$
\begin{array}{ccc}
A_1 & \xrightarrow{\;\;f\;\;} & A_2 \\
 & \searrow{\scriptstyle h} \quad \nearrow{\scriptstyle g} & \\
 & B_1 &
\end{array}
$$

instead of diagram (3.2.1).

Considering representations of $\Omega_1$-algebra in $\Omega_2$-algebras $A_2$ and $B_2$, we are interested in a map $A_2 \to B_2$ that preserves the structure of representation.

DEFINITION 3.2.2. *Let*

$$f : A_1 \xrightarrow{\;*\;} A_2$$

*be representation of $\Omega_1$-algebra $A_1$ in $\Omega_2$-algebra $A_2$ and*

$$g : B_1 \xrightarrow{\;*\;} B_2$$

*be representation of $\Omega_1$-algebra $B_1$ in $\Omega_2$-algebra $B_2$. For $i = 1, 2$, let the map*

$$r_i : A_i \to B_i$$

*be homomorphism of $\Omega_i$-algebra. The tuple of maps $r = (r_1, r_2)$ such, that*

(3.2.2)

$$r_2 \circ f(a) = g(r_1(a)) \circ r_2$$

*is called* **morphism of representations from $f$ into $g$.** *We also say that* **morphism of representations of $\Omega_1$-algebra in $\Omega_2$-algebra** *is defined.* ☐

REMARK 3.2.3. *We may consider a pair of maps $r_1$, $r_2$ as map*

$$F : A_1 \cup A_2 \to B_1 \cup B_2$$

*such that*

$$F(A_1) = B_1 \qquad F(A_2) = B_2$$

*Therefore, hereinafter the tuple of maps $r = (r_1, r_2)$ also is called map and we will use map*

$$r : f \to g$$

Let $a = (a_1, a_2)$ be tuple of $A$-numbers. We will use notation

$$r(a) = (r_1(a_1), r_2(a_2))$$

for image of tuple of $A$-numbers with respect to morphism of representations $r$.  $\square$

**Definition 3.2.4.** *If representation $f$ and $g$ coincide, then morphism of representations $r = (r_1, r_2)$ is called* **morphism of representation $f$**.  $\square$

**Theorem 3.2.5.** *Let*

$$f : A_1 \longrightarrow\!\!\!* \longrightarrow A_2$$

*be representation of $\Omega_1$-algebra $A_1$ in $\Omega_2$-algebra $A_2$ and*

$$g : B_1 \longrightarrow\!\!\!* \longrightarrow B_2$$

*be representation of $\Omega_1$-algebra $B_1$ in $\Omega_2$-algebra $B_2$. The map*

$$(r_1 : A_1 \to B_1, \ r_2 : A_2 \to B_2)$$

*is morphism of representations iff*

$$(3.2.3) \qquad r_2(f(a)(m)) = g(r_1(a))(r_2(m))$$

**Proof.** For any $m \in A_2$, equality (3.2.3) follows from (3.2.2).  $\square$

**Remark 3.2.6.** *Consider morphism of representations*

$$(r_1 : A_1 \to B_1, \ r_2 : A_2 \to B_2)$$

*We denote elements of the set $B_1$ by letter using pattern $b \in B_1$. However if we want to show that $b$ is image of element $a \in A_1$, we use notation $r_1(a)$. Thus equation*

$$r_1(a) = r_1(a)$$

*means that $r_1(a)$ (in left part of equation) is image $a \in A_1$ (in right part of equation). Using such considerations, we denote element of set $B_2$ as $r_2(m)$. We will follow this convention when we consider correspondences between homomorphisms of $\Omega_1$-algebra and maps between sets where we defined corresponding representations.*  $\square$

**Remark 3.2.7.** *There are two ways to interpret (3.2.3)*

- *Let transformation $f(a)$ map $m \in A_2$ into $f(a)(m)$. Then transformation $g(r_1(a))$ maps $r_2(m) \in B_2$ into $r_2(f(a)(m))$.*
- *We represent morphism of representations from $f$ into $g$ using diagram*

$$(3.2.4)$$

*From (3.2.2), it follows that diagram (1) is commutative.*

*We also use diagram*

(3.2.5)

$$\begin{array}{ccc} A_2 & \xrightarrow{\ r_2\ } & B_2 \\ {\scriptstyle f}\big\uparrow & & \big\uparrow{\scriptstyle g} \\ A_1 & \xrightarrow{\ r_1\ } & B_1 \end{array}$$

*instead of diagram* (3.2.4).                                                                  □

Theorem 3.2.8. *Consider representation*

$$f : A_1 \xrightarrow{\quad*\quad} A_2$$

*of $\Omega_1$-algebra $A_1$ and representation*

$$g : B_1 \xrightarrow{\quad*\quad} B_2$$

*of $\Omega_1$-algebra $B_1$. Morphism*

$$(r_1 : A_1 \to B_1,\ r_2 : A_2 \to B_2)$$

*of representations from $f$ into $g$ satisfies equation*

(3.2.6)          $r_2 \circ (f(a_1)...f(a_n)\omega) = (g(r_1(a_1))...g(r_1(a_n))\omega) \circ r_2$

*for any operation $\omega \in \Omega_1(n)$.*

Proof. Since $f$ is homomorphism, we have

(3.2.7)          $r_2 \circ (f(a_1)...f(a_n)\omega) = r_2 \circ f(a_1...a_n\omega)$

From (3.2.2) and (3.2.7) it follows that

(3.2.8)          $r_2 \circ (f(a_1)...f(a_n)\omega) = g(r_1(a_1...a_n\omega)) \circ r_2$

Since $r_1$ is homomorphism, from (3.2.8) it follows that

(3.2.9)          $r_2 \circ (f(a_1)...f(a_n)\omega) = g(r_1(a_1)...r_1(a_n)\omega) \circ r_2$

Since $g$ is homomorphism, (3.2.6) follows from (3.2.9).                          □

Theorem 3.2.9. *Let the map*

$$(r_1 : A_1 \to B_1,\ r_2 : A_2 \to B_2)$$

*be morphism from representation*

$$f : A_1 \xrightarrow{\quad*\quad} A_2$$

*of $\Omega_1$-algebra $A_1$ into representation*

$$g : B_1 \xrightarrow{\quad*\quad} B_2$$

*of $\Omega_1$-algebra $B_1$. If representation $f$ is effective, then the map*

$$r_2^* : \operatorname{End}(\Omega_2; A_2) \to \operatorname{End}(\Omega_2; B_2)$$

*defined by equation*

(3.2.10)          $r_2^*(f(a)) = g(r_1(a))$

*is homomorphism of $\Omega_1$-algebra.*

Proof. Because representation $f$ is effective, then for given transformation $f(a)$ element $a$ is determined uniquely. Therefore, transformation $g(r_1(a))$ is properly defined in equation (3.2.10).

Since $f$ is homomorphism, we have

$$(3.2.11) \qquad r_2^*(f(a_1)...f(a_n)\omega) = r_2^*(f(a_1...a_n\omega))$$

From (3.2.10) and (3.2.11) it follows that

$$(3.2.12) \qquad r_2^*(f(a_1)...f(a_n)\omega) = g(r_1(a_1...a_n\omega))$$

Since $h$ is homomorphism, from (3.2.12) it follows that

$$(3.2.13) \qquad r_2^*(f(a_1)...f(a_n)\omega) = g(r_1(a_1)...r_1(a_n)\omega)$$

Since $g$ is homomorphism,

$$r_2^*(f(a_1)...f(a_n)\omega) = g(r_1(a_1))...g(r_1(a_n))\omega = r_2^*(f(a_1))...r_2^*(f(a_n))\omega$$

follows from (3.2.13). Therefore, the map $r_2^*$ is homomorphism of $\Omega_1$-algebra. □

Theorem 3.2.10. Let

$$f : A_1 \longrightarrow\!\!* A_2$$

be single transitive representation of $\Omega_1$-algebra $A_1$ and

$$g : B_1 \longrightarrow\!\!* B_2$$

be single transitive representation of $\Omega_1$-algebra $B_1$. Given homomorphism of $\Omega_1$ algebra

$$r_1 : A_1 \to B_1$$

there exists morphism of representations from $f$ into $g$

$$(r_1 : A_1 \to B_1, \ r_2 : A_2 \to B_2)$$

Proof. Let us choose homomorphism $r_1$. Let us choose element $m \in A_2$ and element $n \in B_2$. To define map $r_2$, consider following diagram

From commutativity of diagram (1), it follows that

$$r_2(f(a)(m)) = g(r_1(a))(r_2(m))$$

For arbitrary $m' \in A_2$, we defined unambiguously $a \in A_1$ such that $m' = f(a)(m)$. Therefore, we defined map $r_2$ which satisfies to equation (3.2.2). □

Theorem 3.2.11. *Let the representation*

$$f : A_1 \xrightarrow{\quad * \quad} A_2$$

*of $\Omega_1$-algebra $A_1$ be single transitive representation and the representation*

$$g : B_1 \xrightarrow{\quad * \quad} B_2$$

*of $\Omega_1$-algebra $B_1$ be single transitive representation. Given homomorphism of $\Omega_1$-algebra*

$$r_1 : A_1 \to B_1$$

*consider a homomorphism of $\Omega_2$-algebra*

$$r_2 : A_2 \to B_2$$

*such that $r = (r_1, r_2)$ is morphism of representations from $f$ into $g$. The map $H$ is unique up to choice of image $n = r_2(m) \in B_2$ of given element $m \in A_2$.*

Proof. From proof of theorem 3.2.10, it follows that choice of homomorphism $r_1$ and elements $m \in A_2$, $n \in B_2$ uniquely defines the map $r_2$. □

Theorem 3.2.12. *Given single transitive representation*

$$f : A_1 \xrightarrow{\quad * \quad} A_2$$

*of $\Omega_1$-algebra $A_1$, for any endomorphism $r_1 \in \text{End}(\Omega_1; A_1)$ there exists morphism of representation $f$*

$$(r_1 : A_1 \to A_1, \ r_2 : A_2 \to A_2)$$

Proof. Consider following diagram

Statement of theorem is corollary of the theorem 3.2.10. □

## 3.3. Decomposition Theorem for Morphisms of Representations

Theorem 3.3.1. *Let*

$$f : A_1 \xrightarrow{\quad * \quad} A_2$$

*be representation of $\Omega_1$-algebra $A_1$,*

$$g : B_1 \xrightarrow{\quad * \quad} B_2$$

*be representation of $\Omega_1$-algebra $B_1$,*

$$h : C_1 \xrightarrow{\quad * \quad} C_2$$

*be representation of $\Omega_1$-algebra $C_1$. Given morphisms of representations of $\Omega_1$-algebra*

$$(p_1 : A_1 \to B_1, \ p_2 : A_2 \to B_2)$$

$$(q_1 : B_1 \to C_1, \ q_2 : B_2 \to C_2)$$

*There exists morphism of representations of $\Omega_1$-algebra*

$$(r_1 : A_1 \to C_1, \ r_2 : A_2 \to C_2)$$

*where $r_1 = q_1 \circ p_1$, $r_2 = q_2 \circ p_2$. We call morphism $r = (r_1, r_2)$ of representations from $f$ into $h$* **product of morphisms** $p = (p_1, p_2)$ **and** $q = (q_1, q_2)$ **of representations of universal algebra.**

PROOF. We represent statement of theorem using diagram

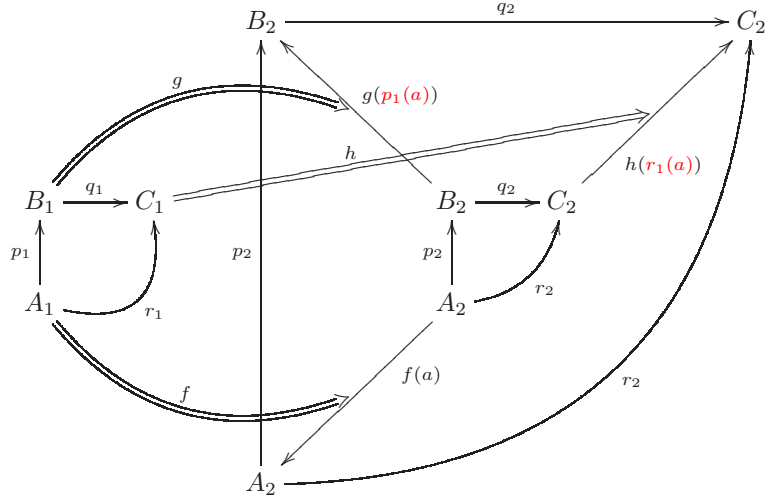

Map $r_1$ is homomorphism of $\Omega_1$-algebra $A_1$ into $\Omega_1$-algebra $C$. We need to show that the map $r = (r_1, r_2)$ satisfies to (3.2.2)

$$r_2(f(a)(m)) = (q_2 \circ p_2)(f(a)(m))$$
$$= q_2(g(p_1(a))(p_2(m)))$$
$$= h((q_1 \circ p_1)(a))((q_2 \circ p_2)(m)))$$
$$= h(r_1(a))(r_2(m))$$

$\square$

DEFINITION 3.3.2. *Let $S$ be equivalence on the set $A_2$. Transformation $f$ is called* **coordinated with equivalence** *$S$, when $f(m_1) \equiv f(m_2)(\mathrm{mod}\ S)$ follows from condition $m_1 \equiv m_2(\mathrm{mod}S)$.* $\square$

THEOREM 3.3.3. *Consider equivalence $S$ on set $A_2$. Consider $\Omega_1$-algebra on the set $\mathrm{End}(\Omega_2, A_2)$. If any transformation $f \in \mathrm{End}(\Omega_2; A_2)$ is coordinated with equivalence $S$, then we can define the structure of $\Omega_1$-algebra on the set $\mathrm{End}(\Omega_2; A_2/S)$.*

PROOF. Let $h = \mathrm{nat}\ S$. If $m_1 \equiv m_2(\mathrm{mod}S)$, then $h(m_1) = h(m_2)$. Since $f \in \mathrm{End}(\Omega_2; A_2)$ is coordinated with equivalence $S$, then $h(f(m_1)) = h(f(m_2))$. This allows us to define transformation $F$ according to rule

$$(3.3.1) \qquad\qquad F([m]) = h(f(m))$$

Let $\omega$ be n-ary operation of $\Omega_1$-algebra. Let $f_1, ..., f_n \in \mathrm{End}(\Omega_2; A_2)$ and

$$F_1([m]) = h(f_1(m)) \quad ... \quad F_n([m]) = h(f_n(m))$$

According to condition of theorem, the transformation

$$f = f_1...f_n\omega \in \operatorname{End}(\Omega_2; A_2)$$

is coordinated with equivalence $S$. Therefore,

(3.3.2)
$$f(m_1) \equiv f(m_2)(\operatorname{mod}S)$$

$$(f_1...f_n\omega)(m_1) \equiv (f_1...f_n\omega)(m_2)(\operatorname{mod}S)$$

follows from condition $m_1 \equiv m_2(\operatorname{mod}S)$ and the definition 3.3.2. Therefore, we can define operation $\omega$ on the set $\operatorname{End}(\Omega_2; A_2/S)$ according to rule

(3.3.3)
$$(F_1...F_n\omega)([m]) = h((f_1...f_n\omega)(m))$$

From the definition (3.3.1) and equation (3.3.2), it follows that we properly defined operation $\omega$ on the set $\operatorname{End}(\Omega_2; A_2/S)$. □

DEFINITION 3.3.4. *Let*

$$f : A_1 \longrightarrow_* A_2$$

*be representation of $\Omega_1$-algebra $A_1$,*

$$g : B_1 \longrightarrow_* B_2$$

*be representation of $\Omega_1$-algebra $B_1$. Let*

$$(r_1 : A_1 \to B_1, \ r_2 : A_2 \to B_2)$$

*be morphism of representations from $f$ into $g$ such that $r_1$ is isomorphism of $\Omega_1$-algebra and $r_2$ is isomorphism of $\Omega_2$-algebra. Then the map $r = (r_1, r_2)$ is called* **isomorphism of repesentations**. □

THEOREM 3.3.5. *Let*

$$f : A_1 \longrightarrow_* A_2$$

*be representation of $\Omega_1$-algebra $A_1$,*

$$g : B_1 \longrightarrow_* B_2$$

*be representation of $\Omega_1$-algebra $B_1$. Let*

$$(t_1 : A_1 \to B_1, \ t_2 : A_2 \to B_2)$$

*be morphism of representations from $f$ into $g$. Then there exist decompositions of $t_1$ and $t_2$, which we describe using diagram*

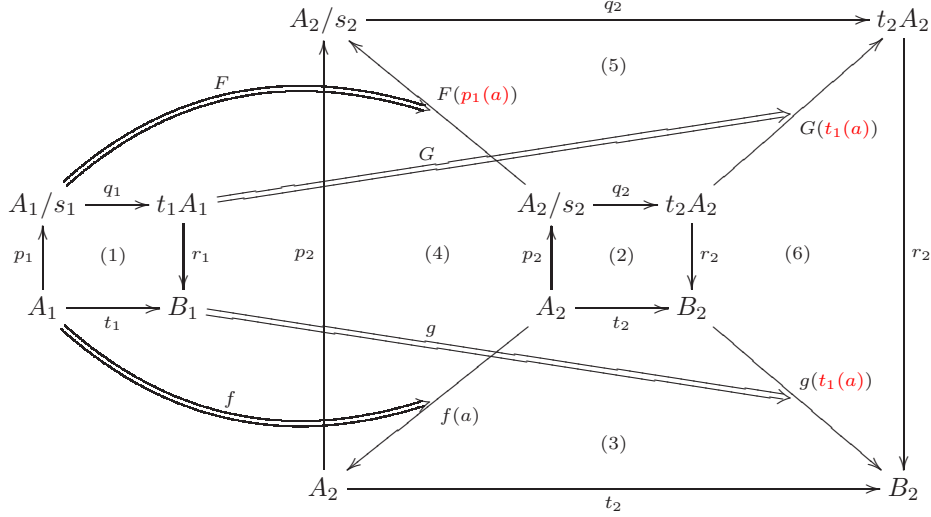

**3.3.5.1:**  *The kernel of homomorphism   $\ker t_i = t_i \circ t_i^{-1}$   is a congruence on $\Omega_i$-algebra $A_i$,  $i = 1$, 2.*

**3.3.5.2:**  *There exists decomposition of homomorphism $t_i$,  $i = 1$, 2,*

$$(3.3.4) \qquad t_i = r_i \circ q_i \circ p_i$$

**3.3.5.3:**  *Maps*

$$p_1(a) = a^{\ker t_1}$$
$$p_2(a) = a^{\ker t_2}$$

*are natural homomorphisms.*

**3.3.5.4:**  *Maps*

$$(3.3.5) \qquad q_1(p_1(a)) = t_1(a)$$

$$(3.3.6) \qquad q_2(p_2(a)) = t_2(a)$$

*are isomorphisms.*

**3.3.5.5:**  *Maps*

$$r_1 : t_1(a) \in f(A_1) \to t_1(a) \in B_1$$
$$r_2 : t_2(a) \in f(A_2) \to t_2(a) \in B_2$$

*are monomorphisms.*

**3.3.5.6:**  *$F$ is representation of $\Omega_1$-algebra $A_1/s$ in $A_2/s_2$*

**3.3.5.7:**  *$G$ is representation of $\Omega_1$-algebra $t_1 A_1$ in $t_2 A_2$*

**3.3.5.8:**  *The map  $p = (p_1, p_2)$  is the morphism of representations $f$ and $F$*

**3.3.5.9:**  *The map  $q = (q_1, q_2)$  is the isomorphism of representations $F$ and $G$*

**3.3.5.10:**  *The map  $r = (r_1, r_2)$  is the morphism of representations $G$ and $g$*

**3.3.5.11:**  *There exists decompositions of morphism of representations*

$$(3.3.7) \qquad (t_1, t_2) = (r_1, r_2) \circ (q_1, q_2) \circ (p_1, p_2)$$

PROOF. Statements 3.3.5.1, 3.3.5.2, 3.3.5.3, 3.3.5.4, 3.3.5.5 follow from the theorem 2.3.8. Therefore, diagrams (1) and (2) are commutative.

We start from diagram (4).

Let $m_1 \equiv m_2 (\mathrm{mod} \ker t_2)$. Then

$$(3.3.8) \qquad t_2(m_1) = t_2(m_2)$$

Since $a_1 \equiv a_2 (\mathrm{mod} \ker t_1)$, then

$$(3.3.9) \qquad t_1(a_1) = t_1(a_2)$$

Therefore, $p_1(a_1) = p_1(a_2)$. Since the map $(t_1, t_2)$ is morphism of representations, then

$$(3.3.10) \qquad t_2(f(a_1)(m_1)) = g(t_1(a_1))(t_2(m_1))$$
$$(3.3.11) \qquad t_2(f(a_2)(m_2)) = g(t_1(a_2))(t_2(m_2))$$

From (3.3.8), (3.3.9), (3.3.10), (3.3.11), it follows that

$$(3.3.12) \qquad t_2(f(a_1)(m_1)) = t_2(f(a_2)(m_2))$$

From (3.3.12) it follows

$$(3.3.13) \qquad f(a_1)(m_1) \equiv f(a_2)(m_2)(\mathrm{mod} \ker t_2)$$

and, therefore,

$$(3.3.14) \qquad p_2(f(a_1)(m_1)) = p_2(f(a_2)(m_2))$$

From (3.3.14) it follows that map

$$(3.3.15) \qquad F(p_1(a))(p_2(m)) = p_2(f(a)(m))$$

is well defined and this map is transformation of set $A_2/\ker t_2$.

From equation (3.3.13) (in case $a_1 = a_2$) it follows that, for any $a$, transformation is coordinated with equivalence $\ker t_2$. From theorem 3.3.3 it follows that we defined structure of $\Omega_1$-algebra on the set $\mathrm{End}(\Omega_2; A_2/\ker t_2)$. Consider $n$-ary operation $\omega$ and $n$ transformations

$$F(p_1(a_i))(p_2(m)) = p_2(f(a_i)(m)) \quad i = 1, ..., n$$

of the $A_2/\ker t_2$. We assume

$$(F(p_1(a_1))...F(p_1(a_n))\omega)(p_2(m)) = p_2((f(a_1)...f(a_n)\omega)(m))$$

Therefore, map $F$ is representations of $\Omega_1$-algebra $A_1/\ker t_1$.

According to the theorem 3.2.5, the statement 3.3.5.8 follows from (3.3.15).

Consider diagram (5).

LEMMA 3.3.6. *The map* $q = (q_1, q_2)$ *is the morphism of representations $F$ and $G$.*

PROOF. Since $q_2$ is bijection, then we identify elements of the set $A_2/\ker t_2$ and the set $t_2(A_2)$, and this identification has form

$$(3.3.16) \qquad q_2(p_2(m)) = t_2(m)$$

We can write transformation $F(p_1(a))$ of the set $A_2/\ker t_2$ as

$$(3.3.17) \qquad F(p_1(a)) : p_2(m) \to F(p_1(a))(p_2(m))$$

Since $q_2$ is bijection, we define transformation

$$(3.3.18) \qquad q_2(p_2(m)) \to q_2(F(p_1(a))(p_2(m)))$$

of the set $t_2(A_2)$. Transformation (3.3.18) depends on $p_1(a) \in A_1/\ker t_1$. Since $q_1$ is bijection, we identify elements of the set $A_1/\ker t_1$ and the set $t_1(A_1)$, and this identification has form

$$(3.3.19) \qquad q_1(p_1(a)) = t_1(a)$$

Therefore, we defined map

$$G : t_1(A_1) \to \text{End}(\Omega_2; t_2(A_2))$$

according to equation

$$(3.3.20) \qquad G(q_1(p_1(a)))(q_2(p_2(m))) = q_2(F(p_1(a))(p_2(m)))$$

Consider $n$-ary operation $\omega$ and $n$ transformations

$$G(t_1(a_i))(t_2(m)) = q_2(F(p_1(a_i))(p_2(m))) \quad i = 1, ..., n$$

of the set $t_2(A_2)$. We assume

$$(3.3.21) \quad (G(t_1(a_1))...G(t_1(a_n))\omega)(t_2(m)) = q_2((F(p_1(a_1))...F(p_1(a_n))\omega)(p_2(m)))$$

According to (3.3.20) operation $\omega$ is well defined on the set $\text{End}(\Omega_2; t_2(A_2))$. Therefore, the map $G$ is representations of $\Omega_1$-algebra.

According to the theorem 3.2.5, the lemma follows from (3.3.20).          ⊙

LEMMA 3.3.7. *The map* $(q_1^{-1}, q_2^{-1})$ *is the morphism of representations $G$ and $F$.*

PROOF. Since $q_2$ is bijection, then from equation (3.3.16), it follows that

$$(3.3.22) \qquad p_2(m) = q_2^{-1}(t_2(m))$$

We can write transformation $G(t_1(a))$ of the set $t_2(A_2)$ as

$$(3.3.23) \qquad G(t_1(a)) : t_2(m) \to G(t_1(a))(t_2(m))$$

Since $q_2$ is bijection, we define transformation

$$(3.3.24) \qquad q_2^{-1}(t_2(m)) \to q_2^{-1}(G(t_1(a))(t_2(m)))$$

of the set $A_2/\ker t_2$. Transformation (3.3.24) depends on $t_1(a) \in t_1(A_1)$. Since $q_1$ is bijection, then from equation (3.3.19) it follows that

$$(3.3.25) \qquad p_1(a) = q_1^{-1}(t_1(a))$$

Since, by construction, diagram (5) is commutative, then transformation (3.3.24) coincides with transformation (3.3.17). We can write the equation (3.3.21) as

$$(3.3.26) \qquad q_2^{-1}((G(t_1(a_1))...G(t_1(a_n))\omega)(t_2(m)))$$
$$= (F(p_1(a_1))...F(p_1(a_n))\omega)(p_2(m))$$

According to the theorem 3.2.5, the lemma follows from (3.3.20), (3.3.22), (3.3.25). ⊙

The statement 3.3.5.9 is corollary of definition 3.3.4 and lemmas 3.3.6 and 3.3.7.

Diagram (6) is the simplest case in our proof. Since map $r_2$ is immersion and diagram (2) is commutative, we identify $n \in B_2$ and $t_2(m)$ when $n \in \text{Im} t_2$. Similarly, we identify corresponding transformations.

$$(3.3.27) \qquad g'(r_1(t_1(a)))(r_2(t_2(m))) = r_2(G(t_1(a))(t_2(m)))$$
$$(g'(t_1(a_1))...g'(t_1(a_n))\omega)(t_2(m)) = r_2((G(t_1(a_1))...G(t_1(a_n))\omega)(t_2(m)))$$

Therefore, $r = (r_1, r_2)$ is morphism of representations $G$ and $g$ (the statement 3.3.5.10).

To prove the statement 3.3.5.11, we need to show that defined in the proof representation $g'$ is congruent with representation $g$, and operations over transformations are congruent with corresponding operations over $\operatorname{End}(\Omega_2, B_2)$.

$$
\begin{aligned}
g'(r_1(t_1(a)))(r_2(t_2(m))) &= r_2(G(t_1(a))(t_2(m))) && \text{by } (3.3.27) \\
&= r_2(G(q_1(p_1(a)))(q_2(p_2(m)))) && \text{by } (3.3.5), (3.3.6), \\
&= r_2 \circ q_2(F(p_1(a))(p_2(m))) && \text{by } (3.3.20) \\
&= r_2 \circ q_2 \circ p_2(f(a)(m)) && \text{by } (3.3.15) \\
&= t_2(f(a)(m)) && \text{by } (3.3.4), \ i = 2 \\
&= g(t_1(a))(t_2(m)) && \text{by } (3.2.2)
\end{aligned}
$$

$$
\begin{aligned}
(G(t_1(a_1))...G(t_1(a_n))\omega)(t_2(m)) &= q_2(F(p_1(a_1))...F(p_1(a_n))\omega)(p_2(m))) \\
&= q_2(F(p_1(a_1)...p_1(a_n)\omega)(p_2(m))) \\
&= q_2(F(p_1(a_1...a_n\omega)(p_2(m))) \\
&= q_2(p_2(f(a_1...a_n\omega)(m)))
\end{aligned}
$$

$\square$

## 3.4. Reduced Morphism of Representations

From theorem 3.3.5, it follows that we can reduce the problem of studying of morphism of representations of $\Omega_1$-algebra to the case described by diagram

(3.4.1)

$$
\begin{array}{ccc}
A_2 & \xrightarrow{p_2} & A_2/\ker t_2 \\
\uparrow{\scriptstyle f} & & \uparrow{\scriptstyle F} \\
A_1 & \xrightarrow{p_1} & A_1/\ker t_1
\end{array}
$$

THEOREM 3.4.1. *We can supplement diagram* (3.4.1) *with representation* $F_1$ *of* $\Omega_1$-*algebra* $A_1$ *into* $\Omega_2$-*algebra* $A_2/\ker t_2$ *such that diagram*

(3.4.2)

$$
\begin{array}{ccc}
A_2 & \xrightarrow{p_2} & A_2/\ker t_2 \\
\uparrow{\scriptstyle f} & \nearrow{\scriptstyle F_1} & \uparrow{\scriptstyle F} \\
A_1 & \xrightarrow{p_1} & A_1/\ker t_1
\end{array}
$$

*is commutative. The set of transformations of representation* $F$ *and the set of transformations of representation* $F_1$ *coincide.*

PROOF. To prove theorem it is enough to assume

$$
F_1(a) = F(p_1(a))
$$

Since map $p_1$ is surjection, then $\operatorname{Im} F_1 = \operatorname{Im} F$. Since $p_1$ and $F$ are homomorphisms of $\Omega_1$-algebra, then $F_1$ is also homomorphism of $\Omega_1$-algebra.          □

Theorem 3.4.1 completes the series of theorems dedicated to the structure of morphism of representations $\Omega_1$-algebra. From these theorems it follows that we can simplify task of studying of morphism of representations $\Omega_1$-algebra and not go beyond morphism of representations of form

$$(\operatorname{id} : A_1 \to A_1, \ r_2 : A_2 \to B_2)$$

Definition 3.4.2. *Let*

$$f : A_1 \longrightarrow_* A_2$$

*be representation of $\Omega_1$-algebra $A_1$ in $\Omega_2$-algebra $A_2$ and*

$$g : A_1 \longrightarrow_* B_2$$

*be representation of $\Omega_1$-algebra $A_1$ in $\Omega_2$-algebra $B_2$. Let*

$$(\operatorname{id} : A_1 \to A_1, \ r_2 : A_2 \to B_2)$$

*be morphism of representations. In this case we identify morphism* $(\operatorname{id}, r_2)$ *of representations of $\Omega_1$-algebra and corresponding homomorphism $r_2$ of $\Omega_2$-algebra and the homomorphism $r_2$ is called* **reduced morphism of representations**. *We will use diagram*

(3.4.3)

*to represent reduced morphism $r_2$ of representations of $\Omega_1$-algebra. From diagram it follows*

(3.4.4)          $$r_2 \circ f(a) = g(a) \circ r_2$$

*We also use diagram*

*instead of diagram* (3.4.3).          □

Theorem 3.4.3. *Let*

$$f : A_1 \longrightarrow_* A_2$$

*be representation of $\Omega_1$-algebra $A_1$ in $\Omega_2$-algebra $A_2$ and*

$$g : A_1 \longrightarrow_* B_2$$

*be representation of $\Omega_1$-algebra $A_1$ in $\Omega_2$-algebra $B_2$. The map*

$$r_2 : A_2 \to B_2$$

*is reduced morphism of representations iff*

(3.4.5) $$r_2(f(a)(m)) = g(a)(r_2(m))$$

Proof. The equality (3.4.5) follows from the equality (3.4.4). □

Theorem 3.4.4. *Let the map*

$$r_2 : A_2 \to B_2$$

*be reduced morphism from representation*

$$f : A_1 \longrightarrow_* A_2$$

*of $\Omega_1$-algebra $A_1$ into representation*

$$g : A_1 \longrightarrow_* B_2$$

*of $\Omega_1$-algebra $A_1$. If representation $f$ is effective, then the map*

$$r_2^* : \mathrm{End}(\Omega_2; A_2) \to \mathrm{End}(\Omega_2; B_2)$$

*defined by equation*

(3.4.6) $$r_2^*(f(a)) = g(a)$$

*is homomorphism of $\Omega_1$-algebra.*

Proof. The theorem follows from the theorem 3.2.9, if we assume $h = id$. □

Theorem 3.4.5. *Let representations*

$$f : A_1 \longrightarrow_* A_2$$

$$g : A_1 \longrightarrow_* B_2$$

*of $\Omega_1$-algebra $A_1$ be single transitive representations. There exists reduced morphism of representations from $f$ into $g$*

$$r_2 : A_2 \to B_2$$

Proof. Let us choose element $m \in A_2$ and element $n \in B_2$. To define map $r_2$, consider following diagram

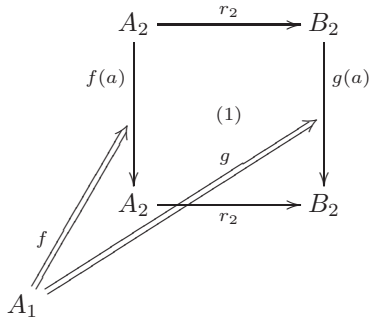

From commutativity of diagram (1), it follows that

$$r_2(f(a)(m)) = g(a)(r_2(m))$$

For arbitrary $m' \in A_2$, we defined unambiguously $a \in A_1$ such that $m' = f(a)(m)$. Therefore, we defined map $r_2$ which satisfies to equation (3.4.4). □

Theorem 3.4.6. *Let representations*

$$f : A_1 \longrightarrow_* \longrightarrow A_2$$

$$g : A_1 \longrightarrow_* \longrightarrow B_2$$

*of $\Omega_1$-algebra $A_1$ be single transitive representations. Reduced morphism of representations from $f$ into $g$*

$$r_2 : A_2 \to B_2$$

*is unique up to choice of image $n = r_2(m) \in B_2$ of given element $m \in A_2$.*

Proof. From proof of theorem 3.4.5, it follows that choice of elements $m \in A_2$, $n \in B_2$ uniquely defines the map $r_2$. □

Theorem 3.4.7. *Let*

$$f : A \longrightarrow_* \longrightarrow B$$

*be representation of $\Omega_1$-algebra $A$ in $\Omega_2$-algebra $B$. Let $N$ be such congruence[3.6] on $\Omega_2$-algebra $B$ that any transformation $h \in \mathrm{End}(\Omega_2, B)$ is coordinated with congruence $N$. There exists representation*

$$f_1 : A \longrightarrow_* \longrightarrow B/N$$

*of $\Omega_1$-algebra $A$ in $\Omega_2$-algebra $B/N$ and the map*

$$\mathrm{nat}\, N : B \to B/N$$

*is reduced morphism of representation $f$ into the representation $f_1$*

$$
\begin{array}{ccc}
B & \xrightarrow{\ \ j\ \ } & B/N \\
\end{array}
\quad j = \mathrm{nat}\, N
$$

Proof. We can represent any element of the set $B/N$ as $j(a)$, $a \in B$.

According to the theorem [14]-II.3.5, there exists a unique $\Omega_2$-algebra structure on the set $B/N$. If $\omega \in \Omega_2(p)$, then we define operation $\omega$ on the set $B/N$ according to the equality (3) on page [14]-59

$$(3.4.7) \qquad j(b_1)...j(b_p)\omega = j(b_1...b_p\omega)$$

As well as in the proof of the theorem 3.3.5, we can define the representation

$$f_1 : A \longrightarrow_* \longrightarrow B/N$$

using equality

$$(3.4.8) \qquad f_1(a) \circ j(b) = j(f(a) \circ b)$$

We can represent the equality (3.4.8) using diagram

$$(3.4.9)
\begin{array}{ccc}
B & \xrightarrow{\ \ j\ \ } & B/N \\
\uparrow{\scriptstyle f(a)} & & \uparrow{\scriptstyle f_1(a)} \\
B & \xrightarrow{\ \ j\ \ } & B/N
\end{array}
$$

─────────

[3.6]See the definition of congruence on p. [14]-57.

Let $\omega \in \Omega_2(p)$. Since the maps $f(a)$ and $j$ are homomorphisms of $\Omega_2$-algebra, then

$$
\begin{aligned}
(3.4.10) \qquad f_1(a) \circ (j(b_1)...j(b_p)\omega) &= f_1(a) \circ j(b_1...b_p\omega) \\
&= j(f(a) \circ (b_1...b_p\omega)) \\
&= j((f(a) \circ b_1)...(f(a) \circ b_p)\omega) \\
&= j(f(a) \circ b_1)...j(f(a) \circ b_p)\omega \\
&= (f_1(a) \circ j(b_1))...(f_1(a) \circ j(b_p))\omega
\end{aligned}
$$

From the equality (3.4.10), it follows that the map $f_1(a)$ is homomorphism of $\Omega_2$-algebra. From the equality (3.4.8), according to the definition 3.4.2, it follows that the map $j$ is reduced morphism of the representation $f$ into the representation $f_1$.                                  □

**Theorem 3.4.8.** *Let*

$$
f : A \longrightarrow\!\!\ast\!\!\longrightarrow B
$$

*be representation of $\Omega_1$-algebra $A$ in $\Omega_2$-algebra $B$. Let $N$ be such congruence on $\Omega_2$-algebra $B$ that any transformation $h \in \mathrm{End}(\Omega_2, B)$ is coordinated with congruence $N$. Consider category $\mathcal{A}$ whose objects are reduced morphisms of representations* [3.7]

$$
R_1 : B \to S_1 \quad \ker R_1 \supseteq N
$$

$$
R_2 : B \to S_2 \quad \ker R_2 \supseteq N
$$

*where $S_1$, $S_2$ are $\Omega_2$-algebras and*

$$
g_1 : A \longrightarrow\!\!\ast\!\!\longrightarrow S_1 \qquad g_2 : A \longrightarrow\!\!\ast\!\!\longrightarrow S_2
$$

*are representations of $\Omega_1$-algebra $A$. We define morphism $R_1 \to R_2$ to be reduced morphism of representations $h : S_1 \to S_2$ making following diagram commutative*

*The reduced morphism* nat $N$ *of representation $f$ into representation $f_1$ (the theorem 3.4.7) is universally repelling in the category $\mathcal{A}$.* [3.8]

---

[3.7]The statement of lemma is similar to the statement on page [2]-119.

[3.8]See definition of universal object of category in definition on p. [2]-57.

Proof. From the theorem 2.1.6, it follows that there exists and unique the map $h$ for which the following diagram is commutative

$j = \text{nat } N \quad \ker R \supseteq N$

Therefore, we can uniquely define the map $h$ using equality

$$(3.4.11) \qquad h(j(b)) = R(b)$$

Let $\omega \in \Omega_2(p)$. Since maps $R$ and $j$ are homomorphisms of $\Omega_2$-algebra, then

$$(3.4.12) \qquad h(j(b_1)...j(b_p)\omega) = h(j(b_1...b_p\omega)) = R(b_1...b_p\omega) = R(b_1)...R(b_p)\omega$$
$$= h(j(b_1))...h(j(b_p))\omega$$

From the equality (3.4.12), it follows that the map $h$ is homomorphism of $\Omega_2$-algebra.

Since the map $R$ is reduced morphism of the representation $f$ into the representation $g$, then the following equality is satisfied

$$(3.4.13) \qquad g(a)(R(b)) = R(f(a)(b))$$

From the equality (3.4.11) it follows that

$$(3.4.14) \qquad g(a)(h(j(b))) = g(a)(R(b))$$

From the equalities (3.4.13), (3.4.14) it follows that

$$(3.4.15) \qquad g(a)(h(j(b))) = R(f(a)(b))$$

From the equalities (3.4.11), (3.4.15) it follows that

$$(3.4.16) \qquad g(a)(h(j(b))) = h(j(f(a)(b)))$$

From the equalities (3.4.8), (3.4.16) it follows that

$$(3.4.17) \qquad g(a)(h(j(b))) = h(f_1(a)(j(b)))$$

From the equality (3.4.17) it follows that the map $h$ is reduced morphism of representation $f_1$ into the representation $g$. □

## 3.5. Automorphism of Representation of Universal Algebra

Definition 3.5.1. Let

$$f : A_1 \longrightarrow A_2$$

be representation of $\Omega_1$-algebra $A_1$ in $\Omega_2$-algebra $A_2$. The reduced morphism of representations of $\Omega_1$-algebra

$$r_2 : A_2 \to A_2$$

such, that $r_2$ is endomorphism of $\Omega_2$-algebra is called **endomorphism of representation** $f$. □

THEOREM 3.5.2. *Given single transitive representation*

$$f : A_1 \longrightarrow_* A_2$$

*of $\Omega_1$-algebra $A_1$, for any $a_{21}, a_{22} \in A_2$ there exists unique endomorphism*

$$r_2 : A_2 \to A_2$$

*of representation $f$ such that $r_2(a_{21}) = a_{22}$.*

PROOF. Consider following diagram

Existence of endomorphism is corollary of the theorem 3.2.10. For given $p, q \in A_2$, uniqueness of endomorphism follows from the theorem 3.2.11 when $r_1 = \mathrm{id}$.    □

THEOREM 3.5.3. *Endomorphisms of representation $f$ form semigroup.*

PROOF. From theorem 3.3.1, it follows that the product of endomorphisms $(id, p_2)$,   $(id, r_2)$ of the representation $f$ is endomorphism $(id, p_2 \circ r_2)$ of the representation $f$.    □

DEFINITION 3.5.4. *Let*

$$f : A_1 \longrightarrow_* A_2$$

*be representation of $\Omega_1$-algebra $A_1$ in $\Omega_2$-algebra $A_2$. The morphism of representations of $\Omega_1$-algebra*

$$r_2 : A_2 \to A_2$$

*such, that $r_2$ is automorphism of $\Omega_2$-algebra is called* **automorphism of representation $f$.**    □

THEOREM 3.5.5. *Let*

$$f : A_1 \longrightarrow_* A_2$$

*be representation of $\Omega_1$-algebra $A_1$ in $\Omega_2$-algebra $A_2$. The set of automorphisms of the representation $f$ forms* **group** $GA(f)$.

PROOF. Let $r_2$, $p_2$ be automorphisms of the representation $f$. According to definition 3.5.4, maps $r_2$, $p_2$ are automorphisms of $\Omega_2$-algebra $A_2$. According to theorem II.3.2, ([14], p. 57), the map $r_2 \circ p_2$ is automorphism of $\Omega_2$-algebra $A_2$. From the theorem 3.3.1 and the definition 3.5.4, it follows that product of automorphisms $r_2 \circ p_2$ of the representation $f$ is automorphism of the representation $f$.

Let $r_2$, $p_2$, $q_2$ be automorphisms of the representation $f$. The associativity of product of maps $r_2$, $p_2$, $q_2$ follows from the chain of equations [3.9]

$$((r_2 \circ p_2) \circ q_2)(a) = (r_2 \circ p_2)(q_2(a)) = r_2(p_2(q_2(a)))$$
$$= r_2((p_2 \circ q_2)(a)) = (r_2 \circ (p_2 \circ q_2))(a)$$

Let $r_2$ be an automorphism of the representation $f$. According to definition 3.5.4 the map $r_2$ is automorphism of $\Omega_2$-algebra $A_2$. Therefore, the map $r_2^{-1}$ is

---

[3.9] To prove the associativity of product I follow to the example of the semigroup from [5], p. 20, 21.

automorphism of $\Omega_2$-algebra $A_2$. The equation (3.2.3) is true for automorphism $r_2$ of representation. Let $m' = r_2(m)$. Since $r_2$ is automorphism of $\Omega_2$-algebra, then $m = r_2^{-1}(m')$ and we can write (3.2.3) in the form

$$(3.5.1) \qquad r_2(f(a')(r_2^{-1}(m'))) = f(a')(m')$$

Since the map $r_2$ is automorphism of $\Omega_2$-algebra $A_2$, then from the equation (3.5.1) it follows that

$$(3.5.2) \qquad f(a')(r_2^{-1}(m')) = r_2^{-1}(f(a')(m'))$$

The equation (3.5.2) corresponds to the equation (3.2.3) for the map $r_2^{-1}$. Therefore, map $r_2^{-1}$ of the representation $f$.                                                 $\square$

CHAPTER 4

# $\Omega$-Group

## 4.1. Set of Homomorphisms of $\Omega$-Algebra

THEOREM 4.1.1. *Let sets $A$, $B$ be $\Omega$-algebras. Then the set* $\operatorname{Hom}(\Omega; A \to B)$ *also is $\Omega$-algebra when for any operations $\omega_1 \in \Omega(m)$, $\omega_2 \in \Omega(n)$, the following equality is true*

$$(4.1.1) \qquad (a_{11}...a_{1n}\omega_2)...(a_{m1}...a_{mn}\omega_2)\omega_1 = (a_{11}...a_{m1}\omega_1)...(a_{1n}...a_{mn}\omega_1)\omega_2$$

PROOF. According to the theorem 2.2.6, the set $B^A$ is $\Omega$-algebra. Let $\omega \in \Omega(n)$. For maps $f_1, ..., f_n \in B^A$, we define the operation $\omega$ by the equality

$$(4.1.2) \qquad (f_1...f_n\omega)(x) = f_1(x)...f_n(x)\omega$$

Let $\omega_1 \in \Omega(m)$, $\omega_2 \in \Omega(n)$. Let maps $f_1, ..., f_m \in \operatorname{Hom}(\Omega; A \to B)$ be homomorphisms from $\Omega$-algebra $A$ to $\Omega$-algebra $B$. In particular, for any $a_1, ..., a_n \in A$

$$(4.1.3) \qquad \begin{aligned} f_1(a_1...a_n\omega_2) &= f_1(a_1)...f_1(a_n)\omega_2 \\ ... &= ... \\ f_m(a_1...a_n\omega_2) &= f_m(a_1)...f_m(a_n)\omega_2 \end{aligned}$$

Since we require that the map $f_1...f_m\omega_1$ is homomorphism from $\Omega$-algebra $A$ to $\Omega$-algebra $B$, then

$$(4.1.4) \qquad (f_1...f_m\omega_1)(a_1...a_n\omega_2) = ((f_1...f_m\omega_1)(a_1))...((f_1...f_m\omega_1)(a_n))\omega_2$$

According to the definition (4.1.2), the equality

$$(4.1.5) \qquad \begin{aligned} &f_1(a_1...a_n\omega_2)...f_m(a_1...a_n\omega_2)\omega_1 \\ &= (f_1(a_1)...f_m(a_1)\omega_1)...(f_1(a_n)...f_m(a_n)\omega_1)\omega_2 \end{aligned}$$

follows from the equality (4.1.4). The equality

$$(4.1.6) \qquad \begin{aligned} &(f_1(a_1)...f_1(a_n)\omega_2)...(f_m(a_1)...f_m(a_n)\omega_2)\omega_1 \\ &= (f_1(a_1)...f_m(a_1)\omega_1)...(f_1(a_n)...f_m(a_n)\omega_1)\omega_2 \end{aligned}$$

follows from the equalities (4.1.3), (4.1.5). Let

$$(4.1.7) \qquad a_{ij} = f_i(a_j)$$

The equality (4.1.1) follows from the equalities (4.1.6), (4.1.7). $\qquad\square$

Not every $\Omega$-algebra satisfies to conditions of the theorem 4.1.1.

THEOREM 4.1.2. *Let $G_1$, $G_2$ be Abelian semigroups. The set* $\operatorname{Hom}(\{+\}; G_1 \to G_2)$ *also is Abelian semigroup.*

PROOF. Since sum in Abelian semigroup is commutative and associative, then the theorem follows from the theorem 4.1.1. $\qquad\square$

THEOREM 4.1.3. *The set* $\operatorname{End}(\{+\}; A)$ *of endomorphism of Abelian group $A$ is Abelian group.*

PROOF. The theorem follows from theorems 2.2.13, 4.1.2 and from the statement that the equation

$$x + a = 0$$

has root in Abelian group.     $\square$

THEOREM 4.1.4. *Let $D_1$, $D_2$ be rings. In general, the set* $\operatorname{Hom}(\{+, *\}; D_1 \to D_2)$ *is not ring.*

PROOF. There are two operations in the ring: sum which is commutative and associative and product which is distributive over sum. According to the theorem 4.1.1, sum and product must satisfy the equality

(4.1.8) $$a_{11}a_{21} + a_{12}a_{22} = (a_{11} + a_{12})(a_{21} + a_{22})$$

However right hand side of the equality (4.1.8) has form

$$(a_{11} + a_{12})(a_{21} + a_{22}) = (a_{11} + a_{12})a_{21} + (a_{11} + a_{12})a_{22}$$
$$= a_{11}a_{21} + a_{12}a_{21} + a_{11}a_{22} + a_{12}a_{22}$$

Therefore, the equality (4.1.8) is not true.     $\square$

Analysis of theorems 4.1.2, 4.1.4 tells us that the set of $\Omega$-algebras which satisfies to conditions of the theorem 4.1.1, is small.

QUESTION 4.1.5. *Is there a universal algebra which is different from the Abelian semigroup and satisfies to conditions of the theorem 4.1.1?*     $\square$

From our experience, it follows that certain $\Omega$-algebras contain an operation which alone generates semigroup. So we change the statement of the theorem 4.1.1.

THEOREM 4.1.6. *Let sets $A$, $B$ be $\Omega$-algebras. Let $\omega \in \Omega(n)$. Then the set* $\operatorname{Hom}(\Omega; A \to B)$ *is closed with respect to operation $\omega$ when the following equality is true*

(4.1.9) $$(a_{11}...a_{1n}\omega)...(a_{n1}...a_{nn}\omega)\omega$$
$$= (a_{11}...a_{n1}\omega)...(a_{1n}...a_{nn}\omega)\omega$$

PROOF. In general, we consider the set $\operatorname{Hom}(\{\omega\}; A \to B)$. The theorem follows from the theorem 4.1.1.     $\square$

THEOREM 4.1.7. *Let the operation $\omega \in \Omega(2)$ be commutative and associative. Then the set $\operatorname{Hom}(\Omega; A \to B)$ is closed with respect to operation $\omega$.*

PROOF. Since the operation $\omega \in \Omega(2)$ is commutative and associative, then

(4.1.10)
$$\begin{aligned}
(a_{11}a_{12}\omega)(a_{21}a_{22}\omega)\omega &= a_{11}(a_{12}(a_{21}a_{22}\omega)\omega)\omega \\
&= a_{11}((a_{12}a_{21}\omega)a_{22}\omega)\omega \\
&= a_{11}((a_{21}a_{12}\omega)a_{22}\omega)\omega \\
&= a_{11}(a_{21}(a_{12}a_{22}\omega)\omega)\omega \\
&= (a_{11}a_{21}\omega)(a_{12}...a_{22}\omega)\omega
\end{aligned}$$

The theorem follows from the equality (4.1.10) and from the theorem 4.1.6.     $\square$

Theorem 4.1.8. *Let the operation  $\omega \in \Omega(2)$   have a neutral element and the set*  $\mathrm{Hom}(\Omega; A \to B)$  *be closed with respect to operation $\omega$. Then the operation $\omega$ is commutative and associative.*

Proof. Equalities

$$(4.1.11) \qquad ab\omega = (ea\omega)(be\omega)\omega = (eb\omega)(ae\omega)\omega = ba\omega$$

$$(4.1.12) \qquad a(bc\omega) = (ae\omega)(bc\omega)\omega = (ab\omega)(ec\omega)\omega = (ab\omega)c\omega$$

follow from equalities (2.4.1), (2.4.2), (4.1.9). Commutativity of the operation $\omega$ follows from the equality (4.1.11). Associativity of the operation $\omega$ follows from equality (4.1.12). $\qquad\square$

Question 4.1.9. *Is there an operator domain $\Omega$, for which following statements are true?*

- *The set*  $\mathrm{Hom}(\Omega; A \to B)$  *is closed with respect to operation*  $\omega \in \Omega(2)$.
- *The operation $\omega$ is not commutative or associative.*

$\qquad\qquad\qquad\qquad\qquad\qquad\qquad\qquad\qquad\qquad\qquad\qquad\qquad\qquad\square$

## 4.2. $\Omega$-Group

Let the operation  $\omega \in \Omega_2(2)$  which is commutative and associative be defined in $\Omega_2$-algebra $A_2$. We identify the operation $\omega$ and sum. We use the symbol + to denote sum. Let

$$\Omega = \Omega_2 \setminus \{+\}$$

Definition 4.2.1. *A map*

$$f : A_2 \to B_2$$

*of $\Omega_2$-algebra $A_2$ into $\Omega_2$-algebra $B_2$ is called* **additive map** *if*

$$f(a + b) = f(a) + f(b)$$

*Let us denote*  $\mathcal{A}(A_2 \to B_2)$  *set of additive maps of $\Omega_2$-algebra $A_2$ into $\Omega_2$-algebra $B_2$.* $\qquad\square$

Theorem 4.2.2. $\mathcal{A}(A_2 \to B_2) = \mathrm{Hom}(\{+\}; A_2 \to B_2)$.

Proof. The theorem follows from definitions 2.2.9, 4.2.1. $\qquad\square$

Definition 4.2.3. *A map*

$$g : A^n \to A$$

*is called* **polyadditive map** *if for any  $i$, $i = 1, ..., n$,*

$$f(a_1, ..., a_i + b_i, ..., a_n) = f(a_1, ..., a_i, ..., a_n) + f(a_1, ..., b_i, ..., a_n)$$

$\qquad\qquad\qquad\qquad\qquad\qquad\qquad\qquad\qquad\qquad\qquad\qquad\qquad\qquad\square$

Theorem 4.2.4. *Let the map*

$$f : A_1 \longrightarrow\!\!\!\!*\!\!\!\!\longrightarrow A_2$$

*be the effective representation of $\Omega$-algebra $A_1$ in Abelian semigroup $A_2$.*

4.2.4.1: *On the set $A_1$ there is a structure of Abelian semigroup*

$$(4.2.1) \qquad f(a_1 + b_1)(a_2) = f(a_1)(a_2) + f(b_1)(a_2)$$

4.2.4.2: *The representation $f$ is additive map.*

4.2.4.3: *The map $f$ is the representation of $\Omega_1$-algebra $A_1$, where $\Omega_1 = \Omega \cup \{+\}$.*

PROOF. According to theorems 2.2.13, 4.1.7, the set $\text{End}(\{+\}, A_2)$ is Abelian semigroup. Since the representation $f$ is effective, then, according to theorems 3.1.3, 4.1.1, for any $A_1$-numbers $a$, $b$, there exists unique $A_1$-number $c$ such that

$$(4.2.2) \qquad f(c)(m) = f(a)(m) + f(b)(m)$$

Based on the equality (4.2.2), we introduce the sum of $A_1$-numbers

$$(4.2.3) \qquad c = a + b$$

The equality (4.2.1) follows from equalities (4.2.2), (4.2.3).

LEMMA 4.2.5. *The sum of $A_1$-numbers is commutative.*

PROOF. Since the sum of $A_2$-numbers is commutative, then the equality

$$(4.2.4) \qquad \begin{aligned} f(a_1 + b_1)(a_2) &= f(a_1)(a_2) + f(b_1)(a_2) = f(b_1)(a_2) + f(a_1)(a_2) \\ &= f(b_1 + a_1)(a_2) \end{aligned}$$

follows from the equality (4.2.1). The lemma follows from the equality (4.2.4). ⊙

LEMMA 4.2.6. *The sum of $A_1$-numbers is associative.*

PROOF. Since the sum of $A_2$-numbers is associative, then the equality

$$(4.2.5) \qquad \begin{aligned} f((a_1 + b_1) + c_1)(a_2) &= f(a_1 + b_1)(a_2) + f(c_1)(a_2) \\ &= (f(a_1)(a_2) + f(b_1)(a_2)) + f(c_1)(a_2) \\ &= f(a_1)(a_2) + (f(b_1)(a_2) + f(c_1)(a_2)) \\ &= f(a_1)(a_2) + f(b_1 + c_1)(a_2) \\ &= f(a_1 + (b_1 + c_1))(a_2) \end{aligned}$$

follows from the equality (4.2.1). The lemma follows from the equality (4.2.5). ⊙

The statement 4.2.4.1 follows from the equality (4.2.3), from lemmas 4.2.5, 4.2.6 and the definition 2.4.4.

The statement 4.2.4.2 follows from the equality (4.2.3). The statement 4.2.4.3 follows from the statement 4.2.4.2, since the map $f$ is homomorphism of $\Omega$-algebra.

□

THEOREM 4.2.7. *Let $\omega \in \Omega(n)$, $\omega_1 \in \Omega(m)$. The map*

$$(4.2.6) \qquad g : a_i \to a_1...a_n\omega$$

*is compatible with the operation $\omega_1$ when the following equality is true*

$$(4.2.7) \qquad a_1...(a_{i1}...a_{im}\omega_1)...a_n\omega = (a_1...a_{i1}...a_n\omega)...(a_1...a_{im}...a_n\omega)\omega_1$$

PROOF. The equality

$$(4.2.8) \qquad \begin{aligned} g(a_{i1}...a_{im}\omega_1) &= a_1...(a_{i1}...a_{im}\omega_1)...a_n\omega \\ &= (a_1...a_{i1}...a_n\omega)...(a_1...a_{im}...a_n\omega)\omega_1 \\ &= g(a_{i1})...g(a_{im})\omega_1 \end{aligned}$$

follows from equalities (4.2.7), (4.2.6). The theorem follows from the definition 2.2.9 and the equality (4.2.8).

□

The equality (4.2.7) is less restrictive than the equality (4.1.1). However, like in the case of the theorem 4.1.1, the majority of operations of universal algebra do not satisfy to the theorem 4.2.7. Since the addition satisfies to the theorem 4.1.1, we expect that there are conditions when addition satisfies to the theorem 4.2.7.

THEOREM 4.2.8. *Let* $\omega \in \Omega(n)$. *Since the map*

(4.2.9) $$g : a_i \to a_1...a_n\omega$$

*is compatible with the addition for any $i$, then the operation $\omega$ is polyadditive map.*

PROOF. According to the theorem 4.2.7, since the map (4.2.9) is compatible with the sum, then following equality is true

(4.2.10) $$a_1...(a_{i1} + a_{i2})...a_n\omega = (a_1...a_{i1}...a_n\omega) + (a_1...a_{i2}...a_n\omega)$$

The theorem follows from the equality (4.2.10) and the definition 4.2.3.    □

THEOREM 4.2.9. *Let* $\omega \in \Omega(n)$ *be polyadditive map. The operation $\omega$ is* **distributive** *over addition*

$$a_1...(a_i + b_i)...a_n\omega = a_1...a_i...a_n\omega + a_1...b_i...a_n\omega \quad i = 1,...,n$$

PROOF. The theorem follows from the theorem 4.2.8.    □

DEFINITION 4.2.10. *Let sum which is not necessarily commutative be defined in $\Omega_1$-algebra $A$. We use the symbol $+$ to denote sum. Let*

$$\Omega = \Omega_1 \setminus \{+\}$$

*If $\Omega_1$-algebra $A$ is group relative to sum and any operation $\omega \in \Omega$ is polyadditive map, then $\Omega_1$-algebra $A$ is called $\Omega$-group. If $\Omega$-group $A$ is associative group relative to sum, then $\Omega_1$-algebra $A$ is called* **associative** *$\Omega$-group. If $\Omega$-group $A$ is Abelian group relative to sum, then $\Omega_1$-algebra $A$ is called* **Abelian** *$\Omega$-group.*    □

EXAMPLE 4.2.11. *The group is the most evident example of $\Omega$-group.*
*A ring is $\Omega$-group.*
*Biring of matrices over division ring ([8]) is $\Omega$-group.*    □

REMARK 4.2.12. *Bourbaki consider similar definition, namely group with operators (see the definition 2 in [16] on page 31).*    □

THEOREM 4.2.13. *Let $A$ be $\Omega$-group. Let $\omega \in \Omega(n)$. The map*

$$g : a_i \to a_1...a_n\omega$$

*is endomorphism of additive group $A$.*

PROOF. The theorem follows from the theorem 4.2.9 and the definition 4.2.10.    □

THEOREM 4.2.14. *Let the map*

$$g : A_1 \longrightarrow A_2$$

*be the representation of $\Omega$-group $A_1$. Then the map*

$$\left( a_i \to a_1...a_n\omega \quad f(a_i) \to f(a_1)...f(a_n)\omega \right)$$

*is morphism of the representation $f$ of additive group $A_1$.*

PROOF. The theorem follows from the theorem 4.2.13 and definitions 3.1.1, 3.2.2.    □

### 4.3. Cartesian Product of Representations

LEMMA 4.3.1. *Let*

$$A = \prod_{i \in I} A_i$$

*be Cartesian product of family of $\Omega_2$-algebras $(A_i, i \in I)$. For each $i \in I$, let the set $\operatorname{End}(\Omega_2, A_i)$ be $\Omega_1$-algebra. Then the set*

(4.3.1)     $${}^\circ A = \{f \in \operatorname{End}(\Omega_2; A) : f(a_i, i \in I) = (f_i(a_i), i \in I)\}$$

*is Cartesian product of $\Omega_1$-algebras $\operatorname{End}(\Omega_2, A_i)$.*

PROOF. According to the definition (4.3.1), we can represent a map $f \in {}^\circ A$ as tuple

$$f = (f_i, i \in I)$$

of maps $f_i \in \operatorname{End}(\Omega_2; A_i)$. According to the definition (4.3.1),

$$(f_i, i \in I)(a_i, i \in I) = (f_i(a_i), i \in I)$$

Let $\omega \in \Omega_2$ be n-ary operation. We define operation $\omega$ on the set ${}^\circ A$ using equality

$$((f_{1i}, i \in I)...(f_{ni}, i \in I)\omega)(a_i, i \in I) = ((f_{1i}(a_i))...(f_{ni}(a_i))\omega, i \in I)$$

$\square$

DEFINITION 4.3.2. *Let $\mathcal{A}_1$ be category of $\Omega_1$-algebras. Let $\mathcal{A}_2$ be category of $\Omega_2$-algebras. We define* **category $\mathcal{A}_1(\mathcal{A}_2)$ of representations**. *Representations of $\Omega_1$-algebra in $\Omega_2$-algebra are objects of this category. Morphisms of corresponding representations are morphisms of this category.* $\square$

THEOREM 4.3.3. *In category $\mathcal{A}_1(\mathcal{A}_2)$ there exists product of single transitive representations of $\Omega_1$-algebra in $\Omega_2$-algebra.*

PROOF. For $j = 1, 2$, let

$$P_j = \prod_{i \in I} B_{ji}$$

be product of family of $\Omega_j$-algebras $\{B_{ji}, i \in I\}$ and for any $i \in I$ the map

$$t_{ji} : P_j \longrightarrow B_{ji}$$

be projection onto factor $i$. For each $i \in I$, let

$$h_i : B_{1i} \dashrightarrow B_{2i}$$

be single transitive $B_{1i}$-representation in $\Omega_2$-algebra $B_{2i}$.

Let $b_1 \in P_1$. According to the statement 2.3.3.3, we can represent $P_1$-number $b_1$ as tuple of $B_{1i}$-numbers

(4.3.2)     $$b_1 = (b_{1i}, i \in I) \quad b_{1i} = t_{1i}(b_1) \in B_{1i}$$

Let $b_2 \in P_2$. According to the statement 2.3.3.3, we can represent $P_2$-number $b_2$ as tuple of $B_{2i}$-numbers

(4.3.3)     $$b_2 = (b_{2i}, i \in I) \quad b_{2i} = t_{2i}(b_2) \in B_{2i}$$

LEMMA 4.3.4. *For each $i \in I$, consider diagram of maps*

(4.3.4)

*Let map*

$$g : P_1 \to \mathrm{End}(\Omega_2; P_2)$$

*be defined by the equality*

(4.3.5)                        $$g(b_1)(b_2) = (h_i(b_{1i})(b_{2i}), i \in I)$$

*Then the map $g$ is single transitive $P_1$-representation in $\Omega_2$-algebra $P_2$*

$$g : P_1 \xrightarrow{\quad * \quad} P_2$$

*The map $(t_{1i}, t_{2i})$ is morphism of representation $g$ into representation $h_i$.*

PROOF.

4.3.4.1: According to definitions 3.1.1, the map $h_i(b_{1i})$ is homomorphism of $\Omega_2$-algebra $B_{2i}$. According to the theorem 2.3.6, from commutativity of the diagram (1) for each $i \in I$, it follows that the map

$$g(b_1) : P_2 \to P_2$$

defined by the equality (4.3.5) is homomorphism of $\Omega_2$-algebra $P_2$.

4.3.4.2: According to the definition 3.1.1, the set $\mathrm{End}(\Omega_2; B_{2i})$ is $\Omega_1$-algebra. According to the lemma 4.3.1, the set $^\circ P_2 \subseteq \mathrm{End}(\Omega_2; P_2)$ is $\Omega_1$-algebra.

4.3.4.3: According to the definition 3.1.1, the map

$$h_i : B_{1i} \to \mathrm{End}(\Omega_2; B_{2i})$$

is homomorphism of $\Omega_1$-algebra. According to the theorem 2.3.6, the map

$$g : P_1 \to \mathrm{End}(\Omega_2; P_2)$$

defined by the equality

$$g(b_1) = (h_i(b_{1i}), i \in I)$$

is homomorphism of $\Omega_1$-algebra.

According to statements 4.3.4.1, 4.3.4.3 and to the definition 3.1.1, the map $g$ is $P_1$-representation in $\Omega_2$-algebra $P_2$.

Let $b_{21}$, $b_{22} \in P_2$. According to the statement 2.3.3.3, we can represent $P_2$-numbers $b_{21}$, $b_{22}$ as tuples of $B_{2i}$-numbers

(4.3.6)          $$b_{21} = (b_{21i}, i \in I) \quad b_{21i} = t_{2i}(b_{21}) \in B_{2i}$$

$$b_{22} = (b_{22i}, i \in I) \quad b_{22i} = t_{2i}(b_{22}) \in B_{2i}$$

According to the theorem 3.1.9, since the representation $h_i$ is single transitive, then there exists unique $B_{1i}$-number $b_{1i}$ such that

$$b_{22i} = h_i(b_{1i})(b_{21i})$$

According to definitions (4.3.2), (4.3.5), (4.3.6), there exists unique $P_1$-number $b_1$ such that

$$b_{22} = g(b_1)(b_{21})$$

According to the theorems 3.1.9, the representation $g$ is single transitive.

From commutativity of diagram (1) and from the definition 3.2.2, it follows that map $(t_{1i}, t_{2i})$ is morphism of representation $g$ into representation $h_i$.    ⊙

Let

(4.3.7) $$d_2 = g(b_1)(b_2) \quad d_2 = (d_{2i}, i \in I)$$

From equalities (4.3.5), (4.3.7), it follows that

(4.3.8) $$d_{2i} = h_i(b_{1i})(b_{2i})$$

For $j = 1, 2$, let $R_j$ be other object of category $\mathcal{A}_j$. For any $i \in I$, let the map

$$r_{1i} : R_1 \longrightarrow B_{1i}$$

be morphism from $\Omega_1$-algebra $R_1$ into $\Omega_1$-algebra $B_{1i}$. According to the definition 2.3.1, there exists a unique morphism of $\Omega_1$-algebra

$$s_1 : R_1 \longrightarrow P_1$$

such that following diagram is commutative

(4.3.9)

$$\begin{array}{ccc} P_1 & \xrightarrow{t_{1i}} & B_{1i} \\ \uparrow{\scriptstyle s_1} & & \\ R_1 & \xrightarrow{r_{1i}} & \end{array} \qquad t_{1i}(s_1) = r_{1i}$$

Let $a_1 \in R_1$. Let

(4.3.10) $$b_1 = s_1(a_1) \in P_1$$

From commutativity of the diagram (4.3.9) and statements (4.3.10), (4.3.2), it follows that

(4.3.11) $$b_{1i} = r_{1i}(a_1)$$

Let

$$f : R_1 \dashrightarrow R_2$$

be single transitive $R_1$-representation in $\Omega_2$-algebra $R_2$. According to the theorem 3.2.11, a morphism of $\Omega_2$-algebra

$$r_{2i} : R_2 \longrightarrow B_{2i}$$

such that map $(r_{1i}, r_{2i})$ is morphism of representations from $f$ into $h_i$ is unique up to choice of image of $R_2$-number $a_2$. According to the remark 3.2.7, in diagram

of maps

(4.3.12)

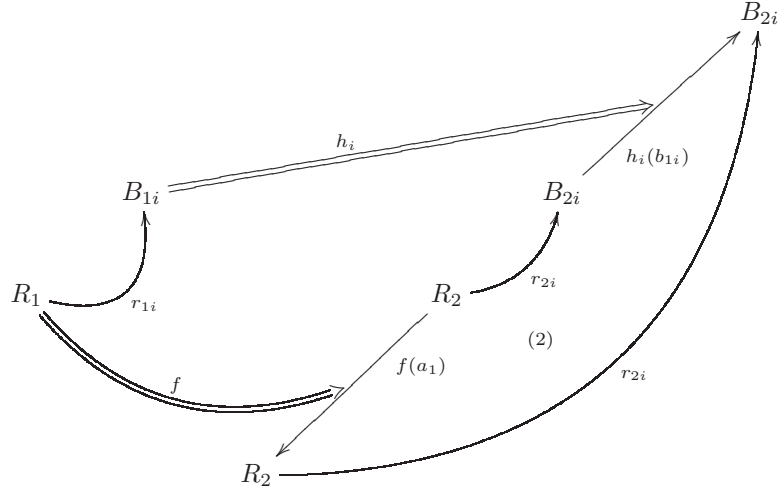

diagram (2) is commutative. According to the definition 2.3.1, there exists a unique morphism of $\Omega_2$-algebra

$$s_2 : R_2 \longrightarrow P_2$$

such that following diagram is commutative

(4.3.13)                    $\begin{array}{c} P_2 \xrightarrow{t_{2i}} B_{2i} \\ \end{array}$        $t_{2i}(s_2) = r_{2i}$

Let $a_2 \in R_2$. Let

(4.3.14)                              $b_2 = s_2(a_2) \in P_2$

From commutativity of the diagram (4.3.13) and statements (4.3.14), (4.3.3), it follows that

(4.3.15)                              $b_{2i} = r_{2i}(a_2)$

Let

(4.3.16)                              $c_2 = f(a_1)(a_2)$

From commutativity of the diagram (2) and equalities (4.3.8), (4.3.15), (4.3.16), it follows that

(4.3.17)                              $d_{2i} = r_{2i}(c_2)$

From equalities (4.3.8), (4.3.17), it follows that

(4.3.18)                              $d_2 = s_2(c_2)$

and this is consistent with commutativity of the diagram (4.3.13).

For each $i \in I$, we join diagrams of maps (4.3.4), (4.3.9), (4.3.13), (4.3.12)

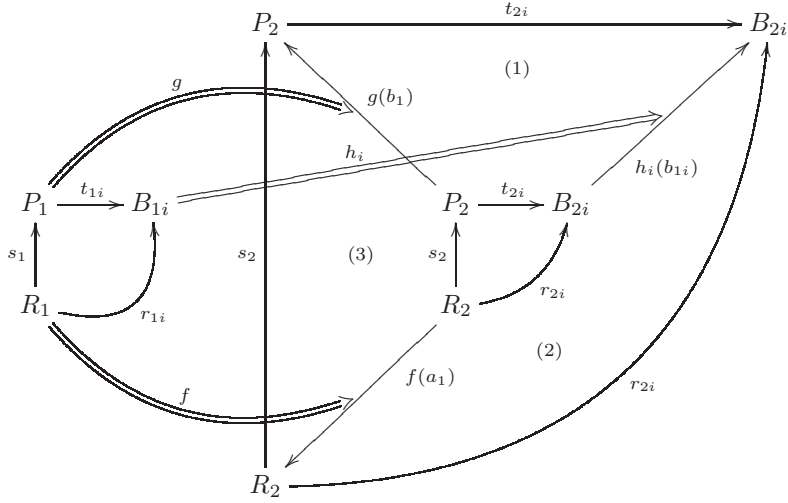

From equalities (4.3.7), (4.3.14) and from equalities (4.3.16), (4.3.18), commutativity of the diagram (3) follows. Therefore, the map $(s_1, s_2)$ is morphism of representations from $f$ into $g$. According to the theorem 3.2.11, the morphism $(s_1, s_2)$ is defined unambiguously, since we require (4.3.18).

According to the definition 2.3.1, the representation $g$ and family of morphisms of representation $((t_{1i}, t_{2i}), i \in I)$ is product in the category $\mathcal{A}_1(\mathcal{A}_2)$. $\square$

DEFINITION 4.3.5. Let $A_1, ..., A_n, A$ be $\Omega_1$-algebras. Let $B_1, ..., B_n, B$ be $\Omega_2$-algebras. Let, for any $k$, $k = 1, ..., n$,

$$f_k : A_k \longrightarrow\!\!\!* B_k$$

be representation of $\Omega_1$-algebra $A_k$ in $\Omega_2$-algebra $B_k$. Let

$$f : A \longrightarrow\!\!\!* B$$

be representation of $\Omega_1$-algebra $A$ in $\Omega_2$-algebra $B$. The map

$$\left( r_{1k} : A_k \to A \quad k = 1, ..., n \quad r_2 : B_1 \times ... \times B_n \to B \right)$$

is called **polymorphism of representations** $f_1, ..., f_n$ into representation $f$, if, for any $k$, $k = 1, ..., n$, provided that all variables except variables $a_k \in A_k$, $b_k \in B_k$ have given value, the map $(r_{1k}, r_2)$ is a morphism of representation $f_k$ into representation $f$.

If $f_1 = ... = f_n$, then we say that the map $((r_{1,k}, k = 1, ..., n) \quad r_2)$ is polymorphism of representation $f_1$ into representation $f$.

If $f_1 = ... = f_n = f$, then we say that the map $((r_{1,k}, k = 1, ..., n) \quad r_2)$ is polymorphism of representation $f$. $\square$

We also say that the map $r = (r_1, r_2)$ is polymorphism of representations in $\Omega_2$-algebras $B_1, ..., B_n$ into representation in $\Omega_2$-algebra $B$.

THEOREM 4.3.6. Let the map $((r_{1,k}, k = 1, ..., n) \quad r_2)$ be polymorphism of representations $f_1, ..., f_n$ into representation $f$. For any $k$, $k = 1, ..., n$, the map

$(r_{1k}, r_2)$ *satisfies to the equality*

(4.3.19) $$r_2(m_1, ..., f_k(a_k)(m_k), ..., m_n) = f(r_{1k}(a_k))(r_2(m_1, ..., m_n))$$

*Let* $\omega_1 \in \Omega_1(p)$. *For any* $k$, $k = 1, ..., n$, *the map* $r_{1k}$ *satisfies to the equality*

(4.3.20) $$r_{1k}(a_{k \cdot 1}...a_{k \cdot p}\omega_1) = r_{1k}(a_{k \cdot 1})...r_{1k}(a_{k \cdot p})\omega_1$$

*Let* $\omega_2 \in \Omega_2(p)$. *For any* $k$, $k = 1, ..., n$, *the map* $r_2$ *satisfies to the equality*

(4.3.21) $$r_2(m_1, ..., m_{k \cdot 1}...m_{k \cdot p}\omega_2, ..., m_n)$$
$$= r_2(m_1, ..., m_{k \cdot 1}, ..., m_n)...r_2(m_1, ..., m_{k \cdot p}, ..., m_n)\omega_2$$

PROOF. The equality (4.3.19) follows from the definition 4.3.5 and the equality (3.2.3). The equality (4.3.20) follows from the statement that, for any $k$, $k = 1, ..., n$, provided that all variables except the variable $x_k \in A_k$ have given value, the map $r_1$ is homomorphism of $\Omega_1$-algebra $A_k$ into $\Omega_1$-algebra $A$. The equality (4.3.21) follows from the statement that, for any $k$, $k = 1, ..., n$, provided that all variables except the variable $m_k \in B_k$ have given value, the map $r_2$ is homomorphism of $\Omega_2$-algebra $B_k$ into $\Omega_2$-algebra $B$. □

## 4.4. Reduced Cartesian Product of Representations

DEFINITION 4.4.1. *Let* $A_1$ *be* $\Omega_1$-*algebra. Let* $\mathcal{A}_2$ *be category of* $\Omega_2$-*algebras. We define* **category** $A_1(\mathcal{A}_2)$ **of representations** *of* $\Omega_1$-*algebra* $A_1$ *in* $\Omega_2$-*algebra. Representations of* $\Omega_1$-*algebra* $A_1$ *in* $\Omega_2$-*algebra are objects of this category. Reduced morphisms of corresponding representations are morphisms of this category.* □

THEOREM 4.4.2. *In category* $A_1(\mathcal{A}_2)$ *there exists product of effective representations of* $\Omega_1$-*algebra* $A_1$ *in* $\Omega_2$-*algebra and the product is effective representation of* $\Omega_1$-*algebra* $A_1$.

PROOF. Let
$$A_2 = \prod_{i \in I} A_{2i}$$
be product of family of $\Omega_2$-algebras $\{A_{2i}, i \in I\}$ and for any $i \in I$ the map
$$t_i : A_2 \longrightarrow A_{2i}$$
be projection onto factor $i$. For each $i \in I$, let
$$h_i : A_1 \dashrightarrow A_{2i}$$
be effective $A_1$-representation in $\Omega_2$-algebra $A_{2i}$.

Let $b_1 \in A_1$. Let $b_2 \in A_2$. According to the statement 2.3.3.3, we can represent $A_2$-number $b_2$ as tuple of $A_{2i}$-numbers

(4.4.1) $$b_2 = (b_{2i}, i \in I) \quad b_{2i} = t_i(b_2) \in A_{2i}$$

LEMMA 4.4.3. *For each $i \in I$, consider diagram of maps*

(4.4.2)

*Let map*

$$g : A_1 \to \operatorname{End}(\Omega_2; A_2)$$

*be defined by the equality*

(4.4.3) 
$$g(b_1)(b_2) = (h_i(b_1)(b_{2i}), i \in I)$$

*Then the map $g$ is effective $A_1$-representation in $\Omega_2$-algebra $A_2$*

$$g : A_1 \xrightarrow{\quad * \quad} A_2$$

*The map $t_i$ is reduced morphism of representation $g$ into representation $h_i$.*

PROOF.

4.4.3.1: According to definitions 3.1.1, the map $h_i(b_1)$ is homomorphism of $\Omega_2$-algebra $A_{2i}$. According to the theorem 2.3.6, from commutativity of the diagram (1) for each $i \in I$, it follows that the map

$$g(b_1) : A_2 \to A_2$$

defined by the equality (4.4.3) is homomorphism of $\Omega_2$-algebra $A_2$.

4.4.3.2: According to the definition 3.1.1, the set $\operatorname{End}(\Omega_2; A_{2i})$ is $\Omega_1$-algebra. According to the lemma 4.3.1, the set $^\circ A_2 \subseteq \operatorname{End}(\Omega_2; A_2)$ is $\Omega_1$-algebra.

4.4.3.3: According to the definition 3.1.1, the map

$$h_i : A_1 \to \operatorname{End}(\Omega_2, A_{2i})$$

is homomorphism of $\Omega_1$-algebra. According to the theorem 2.3.6, the map

$$g : A_1 \to \operatorname{End}(\Omega_2; A_2)$$

defined by the equality

$$g(b_1) = (h_i(b_1), i \in I)$$

is homomorphism of $\Omega_1$-algebra.

According to statements 4.4.3.1, 4.4.3.3 and to the definition 3.1.1, the map $g$ is $A_1$-representation in $\Omega_2$-algebra $A_2$.

For any $i \in I$, according to the definition 3.1.2, $A_1$-number $a_1$ generates unique transformation

(4.4.4) 
$$b_{22i} = h_i(b_1)(b_{21i})$$

Let $b_{21}$, $b_{22} \in A_2$. According to the statement 2.3.3.3, we can represent $A_2$-numbers $b_{21}$, $b_{22}$ as tuples of $A_{2i}$-numbers

(4.4.5) 
$$b_{21} = (b_{21i}, i \in I) \quad b_{21i} = t_i(b_{21}) \in A_{2i}$$
$$b_{22} = (b_{22i}, i \in I) \quad b_{22i} = t_i(b_{22}) \in A_{2i}$$

According to the definition (4.4.3) of the representation $g$, from equalities (4.4.4), (4.4.5), it follows that $A_1$-number $a_1$ generates unique transformation

$$(4.4.6) \qquad b_{22} = (h_i(b_1)(b_{21i}), i \in I) = g(b_1)(b_{21})$$

According to the definition 3.1.2, the representation $g$ is effective.

From commutativity of diagram (1) and from the definition 3.2.2, it follows that map $t_i$ is reduced morphism of representation $g$ into representation $h_i$.     ⊙

Let

$$(4.4.7) \qquad d_2 = g(b_1)(b_2) \quad d_2 = (d_{2i}, i \in I)$$

From equalities (4.4.3), (4.4.7), it follows that

$$(4.4.8) \qquad d_{2i} = h_i(b_1)(b_{2i})$$

Let $R_2$ be other object of category $\mathcal{A}_2$. Let

$$f : A_1 \rightarrowtail R_2$$

be effective $A_1$-representation in $\Omega_2$-algebra $R_2$. For any $i \in I$, let there exist morphism

$$r_i : R_2 \longrightarrow A_{2i}$$

of representations from $f$ into $h_i$. According to the remark 3.2.7, in diagram of maps

$$(4.4.9)$$

diagram (2) is commutative. According to the definition 2.3.1, there exists a unique morphism of $\Omega_2$-algebra

$$s : R_2 \longrightarrow A_2$$

such that following diagram is commutative

$$(4.4.10) \qquad \begin{array}{ccc} A_2 & \xrightarrow{t_i} & A_{2i} \\ {\scriptstyle s}\uparrow & & \\ R_2 & \xrightarrow{r_i} & \end{array} \qquad t_i(s) = r_i$$

Let $a_2 \in R_2$. Let

$$(4.4.11) \qquad b_2 = s(a_2) \in A_2$$

From commutativity of the diagram (4.4.10) and statements (4.4.11), (4.4.1), it follows that

(4.4.12) $$b_{2i} = r_i(a_2)$$

Let

(4.4.13) $$c_2 = f(a_1)(a_2)$$

From commutativity of the diagram (2) and equalities (4.4.8), (4.4.12), (4.4.13), it follows that

(4.4.14) $$d_{2i} = r_i(c_2)$$

From equalities (4.4.8), (4.4.14), it follows that

(4.4.15) $$d_2 = s(c_2)$$

and this is consistent with commutativity of the diagram (4.4.10).

For each $i \in I$, we join diagrams of maps (4.4.2), (4.4.10), (4.4.9)

From equalities (4.4.7), (4.4.11) and from equalities (4.4.13), (4.4.15), commutativity of the diagram (3) follows. Therefore, the map $s$ is reduced morphism of representations from $f$ into $g$. According to the definition 3.4.2, the map $s$ is homomorphism of $\Omega_2$ algebra. According to the theorem 2.3.3 and to the definition 2.3.1, the reduced morphism $s$ is defined unambiguously.

According to the definition 2.3.1, the representation $g$ and family of morphisms of representation $(t_i, i \in I)$ is product in the category $A_1(\mathcal{A}_2)$.          $\square$

DEFINITION 4.4.4. *Let $A$, $B_1$, ..., $B_n$, $B$ be universal algebras. Let, for any $k$, $k = 1$, ..., $n$,*

$$f_k : A \longrightarrow_* B_k$$

*be effective representation of $\Omega_1$-algebra $A$ in $\Omega_2$-algebra $B_k$. Let*

$$f : A \longrightarrow_* B$$

*be effective representation of $\Omega_1$-algebra $A$ in $\Omega_2$-algebra $B$. The map*

$$r_2 : B_1 \times ... \times B_n \to B$$

*is called* **reduced polymorphism of representations** $f_1$, ..., $f_n$ *into representation* $f$, *if, for any* $k$, $k = 1$, ..., $n$, *provided that all variables except the variable* $x_k \in B_k$ *have given value, the map* $r_2$ *is a reduced morphism of representation* $f_k$ *into representation* $f$.

If $f_1 = ... = f_n$, *then we say that the map* $r_2$ *is reduced polymorphism of representation* $f_1$ *into representation* $f$.

If $f_1 = ... = f_n = f$, *then we say that the map* $r_2$ *is reduced polymorphism of representation* $f$.                                                                  □

THEOREM 4.4.5. *Let the map* $r_2$ *be reduced polymorphism of effective representations* $f_1$, ..., $f_n$ *into effective representation* $f$.

- *For any* $k$, $k = 1$, ..., $n$, *the map* $r_2$ *satisfies to the equality*

(4.4.16)        $r_2(m_1, ..., f_k(a)(m_k), ..., m_n) = f(a)(r_2(m_1, ..., m_n))$

- *For any* $k$, $l$, $k = 1$, ..., $n$, $l = 1$, ..., $n$, *the map* $r_2$ *satisfies to the equality*

(4.4.17)
$$r_2(m_1, ..., f_k(a)(m_k), ..., m_l, ..., m_n)$$
$$= r_2(m_1, ..., m_k, ..., f_l(a)(m_l), ..., m_n)$$

- *Let* $\omega_2 \in \Omega_2(p)$. *For any* $k$, $k = 1$, ..., $n$, *the map* $r_2$ *satisfies to the equality*

(4.4.18)
$$r_2(m_1, ..., m_{k \cdot 1}...m_{k \cdot p}\omega_2, ..., m_n)$$
$$= r_2(m_1, ..., m_{k \cdot 1}, ..., m_n)...r_2(m_1, ..., m_{k \cdot p}, ..., m_n)\omega_2$$

PROOF. The equality (4.4.16) follows from the definition 4.4.4 and the equality (3.4.4). The equality (4.4.17) follows from the equality (4.4.16). The equality (4.4.18) follows from the statement that, for any $k$, $k = 1$, ..., $n$, provided that all variables except the variable $m_k \in B_k$ have given value, the map $r_2$ is homomorphism of $\Omega_2$-algebra $B_k$ into $\Omega_2$-algebra $B$.                                                                  □

We also say that the map $r_2$ is reduced polymorphism of representations in $\Omega_2$-algebras $B_1$, ..., $B_n$ into representation in $\Omega_2$-algebra $B$.

## 4.5. Multiplicative Ω-Group

Let the map

$$f : A \longrightarrow\!\!* B$$

be the representation of $\Omega_1$-algebra $A$ in $\Omega_2$ algebra $B$. According to the theorem 3.5.3, the set $\text{End}(A(\Omega_2); B)$ is semigroup. At the same time [4.1]

(4.5.1)                          $\text{End}(A(\Omega_2); B) \subseteq \text{End}(\Omega_2; B)$

According to the definition 3.1.1, the set $\text{End}(\Omega_2, B)$ is $\Omega_2$-algebra. However, the statement (4.5.1) does not imply that the set $\text{End}(A(\Omega_2); B)$ is $\Omega_2$-algebra.

To understand what is condition when the set $\text{End}(A(\Omega_2); B)$ is $\Omega_2$-algebra, we consider connection between the set of representations of $\Omega_1$-algebra $A$ in $\Omega_2$-algebra $B$ and the set of reduced morphisms of these representations.

---

[4.1] In the statement (4.5.1), I designated $\Omega_2$ category of $\Omega_2$-algebras and $A(\Omega_2)$ category of representations of $\Omega_1$-algebra $A$ in $\Omega_2$-algebra.

Theorem 4.5.1. *Let the map*

$$r : B \to B$$

*be reduced endomorphism of the representation*

$$f : A \overset{*}{\longrightarrow} B$$

*of $\Omega_1$-algebra $A$ in $\Omega_2$ algebra $B$. The map*

(4.5.2)          $rf : a \in A \to r \circ f(a) \in \text{End}(\Omega_2; B)$

*is representation of $\Omega_1$-algebra $A$ in $\Omega_2$ algebra $B$ iff, on the set $f(A) \subseteq \text{End}(\Omega_2, B)$, the product $\circ$ of maps is distributive on the left over any operation $\omega \in \Omega_1$*

(4.5.3)          $r \circ (f(a_1)...f(a_p)\omega) = (r \circ f(a_1))...(r \circ f(a_p))\omega$

Proof. According to the definition 3.1.1, the map $f(a)$ is emdomorphism of $\Omega_2$-algebra $B$. According to definitions 3.2.2, 3.4.2, the map $r$ is emdomorphism of $\Omega_2$-algebra $B$. Therefore, the map $r \circ f(a)$ is emdomorphism of $\Omega_2$-algebra $B$.

4.5.1.1: According to the definition 3.1.1, map $rf$ is representation of $\Omega_1$-algebra $A$ in $\Omega_2$ algebra $B$ iff the map $rf$ is homomorphism of $\Omega_1$-algebra.

4.5.1.2: The statement 4.5.1.1 means that, for any operation $\omega \in \Omega_1$, the following equality is true

(4.5.4)          $r \circ f(a_1...a_p\omega) = (rf)(a_1...a_p\omega) = ((rf)(a_1))...((rf)(a_p))\omega$
$$= (r \circ f(a_1))...(r \circ f(a_p))\omega$$

Since the map $f$ is representation of $\Omega_1$-algebra $A$ in $\Omega_2$ algebra $B$, then, according to the definition 3.1.1, the map $f$ is homomorphism of $\Omega_1$-algebra

(4.5.5)          $r \circ f(a_1...a_p\omega) = r \circ (f(a_1)...f(a_p)\omega)$

The equality (4.5.3) follows from equalities (4.5.4), (4.5.5).

The theorem follows from the statement 4.5.1.2.                                    □

Theorem 4.5.2. *Let the map*

$$f : A \overset{*}{\longrightarrow} B$$

*be the representation of $\Omega_1$-algebra $A$ in $\Omega_2$ algebra $B$. Let*

(4.5.6)          $f(A) = \text{End}(A(\Omega_2); B)$

4.5.2.1: *The product in semigroup $\text{End}(A(\Omega_2); B)$ is commutative.*

4.5.2.2: *The product $\circ$ in semigroup $\text{End}(A(\Omega_2); B)$ generates the product $*$ in $\Omega_1$-algebra $A$ such that*

(4.5.7)          $f(a * b) = f(a) \circ f(b)$

4.5.2.3: *The semigroup $\text{End}(A(\Omega_2); B)$ is $\Omega_1$-algebra.*

Proof. Let a map $h$ be endomorphism of the representation $f$. According to the statement (4.5.6), there exists $b \in A$ such that $h = f(b)$. Therefore, the equality

(4.5.8)          $f(a) \circ f(b) = f(b) \circ f(a)$

follows from the equality (3.4.4). According to the statement (4.5.6), maps $f(a)$, $f(b)$ are endomorphisms of the representation $f$. Therefore, the product $\circ$ in semigroup $\text{End}(A(\Omega_2); B)$ is commutative.

According to the theorem 3.5.3, the product of endomorphisms $f(a)$, $f(b)$ of
the representation $f$ is endomorphism $h$ of the representation $f$. According to the
statement (4.5.6), there exists $c \in A$ such that $h = f(c)$. Binary operation $*$ on the
set $A$ is defined by the equality

$$c = a * b$$

Therefore, the statement 4.5.2.2 is true.

Let maps $h_1, ..., h_p$ be endomorphism of the representation $f$. According to
the statement (4.5.6), there exist $A$-numbers $a_1, ..., a_p$ such that

$$h_1 = f(a_1) \quad ... \quad h_n = f(a_n)$$

Since the map $f$ is representation of $\Omega_1$-algebra $A$ in $\Omega_2$ algebra $B$, then, according
to the definition 3.1.1, the map $f$ is homomorphism of $\Omega_1$-algebra $A$

(4.5.9)                $h_1...h_p\omega = f(a_1)...f(a_p)\omega = f(a_1...a_p\omega)$

According to the statement (4.5.6), $h_1...h_p\omega \in \operatorname{End}(A(\Omega_2); B)$. Therefore, the
statement 4.5.2.3 is true.                                                                                      □

According to the theorem 4.5.2, if the statement (4.5.6) is satisfied, then the
set $\operatorname{End}(A(\Omega_2); B)$ is equiped by two algebraic structures. Namely, the set
$\operatorname{End}(A(\Omega_2); B)$ is semigroup and at the same time this set is $\Omega_1$-algebra. Sim-
ilar statement is true for $\Omega_1$-algebra $A$. However, we cannot say that product in
$\Omega_1$-algebra $A$ distributive over any operation $\omega \in \Omega_1$ (see the theorem 4.5.1).

Theorem 4.5.3. *Let the map*

$$f : A \longrightarrow\!\!\!* B$$

*be the representation of $\Omega_1$-algebra $A$ in $\Omega_2$ algebra $B$. and*

$$f(A) = \operatorname{End}(A(\Omega_2); B)$$

*The product $*$ defined in $\Omega_1$-algebra $A$ is distributive over any operation $\omega \in \Omega_1$
iff the map*

(4.5.10)                $f(b * a) : a \in A \to f(b * a) \in \operatorname{End}(\Omega_2; B)$

*is representation of $\Omega_1$-algebra $A$ in $\Omega_2$ algebra $B$*

Proof. According to the statement 4.5.2.2, it does not matter for us whether
we are considering $\Omega_1$-algebra $A$ or we are considering $\Omega_1$-algebra $\operatorname{End}(A(\Omega_2); B)$.
The theorem follows from the definition (4.5.7) of product $*$ in $\Omega_1$-algebra $A$, as
well it follows from the theorem 4.5.1 and statements 4.5.2.1, 4.5.2.3.                     □

In the theorem 4.5.3, we see universal algebra similar to $\Omega$-group, however this
algebra is little different. Since this universal algebra plays an important role in
representation theory, I introduce definitions 4.5.4, 4.5.5.

Definition 4.5.4. *Let product*

$$c_1 = a_1 * b_1$$

*be operation of $\Omega_1$-algebra $A$. Let $\Omega = \Omega_1 \setminus \{*\}$. If $\Omega_1$-algebra $A$ is group with
respect to product and, for any operation $\omega \in \Omega(n)$, the product is distributive
over the operation $\omega$*

$$a * (b_1...b_n\omega) = (a * b_1)...(a * b_n)\omega$$

$$(b_1...b_n\omega) * a = (b_1 * a)...(b_n * a)\omega$$

then $\Omega_1$-algebra $A$ is called **multiplicative $\Omega$-group**.                □

DEFINITION 4.5.5. *If*

$$(4.5.11) \qquad\qquad a * b = b * a$$

then multiplicative $\Omega$-group is called **Abelian**.                □

DEFINITION 4.5.6. *If*

$$(4.5.12) \qquad\qquad a * (b * c) = (a * b) * c$$

then multiplicative $\Omega$-group is called **associative**.                □

THEOREM 4.5.7. *Let $A$, $B_1$, ..., $B_n$, $B$ be universal algebras. Let, for any $k$, $k = 1$, ..., $n$,*

$$f_k : A \longrightarrow\!\!\!\!\ast\!\!\longrightarrow B_k$$

*be representation of $\Omega_1$-algebra $A$ in $\Omega_2$-algebra $B_k$. Let*

$$f : A \longrightarrow\!\!\!\!\ast\!\!\longrightarrow B$$

*be representation of $\Omega_1$-algebra $A$ in $\Omega_2$-algebra $B$. Let the map*

$$r_2 : B_1 \times ... \times B_n \to B$$

*be reduced polymorphism of representations $f_1$, ..., $f_n$ into representation $f$. The product $\circ$ defined in $\Omega_1$-algebra $f(A)$ is commutative.*

*The representation*

$$f : A \longrightarrow\!\!\!\!\ast\!\!\longrightarrow B$$

*permits reduced polymorphism of representations iff following statements are satisfied*

4.5.7.1: *The product $\circ$ defined in $\Omega_1$-algebra $\operatorname{End}(A(\Omega_2); B)$ is distributive over any operation $\omega \in \Omega_1$*

4.5.7.2: $f(a * b) = f(a) \circ f(b)$

PROOF. Using the equality (4.4.16), we can write an expression

$$(4.5.13) \qquad r_2(m_1, ..., f_k(a_k)(m_k), ..., f_l(a_l)(m_l), ..., m_n)$$

either in the following form

$$
\begin{aligned}
(4.5.14) \quad & r_2(m_1, ..., f_k(a_k)(m_k), ..., f_l(a_l)(m_l), ..., m_n) \\
& = f(a_k)(r_2(m_1, ..., m_k, ..., f_l(a_l)(m_l), ..., m_n)) \\
& = f(a_k)(f(a_l)(r_2(m_1, ..., m_k, ..., m_l, ..., m_n))) \\
& = (f(a_k) \circ f(a_l))(r_2(m_1, ..., m_k, ..., m_l, ..., m_n))
\end{aligned}
$$

or in the following form

$$
\begin{aligned}
(4.5.15) \quad & r_2(m_1, ..., f_k(a_k)(m_k), ..., f_l(a_l)(m_l), ..., m_n) \\
& = f(a_l)(r_2(m_1, ..., f_k(a_k)(m_k), ..., m_l, ..., m_n)) \\
& = f(a_l)(f(a_k)(r_2(m_1, ..., m_k, ..., m_l, ..., m_n))) \\
& = (f(a_l) \circ f(a_k))(r_2(m_1, ..., m_k, ..., m_l, ..., m_n))
\end{aligned}
$$

Commutativity of the product $\circ$ follows from the equalities (4.5.14), (4.5.15).                □

Theorem 4.5.8. *Let*

$$f : A \longrightarrow_* \longrightarrow B$$

*be representation of $\Omega_1$-algebra $A$ in $\Omega_2$-algebra $B$ and*

(4.5.16)                        $f(A) = \text{End}(A(\Omega_2); B)$

*Then the representation $f$ permits reduced polymorphism of representations.*

Let $\Omega = \Omega_1 \setminus \{*\}$. The representation

$$h : A_1 \to \text{End}(\Omega; A_1) \quad h(a) : b \in A_1 \to a * b \in A_1$$

*of semigroup $A_1$ in $\Omega$-algebra $A_1$ exists iff, for any operation $\omega \in \Omega(n)$, the product is **distributive** over the operation $\omega$*

(4.5.17)                   $a * (b_1 ... b_n \omega) = (a * b_1) ... (a * b_n) \omega$

(4.5.18)                   $(b_1 ... b_n \omega) * a = (b_1 * a) ... (b_n * a) \omega$

Proof. According to the definition 3.1.1, equalities (4.5.17), (4.5.18) are true iff the map $h$ is representation of semigroup $A_1$ in $\Omega$-algebra $A_1$. The same time equalities (4.5.17), (4.5.18) express distributive law of the product over the operation $\omega$.                                                                      □

In $\Omega_1$-algebra $A_1$, we defined the product coordinated with single transitive representation in $\Omega_2$-algebra $A_2$. We can do such construction in case of any representation with request that a product in $\Omega_1$-algebra $A_1$ is defined uniquely. However, in general case, a product may be non commutative.

Theorem 4.5.9. *Let*

$$A \longrightarrow_* \longrightarrow B_1 \qquad A \longrightarrow_* \longrightarrow B_2 \qquad A \longrightarrow_* \longrightarrow B$$

*be effective representations of Abelian multiplicative $\Omega_1$-group $A$ in $\Omega_2$-algebras $B_1$, $B_2$, $B$. Let $\Omega_2$-algebra have 2 operations, namely $\omega_1 \in \Omega(m)$, $\omega_2 \in \Omega(n)$. The equality*

(4.5.19)   $(a_{11} ... a_{1n} \omega_2) ... (a_{m1} ... a_{mn} \omega_2) \omega_1 = (a_{11} ... a_{m1} \omega_1) ... (a_{1n} ... a_{mn} \omega_1) \omega_2$

*is necessary condition of existence of reduced polymorphism*

$$R : B_1 \times B_2 \to B$$

Proof. Let $a_1, ..., a_p \in B_1$, $b_1, ..., b_q \in B_2$. According to the equality (4.4.18), the expression

(4.5.20)                        $r_2(a_1 ... a_p \omega_1, b_1 ... b_q \omega_2)$

can have 2 values

(4.5.21)
$$\begin{aligned}
&r_2(a_1 ... a_m \omega_1, b_1 ... b_n \omega_2) \\
&= r_2(a_1, b_1 ... b_n \omega_2) ... r_2(a_m, b_1 ... b_n \omega_2) \omega_1 \\
&= (r_2(a_1, b_1) ... r_2(a_1, b_n) \omega_2) ... (r_2(a_m, b_1) ... r_2(a_m, b_n) \omega_2) \omega_1
\end{aligned}$$

(4.5.22)
$$\begin{aligned}
&r_2(a_1 ... a_m \omega_1, b_1 ... b_n \omega_2) \\
&= r_2(a_1 ... a_m \omega_1, b_1) ... r_2(a_1 ... a_m \omega_1, b_n) \omega_2 \\
&= (r_2(a_1, b_1) ... r_2(a_m, b_1) \omega_1) ... (r_2(a_1, b_n) ... r_2(a_m, b_n) \omega_1) \omega_2
\end{aligned}$$

From equalities (4.5.21), (4.5.22), it follows that

$$(4.5.23) \quad (r_2(a_1, b_1)...r_2(a_1, b_n)\omega_2)...(r_2(a_m, b_1)...r_2(a_m, b_n)\omega_2)\omega_1$$
$$= (r_2(a_1, b_1)...r_2(a_m, b_1)\omega_1)...(r_2(a_1, b_n)...r_2(a_m, b_n)\omega_1)\omega_2$$

Therefore, the expression (4.5.20) is properly defined iff the equality (4.5.23) is true. Let

$$(4.5.24) \quad a_{i \cdot j} = r_2(a_i, b_j) \in A$$

The equality (4.5.19) follows from equalities (4.5.23), (4.5.24). $\qquad\square$

THEOREM 4.5.10. *There exists reduced polymorphism of effective representations of Abelian multiplicative $\Omega$-group in Abelian group.*

PROOF. Since sum in Abelian group is commutative and associative, then the theorem follows from the theorem 4.5.9. $\qquad\square$

THEOREM 4.5.11. *There is no reduced polymorphism of effective representations of Abelian multiplicative $\Omega$-group in ring.*

PROOF. There are two operations in the ring: sum which is commutative and associative and product which is distributive over sum. According to the theorem 4.5.9, the existence of polymorphism of effective representation in the ring implies that sum and product must satisfy the equality

$$(4.5.25) \quad a_{11}a_{21} + a_{12}a_{22} = (a_{11} + a_{12})(a_{21} + a_{22})$$

However right hand side of the equality (4.5.25) has form

$$(a_{11} + a_{12})(a_{21} + a_{22}) = (a_{11} + a_{12})a_{21} + (a_{11} + a_{12})a_{22}$$
$$= a_{11}a_{21} + a_{12}a_{21} + a_{11}a_{22} + a_{12}a_{22}$$

Therefore, the equality (4.5.25) is not true. $\qquad\square$

QUESTION 4.5.12. *It is possible that polymorphism of representations exists only for effective representation in Abelian group. However, this statement has not been proved.* $\qquad\square$

## 4.6. $\Omega$-ring

DEFINITION 4.6.1. *Let sum*

$$c_1 = a_1 + b_1$$

*which is not necessarily commutative and product*

$$c_1 = a_1 * b_1$$

*be operations of $\Omega_1$-algebra $A$. Let $\Omega = \Omega_1 \setminus \{+, *\}$. If $\Omega_1$-algebra $A$ is $\Omega \cup \{*\}$-group and multiplicative $\Omega \cup \{+\}$-group, then $\Omega_1$-algebra $A$ is called $\Omega$-ring.* $\qquad\square$

THEOREM 4.6.2. *The product in $\Omega$-ring is distributive over addition*

$$a * (b_1 + b_2) = a * b_1 + a * b_2$$
$$(b_1 + b_2) * a = b_1 * a + b_2 * a$$

PROOF. The theorem follows from the definitions 4.2.10, 4.5.4, 4.6.1. $\qquad\square$

DEFINITION 4.6.3. *Let $A$ be $\Omega$-ring. The* **matrix** *over $\Omega$-ring $A$ is a table of $A$-numbers $a_j^i$ where the index $i$ is the number of row and the index $j$ is the number of column.* □

CONVENTION 4.6.4. *We will use Einstein summation convention. When an index is present in an expression twice (one above and one below) and a set of index is known, we have the sum with respect to repeated index. In this case we assume that we know the set of summation index and do not use summation symbol*

$$a^i v_i = \sum_{i \in I} a^i v_i$$

*If needed to clearly show set of index, I will do it.* □

The product of matrices is associated with the product of homomorphisms of vector spaces over field. According to the custom the product of matrices $a$ and $b$ is defined as product of rows of the matrix $a$ and columns of the matrix $b$.

EXAMPLE 4.6.5. *Let $\overline{\overline{e}}$ be basis of right vector space $V$ over $D$-algebra $A$ (see the definition 9.6.2 and the theorem 9.6.15). We represent the basis $\overline{\overline{e}}$ as row of matrix*

$$e = \begin{pmatrix} e_1 & ... & e_n \end{pmatrix}$$

*We represent coordinates of vector $v$ as vector column*

$$v = \begin{pmatrix} v^1 \\ ... \\ v^n \end{pmatrix}$$

*Therefore, we can represent the vector $v$ as product of matrices*

$$v = \begin{pmatrix} e_1 & ... & e_n \end{pmatrix} \begin{pmatrix} v^1 \\ ... \\ v^n \end{pmatrix} = e_i v^i$$

*We represent homomorphism of right vector space $V$ using matrix*

(4.6.1)                              $$v'^i = f_j^i v^j$$

*The equality (4.6.1) expresses a traditional product of matrices $f$ and $v$.* □

EXAMPLE 4.6.6. *Let $\overline{\overline{e}}$ be basis of left vector space $V$ over $D$-algebra $A$ (see the definition 9.5.2 and the theorem 9.5.15). We represent the basis $\overline{\overline{e}}$ as row of matrix*

$$e = \begin{pmatrix} e_1 & ... & e_n \end{pmatrix}$$

*We represent coordinates of vector $v$ as vector column*

$$v = \begin{pmatrix} v^1 \\ ... \\ v^n \end{pmatrix}$$

*However, we cannot represent the vector*

$$v = v^i e_i$$

*as product of matrices*

$$v = \begin{pmatrix} v^1 \\ ... \\ v^n \end{pmatrix} \qquad e = \begin{pmatrix} e_1 & ... & e_n \end{pmatrix}$$

*because this product is not defined. We represent homomorphism of left vector space $V$ using matrix*

(4.6.2) $$v'^i = v^j f_j^i$$

*We cannot express the equality (4.6.2) as traditional product of matrices $v$ and $f$.*

$\square$

From examples 4.6.5, 4.6.6, it follows that we cannot confine ourselves to traditional product of matrices and we need to define two products of matrices. To distinguish between these products we introduced a new notation. In order to keep this notation consistent with the existing one we assume that we have in mind $_*{}^*$-product when no clear notation is present.

DEFINITION 4.6.7. *Let the nubmer of columns of the matrix $a$ equal the number of rows of the matrix $b$.* $_*{}^*$-**product** *of matrices $a$ and $b$ has form*

(4.6.3) $$\begin{cases} a_*{}^*b = \left( a_k^i b_j^k \right) \\ (a_*{}^*b)_j^i = a_k^i b_j^k \end{cases}$$

*and can be expressed as product of a row of matrix $a$ over a column of matrix $b$.* [4.2]

$\square$

DEFINITION 4.6.8. *Let the nubmer of rows of the matrix $a$ equal the number of columns of the matrix $b$.* $^*{}_*$-**product** *of matrices $a$ and $b$ has form*

(4.6.4) $$\begin{cases} a^*{}_*b = \left( a_i^k b_k^j \right) \\ (a^*{}_*b)_i^j = a_i^k b_k^j \end{cases}$$

*and can be expressed as product of a column of matrix $a$ over a row of matrix $b$.* [4.3]

$\square$

We also consider following operations on the set of matrices.

DEFINITION 4.6.9. *The transpose $a^T$ of the matrix $a$ exchanges rows and columns*

(4.6.5) $$(a^T)_j^i = a_i^j$$

$\square$

---

[4.2] We will use symbol $_*{}^*$- in following terminology and notation. We will read symbol $_*{}^*$ as $rc$-product or product of row over column. To draw symbol of product of row over column, we put two symbols of product in the place of index which participate in sum. For instance, if product of $A$-numbers has form $a \circ b$, then $_*{}^*$-product of matrices $a$ and $b$ has form $a_\circ{}^\circ b$.

[4.3] We will use symbol $^*{}_*$ in following terminology and notation. We will read symbol $^*{}_*$ as $cr$-product or product of column over row. To draw symbol of product of column over row, we put two symbols of product in the place of index which participate in sum. For instance, if product of $A$-numbers has form $a \circ b$, then $^*{}_*$-product of matrices $a$ and $b$ has form $a^\circ{}_\circ b$.

Definition 4.6.10. *The sum of matrices a and b is defined by the equality*

$$(a+b)^i_j = a^i_j + b^i_j$$

<div align="right">□</div>

Remark 4.6.11. *We will use symbol $_*{}^*$- or $^*{}_*$- in name of properties of each product and in the notation. We can read symbols $_*{}^*$ and $^*{}_*$ as rc-product and cr-product. This rule we extend to following terminology.*

<div align="right">□</div>

Theorem 4.6.12.

$$(4.6.6) \qquad (a_*{}^*b)^T = a^T{}^*{}_*b^T$$

Proof. The chain of equalities

$$(4.6.7) \qquad ((a_*{}^*b)^T)^j_i = (a_*{}^*b)^i_j = a^i_k b^k_j = (a^T)^k_i (b^T)^j_k = ((a^T)^*{}_*(b^T))^j_i$$

follows from (4.6.5), (4.6.3) and (4.6.4). The equality (4.6.6) follows from (4.6.7).

<div align="right">□</div>

Definition 4.6.13. *The set $\mathcal{A}$ is a* **biring** *if we defined on $\mathcal{A}$ an unary operation, say transpose, and three binary operations, say $_*{}^*$-product, $^*{}_*$-product and sum, such that*

- $_*{}^*$-product and sum define structure of ring on $\mathcal{A}$
- $^*{}_*$-product and sum define structure of ring on $\mathcal{A}$
- both products have common identity $\delta$
- products satisfy equation

$$(a_*{}^*b)^T = a^T{}^*{}_*b^T$$

- transpose of identity is identity

$$(4.6.8) \qquad \delta^T = \delta$$

- double transpose is original element

$$(4.6.9) \qquad (a^T)^T = a$$

<div align="right">□</div>

Theorem 4.6.14 (duality principle for biring). *Let $\mathcal{A}$ be true statement about biring $A$. If we exchange the same time*

- $a \in A$ and $a^T$
- $_*{}^*$-product and $^*{}_*$-product

*then we soon get true statement.*

Theorem 4.6.15 (duality principle for biring of matrices). *Let $A$ be biring of matrices. Let $\mathcal{A}$ be true statement about matrices. If we exchange the same time*

- rows and columns of all matrices
- $_*{}^*$-product and $^*{}_*$-product

*then we soon get true statement.*

Proof. This is the immediate consequence of the theorem 4.6.14.

<div align="right">□</div>

REMARK 4.6.16. *If product in $\Omega$-$A$ ring is commutative, then*

(4.6.10)                    $$a_{*}{}^{*}b = (a_i^k b_k^j) = (b_k^j a_i^k) = b^{*}{}_{*}a$$

**Reducible biring** *is the biring which holds* **condition of reducibility of products** (4.6.10). *So, in reducible biring, it is enough to consider only $_{*}{}^{*}$-product. However in case when the order of factors is essential we will use $^{*}{}_{*}$-product also.* □

## 4.7. Tensor Product of Representations

DEFINITION 4.7.1. *Let $A$ be Abelian multiplicative $\Omega_1$-group. Let $A_1$, ..., $A_n$ be $\Omega_2$-algebras.* [4.4] *Let, for any $k$, $k = 1, ..., n$,*

$$f_k : A \dashrightarrow A_k$$

*be effective representation of multiplicative $\Omega_1$-group $A$ in $\Omega_2$-algebra $A_k$. Consider category $\mathcal{A}$ whose objects are reduced polymorphisms of representations $f_1$, ..., $f_n$*

$$r_1 : B_1 \times ... \times B_n \longrightarrow S_1 \qquad r_2 : B_1 \times ... \times B_n \longrightarrow S_2$$

*where $S_1$, $S_2$ are $\Omega_2$-algebras and*

$$g_1 : A \dashrightarrow S_1 \qquad g_2 : A \dashrightarrow S_2$$

*are effective representations of multiplicative $\Omega_1$-group $A$. We define morphism $r_1 \to r_2$ to be reduced morphism of representations $h : S_1 \to S_2$ making following diagram commutative*

*Universal object $B_1 \otimes ... \otimes B_n$ of category $\mathcal{A}$ is called* **tensor product** *of representations $A_1$, ..., $A_n$.* □

THEOREM 4.7.2. *Since there exists tensor product of effective representations, then tensor product is unique up to isomorphism of representations.*

PROOF. Let $A$ be Abelian multiplicative $\Omega_1$-group. Let $A_1$, ..., $A_n$ be $\Omega_2$-algebras. Let, for any $k$, $k = 1, ..., n$,

$$f_k : A \longrightarrow B_k$$

be effective representation of multiplicative $\Omega_1$-group $A$ in $\Omega_2$-algebra $B_k$. Let effective representations

$$g_1 : A \dashrightarrow S_1 \qquad g_2 : A \dashrightarrow S_2$$

_________________________

[4.4] I give definition of tensor product of representations of universal algebra following to definition in [2], p. 601 - 603.

be tensor product of representations $B_1$, ..., $B_n$. From commutativity of the diagram

(4.7.1)

it follows that

(4.7.2)
$$R_1 = h_2 \circ h_1 \circ R_1$$
$$R_2 = h_1 \circ h_2 \circ R_2$$

From equalities (4.7.2), it follows that morphisms of representation $h_1 \circ h_2$, $h_2 \circ h_1$ are identities. Therefore, morphisms of representation $h_1$, $h_2$ are isomorphisms. □

CONVENTION 4.7.3. *Algebras $S_1$, $S_2$ may be different sets. However they are indistinguishable for us when we consider them as isomorphic representations. In such case, we write the statement $S_1 = S_2$.* □

DEFINITION 4.7.4. *Tensor product*

$$B^{\otimes n} = B_1 \otimes ... \otimes B_n \quad B_1 = ... = B_n = B$$

*is called* **tensor power** *of representation $B$.* □

THEOREM 4.7.5. *Since there exists polymorphism of representations, then there exists tensor product of representations.*

PROOF. Let

$$f : A \overset{*}{\longrightarrow} M$$

be representation of $\Omega_1$-algebra $A$ generated by Cartesian product $B_1 \times ... \times B_n$ of sets $B_1$, ..., $B_n$.[4.5] Injection

$$i : B_1 \times ... \times B_n \longrightarrow M$$

is defined according to rule [4.6]

(4.7.3)                          $$i \circ (b_1, ..., b_n) = (b_1, ..., b_n)$$

Let $N$ be equivalence generated by following equalities [4.7]

(4.7.4)         $$(b_1, ..., b_{i \cdot 1} ... b_{i \cdot p} \omega, ..., b_n) = (b_1, ..., b_{i \cdot 1}, ..., b_n)...(b_1, ..., b_{i \cdot p}, ..., b_n)\omega$$

(4.7.5)         $$(b_1, ..., f_i(a) \circ b_i, ..., b_n) = f(a) \circ (b_1, ..., b_i, ..., b_n)$$

$$b_k \in B_k \quad k = 1, ..., n \quad b_{i \cdot 1}, ..., b_{i \cdot p} \in B_i \quad \omega \in \Omega_2(p) \quad a \in A$$

---

[4.5]According to theorems 2.3.3, 4.4.2, the set generated by reduced Cartesian product of representations $B_1$, ..., $B_n$ coincides with Cartesian product $B_1 \times ... \times B_n$ of sets $B_1$, ..., $B_n$. At this point of the proof, we do not consider any algebra structure on the set $B_1 \times ... \times B_n$.

[4.6]The equality (4.7.3) states that we identify the basis of the representation $M$ with the set $B_1 \times ... \times B_n$.

[4.7] I considered generating of elements of representation according to the theorem 6.1.4. The theorem 4.7.11 requires the fulfillment of conditions (4.7.4), (4.7.5).

Lemma 4.7.6. *Let $\omega \in \Omega_2(p)$. Then*

$$(4.7.6) \qquad \begin{aligned} &f(c) \circ (b_1, ..., b_{i\cdot 1}...b_{i\cdot p}\omega, ..., b_n) \\ &= f(c) \circ ((b_1, ..., b_{i\cdot 1}, ..., b_n)...(b_1, ..., b_{i\cdot p}, ..., b_n)\omega) \end{aligned}$$

Proof. From the equality (4.7.5), it follows that

$$(4.7.7) \qquad f(c) \circ (b_1, ..., b_{i\cdot 1}...b_{i\cdot p}\omega, ..., b_n) = (b_1, ..., f_i(c) \circ (b_{i\cdot 1}...b_{i\cdot p}\omega), ..., b_n)$$

Since $f_i(c)$ is endomorphism of $\Omega_2$-algebra $B_i$, then from the equality (4.7.7), it follows that

$$(4.7.8) \quad f(c) \circ (b_1, ..., b_{i\cdot 1}...b_{i\cdot p}\omega, ..., b_n) = (b_1, ..., (f_i(c) \circ b_{i\cdot 1})...(f_i(c) \circ b_{i\cdot p})\omega, ..., b_n)$$

From equalities (4.7.8), (4.7.4), it follows that

$$(4.7.9) \qquad \begin{aligned} &f(c) \circ (b_1, ..., b_{i\cdot 1}...b_{i\cdot p}\omega, ..., b_n) \\ &= (b_1, ..., f_i(c) \circ b_{i\cdot 1}, ..., b_n)...(b_1, ..., f_i(c) \circ b_{i\cdot p}, ..., b_n)\omega \end{aligned}$$

From equalities (4.7.9), (4.7.5), it follows that

$$(4.7.10) \qquad \begin{aligned} &f(c) \circ (b_1, ..., b_{i\cdot 1}...b_{i\cdot p}\omega, ..., b_n) \\ &= (f(c) \circ (b_1, ..., b_{i\cdot 1}, ..., b_n))...(f(c) \circ (b_1, ..., b_{i\cdot p}, ..., b_n))\omega \end{aligned}$$

Since $f(c)$ is endomorphism of $\Omega_2$-algebra $B$, then the equality (4.7.6) follows from the equality (4.7.10).                    $\odot$

Lemma 4.7.7.

$$(4.7.11) \qquad f(c) \circ (b_1, ..., f_i(a) \circ b_i, ..., b_n) = f(c) \circ (f(a) \circ (b_1, ..., b_i, ..., b_n))$$

Proof. From the equality (4.7.5), it follows that

$$(4.7.12) \qquad \begin{aligned} f(c) \circ (b_1, ..., f_i(a) \circ b_i, ..., b_n) &= (b_1, ..., f_i(c) \circ (f_i(a) \circ b_i), ..., b_n) \\ &= (b_1, ..., (f_i(c) \circ f_i(a)) \circ b_i, ..., b_n) \\ &= (f(c) \circ f(a)) \circ (b_1, ..., b_i, ..., b_n) \\ &= f(c) \circ (f(a) \circ (b_1, ..., b_i, ..., b_n)) \end{aligned}$$

The equality (4.7.11) follows from the equality (4.7.12).                    $\odot$

Lemma 4.7.8. *For any $c \in A$, endomorphism $f(c)$ of $\Omega_2$-algebra $M$ is coordinated with equivalence $N$.*

Proof. The lemma follows from lemmas 4.7.6, 4.7.7 and from the definition 3.3.2.                    $\odot$

From the lemma 4.7.8 and the theorem 3.3.3, it follows that $\Omega_1$-algebra is defined on the set $^*M/N$. Consider diagram

According to lemma 4.7.8, from the condition

$$j \circ b_1 = j \circ b_2$$

it follows that

$$j \circ (f(a) \circ b_1) = j \circ (f(a) \circ b_2)$$

Therefore, transformation $F(a)$ is well defined and

(4.7.13)                         $$F(a) \circ j = j \circ f(a)$$

If $\omega \in \Omega_1(p)$, then we assume

$$(F(a_1)...F(a_p)\omega) \circ (J \circ b) = J \circ ((f(a_1)...f(a_p)\omega) \circ b)$$

Therefore, map $F$ is representations of $\Omega_1$-algebra $A$. From (4.7.13) it follows that $j$ is reduced morphism of representations $f$ and $F$.

Consider commutative diagram

(4.7.14)

From commutativity of the diagram (4.7.14) and from the equality (4.7.3), it follows that

(4.7.15)                 $$g_1 \circ (b_1, ..., b_n) = j \circ (b_1, ..., b_n)$$

From equalities (4.7.3), (4.7.4), (4.7.5), it follows that

(4.7.16)
$$g_1 \circ (b_1, ..., b_{i \cdot 1}...b_{i \cdot p}\omega, ..., b_n)$$
$$= (g_1 \circ (b_1, ..., b_{i \cdot 1}, ..., b_n))...(g_1 \circ (b_1, ..., b_{i \cdot p}, ..., b_n))\omega$$

(4.7.17)        $$g_1 \circ (b_1, ..., f_i(a) \circ b_i, ..., b_n) = f(a) \circ (g_1 \circ (b_1, ..., b_i, ..., b_n))$$

From equalities (4.7.16) and (4.7.17) it follows that map $g_1$ is reduced polymorphism of representations $f_1, ..., f_n$.

Since $B_1 \times ... \times B_n$ is the basis of representation $M$ of $\Omega_1$ algebra $A$, then, according to the theorem 6.2.10, for any representation

$$A \longrightarrow V$$

and any reduced polymorphism

$$g_2 : B_1 \times ... \times B_n \longrightarrow V$$

there exists a unique morphism of representations $k : M \to V$, for which following diagram is commutative

(4.7.18)

$$\begin{array}{ccc} B_1 \times ... \times B_n & \xrightarrow{\quad} & M \\ & \searrow{\scriptstyle g_2} & \downarrow{\scriptstyle k} \\ & & V \end{array}$$

Since $g_2$ is reduced polymorphism, then $\ker k \supseteq N$.

According to the theorem 3.4.8, map $j$ is universal in the category of morphisms of representation $f$ whose kernel contains $N$. Therefore, we have morphism of representations

$$h : M/N \to V$$

which makes the following diagram commutative

(4.7.19)

$$\begin{array}{ccc} & & M/N \\ & \nearrow{\scriptstyle j} & \downarrow{\scriptstyle h} \\ M & \xrightarrow{\scriptstyle k} & V \end{array}$$

We join diagrams (4.7.14), (4.7.18), (4.7.19), and get commutative diagram

$$\begin{array}{ccc} & & M/N \\ & \nearrow{\scriptstyle g_1} & \\ B_1 \times ... \times B_n & \xrightarrow{\scriptstyle i} M & \downarrow{\scriptstyle h} \\ & \searrow{\scriptstyle g_2} & \searrow{\scriptstyle k} \\ & & V \end{array}$$

Since $\operatorname{Im} g_1$ generates $M/N$, than map $h$ is uniquely determined.                    $\square$

According to proof of theorem 4.7.5

$$B_1 \otimes ... \otimes B_n = M/N$$

If $d_i \in A_i$, we write

(4.7.20)                    $$j \circ (d_1, ..., d_n) = d_1 \otimes ... \otimes d_n$$

From equalities (4.7.15), (4.7.20), it follows that

(4.7.21)                    $$g_1 \circ (d_1, ..., d_n) = d_1 \otimes ... \otimes d_n$$

THEOREM 4.7.9. *The map*

$$(x_1, ..., x_n) \in B_1 \times ... \times B_n \to x_1 \otimes ... \otimes x_n \in B_1 \otimes ... \otimes B_n$$

*is polymorphism.*

PROOF. The theorem follows from definitions 4.4.4, 4.7.1.                    □

THEOREM 4.7.10. *Let $B_1$, ..., $B_n$ be $\Omega_2$-algebras. Let*

$$f : B_1 \times ... \times B_n \to B_1 \otimes ... \otimes B_n$$

*be reduced polymorphism defined by equality*

(4.7.22)                     $$f \circ (b_1, ..., b_n) = b_1 \otimes ... \otimes b_n$$

*Let*

$$g : B_1 \times ... \times B_n \to V$$

*be reduced polymorphism into $\Omega$-algebra $V$. There exists morphism of representations*

$$h : B_1 \otimes ... \otimes B_n \to V$$

*such that the diagram*

*is commutative.*

PROOF. equality (4.7.22) follows from equalities (4.7.3) and (4.7.20). An existence of the map $h$ follows from the definition 4.7.1 and constructions made in the proof of the theorem 4.7.5.                    □

THEOREM 4.7.11. *Let*

$$b_k \in B_k \quad k = 1, ..., n \quad b_{i \cdot 1}, ..., b_{i \cdot p} \in B_i \quad \omega \in \Omega_2(p) \quad a \in A$$

*Tensor product is distributive over operation $\omega$*

(4.7.23)      $$\begin{aligned} b_1 \otimes ... \otimes (b_{i \cdot 1} ... b_{i \cdot p} \omega) \otimes ... \otimes b_n \\ = (b_1 \otimes ... \otimes b_{i \cdot 1} \otimes ... \otimes b_n) ... (b_1 \otimes ... \otimes b_{i \cdot p} \otimes ... \otimes b_n) \omega \end{aligned}$$

*The representation of multiplicative $\Omega_1$-group $A$ in tensor product is defined by equality*

(4.7.24)      $$b_1 \otimes ... \otimes (f_i(a) \circ b_i) \otimes ... \otimes b_n = f(a) \circ (b_1 \otimes ... \otimes b_i \otimes ... \otimes b_n)$$

PROOF. The equality (4.7.23) follows from the equality (4.7.16) and from the definition (4.7.21). The equality (4.7.24) follows from the equality (4.7.17) and from the definition (4.7.21).                    □

## 4.8. Associativity of Tensor Product

Let $A$ be multiplicative $\Omega_1$-group. Let $B_1$, $B_2$, $B_3$ be $\Omega_2$-algebras. Let, for $k = 1, 2, 3$,

$$f_k : A \longrightarrow\!\!\!\!\!*\longrightarrow B_k$$

be effective representation of multiplicative $\Omega_1$-group $A$ in $\Omega_2$-algebra $B_k$.

LEMMA 4.8.1. *For given value of $x_3 \in B_3$, the map*

(4.8.1)          $h_{12} : (B_1 \otimes B_2) \times B_3 \to B_1 \otimes B_2 \otimes B_3$

*defined by equality*

(4.8.2)          $h_{12}(x_1 \otimes x_2, x_3) = x_1 \otimes x_2 \otimes x_3$

*is reduced morphism of the representation $B_1 \otimes B_2$ into the representation $B_1 \otimes B_2 \otimes B_3$ .*

PROOF. According to the theorem 4.7.9, for given value of $x_3 \in B_3$, the map

(4.8.3)          $(x_1, x_2, x_3) \in B_1 \times B_2 \times B_3 \to x_1 \otimes x_2 \otimes x_3 \in B_1 \otimes B_2 \otimes B_3$

is polymorphism with respect to $x_1 \in B_1$, $x_2 \in B_2$. Therefore, for given value of $x_3 \in B_3$, the lemma follows from the theorem 4.7.10.                    □

LEMMA 4.8.2. *For given value of $x_{12} \in B_1 \otimes B_2$ the map $h_{12}$ is reduced morphism of the representation $B_3$ into the representation $B_1 \otimes B_2 \otimes B_3$ .*

PROOF. According to the theorem 4.7.9 and the equality (4.7.21), for given value of $x_1 \in B_1$, $x_2 \in B_2$, the map

(4.8.4)          $(x_1 \otimes x_2, x_3) \in B_1 \times B_2 \times B_3 \to x_1 \otimes x_2 \otimes x_3 \in B_1 \otimes B_2 \otimes B_3$

is morphism with respect to $x_3 \in B_3$. Therefore, the theorem follows from the equality (4.4.16) and from the theorem 4.5.9.                    □

LEMMA 4.8.3. *There exists reduced morphism of representations*

$$h : (B_1 \otimes B_2) \otimes B_3 \to B_1 \otimes B_2 \otimes B_3$$

PROOF. According to lemmas 4.8.1, 4.8.2 and to the definition 4.4.4, the map $h_{12}$ is reduced polymorphism of representations. The lemma follows from the theorem 4.7.10.                    □

LEMMA 4.8.4. *There exists reduced morphism of representations*

$$g : B_1 \otimes B_2 \otimes B_3 \to (B_1 \otimes B_2) \otimes B_3$$

PROOF. The map

$$(x_1, x_2, x_3) \in B_1 \times B_2 \times B_3 \to (x_1 \otimes x_2) \otimes x_3 \in (B_1 \otimes B_2) \otimes B_3$$

is polymorphism with respect to $x_1 \in B_1$, $x_2 \in B_2$, $x_3 \in B_3$. Therefore, the lemma follows from the theorem 4.7.10.                    □

THEOREM 4.8.5.

(4.8.5)          $(A_1 \otimes A_2) \otimes A_3 = A_1 \otimes (A_2 \otimes A_3) = A_1 \otimes A_2 \otimes A_3$

Proof. According to lemma 4.8.3, there exists reduced morphism of representations

$$h : (B_1 \otimes B_2) \otimes B_3 \to B_1 \otimes B_2 \otimes B_3$$

According to lemma 4.8.4, there exists reduced morphism of representations

$$g : B_1 \otimes B_2 \otimes B_3 \to (B_1 \otimes B_2) \otimes B_3$$

Therefore, reduced morphisms of representations $h$, $g$ are isomorphisms. Therefore, the following equality is true

(4.8.6)                   $(B_1 \otimes B_2) \otimes B_3 = B_1 \otimes B_2 \otimes B_3$

We prove similarly the equality

$$B_1 \otimes (B_2 \otimes B_3) = B_1 \otimes B_2 \otimes B_3$$

$\square$

Remark 4.8.6. *It is evident that structures of $\Omega_2$-algebras $(B_1 \otimes B_2) \otimes B_3$, $B_1 \otimes B_2 \otimes B_3$ are little different. We write down the equality (4.8.6) based on the convention 4.7.3 and this allows us to speak about associativity of tensor product of representations.* $\square$

# Representation of Multiplicative $\Omega$-Group

## 5.1. Representation of Multiplicative $\Omega$-Group

Consistency of product in multiplicative $\Omega$-group $G$ and corresponding transformations of the representation $f$ allows us to consider more details of the representation $f$. However, the construction considered in the theorem 4.5.7 is not complete in case of non commutative product.

If for given representation

$$g : A_1 \longrightarrow A_2$$

for any $A_1$-numbers $a_1$, $b_1$, there is unique $A_1$-number $c_1$ such that

$$f(c_1) = f(a_1) \circ f(b_1)$$

then what format of the product we should choose:

$$(5.1.1) \qquad c_1 = a_1 * b_1$$

or

$$(5.1.2) \qquad c_1 = b_1 * a_1$$

EXAMPLE 5.1.1. *Let*

$$\overline{\overline{e}} = \begin{pmatrix} e_1 & \ldots & e_n \end{pmatrix}$$

*be basis of left vector space $V$ over associative division algebra $A$. We can represent any vector $\overline{v} \in V$ as $*_*$-product of matrices*

$$(5.1.3) \qquad \overline{v} = v *_* \overline{\overline{e}} = \begin{pmatrix} v^1 \\ \ldots \\ v^n \end{pmatrix} *_* \begin{pmatrix} e_1 & \ldots & e_n \end{pmatrix}$$

*where*

$$v = \begin{pmatrix} v^1 \\ \ldots \\ v^n \end{pmatrix}$$

*is matrix of coordinates of the vector $\overline{v}$ with respect to the basis $\overline{\overline{e}}$.*

*We introduce single transitive action of the group $G$ on the basis manifold by the equality*

$$(5.1.4) \qquad g *_* \overline{\overline{e}} = \begin{pmatrix} g_1^1 & \ldots & g_n^1 \\ \ldots & \ldots & \ldots \\ g_1^n & \ldots & g_n^n \end{pmatrix} *_* \begin{pmatrix} e_1 & \ldots & e_n \end{pmatrix}$$

*where we identify G-number g and non-singular matrix*

$$\begin{pmatrix} g_1^1 & ... & g_n^1 \\ ... & ... & ... \\ g_1^n & ... & g_n^n \end{pmatrix}$$

*Action of the group G on the basis manifold is representation, because the following equality is true*

(5.1.5)                    $g_1{}^*{}_*(g_2{}^*{}_*\overline{\overline{e}}) = (g_1{}^*{}_*g_2){}^*{}_*\overline{\overline{e}}$

   *Let*

(5.1.6)                    $v_i = \begin{pmatrix} v_i^1 \\ ... \\ v_i^n \end{pmatrix}$

*be matrix of coordinates of the vector $\overline{v}$ with respect to the basis $\overline{\overline{e}}_i$, $i = 1$, 2, 3. Then*

(5.1.7)                    $\overline{v} = v_1{}^*{}_*\overline{\overline{e}}_1 = v_2{}^*{}_*\overline{\overline{e}}_2 = v_3{}^*{}_*\overline{\overline{e}}_3$

*Let G-number $g_1$ map the basis $\overline{\overline{e}}_1$ into the basis $\overline{\overline{e}}_2$*

(5.1.8)                    $\overline{\overline{e}}_2 = g_1{}^*{}_*\overline{\overline{e}}_1$

*Let G-number $g_2$ map the basis $\overline{\overline{e}}_2$ into the basis $\overline{\overline{e}}_3$*

(5.1.9)                    $\overline{\overline{e}}_3 = g_2{}^*{}_*\overline{\overline{e}}_2$

*The equality*

(5.1.10)                    $\overline{\overline{e}}_3 = (g_2{}^*{}_*g_1){}^*{}_*\overline{\overline{e}}_1$

*follows from the equalities* (5.1.8), (5.1.9). *The equality*

(5.1.11)                    $v_1{}^*{}_*\overline{\overline{e}}_1 = v_2{}^*{}_*g_1{}^*{}_*\overline{\overline{e}}_1$

*follows from the equalities* (5.1.7), (5.1.8). *The equality*

(5.1.12)                    $v_1 = v_2{}^*{}_*g_1$

*follows from the equality* (5.1.11) *because coordinates of vector $\overline{v}$ are unique with respect to basis $\overline{\overline{e}}_1$. The equality*

(5.1.13)                    $v_2 = v_1{}^*{}_*g_1^{-1}$

*follows from the equality* (5.1.12). *Similarly, the equality*

(5.1.14)                    $v_3 = v_2{}^*{}_*g_2^{-1}$

*follows from equalities* (5.1.7), (5.1.9) *and the equality*

(5.1.15)                    $v_3 = v_1{}^*{}_*(g_2{}^*{}_*g_1)^{-1}$

*follows from equalities* (5.1.7), (5.1.10). *The equality*

(5.1.16)                    $v_3 = v_1{}^*{}_*g_1^{-1}{}^*{}_*g_2^{-1}$

*follows from equalities* (5.1.13), (5.1.14).                    $\square$

EXAMPLE 5.1.2. *Let $V$ be left module over ring $D$. It means that we defined representation*

$$f : D \dashrightarrow V \qquad f(d) : v \to d\,v$$

*such that*

$$(d_1 + d_2)v = d_1 v + d_2 v$$
$$d(v_1 + v_2) = dv_1 + dv_2$$
$$d_1(d_2 v) = (d_1 d_2)v$$

*The map*

$$w : V \to D$$

*is called additive if*

$$w(v_1 + v_2) = w(v_1) + w(v_2)$$

*We use notation*

$$(w, v) = w(v)$$

*for image of additive map. We define sum of additive maps by the equality*

$$(w_1 + w_2, v) = (w_1, v) + (w_2, v)$$

*It is easy to show that the set $W$ of additive maps is Abelian group.*

*We define the map*

$$h : D \dashrightarrow W \qquad h(d) : w \to w\,d$$

*using the equality*

$$(wd, v) = (w, dv)$$

*From equalities*

$$\begin{aligned}
((w_1 + w_2)d, v) &= (w_1 + w_2, dv) = (w_1, dv) + (w_2, dv) \\
&= (w_1 d, v) + (w_2 d, v) \\
&= (w_1 d + w_2 d, v) \\
(w(d_1 + d_2), v) &= (w, (d_1 + d_2)v) = (w, d_1 v + d_2 v) \\
&= (w, d_1 v) + (w, d_2 v) = (wd_1, v) + (wd_2, v) \\
&= (wd_1 + wd_2, v) \\
((wd_1)d_2, v) &= (wd_1, d_2 v) = (w, d_1(d_2 v)) = (w, (d_1 d_2)v) \\
&= (w(d_1 d_2), v)
\end{aligned}$$

(5.1.17)

*it follows that the map $h$ is representation of group $G$. However we can write the equality* (5.1.17) *in the following form*

$$((h(d_2) \circ h(d_1)(w)), v) = ((h(d_2)h(d_1)(w)), v) = (h(d_1 d_2)(w), v)$$

*which implies that the map $h$ is not homomorphism of group $G$.* $\qquad\square$

We assume that transformations of representation of multiplicative $\Omega$-group $A_1$ may act on $A_2$-numbers either on the left or on the right. In this case it is sufficient to restrict ourself to the product (5.1.1) in multiplicative $\Omega$-group $A_1$. Thus, the idea of representation of multiplicative $\Omega$-group is that we multiply elements of multiplicative $\Omega$-group in the same order as we multiply transformations of representation. This point of view is reflected in the example 5.1.2. We also see that we need to change notation before we can use this point of view. Instead of considering $f \in \text{End}(\Omega_2; A_2)$ as map

$$f : a_2 \in A_2 \to f(a_2) \in A_2$$

we must consider an endomorphism $f$ as operator.

DEFINITION 5.1.3. *Let* $\text{End}(\Omega_2, A_2)$ *be a multiplicative $\Omega$-group with product*[5.1]

$$(f, g) \to f \bullet g$$

*Let an endomorphism $f$ act on $A_2$-number $a$ on the left. We will use notation*

(5.1.18)                                  $$f(a_2) = f \bullet a_2$$

*Let $A_1$ be multiplicative $\Omega$-group with product*

$$(a, b) \to a * b$$

*We call a homomorphism of multiplicative $\Omega$-group*

(5.1.19)                                  $$f : A_1 \to \text{End}(\Omega_2, A_2)$$

**left-side representation** *of multiplicative $\Omega$-group $A_1$ or* **left-side $A_1$-representation** *in $\Omega_2$-algebra $A_2$ if the map $f$ holds*

(5.1.20)                          $$f(a_1 * b_1) \bullet a_2 = (f(a_1) \bullet f(b_1)) \bullet a_2$$

*We identify an $A_1$-number $a_1$ and its image $f(a_1)$ and write left-side transformation caused by $A_1$-number $a_1$ as*

$$a_2' = f(a_1) \bullet a_2 = a_1 * a_2$$

*In this case, the equality* (5.1.20) *gets following form*

(5.1.21)                              $$f(a_1 * b_1) \bullet a_2 = (a_1 * b_1) * a_2$$

*The map*

$$(a_1, a_2) \in A_1 \times A_2 \to a_1 * a_2 \in A_2$$

*generated by left-side representation $f$ is called* **left-side product** *of $A_2$-number $a_2$ over $A_1$-number $a_1$.*                                                                    $\square$

Let

$$f : A_2 \to A_2$$
$$g : A_2 \to A_2$$

be endomorphisms of $\Omega_2$-algebra $A_2$. Let product in multiplicative $\Omega$-group $\text{End}(\Omega_2, A_2)$ is composition of endomorphisms. Since the product of maps $f$ and $g$ is defined in the same order as these maps act on $A_2$-number, then we consider the equality

(5.1.22)                                  $$(f \circ g) \circ a = f \circ (g \circ a)$$

as **associative law**. This allows writing of equality (5.1.22) without using of brackets

$$f \circ g \circ a = f \circ (g \circ a) = (f \circ g) \circ a$$

---

[5.1] Very often a product in multiplicative $\Omega$-group $\text{End}(\Omega_2, A_2)$ is superposition of endomorphisms

$$f \bullet g = f \circ g$$

However, as we see in the example 5.2.5, a product in multiplicative $\Omega$-group $\text{End}(\Omega_2, A_2)$ may be different from superposition of endomorphisms. According to the definition 4.6.13, we can consider two products in universal algebra $A$.

as well it allows writing of equality (5.1.20) in the following form

$$(5.1.23) \qquad f(a_1 * b_1) \circ a_2 = f(a_1) \circ f(b_1) \circ a_2$$

From the equality (5.1.21), it follows that

$$(5.1.24) \qquad (a_1 * b_1) * a_2 = a_1 * (b_1 * a_2)$$

We consider the equality (5.1.24) as **associative law**.

REMARK 5.1.4. *Let the map*

$$f : A_1 \longrightarrow\!\ast\!\longrightarrow A_2$$

*be the left-side representation of multiplicative $\Omega$-group $A_1$ in $\Omega_2$-algebra $A_2$. Let the map*

$$g : B_1 \longrightarrow\!\ast\!\longrightarrow B_2$$

*be the left-side representation of multiplicative $\Omega$-group $B_1$ in $\Omega_2$-algebra $B_2$. Let the map*

$$(r_1 : A_1 \to B_1, \ r_2 : A_2 \to B_2)$$

*be morphism of representations. We use notation*

$$r_2(a_2) = r_2 \circ a_2$$

*for image of $A_2$-number $a_2$ with respect to the map $r_2$. Then we can write the equality (3.2.3) in the following form*

$$r_2 \circ (a_1 * a_2) = r_1(a_1) * (r_2 \circ a_2)$$

$\square$

DEFINITION 5.1.5. *Let* $\mathrm{End}(\Omega_2, A_2)$ *be a multiplicative $\Omega$-group with product*[5.2]

$$(f, g) \to f \bullet g$$

*Let an endomorphism $f$ act on $A_2$-number $a$ on the right. We will use notation*

$$(5.1.25) \qquad f(a_2) = a_2 \bullet f$$

*Let $A_1$ be multiplicative $\Omega$-group with product*

$$(a, b) \to a * b$$

*We call a homomorphism of multiplicative $\Omega$-group*

$$(5.1.26) \qquad f : A_1 \to \mathrm{End}(\Omega_2, A_2)$$

**right-side representation** *of multiplicative $\Omega$-group $A_1$ or* **right-side $A_1$-representation** *in $\Omega_2$-algebra $A_2$ if the map $f$ holds*

$$(5.1.27) \qquad a_2 \bullet f(a_1 * b_1) = a_2 \bullet (f(a_1) \bullet f(b_1))$$

*We identify an $A_1$-number $a_1$ and its image $f(a_1)$ and write right-side transformation caused by $A_1$-number $a_1$ as*

$$a_2' = a_2 \bullet f(a_1) = a_2 * a_1$$

---

[5.2] Very often a product in multiplicative $\Omega$-group $\mathrm{End}(\Omega_2, A_2)$ is superposition of endomorphisms

$$f \bullet g = f \circ g$$

However, as we see in the example 5.2.5, a product in multiplicative $\Omega$-group $\mathrm{End}(\Omega_2, A_2)$ may be different from superposition of endomorphisms. According to the definition 4.6.13, we can consider two products in universal algebra $A$.

*In this case, the equality* (5.1.27) *gets following form*

(5.1.28) $$a_2 \bullet f(a_1 * b_1) = a_2 * (a_1 * b_1)$$

*The map*

$$(a_1, a_2) \in A_1 \times A_2 \to a_2 * a_1 \in A_2$$

*generated by right-side representation* $f$ *is called* **right-side product** *of* $A_2$-*number* $a_2$ *over* $A_1$-*number* $a_1$. $\qquad\square$

Let

$$f : A_2 \to A_2$$
$$g : A_2 \to A_2$$

be endomorphisms of $\Omega_2$-algebra $A_2$. Let product in multiplicative $\Omega$-group $\operatorname{End}(\Omega_2, A_2)$ is composition of endomorphisms. Since the product of maps $f$ and $g$ is defined in the same order as these maps act on $A_2$-number, then we consider the equality

(5.1.29) $$a \circ (g \circ f) = (a \circ g) \circ f$$

as **associative law**. This allows writing of equality (5.1.29) without using of brackets

$$a \circ g \circ f = (a \circ g) \circ f = a \circ (g \circ f)$$

as well it allows writing of equality (5.1.27) in the following form

(5.1.30) $$a_2 \circ f(a_1 * b_1) = a_2 \circ f(a_1) \circ f(b_1)$$

From the equality (5.1.28), it follows that

(5.1.31) $$a_2 * (a_1 * b_1) = (a_2 * a_1) * b_1$$

We consider the equality (5.1.31) as **associative law**.

REMARK 5.1.6. *Let the map*

$$f : A_1 \longrightarrow\!\!\!* \longrightarrow A_2$$

*be the left-side representation of multiplicative* $\Omega$-*group* $A_1$ *in* $\Omega_2$-*algebra* $A_2$. *Let the map*

$$g : B_1 \longrightarrow\!\!\!* \longrightarrow B_2$$

*be the left-side representation of multiplicative* $\Omega$-*group* $B_1$ *in* $\Omega_2$-*algebra* $B_2$. *Let the map*

$$(r_1 : A_1 \to B_1, \ r_2 : A_2 \to B_2)$$

*be morphism of representations. We use notation*

$$r_2(a_2) = r_2 \circ a_2$$

*for image of* $A_2$-*number* $a_2$ *with respect to the map* $r_2$. *Then we can write the equality* (3.2.3) *in the following form*

$$r_2 \circ (a_2 * a_1) = (r_2 \circ a_2) * r_1(a_1)$$

$\qquad\square$

If multiplicative $\Omega$-group $A_1$ is Abelian, then there is no difference between left-side and right-side representations.

DEFINITION 5.1.7. *Let $A_1$ be Abelian multiplicative $\Omega$-group. We call a homomorphism of multiplicative $\Omega$-group*

$$(5.1.32) \qquad\qquad f : A_1 \to \operatorname{End}(\Omega_2, A_2)$$

**representation** *of multiplicative $\Omega$-group $A_1$ or $A_1$-**representation** in $\Omega_2$-algebra $A_2$ if the map $f$ holds*

$$(5.1.33) \qquad\qquad f(a_1 * b_1) \bullet a_2 = (f(a_1) \bullet f(b_1)) \bullet a_2$$

$\square$

Usually we identify a representation of the Abelian multiplicative $\Omega$-group $A_1$ and a left-side representation of the multiplicative $\Omega$-group $A_1$. However, if it is necessary for us, we identify a representation of the Abelian multiplicative $\Omega$-group $A_1$ and a right-side representation of the multiplicative $\Omega$-group $A_1$.

From the analysis of the example 5.1.2, it follows that choice between left-side and right-side representation depends from considered model. Since left-side representation and right-side representation are based on homomorphism of $\Omega$-group, then the following statement is true.

THEOREM 5.1.8 (duality principle for representation of multiplicative $\Omega$-group). *Any statement which holds for left-side representation of multiplicative $\Omega$-group $A_1$ holds also for right-side representation of multiplicative $\Omega$-group $A_1$, if we will use right-side product over $A_1$-number $a_1$ instead of left-side product over $A_1$-number $a_1$.* $\square$

REMARK 5.1.9. *If $\Omega_1$-algebra is not multiplicative $\Omega$-group, then we cannot tell whether representation acts on left or on right. In this case we continue to us functional notation $f(a_1)(a_2)$ for representation of $\Omega_1$-algebra.* $\square$

From the analysis of equalities (5.1.15), (5.1.16), it follows that the action of the group $G$ on the set of coordinates of the vector $\overline{v}$ (the example 5.1.1) does not correspond to either left-side or right-side representation. it follows that we have two choices. We accept that in multiplicative $\Omega$-group $A_1$ we can define both types of product ((5.1.1) and (5.1.2)) in order to coordinate product in multiplicative $\Omega$-group $A_1$ and product of transformations of representation of multiplicative $\Omega$-group $A_1$. This point of view is reflected in definitions 5.1.10, 5.1.11.

DEFINITION 5.1.10. *Left-side representation*

$$f : A_1 \longrightarrow_* A_2$$

*is called **covariant** if the equality*

$$a_1 * (b_1 * a_2) = (a_1 * b_1) * a_2$$

*is true.* $\square$

DEFINITION 5.1.11. *Left-side representation*

$$f : A_1 \longrightarrow_* A_2$$

*is called **contravariant** if the equality*

$$(5.1.34) \qquad\qquad a_1^{-1} * (b_1^{-1} * a_2) = (b_1 * a_1)^{-1} * a_2$$

*is true.* $\square$

If type of representation is not specified, then we assume that the representation is covariant. From equalities (5.1.15), (5.1.16), it follows that the action of the group $G$ on the set of coordinates of the vector $\overline{v}$ (the example 5.1.1) is contravariant right-side representation.

How big is the difference between covariant and contravariant representations. Since

$$(b_1 * a_1)^{-1} = a_1^{-1} * b_1^{-1}$$

then the equality

(5.1.35)           $$a_1^{-1} * (b_1^{-1} * a_2) = (a_1^{-1} * b_1^{-1}) * a_2$$

follows from the equality (5.1.34). From the equality (5.1.35), it follows that we can consider contravariant representation of the group $G$ as covariant representation of the group $G$, generated by $G$-numbers of the form $a^{-1}$. The same way as in the example 5.1.1, we consider two coordinated representations of the group $G$

$$f : G \longmapsto A_2$$

$$h : G \longmapsto B_2$$

moreover $G$-number $g$ generates the transformation

$$a_1 \in G : a_2 \in A_2 \to a_1 * a_2 \in A_2$$

in the universal algebra $A_2$ and the transformation

$$a_1 \in G : b_2 \in B_2 \to a_1^{-1} * b_2 \in B_2$$

in the universal algebra $B_2$.

## 5.2. Left and Right Shifts

THEOREM 5.2.1. *The product*

$$(a, b) \to a * b$$

*in multiplicative $\Omega$-group $A$ determines two different representations.*

- *the* **left shift**

$$a' = L(b) \circ a = b * a$$

  *is left-side representation of multiplicative $\Omega$-group $A$ in $\Omega$-algebra $A$*

(5.2.1)           $$L(c * b) = L(c) \circ L(b)$$

- *the* **right shift**

$$a' = a \circ R(b) = a * b$$

  *is right-side representation of multiplicative $\Omega$-group $A$ in $\Omega$-algebra $A$*

(5.2.2)           $$R(b * c) = R(b) \circ R(c)$$

PROOF. According to the definition 4.5.4, left and right shifts are endomorphisms of $\Omega$-algebra $A$. According to the definition 4.5.4, we can define $\Omega$-algebra on the set of left shifts. According to the definition of multiplicative group, [5.3] the equality $a_1 = a_2$ follows from the equality

$$L(a_1) \circ x = a_1 * x = a_2 * x = L(a_2) \circ x$$

---

[5.3] See, for instance, definition on pages [2]-3, [2]-7.

for any $x$. In particular, the equality (5.2.1) follows from the equality

$$L(c * b) \circ a = (c * b) * a = c * (b * a) = L(c) \circ (L(b) \circ a) = L(c) \circ L(b) \circ a$$

Therefore, the map

$$a \in A \to L(a)$$

is left-side representation of multiplicative $\Omega$-group $A$ in $\Omega$-algebra $A$. Similar reasoning is true for right shift. $\qquad \square$

Associative $D$-algebra is multiplicative $\Omega$-group. Non associative $D$-algebra $A$ is not $\Omega$-group, because $A$ is groupoid with respect to product. However we also study representation of non associative $D$-algebra.

Definition 5.2.2. *Let product*

$$c_1 = a_1 * b_1$$

*be operation of $\Omega_1$-algebra $A$. Let $\Omega = \Omega_1 \setminus \{*\}$. If $\Omega_1$-algebra $A$ is groupoid with respect to product and, for any operation $\omega \in \Omega(n)$, the product is distributive over the operation $\omega$*

$$a * (b_1...b_n\omega) = (a * b_1)...(a * b_n)\omega$$

$$(b_1...b_n\omega) * a = (b_1 * a)...(b_n * a)\omega$$

*then $\Omega_1$-algebra $A$ is called $\Omega$-**groupoid**.* $\qquad \square$

We will use the same notation for representation of $\Omega$-groupoid as we use for representation of multiplicative $\Omega$-group.

Theorem 5.2.3. *The product in non associative $\Omega$-groupoid $A$ determines two different representations.*

- *The left shift*

$$a' = L(b) \circ a = b * a$$

 *is representation of $\Omega$-algebra $A$ in $\Omega$-algebra $A$.*

- *The right shift*

$$a' = a \circ R(b) = a * b$$

 *is representation of $\Omega$-algebra $A$ in $\Omega$-algebra $A$.*

Proof. According to the definition 4.5.4, left and right shifts are endomorphisms of $\Omega$-algebra $A$. According to the definition 5.2.2, we can define $\Omega$-algebra on the set of left shifts. Therefore, the map

$$a \in A \to L(a)$$

is representation of $\Omega$-algebra $A$ in $\Omega$-algebra $A$. $\qquad \square$

Theorem 5.2.4. *Let*

$$L : A \longrightarrow\!\!\!* A$$

*be representation of non associative $\Omega$-groupoid $A$ in $\Omega$-algebra $A$. Then, on the set $\mathrm{End}(\Omega, A)$, there exists product which is different from the superposition of endomorphisms.*

PROOF. Consider the map

$$L : A \to \mathrm{End}(\Omega, A) \quad L(a) : b \to ab$$

Since the product in Ω-groupoid $A$ is not associative, then, in general

$$L(a) \circ (L(b) \circ c) = a * (b * c) \neq (a * b) * c = L(a * b) \circ c$$

Therefore, $L(ab) \neq L(a) \circ L(b)$.  □

According to the theorem 5.2.1, if $A$ is multiplicative Ω-group, then the equality (5.2.1) ensures that left shift generates is left-side representation of multiplicative Ω-group $A$ in Ω-algebra $A$. According to the theorem 5.2.4 this equality is not true in non associative Ω-groupoid $A$. However theorems 5.2.3, 5.2.4 do not answer the question about the possibility of consideration of left-side representation of non associative Ω-groupoid $A$ in Ω-algebra $A$. According to the example 5.2.5, there exists posibility of such representation, even product in Ω-groupoid is non associative.

EXAMPLE 5.2.5. *Let $A$ be Lie algebra. The product*[5.4] *$[a, b]$ of $A$-numbers $a$, $b$ satisfies to the equality*

$$(5.2.3) \qquad\qquad [a, b] = -[b, a]$$

*and to the Lee identity*

$$(5.2.4) \qquad\qquad [c, [b, a]] + [b, [a, c]] + [a, [c, b]] = 0$$

*We define left shift on Lie algebra $A$ by the equality*

$$(5.2.5) \qquad\qquad L(b) \circ a = [b, a]$$

*From the equality* (5.2.5), *it follows that*

$$(5.2.6) \qquad L(c) \circ L(b) \circ a = L(c) \circ (L(b) \circ a) = [c, [b, a]]$$

*The equality*

$$(5.2.7) \qquad \begin{aligned} L(c) \circ L(b) \circ a - L(b) \circ L(c) \circ a &= [c, [b, a]] - [b, [c, a]] \\ &= [c, [b, a]] + [b, [a, c]] \end{aligned}$$

*follows from equalities* (5.2.3), (5.2.6). *The equality*

$$(5.2.8) \qquad\qquad [c, [b, a]] + [b, [a, c]] = -[a, [c, b]] = [[c, b], a]$$

*follows from equalities* (5.2.3), (5.2.4). *The equality*

$$(5.2.9) \qquad L(c) \circ L(b) \circ a - L(b) \circ L(c) \circ a = L([c, b]) \circ a$$

*follows from equalities* (5.2.5), (5.2.7), (5.2.8).

*If I define Lie product*

$$[L(c), L(b)] \circ a = L(c) \circ L(b) \circ a - L(b) \circ L(c) \circ a$$

*on the set of left shifts then the equality* (5.2.9) *gets the form*

$$(5.2.10) \qquad\qquad [L(c), L(b)] \circ a = L([c, b]) \circ a$$

*Therefore, Lie algebra $A$ with product $[a, b]$ generates representation in vector space $A$.*  □

---

[5.4] See definition [17]-1 on the page 3.

### 5.3. Orbit of Representation of Multiplicative $\Omega$-Group

THEOREM 5.3.1. *Let the map*

$$f : A_1 \longrightarrow\!\!*\!\!\longrightarrow A_2$$

*be the left-side representation of multiplicative $\Omega_1$-group $A_1$. and $e$ be unit of multiplicative $\Omega_1$-group $A_1$. Then*

$$f(e) = \delta$$

*where $\delta$ is identity transformation of $\Omega_2$-algebra $A_2$.*

PROOF. The theorem follows from the equality

$$f(a) = f(a * e) = f(a) \circ f(e)$$

for any $A_1$-number. $\qquad\square$

THEOREM 5.3.2. *Let the map*

$$g : A_1 \longrightarrow\!\!*\!\!\longrightarrow A_2$$

*be the left-side representation of multiplicative $\Omega$-group $A_1$. For any $g \in A_1$ transformation has inverse map and satisfies the equality*

(5.3.1) $$f(g^{-1}) = f(g)^{-1}$$

PROOF. Let $e$ be unit of multiplicative $\Omega$-group $A_1$ and $\delta$ be identity transformation of the set $A_2$. Based on (5.1.20) and the theorem 5.3.1, we have

$$u = \delta \circ u = f(gg^{-1}) \circ u = f(g) \circ f(g^{-1}) \circ u$$

This completes the proof. $\qquad\square$

DEFINITION 5.3.3. *Let $A_1$ be $\Omega$-groupoid with product*

$$(a, b) \to a * b$$

*Let the map*

$$f : A_1 \longrightarrow\!\!*\!\!\longrightarrow A_2$$

*be the left-side representation of $\Omega$-groupoid $A_1$ in $\Omega_2$-algebra $A_2$. For any $a_2 \in A_2$, we define* **orbit of representation** *of the $\Omega$-groupoid $A_1$ as set*

$$A_1 * a_2 = \{b_2 = a_1 * a_2 : a_1 \in A_1\}$$

$\qquad\square$

DEFINITION 5.3.4. *Let $A_1$ be $\Omega$-groupoid with product*

$$(a, b) \to a * b$$

*Let the map*

$$f : A_1 \longrightarrow\!\!*\!\!\longrightarrow A_2$$

*be the left-side representation of $\Omega$-groupoid $A_1$ in $\Omega_2$-algebra $A_2$. For any $a_2 \in A_2$, we define* **orbit of representation** *of the $\Omega$-groupoid $A_1$ as set*

$$a_2 * A_1 = \{b_2 = a_2 * a_1 : a_1 \in A_1\}$$

$\qquad\square$

Theorem 5.3.5. *Let the map*

$$f : A_1 \longrightarrow_* A_2$$

*be the left-side representation of multiplicative $\Omega$-group $A_1$. Then $a_2 \in A_1 * a_2$.*

Proof. According to the theorem 5.3.1,

$$a_2 = e * a_2 = f(e) \circ a_2$$

$\square$

Theorem 5.3.6. *Let*

$$L : A \longrightarrow_* A$$

*be representation of Lie algebra generated by the set of left shifts. Then $a \notin [A, a]$.*

Proof. The theorem follows from absence of unit in Lie algebra. Besides, the set of vectors of three dimensional space where we defined cross product is the most simple example of Lie algebra. It is evident that there no exist vector $b$ such that

$$a = b \times a$$

$\square$

Theorem 5.3.7. *Let the map*

$$f : A_1 \longrightarrow_* A_2$$

*be the left-side representation of multiplicative $\Omega$-group $A_1$. Let*

(5.3.2)                              $b_2 \in A_1 * a_2$

*Then*

(5.3.3)                              $A_1 * a_2 = A_1 * b_2$

Proof. From (5.3.2) it follows that there exists $a_1 \in A_1$ such that

(5.3.4)                              $b_2 = a_1 * a_2$

Let $c_2 \in A_1 * b_2$. Then there exists $b_1 \in A_1$ such that

(5.3.5)                              $c_2 = b_1 * b_2$

If we substitute (5.3.4) into (5.3.5) we get

(5.3.6)                              $c_2 = b_1 * a_1 * a_2$

Based (5.1.20), we see that from (5.3.6) it follows that $c_2 \in A_1 * a_2$. Thus

(5.3.7)                              $A_1 * b_2 \subseteq A_1 * a_2$

Based (5.3.1), we see that, from (5.3.4), it follows that

(5.3.8)                              $a_2 = a_1^{-1} * b_2$

From (5.3.8) it follows that $a_2 \in A_1 * b_2$ and therefore

(5.3.9)                              $A_1 * a_2 \subseteq A_1 * b_2$

The equality (5.3.3) follows from statements (5.3.7), (5.3.9).           $\square$

Thus, the left-side representation of multiplicative $\Omega$-group $A_1$ in $\Omega_2$-algebra $A_2$ forms equivalence $S$ and the orbit $A_1 * a_2$ is equivalence class. We will use notation $A_2/A_1$ for quotient set $A_2/S$ and this set is called **space of orbits of left-side representation** $f$.

## 5.4. Representation in $\Omega$-Group

THEOREM 5.4.1. *We call kernel of inefficiency of left-side representation of multiplicative $\Omega$-group $A_1$ in $\Omega_2$-algebra $A_2$ a set*

$$K_f = \{a_1 \in A_1 : f(a_1) = \delta\}$$

*A kernel of inefficiency of left-side representation is a subgroup of the multiplicative group $A_1$.*

PROOF. Assume $f(a_1) = \delta$ and $f(a_2) = \delta$. Then

$$f(a_1 * a_2) = f(a_1) \bullet f(a_2) = \delta$$
$$f(a_1^{-1}) = (f(a_1))^{-1} = \delta$$

$\square$

THEOREM 5.4.2. *Left-side representation of multiplicative $\Omega$-group $A_1$ in $\Omega_2$-algebra $A_2$ is* **effective** *iff kernel of inefficiency $K_f = \{e\}$.*

PROOF. Statement follows from the definitions 3.1.2 and from the theorem 5.4.1. $\square$

THEOREM 5.4.3. *If a representation*

$$f : A_1 \longrightarrow\!\!\!* \longrightarrow A_2$$

*of multiplicative $\Omega$-group $A_1$ in $\Omega_2$-algebra $A_2$ is not effective we can switch to the effective representation replacing the multiplicative $\Omega$-group $A_1$ by the multiplicative $\Omega$-group $A_1' = A_1/K_f$.*

PROOF. Let the operation $\omega \in \Omega(n)$. To prove the theorem, we need to show that the equality

(5.4.1)                    $f(a_1...a_n\omega) = f(b_1...b_n\omega)$

follows from the statement $f(a_1) = f(b_1), ..., f(a_n) = f(b_n)$. Indeed, the equality (5.4.1) follows from the equality

$$f(a_1...a_n\omega) = f(a_1)...f(a_n)\omega = f(b_1)...f(b_n)\omega = f(b_1...b_n\omega)$$

$\square$

The theorem 5.4.3 means that we can study only an effective action.

## 5.5. Single Transitive Right-Side Representation of Group

THEOREM 5.5.1. *Let the map*

$$g : A_1 \longrightarrow\!\!\!* \longrightarrow A_2$$

*be the left-side representation of multiplicative $\Omega$-group $A_1$ in $\Omega_2$-algebra $A_2$. A* **little group** *or* **stability group** *of $a_2 \in A_2$ is the set*

$$A_{1a_2} = \{a_1 \in A_1 : a_1 * a_2 = a_2\}$$

*The representation $f$ is* **free***, iff, for any $a_2 \in A_2$, stability group $A_{1a_2} = \{e\}$.*

PROOF. According to the definition 3.1.4, the representation $f$ is free iff the statement

$$(5.5.1) \qquad\qquad f(a_1) = f(b_1)$$

implies the equality $a_1 = b_1$. The equality (5.5.1) is equivalent to the equality

$$(5.5.2) \qquad\qquad f(b_1^{-1} * a_1) = \delta$$

The statement (5.5.2) implies the equality $a_1 = b_1$ iff, for any $a_2 \in A_2$, stability group $A_{1a_2} = \{e\}$. $\qquad\square$

THEOREM 5.5.2. *Let the map*

$$f : A_1 \longrightarrow\!\!\!* A_2$$

*be the free left-side representation of multiplicative $\Omega$-group $A_1$ in $\Omega_2$-algebra $A_2$. There exist $1 - 1$ correspondence between any two orbits of representation, as well between any orbit of representation and multiplicative $\Omega$-group $A_1$.*

PROOF. Given $a_2 \in A_2$ there exist $a_1, b_1 \in A_1$

$$(5.5.3) \qquad\qquad a_1 * a_2 = b_1 * a_2$$

We multiply both parts of equation (5.5.3) by $a_1^{-1}$

$$a_2 = a_1^{-1} * b_1 * a_2$$

Since the representation is free, $a_1 = b_1$. Since we established $1 - 1$ correspondence between orbit and multiplicative $\Omega$-group $A_1$, we proved the statement of the theorem. $\qquad\square$

THEOREM 5.5.3. *Left-side representation*

$$g : A_1 \longrightarrow\!\!\!* A_2$$

*of multiplicative $\Omega$-group $A_1$ in $\Omega_2$-algebra $A_2$ is* **single transitive** *iff, for any $a_2$, $b_2 \in A_2$, exists one and only one $a_1 \in A_1$ such that $a_2 = a_1 * b_2$.*

PROOF. Corollary of definitions 3.1.2 and 3.1.8. $\qquad\square$

THEOREM 5.5.4. *If there exists single transitive representation*

$$f : A_1 \longrightarrow\!\!\!* A_2$$

*of multiplicative $\Omega$-group $A_1$ in $\Omega_2$-algebra $A_2$, then we can uniquely define coordinates on $A_2$ using $A_1$-numbers.*

*If $f$ is left-side single transitive representation then $f(a)$ is equivalent to the left shift $L(a)$ on the group $A_1$. If $f$ is right-side single transitive representation then $f(a)$ is equivalent to the right shift $R(a)$ on the group $A_1$.*

PROOF. Let $f$ be left-side single transitive representation. We select $A_2$-number $a_2$ and define coordinates of $A_2$-number $b_2$ as $A_1$-number $a_1$ such that

$$b_2 = a_1 * a_2 = (a_1 * e) * a_2 = (L(a_1) \circ e) * a_2$$

Coordinates defined this way are unique up to choice of $A_2$-number $a_2$ because the action is effective. For left-side single transitive representation, we also use notation

$$f(a_1) \bullet a_2 = L(a_1) \circ a_2 = (L(a_1) \circ e) * a_2$$

We use notation $L(a_1) \circ a_2$ for left-side single transitive representation $f$ because, according to the theorem 5.2.1, product of left shifts equals their composition.

Let $f$ be right-side single transitive representation. We select $A_2$-number $a_2$ and define coordinates of $A_2$-number $b_2$ as $A_1$-number $a_1$ such that

$$b_2 = a_2 * a_1 = a_2 * (e * a_1) = a_2 * (e \circ R(a_1))$$

Coordinates defined this way are unique up to choice of $A_2$-number $a_2$ because the action is effective. For rigt-side single transitive representation, we also use notation

$$a_2 \bullet f(a_1) = a_2 \circ R(a_1) = a_2 * (e \circ R(a_1))$$

We use notation $a_2 \circ R(a_1)$ for rigt-side single transitive representation $f$ because, according to the theorem 5.2.1, product of rigt shifts equals their composition.   □

DEFINITION 5.5.5. *We call* $\Omega_2$-*algebra* $A_2$ **homogeneous space** *of multiplicative* $\Omega$-*group* $A_1$ *if there exists single transitive left-side representation*

$$f : A_1 \xrightarrow{\quad * \quad} A_2$$

   □

THEOREM 5.5.6. *Free left-side representation of multiplicative* $\Omega$-*group* $A_1$ *in* $\Omega_2$-*algebra* $A_2$ *is single transitive representation on orbit.*

PROOF. The theorem follows from the theorem 5.5.2.   □

THEOREM 5.5.7. *Left and right shifts on multiplicative* $\Omega$-*group* $A_1$ *are commuting.*

PROOF. The theorem follows from the associativity of product on multiplicative $\Omega$-group $A_1$

$$(L(a) \circ c) \circ R(b) = (a * c) * b = a * (c * b) = L(a) \circ (c \circ R(b))$$

   □

Theorem 5.5.7 can be phrased n the following way.

THEOREM 5.5.8. *Let* $A_1$ *be multiplicative* $\Omega$-*group. For any* $a_1 \in A_1$, *the map* $L(a_1)$ *is automorphism of representation* $R$.

PROOF. According to theorem 5.5.7

$$(5.5.4) \qquad\qquad L(a_1) \circ R(b_1) = R(b_1) \circ L(a_1)$$

Equation (5.5.4) coincides with equation (3.2.2) from definition 3.2.2 when $r_1 = id$, $r_2 = L(a_1)$.   □

THEOREM 5.5.9. *Let left-side* $A_1$-*representation* $f$ *on* $\Omega_2$-*algebra* $A_2$ *be single transitive. Then we can uniquely define a single transitive right-side* $A_1$-*representation* $h$ *on* $\Omega_2$-*algebra* $A_2$ *such that diagram*

$$
\begin{array}{ccc}
A_2 & \xrightarrow{\ h(a_1)\ } & A_2 \\
{\scriptstyle f(b_1)}\downarrow & & \downarrow{\scriptstyle f(b_1)} \\
A_2 & \xrightarrow{\ h(a_1)\ } & A_2
\end{array}
$$

*is commutative for any* $a_1$, $b_1 \in A_1$.[5.5]

---

[5.5]You can see this statement in [4].

PROOF. We use group coordinates for $A_2$-numbers $a_2$. Then according to theorem 5.5.4 we can write the left shift $L(a_1)$ instead of the transformation $f(a_1)$.

Let $a_2, b_2 \in A_2$. Then we can find one and only one $a_1 \in A_1$ such that

$$b_2 = a_2 * a_1 = a_2 \circ R(a_1)$$

We assume

$$h(a) = R(a)$$

For some $b_1 \in A_1$, we have

$$c_2 = f(b_1) \bullet a_2 = L(b_1) \circ a_2 \quad d_2 = f(b_1) \bullet b_2 = L(b_1) \circ b_2$$

According to the theorem 5.5.7, the diagram

(5.5.5)
$$
\begin{array}{ccc}
a_2 & \xrightarrow{\ h(a_1)=R(a_1)\ } & b_2 \\
{\scriptstyle f(b_1)=L(b_1)}\big\downarrow & & \big\downarrow{\scriptstyle f(b_1)=L(b_1)} \\
c_2 & \xrightarrow{\ h(a_1)=R(a_1)\ } & d_2
\end{array}
$$

is commutative.

Changing $b_1$ we get that $c_2$ is an arbitrary $A_2$-number.

We see from the diagram that if $a_2 = b_2$ then $c_2 = d_2$ and therefore $h(e) = \delta$. On other hand if $a_2 \neq b_2$ then $c_2 \neq d_2$ because the left-side $A_1$-representation $f$ is single transitive. Therefore the right-side $A_1$-representation $h$ is effective.

In the same way we can show that for given $c_2$ we can find $a_1$ such that $d_2 = c_2 \bullet h(a_1)$. Therefore the right-side $A_1$-representation $h$ is single transitive.

In general the product of transformations of the left-side $A_1$-representation $f$ is not commutative and therefore the right-side $A_1$-representation $h$ is different from the left-side $A_1$-representation $f$. In the same way we can create a left-side $A_1$-representation $f$ using the right-side $A_1$-representation $h$. $\qquad \square$

Representations $f$ and $h$ are called **twin representations** of the multiplicative $\Omega$-group $A_1$.

REMARK 5.5.10. *It is clear that transformations $L(a)$ and $R(a)$ are different until the multiplicative $\Omega$-group $A_1$ is nonabelian. However they both are maps onto. Theorem 5.5.9 states that if both right and left shift presentations exist on the set $A_2$, then we can define two commuting representations on the set $A_2$. The right shift or the left shift only cannot represent both types of representation. To understand why it is so let us change diagram (5.5.5) and assume*

$$h(a_1) \bullet a_2 = L(a_1) \circ a_2 = b_2$$

*instead of*

$$a_2 \bullet h(a_1) = a_2 \circ R(a_1) = b_2$$

*and let us see what expression $h(a_1)$ has at the point $c_2$. The diagram*

$$
\begin{array}{ccc}
a_2 & \xrightarrow{\ h(a_1)=L(a_1)\ } & b_2 \\
{\scriptstyle f(b_1)=L(b_1)}\big\downarrow & & \big\downarrow{\scriptstyle f(b_1)=L(b_1)} \\
c_2 & \xrightarrow{\quad h(a_1)\quad } & d_2
\end{array}
$$

*is equivalent to the diagram*

$$
\begin{array}{ccc}
a_2 & \xrightarrow{\;\;h(a_1)=L(a_1)\;\;} & b_2 \\[2pt]
\Big\uparrow{\scriptstyle (f(b_1))^{-1}=L(b_1^{-1})} & & \Big\downarrow{\scriptstyle f(b_1)=L(b_1)} \\[2pt]
c_2 & \xrightarrow{\;\;h(a_1)\;\;} & d_2
\end{array}
$$

*and we have*  $d_2 = b_1 b_2 = b_1 a_1 a_2 = b_1 a_1 b_1^{-1} c_2$.   *Therefore*

$$h(a_1) \bullet c_2 = (b_1 a_1 b_1^{-1}) c_2$$

*We see that the representation of h depends on its argument.* □

Theorem 5.5.11. *Let f and h be twin representations of the multiplicative Ω-group $A_1$. For any  $a_1 \in A_1$  the map $h(a_1)$ is automorphism of representation f.*

Proof. The statement of theorem is corollary of theorems 5.5.8 and 5.5.9.   □

Question 5.5.12. *Is there a morphism of representations from L to L different from automorphism $R(a_1)$? If we assume*

$$r_1(a_1) = c_1 a_1 c_1^{-1}$$

$$r_2(a_1) \circ a_2 = c_1 a_2 a_1 c_1^{-1}$$

*then it is easy to see that the map  $(r_1 \;\; r_2(a_1))$  is morphism of the representations from L to L.  However this map is not automorphism of the representation L, because $r_1 \neq \mathrm{id}$.* □

# Basis of Representation of Universal Algebra

## 6.1. Generating Set of Representation

DEFINITION 6.1.1. *Let*

$$f : A_1 \longrightarrow_* A_2$$

*be representation of $\Omega_1$-algebra $A_1$ in $\Omega_2$-algebra $A_2$. The set $B_2 \subset A_2$ is called* **stable set of representation** *$f$, if $f(a)(m) \in B_2$ for each $a \in A_1$, $m \in B_2$.* □

We also say that the set $A_2$ is stable with respect to the representation $f$.

THEOREM 6.1.2. *Let*

$$f : A_1 \longrightarrow_* A_2$$

*be representation of $\Omega_1$-algebra $A_1$ in $\Omega_2$-algebra $A_2$. Let set $B_2 \subset A_2$ be subalgebra of $\Omega_2$-algebra $A_2$ and stable set of representation $f$. Then there exists representation*

$$f_{B_2} : A_1 \relbar_* B_2$$

*such that $f_{B_2}(a) = f(a)|_{B_2}$. Representation $f_{B_2}$ is called* **subrepresentation** *of representation $f$.*

PROOF. Let $\omega_1$ be $n$-ary operation of $\Omega_1$-algebra $A_1$. Then for each $a_1$, ..., $a_n \in A_1$ and each $b \in B_2$

$$(f_{B_2}(a_1)...f_{B_2}(a_n)\omega_1)(b) = (f(a_1)...f(a_n)\omega_1)(b) = f(a_1...a_n\omega_1)(b)$$
$$= f_{B_2}(a_1...a_n\omega_1)(b)$$

Let $\omega_2$ be $n$-ary operation of $\Omega_2$-algebra $A_2$. Then for each $b_1$, ..., $b_n \in B_2$ and each $a \in A_1$

$$f_{B_2}(a)(b_1)...f_{B_2}(a)(b_n)\omega_2 = f(a)(b_1)...f(a)(b_n)\omega_2 = f(a)(b_1...b_n\omega_2)$$
$$= f_{B_2}(a)(b_1...b_n\omega_2)$$

We proved the statement of theorem. □

From the theorem 6.1.2, it follows that if $f_{B_2}$ is subrepresentation of representation $f$, then the map

$$(\text{id} : A \to A, \text{id}_{B_2} : B_2 \to A_2)$$

is morphism of representations.

THEOREM 6.1.3. *The set* [6.1] $\mathcal{B}_f$ *of all subrepresentations of representation* $f$ *generates a closure system on* $\Omega_2$-*algebra* $A_2$ *and therefore is a complete lattice.*

PROOF. Let $(K_\lambda)_{\lambda \in \Lambda}$ be the set off subalgebras of $\Omega_2$-algebra $A_2$ that are stable with respect to representation $f$. We define the operation of intersection on the set $\mathcal{B}_f$ according to rule

$$\bigcap f_{K_\lambda} = f_{\cap K_\lambda}$$

We defined the operation of intersection of subrepresentations properly. $\cap K_\lambda$ is subalgebra of $\Omega_2$-algebra $A_2$. Let $m \in \cap K_\lambda$. For each $\lambda \in \Lambda$ and for each $a \in A_1$, $f(a)(m) \in K_\lambda$. Therefore, $f(a)(m) \in \cap K_\lambda$. Therefore, $\cap K_\lambda$ is the stable set of representation $f$.                                                                                                        $\square$

We denote the corresponding closure operator by $J[f]$. Thus $J[f, X]$ is the intersection of all subalgebras of $\Omega_2$-algebra $A_2$ containing $X$ and stable with respect to representation $f$.

THEOREM 6.1.4. *Let* [6.2]

$$g : A_1 \longrightarrow\!\!\!* A_2$$

*be representation of* $\Omega_1$-*algebra* $A_1$ *in* $\Omega_2$-*algebra* $A_2$. *Let* $X \subset A_2$. *Define a subset* $X_k \subset A_2$ *by induction on* $k$.

6.1.4.1: $X_0 = X$

6.1.4.2: $x \in X_k => x \in X_{k+1}$

6.1.4.3: $x_1 \in X_k, ..., x_n \in X_k, \omega \in \Omega_2(n) => x_1...x_n\omega \in X_{k+1}$

6.1.4.4: $x \in X_k, a \in A => f(a)(x) \in X_{k+1}$

*Then*

$$(6.1.1) \qquad \bigcup_{k=0}^{\infty} X_k = J[f, X]$$

PROOF. If we put $U = \cup X_k$, then by definition of $X_k$, we have $X_0 \subset J[f, X]$, and if $X_k \subset J[f, X]$, then $X_{k+1} \subset J[f, X]$. By induction it follows that $X_k \subset J[f, X]$ for all $k$. Therefore,

$$(6.1.2) \qquad U \subset J[f, X]$$

If $a \in U^n$, $a = (a_1, ..., a_n)$, where $a_i \in X_{k_i}$, and if $k = \max\{k_1, ..., k_n\}$, then $a_1...a_n\omega \in X_{k+1} \subset U$. Therefore, $U$ is subalgebra of $\Omega_2$-algebra $A_2$.

If $m \in U$, then there exists such $k$ that $m \in X_k$. Therefore, $f(a)(m) \in X_{k+1} \subset U$ for any $a \in A_1$. Therefore, $U$ is stable set of the representation $f$.

Since $U$ is subalgebra of $\Omega_2$-algebra $A_2$ and is a stable set of the representation $f$, then subrepresentation $f_U$ is defined. Therefore,

$$(6.1.3) \qquad J[f, X] \subset U$$

From (6.1.2), (6.1.3), it follows that $J[f, X] = U$.                                        $\square$

---

[6.1] This definition is similar to definition of the lattice of subalgebras ([14], p. 79, 80). In general, In this and subsequent theorems of this chapter, it is necessary to consider the structure of universal algebras $A_1$ and $A_2$. Because the main task of this chapter is is the study of the structure of the representation, I deliberately simplified the theorems so that the details do not obscure the basic statements. This topic will be discussed in more details in the chapter 8, where theorems will be formulated in general form.

[6.2] The statement of theorem is similar to the statement of theorem 5.1, [14], page 79.

DEFINITION 6.1.5. $J[f, X]$ *is called* **subrepresentation** *generated by set* $X$, *and* $X$ *is a generating set of subrepresentation* $J[f, X]$. *In particular, a* **generating set** *of representation* $f$ *is a subset* $X \subset A_2$ *such that* $J[f, X] = A_2$.          □

The next definition follows from the theorem 6.1.4.

DEFINITION 6.1.6. *Let* $X \subset A_2$. *For each* $m \in J[f, X]$ *there exists* $\Omega_2$-**word** *defined according to following rules.* $w[f, X, m]$

6.1.6.1: *If* $m \in X$, *then* $m$ *is* $\Omega_2$-word.

6.1.6.2: *If* $m_1$, ..., $m_n$ *are* $\Omega_2$-words and $\omega \in \Omega_2(n)$, *then* $m_1 ... m_n \omega$ *is* $\Omega_2$-word.

6.1.6.3: *If* $m$ *is* $\Omega_2$-word and $a \in A_1$, *then* $f(a)(m)$ *is* $\Omega_2$-word.

*We will identify an element* $m \in J[f, X]$ *and corresponding it* $\Omega_2$-word using equation*

$$m = w[f, X, m]$$

*Similarly, for an arbitrary set* $B \subset J[f, X]$ *we consider the set of* $\Omega_2$-words[6.3]

$$w[f, X, B] = \{w[f, X, m] : m \in B\}$$

*We also use notation*

$$w[f, X, B] = (w[f, X, m], m \in B)$$

*Denote* $w[f, X]$ **the set of** $\Omega_2$-**words of representation** $J[f, X]$.          □

THEOREM 6.1.7. $w[f, X, X] = X$.

PROOF. The theorem follows from the statement 6.1.6.1.          □

THEOREM 6.1.8. *Let* $X$, $Y$ *be generating sets of representation*

$$f : A_1 \dashrightarrow A_2$$

*Let* $w[f, X, m]$ *be* $\Omega_2$-*word of* $A_2$-*number* $m$ *relative generating set* $X$. *Let* $w[f, Y, X]$ *be the set of* $\Omega_2$-*words of the set* $X$ *relative generating set* $Y$. *If, in the word* $w[f, X, m]$, *we substitute image* $w[f, Y, x]$ *of each* $x \in X$, *then we get* $\Omega_2$-*word* $w[f, Y, m]$ *of* $A_2$-*number* $m$ *relative generating set* $Y$.

*Transformation of* $\Omega_2$-*words*

$$w[f, X, m] \to w[f, Y, m]$$

$$w[f, Y, m] = w[f, Y, X] \circ w[f, X, m]$$

*is called supperposition of coordinates.*

PROOF. We prove the theorem by induction over complexity of $\Omega_2$-word.

If $m \in X$, then $w[f, X, m] = m$. If we substitute image $w[f, Y, x]$ of $m$, then we get $\Omega_2$-word $w[f, Y, m]$ of $A_2$-number $m$ relative generating set $Y$.

Let $\Omega_2$-word $w[f, X, m]$ of $A_2$-number $m$ has form

(6.1.4)          $$w[f, X, m] = w[f, X, m_1] ... w[f, X, m_n]\omega$$

where $\omega \in \Omega_2(n)$ and, for each $A_2$-number $m_i$, we defined map

$$w[f, X, m_i] \to w[f, Y, m_i]$$

---

[6.3]The expression $w[f, X, m]$ is a special case of the expression $w[f, X, B]$, namely

$$w[f, X, \{m\}] = \{w[f, X, m]\}$$

According to the statement 6.1.6.2, the expression

$$w[f, Y, m_1]...w[f, Y, m_n]\omega$$

is $\Omega_2$-word $w[f, Y, m]$ of $A_2$-number $m$ relative generating set $Y$. Therefore, we defined map

$$w[f, X, m] \to w[f, Y, m]$$

for $A_2$-number $m$.

Let $\Omega_2$-word $w[f, X, m]$ of $A_2$-number $m$ has form

(6.1.5)          $$w[f, X, m] = f(a)(w[f, X, m_1])$$

where, for $A_2$-number $m_1$, we defined map

$$w[f, X, m_1] \to w[f, Y, m_1]$$

According to the statement 6.1.6.3, the expression

$$f(a)(w[f, Y, m_1])$$

is $\Omega_2$-word $w[f, Y, m]$ of $A_2$-number $m$ relative generating set $Y$. Therefore, we defined map

$$w[f, X, m] \to w[f, Y, m]$$

for $A_2$-number $m$.          $\square$

Choice of $\Omega_2$-word relative generating set $X$ is ambiguous. Therefore, if $\Omega_2$-number has different $\Omega_2$-words, then we will use indexes to distinguish them: $w[f, X, m]$, $w_1[f, X, m]$, $w_2[f, X, m]$.

**DEFINITION 6.1.9.** *Generating set $X$ of representation $f$ generates equivalence*

$$\rho[f, X] = \{(w[f, X, m], w_1[f, X, m]) : m \in A_2\}$$

*on the set of $\Omega_2$-words.*          $\square$

According to the definition 6.1.9, two $\Omega_2$-words with respect to the generating set $X$ of representation $f$ are equivalent iff they correspond to the same $A_2$-number. When we write equality of two $\Omega_2$-words with respect to the generating set $X$ of representation $f$, we will keep in mind that this equality is true up to equivalence $\rho[f, X]$.

**THEOREM 6.1.10.** *Let*

$$f : A_1 \longrightarrow_* A_2$$

*be representation of $\Omega_1$-algebra $A_1$ in $\Omega_2$-algebra $A_2$. Let*

$$g : A_1 \longrightarrow_* B_2$$

*be representation of $\Omega_1$-algebra $A_1$ in $\Omega_2$-algebra $B_2$. Let $X$ be the generating set of representation $f$. Let*

$$R : A_2 \to B_2$$

*be reduced morphism of representation[6.4] and $X' = R(X)$. Reduced morphism $R$ of representation generates the map of $\Omega_2$-words*

$$w[f \to g, X, R] : w[f, X] \to w[g, X']$$

*such that*

---

[6.4] I considered morphism of representation in the theorem 8.1.7.

6.1.10.1: *If* $m \in X$, $m' = R(m)$, *then*

$$w[f \rightarrow g, X, R](m) = m'$$

6.1.10.2: *If*

$$m_1, ..., m_n \in w[f, X]$$

$$m'_1 = w[f \rightarrow g, X, R](m_1) \quad ... \quad m'_n = w[f \rightarrow g, X, R](m_n)$$

*then for operation* $\omega \in \Omega_2(n)$ *holds*

$$w[f \rightarrow g, X, R](m_1...m_n\omega) = m'_1...m'_n\omega$$

6.1.10.3: *If*

$$m \in w[f, X] \quad m' = w[f \rightarrow g, X, R](m) \quad a \in A_1$$

*then*

$$w[f \rightarrow g, X, R](f(a)(m)) = g(a)(m')$$

Proof. Statements 6.1.10.1, 6.1.10.2 are true by definition of the reduced morphism $R$. The statement 6.1.10.3 follows from the equality (3.4.5). □

Remark 6.1.11. *Let*

$$R : A_2 \rightarrow B_2$$

*be reduced morphism of representation. Let*

$$m \in J[f, X] \quad m' = R(m) \quad X' = R(X)$$

*The theorem 6.1.10 states that* $m' \in J[g, X']$. *The theorem 6.1.10 also states that* $\Omega_2$-word representing $m$ relative $X$ and $\Omega_2$-word representing $m'$ relative $X'$ are generated according to the same algorithm. This allows considering of the set of $\Omega_2$-words $w[g, X', m']$ as map

(6.1.6)           $W[f, X, m] : (g, X') \rightarrow (g, X') \circ W[f, X, m] = w[g, X', m']$

*where*

$$X' = R(X) \quad m' = R(m)$$

*for certain reduced morphism* $R$.

*If* $f = g$, *then, instead of the map* (6.1.6), *we consider the map*

$$W[f, X, m] : X' \rightarrow X' \circ W[f, X, m] = w[f, X', m']$$

*such that, if for certain endomorphism* $R$

$$X' = R(X) \quad m' = R(m)$$

*then*

$$W[f, X, m](X') = X' \circ W[f, X, m] = w[f, X', m'] = m'$$

*The map* $W[f, X, m]$ *is called* **coordinates of** $A_2$**-number** $m$ *relative to the set* $X$. *Similarly, we consider coordinates of a set* $B \subset J[f, X]$ *relative to the set* $X$

$$W[f, X, B] = \{W[f, X, m] : m \in B\} = (W[f, X, m], m \in B)$$

*Denote*

$$W[f, X] = \{W[f, X, m] : m \in J[f, X]\} = (W[f, X, m], m \in J[f, X])$$

**the set of coordinates of representation** $J[f, X]$.                               $\square$

THEOREM 6.1.12. *There is a structure of $\Omega_2$-algebra on the set of coordinates $W[f, X]$.*

PROOF. Let $\omega \in \Omega_2(n)$. Then for any $m_1, ..., m_n \in J[f, X]$ , we assume

$$(6.1.7) \qquad W[f, X, m_1]...W[f, X, m_n]\omega = W[f, X, m_1...m_n\omega]$$

According to the remark 6.1.11,

$$(6.1.8) \qquad \begin{aligned} X \circ (W[f, X, m_1]...W[f, X, m_n]\omega) &= X \circ W[f, X, m_1...m_n\omega] \\ &= w[f, X, m_1...m_n\omega] \end{aligned}$$

follows from the equation (6.1.7). According to rule 6.1.6.2, from the equation (6.1.8), it follows that

$$(6.1.9) \qquad \begin{aligned} &X \circ (W[f, X, m_1]...W[f, X, m_n]\omega) \\ &= w[f, X, m_1]...w[f, X, m_n]\omega \\ &= (X \circ W[f, X, m_1])...(X \circ W[f, X, m_n])\omega \end{aligned}$$

From the equation (6.1.9), it follows that the operation $\omega$ defined by the equation (6.1.7) on the set of coordinates is defined properly.                               $\square$

THEOREM 6.1.13. *There exists the representation of $\Omega_1$-algebra $A_1$ in $\Omega_2$-algebra $W[f, X]$.*

PROOF. Let $a \in A_1$. Then for any $m \in J[f, X]$, we assume

$$(6.1.10) \qquad f(a)(W[f, X, m]) = W[f, X, f(a)(m)]$$

According to the remark 6.1.11,

$$(6.1.11) \qquad X \circ (f(a)(W[f, X, m])) = X \circ W[f, X, f(a)(m)] = w[f, X, f(a)(m)]$$

follows from the equation (6.1.10). According to rule 6.1.6.3, from the equation (6.1.11), it follows that

$$(6.1.12) \qquad X \circ (f(a)(W[f, X, m])) = f(a)(w[f, X, m]) = f(a)(X \circ W[f, X, m])$$

From the equation (6.1.12), it follows that the representation (6.1.10) of $\Omega_1$-algebra $A_1$ in $\Omega_2$-algebra $W[f, X]$ is defined properly.                               $\square$

THEOREM 6.1.14. *Let*

$$f : A_1 \longrightarrow_* A_2$$

*be representation of $\Omega_1$-algebra $A_1$ in $\Omega_2$-algebra $A_2$. Let*

$$g : A_1 \longrightarrow_* B_2$$

*be representation of $\Omega_1$-algebra $A_1$ in $\Omega_2$-algebra $B_2$. For given sets $X \subset A_2$, $X' \subset B_2$, let map*

$$R_1 : X \to X'$$

*agree with the structure of representation $f$, i. e.*

$$\omega \in \Omega_2(n) \ x_1, \ ..., \ x_n, \ x_1...x_n\omega \in X, \ R_1(x_1...x_n\omega) \in X'$$

$$=>R_1(x_1...x_n\omega) = R_1(x_1)...R_1(x_n)\omega$$

$$x \in X, \ a \in A, \ R_1(f(a)(x)) \in X'$$

$$=>R_1(f(a)(x)) = g(a)(R_1(x))$$

*Consider the map of $\Omega_2$-words*

$$w[f \to g, X, X', R_1] : w[f, X] \to w[g, X']$$

*that satisfies conditions 6.1.10.1, 6.1.10.2, 6.1.10.3 and such that*

$$x \in X => w[f \to g, X, X', R_1](x) = R_1(x)$$

*There exists unique map*

$$R : A_2 \to B_2$$

*defined by rule*

$$R(m) = w[f \to g, X, X', R_1](w[f, X, m])$$

*which is reduced morphism of representations $J[f, X]$ and $J[g, X']$.*

Proof. We prove the theorem by induction over complexity of $\Omega_2$-word.
If $w[f, X, m] = m$, then $m \in X$. According to condition 6.1.10.1,

$$R(m) = w[f \to g, X, X', R_1](w[f, X, m]) = w[f, X, R_1](m) = R_1(m)$$

Therefore, maps $R$ and $R_1$ coinside on the set $X$, and the map $R$ agrees with structure of representation $f$.

Let $\omega \in \Omega_2(n)$. Let the map $R$ be defined for $m_1, ..., m_n \in J[f, X]$. Let

$$w_1 = w[f, X, m_1] \quad ... \quad w_n = w[f, X, m_n]$$

f $m = m_1...m_n\omega$, then according to rule 6.1.6.2,

$$w[f, X, m] = w_1...w_n\omega$$

According to condition 6.1.10.2,

$$R(m) = w[f \to g, X, X', R_1](w[f, X, m]) = w[f \to g, X, X', R_1](w_1...w_n\omega)$$

$$= w[f \to g, X, X', R_1](w_1)...w[f \to g, X, X', R_1](w_n)\omega$$

$$= R(m_1)...R(m_n)\omega$$

Therefore, the map $R$ is endomorphism of $\Omega_2$-algebra $A_2$.

Let the map $R$ be defined for $m_1 \in J[f, X]$, $w_1 = w[f, X, m_1]$. Let $a \in A_1$.
If $m = f(a)(m_1)$, then according to rule 6.1.6.3,

$$w[f, X, f(a)(m_1)] = f(a)(w_1)$$

According to condition 6.1.10.3,

$$R(m) = w[f \to g, X, X', R_1](w[f, X, m]) = w[f \to g, X, X', R_1](f(a)(w_1))$$

$$= f(a)(w[f \to g, X, X', R_1](w_1)) = f(a)(R(m_1))$$

From equation (3.2.3), it follows that the map $R$ is morphism of the representation $f$.

The statement that the endomorphism $R$ is unique and therefore this endomorphism is defined properly follows from the following argument. Let $m \in A_2$ have different $\Omega_2$-words relative the set $X$, for instance

(6.1.13)                    $m = x_1...x_n\omega = f(a)(x)$

Because $R$ is endomorphism of representation, then, from the equation (6.1.13), it follows that

(6.1.14)        $R(m) = R(x_1...x_n\omega) = R(x_1)...R(x_n)\omega = R(f(a)(x)) = f(a)(R(x))$

From the equation (6.1.14), it follows that

(6.1.15)              $R(m) = R(x_1)...R(x_n)\omega = f(a)(R(x))$

From equations (6.1.13), (6.1.15), it follows that the equation (6.1.13) is preserved under the map. Therefore, the image of $A_2$ does not depend on the choice of coordinates.                    □

REMARK 6.1.15. *The theorem 6.1.14 is the theorem of extension of map. The only statement we know about the set $X$ is the statement that $X$ is generating set of the representation $f$. However, between the elements of the set $X$ there may be relationships generated by either operations of $\Omega_2$-algebra $A_2$, or by transformation of representation $f$. Therefore, any map of set $X$, in general, cannot be extended to a reduced morphism of representation $f$.[6.5] However, if the map $R_1$ is coordinated with the structure of representation on the set $X$, then we can construct an extension of this map and this extension is reduced morphism of representation $f$.*                    □

DEFINITION 6.1.16. *Let $X$ be the generating set of the representation*

$$f : A_1 \longrightarrow_* A_2$$

*of $\Omega_1$-algebra $A_1$ in $\Omega_2$-algebra $A_2$. Let $Y$ be the generating set of the representation*

$$g : A_1 \longrightarrow_* B_2$$

*of $\Omega_1$-algebra $A_1$ in $\Omega_2$-algebra $B_2$. Let*

$$R : A_2 \to B_2$$

*be the reduced morphism of the representation $f$. The set of coordinates $W[g, Y, R(X)]$ is called* **coordinates of reduced morphism of representation**.                    □

From definitions 6.1.6, 6.1.16, it follows that

$$W[g, Y, R(X)] = (W[g, Y, R(x)], x \in X)$$

Let $m \in A_2$. If, in the word $w[f, X, m]$, we substitute image $w[g, Y, R(x)]$ of each $x \in X$, then, according to the theorem 6.1.14, we get $\Omega_2$-word $w[g, Y, R(m)]$. The definition 6.1.17 follows from this statement.

---

[6.5]In the theorem 6.2.10, requirements to generating set are more stringent. Therefore, the theorem 6.2.10 says about extension of arbitrary map. A more detailed analysis is given in the remark 6.2.12.

DEFINITION 6.1.17. *Let $X$ be the generating set of the representation*

$$f : A_1 \longrightarrow_* A_2$$

*of $\Omega_1$-algebra $A_1$ in $\Omega_2$-algebra $A_2$. Let $Y$ be the generating set of the representation*

$$g : A_1 \longrightarrow_* B_2$$

*of $\Omega_1$-algebra $A_1$ in $\Omega_2$-algebra $B_2$. Let $R$*

$$R : A_2 \to B_2$$

*be the reduced morphism of the representation $f$. Let $m \in A_2$. We define* **superposition of coordinates** *of the reduced morphism $R$ of the representation $f$ and $A_2$-number $m$ as coordinates defined according to rule*

(6.1.16)        $W[g, Y, R(X)] \circ W[f, X, m] = W[g, Y, R(m)]$

*We define superposition of coordinates of the reduced morphism $R$ of the representation $f$ and the set $B \subseteq A_2$ according to rule*

(6.1.17)        $W[g, Y, R(X)] \circ W[f, X, B] = (W[g, Y, R(X)] \circ W[f, X, m], m \in B)$

$W[g, Y, R(X)] \circ w[f, X, B] = w[g, Y, R(X)] \circ W[f, X, B] = w[g, Y, R(B)]$

□

THEOREM 6.1.18. *Let $X$ be the generating set of the representation*

$$f : A_1 \longrightarrow_* A_2$$

*of $\Omega_1$-algebra $A_1$ in $\Omega_2$-algebra $A_2$. Let $Y$ be the generating set of the representation*

$$g : A_1 \longrightarrow_* B_2$$

*of $\Omega_1$-algebra $A_1$ in $\Omega_2$-algebra $B_2$. Reduced morphism of representation*

$$R : A_2 \to B_2$$

*generates the map of coordinates of representation*

(6.1.18)        $W[f \to g, X, Y, R] : W[f, X] \to W[g, Y]$

*such that*

(6.1.19)        $W[f, X, m] \to W[f \to g, X, Y, R] \circ W[f, X, m] = W[g, Y, R(m)]$

PROOF. According to the remark 6.1.11, we consider equations (6.1.16), (6.1.18) relative to given generating sets $X$, $Y$. The word

(6.1.20)        $X \circ W[f, X, m] = w[f, X, m]$

corresponds to coordinates $W[f, X, m]$; the word

(6.1.21)        $Y \circ W[g, Y, R(m)] = w[g, Y, R(m)]$

corresponds to coordinates $W[g, Y, R(m)]$. Therefore, in order to prove the theorem, it is sufficient to show that the map $W[f, X, R]$ corresponds to map $w[f, X, R]$. We prove this statement by induction over complexity of $\Omega_2$-word.

If $m \in X$, $m' = R(m)$, then, according to equations (6.1.20), (6.1.21), maps $W[f, X, R]$ and $w[f, X, R]$ are coordinated.

Let for $m_1, ..., m_n \in X$ maps $W[f, X, R]$ and $w[f, X, R]$ be coordinated. Let $\omega \in \Omega_2(n)$. According to the theorem 6.1.12

(6.1.22)        $W[f, X, m_1...m_n\omega] = W[f, X, m_1]...W[f, X, m_n]\omega$

Because $R$ is endomorphism of $\Omega_2$-algebra $A_2$, then from the equation (6.1.22), it follows that

$$
(6.1.23) \qquad \begin{aligned}
W[f, X, R \circ (m_1...m_n\omega)] &= W[f, X, (R \circ m_1)...(R \circ m_n)\omega] \\
&= W[f, X, R \circ m_1]...W[f, X, R \circ m_n]\omega
\end{aligned}
$$

From equations (6.1.22), (6.1.23) and the statement of induction, it follows that the maps $W[f, X, R]$ and $w[f, X, R]$ are coordinated for $m = m_1...m_n\omega$.

Let for $m_1 \in A_2$ maps $W[f, X, R]$ and $w[f, X, R]$ are coordinated. Let $a \in A_1$. According to the theorem 6.1.13

$$
(6.1.24) \qquad W[f, X, f(a)(m_1)] = f(a)(W[f, X, m_1])
$$

Because $R$ is endomorphism of representation $f$, then, from the equation (6.1.24), it follows that

$$
(6.1.25) \quad W[f, X, R \circ f(a)(m_1)] = W[f, X, f(a)(R \circ m_1)] = f(a)(W[f, X, R \circ m_1])
$$

From equations (6.1.24), (6.1.25) and the statement of induction, it follows that maps $W[f, X, R]$ and $w[f, X, R]$ are coordinated for $m = f(a)(m_1)$. $\qquad\square$

**Corollary 6.1.19.** *Let $X$ be the generating set of the representation $f$. Let $R$ be the endomorphism of the representation $f$. The map $W[f, X, R]$ is endomorphism of representation of $\Omega_1$-algebra $A_1$ in $\Omega_2$-algebra $W[f, X]$.* $\qquad\square$

Hereinafter we will identify map $W[f, X, R]$ and the set of coordinates $W[f, X, R \circ X]$.

**Theorem 6.1.20.** *Let $X$ be the generating set of the representation $f$. Let $R$ be the endomorphism of the representation $f$. Let $Y \subset A_2$. Then*

$$
(6.1.26) \qquad W[f, X, R(X)] \circ W[f, X, Y] = W[f, X, R(Y)]
$$

**Proof.** The equation (6.1.26) follows from the equation

$$
R \circ Y = (R \circ m, m \in Y)
$$

as well from equations (6.1.16), (6.1.17). $\qquad\square$

**Theorem 6.1.21.** *Let $X$ be the generating set of the representation $f$. Let $R$, $S$ be the endomorphisms of the representation $f$. Then*

$$
(6.1.27) \qquad W[f, X, R] \circ W[f, X, S] = W[f, X, R \circ S]
$$

**Proof.** The equation (6.1.27) follows from the equation (6.1.26), if we assume $Y = S \circ X$. $\qquad\square$

The concept of superposition of coordinates is very simple and resembles a kind of Turing machine. If element $m \in A_2$ has form either

$$
m = m_1...m_n\omega
$$

or

$$
m = f(a)(m_1)
$$

then we are looking for coordinates of elements $m_i$ to substitute them in an appropriate expression. As soon as an element $m \in A_2$ belongs to the generating set of $\Omega_2$-algebra $A_2$, we choose coordinates of the corresponding element of the second factor. Therefore, we require that the second factor in the superposition has been the set of coordinates of the image of the generating set $X$.

The following forms of writing an image of the set $Y$ under endomorphism $R$ are equivalent.

(6.1.28)        $R \circ Y = (R(X)) \circ W[f, X, Y] = (X \circ W[f, X, R]) \circ W[f, X, Y]$

From equations (6.1.26), (6.1.28), it follows that

(6.1.29)        $X \circ (W[f, X, R] \circ W[f, X, Y]) = (X \circ W[f, X, R]) \circ W[f, X, Y]$

The equation (6.1.29) is associative law for composition and allows us to write expression

$$X \circ W[f, X, R] \circ W[f, X, Y]$$

without brackets.

DEFINITION 6.1.22. *Let $X \subset A_2$ be generating set of representation*

$$f : A_1 \longrightarrow_* A_2$$

*Let the map*

$$H : A_2 \to A_2$$

*be endomorphism of the representation $f$. Let the set $X' = H \circ X$ be the image of the set $X$ under the map $H$. Endomorphism $H$ of representation $f$ is called regular on the generating set $X$, if the set $X'$ is the generating set of representation $f$. Otherwise, endomorphism $H$ of representation $f$ is called singular on the generating set $X$.* □

DEFINITION 6.1.23. *Endomorphism of representation $f$ is called **regular**, if it is regular on every generating set. Otherwise, endomorphism $H$ of representation $f$ is called **singular**.* □

THEOREM 6.1.24. *Automorphism $R$ of representation*

$$f : A_1 \longrightarrow_* A_2$$

*is regular endomorphism.*

PROOF. Let $X$ be generating set of representation $f$. Let $X' = R(X)$.

According to theorem 6.1.10 endomorphism $R$ forms the map of $\Omega_2$-words $w[f \to g, X, R]$.

Let $m' \in A_2$. Since $R$ is automorphism, then there exists $m \in A_2$, $R \circ m = m'$. According to definition 6.1.6, $w[f, X, m]$ is $\Omega_2$-word, representing $A_2$ relative to generating set $X$. According to theorem 6.1.10, $w[f, X', m']$ is $\Omega_2$-word, representing of $m'$ relative to generating set $X'$

$$w[f, X', m'] = w[f \to g, X, R](w[f, X, m])$$

Therefore, $X'$ is generating set of representation $f$. According to definition 6.1.23, automorphism $R$ is regular. □

## 6.2. Basis of representation

DEFINITION 6.2.1. *Let*

$$f : A_1 \longrightarrow_* A_2$$

*be representation of $\Omega_1$-algebra $A_1$ in $\Omega_2$-algebra $A_2$ and*

$$Gen[f] = \{X \subseteq A_2 : J[f, X] = A_2\}$$

*If, for the set $X \subset A_2$, it is true that $X \in Gen[f]$, then for any set $Y$, $X \subset Y \subset A_2$, also it is true that $Y \in Gen[f]$. If there exists minimal set $X \in Gen[f]$, then the set $X$ is called* **quasibasis of representation** $f$.                     $\square$

THEOREM 6.2.2. *If the set $X$ is the quasibasis of the representation $f$, then, for any $m \in X$, the set $X \setminus \{m\}$ is not generating set of the representation $f$.*

PROOF. Let $X$ be quasibasis of the representation $f$. Assume that for some $m \in X$ there exist $\Omega_2$-word

$$w = w[f, X \setminus \{m\}, m]$$

Consider $A_2$-number $m'$ such that it has $\Omega_2$-word $w' = w[f, X, m']$ that depends on $m$. According to the definition 6.1.6, any occurrence of $A_2$-number $m$ into $\Omega_2$-word $w'$ can be substituted by the $\Omega_2$-word $w$. Therefore, the $\Omega_2$-word $w'$ does not depend on $m$, and the set $X \setminus \{m\}$ is generating set of representation $f$. Therefore, $X$ is not quasibasis of representation $f$.                                       $\square$

REMARK 6.2.3. *The proof of the theorem 6.2.2 gives us effective method for constructing the quasibasis of the representation $f$. Choosing an arbitrary generating set, step by step, we remove from set those elements which have coordinates relative to other elements of the set. If the generating set of the representation is infinite, then this construction may not have the last step. If the representation has finite generating set, then we need a finite number of steps to construct a quasibasis of this representation.*                                                                     $\square$

We introduced $\Omega_2$-word of $x \in A_2$ relative generating set $X$ in the definition 6.1.6. From the theorem 6.2.2, it follows that if the generating set $X$ is not an quasibasis, then a choice of $\Omega_2$-word relative generating set $X$ is ambiguous. However, even if the generating set $X$ is an quasibasis, then a representation of $m \in A_2$ in form of $\Omega_2$-word is ambiguous.

REMARK 6.2.4. *There are three reasons of ambiguity in notation of $\Omega_2$-word.*

6.2.4.1: *In $\Omega_i$-algebra $A_i$, $i = 1, 2$, equalities may be defined. For instance, if $e$ is unit of multiplicative group $A_i$, then the equality*

$$ae = a$$

*is true for any $a \in A_i$.*

6.2.4.2: *Ambiguity of choice of $\Omega_2$-word may be associated with properties of representation. For instance, if $m_1, ..., m_n$ are $\Omega_2$-words, $\omega \in \Omega_2(n)$ and $a \in A_1$, then*[6.6]

(6.2.1)                $f(a)(m_1...m_n\omega) = (f(a)(m_1))...(f(a)(m_n))\omega$

*At the same time, if $\omega$ is operation of $\Omega_1$-algebra $A_1$ and operation of $\Omega_2$-algebra $A_2$, then we require that $\Omega_2$-words $f(a_1...a_n\omega)(x)$ and*

---

[6.6] For instance, let $\{e_1, e_2\}$ be the basis of vector space over field $k$. The equation (6.2.1) has the form of distributive law

$$a(b^1e_1 + b^2e_2) = (ab^1)e_1 + (ab^2)e_2$$

$(f(a_1)(x))...(f(a_n)(x))\omega$   *describe the same element of $\Omega_2$-algebra $A_2$.* [6.7]

(6.2.2) $$f(a_1...a_n\omega)(x) = (f(a_1)(x))...(f(a_n)(x))\omega$$

6.2.4.3: *Equalities like* (6.2.1), (6.2.2) *persist under morphism of representation. Therefore we can ignore this form of ambiguity of $\Omega_2$-word. However, a fundamentally different form of ambiguity is possible. We can see an example of such ambiguity in theorems 9.3.15, 9.3.16.*

*So we see that we can define different equivalence relations on the set of $\Omega_2$-words.* [6.8] *Our goal is to find a maximum equivalence on the set of $\Omega_2$-words which persist under morphism of representation.*

*A similar remark concerns the map $W[f, X, m]$ defined in the remark 6.1.11.* [6.9]

<div align="right">□</div>

THEOREM 6.2.5. *Let $X$ be quasibasis of the representation*

$$f : A_1 \longrightarrow_* A_2$$

*Consider equivalence*

$$\lambda[f, X] \subseteq w[f, X] \times w[f, X]$$

*which is generated exclusively by the following statements.*

6.2.5.1: *If in $\Omega_2$-algebra $A_2$ there is an equality*

$$w_1[f, X, m] = w_2[f, X, m]$$

*defining structure of $\Omega_2$-algebra, then*

$$(w_1[f, X, m], w_2[f, X, m]) \in \lambda[f, X]$$

6.2.5.2: *If in $\Omega_1$-algebra $A_1$ there is an equality*

$$w_1[f, X, m] = w_2[f, X, m]$$

*defining structure of $\Omega_1$-algebra, then*

$$(f(w_1)(w[f, X, m]), f(w_2)(w[f, X, m])) \in \lambda[f, X]$$

---

[6.7]For vector space, this requirement has the form of distributive law

$$(a + b)e_1 = ae_1 + be_1$$

[6.8] Evidently each of the equalities (6.2.1), (6.2.2) generates some equivalence relation.

[6.9]If vector space has finite basis, then we represent the basis as matrix

$$\overline{\overline{e}} = \begin{pmatrix} e_1 & ... & e_2 \end{pmatrix}$$

We present the map $W[f, \overline{\overline{e}}](v)$ as matrix

$$W[f, \overline{\overline{e}}, v] = \begin{pmatrix} v^1 \\ ... \\ v^n \end{pmatrix}$$

Then

$$W[f, \overline{\overline{e}}, v](\overline{\overline{e}}') = W[f, \overline{\overline{e}}, v]\begin{pmatrix} e_1' & ... & e_n' \end{pmatrix} = \begin{pmatrix} v^1 \\ ... \\ v^n \end{pmatrix}\begin{pmatrix} e_1' & ... & e_n' \end{pmatrix}$$

has form of matrix product.

6.2.5.3: *For any operation* $\omega \in \Omega_1(n)$,

$$(f(a_{11}...a_{1n}\omega)(a_2), (f(a_{11})...f(a_{1n})\omega)(a_2)) \in \lambda[f, X]$$

6.2.5.4: *For any operation* $\omega \in \Omega_2(n)$,

$$(f(a_1)(a_{21}...a_{2n}\omega), f(a_1)(a_{21})...f(a_1)(a_{2n})\omega) \in \lambda[f, X]$$

6.2.5.5: *Let* $\omega \in \Omega_1(n) \cap \Omega_2(n)$. *If the representation* $f$ *satisfies equality*[6.10]

$$f(a_{11}...a_{1n}\omega)(a_2) = (f(a_{11})(a_2))...(f(a_{1n})(a_2))\omega$$

*then we can assume that the following equality is true*

$$(f(a_{11}...a_{1n}\omega)(a_2), (f(a_{11})(a_2))...(f(a_{1n})(a_2))\omega) \in \lambda[f, X]$$

Proof. The theorem is true because considered equalities are preserved under homomorphisms of universal algebras $A_1$ and $A_2$. □

Definition 6.2.6. *Quasibasis* $\overline{\overline{e}}$ *of the representation* $f$ *such that*

$$\rho[f, \overline{\overline{e}}] = \lambda[f, \overline{\overline{e}}]$$

*is called* **basis of representation** $f$. □

Remark 6.2.7. *As noted by Paul Cohn in [14], p. 82, 83, the representation may have inequivalent bases. For instance, the cyclic group of order six has bases* $\{a\}$ *and* $\{a^2, a^3\}$ *which we cannot map one into another by endomorphism of the representation.* □

Remark 6.2.8. *We write a basis also in following form*

$$\overline{\overline{e}} = (e, e \in \overline{\overline{e}})$$

*If basis is finite, then we also use notation*

$$\overline{\overline{e}} = (e_i, i \in I) = (e_1, ..., e_n)$$

□

Theorem 6.2.9. *Automorphism of the representation* $f$ *maps a basis of the representation* $f$ *into basis.*

Proof. Let the map $R$ be automorphism of the representation $f$. Let the set $\overline{\overline{e}}$ be a basis of the representation $f$. Let[6.11] $\overline{\overline{e}}' = R \circ \overline{\overline{e}}$. Assume that the set $\overline{\overline{e}}'$ is not basis. According to the theorem 6.2.2 there exists such $e' \in \overline{\overline{e}}'$ that $\overline{\overline{e}}' \setminus \{e'\}$ is generating set of the representation $f$. According to the theorem 3.5.5, the map $R^{-1}$ is automorphism of the representation $f$. According to the theorem 6.1.24 and definition 6.1.23, the set $\overline{\overline{e}} \setminus \{e\}$ is generating set of the representation $f$. The contradiction completes the proof of the theorem. □

---

[6.10] Consider a representation of commutative ring $D$ in $D$-algebra $A$. We will use notation

$$f(a)(v) = av$$

The operations of addition and multiplication are defined in both algebras. However the equality

$$f(a + b)(v) = f(a)(v) + f(b)(v)$$

is true, and the equality

$$f(ab)(v) = f(a)(v)f(b)(v)$$

is wrong.

[6.11] According to definitions 5.1.3, 6.4.1, we will use notation $R(\overline{\overline{e}}) = R \circ \overline{\overline{e}}$.

THEOREM 6.2.10. *Let $\overline{\overline{e}}$ be the basis of the representation $f$. Let*

$$R_1 : \overline{\overline{e}} \to \overline{\overline{e}}'$$

*be arbitrary map of the set $X$. Consider the map of $\Omega_2$-words*

$$w[f \to g, \overline{\overline{e}}, \overline{\overline{e}}', R_1] : w[f, \overline{\overline{e}}] \to w[g, \overline{\overline{e}}']$$

*that satisfies conditions 6.1.10.1, 6.1.10.2, 6.1.10.3 and such that*

$$e \in \overline{\overline{e}} \Longrightarrow w[f \to g, \overline{\overline{e}}, \overline{\overline{e}}', R_1](e) = R_1(e)$$

*There exists unique endomorphism of representation $f$* [6.12]

$$r_2 : A_2 \to A_2$$

*defined by rule*

$$R(m) = w[f \to g, \overline{\overline{e}}, \overline{\overline{e}}', R_1](w[f, \overline{\overline{e}}, m])$$

PROOF. The statement of theorem is corollary theorems 6.1.10, 6.1.14. □

COROLLARY 6.2.11. *Let $\overline{\overline{e}}$, $\overline{\overline{e}}'$ be the bases of the representation $f$. Let $R$ be the automorphism of the representation $f$ such that $\overline{\overline{e}}' = R \circ \overline{\overline{e}}$. Automorphism $R$ is uniquely defined.* □

REMARK 6.2.12. *The theorem 6.2.10, as well as the theorem 6.1.14, is the theorem of extension of map. However in this theorem, $\overline{\overline{e}}$ is not arbitrary generating set of the representation, but basis. According to remark 6.2.3, we cannot determine coordinates of any element of basis through the remaining elements of the same basis. Therefore, we do not need to coordinate the map of the basis with representation.* □

THEOREM 6.2.13. *The set of coordinates $W[f, \overline{\overline{e}}, \overline{\overline{e}}]$ corresponds to identity transformation*

$$W[f, \overline{\overline{e}}, E] = W[f, \overline{\overline{e}}, \overline{\overline{e}}]$$

PROOF. The statement of the theorem follows from the equation

$$m = \overline{\overline{e}} \circ W[f, \overline{\overline{e}}, m] = \overline{\overline{e}} \circ W[f, \overline{\overline{e}}, \overline{\overline{e}}] \circ W[f, \overline{\overline{e}}, m]$$

□

THEOREM 6.2.14. *Let $W[f, \overline{\overline{e}}, R \circ \overline{\overline{e}}]$ be the set of coordinates of automorphism $R$. There exists set of coordinates $W[f, R \circ \overline{\overline{e}}, \overline{\overline{e}}]$ corresponding to automorphism $R^{-1}$. The set of coordinates $W[f, R \circ \overline{\overline{e}}, \overline{\overline{e}}]$ satisfies to equation*

(6.2.3)       $$W[f, \overline{\overline{e}}, R \circ \overline{\overline{e}}] \circ W[f, R \circ \overline{\overline{e}}, \overline{\overline{e}}] = W[f, \overline{\overline{e}}, \overline{\overline{e}}]$$

$$W[f \to f, \overline{\overline{e}}, \overline{\overline{e}}, R^{-1}] = W[f \to f, \overline{\overline{e}}, \overline{\overline{e}}, R]^{-1} = W[f, R \circ \overline{\overline{e}}, \overline{\overline{e}}]$$

PROOF. Since $R$ is automorphism of the representation $f$, then, according to the theorem 6.2.9, the set $R \circ \overline{\overline{e}}$ is a basis of the representation $f$. Therefore, there exists the set of coordinates $W[f, R \circ \overline{\overline{e}}, \overline{\overline{e}}]$. The equation (6.2.3) follows from the chain of equations

$$W[f, \overline{\overline{e}}, R \circ \overline{\overline{e}}] \circ W[f, R \circ \overline{\overline{e}}, \overline{\overline{e}}] = W[f, \overline{\overline{e}}, R \circ \overline{\overline{e}}] \circ W[f, \overline{\overline{e}}, R^{-1} \circ \overline{\overline{e}}]$$
$$= W[f, \overline{\overline{e}}, R \circ R^{-1} \circ \overline{\overline{e}}] = W[f, \overline{\overline{e}}, \overline{\overline{e}}]$$

□

---

[6.12]This statement is similar to the theorem [2]-4.1, p. 135.

Remark 6.2.15. *In $\Omega_2$-algebra $A_2$ there is no universal algorithm for determining the set of coordinates $W[f, R \circ \overline{\overline{e}}, \overline{\overline{e}}]$ for given set $W[f, \overline{\overline{e}}, R \circ \overline{\overline{e}}]$.* [6.13] *We assume that in the theorem 6.2.14 this algorithm is given implicitly. It is evident also that the set of $\Omega_2$-words*

$$(6.2.4) \qquad \overline{\overline{e}} \circ W[f, R \circ \overline{\overline{e}}, \overline{\overline{e}}] \circ W[f, \overline{\overline{e}}, R \circ \overline{\overline{e}}]$$

*in general, does not coincide with the set of $\Omega_2$-words*

$$(6.2.5) \qquad \overline{\overline{e}} \circ W[f, \overline{\overline{e}}, \overline{\overline{e}}]$$

*The theorem 6.2.14 states that sets of $\Omega_2$-words (6.2.4) and (6.2.5) coincide up to equivalence generated by the representation $f$.*  □

Theorem 6.2.16. *Let $W[f, \overline{\overline{e}}, R \circ \overline{\overline{e}}]$ be the set of coordinates of automorphism $R$. Let $W[f, \overline{\overline{e}}, S \circ \overline{\overline{e}}]$ be the set of coordinates of automorphism $S$. The set of coordinates of automorphism $(R \circ S)^{-1}$ satisfies to the equality*

$$(6.2.6) \quad W[f, (R \circ S) \circ \overline{\overline{e}}, \overline{\overline{e}}] = W[f, S \circ (R \circ \overline{\overline{e}}), \overline{\overline{e}}] = W[f, S \circ \overline{\overline{e}}, \overline{\overline{e}}] \circ W[f, R \circ \overline{\overline{e}}, \overline{\overline{e}}]$$

Proof. The equality

$$(6.2.7) \quad \begin{aligned} W[f, (R \circ S) \circ \overline{\overline{e}}, \overline{\overline{e}}] &= W[f, \overline{\overline{e}}, (R \circ S)^{-1} \circ \overline{\overline{e}}] = W[f, \overline{\overline{e}}, S^{-1} \circ R^{-1} \circ \overline{\overline{e}}] \\ &= W[f, \overline{\overline{e}}, S^{-1} \circ \ e] \circ W[f, \overline{\overline{e}}, R^{-1} \circ \overline{\overline{e}}] \\ &= W[f, S \circ \overline{\overline{e}}, \overline{\overline{e}}] \circ W[f, R \circ \overline{\overline{e}}, \overline{\overline{e}}] \\ &= W[f, S \circ (R \circ \overline{\overline{e}}), \overline{\overline{e}}] \end{aligned}$$

follows from theorems 6.1.21, 6.2.14. The equality (6.2.6) follows from the equality (6.2.7).  □

Theorem 6.2.17. *The group of automorphisms $GA(f)$ of effective representation $f$ in $\Omega_2$-algebra $A_2$ generates effective left-side representation in $\Omega_2$-algebra $A_2$.*

Proof. From the corollary 6.2.11, it follows that if automorphism $R$ maps a basis $\overline{\overline{e}}$ into a basis $\overline{\overline{e}}'$, then the set of coordinates $W[f, \overline{\overline{e}}, \overline{\overline{e}}']$ uniquely determines an automorphism $R$. From the theorem 6.1.18, it follows that the set of coordinates $W[f, \overline{\overline{e}}, \overline{\overline{e}}']$ determines the map of coordinates relative to the basis $\overline{\overline{e}}$ under automorphism of the representation $f$. From the equation (6.1.28), it follows that automorphism $R$ acts from the left on elements of $\Omega_2$-algebra $A_2$. From the equation (6.1.27), it follows that the representation of group is left-side representation. According to the theorem 6.2.13, the set of coordinates $W[f, \overline{\overline{e}}, \overline{\overline{e}}]$ corresponds to identity transformation. From the theorem 6.2.14, it follows that the set of coordinates $W[f, R \circ \overline{\overline{e}}, \overline{\overline{e}}]$ corresponds to transformation, inverse to transformation $W[f, \overline{\overline{e}}, R \circ \overline{\overline{e}}]$.  □

## 6.3. Free Representation

In the section 3.1, we considered the definition 3.1.4 of free representation. However we may consider another definition which is similar to the definition of free module.

---

[6.13]In vector space, the matrix of numbers corresponds to linear transformation. Accordingly, the inverse matrix corresponds to inverse transformation.

DEFINITION 6.3.1. *The representation*

$$f : A_1 \longrightarrow_* A_2$$

*is called free representation if this representation has a basis.* ∎

THEOREM 6.3.2. *Let*

$$f : A_1 \longrightarrow_* A_2$$

*be free representation according to the definition 6.3.1. Then the representation $f$ is free according to the definition 3.1.4.*

PROOF. Let $\overline{\overline{e}}$ be basis of representation $f$ and $m \in \overline{\overline{e}}$. Let there exist $A_1$-numbers $a_1$, $b_1$ such that $f(a_1) = f(b_1)$. According to the assumption, $f(a_1)(m) = f(b_1)(m)$. However, if $a_1 \neq b_1$, then $f(a_1)(m)$ and $f(b_1)(m)$ are different $\Omega_2$-words. Therefore, $\overline{\overline{e}}$ is not a basis. From this contradiction, it follows that $a_1 = b_1$. Therefore, the representation $f$ is free according to the definition 3.1.4. ∎

THEOREM 6.3.3. *Let*

$$f : A_1 \longrightarrow_* A_2$$

*be free representation according to the definition 3.1.4. Then the representation $f$ is free according to the definition 6.3.1.*

QUESTION 6.3.4. *It is very important to find a proof of the theorem 6.3.3 or to find an example when this theorem is wrong. We will see in the chapter 7 how is important a free presentation defined in the definition 6.3.1. Because in the future I will assume that a representation always has a basis, then within the framework of this book I can use the theorem 6.3.2.* ∎

## 6.4. Basis Manifold of Representation

The set $\mathcal{B}[f]$ of bases of representation $f$ is called **basis manifold** of representation $f$.

DEFINITION 6.4.1. *According to theorems 6.1.20, 6.2.9, automorphism $R$ of the representation $f$ generates transformation*

(6.4.1)
$$R : \overline{\overline{h}} \to R \circ \overline{\overline{h}}$$
$$R \circ \overline{\overline{h}} = W[f, \overline{\overline{e}}, R \circ \overline{\overline{e}}] \circ \overline{\overline{h}}$$

*of the basis manifold of representation. This transformation is called **active**. According to the theorem 3.5.5, we defined left-side representation*

$$A(f) : GA(f) \longrightarrow_* \mathcal{B}[f]$$

*of group $GA(f)$ in basis manifold $\mathcal{B}[f]$. Representation $A(f)$ is called **active representation**. According to the corollary 6.2.11, this representation is single transitive.* ∎

REMARK 6.4.2. *According to remark 6.2.3, it is possible that there exist bases of representation $f$ such that there is no active transformation between them. Then we consider the orbit of selected basis as basis manifold. Therefore, it is possible that the representation $f$ has different basis manifolds. We will assume that we have chosen a basis manifold.*

Theorem 6.4.3. *There exists single transitive right-side representation*

$$P(f) : GA(f) \longmapsto \mathcal{B}[f]$$

*of group $GA(f)$ in basis manifold $\mathcal{B}[f]$. Representation $P(f)$ is called* **passive representation**.

Proof. Since $A(f)$ is single transitive left-side representation of group $GA(f)$, then single transitive right-side representation $P(f)$ is uniquely defined according to the theorem 5.5.9.  □

Theorem 6.4.4. *Transformation of representation $P(f)$ is called* **passive transformation of the basis manifold** *of representation. We also use notation*

$$S(\overline{\overline{e}}) = \overline{\overline{e}} \circ S$$

*to denote the image of basis $\overline{\overline{e}}$ under passive transformation $S$. Passive transformation of basis has form*

(6.4.2)
$$S : \overline{\overline{h}} \to \overline{\overline{h}} \circ S$$
$$\overline{\overline{h}} \circ S = \overline{\overline{h}} \circ W[f, \overline{\overline{e}}, \overline{\overline{e}} \circ S]$$

Proof. According to the equality (6.4.1), active transformation acts from left on coordinates of basis. The equality (6.4.2) follows from theorems 5.5.8, 5.5.9, 5.5.11; according to these theorems, passive transformation acts from right on coordinates of basis.  □

Theorem 6.4.5. *Passive transformation of the basis manifold is automorphism of representation $A(f)$.*

Proof. The theorem follows from the theorem 5.5.11.  □

Theorem 6.4.6. *Let $S$ be passive transformation of the basis manifold of the representation $f$. Let $\overline{\overline{e}}_1$ be the basis of the representation $f$, $\overline{\overline{e}}_2 = \overline{\overline{e}}_1 \circ S$. For basis $\overline{\overline{e}}_3$, let there exists an active transformation $R$ such that $\overline{\overline{e}}_3 = R \circ \overline{\overline{e}}_1$. Let $\overline{\overline{e}}_4 = R \circ \overline{\overline{e}}_2$. Then $\overline{\overline{e}}_4 = \overline{\overline{e}}_3 \circ S$.*

Proof. According to the equality (6.4.1), active transformation of coordinates of basis $\overline{\overline{e}}_3$ has form

(6.4.3)  $$\overline{\overline{e}}_4 = W[f, \overline{\overline{e}}_1, \overline{\overline{e}}_3] \circ \overline{\overline{e}}_2 = W[f, \overline{\overline{e}}_1, \overline{\overline{e}}_3] \circ \overline{\overline{e}}_1 \circ W[f, \overline{\overline{e}}_1, \overline{\overline{e}}_2]$$

Let $\overline{\overline{e}}_5 = \overline{\overline{e}}_3 \circ S$. From the equality (6.4.2), it follows that

(6.4.4)  $$\overline{\overline{e}}_5 = \overline{\overline{e}}_3 \circ W[f, \overline{\overline{e}}_1, \overline{\overline{e}}_2] = W[f, \overline{\overline{e}}_1, \overline{\overline{e}}_3] \circ \overline{\overline{e}}_1 \circ W[f, \overline{\overline{e}}_1, \overline{\overline{e}}_2]$$

From match of expressions in equalities (6.4.3), (6.4.4), it follows that $\overline{\overline{e}}_4 = \overline{\overline{e}}_5$. Therefore, the diagram

$$
\begin{array}{ccc}
\overline{\overline{e}}_1 \in \mathcal{B}[f] & \xrightarrow{\quad R \quad} & \overline{\overline{e}}_3 \in \mathcal{B}[f] \\
{\scriptstyle S} \downarrow & & \downarrow {\scriptstyle S} \\
\overline{\overline{e}}_2 \in \mathcal{B}[f] & \xrightarrow{\quad R \quad} & \overline{\overline{e}}_4 \in \mathcal{B}[f]
\end{array}
$$

is commutative.  □

### 6.5. Geometric Object of Representation of Universal Algebra

An active transformation changes a basis of the representation and $\Omega_2$-number uniformly and coordinates of $\Omega_2$-number relative basis do not change. A passive transformation changes only the basis and it leads to change of coordinates of $\Omega_2$-number relative to the basis.

THEOREM 6.5.1. *Let passive transformation $S \in GA(f)$ maps basis $\overline{\overline{e}}_1 \in \mathcal{B}[f]$ into basis $\overline{\overline{e}}_2 \in \mathcal{B}[f]$*

$$\overline{\overline{e}}_2 = \overline{\overline{e}}_1 \circ S = \overline{\overline{e}}_1 \circ W[f, \overline{\overline{e}}_1, \overline{\overline{e}}_1 \circ S] \tag{6.5.1}$$

*Let $A_2$-number $m$ has $\Omega_2$-word*

$$m = \overline{\overline{e}}_1 \circ W[f, \overline{\overline{e}}_1, m] \tag{6.5.2}$$

*relative to basis $\overline{\overline{e}}_1$ and has $\Omega_2$-word*

$$m = \overline{\overline{e}}_2 \circ W[f, \overline{\overline{e}}_2, m] \tag{6.5.3}$$

*relative to basis $\overline{\overline{e}}_2$. Coordinate transformation*

$$W[f, \overline{\overline{e}}_2, m] = W[f, \overline{\overline{e}}_1 \circ S, \overline{\overline{e}}_1] \circ W[f, \overline{\overline{e}}_1, m] \tag{6.5.4}$$

*does not depend on $A_2$-number $m$ or basis $\overline{\overline{e}}_1$, but is defined only by coordinates of $A_2$-number $m$ relative to basis $\overline{\overline{e}}_1$.*

PROOF. From (6.5.1) and (6.5.3), it follows that

$$\begin{aligned} \overline{\overline{e}}_1 \circ W[f, \overline{\overline{e}}_1, m] &= \overline{\overline{e}}_2 \circ W[f, \overline{\overline{e}}_2, m] = \overline{\overline{e}}_1 \circ W[f, \overline{\overline{e}}_1, \overline{\overline{e}}_2] \circ W[f, \overline{\overline{e}}_2, m] \\ &= \overline{\overline{e}}_1 \circ W[f, \overline{\overline{e}}_1, \overline{\overline{e}}_1 \circ S] \circ W[f, \overline{\overline{e}}_2, m] \end{aligned} \tag{6.5.5}$$

Comparing (6.5.2) and (6.5.5) we get

$$W[f, \overline{\overline{e}}_1, m] = W[f, \overline{\overline{e}}_1, \overline{\overline{e}}_1 \circ S] \circ W[f, \overline{\overline{e}}_2, m] \tag{6.5.6}$$

Since $S$ is automorphism, then the equality (6.5.4) follows from (6.5.6) and the theorem 6.2.14. $\square$

THEOREM 6.5.2. *Coordinate transformations (6.5.4) form effective contravariant right-side representation of group $GA(f)$ which is called* **coordinate representation** *in $\Omega_2$-algebra.*

PROOF. According to corollary 6.1.19, the transformation (6.5.4) is the endomorphism of representation [6.14]

$$f : A_1 \longrightarrow\!\!\!* \longrightarrow W[f, \overline{\overline{e}}_1]$$

Suppose we have two consecutive passive transformations $S$ and $T$. Coordinate transformation

$$W[f, \overline{\overline{e}}_2, m] = W[f, \overline{\overline{e}}_1 \circ S, \overline{\overline{e}}_1] \circ W[f, \overline{\overline{e}}_1, m] \tag{6.5.7}$$

corresponds to passive transformation $S$. Coordinate transformation

$$W[f, \overline{\overline{e}}_2, m] = W[f, \overline{\overline{e}}_1 \circ T, \overline{\overline{e}}_1] \circ W[f, \overline{\overline{e}}_1, m] \tag{6.5.8}$$

---

[6.14]This transformation does not generate an endomorphism of the representation $f$. Coordinates change because basis relative which we determinate coordinates changes. However, $A_2$-number, coordinates of which we are considering, does not change.

corresponds to passive transformation $T$. According to the theorem 8.3.3, product of coordinate transformations (6.5.7) and (6.5.8) has form

$$W[f, \overline{\overline{e}}_3, m] = W[f, \overline{\overline{e}}_1 \circ T, \overline{\overline{e}}_1] \circ W[f, \overline{\overline{e}}_1 \circ S, \overline{\overline{e}}_1] \circ W[f, \overline{\overline{e}}_1, m]$$
$$= W[f, \overline{\overline{e}}_1 \circ T \circ S, \overline{\overline{e}}_1] \circ W[f, \overline{\overline{e}}_1, m]$$

and is coordinate transformation corresponding to passive transformation $S \circ T$. According to theorems 6.2.14, 6.2.16 and to the definition 5.1.11, coordinate transformations form right-side contravariant representation of group $GA(f)$.

Suppose coordinate transformation does not change coordinates of selected basis. Then unit of group $GA(f)$ corresponds to it because representation is single transitive. Therefore, coordinate representation is effective.      □

Let $f$ be representation of $\Omega_1$-algebra $A_1$ in $\Omega_2$-algebra $A_2$. Let $g$ be representation of $\Omega_1$-algebra $A_1$ in $\Omega_3$-algebra $A_3$. Passive representation $P(g)$ is coordinated with passive representation $P(f)$, if there exists homomorphism $h$ of group $GA(f)$ into group $GA(g)$. Consider diagram

$$\begin{array}{ccc} \text{End}(B[f]) & \xrightarrow{\ H\ } & \text{End}(B[g]) \\ {\scriptstyle P(f)}\Big\uparrow & {\scriptstyle f'} \nearrow & \Big\uparrow{\scriptstyle P(g)} \\ GA(f) & \xrightarrow{\ h\ } & GA(g) \end{array}$$

Since maps $P(f)$, $P(g)$ are isomorphisms of group, then map $H$ is homomorphism of groups. Therefore, map $f'$ is representation of group $GA(f)$ in basis manifold $\mathcal{B}(g)$. According to design, passive transformation $H(S)$ of basis manifold $\mathcal{B}(g)$ corresponds to passive transformation $S$ of basis manifold $\mathcal{B}(f)$

$$(6.5.9) \qquad\qquad \overline{\overline{e}}_{g1} = \overline{\overline{e}}_g \circ H(S)$$

Then coordinate transformation in representation $g$ gets form

$$(6.5.10) \qquad W[g, \overline{\overline{e}}_{g1}, m] = W[g, \overline{\overline{e}}_g \circ H(S), \overline{\overline{e}}_g] \circ W[g, \overline{\overline{e}}_g, m]$$

**Definition 6.5.3.** *Orbit*

$$\mathcal{O}(f, g, \overline{\overline{e}}_g, m) = H(GA(f)) \circ W[g, \overline{\overline{e}}_g, m]$$
$$= (W[g, \overline{\overline{e}}_g \circ H(S), \overline{\overline{e}}_g] \circ W[g, \overline{\overline{e}}_g, m], \overline{\overline{e}}_f \circ S, S \in GA(f))$$

*is called* **geometric object in coordinate representation** *defined in the representation $f$. For any basis $\overline{\overline{e}}_{f1} = \overline{\overline{e}}_f \circ S$ corresponding point (6.5.10) of orbit defines* **coordinates of geometric object** *relative basis $\overline{\overline{e}}_{f1}$.*      □

**Definition 6.5.4.** *Orbit*

$$\mathcal{O}(f, g, m) = (W[g, \overline{\overline{e}}_g \circ H(S), \overline{\overline{e}}_g] \circ W[g, \overline{\overline{e}}_g, m], \overline{\overline{e}}_g \circ H(S), \overline{\overline{e}}_f \circ S, S \in GA(f))$$

*is called* **geometric object** *defined in the representation $f$. We also say that $m$ is a* **geometric object of type** $H$. *For any basis $\overline{\overline{e}}_{f1} = \overline{\overline{e}}_f \circ S$ corresponding point (6.5.10) of orbit defines $A_2$-number*

$$m = \overline{\overline{e}}_g \circ W[g, \overline{\overline{e}}_g, m]$$

*called* **representative of geometric object** *in the representation $f$.*      □

Since a geometric object is an orbit of representation, we see that according to the theorem 5.3.7 the definition of the geometric object is a proper definition.

Definition 6.5.3 introduces a geometric object in coordinate space. We assume in definition 6.5.4 that we selected a basis of representation $g$. This allows using a representative of the geometric object instead of its coordinates.

THEOREM 6.5.5 (invariance principle). *Representative of geometric object does not depend on selection of basis $\overline{\overline{e}}_f$.*

PROOF. To define representative of geometric object, we need to select basis $\overline{\overline{e}}_f$ of representation $f$, basis $\overline{\overline{e}}_g$ of representation $g$ and coordinates of geometric object $W[g, \overline{\overline{e}}_g, n]$. Corresponding representative of geometric object has form

$$n = \overline{\overline{e}}_g \circ W[g, \overline{\overline{e}}_g, n]$$

Suppose we map basis $\overline{\overline{e}}_f$ to basis $\overline{\overline{e}}_{f1}$ by passive transformation

$$\overline{\overline{e}}_{f1} = \overline{\overline{e}}_f \circ S$$

According building this forms passive transformation (6.5.9) and coordinate transformation (6.5.10). Corresponding representative of geometric object has form

$$\begin{aligned}
n' &= \overline{\overline{e}}_{g1} \circ W[g, \overline{\overline{e}}_{g1}, n'] \\
&= \overline{\overline{e}}_g \circ W[g, \overline{\overline{e}}_g, \overline{\overline{e}}_g \circ H(S)] \circ W[g, \overline{\overline{e}}_g \circ H(S), \overline{\overline{e}}_g] \circ W[g, \overline{\overline{e}}_g, n] \\
&= \overline{\overline{e}}_g \circ W[g, \overline{\overline{e}}_g, n] = n
\end{aligned}$$

Therefore representative of geometric object is invariant relative selection of basis. □

THEOREM 6.5.6. *The set of geometric objects of type $H$ is $\Omega_3$-algebra.*

PROOF. Let

$$m_i = \overline{\overline{e}}_g \circ W[g, \overline{\overline{e}}_g, m_i] \quad i = 1, ..., n$$

For operation $\omega \in \Omega_3(n)$ we assume

(6.5.11)         $$m_1...m_n\omega = \overline{\overline{e}}_g \circ (W[g, \overline{\overline{e}}_g, m_1]...W[g, \overline{\overline{e}}_g, m_n]\omega)$$

Since for arbitrary endomorphism $S$ of $\Omega_2$-algebra $A_2$, the map $W[g, \overline{\overline{e}}_g, \overline{\overline{e}}_g \circ H(S)]$ is endomorphism of $\Omega_3$-algebra $A_3$, then the definition (6.5.11) is correct. □

THEOREM 6.5.7. *There exists the representation of $\Omega_1$-algebra $A_1$ in $\Omega_3$-algebra $N$ of geometric objects of type $H$.*

PROOF. Let

$$m = \overline{\overline{e}}_g \circ W[g, \overline{\overline{e}}_g, m]$$

For $a \in A_1$, we assume

(6.5.12)         $$f(a)(m) = \overline{\overline{e}}_g \circ f(a)(W[g, \overline{\overline{e}}_g, m])$$

Since for arbitrary endomorphism $S$ of $\Omega_2$-algebra $A_2$, the map $W[g, \overline{\overline{e}}_g, \overline{\overline{e}}_g \circ H(S)]$ is endomorphism of representation $g$, then the definition (6.5.12) is correct. □

# Diagram of Representations of Universal Algebras

## 7.1. Diagram of Representations of Universal Algebras

From a comparison of theorems 6.1.4 and [14]-5.1, it follows that there is no rigid boundary between universal algebra and representation of universal algebra. This implies the possibility of a generalization of representation of universal algebra.

The simplest construction arises as follows. Let

$$f_{12} : A_1 \xrightarrow{\quad\ast\quad} A_2$$

be representation of $\Omega_1$-algebra $A_1$ in $\Omega_2$-algebra $A_2$. If, instead of $\Omega_2$-algebra $A_2$, we consider representation

$$f_{23} : A_2 \xrightarrow{\quad\ast\quad} A_3$$

of $\Omega_2$-algebra $A_2$ in $\Omega_3$-algebra $A_3$, then we get a diagram of the following form

$$(7.1.1) \qquad A_1 \xrightarrow{f_{12}} A_2 \xrightarrow{f_{23}} A_3$$

It is evident that, in the diagram (7.1.1), we assume that $A_3$ is representation

$$f_{34} : A_3 \xrightarrow{\quad\ast\quad} A_4$$

We can make a chain of representations of universal algebras as long as we wish. Thus we obtain the following definition.

DEFINITION 7.1.1. *Consider set of $\Omega_k$-algebras $A_k$, $k = 1, ..., n$. Let $A = (A_1, ..., A_n)$. Let $f = (f_{1\,2}, ..., f_{n-1\,n})$. Set of representations $f_{k\,k+1}$, $k = 1, ..., n$, of $\Omega_k$-algebra $A_k$ in $\Omega_{k+1}$-algebra $A_{k+1}$ is called* **tower** $(f, A)$ **of representations of $\Omega$-algebras.** $\qquad\square$

We can represent tower of representations $(f, A)$ using the following diagram

$$A_1 \xrightarrow{f_{1\,2}} A_2 \xrightarrow{f_{2\,3}} ... \xrightarrow{f_{n-1\,n}} A_n$$

When we consider the tower of representations, we again consider that $A_2$ or $A_3$ are representations of universal algebras or towers of representations. In this case, the diagram (7.1.1) gets form

$$A_1 \xrightarrow{f_{12}} A_2 \xrightarrow{f_{23}} A_3$$
$$\uparrow f_{42} \qquad \uparrow f_{53}$$
$$A_4 \qquad\qquad A_5$$
$$\uparrow f_{64}$$
$$A_6$$

or

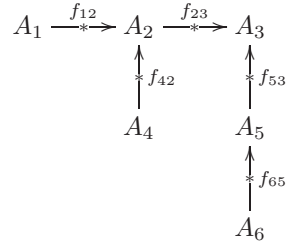

We also assume that some algebras or maps on the diagram coincide. Thus, we say that diagrams

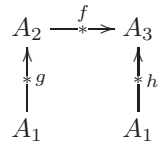

and

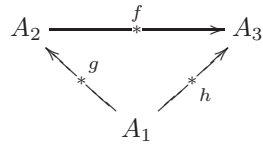

are equivalent.

DEFINITION 7.1.2. **Diagram** $(f, A)$ **of representations of universal algebras** *is oriented graph such that*

7.1.2.1: *the vertex of $A_k$, $k = 1, ..., n$, is $\Omega_k$-algebra;*

7.1.2.2: *the edge $f_{kl}$ is representation of $\Omega_k$-algebra $A_k$ in $\Omega_l$-algebra $A_l$;*

*We require that this graph is connected graph and does not have loops. Let $A_{[0]}$ be set of initial vertices of the graph. Let $A_{[k]}$ be set of vertices of the graph for which the maximum path from the initial vertices is $k$.* □

REMARK 7.1.3. *Since different vertices of the graph can be the same algebra, then we denote $A = (A_{(1)} \ ... \ A_{(n)})$ the set of universal algebras which are distinct. From the equality*

$$A = (A_{(1)} \ ... \ A_{(n)}) = (A_1 \ ... \ A_n)$$

*it follows that, for any index $(i)$, there exists at least one index $i$ such that $A_{(i)} = A_i$. If there are two sets of sets $A = (A_{(1)} \ ... \ A_{(n)})$, $B = (B_{(1)} \ ... \ B_{(n)})$ and there is a map*

$$h_{(i)} : A_{(i)} \to B_{(i)}$$

*for an index $(i)$, then also there is a map*

$$h_i : A_i \to B_i$$

*for any index $i$ such that $A_{(i)} = A_i$ and in this case $h_i = h_{(i)}$.* □

THEOREM 7.1.4 (induction over diagram of representations). *Let the theorem $\mathcal{T}$ be true for the set of universal algebras $A_{[0]}$ of diagram $(f, A)$ of representations of universal algebras. Let the statement that the theorem $\mathcal{T}$ is true for the set of universal algebras $A_{[k]}$ of diagram $(f, A)$ of representations imply the statement that the theorem $\mathcal{T}$ is true for the set of universal algebras $A_{[k+1]}$ of diagram $(f, A)$ of* *representations. Then the theorem $\mathcal{T}$ is true for the set of universal algebras of diagram $(f, A)$ of representations.*

Proof. The theorem follows from the principle of mathematical induction. $\square$

Definition 7.1.5. *Diagram $(f, A)$ of representations of universal algebras is called* **commutative** *when diagram meets the following requirement. for each pair of representations*

$$f_{ik} : A_i \longrightarrow A_k$$
$$f_{jk} : A_j \longrightarrow A_k$$

*the following equality is true* [7.1]

$$(7.1.2) \qquad f_{ik}(a_i)(f_{jk}(a_j)(a_k)) = f_{jk}(a_j)(f_{ik}(a_i)(a_k))$$

$\square$

Theorem 7.1.6. *Let*

$$f_{ij} : A_i \longrightarrow A_j$$

*be representation of $\Omega_i$-algebra $A_i$ in $\Omega_j$-algebra $A_j$. Let*

$$f_{jk} : A_j \longrightarrow A_k$$

*be representation of $\Omega_j$-algebra $A_j$ in $\Omega_k$-algebra $A_k$. We represent the fragment* [7.2]

$$A_i \overset{f_{ij}}{\longrightarrow} A_j \overset{f_{jk}}{\longrightarrow} A_k$$

*of the diagram of representations using the diagram*

(7.1.3)

*The map*

$$f_{ijk} : A_i \to \operatorname{End}(\Omega_j, \operatorname{End}(\Omega_k, A_k))$$

*is defined by the equality*

$$(7.1.4) \qquad f_{ijk}(a_i)(f_{jk}(a_j)) = f_{jk}(f_{ij}(a_i)(a_j))$$

---

[7.1] Metaphorically speaking, representations $f_{ik}$ and $f_{jk}$ are transparent to each other.

[7.2] The theorem 7.1.6 states that transformations in diagram of representations are coordinated.

*where $a_i \in A_i$,  $a_j \in A_j$. If the representation $f_{jk}$ is effective and the representation $f_{ij}$ is free, then the map $f_{ijk}$ is free representation*

$$f_{ijk} : A_i \longrightarrow \text{End}(\Omega_k, A_k)$$

*of $\Omega_i$-algebra $A_i$ in $\Omega_j$-algebra  $\text{End}(\Omega_k, A_k)$.*

Proof.

Lemma 7.1.7. *The map $f_{ijk}$ is injection.*

Proof. Let  $(a_i, b_i) \in \ker f_{ijk}$. Then

$$(7.1.5) \qquad f_{jk}(f_{ij}(a_i)(a_j)) = f_{ijk}(a_i)(f_{jk}(a_j)) = f_{ijk}(b_i)(f_{jk}(a_j))$$
$$= f_{jk}(f_{ij}(b_i)(a_j))$$

If the representation $f_{jk}$ is effective, then the equality

$$(7.1.6) \qquad f_{ij}(a_i)(a_j) = f_{ij}(b_i)(a_j)$$

follows from the definition 3.1.2 and from the equality (7.1.5) for any  $a_j \in A_j$. The statement $a_i = b_i$ follows from the definition 3.1.4.               ⊙

Lemma 7.1.8. *There is the structure of $\Omega_i$-algebra on the set $\text{End}(\Omega_j, \text{End}(\Omega_k, A_k))$.*

Proof. Let  $\omega \in \Omega_i$. Let $a_1, ..., a_m \in A_i$. We define operation $\omega$ on the set $\text{End}(\Omega_j, \text{End}(\Omega_k, A_k))$  using the following equality

$$(7.1.7) \qquad f_{ijk}(a_1)...f_{ijk}(a_m)\omega = f_{ijk}(a_1...a_m\omega)$$

According to the lemma 7.1.7, the equality (7.1.7) properly defines operation $\omega$. ⊙

Corollary 7.1.9. *The map $f_{ijk}$ is homomorphism of $\Omega_i$-algebra.*      ⊙

Lemma 7.1.10. *The map $f_{ijk}(a)$ is homomorphism of $\Omega_j$-algebra.*

Proof. Let $b_1, ..., b_m \in A_j$. Then the equality

$$(7.1.8) \quad f_{ijk}(a)(f_{jk}(b_1))...f_{ijk}(a)(f_{jk}(b_m))\omega = f_{jk}(f_{ij}(a)(b_1))...f_{jk}(f_{ij}(a)(b_m))\omega$$

follows from the equality (7.1.4). Since maps  $f_{ij}(a)$, $f_{jk}$  are homomorphisms of $\Omega_j$-algebra, then the equality

$$(7.1.9) \qquad f_{ijk}(a)(f_{jk}(b_1))...f_{ijk}(a)(f_{jk}(b_m))\omega$$
$$= f_{jk}(f_{ij}(a)(b_1)...f_{ij}(a)(b_m)\omega)$$
$$= f_{jk}(f_{ij}(a)(b_1...b_m\omega))$$

follows from the equality (7.1.8). The equality

$$(7.1.10) \qquad f_{ijk}(a)(f_{jk}(b_1))...f_{ijk}(a)(f_{jk}(b_m))\omega = f_{ijk}(a)(f_{jk}(b_1...b_m\omega))$$

follows from the equalities (7.1.4), (7.1.9). Since the map $f_{jk}$ is homomorphism of $\Omega_2$-algebra, then the equality

$$(7.1.11) \qquad f_{ijk}(a)(f_{jk}(b_1))...f_{ijk}(a)(f_{jk}(b_m))\omega = f_{ijk}(a)(f_{jk}(b_1)...f_{jk}(b_m)\omega)$$

follows from the equality (7.1.10).                                ⊙

The theorem follows from corollary 7.1.9 and from the lemma 7.1.10.       □

Theorem 7.1.11. *The map $f_{jk}$ is reduced morphism of representations from $f_{ij}$ into $f_{ijk}$.*

Proof. Consider diagram (7.1.3) in more detail.

(7.1.12)

$$A_j \xrightarrow{\quad f_{jk} \quad} \mathrm{End}(\Omega_k, A_k)$$

$$f_{ij} \qquad A_i \qquad f_{ijk}$$

The statement of theorem follows from the equality (7.1.4) and from the definition 3.4.2. □

Theorem 7.1.12. *Let*

$$f_{ij} : A_i \longrightarrow A_j$$

*be representation of* $\Omega_i$*-algebra* $A_i$ *in* $\Omega_j$*-algebra* $A_j$. *Let*

$$f_{jk} : A_j \longrightarrow A_k$$

*be representation of* $\Omega_j$*-algebra* $A_j$ *in* $\Omega_k$*-algebra* $A_k$. *Then there exists representation*

$$f_{ij,k} : A_i \times A_j \longrightarrow A_k$$

*of the set*[7.3] $A_i \times A_j$ *in* $\Omega_k$*-algebra* $A_k$.

Proof. We represent the fragment

$$A_i \xrightarrow{\quad f_{ij} \quad} A_j \xrightarrow{\quad f_{jk} \quad} A_k$$

of the diagram of representations using the diagram

(7.1.13)

$$A_i \times A_j$$

$$f_{ij,k}$$

$$A_k \xrightarrow{\hspace{5cm}} A_k$$

$$f_{jk}(f_{ij}(a_i)(a_j))$$

$$f_{jk}$$

$$A_j \xrightarrow{\quad f_{ij}(a_i) \quad} A_j$$

$$f_{ij}$$

$$A_i$$

From the diagram (7.1.13), it follows that the map $f_{ij,k}$ is defined by the equality

$$f_{ij,k}(a_i, a_j) = f_{jk}(f_{ij}(a_i)(a_j))$$

□

---

[7.3] Since $\Omega_i$-algebra $A_i$ and $\Omega_j$-algebra $A_j$ have different set of operations, we cannot define the structure of universal algebra on the set $A_i \times A_j$.

## 7.2. Morphism of Diagram of Representations

Definition 7.2.1. *Let $(f, A)$ be the diagram of representations where $A = (A_{(1)} \ ... \ A_{(n)})$ is the set of universal algebras. Let $(B, g)$ be the diagram of representations where $B = (B_{(1)} \ ... \ B_{(n)})$ is the set of universal algebras. The set of maps $h = (h_{(1)} \ ... \ h_{(n)})$*

$$h_{(i)} : A_{(i)} \to B_{(i)}$$

*is called* **morphism from diagram of representations** $(f, A)$ **into diagram of representations** $(B, g)$, *if for any indexes $(i)$, $(j)$, $i$, $j$ such that $A_{(i)} = A_i$, $A_{(j)} = A_j$ and for any representation*

$$f_{ji} : A_j \longrightarrow{}^* A_i$$

*the tuple of maps $(h_j \ \ h_i)$ is morphism of representations from $f_{ji}$ into $g_j i$.* □

We will use notation

$$h : A \to B$$

if tuple of maps $h$ is morphism from diagram of representations $(f, A)$ into diagram of representations $(B, g)$.

When studying morphism of the representation of universal algebra, we very often assume that algebra generating representation is given. Therefore, we are not interested in the map of this algebra; this convention simplifies the structure of morphism. Such morphism of representation we call reduced morphism of representation.

We see a similar problem when we study morphism of diagram of representations. For each universal algebra from diagram of representations, there exists the set of algebras preceding this algebra in corresponding graph. We can assume that some of these algebras are given and we will not consider corresponding homomorphisms. Corresponding morphism of diagram of representations also is called reduced. However, because diagram of representations is complicated structure, we will not consider reduced morphism of diagram of representations.

For any representation $f_{ij}$, $i = 1, ..., n$, $j = 1, ..., n$, we have diagram

(7.2.1)

Equalities

(7.2.2)
$$h_j \circ f_{ij}(a_i) = g_{ij}(h_i(a_i)) \circ h_j$$

(7.2.3)
$$h_j(f_{ij}(a_i)(a_j)) = g_{ij}(h_i(a_i))(h_j(a_j))$$

express commutativity of diagram (1).

Let representations $f_{ij}$ and $f_{jk}$ of universal algebras be defined. Assuming diagram (7.2.1) for representations $f_{ij}$ and $f_{jk}$ we will get the following diagram

(7.2.4)

It is evident that there exists morphism from $\operatorname{End}(\Omega_k, A_k)$ into $\operatorname{End}(\Omega_k, B_k)$, which maps $f_{ijk}(a_i)$ into $g_{ijk}(h_k(a_i))$.

THEOREM 7.2.2. *If the representation $f_{jk}$ is effective and the representation $f_{ij}$ is free, then*[7.4] *$(h_i, h_k^*)$ is morphism of representations from representation $f_{ijk}$ into representation $g_{ijk}$ of $\Omega_i$-algebra.*

PROOF. Consider the diagram

The existence of map $h_k^*$ and commutativity of the diagram (2) and (3) follows from effectiveness of map $f_{jk}$ and theorem 3.2.9. Commutativity of diagrams (4) and (5) follows from theorem 7.1.11.

From commutativity of the diagram (4) it follows that

(7.2.5) $$f_{jk} \circ f_{ij}(a_i) = f_{ijk}(a_i) \circ f_{jk}$$

From the equalitiy (7.2.5) it follows that

(7.2.6) $$h_k^* \circ f_{jk} \circ f_{ij}(a_i) = h_k^* \circ f_{ijk}(a_i) \circ f_{jk}$$

---

[7.4] See the definition of the map $h^*$ in the theorem 3.2.9.

From commutativity of diagram (3) it follows that

$$(7.2.7) \qquad h_k^* \circ f_{jk} = g_{jk} \circ h_j$$

From the equalitiy (7.2.7) it follows

$$(7.2.8) \qquad h_k^* \circ f_{jk} \circ f_{ij}(a_i) = g_{jk} \circ h_j \circ f_{ij}(a_i)$$

From equalities (7.2.6) and (7.2.8) it follows that

$$(7.2.9) \qquad h_k^* \circ f_{ijk}(a_i) \circ f_{jk} = g_{jk} \circ h_j \circ f_{ij}(a_i)$$

From commutativity of the diagram (5) it follows that

$$(7.2.10) \qquad g_{jk} \circ g_{ij}(h_i(a_i)) = g_{ijk}(h_i(a_i)) \circ g_{jk}$$

From the equalitiy (7.2.10) it follows that

$$(7.2.11) \qquad g_{jk} \circ g_{ij}(h_i(a_i)) \circ h_j = g_{ijk}(h_i(a_i)) \circ g_{jk} \circ h_j$$

From commutativity of the diagram (2) it follows that

$$(7.2.12) \qquad h_k^* \circ f_{jk} = g_{jk} \circ h_j$$

From the equalitiy (7.2.12) it follows that

$$(7.2.13) \qquad g_{ijk}(h_i(a_i)) \circ h_k^* \circ f_{jk} = g_{ijk}(h_i(a_i)) \circ g_{jk} \circ h_j$$

From equalities (7.2.11) and (7.2.13) it follows that

$$(7.2.14) \qquad g_{jk} \circ g_{ij}(h_i(a_i)) \circ h_j = g_{ij}(h_i(a_i)) \circ h_k^* \circ f_{jk}$$

External diagram is diagram (7.2.1) when $i = 1$. Therefore, external diagram is commutative

$$(7.2.15) \qquad h_j \circ f_{ij}(a_i) = g_{ij}(h_i(a_i)) \circ h_j$$

From the equalitiy (7.2.15) it follows that

$$(7.2.16) \qquad g_{jk} \circ h_J \circ f_{ij}(a_i) = g_{jk} \circ g_{ij}(h_i(a_i)) \circ h_j(a_j)$$

From equalities (7.2.9), (7.2.14) and (7.2.16) it follows that

$$(7.2.17) \qquad h_k^* \circ f_{ijk}(a_i) \circ f_{jk} = g_{ijk}(h_i(a_i)) \circ h_k^* \circ f_{jk}$$

Because the map $f_{i+1,i+2}$ is injection, then from the equalitiy (7.2.17) it follows that

$$(7.2.18) \qquad h_k^* \circ f_{ijk}(a_i) = g_{ijk}(h_i(a_i)) \circ h_k^*$$

From the equalitiy (7.2.18) commutativity of the diagram (1) follows. This proves the statement of theorem.    $\square$

The theorem 7.2.2 states that unknown map on the diagram (7.2.4) is the map $h_k^*$. Meaning of theorems 7.1.11 and 7.2.2 is that all maps in diagram of representations act coherently when all representations are free.

THEOREM 7.2.3. *Consider the set of $\Omega_i$-algebras $A_{(i)}$, $B_{(i)}$, $C_{(i)}$, $(i) = (1), ..., (n)$. Let*

$$p : (f, A) \to (g, B)$$
$$q : (g, B) \to (h, C)$$

*be morphisms of diagrams of representations. There exists morphism of representations of $\Omega$-algebra*

$$r : (f, A) \to (h, C)$$

*where $r_{(k)} = q_{(k)} \circ p_{(k)}$, $(k) = (1), ..., (n)$. We call morphism $r$ of diagram of representations from $f$ into $h$* **product of morphisms $p$ and $q$ of diagram of representations**.

PROOF. For any $i$, $j$ such that $A_{(j)} = A_j$, if there exists the representation $f_{ij}$, we represent statement of theorem using diagram

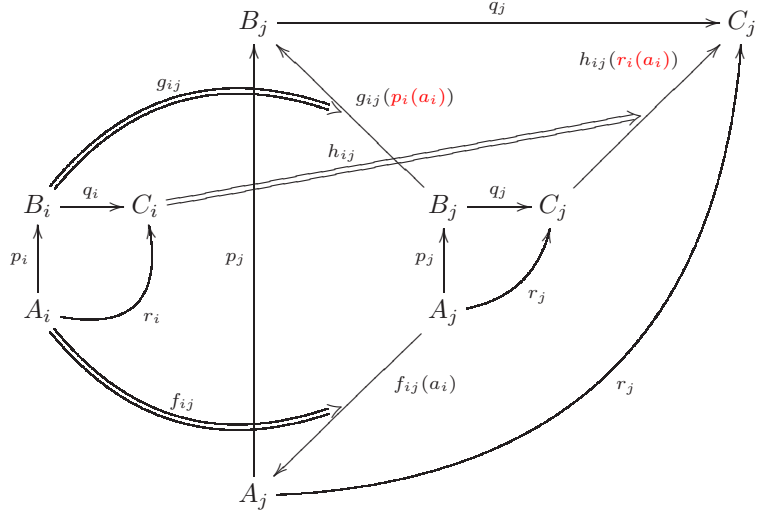

Map $r_i$ is homomorphism of $\Omega_i$-algebra $A_i$ into $\Omega_i$-algebra $C_i$. We need to show that tuple of maps $(r_i, r_j)$ satisfies to (7.2.2):

$$\begin{aligned}
r_k(f_{ij}(a_i)(a_j)) &= q_j \circ p_j(f_{ij}(a_i)(a_j)) \\
&= q_k(g_{ij}(p_i(a_i))(p_j(a_j))) \\
&= h_{ij}(q_i \circ p_i(a_i))(q_j \circ p_j(a_j)) \\
&= h_{ij}(r(a_i))(r_j(a_j))
\end{aligned}$$

$\square$

## 7.3. Automorphism of Diagram of Representations

DEFINITION 7.3.1. *Let $(f, A)$ be diagram of representations of universal algebras. The morphism of diagram of representations $(h_{(1)}, ..., h_{(n)})$ such, that for each $(k)$, $(k) = (1), ..., (n)$, $h_{(k)}$ is endomorphism of $\Omega_{(k)}$-algebra $A_{(k)}$ is called* **endomorphism of diagram of representations**. $\square$

DEFINITION 7.3.2. *Let $(f, A)$ be diagram of representations of universal algebras. The morphism of diagram of representations $(h_{(1)}, ..., h_{(n)})$ such, that for each $(k)$, $(k) = (1), ..., (n)$, $h_{(k)}$ is automorphism of $\Omega_{(k)}$-algebra $A_{(k)}$ is called* **automorphism of diagram of representations**. $\square$

THEOREM 7.3.3. *Let $(f, A)$ be diagram of representations of universal algebras. The set of automorphisms of the diagram of representations $(f, A)$ forms* **group** *$GA(f)$.*

PROOF. Let $r = (r_{(1)}, ..., r_{(n)})$, $p = (p_{(1)}, ..., p_{(n)})$ be automorphisms of the diagram of representations $(f, A)$. According to definition 7.3.2, for each $(k)$, $(k) = (1), ..., (n)$, maps $r_{(k)}$, $p_{(k)}$ are automorphisms of $\Omega_{(k)}$-algebra $A_{(k)}$.

According to theorem II.3.2 ([14], p. 57), for each $(k)$, $(k) = (1)$, ..., $(n)$, the map $r_{(k)} \circ p_{(k)}$ is automorphism of $\Omega_{(k)}$-algebra $A_{(k)}$. From the theorem 7.2.3 and the definition 7.3.2, it follows that product of automorphisms $r \circ p$ of the diagram of representations $(f, A)$ is automorphism of the diagram of representations $(f, A)$.

According to proof of the theorem 3.5.5, for any $(k)$, $(k) = (1)$, ..., $(n)$, the product of automorphisms of $\Omega_{(k)}$-algebra is associative. Therefore, the product of automorphisms of diagram of representations is associative.

Let $r = (r_{(1)}, ..., r_{(n)})$ be an automorphism of the diagram of representations $(f, A)$. According to definition 7.3.2 for each $(k)$, $(k) = (1)$, ..., $(n)$, the map $r_{(k)}$ is automorphism of $\Omega_{(k)}$-algebra $A_{(k)}$. Therefore, for each $(k)$, $(k) = (1)$, ..., $(n)$, the map $r_{(k)}^{-1}$ is automorphism of $\Omega_{(k)}$-algebra $A_{(k)}$. Then the equalitiy (7.2.3) is true for automorphism $r = (r_{(1)}, ..., r_{(n)})$. Let $j$ be index such that $j = (k)$. Let $a_j' = r_j(a_j)$. Since $r_j$ is automorphism then $a_j = r_j^{-1}(a_j')$ and, for any $i$, $j$, in case that there exists representation $f_{ij}$, we can write (7.2.3) in the form

$$(7.3.1) \qquad h_j(f_{ij}(h_i^{-1}(a_i'))(h_j^{-1}(a_j'))) = g_{ij}(a_i')(a_j')$$

Similarly, from the equalitiy (7.3.1) it follows that

$$(7.3.2) \qquad f_{ij}(h_i^{-1}(a_i'))(h_j^{-1}(a_j')) = h_j^{-1}(g_{ij}(a_i')(a_j'))$$

The equalitiy (7.3.2) corresponds to the equalitiy (7.2.3) for the map $r^{-1}$. Therefore, map $r^{-1}$ is automorphism of the diagram of representations $(f, A)$.    $\square$

# Basis of Diagram of Representations of Universal Algebra

## 8.1. Generating Set of Diagram of Representations

We construct the basis of the diagram of representations in a similar way that we constructed the basis of representation in the section 6.2.

DEFINITION 8.1.1. *Let $(f, A)$ be diagram of representations. The tuple of sets*

$$N = (N_{(1)} \subset A_{(1)}, ..., N_{(n)} \subset A_{(n)}) = (N_1 \subset A_1, ..., N_n \subset A_n)$$

*is called* **tuple of stable sets of diagram of representations** $(f, A)$, *if*

$$f_{ij}(a_i)(a_j) \in N_j \quad i, j = 1, ..., n$$

*for every* $a_1 \in N_1,$ ..., $a_n \in N_n,$ *in case that there exists the representation $f_{ij}$. We also will say that tuple of sets*

$$N = (N_{(1)} \subset A_{(1)}, ..., N_{(n)} \subset A_{(n)}) = (N_1 \subset A_1, ..., N_n \subset A_n)$$

*is stable relative to diagram of representations* $(f, A)$. □

THEOREM 8.1.2. *Let $(f, A)$ be diagram of representations. Let set $N_{(i)} \subset A_{(i)}$ be subalgebra of $\Omega_i$-algebra $A_{(i)}$, $(i) = (1),$ ..., $(n)$. Let tuple of sets*

$$N = (N_{(1)} \subset A_{(1)}, ..., N_{(n)} \subset A_{(n)}) = (N_1 \subset A_1, ..., N_n \subset A_n)$$

*be stable relative to diagram of representations* $(f, A)$. *Then there exists diagram of representations*

$$(8.1.1) \qquad\qquad (N, f_N = (f_{Nij}))$$

*such that*

$$f_{Nij}(a_i) = f_{ij}(a_i)|_{N_j} \quad i = 1, ..., n \quad j = 1, ..., n$$

*The diagram of representations* (8.1.1) *is called* **subrepresentation** *of diagram of representations* $(f, A)$.

PROOF. Let $\omega$ be $m$-ary operation of $\Omega_i$-algebra $A_i$, $i = 1,$ ..., $n$. Then for any $a_{i,1},$ ..., $a_{i,m} \in N_i$ and any $a_j \in N_j$

$$\begin{aligned}
(f_{Nij}(a_{i,1})...f_{Nij}(a_{i,m})\omega)(a_j) &= (f_{ij}(a_{i,1})...f_{ij}(a_{i,m})\omega)(a_j) \\
&= f_{ij}(a_{i,1}...a_{i,m}\omega)(a_j) \\
&= f_{Nij}(a_{i,1}...a_{i,m}\omega)(a_i)
\end{aligned}$$

Let $\omega$ be $m$-ary operation of $\Omega_j$-algebra $A_j$, $j = 1, ..., n$. Then for any $a_{j,1}, ..., a_{j,m} \in N_j$ and any $a_i \in N_i$

$$
\begin{aligned}
f_{Nij}(a_i)(a_{j,1})...f_{Nij}(a_i)(a_{i,m})\omega &= f_{ij}(a_i)(a_{j,1})...f_{ij}(a_i)(a_{j,m})\omega \\
&= f_{ij}(a_i)(a_{j,1}...a_{j,m}\omega) \\
&= f_{Nij}(a_i)(a_{j,1}...a_{j,m}\omega)
\end{aligned}
$$

We proved the statement of theorem.                                                    □

From theorem 8.1.2, it follows that if diagram of representations (8.1.1) is diagram of subrepresentations of diagram of representations $(f, A)$, then map

$$(id_{(1)} : N_{(1)} \to A_{(1)}, ..., id_{(n)} : N_{(n)} \to A_{(n)})$$

is morphism of diagrams of representations.

THEOREM 8.1.3. *The set[8.1] $\mathcal{B}[f, A]$ of all diagrams of subrepresentations of diagram of representations $(f, A)$ generates a closure system on diagram of representations $(f, A)$ and therefore is a complete lattice.*

PROOF. Let for given $\lambda \in \Lambda$,

$$K_\lambda = (K_{\lambda,(1)} \subset A_{(1)}, ..., K_{\lambda,(n)} \subset A_{(n)})$$

be tuple of sets stable relative to diagram of representations $(f, A)$. We determine the operation of intersection on the set $\mathcal{B}[f, A]$ according to rule

$$\bigcap f_{K_\lambda ij} = f_{\cap K_\lambda ij} \quad i, j = 1, ..., n$$

$$\bigcap K_\lambda = \left( K_{(1)} = \bigcap K_{\lambda,(1)}, ..., K_{(n)} = \bigcap K_{\lambda,(n)} \right)$$

$\cap K_{\lambda,(i)}$ is subalgebra of $\Omega_{(i)}$-algebra $A_{(i)}$. Let

$$a_j \in \bigcap K_{\lambda,j} = \bigcap K_{\lambda,(j)}$$

For any $\lambda \in \Lambda$ and for any $a_i \in K_i = K_{(i)}$

$$f_{ij}(a_i)(a_j) \in K_{\lambda,j} = K_{\lambda,(j)}$$

Therefore,

$$f_{ij}(a_i)(a_j) \in K_j = K_{(j)}$$

Therefore, we determined the operation of intersection of diagrams of subrepresentations properly.                                                    □

We denote the corresponding closure operator by $J[f]$. If

$$X = (X_{(1)} \subset A_{(1)}, ..., X_{(n)} \subset A_{(n)}) = (X_1 \subset A_1, ..., X_n \subset A_n)$$

is the tuple of sets, then

$$J[f, X] = (J_{(1)}[f, X], ..., J_{(n)}[f, X]) = (J_1[f, X], ..., J_n[f, X])$$

is the intersection of all tuples

$$K = (K_{(1)} \subset A_{(1)}, ..., K_{(n)} \subset A_{(n)}) = (K_1 \subset A_1, ..., K_n \subset A_n)$$

stable with respect to diagram of representations $(f, A)$ and such that for $(i) = (1), ..., (n)$, $K_{(i)}$ is subalgebra of $\Omega_{(i)}$-algebra $A_{(i)}$ containing $X_{(i)}$.

---

[8.1]This theorem is similar to definition of the lattice of subalgebras ([14], p. 79, 80)

THEOREM 8.1.4. *Let*[8.2] $(f, A)$ *be the diagram of representations. Let*

$$X = (X_{(1)} \subset A_{(1)}, ..., X_{(n)} \subset A_{(n)}) = (X_1 \subset A_1, ..., X_n \subset A_n)$$

*For every value of* $(i)$, $(i) = (1)$, ..., $(n)$, *we define a subsets* $X_{(i)k} \subset A_{(i)}$ *by induction on* $k$.

8.1.4.1: $X_{(i)0} = X_{(i)}$

8.1.4.2: $x \in X_{(i)k} => x \in X_{(i)k+1}$

8.1.4.3: $x_1 \in X_{(i)k}$, ..., $x_p \in X_{(i)k}$, $\omega \in \Omega_{(i)}(p) => x_1...x_p\omega \in X_{(i)k+1}$

8.1.4.4: $x_i \in X_{ik} = X_{(i)k}$, $x_j \in X_{jk} = X_{(j)k} => f_{ji}(x_j)(x_i) \in X_{(i)k+1}$

*For each value of* $(i)$, *we assume*

$$Y_{(i)} = \bigcup_{m=0}^{\infty} X_{(i)m}$$

*Then*

$$J_{(i)}[f, X] = Y_{(i)} \qquad (i) = (1), ..., (n)$$

PROOF. For each value of $(i)$ the proof of the theorem coincides with the proof of theorem 6.1.4. □

$J[f, X]$ is called **subrepresentation** of diagram of representations $(f, A)$ generated by tuple of sets $X$ and $X$ is a **generating set** of diagram of representations J[f,X]. In particular, a generating set of diagram of representations $(f, A)$ is a tuple

$$X = (X_{(1)} \subset A_{(1)}, ..., X_{(n)} \subset A_{(n)}) = (X_1 \subset A_1, ..., X_n \subset A_n)$$

such that $J[f, X] = A$.

From theorem 8.1.4, it follows next definition.

DEFINITION 8.1.5. *Let*

$$X = (X_{(1)} \subset A_{(1)}, ..., X_{(n)} \subset A_{(n)}) = (X_1 \subset A_1, ..., X_n \subset A_n)$$

*be tuple of sets. For each tuple of* $A$*-numbers* $a \in J[f, X]$

$$a = (a_{(1)} \quad ... \quad a_{(n)}) = (a_1 \quad ... \quad a_n)$$

*there exists* **tuple of** $\Omega$**-words**

$$w[f, X, a] = (w_{(1)}[f, X, a_{(1)}], ..., w_{(n)}[f, X, a_{(n)}])$$
$$= (w_1[f, X, a_1], ..., w_n[f, X, a_n])$$

*defined according to following rule.*

8.1.5.1: *If* $a_{(i)} \in X_{(i)}$, $(i) = (1)$, ..., $(n)$, *then* $a_{(i)}$ *is* $\Omega_{(i)}$*-word*

$$w_{(i)}[f, X, a_{(i)}] = a_{(i)}$$

8.1.5.2: *If* $a_{(i)1}$, ..., $a_{(i)p}$ *are* $\Omega_{(i)}$*-words,* $(i) = (1)$, ..., $(n)$, *and* $\omega \in \Omega_{(i)}(p)$, *then* $a_{(i)1}...a_{(i)p}\omega$ *is* $\Omega_{(i)}$*-word.*

8.1.5.3: *Let* $a_i = a_{(i)}$ *be* $\Omega_{(i)}$*-word,* $a_j = a_{(j)}$ *be* $\Omega_{(j)}$*-word. Let there exist the representation* $f_{ij}$. *Then* $f_{ij}(a_i)(a_j)$ *is* $\Omega_{(j)}$*-word.*

*Denote* $w[f, X]$ **the set of tuples of** $\Omega$**-words** *of diagram of representations* J[f,X]. □

---

We consider tuple of $A$-numbers in the definition 8.1.5 because we need an algorithm of generation of tuple of $\Omega$-words. However, to solve specific problem we need only some subset of tuples of $A$-numbers. For instance, in affine space we can consider or a set of points, either a set of vectors.

Choice of $\Omega_{(i)}$-word relative generating set $X$ is ambiguous. Therefore, if $\Omega_{(i)}$-number has different $\Omega_{(i)}$-words, then we will use indexes to distinguish them: $w_{(i)}[f, X, m]$, $w_{(i)1}[f, X, m]$, $w_{(i)2}[f, X, m]$.

**DEFINITION 8.1.6.** *Generating set $X$ of diagram of representations $(f, A)$ generates tuple of equivalences*

$$\rho[f, X] = (\rho_{(1)}[f, X], ..., \rho_{(n)}[f, X])$$

$$\rho_{(i)}[f, X] = \{(w_{(i)}[f, X, m_{(i)}], w_{(i)1}[f, X, m_{(i)}]) : m_{(i)} \in A_{(i)}\}$$

*on the set of tuples of $\Omega$-words.*                                              □

According to the definition 8.1.6, two $\Omega_{(i)}$-words with respect to the generating set $X$ of diagram of representations $(f, A)$ are equivalent iff they correspond to the same $A_{(i)}$-number. When we write equality of two $\Omega_{(i)}$-words with respect to the generating set $X$ of diagram of representations $(f, A)$, we will keep in mind that this equality is true up to equivalence $\rho_{(i)}[f, X]$.

We will use notation

$$r(a) = (r_{(1)}(a_{(1)}), ..., r_{(n)}(a_{(n)}))$$

for image of tuple of elements $a = (a_{(1)}, ..., a_{(n)})$ under the morphism of diagram of representations.

**THEOREM 8.1.7.** *Let $X$ be the generating set of diagram of representations $(f, A)$. Let $Y$ be the generating set of diagram of representations $(g, B)$. Morphism $r$ of diagram of representations $(f, A)$ forms the map of $\Omega$-words*

$$w[f \to g, X, Y, r] : w[f, X] \to w[g, Y]$$

$$X_{(i)} \subset A_{(i)} \quad Y_{(i)} = r_{(i)}(X_{(i)}) \quad (i) = (1), ..., (n)$$

*such that for any $(i)$, $(i) = (1), ..., (n)$,*

8.1.7.1: *If $a_{(i)} \in X_{(i)}$, $a'_{(i)} = r_{(i)}(a_{(i)})$, then*

$$w_{(i)}[f \to g, X, Y, r](a_{(i)}) = a'_{(i)}$$

8.1.7.2: *If*

$$a_{(i)1}, ..., a_{(i)p} \in w_{(i)}[f, X]$$

$$a'_{(i)1} = w_{(i)}[f \to g, X, Y, r](a_{(i)1}) \quad ... \quad a'_{(i)p} = w_{(i)}[f \to g, X, Y, r](a_{(i)p})$$

*then for operation $\omega \in \Omega_{(i)}(p)$ holds*

$$w_{(i)}[f \to g, X, Y, r](a_{(i),1}...a_{(i),p}\omega) = a'_{(i),1}...a'_{(i),p}\omega$$

8.1.7.3: *If*

$$a_i = a_{(i)} \in w_{(i)}[f, X] \quad a'_{(i)} = w_{(i)}[f \to g, X, Y, r](a_{(i)})$$

$$a_j = a_{(j)} \in w_{(j)}[f, X] \quad a'_j = a'_{(j)} = w_{(j)}[f, X, r](a_{(j)})$$

*then*

$$w_{(i)}[f \to g, X, Y, r](f_{ji}(a_j)(a_i)) = g_{ji}(a'_j)(a'_i)$$

PROOF. Statements 8.1.7.1, 8.1.7.2 are true by definition of the morphism $r$. The statement 8.1.7.3 follows from the equality (7.2.3). □

REMARK 8.1.8. *Let $r$ be morphism of diagram of representations $(f, A)$ into diagram of representations $(g, B)$. Let*

$$a \in J[f, X] \quad a' = r(a) \quad Y = r(X)$$

*The theorem 8.1.7 states that $a' \in J[g, Y]$. The theorem 8.1.7 also states that the tuple of $\Omega$-words representing a relative $X$ and the tuple of $\Omega$-words representing $a'$ relative $Y$ are generated according to the same algorithm. This allows considering of the tuple of $\Omega$-words $w[g, Y, a']$ as tuple of maps*

$$W[f, X, a] = (W_{(1)}[f, X, a], ..., W_{(n)}[f, X, a]) = (W_1[f, X, a], ..., W_n[f, X, a])$$

(8.1.2) $\qquad W_{(k)}[f, X, a] : (g, X') \to (g, X') \circ W_{(k)}[f, X, a] = w_{(k)}[g, X', a']$

*If $f = g$, then, instead of the map (8.1.2), we consider the map*

$$W_{(k)}[f, X, a] : Y \to Y \circ W_{(k)}[f, X, a] = w_{(k)}[f, Y, a']$$

$$W_{(k)}[f, X, a](Y) = Y \circ W_{(k)}[f, X, a]$$

*such that, if for certain morphism $r$*

$$Y = r(X) \quad a' = r(a)$$

*then*

$$W_{(k)}[f, X, a](Y) = Y \circ W[f, X, a] = w[f, Y, a'] = a'$$

*The map $W_{(k)}[f, X, a]$ is called* **coordinates of $A_{(k)}$-number** $a_{(k)}$ *relative to the tuples of sets $X$. Similarly, we consider coordinates of a set $B \subset J_{(k)}[f, X]$ relative to the set $X$*

$$W_{(k)}[f, X, B] = \{W_{(k)}[f, X, a] : a \in B\} = (W_{(k)}[f, X, a], a \in B)$$

*Denote*

$$W[f, X] = (W_{(1)}[f, X], ..., W_{(n)}[f, X]) = (W_1[f, X], ..., W_n[f, X])$$

$$W_{(k)}[f, X] = \{W_{(k)}[f, X, a] : a \in J_{(k)}[f, X]\} = (W_{(k)}[f, X, a], a \in J_{(k)}[f, X])$$

**the set of coordinates of representation** $J[f, X]$. □

THEOREM 8.1.9. *There is a structure of $\Omega_{(k)}$-algebra on the set of coordinates $W_{(k)}[f, X]$.*

PROOF. Let $\omega \in \Omega_{(k)}(n)$. Then for any $m_1, ..., m_n \in J_{(k)}[f, X]$ , we assume

(8.1.3) $\qquad W_{(k)}(f, X, m_1)...W_{(k)}(f, X, m_n)\omega = W_{(k)}(f, X, m_1...m_n\omega)$

According to the remark 8.1.8,

(8.1.4) $\qquad X \circ (W_{(k)}[f, X, m_1]...W_{(k)}[f, X, m_n]\omega) = X \circ W_{(k)}[f, X, m_1...m_n\omega]$
$$= w_{(k)}[f, X, m_1...m_n\omega]$$

follows from the equality (8.1.3). According to rule 8.1.5.2, from the equality (8.1.4), it follows that

$$
\begin{aligned}
(8.1.5) \quad & X \circ (W_{(k)}[f, X, m_1]...W_{(k)}[f, X, m_n]\omega) \\
& = w_{(k)}[f, X, m_1]...w_{(k)}[f, X, m_n]\omega \\
& = (X \circ W_{(k)}[f, X, m_1])...(X \circ W_{(k)}[f, X, m_n])\omega
\end{aligned}
$$

From the equality (8.1.5), it follows that the operation $\omega$ defined by the equality (8.1.3) on the set of coordinates $W_{(k)}[f, X]$ is defined properly.  □

THEOREM 8.1.10. *If there exists the representation $f_{jk}$ of $\Omega_j$-algebra $A_j$ in $\Omega_k$-algebra $A_k$, then there exists the representation $F_{jk}$ of $\Omega_j$-algebra $W_j[f, X]$ in $\Omega_k$-algebra $W_k[f, X]$.*

PROOF. Let $a_j \in J_j[f, X]$. Then for any $a_k \in J_k[f, X]$, we assume

$$(8.1.6) \qquad F_{jk}(W_j[f, X, a_j])(W_k[f, X, a_k]) = W_k[f, X, f_{jk}(a_j)(a_k)]$$

According to the remark 8.1.8,

$$
\begin{aligned}
(8.1.7) \quad & X \circ (F_{jk}(W_j[f, X, a_j])(W_k[f, X, a_k])) = X \circ W_k[f, X, f_{jk}(a_j)(a_k)] \\
& = w_k[f, X, f_{jk}(a_j)(a_k)]
\end{aligned}
$$

follows from the equality (8.1.6). According to rule 8.1.5.3, from the equality (8.1.7), it follows that

$$
\begin{aligned}
(8.1.8) \quad & X \circ (F_{jk}(W_j[f, X, a_j])(W_k[f, X, a_k])) \\
& = f_{jk}(w_j[f, X, a_j])(w_k[f, X, a_k]) \\
& = f_{jk}(X \circ W_j(f, X, a_j))(X \circ W_k(f, X, a_k))
\end{aligned}
$$

From the equality (8.1.8), it follows that the representation (8.1.6) of $\Omega_{k-1}$-algebra $W_{k-1}[f, X]$ in $\Omega_k$-algebra $W_k[f, X]$ is defined properly.  □

COROLLARY 8.1.11. *Tuple of $\Omega$-algebras*

$$W[f, X] = (W_{(1)}[f, X], ..., W_{(n)}[f, X])$$

*and the set of representations $F$ forms the diagram of representations $(F, W[f, X])$.*  □

THEOREM 8.1.12. *Let $(f, A)$, $(g, B)$ be diagrams of representations. For given sets $X_{(k)} \subset A_{(k)}$, $Y_{(k)} \subset B_{(k)}$, $(k) = (1)$, ..., $(n)$, consider tuple of maps*

$$R = (R_{(1)}, ..., R_{(n)})$$

*such that for any $(k)$, $(k) = (1)$, ..., $(n)$, the map*

$$R_{(k)} : X_{(k)} \to Y_{(k)}$$

*agree with the structure of diagram of representations, i. e.*

$$(8.1.9) \qquad
\begin{cases}
\omega \in \Omega_{(k)}(p), \ x_{(k)1}, \ ..., \ x_{(k)p}, \ x_{(k)1}...x_{(k)p}\omega \in X_{(k)}, \\
R_{(k)}(x_{(k)1}...x_{(k)p}\omega) \in Y_{(k)} \\
=> R_{(k)}(x_{(k)1}...x_{(k)p}\omega) = R_{(k)}(x_{(k)1})...R_{(k)}(x_{(k)p})\omega
\end{cases}
$$

$$(8.1.10) \qquad
\begin{cases}
a_j \in X_j, \ a_k \in X_k, \ R_k(f_{jk}(a_j)(a_k)) \in Y_k \\
=> R_k(f_{jk}(a_j)(a_k)) = g_{jk}(R_j(a_j))(R_k(a_k))
\end{cases}
$$

*Consider the tuple of maps of $\Omega$-words*

$$w_{(k)}[f \to g, \overline{\overline{e}}, Y, R] : w_{(k)}[f, \overline{\overline{e}}] \to w_{(k)}[g, Y]$$

*that satisfies conditions 8.1.7.1, 8.1.7.2, 8.1.7.3 and such that*

$$e_{(k)i} \in \overline{\overline{e}}_{(k)} => w_{(k)}[f \to g, \overline{\overline{e}}, Y, R](e_{(k)i}) = R_{(k)}(e_{(k)i})$$

*For each $(k)$, $(k)$, $(k) = (1), ..., (n)$, there exists homomorphism of $\Omega_{(k)}$-algebra*

$$r_{(k)} : A_{(k)} \to B_{(k)}$$

*defined by rule*

(8.1.11)          $$r_{(k)}(a_{(k)}) = w_{(k)}[f \to g, X, Y, R](w_{(k)}[f, X, a_{(k)}])$$

*Tuple of homomorphisms*

$$r = (r_{(1)} \quad ... \quad r_{(n)}) = (r_1 \quad ... \quad r_n)$$

*is morphism of diagrams of representations $J[f, X]$ and $J[g, Y]$.*

PROOF. For any $(k)$, $(k) = (1), ..., (n)$, consider the map

$$r_{(k)} : A_{(k)} \to B_{(k)}$$

LEMMA 8.1.13. *For any $(k)$, $(k) = (1), ..., (n)$, maps $r_{(k)}$ and $R_{(k)}$ coinside on the set $X_{(k)}$, and the map $r_{(k)}$ agrees with structure of $\Omega_{(k)}$-algebra.*

PROOF. If

(8.1.12)          $$w_{(k)}[f, X, a_{(k)}] = a_{(k)}$$

then $a_{(k)} \in X_{(k)}$. According to condition 8.1.7.1, the equality

(8.1.13)
$$r_{(k)}(a_{(k)}) = w_{(k)}[f \to g, X, Y, R](w_{(k)}[f, X, a_{(k)}]) = w_{(k)}[f \to g, X, Y, R](a_{(k)})$$
$$= R_{(k)}(a_{(k)})$$

follows from equalities (8.1.11), (8.1.12). The lemma follows from the equality (8.1.13). ⊙

LEMMA 8.1.14. *Let $\omega \in \Omega_{(k)}(p)$.*

(8.1.14)          $$r_{(k)}(x_{(k)1}...x_{(k)p}\omega) = r_{(k)}(x_{(k)1})...r_{(k)}(x_{(k)p})\omega$$

PROOF. We prove the lemma by induction over complexity of $\Omega_{(k)}$-word. If

$$x_{(k)1}, ..., x_{(k)p}, x_{(k)1}...x_{(k)p}\omega \in X_{(k)}$$

then the equality (8.1.14) follows from the statement (8.1.9).

Let the statement of induction be true for

$$a_{(k)1}, ..., a_{(k)p} \in J_{(k)}[f, X]$$

Let

(8.1.15)          $$w_{(k)1} = w_{(k)}[f, X, a_{(k)1}] \quad ... \quad w_{(k)p} = w_{(k)}[f, X, a_{(k)p}]$$

According to the statement of induction, the equality

(8.1.16)
$$r_{(k)}(a_{(k)1}) = w_{(k)}[f \to g, X, Y, R](w_{(k)1})$$
$$... = ...$$
$$r_{(k)}(a_{(k)p}) = w_{(k)}[f \to g, X, Y, R](w_{(k)p})$$

follows from equalities (8.1.11), (8.1.15). If

$$(8.1.17) \qquad a_{(k)} = a_{(k)1}...a_{(k)p}\omega$$

then according to condition 8.1.5.2,

$$w_{(k)}[f, X, a_{(k)}] = w_{(k)1}...w_{(k)p}\omega$$

According to condition 8.1.7.2, the equality

$$
\begin{aligned}
(8.1.18) \quad r_{(k)}(a_{(k)}) &= w_{(k)}[f \to g, X, Y, R](w_{(k)}[f, X, a_{(k)}]) \\
&= w_{(k)}[f \to g, X, Y, R](w_{(k)1}...w_{(k)p}\omega) \\
&= w_{(k)}[f \to g, X, Y, R](w_{(k)1})...w_{(k)}[f \to g, X, Y, R](w_{(k)p})\omega \\
&= (r_{(k)}(a_{(k)1}))...(r_{(k)}(a_{(k)p}))\omega
\end{aligned}
$$

follows from equalities (8.1.11), (8.1.17), (8.1.16). The equality (8.1.14) follows from the equality (8.1.13).                              ⊙

According to the lemma 8.1.13, maps $r_{(k)}$ and $R_{(k)}$ coinside on the set $X_{(k)}$. According to the lemma 8.1.14, the map $r_{(k)}$ is homomorphism of $\Omega_{(k)}$-algebra $A_{(k)}$ into $\Omega_{(k)}$-algebra $B_{(k)}$. To prove the theorem, it suffices to show that existence of the representation

$$f_{ji} : A_j \longrightarrow \ast \longrightarrow A_i$$

implies that the tuple of maps $(r_j \ \ r_i)$ is morphism of representations from $f_{ji}$ into $g_{ji}$ (the definition 7.2.1).

We prove the theorem by induction over complexity of $\Omega_i$-word.

If $a_i \in X_i$, $a_j \in X_j$, then the statement of induction follows from the statement (8.1.10)

Let the statement of induction be true for

$$a_j \in J_j[f, X] \quad w_j[f, X, a_j] = m_j$$

$$a_i \in J_i[f, X] \quad w_i[f, X, a_i] = m_i$$

According to condition 8.1.5.3,

$$(8.1.19) \qquad w_i(f, X, f_{ji}(a_j)(a_i)) = f_{ji}(m_j)(m_i)$$

According to condition 8.1.7.3, the equality

$$
\begin{aligned}
(8.1.20) \quad r_i(f_{ji}(a_j)(a_i)) &= w_i[f \to g, X, Y, R](w_i[f, X, f_{ji}(a_j)(a_i)]) \\
&= w_i[f \to g, X, Y, R](f_{ji}(m_j)(m_i)) \\
&= g_{ji}(w_j[g, Y, r_j(a_j)])(w_i[g, Y, r_i(a_i)]) \\
&= g_{ji}(r_j(a_j))(r_i(a_i))
\end{aligned}
$$

follows from equalities (8.1.11), (8.1.19), From equalities (7.2.3), (8.1.20), it follows that the map $r$ is morphism of the diagram of representations $(f, A)$.         □

REMARK 8.1.15. *The theorem 8.1.12 is the theorem of extension of map. The only statement we know about the tuple of sets $X$ is the statement that $X$ is the tuple of generating sets of the diagram of representations $(f, A)$. However, between the elements of the set $X_{(k)}$, $(k) = (1), ..., (n)$, there may be relationships generated by either operations of $\Omega_{(k)}$-algebra $A_{(k)}$, or by transformation of representation $f_{jk}$. Therefore, any map of tuple of sets $X$, in general, cannot be extended to an* *endomorphism of diagram of representations* $(f, A)$.[8.3] *However, if for any* $(k)$, $(k) = (1), ..., (n)$, *the map* $R_{(k)}$ *is coordinated with the structure of diagram of representations, then we can construct an extension of this map and this extension is morphism of diagram of representations* $(f, A)$.          □

DEFINITION 8.1.16. *Let* $X$ *be the tuple of generating sets of diagram of representations* $(f, A)$. *Let* $Y$ *be the tuple of generating sets of diagram of representations* $(g, B)$. *Let* $r$ *be the morphism of the diagram of representations* $(f, A)$ *into the diagram of representations* $(g, B)$. *The set of coordinates* $W[g, Y, r(X)]$ *is called* **coordinates of morphism of diagram of representations**.          □

DEFINITION 8.1.17. *Let* $X$ *be the tuple of generating sets of diagram of representations* $(f, A)$. *Let* $Y$ *be the tuple of generating sets of diagram of representations* $(g, B)$. *Let* $r$ *be the morphism of the diagram of representations* $(f, A)$ *into the diagram of representations* $(g, B)$. *Let for* $(k) = (1), ..., (n)$, $a_{(k)} \in A_{(k)}$. *We define* **superposition of coordinates** *of morphism* $r$ *of the diagram of representations and* $A_{(k)}$-*number* $a_{(k)}$ *as coordinates defined according to rule*

(8.1.21)          $W_{(k)}[g, Y, r_{(k)}(X_{(k)})] \circ W_{(k)}[f, X, a_{(k)}] = W_{(k)}[g, Y, r_{(k)}(a_{(k)})]$

*Let* $Y_{(k)} \subset A_{(k)}$. *We define superposition of coordinates of morphism* $r$ *of the diagram of representations and the set* $Y_{(k)}$ *according to rule*

$$(8.1.22) \qquad W_{(k)}[g, Y, r_{(k)}(X_{(k)})] \circ W_{(k)}[f, X, Y_{(k)}]$$
$$= (W_{(k)}[g, Y, r_{(k)}(X_{(k)})] \circ W_{(k)}[f, X, a_{(k)}], a_{(k)} \in Y_{(k)})$$

□

THEOREM 8.1.18. *Morphism* $r$ *of diagram of representations* $(f, A)$ *into diagram of representations* $(g, B)$ *generates the map of coordinates of diagram of representations*

(8.1.23)          $W_{(k)}[f \rightarrow g, X, Y, r] : W_{(k)}[f, X] \rightarrow W_{(k)}[g, Y]$

$(k) = (1), ..., (n)$, *such that*

$$(8.1.24) \qquad W_{(k)}[f, X, a] \rightarrow W_{(k)}[f \rightarrow g, X, Y, r] \circ W_{(k)}[f, X, a_{(k)}]$$
$$= W_{(k)}[g, Y, r_{(k)}(a_{(k)})]$$

$$(8.1.25) \qquad W_{(k)}[f \rightarrow g, X, Y, r] \circ W_{(k)}[f, X, a_{(k)}]$$
$$= W_{(k)}[g, Y, r_{(k)}(X_{(k)})] \circ W_{(k)}[f, X, a_{(k)}]$$

PROOF. According to the remark 8.1.8, we consider equalities (8.1.21), (8.1.23) relative to given tuple of generating sets $X$. The tuple of words

(8.1.26)          $X \circ W_{(k)}[f, X, a_{(k)}] = w_{(k)}[f, X, a_{(k)}]$

corresponds to coordinates $W_{(k)}[f, X, a_{(k)}]$; the tuple of words

(8.1.27)          $Y \circ W_{(k)}[g, Y, r_{(k)}(a_{(k)})] = w_{(k)}[g, Y, r_{(k)}(a_{(k)})]$

corresponds to coordinates $W_{(k)}[g, Y, r_{(k)}(a_{(k)})]$. Therefore, in order to prove the theorem, it is sufficient to show that the map $W_{(k)}[f \rightarrow g, X, Y, r]$ corresponds to map $w_{(k)}[f \rightarrow g, X, Y, r]$.

---

[8.3] In the theorem 8.2.9, requirements to tuple of generating sets are more stringent. Therefore, the theorem 8.2.9 says about extension of arbitrary map. A more detailed analysis is given in the remark 8.2.11.

We prove the theorem by induction over complexity of $\Omega_{(k)}$-word.

If $a_{(k)} \in X_{(k)}$, $a'_{(k)} = r_{(k)}(a_{(k)})$, then, according to equalities (8.1.26), (8.1.27), maps $W_{(k)}[f \to g, X, Y, r]$ and $w_{(k)}[f \to g, X, Y, r]$ are coordinated.

Let for $a_{(k)1}, \ldots, a_{(k)p} \in X_{(k)}$ maps $W_{(k)}[f \to g, X, Y, r]$ and $w_{(k)}[f \to g, X, Y, r]$ are coordinated. Let $\omega \in \Omega_{(k)}(p)$. According to the theorem 6.1.12

$$(8.1.28) \qquad W_{(k)}[f, X, a_{(k)1}\ldots a_{(k)p}\omega] = W_{(k)}[f, X, a_{(k)1}]\ldots W_{(k)}[f, X, a_{(k)p}]\omega$$

Because the map

$$r_{(k)} : A_{(k)} \to B_{(k)}$$

is homomorphism of $\Omega_{(k)}$-algebra, then from the equality (8.1.28) it follows that

$$
\begin{aligned}
(8.1.29) \qquad & W_{(k)}[g, Y, r_{(k)}(a_{(k)1}\ldots a_{(k)p}\omega)] \\
& = W_{(k)}[g, Y, (r(a_{(k)1})\ldots(r_{(k)}(a_{(k)p}))\omega] \\
& = W_{(k)}[g, Y, r_{(k)}(a_{(k)1})]\ldots W_{(k)}[g, Y, r_{(k)}(a_{(k)p})]\omega
\end{aligned}
$$

From equalities (8.1.28), (8.1.29) and the statement of induction, it follows that the maps $W_{(k)}[f \to g, X, Y, r]$ and $w_{(k)}[f \to g, X, Y, r]$ are coordinated for $a_{(k)} = a_{(k)1}\ldots a_{(k)p}\omega$.

Let for $a_{j1} \in A_j$ maps $W_j[f \to g, X, Y, r]$ and $w_j[f \to g, X, Y, r]$ are coordinated. Let for $a_{i1} \in A_i$ maps $W_i[f \to g, X, Y, r]$ and $w_i[f \to g, X, Y, r]$ are coordinated. According to the theorem 8.1.10

$$(8.1.30) \qquad W_i(f, X, f_{ji}(a_j)(a_i)) = F_{ji}(W_j(f, X, a_j))(W_i(f, X, a_i))$$

Because the map $(r_j, r_i)$ is morphism of representations $f_{ji}$ into representation $F_{ji}$, then from the equality (8.1.30) it follows that

$$
\begin{aligned}
(8.1.31) \qquad W_i[g, Y, r_i(f_{ji}(a_j)(a_i))] & = W_i[g, Y, g_{ji}(r_j(a_j))(r_i(a_i))] \\
& = G_{ji}(W_j[g, Y, r_j(a_j)])(W_i[g, Y, r_i(a_{n,1})])
\end{aligned}
$$

From equalities (8.1.30), (8.1.31) and the statement of induction, it follows that maps $W_i[f \to g, X, Y, r]$ and $w_i[f \to g, X, Y, r]$ are coordinated for $b_i = f_{ji}(a_j)(a_i)$. $\qquad\square$

COROLLARY 8.1.19. *Let $X$ be the tuple of generating sets of the diagram of representations $(f, A)$. Let $Y$ be the tuple of generating sets of the diagram of representations $(g, B)$. Let $r$ be the morphism of the diagram of representations $(f, A)$ into the diagram of representations $(g, B)$. The map*

$$W[f \to g, X, Y, r] = (W_{(1)}[f \to g, X, Y, r], \ldots, W_{(n)}[f \to g, X, Y, r])$$

*is morphism of diagram of representations $(F, W[f, X])$ into diagram of representations $(G, W[g, Y])$.* $\qquad\square$

Hereinafter we will identify the map $W[f \to g, X, Y, r]$ and the set of coordinates $W[g, Y, r(X)]$.

THEOREM 8.1.20. *Let $X$ be the tuple of generating sets of the diagram of representations $(f, A)$. Let $Y$ be the tuple of generating sets of the diagram of representations $(g, B)$. Let $r$ be the morphism of the diagram of representations $(f, A)$ into the diagram of representations $(g, B)$. Let $Y \subset A$. Then*

$$(8.1.32) \qquad W[g, Y, r(X)] \circ W[f, X, X'] = W[g, Y, r(X')]$$

$$(8.1.33) \qquad W[f \to g, X, Y, r] \circ W[f, X, X'] = W[g, Y, r(X')]$$

Proof. The equality (8.1.32) follows from the equality

$$r(X') = (r(a), a \in X')$$

as well from equalities (8.1.21), (8.1.22). The equality (8.1.33) is corollary of equalities (8.1.32), (8.1.24). □

Theorem 8.1.21. *Let $X$ be the tuple of generating sets of the diagram of representations $(f, A)$. Let $Y$ be the tuple of generating sets of the diagram of representations $(g, B)$. Let $Z$ be the tuple of generating sets of the diagram of representations $(h, C)$. Let $r$ be the morphism of the diagram of representations $(f, A)$ into the diagram of representations $(g, B)$. Let $s$ be the morphism of the diagram of representations $(g, B)$ into the diagram of representations $(h, C)$. Then*

$$(8.1.34) \qquad W[h, Z, s(Y)] \circ W[g, Y, r(X)] = W[h, Z, (s \circ r)(X)]$$

$$(8.1.35) \qquad W[g \to h, Y, Z, s] \circ W[f \to g, X, Y, r] = W[f \to h, X, Z, s \circ r]$$

Proof. The equality

$$(8.1.36) \qquad W[h, Z, s(Y')] \circ W[g, Y, Y'] = W[h, z, s(Y')]$$

follows from the equality (8.1.32). The equality (8.1.34) follows from the equality (8.1.36), if we assume $Y' = r(X)$. The equality (8.1.35) follows from the equality (8.1.34). □

Definition 8.1.22. *We can generalize the definition of the superposition of coordinates and assume that one of the factors is a tuple of sets of $\Omega$-words. Accordingly, the definition of the superposition of coordinates has the form*

$$W[g, Y, r(X)] \circ w[f, X, X'] = w[g, Y, r(X)] \circ W[f, X, X'] = w[g, Y, r(X')]$$

□

The following forms of writing an image of a tuple of sets $X'$ under morphism $r$ of diagram of representations are equivalent

$$(8.1.37) \qquad \begin{aligned} r(X') &= r(X) \circ W[f, X, X'] \\ &= (Y \circ W[g, Y, r(X)]) \circ W[f, X, X'] \end{aligned}$$

From equalities (8.1.32), (8.1.37), it follows that

$$(8.1.38) \qquad \begin{aligned} &Y \circ (W[g, Y, r(X) \circ W[f, X, X']) \\ =&(Y \circ W[g, Y, r(X)) \circ W[f, X, X'] \end{aligned}$$

The equality (8.1.38) is associative law for composition and allows us to write expression

$$Y \circ W[g, Y, r(X) \circ W[f, X, X']$$

without brackets.

Definition 8.1.23. *Let*

$$X = (X_{(1)} \subset A_{(1)}, ..., X_{(n)} \subset A_{(n)}) = (X_1 \subset A_1, ..., X_n \subset A_n)$$

*be generating set of diagram of representations $(f, A)$. Let map $r$ be endomorphism of diagram of representations $(f, A)$. Let the tuple of sets $Y = r(X)$ be image of tuple of sets $X$ under the map $r$. Endomorphism $r$ of diagram of representations $(f, A)$ is called regular on the tuple of generating set $X$, if the tuple of sets $Y$ is tuple*

*of generating sets of diagram of representations* $(f, A)$. *Otherwise, endomorphism* $r$ *is called singular on the tuple of generating sets* $X$.                    □

DEFINITION 8.1.24. *Endomorphism* $r$ *of diagram of representations* $(f, A)$ *is called* **regular**, *if it is regular on any tuple of generating sets. Otherwise, endomorphism* $r$ *is called* **singular**.                    □

THEOREM 8.1.25. *Automorphism* $r$ *of diagram of representations* $(f, A)$ *is regular endomorphism.*

PROOF. Let $X$ be tuple of generating sets of diagram of representations $(f, A)$. Let $Y = r(X)$. According to theorem 8.1.18 endomorphism $r$ forms the map of $\Omega$-words $w[f \to g, X, Y, r]$. Let $a' \in A$. Since $r$ is automorphism, then there exists $a \in A$, $r(a) = a'$. According to definition 8.1.5, $w[f, X, a]$ is tuple of $\Omega$-words representing $a$ relative to tuple of generating sets $X$. According to theorem 8.1.18, $w[f, X', a']$ is tuple of $\Omega$-words representing $a'$ relative to tuple of sets $Y$

$$w[f, Y, a'] = w[f \to g, X, Y, r](w[f, X, a])$$

Therefore, $Y$ is generating set of diagram of representations $(f, A)$. According to definition 8.1.24, automorphism $r$ is regular.                    □

## 8.2. Basis of Diagram of Representations

DEFINITION 8.2.1. *Let* $(f, A)$ *be diagram of representations and*

$$Gen[f, A] = \{X = (X_{(1)}, ..., X_{(n)}) : X_{(k)} \subseteq A_{(k)}, J_{(k)}[f, X] = A_{(k)}\}$$

*If, for the tuple of sets* $X \subset A_2$, *it is true that* $X \in Gen[f, A]$, *then for any tuple of sets* $Y$, $X_{(k)} \subset Y_{(k)} \subset A_{(k)}$, $(k) = (1), ..., (n)$ *also it is true that* $Y \in Gen[f, A]$. *If there exists minimal tuple of sets* $X \in Gen[f, A]$, *then the tuple of sets* $X$ *is called* **quasibasis of diagram of representations** $(f, A)$.                    □

THEOREM 8.2.2. *If the tuple of sets* $X$ *is the quasibasis of diagram of representations* $(f, A)$, *then for any* $(k)$, $(k) = (1), ..., (n)$, *and for any* $m \in X_{(k)}$ *the tuple of sets*

$$X' = (X_{(1)}, ..., X'_{(k)} = X_{(k)} \setminus \{m\}, ..., X_{(n)})$$

*is not generating set of diagram of representations* $(f, A)$.

PROOF. Let $X$ be quasibasis of diagram of representations $(f, A)$. Assume that for some $m \in X_{(k)}$ there exist $\Omega_{(k)}$-word

$$w = w_{(k)}[f, X', m]$$

Consider $A_{(k)}$-number $m'$ such that it has $\Omega_{(k)}$-word $w' = w_{(k)}[f, X, m']$ that depends on $m$. According to the definition 8.1.5, any occurrence of $A_{(k)}$-number $m$ into $\Omega_{(k)}$-word $w'$ can be substituted by the $\Omega_{(k)}$-word $w$. Therefore, the $\Omega_{(k)}$-word $w'$ does not depend on $m$, and the tuple of sets $X'$ is generating set of diagram of representations $(f, A)$. Therefore, $X$ is not quasibasis of diagram of representations $(f, A)$.                    □

REMARK 8.2.3. *The proof of the theorem 8.2.2 gives us effective method for constructing the quasibasis of diagram of representations* $(f, A)$. *We define a quasibasis of diagram of representations by induction over diagram of representations. We start to build a quasibasis in* $\Omega$-algebras *from the set* $A_{[0]}$. *When the quasibasis is constructed in* $\Omega$-algebras *from the set* $A_{[i]}$, *we can proceed to the construction of quasibasis in* $\Omega$-algebras *from the set* $A_{[i+1]}$.                    □

For each $(k)$, $(k) = (1), \ldots, (n)$, we introduced $\Omega_{(k)}$-word of $A_{(k)}$-number $x$ relative generating set $X$ in the definition [8.4] 8.1.5. From the theorem 8.2.2, it follows that if the generating set $X$ is not an quasibasis, then a choice of $\Omega_{(k)}$-word relative generating set $X$ is ambiguous. However, even if the generating set $X$ is an quasibasis, then a representation of $m \in A_{(k)}$ in form of $\Omega_{(k)}$-word is ambiguous.

REMARK 8.2.4. *There are three reasons of ambiguity in notation of $\Omega_{(k)}$-word.*

8.2.4.1: *In $\Omega_{(k)}$-algebra $A_{(k)}$, $(k) = (1), \ldots, (n)$, equalities may be defined. For instance, if $e$ is unit of multiplicative group $A_{(k)}$, then the equality*

$$ae = a$$

*is true for any $a \in A_{(k)}$.*

8.2.4.2: *Ambiguity of choice of $\Omega_{(k)}$-word may be associated with properties of representation. For instance, let there exist representation $f_{ik}$ of $\Omega_i$-algebra $A_i$ in $\Omega_k$-algebra $A_k$. If $m_1, \ldots, m_n$ are $\Omega_k$-words, $\omega \in \Omega_k(n)$ and $a$ is $\Omega_i$-word, then* [8.5]

$$(8.2.1) \qquad f_{ik}(a)(m_1 \ldots m_n \omega) = (f_{ik}(a)(m_1)) \ldots (f_{ik}(a)(m_n)) \omega$$

*At the same time, if $\omega$ is operation of $\Omega_i$-algebra $A_i$ and operation of $\Omega_k$-algebra $A_k$, then we require that $\Omega_k$-words $f(a_1 \ldots a_n \omega)(x)$ and $(f(a_1)(x)) \ldots (f(a_n)(x)) \omega$ describe the same element of $\Omega_k$-algebra $A_k$.* [8.6]

$$(8.2.2) \qquad f(a_1 \ldots a_n \omega)(x) = (f(a_1)(x)) \ldots (f(a_n)(x)) \omega$$

8.2.4.3: *Equalities like (8.2.1), (8.2.2) persist under morphism of diagram of representations. Therefore we can ignore this form of ambiguity of $\Omega_{(k)}$-word. However, a fundamentally different form of ambiguity is possible. We can see an example of such ambiguity in theorems 9.3.15, 9.3.16.*

*So we see that we can define different equivalence relations on the set of $\Omega_{(k)}$-words.* [8.7] *Our goal is to find a maximum equivalence on the set of $\Omega_{(k)}$-words which persist under morphism of representation.*

---

[8.4] Naturally, arguments at the beginning of this section are the same as arguments at the beginning of the section 6.2 and I saved these arguments for the completeness of the text.

[8.5] For instance, let $\{e_1, e_2\}$ be the basis of vector space over field $k$. The equation (8.2.1) has the form of distributive law

$$a(b^1 e_1 + b^2 e_2) = (ab^1)e_1 + (ab^2)e_2$$

[8.6] For vector space, this requirement has the form of distributive law

$$(a + b)e_1 = ae_1 + be_1$$

[8.7] Evidently each of the equalities (8.2.1), (8.2.2) generates some equivalence relation.

*A similar remark concerns the map $W[f, X, m]$ defined in the remark 8.1.8.*[8.8]

$$\square$$

THEOREM 8.2.5. *Let $X$ be quasibasis of the diagram of representations $(f, A)$. Consider tuple of equivalences*

$$\lambda[f, X] = (\lambda_{(1)}[f, X], ..., \lambda_{(n)}[f, X])$$

$$\lambda_{(k)}[f, X] \subseteq w_{(k)}[f, X] \times w_{(k)}[f, X]$$

*which is generated exclusively by the following statements.*

8.2.5.1: *If in $\Omega_{(k)}$-algebra $A_{(k)}$ there is an equality*

$$w_{(k)1}[f, X, m] = w_{(k)2}[f, X, m]$$

*defining structure of $\Omega_{(k)}$-algebra, then*

$$(w_{(k)1}[f, X, m], w_{(k)2}[f, X, m]) \in \lambda_{(k)}[f, X]$$

8.2.5.2: *If there exists representation $f_{ik}$ and in $\Omega_i$-algebra $A_i$ there is an equality*

$$w_{i1}[f, X, m] = w_{i2}[f, X, m]$$

*defining structure of $\Omega_i$-algebra, then*

$$(f_{ik}(w_{i1})(w_k[f, X, m]), f_{ik}(w_{i2})(w_k[f, X, m])) \in \lambda_k[f, X]$$

8.2.5.3: *If there exists representation $f_{ik}$, then for any operation $\omega \in \Omega_i(n)$,*

$$(f_{ik}(a_{i1}...a_{in}\omega)(a_2), (f_{ik}(a_{i1})...f_{ik}(a_{in})\omega)(a_2)) \in \lambda_k[f, X]$$

8.2.5.4: *If there exists representation $f_{ik}$, then for any operation $\omega \in \Omega_k(n)$,*

$$(f_{ik}(a_i)(a_{k1}...a_{kn}\omega), f_{ik}(a_i)(a_{k1})...f_{ik}(a_i)(a_{kn})\omega) \in \lambda_k[f, X]$$

---

[8.8]If vector space has finite basis, then we represent the basis as matrix

$$\overline{\overline{e}} = \begin{pmatrix} e_1 & ... & e_2 \end{pmatrix}$$

We present the map $W[f, \overline{\overline{e}}](v)$ as matrix

$$W[f, \overline{\overline{e}}, v] = \begin{pmatrix} v^1 \\ ... \\ v^n \end{pmatrix}$$

Then

$$W[f, \overline{\overline{e}}, v](\overline{\overline{e}}') = W[f, \overline{\overline{e}}, v] \begin{pmatrix} e'_1 & ... & e'_n \end{pmatrix} = \begin{pmatrix} v^1 \\ ... \\ v^n \end{pmatrix} \begin{pmatrix} e'_1 & ... & e'_n \end{pmatrix}$$

has form of matrix product.

**8.2.5.5:** *If there exists representation $f_{ik}$, $\omega \in \Omega_i(n) \cap \Omega_k(n)$ and the representation $f_{ik}$ satisfies equality*[8.9]

$$f(a_{i1}...a_{in}\omega)(a_k) = (f(a_{i1})(a_k))...(f(a_{in})(a_k))\omega$$

*then we can assume that the following equality is true*

$$(f(a_{i1}...a_{in}\omega)(a_k), (f(a_{i1})(a_k))...(f(a_{in})(a_k))\omega) \in \lambda_k[f, X]$$

PROOF. The theorem is true because considered equalities are preserved under homomorphisms of universal algebras $A_{(k)}$. □

DEFINITION 8.2.6. *Quasibasis $\overline{\overline{e}}$ of the diagram of representations $(f, A)$ such that*

$$\rho[f, \overline{\overline{e}}] = \lambda[f, \overline{\overline{e}}]$$

*is called* **basis of diagram of representations** *$(f, A)$.* □

REMARK 8.2.7. *We write a basis also in following form*

$$\overline{\overline{e}} = (\overline{\overline{e}}_{(1)}, ..., \overline{\overline{e}}_{(n)})$$

$$\overline{\overline{e}}_{(k)} = (e_{(k)l}, e_{(k)l} \in \overline{\overline{e}}_{(k)}) \quad (k) = (1), ..., (n)$$

*If basis is finite, then we also use notation*

$$\overline{\overline{e}}_{(k)} = (e_{(k)i}, i \in I_{(k)}) = (e_{(k)1}, ..., e_{(k)p_{(k)}}) \quad (k) = (1), ..., (n)$$

□

THEOREM 8.2.8. *Automorphism of the diagram of representations $(f, A)$ maps a basis of the diagram of representations $(f, A)$ into basis.*

PROOF. Let the map $r$ be automorphism of the diagram of representations $(f, A)$. Let the tuple of sets $\overline{\overline{e}}$ be a basis of the diagram of representations $(f, A)$. Let [8.10] $\overline{\overline{e}}' = r \circ \overline{\overline{e}}$. Assume that the tuple of sets $\overline{\overline{e}}'$ is not basis. According to the theorem 8.2.2, there exist $(k)$, $(k) = (1), ..., (n)$, and $e'_{(k)i} \in \overline{\overline{e}}'_{(k)}$ such that the tuple of sets

$$Z = (\overline{\overline{e}}'_{(1)}, ..., Z_{(k)} = \overline{\overline{e}}'_{(k)} \setminus \{e'_{(k)i}\}, ..., \overline{\overline{e}}'_{(n)})$$

is generating set of the diagram of representations $(f, A)$. According to the theorem 7.3.3 the map $r^{-1}$ is automorphism of the diagram of representations $(f, A)$. According to the theorem 8.1.25 and definition 8.1.24, the tuple of sets

$$X = (\overline{\overline{e}}_{(1)}, ..., X_{(k)} = \overline{\overline{e}}_{(k)} \setminus \{r^{-1}_{(k)}(e'_{(k)i})\}, ..., \overline{\overline{e}}_{(n)})$$

is generating set of the diagram of representations $(f, A)$. The contradiction completes the proof of the theorem. □

---

[8.9] Consider a representation of commutative ring $D$ in $D$-algebra $A$. We will use notation

$$f(a)(v) = av$$

The operations of addition and multiplication are defined in both algebras. However the equality

$$f(a + b)(v) = f(a)(v) + f(b)(v)$$

is true, and the equality

$$f(ab)(v) = f(a)(v)f(b)(v)$$

is wrong.

[8.10] According to definitions 5.1.3, 8.3.1, we will use notation $r(\overline{\overline{e}}) = r \circ \overline{\overline{e}}$.

Theorem 8.2.9. *Let $\overline{\overline{e}}$ be the basis of the diagram of representations $(f, A)$. Let $(g, B)$ be diagram of representations. Let*

$$R : \overline{\overline{e}} \to Y$$

*be arbitrary map of the tuple of sets $\overline{\overline{e}}$, $Y_{(k)} \subseteq B_{(k)}$, $(k) = (1), ..., (n)$. Consider the tuple of maps*

$$w_{(k)}[f \to g, \overline{\overline{e}}, Y, R] : w_{(k)}[f, \overline{\overline{e}}] \to w_{(k)}[g, Y]$$

*that satisfy conditions 8.1.7.1, 8.1.7.2, 8.1.7.3 and such that*

$$e_{(k)i} \in \overline{\overline{e}}_{(k)} => w_{(k)}[f \to g, \overline{\overline{e}}, Y, R](e_{(k)i}) = R_{(k)}(e_{(k)i})$$

*There exists unique morphism of diagram of representations*[8.11]

$$r : A \to B$$

*defined by rule*

$$r(a) = w[f \to g, \overline{\overline{e}}, Y, R](w[f, \overline{\overline{e}}, a])$$

Proof. The statement of theorem is corollary theorems 6.1.10, 6.1.14. □

Corollary 8.2.10. *Let $\overline{\overline{e}}$, $\overline{\overline{e}}'$ be bases of the representation $(f, A)$. Let $r$ be the automorphism of the representation $(f, A)$ such that $\overline{\overline{e}}' = r \circ \overline{\overline{e}}$. Automorphism $r$ is uniquely defined.* □

Remark 8.2.11. *The theorem 8.2.9, as well as the theorem 8.1.12, is the theorem of extension of map. However in this theorem, $\overline{\overline{e}}$ is not arbitrary generating set of the diagram of representations, but basis. According to the remark 8.2.3, we cannot determine coordinates of any element of basis through the remaining elements of the same basis. Therefore, we do not need to coordinate the map of the basis with representation.* □

Theorem 8.2.12. *The set of coordinates $W[f, \overline{\overline{e}}, \overline{\overline{e}}]$ corresponds to identity transformation*

$$W[f, \overline{\overline{e}}, E] = W[f, \overline{\overline{e}}, \overline{\overline{e}}]$$

Proof. The statement of the theorem follows from the equality

$$a = \overline{\overline{e}} \circ W[f, \overline{\overline{e}}, a] = \overline{\overline{e}} \circ W[f, \overline{\overline{e}}, \overline{\overline{e}}] \circ W[f, \overline{\overline{e}}, a]$$

□

Theorem 8.2.13. *Let $W[f, \overline{\overline{e}}, r \circ \overline{\overline{e}}]$ be the set of coordinates of automorphism $r$. There exists set of coordinates $W[f, r \circ \overline{\overline{e}}, \overline{\overline{e}}]$, corresponding to automorphism $r^{-1}$. The set of coordinates $W[f, r \circ \overline{\overline{e}}, \overline{\overline{e}}]$ satisfy to equalities*[8.12]

$$(8.2.3) \qquad W[f, \overline{\overline{e}}, r \circ \overline{\overline{e}}] \circ W[f, r \circ \overline{\overline{e}}, \overline{\overline{e}}] = W[f, \overline{\overline{e}}, \overline{\overline{e}}]$$

$$W[f \to f, \overline{\overline{e}}, \overline{\overline{e}}, r^{-1}] = W[f \to f, \overline{\overline{e}}, \overline{\overline{e}}, r]^{-1} = W[f, r \circ \overline{\overline{e}}, \overline{\overline{e}}]$$

---

[8.11] This statement is similar to the theorem [2]-4.1, p. 135.
[8.12] See also remark 6.2.15.

Proof. Since $r$ is automorphism of the diagram of representations $(f, A)$, then, according to the theorem 8.2.8, the set $r \circ \overline{\overline{e}}$ is a basis of the diagram of representations $(f, A)$. Therefore, there exists the set of coordinates $W[f, r \circ \overline{\overline{e}}, \overline{\overline{e}}]$. The equality (8.2.3) follows from the chain of equalities

$$W[f, \overline{\overline{e}}, r \circ \overline{\overline{e}}] \circ W[f, r \circ \overline{\overline{e}}, \overline{\overline{e}}] = W[f, \overline{\overline{e}}, r \circ \overline{\overline{e}}] \circ W[f, \overline{\overline{e}}, r^{-1} \circ \overline{\overline{e}}]$$
$$= W[f, \overline{\overline{e}}, r \circ r^{-1} \circ \overline{\overline{e}}] = W[f, \overline{\overline{e}}, \overline{\overline{e}}]$$

□

Theorem 8.2.14. *Let* $W[f, \overline{\overline{e}}, r \circ \overline{\overline{e}}]$ *be the set of coordinates of automorphism* $r$. *Let* $W[f, \overline{\overline{e}}, s \circ \overline{\overline{e}}]$ *be the set of coordinates of automorphism* $s$. *The set of coordinates of automorphism* $(r \circ s)^{-1}$ *satisfies to the equality*

$$(8.2.4) \qquad W[f, (r \circ s) \circ \overline{\overline{e}}, \overline{\overline{e}}] = W[f, s \circ (r \circ \overline{\overline{e}}), \overline{\overline{e}}] = W[f, s \circ \overline{\overline{e}}, \overline{\overline{e}}] \circ W[f, r \circ \overline{\overline{e}}, \overline{\overline{e}}]$$

Proof. The equality

$$(8.2.5) \qquad \begin{aligned} W[f, (r \circ s) \circ \overline{\overline{e}}, \overline{\overline{e}}] &= W[f, \overline{\overline{e}}, (r \circ s)^{-1} \circ \overline{\overline{e}}] = W[f, \overline{\overline{e}}, s^{-1} \circ r^{-1} \circ \overline{\overline{e}}] \\ &= W[f, \overline{\overline{e}}, s^{-1} \circ e] \circ W[f, \overline{\overline{e}}, r^{-1} \circ \overline{\overline{e}}] \\ &= W[f, s \circ \overline{\overline{e}}, \overline{\overline{e}}] \circ W[f, r \circ \overline{\overline{e}}, \overline{\overline{e}}] \\ &= W[f, s \circ (r \circ \overline{\overline{e}}), \overline{\overline{e}}] \end{aligned}$$

follows from theorems 8.1.21, 8.2.13. The equality (8.2.4) follows from the equality (8.2.5).                □

Theorem 8.2.15. *The group of automorphisms* $GA(f)$ *of the diagram of effective representations* $(f, A)$ *generates effective left-side representation in the diagram of representations* $(f, A)$.

Proof. From the corollary 8.2.10, it follows that if automorphism $r$ maps a basis $\overline{\overline{e}}$ into a basis $\overline{\overline{e}}'$, then the set of coordinates $W[f, \overline{\overline{e}}, \overline{\overline{e}}']$ uniquely determines an automorphism $r$. From the theorem 8.1.18, it follows that the set of coordinates $W[f, \overline{\overline{e}}, \overline{\overline{e}}']$ determines the map of coordinates relative to the basis $\overline{\overline{e}}$ under automorphism of the diagram of representations $(f, A)$. From the equality (8.1.37), it follows that automorphism $r$ acts from the left on elements of $\Omega_{(k)}$-algebra $A_{(k)}$, $(k) = (1), ..., (n)$. From the equality (8.1.34), it follows that the representation of group is left-side representation. According to the theorem 8.2.12 the set of coordinates $W[f, \overline{\overline{e}}, \overline{\overline{e}}]$ corresponds to identity transformation. From the theorem 8.2.13, it follows that the set of coordinates $W[f, r \circ \overline{\overline{e}}, \overline{\overline{e}}]$ corresponds to transformation, inverse to transformation $W[f, \overline{\overline{e}}, r \circ \overline{\overline{e}}]$.                □

## 8.3. Basis Manifold of Diagram of Representations

The set $\mathcal{B}[f]$ of bases of diagram of representations $(f, A)$ is called **basis manifold** of diagram of representations $(f, A)$.

Definition 8.3.1. *According to the theorem 8.2.8 and to the definition 8.1.22, automorphism* $r$ *of the diagram of representations* $(f, A)$ *generates transformation*

$$(8.3.1) \qquad \begin{aligned} r &: \overline{\overline{h}} \to r \circ \overline{\overline{h}} \\ r \circ \overline{\overline{h}} &= W[f, \overline{\overline{e}}, r \circ \overline{\overline{e}}] \circ \overline{\overline{h}} \end{aligned}$$

*of the basis manifold of diagram of representations. This transformation is called* **active**. *According to the theorem 7.3.3, we defined left-side representation*

$$A(f) : GA(f) \longrightarrow \mathcal{B}[f]$$

*of group $GA(f)$ in basis manifold $\mathcal{B}[f]$. Representation $A(f)$ is called* **active representation**. *According to the corollary 8.2.10, this representation is single transitive.* □

REMARK 8.3.2. *According to remark 8.2.3, it is possible that there exist bases of diagram of representations $(f, A)$ such that there is no active transformation between them. Then we consider the orbit of selected basis as basis manifold. Therefore, it is possible that the diagram of representations $(f, A)$ has different basis manifolds. We will assume that we have chosen a basis manifold.*

THEOREM 8.3.3. *There exists single transitive right-side representation*

$$P(f) : GA(f) \longrightarrow \mathcal{B}[f]$$

*of group $GA(f)$ in basis manifold $\mathcal{B}[f]$. Representation $P(f)$ is called* **passive representation**.

PROOF. Since $A(f)$ is single transitive left-side representation of group $GA(f)$, then single transitive right-side representation $P(f)$ is uniquely defined according to the theorem 5.5.9. □

THEOREM 8.3.4. *Transformation of representation $P(f)$ is called* **passive transformation of the basis manifold** *of diagram of representations. We also use notation*

$$s(\overline{\overline{e}}) = \overline{\overline{e}} \circ s$$

*to denote the image of basis $\overline{\overline{e}}$ under passive transformation $s$. Passive transformation of basis has form*

(8.3.2)
$$s : \overline{\overline{h}} \to \overline{\overline{h}} \circ s$$
$$\overline{\overline{h}} \circ s = \overline{\overline{h}} \circ W[f, \overline{\overline{e}}, \overline{\overline{e}} \circ s]$$

PROOF. According to the equality (8.3.1), active transformation acts from left on coordinates of basis. The equality (8.3.2) follows from theorems 5.5.8, 5.5.9, 5.5.11; according to these theorems, passive transformation acts from right on coordinates of basis. □

THEOREM 8.3.5. *Passive transformation of the basis manifold is automorphism of representation $A(f)$.*

PROOF. The theorem follows from the theorem 5.5.11. □

THEOREM 8.3.6. *Let $s$ be passive transformation of the basis manifold of the diagram of representations $(f, A)$. Let $\overline{\overline{e}}_1$ be the basis of the diagram of representations $(f, A)$, $\overline{\overline{e}}_2 = \overline{\overline{e}}_1 \circ s$. For basis $\overline{\overline{e}}_3$, let there exists an active transformation $r$ such that $\overline{\overline{e}}_3 = r \circ \overline{\overline{e}}_1$. Let $\overline{\overline{e}}_4 = r \circ \overline{\overline{e}}_2$. Then $\overline{\overline{e}}_4 = \overline{\overline{e}}_3 \circ s$.*

PROOF. According to the equality (8.3.1), active transformation of coordinates of basis $\overline{\overline{e}}_3$ has form

(8.3.3)      $$\overline{\overline{e}}_4 = W[f, \overline{\overline{e}}_1, \overline{\overline{e}}_3] \circ \overline{\overline{e}}_2 = W[f, \overline{\overline{e}}_1, \overline{\overline{e}}_3] \circ \overline{\overline{e}}_1 \circ W[f, \overline{\overline{e}}_1, \overline{\overline{e}}_2]$$

Let $\overline{\overline{e}}_5 = \overline{\overline{e}}_3 \circ s$. From the equality (8.3.2), it follows that

$$(8.3.4) \qquad \overline{\overline{e}}_5 = \overline{\overline{e}}_3 \circ W[f, \overline{\overline{e}}_1, \overline{\overline{e}}_2] = W[f, \overline{\overline{e}}_1, \overline{\overline{e}}_3] \circ \overline{\overline{e}}_1 \circ W[f, \overline{\overline{e}}_1, \overline{\overline{e}}_2]$$

From match of expressions in equalities (8.3.3), (8.3.4), it follows that $\overline{\overline{e}}_4 = \overline{\overline{e}}_5$. Therefore, the diagram

$$
\begin{array}{ccc}
\overline{\overline{e}}_1 \in \mathcal{B}[f] & \xrightarrow{\ \ r\ \ } & \overline{\overline{e}}_3 \in \mathcal{B}[f] \\
\downarrow{\scriptstyle s} & & \downarrow{\scriptstyle s} \\
\overline{\overline{e}}_2 \in \mathcal{B}[f] & \xrightarrow{\ \ r\ \ } & \overline{\overline{e}}_4 \in \mathcal{B}[f]
\end{array}
$$

is commutative. $\qquad\qquad\square$

## 8.4. Geometric Object of Diagram of Representations

An active transformation changes a basis of the diagram of representations and tuple of $\Omega$-numbers uniformly and coordinates of $\Omega_2$-number relative basis do not change. A passive transformation changes only the basis and it leads to change of coordinates of tuple of $\Omega$-numbers relative to the basis.

THEOREM 8.4.1. *Let passive transformation $s \in GA(f)$ maps basis $\overline{\overline{e}}_1 \in \mathcal{B}[f]$ into basis $\overline{\overline{e}}_2 \in \mathcal{B}[f]$*

$$(8.4.1) \qquad \overline{\overline{e}}_2 = \overline{\overline{e}}_1 \circ s = \overline{\overline{e}}_1 \circ W[f, \overline{\overline{e}}_1, \overline{\overline{e}}_1 \circ s]$$

*Let tuple of A-numbers $a$ has tuple of $\Omega$-words*

$$(8.4.2) \qquad a = \overline{\overline{e}}_1 \circ W[f, \overline{\overline{e}}_1, a]$$

*relative to basis $\overline{\overline{e}}_1$ and has tuple of $\Omega$-words*

$$(8.4.3) \qquad a = \overline{\overline{e}}_2 \circ W[f, \overline{\overline{e}}_2, a]$$

*relative to basis $\overline{\overline{e}}_2$. Coordinate transformation*

$$(8.4.4) \qquad W[f, \overline{\overline{e}}_2, a] = W[f, \overline{\overline{e}}_1 \circ s, \overline{\overline{e}}_1] \circ W[f, \overline{\overline{e}}_1, a]$$

*does not depend on tuple of A-numbers $a$ or basis $\overline{\overline{e}}_1$, but is defined only by coordinates of tuple of A-numbers $a$ relative to basis $\overline{\overline{e}}_1$.*

PROOF. From (8.4.1) and (8.4.3), it follows that

$$(8.4.5) \qquad
\begin{aligned}
\overline{\overline{e}}_1 \circ W[f, \overline{\overline{e}}_1, a] &= \overline{\overline{e}}_2 \circ W[f, \overline{\overline{e}}_2, a] = \overline{\overline{e}}_1 \circ W[f, \overline{\overline{e}}_1, \overline{\overline{e}}_2] \circ W[f, \overline{\overline{e}}_2, a] \\
&= \overline{\overline{e}}_1 \circ W[f, \overline{\overline{e}}_1, \overline{\overline{e}}_1 \circ s] \circ W[f, \overline{\overline{e}}_2, a]
\end{aligned}$$

Comparing (8.4.2) and (8.4.5) we get

$$(8.4.6) \qquad W[f, \overline{\overline{e}}_1, a] = W[f, \overline{\overline{e}}_1, \overline{\overline{e}}_1 \circ s] \circ W[f, \overline{\overline{e}}_2, a]$$

Since $s$ is automorphism, then the equality (8.4.4) follows from (8.4.6) and the theorem 8.2.13. $\qquad\square$

THEOREM 8.4.2. *Coordinate transformations (8.4.4) form effective contravariant right-side representation of group $GA(f)$ which is called* **coordinate representation** *in tuple of $\Omega$-algebras.*

Proof. According to corollary 8.1.19, the transformation (8.4.4) is the endomorphism of diagram of representations [8.13] $(F, W[f, \overline{\overline{e}}_1])$.

Suppose we have two consecutive passive transformations $s$ and $t$. Coordinate transformation

$$(8.4.7) \qquad W[f, \overline{\overline{e}}_2, a] = W[f, \overline{\overline{e}}_1 \circ s, \overline{\overline{e}}_1] \circ W[f, \overline{\overline{e}}_1, a]$$

corresponds to passive transformation $s$. Coordinate transformation

$$(8.4.8) \qquad W[f, \overline{\overline{e}}_2, a] = W[f, \overline{\overline{e}}_1 \circ t, \overline{\overline{e}}_1] \circ W[f, \overline{\overline{e}}_1, a]$$

corresponds to passive transformation $t$. According to the theorem 8.3.3, product of coordinate transformations (8.4.7) and (8.4.8) has form

$$W[f, \overline{\overline{e}}_3, a] = W[f, \overline{\overline{e}}_1 \circ t, \overline{\overline{e}}_1] \circ W[f, \overline{\overline{e}}_1 \circ s, \overline{\overline{e}}_1] \circ W[f, \overline{\overline{e}}_1, a]$$
$$= W[f, \overline{\overline{e}}_1 \circ t \circ s, \overline{\overline{e}}_1] \circ W[f, \overline{\overline{e}}_1, a]$$

and is coordinate transformation corresponding to passive transformation $s \circ t$. According to theorems 8.2.13, 8.2.14 and to the definition 5.1.11, coordinate transformations form right-side contravariant representation of group $GA(f)$.

Suppose coordinate transformation does not change coordinates of selected basis. Then unit of group $GA(f)$ corresponds to it because representation is single transitive. Therefore, coordinate representation is effective.                    □

Let $(f, A)$, $(g, B)$ be diagrams of representations. Passive representation $P(g)$ is coordinated with passive representation $P(f)$, if there exists homomorphism $h$ of group $GA(f)$ into group $GA(g)$. Consider diagram

$$
\begin{array}{ccc}
\mathrm{End}(B[f]) & \xrightarrow{\;\;H\;\;} & \mathrm{End}(B[g]) \\
{\scriptstyle P(f)}\big\uparrow & \nearrow{\scriptstyle f'} & \big\uparrow{\scriptstyle P(g)} \\
GA(f) & \xrightarrow{\;\;h\;\;} & GA(g)
\end{array}
$$

Since maps $P(f)$, $P(g)$ are isomorphisms of group, then map $H$ is homomorphism of groups. Therefore, map $f'$ is representation of group $GA(f)$ in basis manifold $\mathcal{B}(g)$. According to design, passive transformation $H(s)$ of basis manifold $\mathcal{B}(g)$ corresponds to passive transformation $s$ of basis manifold $\mathcal{B}(f)$

$$(8.4.9) \qquad\qquad \overline{\overline{e}}_{g1} = \overline{\overline{e}}_g \circ H(s)$$

Then coordinate transformation in diagram of representations $(B, g)$ gets form

$$(8.4.10) \qquad W[g, \overline{\overline{e}}_{g1}, a] = W[g, \overline{\overline{e}}_g \circ H(s), \overline{\overline{e}}_g] \circ W[g, \overline{\overline{e}}_g, a]$$

Definition 8.4.3. *Orbit*

$$\mathcal{O}(f, g, \overline{\overline{e}}_g, a) = H(GA(f)) \circ W[g, \overline{\overline{e}}_g, a]$$
$$= (W[g, \overline{\overline{e}}_g \circ H(s), \overline{\overline{e}}_g] \circ W[g, \overline{\overline{e}}_g, a], \overline{\overline{e}}_f \circ s, s \in GA(f))$$

*is called* **geometric object in coordinate representation** *defined in the diagram of representations* $(f, A)$. *For any basis* $\overline{\overline{e}}_{f1} = \overline{\overline{e}}_f \circ s$ *corresponding point* (8.4.10) *of orbit defines* **coordinates of geometric object** *relative basis* $\overline{\overline{e}}_{f1}$.                    □

---

[8.13]This transformation does not generate an endomorphism of the diagram of representations $(f, A)$. Coordinates change because basis relative which we determinate coordinates changes. However, tuple of $A$-numbers, coordinates of which we are considering, does not change.

DEFINITION 8.4.4. *Orbit*

$$\mathcal{O}(f,g,a) = (W[g,\overline{\overline{e}}_g \circ H(s),\overline{\overline{e}}_g] \circ W[g,\overline{\overline{e}}_g,a], \overline{\overline{e}}_g \circ H(s),\overline{\overline{e}}_f \circ s, s \in GA(f))$$

*is called* **geometric object** *defined in the diagram of representations* $(f,A)$. *We also say that* $a$ *is a* **geometric object of type** $H$. *For any basis* $\overline{\overline{e}}_{f1} = \overline{\overline{e}}_f \circ s$ *corresponding point* (8.4.10) *of orbit defines tuple of* $A$-*numbers*

$$a = \overline{\overline{e}}_g \circ W[g,\overline{\overline{e}}_g,a]$$

*called* **representative of geometric object** *in the diagram of representations* $(f,A)$.                                                                   □

Since a geometric object is an orbit of representation, we see that according to the theorem 5.3.7 the definition of the geometric object is a proper definition.

Definition 8.4.3 introduces a geometric object in coordinate space. We assume in definition 8.4.4 that we selected a basis of representation $g$. This allows using a representative of the geometric object instead of its coordinates.

THEOREM 8.4.5 (invariance principle). *Representative of geometric object does not depend on selection of basis* $\overline{\overline{e}}_f$.

PROOF. To define representative of geometric object, we need to select basis $\overline{\overline{e}}_f$ of diagram of representations $(f,A)$, basis $\overline{\overline{e}}_g$ of diagram of representations $(B,g)$ and coordinates of geometric object $W[g,\overline{\overline{e}}_g,b]$. Corresponding representative of geometric object has form

$$b = \overline{\overline{e}}_g \circ W[g,\overline{\overline{e}}_g,b]$$

Suppose we map basis $\overline{\overline{e}}_f$ to basis $\overline{\overline{e}}_{f1}$ by passive transformation

$$\overline{\overline{e}}_{f1} = \overline{\overline{e}}_f \circ s$$

According building this forms passive transformation (8.4.9) and coordinate transformation (8.4.10). Corresponding representative of geometric object has form

$$\begin{aligned}
b' &= \overline{\overline{e}}_{g1} \circ W[g,\overline{\overline{e}}_{g1},b'] \\
&= \overline{\overline{e}}_g \circ W[g,\overline{\overline{e}}_g,\overline{\overline{e}}_g \circ H(s)] \circ W[g,\overline{\overline{e}}_g \circ H(s),\overline{\overline{e}}_g] \circ W[g,\overline{\overline{e}}_g,b] \\
&= \overline{\overline{e}}_g \circ W[g,\overline{\overline{e}}_g,b] = b
\end{aligned}$$

Therefore representative of geometric object is invariant relative selection of basis.

                                                                   □

CHAPTER 9

# Examples of Diagram of Representations: Module

## 9.1. About this Chapter

Representation of universal algebra is very important tool which I have been using for many years to study algebra, geometry, calculus. The main goal of this and following chapters is to show how representation of universal algebra works in different fields of mathematics.

Examples in this chapter are related to different structures associated with module over ring.

Abelian group is first example. Module is effective representation of ring in Abelian group. So there is some similarity between Abelian group and module. I consider this similarity in the section 9.2.

Module over commutative ring is relatively simple theory. On the other hand, definitions of representation theory (basis of representation, morphism of representation, free representation) are based on similar definitions in module. So I dedicated the section 9.3 to thorough consideration of module over commutative ring.

I consider algebra over commutative ring in the section 9.4 and left module over $D$-algebra in the section 9.5. We can consider module over non commutative ring the same way as we considered module over commutative ring. However, we encounter serious problems when studying linear map.

Consideration a non commutative ring as algebra over center of the ring significantly changes the picture. Analysis of diagram of representations describing module $V$ over $D$-algebra $A$, allows us to consider different groups of maps which preserve the structure of algebra. Among these maps, we distinguish linear maps of $A$-module $V$ (reduced morphism of $D$-module $V$) and homomorphisms of $A$-module $V$ (reduced morphism of diagram of representations). Such definition of linear map allows us to consider polylinear map of module over $D$-algebra $A$.

If $D$-algebra $A$ is Banach algebra, then we get the tool to study multivariable calculus. Unfortunately, the structure of linear map of noncommutative algebra is outside the scope of this chapter. The reader can study more on this topic in the book [12].

## 9.2. Abelian Group

DEFINITION 9.2.1. *The action of ring of rational integers $Z$ in Abelian group $G$ is defined using following rules*

$$(9.2.1) \qquad 0g = 0$$
$$(9.2.2) \qquad (n+1)g = ng + g$$
$$(9.2.3) \qquad (n-1)g = ng - g$$

$\square$

Theorem 9.2.2. *The action of ring of rational integers $Z$ in Abelian group $G$ defined in the definition 9.2.1 is representation. The following equalities are true*

(9.2.4)
$$1a = a$$

(9.2.5)
$$(nm)a = n(ma)$$

(9.2.6)
$$(m + n)a = ma + na$$

(9.2.7)
$$(m - n)a = ma - na$$

(9.2.8)
$$n(a + b) = na + nb$$

Proof. The equality (9.2.4) follows from the equality (9.2.1) and from the equality (9.2.2) when $n = 0$.

From the equality (9.2.1), it follows that the equality (9.2.6) is true when $n = 0$.

- Let the equality (9.2.6) is true when $n = k \geq 0$. Then
$$(m + k)a = ma + ka$$
The equality
$$(m + (k + 1))a = ((m + k) + 1)a = (m + k)a + a = ma + ka + a$$
$$= ma + (k + 1)a$$
follows from the equality (9.2.2). Therefore, the equality (9.2.6) is true when $n = k + 1$. According to mathematical induction, the equality (9.2.6) is true for any $n \geq 0$.

- Let the equality (9.2.6) is true when $n = k \leq 0$. Then
$$(m + k)a = ma + ka$$
The equality
$$(m + (k - 1))a = ((m + k) - 1)a = (m + k)a - a = ma + ka - a$$
$$= ma + (k - 1)a$$
follows from the equality (9.2.3). Therefore, the equality (9.2.6) is true when $n = k - 1$. According to mathematical induction, the equality (9.2.6) is true for any $n \leq 0$.

- Therefore, the equality (9.2.6) is true for any $n \in Z$.

The equality

(9.2.9)
$$(k + n)a - na = ka$$

follows from the equality (9.2.6). The equality (9.2.7) follows from the equality (9.2.9), if we assume $m = k + n$, $k = m - n$.

From the equality (9.2.1), it follows that the equality (9.2.5) is true when $n = 0$.

- Let the equality (9.2.5) is true when $n = k \geq 0$. Then
$$(km)a = k(ma)$$
The equality
$$((k + 1)m)a = (km + m)a = (km)a + ma = k(ma) + ma$$
$$= (k + 1)(ma)$$

follows from equalities (9.2.2), (9.2.6). Therefore, the equality (9.2.5) is true when $n = k + 1$. According to mathematical induction, the equality (9.2.5) is true for any $n \geq 0$.

- Let the equality (9.2.6) is true when $n = k \leq 0$. Then

$$(km)a = k(ma)$$

The equality

$$((k-1)m)a = (km - m)a = (km)a - ma = k(ma) - ma$$
$$= (k-1)(ma)$$

follows from equalities (9.2.3), (9.2.7). Therefore, the equality (9.2.5) is true when $n = k - 1$. According to mathematical induction, the equality (9.2.5) is true for any $n \leq 0$.

- Therefore, the equality (9.2.5) is true for any $n \in Z$.

From the equality (9.2.1), it follows that the equality (9.2.8) is true when $n = 0$.

- Let the equality (9.2.8) is true when $n = k \geq 0$. Then

$$k(a + b) = ka + kb$$

The equality

$$(k+1)(a+b) = k(a+b) + a + b = ka + kb + a + b$$
$$= ka + a + kb + b$$
$$= (k+1)a + (k+1)b$$

follows from the equality (9.2.2). Therefore, the equality (9.2.8) is true when $n = k + 1$. According to mathematical induction, the equality (9.2.8) is true for any $n \geq 0$.

- Let the equality (9.2.6) is true when $n = k \leq 0$. Then

$$k(a + b) = ka + kb$$

The equality

$$(k-1)(a+b) = k(a+b) - (a+b) = ka + kb - a - b$$
$$= ka - a + kb - b$$
$$= (k-1)a + (k-1)b$$

follows from the equality (9.2.3). Therefore, the equality (9.2.8) is true when $n = k - 1$. According to mathematical induction, the equality (9.2.8) is true for any $n \leq 0$.

- Therefore, the equality (9.2.8) is true for any $n \in Z$.

From the equality (9.2.8), it follows that the map

$$\varphi(n) : a \in G \to na \in G$$

is an endomorphism of Abelian group $G$. From equalities (9.2.6), (9.2.5), it follows that the map

$$\varphi : Z \to \text{End}(Ab, G)$$

is a homomorphism of the ring $Z$. According to the definition 3.1.1, the map $\varphi$ is representation of ring of rational integers $Z$ in Abelian group $G$.                □

THEOREM 9.2.3. *Let $G$ be Abelian group. The set of $G$-numbers generated by the set $S = \{s_i : i \in I\}$ has form*

$$(9.2.10) \qquad J(S) = \left\{ g : g = \sum_{i \in I} g^i s_i, g^i \in Z \right\}$$

*where the set $\{i \in I : g^i \neq 0\}$ is finite.*

PROOF. We prove the theorem by induction based on the theorems [14]-5.1, page 79 and 6.1.4.

For any $s_k \in S$, let $g^i = \delta^i_k$. Then

$$(9.2.11) \qquad s_k = \sum_{i \in I} g^i s_i$$

$s_k \in J(S)$ follows from (9.2.10), (9.2.11).

Let $g_1, g_2 \in X_k \subseteq J(S)$. Since $G$ is Abelian group, then, according to the statement 6.1.4.3, $g_1 + g_2 \in J(S)$. According to the equality (9.2.10), there exist $Z$-numbers $g_1^i, g_2^i, i \in I$, such that

$$(9.2.12) \qquad g_1 = \sum_{i \in I} g_1^i v_i \quad g_2 = \sum_{i \in I} g_1^i v_i$$

where sets

$$(9.2.13) \qquad H_1 = \{i \in I : g_1^i \neq 0\} \quad H_2 = \{i \in I : g_2^i \neq 0\}$$

are finite. From the equality (9.2.12), it follows that

$$(9.2.14) \qquad g_1 + g_2 = \sum_{i \in I} g_1^i v_i + \sum_{i \in I} g_2^i v_i = \sum_{i \in I} (g_1^i v_i + g_2^i v_i)$$

The equality

$$(9.2.15) \qquad g_1 + g_2 = \sum_{i \in I} (g_1^i + g_2^i) v_i$$

follows from equalities (9.2.6), (9.2.14). From the equality (9.2.13), it follows that the set

$$\{i \in I : g_1^i + g_2^i \neq 0\} \subseteq H_1 \cup H_2$$

is finite. From the equality (9.2.15), it follows that $g_1 + g_2 \in J(S)$. $\qquad \square$

## 9.3. Vector Space

### 9.3.1. Module over Commutative Ring.

DEFINITION 9.3.1. *Effective representation of commutative ring $D$ in an Abelian group $V$*

$$(9.3.1) \qquad f : D \overset{*}{\longrightarrow} V \quad f(d) : v \to dv$$

*is called* **module over ring** $D$ *or* $D$**-module**. $V$*-number is called* **vector**. $\qquad \square$

THEOREM 9.3.2. *The following diagram of representations describes $D$-module $V$*

$$(9.3.2) \qquad \begin{array}{ccc} D & \overset{g_2}{\longrightarrow} & V \\ & & \uparrow g_1 \\ & & Z \end{array}$$

*The diagram of representations* (9.3.2) *holds* **commutativity of representations of ring of rational integers** $Z$ *and commutative ring* $D$ *in Abelian group* $V$

$$(9.3.3) \qquad\qquad a(nv) = n(av)$$

PROOF. The diagram of representations (9.3.2) follows from the definition 9.3.1 and from the theorem 9.2.2. Since transformation $g_2(a)$ is endomorphism of $Z$-module $V$, we obtain the equality (9.3.3). $\qquad\square$

THEOREM 9.3.3. *Let* $V$ *be* $D$-*module. For any vector* $v \in V$, *vector generated by the diagram of representations* (9.3.2) *has the following form*

$$(9.3.4) \qquad\qquad (a+n)v = av + nv \quad a \in D \quad n \in Z$$

9.3.3.1: *The set of maps*

$$(9.3.5) \qquad\qquad a+n : v \in V \to (a+n)v \in V$$

*generates*[9.1] *ring* $D_{(1)}$ *where the sum is defined by the equality*

$$(9.3.6) \qquad\qquad (a+n) + (b+m) = (a+b) + (n+m)$$

*and the product is defined by the equality*

$$(9.3.7) \qquad\qquad (a+n)(b+m) = (ab + ma + nb) + (nm)$$

*The ring* $D_{(1)}$ *is called* **unital extension** *of the ring* $D$.

| *If ring* $D$ *has unit, then* | $Z \subseteq D$ | $D_{(1)} = D$ |
|---|---|---|
| *If ring* $D$ *is ideal of* $Z$, *then* | $D \subseteq Z$ | $D_{(1)} = Z$ |
| *Otherwise* | | $D_{(1)} = D \oplus Z$ |

9.3.3.2: *The ring* $D$ *is ideal of ring* $D_{(1)}$.
9.3.3.3: *The set of transormations* (9.3.4) *is representation of ring* $D_{(1)}$ *in Abelian group* $V$.

*We use the notation* $D_{(1)}v$ *for the set of vectors generated by vector* $v$.

THEOREM 9.3.4. *Following conditions hold for* $D$-*module* $V$:

9.3.4.1: **associative law**

$$(9.3.8) \qquad\qquad (pq)v = p(qv)$$

9.3.4.2: **distributive law**

$$(9.3.9) \qquad\qquad p(v+w) = pv + pw$$
$$(9.3.10) \qquad\qquad (p+q)v = pv + qv$$

9.3.4.3: **unitarity law**

$$(9.3.11) \qquad\qquad\qquad 1v = v$$

*for any* $p, q \in D_{(1)}, v, w \in V$.

PROOF OF THEOREMS 9.3.3, 9.3.4. Let $v \in V$.

LEMMA 9.3.5. *Let* $n \in Z, \quad a \in D$. *The map* (9.3.5) *is endomorphism of Abelian group* $V$.

---

[9.1] See the definition of unital extension also on the pages [6]-52, [7]-64.

PROOF. Statements $nv \in V$, $av \in V$ follow from the theorems 6.1.4, 9.3.2. Since $V$ is Abelian group, then

$$nv + av \in V \quad n \in Z \quad a \in D$$

Therefore, for any $Z$-number $n$ and for any $D$-number $a$, we defined the map (9.3.5). Since transformation $g_1(n)$ and transformation $g_2(a)$ are endomorphisms of Abelian group $V$, then the map (9.3.5) is endomorphism of Abelian group $V$. $\odot$

Let $D_{(1)}$ be the set of maps (9.3.5). The equality (9.3.9) follows from the lemma 9.3.5.

Let $p = a + n \in D_{(1)}$, $q = b + m \in D_{(1)}$. According to the statement 9.3.3.3, we define the sum of $D_{(1)}$-numbers $p$ and $q$ by the equality (9.3.10). The equality

$$(9.3.12) \qquad ((a + n) + (b + m))v = (a + n)v + (b + m)v$$

follows from the equality (9.3.10). Since representation $g_1$ is homomorphism of the aditive group of ring $Z$, we obtain the equality

$$(9.3.13) \qquad (n + m)v = nv + mv$$

Since representation $g_2$ is homomorphism of the aditive group of ring $D$, we obtain the equality

$$(9.3.14) \qquad (a + b)v = av + bv$$

Since $V$ is Abelian group, then the equality

$$(9.3.15) \qquad \begin{aligned} ((a + n) + (b + m))v &= av + nv + bv + mv = av + bv + nv + mv \\ &= (a + b)v + (n + m)v = ((a + b) + (n + m))v \end{aligned}$$

follows from equalities (9.3.12), (9.3.13), (9.3.14). From the equality (9.3.15), it follows that the definition (9.3.6) of sum on the set $D_{(1)}$ does not depend on vector $v$.

Equalities (9.3.8), (9.3.11) follow from the statement 9.3.3.3. Let $p = a + n \in D_{(1)}$, $q = b + m \in D_{(1)}$. Since representation $g_1$ is representation of the multiplicative group of ring $Z$, we obtain the equality

$$(9.3.16) \qquad (mn)v = m(nv)$$

Since representation $g_2$ is representation of the multiplicative group of ring $D$, we obtain the equality

$$(9.3.17) \qquad (ab)v = a(bv)$$

Since the ring $D$ is Abelian group, we obtain the equality

$$(9.3.18) \qquad (md)v = m(dv)$$

The equality

$$(9.3.19) \qquad \begin{aligned} ((a + n)(b + m))v &= (a + n)((b + m)v) = (a + n)(bv + mv) \\ &= a(bv + mv) + n(bv + mv) \\ &= a(bv) + a(mv) + n(bv) + n(mv) \\ &= (ab)v + m(av) + + n(bv) + (nm)v \\ &= (ab)v + (ma)v + + (nb)v + (nm)v \\ &= ((ab + ma + nb) + nm)v \end{aligned}$$

follows from equalities (9.3.3), (9.3.4), (9.3.8), (9.3.16), (9.3.17), (9.3.18). The equality (9.3.7) follows from the equality (9.3.19).

The statement 9.3.3.2 follows from the equality (9.3.7). □

Theorem 9.3.6. *Let $V$ be $D$-module. The set of vectors generated by the set of vectors $v = (v_i \in V, i \in I)$ has form* [9.2]

$$(9.3.20) \qquad J(v) = \left\{ w : w = \sum_{i \in I} c^i v_i, c^i \in D_{(1)}, |\{i : c^i \neq 0\}| < \infty \right\}$$

Proof. We prove the theorem by induction based on the theorem 6.1.4, Acording to the theorem 6.1.4, we need to prove following statements:

9.3.6.1: $v_k \in X_0 \subseteq J(v)$

9.3.6.2: $c^k v_k \in J(v)$, $c^k \in D_{(1)}$, $k \in I$

9.3.6.3: $\sum_{k \in I} c^k v_k \in J(v)$, $c^k \in D_{(1)}$, $|\{i : c^i \neq 0\}| < \infty$

9.3.6.4: $w_1, w_2 \in J(v) \Rightarrow w_1 + w_2 \in J(v)$

9.3.6.5: $a \in D$, $w \in J(v) \Rightarrow aw \in J(v)$

- For any $v_k \in v$, let $c^i = \delta_k^i \in D_{(1)}$. Then

$$(9.3.21) \qquad v_k = \sum_{i \in I} c^i v_i$$

  The statement 9.3.6.1 follows from (9.3.20), (9.3.21).

- The statement 9.3.6.2 follow from the theorems 6.1.4, 9.3.3 and from the statement 9.3.6.1.

- Since $V$ is Abelian group, then the statement 9.3.6.3 follows from the statement 9.3.6.2 and from theorems 6.1.4, 9.2.3.

- Let $w_1, w_2 \in X_k \subseteq J(v)$. Since $V$ is Abelian group, then, according to the statement 6.1.4.3,

$$(9.3.22) \qquad w_1 + w_2 \in X_{k+1}$$

  According to the equality (9.3.20), there exist $D_{(1)}$-numbers $w_1^i$, $w_2^i$, $i \in I$, such that

$$(9.3.23) \qquad w_1 = \sum_{i \in I} w_1^i v_i \qquad w_2 = \sum_{i \in I} w_1^i v_i$$

  where sets

$$(9.3.24) \qquad H_1 = \{i \in I : w_1^i \neq 0\} \qquad H_2 = \{i \in I : w_2^i \neq 0\}$$

  are finite. Since $V$ is Abelian group, then from the equality (9.3.23) it follows that

$$(9.3.25) \qquad w_1 + w_2 = \sum_{i \in I} w_1^i v_i + \sum_{i \in I} w_2^i v_i = \sum_{i \in I} (w_1^i v_i + w_2^i v_i)$$

  The equality

$$(9.3.26) \qquad w_1 + w_2 = \sum_{i \in I} (w_1^i + w_2^i) v_i$$

---

[9.2] For a set $A$, we denote by $|A|$ the cardinal number of the set $A$. The notation $|A| < \infty$ means that the set $A$ is finite.

follows from equalities (9.3.10), (9.3.25). From the equality (9.3.24), it follows that the set

$$\{i \in I : w_1^i + w_2^i \neq 0\} \subseteq H_1 \cup H_2$$

is finite.

- Let $w \in X_k \subseteq J(v)$. According to the statement 6.1.4.4, for any $D_{(1)}$-number $a$,

$$(9.3.27) \qquad\qquad aw \in X_{k+1}$$

According to the equality (9.3.20), there exist $D_{(1)}$-numbers $w^i$, $i \in I$, such that

$$(9.3.28) \qquad\qquad w = \sum_{i \in I} w^i v_i$$

where

$$(9.3.29) \qquad\qquad |\{i \in I : w^i \neq 0\}| < \infty$$

From the equality (9.3.28) it follows that

$$(9.3.30) \qquad aw = a\sum_{i \in I} w^i v_i = \sum_{i \in I} a(w^i v_i) = \sum_{i \in I}(aw^i)v_i$$

From the statement (9.3.29), it follows that the set $\{i \in I : aw^i \neq 0\}$ is finite.

From equalities (9.3.22), (9.3.26), (9.3.27), (9.3.30), it follows that $X_{k+1} \subseteq J(v)$.
$\square$

DEFINITION 9.3.7. *Let* $v = (v_i \in V, i \in I)$ *be set of vectors. The expression* $w^i v_i$ *is called* **linear combination** *of vectors* $v_i$. *A vector* $\overline{w} = w^i v_i$ *is called* **linearly dependent** *on vectors* $v_i$. $\square$

We represent the set of $D_{(1)}$-numbers $w^i$, $i \in I$, as matrix

$$w = \begin{pmatrix} w^1 \\ ... \\ w^n \end{pmatrix}$$

We represent the set of vectors $v_i$, $i \in I$, as matrix

$$v = \begin{pmatrix} v_1 & ... & v_n \end{pmatrix}$$

Then we can represent linear combination of vectors $\overline{w} = w^i v_i$ as

$$\overline{w} = w^*_* v$$

THEOREM 9.3.8. *Let $D$ be field. Since the equation*

$$w^i v_i = 0$$

*implies existence of index $i = j$ such that $w^j \neq 0$, then the vector $v_j$ linearly depends on rest of vectors $v$.*

Proof. The theorem follows from the equality

$$v_j = \sum_{i \in I \setminus \{j\}} \frac{w^i}{w^j} v_i$$

and from the definition 9.3.7.                                                                    □

It is evident that for any set of vectors $v_i$

$$w^i = 0 \Rightarrow w^*{}_* v = 0$$

Definition 9.3.9. *The set of vectors*[9.3] $v_i$, $i \in I$, *of $D$-module $V$ is* **linearly independent** *if $w = 0$ follows from the equation*

$$w^i v_i = 0$$

*Otherwise the set of vectors $v_i$, $i \in I$, is* **linearly dependent**.                    □

The following definition follows from the theorems 9.3.6, 6.1.4 and from the definition 6.1.5.

Definition 9.3.10. *$J(v)$ is called* **submodule generated by set** *$v$, and $v$ is a* **generating set** *of submodule $J(v)$. In particular, a* **generating set** *of $D$-module $V$ is a subset $X \subset V$ such that $J(X) = V$.*                                □

The following definition follows from the theorems 9.3.6, 6.1.4 and from the definition 6.2.6.

Definition 9.3.11. *If the set $X \subset V$ is generating set of $D$-module $V$, then any set $Y$, $X \subset Y \subset V$ also is generating set of $D$-module $V$. If there exists minimal set $X$ generating the $D$-module $V$, then the set $X$ is called* **basis** *of $D$-module $V$.*                                                              □

Theorem 9.3.12. *The set of vectors $\overline{\overline{e}} = (e_i, i \in I)$ is basis of $D$-module $V$, if following statements are true.*

9.3.12.1: *Arbitrary vector $v \in V$ is linear combination of vectors of the set $\overline{\overline{e}}$.*

9.3.12.2: *Vector $e_i$ cannot be represented as a linear combination of the remaining vectors of the set $\overline{\overline{e}}$.*

Proof. According to the statement 9.3.12.1, the theorem 9.3.6 and the definition 9.3.7, the set $\overline{\overline{e}}$ generates $D$-module $V$ (the definition 9.3.10). According to the statement 9.3.12.2, the set $\overline{\overline{e}}$ is minimal set generating $D$-module $V$. According to the definitions 9.3.11, the set $\overline{\overline{e}}$ is a basis of $D$-module $V$.                        □

Theorem 9.3.13. *Let $D$ be field. The set of vectors $\overline{\overline{e}} = (e_i, i \in I)$ is a* **basis of $D$-vector space** *$V$ if vectors $e_i$ are linearly independent and any vector $v \in V$ linearly depends on vectors $e_i$.*

Proof. Let the set of vectors $e_i$, $i \in I$, be linear dependent. Then the equation

$$w^i e_i = 0$$

implies existence of index $i = j$ such that $w^j \neq 0$. According to the theorem 9.3.8, the vector $e_j$ linearly depends on rest of vectors of the set $\overline{\overline{e}}$. According to the definition 9.3.11, the set of vectors $e_i$, $i \in I$, is not a basis for $D$-vector space $V$.

---

[9.3] I follow to the definition in [2], page 130.

Therefore, if the set of vectors $e_i$, $i \in I$, is a basis, then these vectors are linearly independent. Since an arbitrary vector $v \in V$ is linear combination of vectors $e_i$, $i \in I$, , then the set of vectors $v$, $e_i$, $i \in I$, is not linearly independent. $\square$

**Definition 9.3.14.** *Let $\overline{\overline{e}}$ be the basis of $D$-module $V$ and vector $\overline{v} \in V$ has expansion*

$$\overline{v} = v^*{}_*e = v^i e_i$$

*with respect to the basis $\overline{\overline{e}}$. $D_{(1)}$-numbers $v^i$ are called* **coordinates** *of vector $\overline{v}$ with respect to the basis $\overline{\overline{e}}$. Matrix of $D_{(1)}$-numbers $v = (v^i, i \in I)$ is called* **coordinate matrix of vector** $\overline{v}$ *in basis* $\overline{\overline{e}}$. $\square$

**Theorem 9.3.15.** *Let $D$ be ring. Let $\overline{\overline{e}}$ be basis of $D$-module $V$. Let*

$$(9.3.31) \qquad\qquad w^i e_i = 0$$

*be linear dependence of vectors of the basis $\overline{\overline{e}}$. Then*

9.3.15.1: $D_{(1)}$-*number $w^i$, $i \in I$, does not have inverse element in ring $D_{(1)}$.*

9.3.15.2: *The set $D'$ of matrices $w = (w^i, i \in I)$ generates $D$-module.*

Proof. Let there exist matrix $w = (w^i, i \in I)$ such that the equality (9.3.31) is true and there exist index $i = j$ such that $w^j \neq 0$. If we assume that $D_{(1)}$-number $c^j$ has inverse one, then the equality

$$e_j = \sum_{i \in I \setminus \{j\}} \frac{w^i}{w^j} e_i$$

follows from the equality (9.3.31). Therefore, the vector $e_j$ is linear combination of other vectors of the set $\overline{\overline{e}}$ and the set $\overline{\overline{e}}$ is not basis. Therefore, our assumption is false, and $D_{(1)}$-number $c^j$ does not have inverse.

Let matrices $b = (b^i, i \in I) \in D'$, $c = (c^i, i \in I) \in D'$. From equalities

$$b^i e_i = 0$$

$$c^i e_i = 0$$

it follows that

$$(b^i + c^i)e_i = 0$$

Therefore, the set $D'$ is Abelian group.

Let matrix $c = (c^i, i \in I) \in D'$ and $a \in D$. From the equality

$$w^i e_i = 0$$

it follows that

$$(ac^i)e_i = 0$$

Therefore, Abelian group $D'$ is $D$-module. $\square$

**Theorem 9.3.16.** *Let $D$-module $V$ have the basis $\overline{\overline{e}}$ such that in the equality*

$$(9.3.32) \qquad\qquad w^i e_i = 0$$

*there exists index $i = j$ such that $w^j \neq 0$. Then*

9.3.16.1: *The matrix $w = (w^i, i \in I)$ determines coordinates of vector $0 \in V$ with respect to basis $\overline{\overline{e}}$.*

9.3.16.2: *Coordinates of vector $\overline{v}$ with respect to basis $\overline{\overline{e}}$ are uniquely determined up to a choice of coordinates of vector $0 \in V$.*

PROOF. The statement 9.3.16.1 follows from the equality (9.3.32) and from the definition 9.3.14.

Let vector $\overline{v}$ have expansion

$$(9.3.33) \qquad \overline{v} = v^*{}_* e = v^i e_i$$

with respect to basis $\overline{\overline{e}}$. The equality

$$(9.3.34) \qquad \overline{v} = \overline{v} + 0 = v^i e_i + c^i e_i = (v^i + c^i) e_i$$

follows from equalities (9.3.32), (9.3.33). The statement 9.3.16.2 follows from equalities (9.3.33), (9.3.34) and from the definition 9.3.14. ∎

DEFINITION 9.3.17. *The $D$-module $V$ is **free** $D$-module,[9.4] if $D$-module $V$ has basis and vectors of the basis are linearly independent.* ∎

THEOREM 9.3.18. *Coordinates of vector $v \in V$ relative to basis $\overline{\overline{e}}$ of free $D$-module $V$ are uniquely defined.*

PROOF. The theorem follows from the theorem 9.3.16 and from definitions 9.3.9, 9.3.17. ∎

EXAMPLE 9.3.19. *From the theorem 9.2.2 and the definition 9.3.1, it follows that Abelian group $G$ is module over ring of integers $Z$.* ∎

### 9.3.2. Linear Map.

DEFINITION 9.3.20. *Morphism of representations*

$$\Big( h : D_1 \to D_2 \quad f : V_1 \to V_2 \Big)$$

*of $D_1$-module $A_1$ into $D_2$-module $A_2$ is called **linear map** of $D_1$-module $A_1$ into $D_2$-module $A_2$. Let us denote $\mathcal{L}(D_1 \to D_2; A_1 \to A_2)$ set of linear maps of $D_1$-module $A_1$ into $D_2$-module $A_2$.* ∎

If the map

$$f : A_1 \to A_2$$

is linear map of $D$-algebra $A_1$ into $D$-algebra $A_2$, then I use notation

$$f \circ a = f(a)$$

for image of the map $f$.

THEOREM 9.3.21. *Linear map*

$$\Big( h : D_1 \to D_2 \quad f : A_1 \to A_2 \Big)$$

*of $D_1$-module $A_1$ into $D_2$-module $A_2$ satisfies to equations[9.5]*

$$(9.3.35) \qquad h(d_1 + d_2) = h(d_1) + h(d_2)$$

$$(9.3.36) \qquad h(d_1 d_2) = h(d_1) h(d_2)$$

$$(9.3.37) \qquad f \circ (a + b) = f \circ a + f \circ b$$

$$(9.3.38) \qquad f \circ (da) = h(d)(f \circ a)$$

$$a, b \in A_1 \quad d, d_1, d_2 \in D_1$$

---

[9.4] I follow to the definition in [2], page 135.

[9.5] In some books (for instance, on page [2]-119) the theorem 9.3.21 is considered as a definition.

PROOF. From definitions 3.2.2, 9.3.20, it follows that the map $h$ is a homomorphism of the ring $D_1$ into the ring $D_2$ (the equalities (9.3.35), (9.3.36)) and the map $f$ is a homomorphism of the Abelian group $A_1$ into the Abelian group $A_2$ (the equality (9.3.37)). The equality (9.3.38) follows from the equality (3.2.3). □

THEOREM 9.3.22. *Let*

$$\overline{\overline{e}}_1 = (e_{1 \cdot i}, i \in I)$$

*be a basis of $D_1$-module $A_1$. Let*

$$\overline{\overline{e}}_2 = (e_{2 \cdot j}, j \in J)$$

*be a basis of $D_2$-module $A_2$. Then linear map*

$$\left( h : D_1 \to D_2 \quad \overline{f} : A_1 \to A_2 \right)$$

*has presentation*

$$(9.3.39) \qquad b = h(a)^*{}_* f$$

*relative to selected bases. Here*

- *$a$ is coordinate matrix of $A_1$-number $\overline{a}$ relative the basis $\overline{\overline{e}}_1$*

$$(9.3.40) \qquad \overline{a} = a^*{}_* e_1$$

- *$h(a) = (h(a_i), i \in I)$ is a matrix of $D_2$-numbers.*
- *$b$ is coordinate matrix of vector*

$$(9.3.41) \qquad \overline{b} = \overline{f} \circ \overline{a}$$

  *relative the basis $\overline{\overline{e}}_2$*

$$(9.3.42) \qquad \overline{b} = b^*{}_* e_2$$

- *$f$ is coordinate matrix of set of vectors $(\overline{f} \circ e_{1 \cdot i}, i \in I)$ relative the basis $\overline{\overline{e}}_2$. The matrix $f$ is called **matrix of linear map** $\overline{f}$ relative bases $\overline{\overline{e}}_1$ and $\overline{\overline{e}}_2$.*

PROOF. Since

$$\left( h : D_1 \to D_2 \quad \overline{f} : A_1 \to A_2 \right)$$

is a linear map, then the equality

$$(9.3.43) \qquad \overline{b} = \overline{f} \circ \overline{a} = \overline{f} \circ (a^*{}_* e_1) = h(a)^*{}_* (\overline{f} \circ e_1)$$

follows from equalities (9.3.38), (9.3.40), (9.3.41). $A_2$-number $\overline{f} \circ e_{1 \cdot i}$ has expansion

$$(9.3.44) \qquad \overline{f} \circ e_{1 \cdot i} = f_i{}^*{}_* e_2 = f_i^j e_{2 \cdot j}$$

relative to basis $\overline{\overline{e}}_2$. Combining (9.3.43) and (9.3.44), we get

$$(9.3.45) \qquad \overline{b} = h(a)^*{}_* f^*{}_* e_2$$

(9.3.39) follows from comparison of (9.3.42) and (9.3.45) and theorem 9.3.18. □

DEFINITION 9.3.23. *Reduced morphism of representations*

$$f : A_1 \to A_2$$

*of $D$-module $A_1$ into $D$-module $A_2$ is called **linear map** of $D$-module $A_1$ into $D$-module $A_2$. Let us denote $\mathcal{L}(D; A_1 \to A_2)$ set of linear maps of $D$-module $A_1$ into $D$-module $A_2$.* □

Theorem 9.3.24. *Linear map*

$$f : A_1 \to A_2$$

*of D-module $A_1$ into D-module $A_2$ satisfies to equations*[9.6]

$$(9.3.46) \qquad\qquad f \circ (a + b) = f \circ a + f \circ b$$

$$(9.3.47) \qquad\qquad f \circ (da) = d(f \circ a)$$

$$a, b \in A_1 \quad d \in D$$

Proof. From definitions 3.4.2, 9.3.23, it follows that the map $f$ is a homomorphism of the Abelian group $A_1$ into the Abelian group $A_2$ (the equality (9.3.46)). The equality (9.3.47) follows from the equality (3.4.4). □

Theorem 9.3.25. *Let*

$$\overline{\overline{e}}_1 = (e_{1 \cdot i}, i \in I)$$

*be a basis of D-module $A_1$. Let*

$$\overline{\overline{e}}_2 = (e_{2 \cdot j}, j \in J)$$

*be a basis of D-module $A_2$. Then linear map*

$$\overline{f} : A_1 \to A_2$$

*has presentation*

$$(9.3.48) \qquad\qquad b = a^*{}_* f$$

*relative to selected bases. Here*

- *a is coordinate matrix of $A_1$-number $\overline{a}$ relative the basis $\overline{\overline{e}}_1$*

$$(9.3.49) \qquad\qquad \overline{a} = a^*{}_* e_1$$

- *b is coordinate matrix of vector*

$$(9.3.50) \qquad\qquad \overline{b} = \overline{f} \circ \overline{a}$$

  *relative the basis $\overline{\overline{e}}_2$*

$$(9.3.51) \qquad\qquad \overline{b} = b^*{}_* e_2$$

- *f is coordinate matrix of set of vectors $(\overline{f} \circ e_{1 \cdot i}, i \in I)$ relative the basis $\overline{\overline{e}}_2$. The matrix f is called **matrix of linear map** $\overline{f}$ relative bases $\overline{\overline{e}}_1$ and $\overline{\overline{e}}_2$.*

Proof. Since

$$\overline{f} : A_1 \to A_2$$

is a linear map, then the equality

$$(9.3.52) \qquad \overline{b} = \overline{f} \circ \overline{a} = \overline{f} \circ (a^*{}_* e_1) = a^*{}_* (\overline{f} \circ e_1)$$

follows from equalities (9.3.47), (9.3.49), (9.3.50). $A_2$-number $\overline{f} \circ e_{1 \cdot i}$ has expansion

$$(9.3.53) \qquad\qquad \overline{f} \circ e_{1 \cdot i} = f_i{}^*{}_* e_2 = f_i^j e_{2 \cdot j}$$

relative to basis $\overline{\overline{e}}_2$. Combining (9.3.52) and (9.3.53), we get

$$(9.3.54) \qquad\qquad \overline{b} = a^*{}_* f^*{}_* e_2$$

(9.3.48) follows from comparison of (9.3.51) and (9.3.54) and theorem 9.3.18. □

---

[9.6] In some books (for instance, on page [2]-119) the theorem 9.3.24 is considered as a definition.

### 9.3.3. Polylinear Map.

DEFINITION 9.3.26. *Let $D$ be the commutative ring. Reduced polymorphism of $D$-modules $A_1$, ..., $A_n$ into $D$-module $S$*

$$f : A_1 \times ... \times A_n \to S$$

*is called* **polylinear map** *of $D$-modules $A_1$, ..., $A_n$ into $D$-module $S$. We denote* $\mathcal{L}(D; A_1 \times ... \times A_n \to S)$ *the set of polylinear maps of $D$-modules $A_1$, ..., $A_n$ into $D$-module $S$. Let us denote* $\mathcal{L}(D; A^n \to S)$ *set of $n$-linear maps of $D$-module $A$ ($A_1 = ... = A_n = A$) into $D$-module $S$.* □

THEOREM 9.3.27. *Let $D$ be the commutative ring. The polylinear map of $D$-modules $A_1$, ..., $A_n$ into $D$-module $S$*

$$f : A_1 \times ... \times A_n \to S$$

*satisfies to equalities*

$$f \circ (a_1, ..., a_i + b_i, ..., a_n) = f \circ (a_1, ..., a_i, ..., a_n) + f \circ (a_1, ..., b_i, ..., a_n)$$

$$f \circ (a_1, ..., pa_i, ..., a_n) = pf \circ (a_1, ..., a_i, ..., a_n)$$

$$f \circ (a_1, ..., a_i + b_i, ..., a_n) = f \circ (a_1, ..., a_i, ..., a_n) + f \circ (a_1, ..., b_i, ..., a_n)$$

$$f \circ (a_1, ..., pa_i, ..., a_n) = pf \circ (a_1, ..., a_i, ..., a_n)$$

$$1 \le i \le n \quad a_i, b_i \in A_i \quad p \in D$$

PROOF. The theorem follows from definitions 4.4.4, 9.3.23, 9.3.26 and from the theorem 9.3.24. □

THEOREM 9.3.28. *Let $D$ be the commutative ring. Let $A_1$, ..., $A_n$, $S$ be $D$-modules. The map*

$$(9.3.55) \qquad f + g : A_1 \times ... \times A_n \to S \quad f, g \in \mathcal{L}(D; A_1 \times ... \times A_n \to S)$$

*defined by the equality*

$$(9.3.56) \qquad (f + g) \circ (a_1, ..., a_n) = f \circ (a_1, ..., a_n) + g \circ (a_1, ..., a_n)$$

*is called* **sum of polylinear maps** *$f$ and $g$ and is polylinear map. The set $\mathcal{L}(D; A_1 \times ... \times A_n \to S)$ is an Abelian group relative sum of maps.*

PROOF. According to the theorem 9.3.27

$$(9.3.57) \quad f \circ (a_1, ..., a_i + b_i, ..., a_n) = f \circ (a_1, ..., a_i, ..., a_n) + f \circ (a_1, ..., b_i, ..., a_n)$$

$$(9.3.58) \qquad f \circ (a_1, ..., pa_i, ..., a_n) = pf \circ (a_1, ..., a_i, ..., a_n)$$

$$(9.3.59) \quad g \circ (a_1, ..., a_i + b_i, ..., a_n) = g \circ (a_1, ..., a_i, ..., a_n) + g \circ (a_1, ..., b_i, ..., a_n)$$

$$(9.3.60) \qquad g \circ (a_1, ..., pa_i, ..., a_n) = pg \circ (a_1, ..., a_i, ..., a_n)$$

The equality

$$
\begin{aligned}
(9.3.61) \quad & (f + g) \circ (x_1, ..., x_i + y_i, ..., x_n) \\
& = f \circ (x_1, ..., x_i + y_i, ..., x_n) + g \circ (x_1, ..., x_i + y_i, ..., x_n) \\
& = f \circ (x_1, ..., x_i, ..., x_n) + f \circ (x_1, ..., y_i, ..., x_n) \\
& \quad + g \circ (x_1, ..., x_i, ..., x_n) + g \circ (x_1, ..., y_i, ..., x_n) \\
& = (f + g) \circ (x_1, ..., x_i, ..., x_n) + (f + g) \circ (x_1, ..., y_i, ..., x_n)
\end{aligned}
$$

follows from the equalities (9.3.56), (9.3.57), (9.3.59). The equality

$$(9.3.62) \quad \begin{aligned} & (f+g) \circ (x_1, ..., px_i, ..., x_n) \\ =& f \circ (x_1, ..., px_i, ..., x_n) + g \circ (x_1, ..., px_i, ..., x_n) \\ =& pf \circ (x_1, ..., x_i, ..., x_n) + pg \circ (x_1, ..., x_i, ..., x_n) \\ =& p(f \circ (x_1, ..., x_i, ..., x_n) + g \circ (x_1, ..., x_i, ..., x_n)) \\ =& p(f+g) \circ (x_1, ..., x_i, ..., x_n) \end{aligned}$$

follows from the equalities (9.3.56), (9.3.58), (9.3.60). From equalities (9.3.61), (9.3.62) and from the theorem 9.3.27, it follows that the map (9.3.55) is linear map of $D$-modules.

Let $f$, $g$, $h \in \mathcal{L}(D; A_1 \times ... \times A_2 \to S)$. For any $a = (a_1, ..., a_n)$, $a_1 \in A_1$, ..., $a_n \in A_n$,

$$\begin{aligned} (f+g) \circ a =& f \circ a + g \circ a = g \circ a + f \circ a \\ =& (g+f) \circ a \\ ((f+g)+h) \circ a =& (f+g) \circ a + h \circ a = (f \circ a + g \circ a) + h \circ a \\ =& f \circ a + (g \circ a + h \circ a) = f \circ a + (g+h) \circ a \\ =& (f+(g+h)) \circ a \end{aligned}$$

Therefore, sum of polylinear maps is commutative and associative.

From the equality (9.3.56), it follows that the map

$$0 : v \in A_1 \times ... \times A_n \to 0 \in S$$

is zero of addition

$$(0+f) \circ (a_1, ..., a_n) = 0 \circ (a_1, ..., a_n) + f \circ (a_1, ..., a_n) = f \circ (a_1, ..., a_n)$$

From the equality (9.3.56), it follows that the map

$$-f : (a_1, ..., a_n) \in A_1 \times ... \times A_n \to -(f \circ (a_1, ..., a_n)) \in S$$

is map inversed to map $f$

$$f + (-f) = 0$$

because

$$\begin{aligned} (f+(-f)) \circ (a_1, ..., a_n) =& f \circ (a_1, ..., a_n) + (-f) \circ (a_1, ..., a_n) \\ =& f \circ (a_1, ..., a_n) - f \circ (a_1, ..., a_n) \\ =& 0 = 0 \circ (a_1, ..., a_n) \end{aligned}$$

From the equality

$$\begin{aligned} (f+g) \circ (a_1, ..., a_n) =& f \circ (a_1, ..., a_n) + g \circ (a_1, ..., a_n) \\ =& g \circ (a_1, ..., a_n) + f \circ (a_1, ..., a_n) \\ =& (g+f) \circ (a_1, ..., a_n) \end{aligned}$$

it follows that sum of maps is commutative. Therefore, the set $\mathcal{L}(D; A_1 \times ... \times A_n \to S)$ is an Abelian group. $\square$

COROLLARY 9.3.29. *Let $A_1$, $A_2$ be $D$-modules. The map*

$$(9.3.63) \qquad f+g : A_1 \to A_2 \quad f, g \in \mathcal{L}(D; A_1 \to A_2)$$

*defined by equation*

(9.3.64)                     $(f + g) \circ x = f \circ x + g \circ x$

*is called* **sum of maps** *$f$ and $g$ and is linear map. The set $\mathcal{L}(D; A_1; A_2)$ is an Abelian group relative sum of maps.*                                                    □

THEOREM 9.3.30. *Let $D$ be the commutative ring. Let $A_1$, ..., $A_n$, $S$ be $D$-modules. The map*

(9.3.65)          $d\,f : A_1 \times ... \times A_n \to S \quad d \in D \quad f \in \mathcal{L}(D; A_1 \times ... \times A_n \to S)$

*defined by equality*

(9.3.66)                  $(d\,f) \circ (a_1, ..., a_n) = d(f \circ (a_1, ..., a_n))$

*is polylinear map and is called* **product of map** *$f$* **over scalar** *$d$. The representation*

(9.3.67)     $a : f \in \mathcal{L}(D; A_1 \times ... \times A_n \to S) \to af \in \mathcal{L}(D; A_1 \times ... \times A_n \to S)$

*of ring $D$ in Abelian group  $\mathcal{L}(D; A_1 \times ... \times A_n \to S)$  generates structure of $D$-module.*

PROOF. According to the theorem 9.3.27

(9.3.68)   $f \circ (a_1, ..., a_i + b_i, ..., a_n) = f \circ (a_1, ..., a_i, ..., a_n) + f \circ (a_1, ..., b_i, ..., a_n)$

(9.3.69)            $f \circ (a_1, ..., pa_i, ..., a_n) = pf \circ (a_1, ..., a_i, ..., a_n)$

The equality

$$
\begin{aligned}
&(pf) \circ (x_1, ..., x_i + y_i, ..., x_n) \\
&= p\, f \circ (x_1, ..., x_i + y_i, ..., x_n) \\
&= p\, (f \circ (x_1, ..., x_i, ..., x_n) + f \circ (x_1, ..., y_i, ..., x_n)) \\
&= p(f \circ (x_1, ..., x_i, ..., x_n)) + p(f \circ (x_1, ..., y_i, ..., x_n)) \\
&= (pf) \circ (x_1, ..., x_i, ..., x_n) + (pf) \circ (x_1, ..., y_i, ..., x_n)
\end{aligned}
$$

(9.3.70)

follows from equalities (9.3.66), (9.3.68). The equality

$$
\begin{aligned}
&(pf) \circ (x_1, ..., qx_i, ..., x_n) \\
&= p(f \circ (x_1, ..., qx_i, ..., x_n)) = pq(f \circ (x_1, ..., x_i, ..., x_n)) \\
&= qp(f \circ (x_1, ..., x_n)) = q(pf) \circ (x_1, ..., x_n)
\end{aligned}
$$

(9.3.71)

follows from equalities (9.3.66), (9.3.69). From equalities (9.3.70), (9.3.71) and from the theorem 9.3.27, it follows that the map (9.3.65) is polylinear map of $D$-modules.

The equality

(9.3.72)                          $(p + q)f = pf + qf$

follows from the equality

$$
\begin{aligned}
((p + q)f) \circ (x_1, ..., x_n) &= (p + q)(f \circ (x_1, ..., x_n)) \\
&= p(f \circ (x_1, ..., x_n)) + q(f \circ (x_1, ..., x_n)) \\
&= (pf) \circ (x_1, ..., x_n) + (qf) \circ (x_1, ..., x_n)
\end{aligned}
$$

The equality

(9.3.73)                          $p(qf) = (pq)f$

follows from the equality

$$(p(qf)) \circ (x_1, ..., x_n) = p\ (qf) \circ (x_1, ..., x_n) = p\ (q\ f \circ (x_1, ..., x_n))$$
$$= (pq)\ f \circ (x_1, ..., x_n) = ((pq)f) \circ (x_1, ..., x_n)$$

From equalities (9.3.72) (9.3.73) it follows that the map (9.3.67) is representation of ring $D$ in Abelian group $\mathcal{L}(D; A_1 \times ... \times A_n \to S)$. Since specified representation is effective, then, according to the definition 9.3.1 and the theorem 9.3.28, Abelian group $\mathcal{L}(D; A_1 \to A_2)$ is $D$-module. $\qquad\square$

COROLLARY 9.3.31. *Let* $A_1$, $A_2$ *be* $D$-*modules. The map*

$$(9.3.74) \qquad\qquad d\,f : A_1 \to A_2 \quad d \in D \quad f \in \mathcal{L}(D; A_1 \to A_2)$$

*defined by equality*

$$(9.3.75) \qquad\qquad (d\,f) \circ x = d(f \circ x)$$

*is linear map and is called* **product of map** $f$ **over scalar** $d$. *The representation*

$$(9.3.76) \qquad a : f \in \mathcal{L}(D; A_1 \to A_2) \to af \in \mathcal{L}(D; A_1 \to A_2)$$

*of ring* $D$ *in Abelian group* $\mathcal{L}(D; A_1 \to A_2)$ *generates structure of* $D$-*module.* $\qquad\square$

## 9.4. Algebra over Commutative Ring

DEFINITION 9.4.1. *Let* $D$ *be commutative ring.* $D$-*module* $A$ *is called* **algebra over ring** $D$ *or* $D$-**algebra**, *if we defined product*[9.7] *in* $A$

$$(9.4.1) \qquad\qquad v\,w = C \circ (v, w)$$

*where* $C$ *is bilinear map*

$$C : A \times A \to A$$

*If* $A$ *is free* $D$-*module, then* $A$ *is called* **free algebra over ring** $D$. $\qquad\square$

THEOREM 9.4.2. *Let* $D$ *be commutative ring and* $A$ *be Abelian group. The diagram of representations*

$$D \xrightarrow{g_{12}} A \xrightarrow{g_{23}} A \qquad g_{12}(d) : v \to d\,v$$
$$\qquad\qquad \uparrow{\scriptstyle g_{12}} \qquad g_{23}(v) : w \to C \circ (v, w)$$
$$\qquad D \qquad\qquad C \in \mathcal{L}(D; A^2 \to A)$$

*generates the structure of* $D$-*algebra* $A$.

PROOF. The structure of $D$-module $A$ is generated by effective representation

$$g_{12} : D \longrightarrow A$$

of ring $D$ in Abelian group $A$.

---

[9.7] I follow the definition given in [20], page 1, [13], page 4. The statement which is true for any $D$-module, is true also for $D$-algebra.

Lemma 9.4.3. *Let the structure of $D$-algebra $A$ defined in $D$-module $A$, be generated by product*

$$v\,w = C \circ (v, w)$$

**Left shift of $D$-module $A$** *defined by equation*

(9.4.2) $$l \circ v : w \in A \to v\,w \in A$$

*generates the representation*

$$A \xrightarrow{g_{23}} A \qquad \begin{aligned} g_{23} &: v \to l \circ v \\ g_{23} \circ v &: w \to (l \circ v) \circ w \end{aligned}$$

*of $D$-module $A$ in $D$-module $A$*

Proof. According to definitions 9.4.1 and 9.3.26, left shift of $D$-module $A$ is linear map. According to the definition 9.3.23, the map $l(v)$ is endomorphism of $D$-module $A$. The equation

(9.4.3) $$(l \circ (v_1 + v_2)) \circ w = (v_1 + v_2)w = v_1 w + v_2 w = (l \circ v_1) \circ w + (l \circ v_2) \circ w$$

follows from the definition 9.3.26 and from the equation (9.4.2). According to the corollary 9.3.29, the equation

(9.4.4) $$l \circ (v_1 + v_2) = l \circ v_1 + l \circ v_2$$

follows from equation (9.4.3). The equation

(9.4.5) $$(l \circ (dv)) \circ w = (dv)w = d(vw) = d((l \circ v) \circ w)$$

follows from the definition 9.3.26 and from the equation (9.4.2). 9.3.29, the equation

(9.4.6) $$l \circ (dv) = d(l \circ v)$$

follows from equation (9.4.5). The lemma follows from equalities (9.4.4), (9.4.6). ⊙

Lemma 9.4.4. *The representation*

$$A \xrightarrow{g_{23}} A \qquad \begin{aligned} g_{23} &: v \to l \circ v \\ g_{23} \circ v &: w \to (l \circ v) \circ w \end{aligned}$$

*of $D$-module $A$ in $D$-module $A$ determines the product in $D$-module $A$ according to rule*

$$ab = (g_{23} \circ a) \circ b$$

Proof. Since map $g_{23} \circ v$ is endomorphism of $D$-module $A$, then

(9.4.7) $$\begin{aligned} (g_{23} \circ v)(w_1 + w_2) &= (g_{23} \circ v) \circ w_1 + (g_{23} \circ v) \circ w_2 \\ (g_{23} \circ v) \circ (dw) &= d((g_{23} \circ v) \circ w) \end{aligned}$$

Since the map $g_{23}$ is linear map

$$g_{23} : A \to \mathcal{L}(D; A \to A)$$

then, according to corollarys 9.3.29, 9.3.31,

(9.4.8) $$(g_{23} \circ (v_1 + v_2)) \circ w = (g_{23} \circ v_1 + g_{23} \circ v_2)(w) = (g_{23} \circ v_1) \circ w + (g_{23} \circ v_2) \circ w$$

(9.4.9) $$(g_{23} \circ (d\,v)) \circ w = (d\,(g_{23} \circ v)) \circ w = d\,((g_{23} \circ v) \circ w)$$

From equations (9.4.7), (9.4.8), (9.4.9) and the definition 9.3.26, it follows that the map $g_{23}$ is bilinear map. Therefore, the map $g_{23}$ determines the product in $D$-module $A$ according to rule

$$ab = (g_{23} \circ a) \circ b$$

⊙

The theorem follows from lemmas 9.4.3, 9.4.4.                                  □

Usually, when we consider the $D$-algebra $A$, we choose a basis $\overline{\overline{e}}$ of corresponding $D$-module $A$. This choice is convenient because if $D$-module $A$ is free $D$-module $A$, then expansion of the vector is unique relative to basis of $D$-module $A$. This, in particular, allows us to define product by specifying structural constants of algebra relative to given basis.

In general, the basis of $R$-module $A$ may appear a generating set. For instance, if in vector space $H$, where we consider quaternion algebra over real field, we consider the basis

(9.4.10)              $e_0 = 1 \quad e_1 = i \quad e_2 = j \quad e_3 = k$

then in the algebra $H$ the following equation is true

(9.4.11)
$$e_0 = -e_1 e_1 = -e_2 e_2$$
$$e_3 = e_1 e_2$$

Therefore, the set $(e_1, e_2)$ is a basis of algebra $H$. Ambiguity of representation of quaternion relative to the given basis is consequence of the equation (9.4.11). Namely, we can present a quaternion $a \in H$ as

$$a = (a^0 - a^4)e_1 e_1 + a^4 e_2 e_2 + a^1 e_1 + a^2 e_2 + a^3 e_1 e_2$$

where $a^4$ is arbitrary.

## 9.5. Left Module over Algebra

DEFINITION 9.5.1. *Effective left-side representation*

(9.5.1)         $f : A \dashrightarrow V \quad f(a) : v \in V \to av \in V \quad a \in A$

*of associative $D$-algebra $A$ in $D$-module $V$ is called* **left module** *over $D$-algebra $A$. We will also say that $D$-module $V$ is* **left $A$-module** *or $A*$-module*. *$V$-number is called* **vector**.                                  □

DEFINITION 9.5.2. *Let $A$ be division algebra. Effective left-side representation*

$$f : A \dashrightarrow V \quad f(a) : v \in V \to av \in V \quad a \in A$$

*of Abelian group $A$ in $D$-module $V$ is called* **left vector space** *over $D$-algebra $A$. We will also say that $D$-module $V$ is* **left $A$-vector space** *or $A*$-vector space*. *$V$-number is called* **vector**.                                  □

Theorem 9.5.3. *The following diagram of representations describes left $A$-module $V$*

$$(9.5.2)$$

$$
\begin{array}{c}
A \xrightarrow{g_{23}} A \xrightarrow{g_{3,4}} V \\
\end{array}
$$

$g_{12}(d) : a \to d\,a$

$g_{23}(v) : w \to C(w,v)$

$\qquad C \in \mathcal{L}(A^2 \to A)$

$g_{3,4}(a) : v \to v\,a$

$g_{1,4}(d) : v \to d\,v$

*The diagram of representations* (9.5.2) *holds* **commutativity of representations** *of commutative ring $D$ and $D$-algebra $A$ in Abelian group $V$*

$$(9.5.3) \qquad a(dv) = d(av)$$

Proof. The diagram of representations (9.5.2) follows from the definition 9.5.1 and the theorem 9.4.2. Since left-side transformation $g_{3,4}(a)$ is endomorphism of $D$-module $V$, we obtain the equality (9.5.3).  $\square$

Theorem 9.5.4. *Let $g$ be effective left-side representation of $D$-algebra $A$ in $D$-module $V$. Then $D$-algebra $A$ is associative.*

Proof. Let $a, b, c \in A$, $v \in V$. Since left-side representation $g$ is left-side representation of the multiplicative group of $D$-algebra $A$, we obtain the equality

$$(9.5.4) \qquad (ab)v = a(bv)$$

The equality

$$(9.5.5) \qquad a(b(cv)) = a((bc)v) = (a(bc))v$$

follows from the equality (9.5.4). Since $cv \in A$, the equality

$$(9.5.6) \qquad a(b(cv)) = (ab)(cv) = ((ab)c)v$$

follows from the equality (9.5.4). The equality

$$(9.5.7) \qquad (a(bc))v = ((ab)c)v$$

follows from equalities (9.5.5), (9.5.7). Since $v$ is any vector of $A$-module $V$, the equality

$$(9.5.8) \qquad a(bc) = (ab)c$$

follows from the equality (9.5.7). Therefore, $D$-algebra $A$ is associative.  $\square$

Theorem 9.5.5. *Let $V$ be left $A$-module. For any vector $v \in V$, vector generated by the diagram of representations* (9.5.2) *has the following form*

$$(9.5.9) \qquad (a+n)v = av + nv \quad a \in A \quad n \in D$$

9.5.5.1: *The set of maps*

$$(9.5.10) \qquad a + n : v \in V \to (a+n)v \in V$$

*generates[9.8] $D$-algebra $A_{(1)}$ where the sum is defined by the equality*

$$(9.5.11) \qquad (a+n) + (b+m) = (a+b) + (n+m)$$

---

[9.8] See the definition of unital extension also on the pages [6]-52, [7]-64.

*and the product is defined by the equality*

(9.5.12) $$(a + n)(b + m) = (ab + ma + nb) + (nm)$$

*The $D$-algebra $A_{(1)}$ is called* **unital extension** *of the $D$-algebra $A$.*

| | | |
|---|---|---|
| *If $D$-algebra $A$ has unit, then* | $D \subseteq A$ | $A_{(1)} = A$ |
| *If $D$-algebra $A$ is ideal of $D$, then* | $A \subseteq D$ | $A_{(1)} = D$ |
| *Otherwise* | | $A_{(1)} = A \oplus D$ |

9.5.5.2: *The $D$-algebra $A$ is left ideal of $D$-algebra $A_{(1)}$.*

9.5.5.3: *The set of transormations* (9.5.9) *is left-side representation of $D$-algebra $A_{(1)}$ in Abelian group $V$.*

*We use the notation* $A_{(1)}v$ *for the set of vectors generated by vector $v$.*

Theorem 9.5.6. *Following conditions hold for left $A$-module $V$:*

9.5.6.1: **associative law**

(9.5.13) $$(pq)v = p(qv)$$

9.5.6.2: **distributive law**

(9.5.14) $$p(v + w) = pv + pw$$

(9.5.15) $$(p + q)v = pv + qv$$

9.5.6.3: **unitarity law**

(9.5.16) $$1v = v$$

*for any $p, q \in A_{(1)}, v, w \in V$.*

Proof of theorems 9.5.5, 9.5.6. Let $v \in V$.

Lemma 9.5.7. *Let $d \in D$, $a \in A$. The map* (9.5.10) *is endomorphism of Abelian group $V$.*

Proof. Statements $dv \in V$, $av \in V$ follow from the theorems 6.1.4, 9.5.3. Since $V$ is Abelian group, then

$$dv + av \in V \quad d \in D \quad a \in A$$

Therefore, for any $D$-number $d$ and for any $A$-number $a$, we defined the map (9.5.10). Since transformation $g_{1,4}(d)$ and left-side transformation $g_{3,4}(a)$ are endomorphisms of Abelian group $V$, then the map (9.5.10) is endomorphism of Abelian group $V$. ⊙

Let $A_{(1)}$ be the set of maps (9.5.10). The equality (9.5.14) follows from the lemma 9.5.7.

Let $p = a + n \in A_{(1)}, q = b + m \in A_{(1)}$. According to the statement 9.3.3.3, we define the sum of $A_{(1)}$-numbers $p$ and $q$ by the equality (9.5.15). The equality

(9.5.17) $$((a + n) + (b + m))v = (a + n)v + (b + m)v$$

follows from the equality (9.5.15). Since representation $g_{1,4}$ is homomorphism of the aditive group of ring $D$, we obtain the equality

(9.5.18) $$(n + m)v = cn + dm$$

Since left-side representation $g_{3,4}$ is homomorphism of the aditive group of $D$-algebra $A$, we obtain the equality

$$(9.5.19) \qquad (a+b)v = av + bv$$

Since $V$ is Abelian group, then the equality

$$(9.5.20) \quad ((a+n)+(b+m))v = av + nv + bv + mv = av + bv + nv + mv$$
$$= (a+b)v + (n+m)v = ((a+b)+(n+m))v$$

follows from equalities (9.5.17), (9.5.18), (9.5.19). From the equality (9.5.20), it follows that the definition (9.5.11) of sum on the set $A_{(1)}$ does not depend on vector $v$.

Equalities (9.5.13), (9.5.16) follow from the statement 9.5.5.3. Let $p = a + n \in A_{(1)}$, $q = b + m \in A_{(1)}$. Since the product in $D$-algebra $A$ can be non associative, then, based on the theorem 9.5.6, we consider product of $A_{(1)}$-numbers $p$ and $q$ as bilinear map

$$f : A_{(1)} \times A_{(1)} \to A_{(1)}$$

such that following equalities are true

$$(9.5.21) \qquad f(a,b) = ab \quad a,b \in A$$

$$(9.5.22) \qquad f(1,p) = f(p,1) = p \quad p \in A_{(1)} \quad 1 \in D_{(1)}$$

The equality

$$(9.5.23) \quad \begin{aligned} (a+n)(b+m) &= f(a+n, b+m) \\ &= f(a,b) + f(a,m) + f(n,b) + f(n,m) \\ &= f(a,b) + mf(a,1) + nf(1,b) + nf(1,m) \\ &= ab + ma + nb + nm \end{aligned}$$

follows from equalities (9.5.21), (9.5.22). The equality (9.5.12) follows from the equality (9.5.23).

The statement 9.5.5.2 follows from the equality (9.5.12). $\qquad \square$

Bilinear map

$$(a,v) \in A \times V \to av \in V$$

generated by left-side representation $g_{2,3}$ is called **left-side product** of vector over scalar.

THEOREM 9.5.8. *Let $V$ be left $A$-module. The set of vectors generated by the set of vectors $v = (v_i \in V, i \in I)$ has form*[9.9]

$$(9.5.24) \qquad J(v) = \left\{ w : w = \sum_{i \in I} c^i v_i, c^i \in A_{(1)}, |\{i : c^i \neq 0\}| < \infty \right\}$$

PROOF. We prove the theorem by induction based on the theorem 6.1.4, According to the theorem 6.1.4, we need to prove following statements:

9.5.8.1: $v_k \in X_0 \subseteq J(v)$

9.5.8.2: $c^k v_k \in J(v)$, $c^k \in A_{(1)}$, $k \in I$

9.5.8.3: $\sum_{k \in I} c^k v_k \in J(v)$, $c^k \in A_{(1)}$, $|\{i : c^i \neq 0\}| < \infty$

---

[9.9] For a set $A$, we denote by $|A|$ the cardinal number of the set $A$. The notation $|A| < \infty$ means that the set $A$ is finite.

9.5.8.4: $w_1$, $w_2 \in J(v) \Rightarrow w_1 + w_2 \in J(v)$

9.5.8.5: $a \in A$, $w \in J(v) \Rightarrow aw \in J(v)$

- For any $v_k \in v$, let $c^i = \delta^i_k \in A_{(1)}$. Then

$$(9.5.25) \qquad v_k = \sum_{i \in I} c^i v_i$$

The statement 9.5.8.1 follows from (9.5.24), (9.5.25).

- The statement 9.5.8.2 follow from the theorems 6.1.4, 9.5.5 and from the statement 9.5.8.1.

- Since $V$ is Abelian group, then the statement 9.5.8.3 follows from the statement 9.5.8.2 and from theorems 6.1.4, 9.2.3.

- Let $w_1$, $w_2 \in X_k \subseteq J(v)$. Since $V$ is Abelian group, then, according to the statement 6.1.4.3,

$$(9.5.26) \qquad w_1 + w_2 \in X_{k+1}$$

According to the equality (9.5.24), there exist $A_{(1)}$-numbers $w^i_1$, $w^i_2$, $i \in I$, such that

$$(9.5.27) \qquad w_1 = \sum_{i \in I} w^i_1 v_i \quad w_2 = \sum_{i \in I} w^i_2 v_i$$

where sets

$$(9.5.28) \qquad H_1 = \{i \in I : w^i_1 \neq 0\} \quad H_2 = \{i \in I : w^i_2 \neq 0\}$$

are finite. Since $V$ is Abelian group, then from the equality (9.5.27) it follows that

$$(9.5.29) \qquad w_1 + w_2 = \sum_{i \in I} w^i_1 v_i + \sum_{i \in I} w^i_2 v_i = \sum_{i \in I} (w^i_1 v_i + w^i_2 v_i)$$

The equality

$$(9.5.30) \qquad w_1 + w_2 = \sum_{i \in I} (w^i_1 + w^i_2) v_i$$

follows from equalities (9.5.15), (9.5.29). From the equality (9.5.28), it follows that the set

$$\{i \in I : w^i_1 + w^i_2 \neq 0\} \subseteq H_1 \cup H_2$$

is finite.

- Let $w \in X_k \subseteq J(v)$. According to the statement 6.1.4.4, for any $A_{(1)}$-number $a$,

$$(9.5.31) \qquad aw \in X_{k+1}$$

According to the equality (9.5.24), there exist $A_{(1)}$-numbers $w^i$, $i \in I$, such that

$$(9.5.32) \qquad w = \sum_{i \in I} w^i v_i$$

where

$$(9.5.33) \qquad |\{i \in I : w^i \neq 0\}| < \infty$$

From the equality (9.5.32) it follows that

$$(9.5.34) \qquad aw = a \sum_{i \in I} w^i v_i = \sum_{i \in I} a(w^i v_i) = \sum_{i \in I} (aw^i) v_i$$

From the statement (9.5.33), it follows that the set $\{i \in I : aw^i \neq 0\}$ is finite.

From equalities (9.5.26), (9.5.30), (9.5.31), (9.5.34), it follows that $X_{k+1} \subseteq J(v)$.

$\square$

DEFINITION 9.5.9. *Let $v = (v_i \in V, i \in I)$ be set of vectors. The expression $w^i v_i$ is called* **linear combination** *of vectors $v_i$. A vector $\overline{w} = w^i v_i$ is called* **linearly dependent** *on vectors $v_i$.* $\square$

We represent the set of $A_{(1)}$-numbers $w^i$, $i \in I$, as matrix

$$w = \begin{pmatrix} w^I \\ ... \\ w^n \end{pmatrix}$$

We represent the set of vectors $v_i$, $i \in I$, as matrix

$$v = \begin{pmatrix} v_I & ... & v_n \end{pmatrix}$$

Then we can represent linear combination of vectors $\overline{w} = w^i v_i$ as

$$\overline{w} = w^*{}_* v$$

THEOREM 9.5.10. *Let $A$ be associative division $D$-algebra. Since the equation*

$$w^i v_i = 0$$

*implies existence of index $i = j$ such that $w^j \neq 0$, then the vector $v_j$ linearly depends on rest of vectors $v$.*

PROOF. The theorem follows from the equality

$$v_j = \sum_{i \in I \setminus \{j\}} (w^j)^{-1} w^i v_i$$

and from the definition 9.5.9. $\square$

It is evident that for any set of vectors $v_i$

$$w^i = 0 \Rightarrow w^*{}_* v = 0$$

DEFINITION 9.5.11. *The set of vectors[9.10] $v_i$, $i \in I$, of left $A$-module $V$ is* **linearly independent** *if $w = 0$ follows from the equation*

$$w^i v_i = 0$$

*Otherwise the set of vectors $v_i$, $i \in I$, is* **linearly dependent**. $\square$

The following definition follows from the theorems 9.5.8, 6.1.4 and from the definition 6.1.5.

DEFINITION 9.5.12. *$J(v)$ is called* **submodule generated by set** *$v$, and $v$ is a* **generating set** *of submodule $J(v)$. In particular, a* **generating set** *of left $D$-module $V$ is a subset $X \subset V$ such that $J(X) = V$.* $\square$

---

[9.10] I follow to the definition in [2], page 130.

The following definition follows from the theorems 9.5.8, 6.1.4 and from the definition 6.2.6.

**Definition 9.5.13.** *If the set* $X \subset V$ *is generating set of left D-module* $V$, *then any set* $Y$, $X \subset Y \subset V$ *also is generating set of left D-module* $V$. *If there exists minimal set* $X$ *generating the left D-module* $V$, *then the set* $X$ *is called* **basis of left D-module** $V$.                    □

**Theorem 9.5.14.** *The set of vectors* $\overline{\overline{e}} = (e_i, i \in I)$ *is basis of left A-module* $V$, *if following statements are true.*

9.5.14.1: *Arbitrary vector* $v \in V$ *is linear combination of vectors of the set* $\overline{\overline{e}}$.

9.5.14.2: *Vector* $e_i$ *cannot be represented as a linear combination of the remaining vectors of the set* $\overline{\overline{e}}$.

**Proof.** According to the statement 9.5.14.1, the theorem 9.5.8 and the definition 9.5.9, the set $\overline{\overline{e}}$ generates left A-module $V$ (the definition 9.5.12). According to the statement 9.5.14.2, the set $\overline{\overline{e}}$ is minimal set generating left A-module $V$. According to the definitions 9.5.13, the set $\overline{\overline{e}}$ is a basis of left A-module $V$.            □

**Theorem 9.5.15.** *Let A be associative division D-algebra. The set of vectors* $\overline{\overline{e}} = (e_i, i \in I)$ *is a* **basis of left A-vector space** $V$ *if vectors* $e_i$ *are linearly independent and any vector* $v \in V$ *linearly depends on vectors* $e_i$.

**Proof.** Let the set of vectors $e_i$, $i \in I$, be linear dependent. Then the equation

$$w^i e_i = 0$$

implies existence of index $i = j$ such that $w^j \neq 0$. According to the theorem 9.5.10, the vector $e_j$ linearly depends on rest of vectors of the set $\overline{\overline{e}}$. According to the definition 9.5.13, the set of vectors $e_i$, $i \in I$, is not a basis for left A-vector space $V$.

Therefore, if the set of vectors $e_i$, $i \in I$, is a basis, then these vectors are linearly independent. Since an arbitrary vector $v \in V$ is linear combination of vectors $e_i, i \in I$, , then the set of vectors $v$, $e_i, i \in I$, is not linearly independent.                    □

**Definition 9.5.16.** *Let* $\overline{\overline{e}}$ *be the basis of left A-module* $V$ *and vector* $\overline{v} \in V$ *has expansion*

$$\overline{v} = v^*{}_* e = v^i e_i$$

*with respect to the basis* $\overline{\overline{e}}$. $A_{(1)}$-*numbers* $v^i$ *are called* **coordinates** *of vector* $\overline{v}$ *with respect to the basis* $\overline{\overline{e}}$. *Matrix of* $A_{(1)}$-*numbers* $v = (v^i, i \in I)$ *is called* **coordinate matrix of vector** $\overline{v}$ **in basis** $\overline{\overline{e}}$.            □

**Theorem 9.5.17.** *Let A be associative D-algebra. Let* $\overline{\overline{e}}$ *be basis of left A-module* $V$. *Let*

$$(9.5.35) \qquad\qquad w^i e_i = 0$$

*be linear dependence of vectors of the basis* $\overline{\overline{e}}$. *Then*

9.5.17.1: $A_{(1)}$-*number* $w^i$, $i \in I$, *does not have inverse element in D-algebra* $A_{(1)}$.

9.5.17.2: *The set* $A'$ *of matrices* $w = (w^i, i \in I)$ *generates left A-module.*

Proof. Let there exist matrix $w = (w^i, i \in I)$ such that the equality (9.5.35) is true and there exist index $i = j$ such that $w^j \neq 0$. If we assume that $A_{(1)}$-number $c^j$ has inverse one, then the equality

$$e_j = \sum_{i \in I \setminus \{j\}} (w^j)^{-1} w^i e_i$$

follows from the equality (9.5.35). Therefore, the vector $e_j$ is linear combination of other vectors of the set $\overline{\overline{e}}$ and the set $\overline{\overline{e}}$ is not basis. Therefore, our assumption is false, and $A_{(1)}$-number $c^j$ does not have inverse.

Let matrices $b = (b^i, i \in I) \in A'$, $c = (c^i, i \in I) \in A'$. From equalities

$$b^i e_i = 0$$

$$c^i e_i = 0$$

it follows that

$$(b^i + c^i) e_i = 0$$

Therefore, the set $A'$ is Abelian group.

Let matrix $c = (c^i, i \in I) \in A'$ and $a \in A$. From the equality

$$w^i e_i = 0$$

it follows that

$$(ac^i) e_i = 0$$

Therefore, Abelian group $A'$ is left $A$-module. $\qquad \square$

Theorem 9.5.18. *Let left $A$-module $V$ have the basis $\overline{\overline{e}}$ such that in the equality*

$$(9.5.36) \qquad\qquad w^i e_i = 0$$

*there exists index $i = j$ such that $w^j \neq 0$. Then*

9.5.18.1: *The matrix $w = (w^i, i \in I)$ determines coordinates of vector $0 \in V$ with respect to basis $\overline{\overline{e}}$.*

9.5.18.2: *Coordinates of vector $\overline{v}$ with respect to basis $\overline{\overline{e}}$ are uniquely determined up to a choice of coordinates of vector $0 \in V$.*

Proof. The statement 9.5.18.1 follows from the equality (9.5.36) and from the definition 9.5.16.

Let vector $\overline{v}$ have expansion

$$(9.5.37) \qquad\qquad \overline{v} = v^*_* e = v^i e_i$$

with respect to basis $\overline{\overline{e}}$. The equality

$$(9.5.38) \qquad\qquad \overline{v} = \overline{v} + 0 = v^i e_i + c^i e_i = (v^i + c^i) e_i$$

follows from equalities (9.5.36), (9.5.37). The statement 9.5.18.2 follows from equalities (9.5.37), (9.5.38) and from the definition 9.5.16. $\qquad \square$

Definition 9.5.19. *The left $A$-module $V$ is* **free left $A$-module,** [9.11] *if left $A$-module $V$ has basis and vectors of the basis are linearly independent.* $\qquad \square$

Theorem 9.5.20. *Coordinates of vector $v \in V$ relative to basis $\overline{\overline{e}}$ of free left $A$-module $V$ are uniquely defined.*

Proof. The theorem follows from the theorem 9.5.18 and from definitions 9.5.11, 9.5.19. $\qquad \square$

---

[9.11] I follow to the definition in [2], page 135.

## 9.6. Right Module over Algebra

Definition 9.6.1. *Effective right-side representation*

$$(9.6.1) \qquad f : A \rightarrowtail V \quad f(a) : v \in V \to va \in V \quad a \in A$$

*of associative D-algebra A in D-module V is called* **right module** *over D-algebra A. We will also say that D-module V is* **right A-module** *or* **∗A-module**. *V-number is called* **vector**. $\qquad\square$

Definition 9.6.2. *Let A be division algebra. Effective right-side representation*

$$f : A \rightarrowtail V \quad f(a) : v \in V \to va \in V \quad a \in A$$

*of Abelian group A in D-module V is called* **right vector space** *over D-algebra A. We will also say that D-module V is* **right A-vector space** *or* **∗A-vector space**. *V-number is called* **vector**. $\qquad\square$

Theorem 9.6.3. *The following diagram of representations describes right A-module V*

$$(9.6.2)$$

$$g_{12}(d) : a \to d\,a$$
$$g_{23}(v) : w \to C(w, v)$$
$$C \in \mathcal{L}(A^2 \to A)$$
$$g_{3,4}(a) : v \to v\,a$$
$$g_{1,4}(d) : v \to v\,d$$

*The diagram of representations* (9.6.2) *holds* **commutativity of representations** *of commutative ring D and D-algebra A in Abelian group V*

$$(9.6.3) \qquad (vd)a = (va)d$$

Proof. The diagram of representations (9.6.2) follows from the definition 9.6.1 and the theorem 9.4.2. Since right-side transformation $g_{3,4}(a)$ is endomorphism of D-module V, we obtain the equality (9.6.3). $\qquad\square$

Theorem 9.6.4. *Let g be effective left-side representation of D-algebra A in D-module V. Then D-algebra A is associative.*

Proof. Let $a, b, c \in A$, $v \in V$. Since right-side representation g is right-side representation of the multiplicative group of D-algebra A, we obtain the equality

$$(9.6.4) \qquad v(ab) = (va)b$$

The equality

$$(9.6.5) \qquad ((vc)b)a = (v(cb))a = v((cb)a)$$

follows from the equality (9.6.4). Since $vc \in A$, the equality

$$(9.6.6) \qquad ((vc)b)a = (vc)(ba) = v(c(ba))$$

follows from the equality (9.6.4). The equality

$$(9.6.7) \qquad v((cb)a) = v(c(ba))$$

follows from equalities (9.6.5), (9.6.7). Since $v$ is any vector of $A$-module $V$, the equality

$$(9.6.8) \qquad (cb)a = c(ba)$$

follows from the equality (9.6.7). Therefore, $D$-algebra $A$ is associative.     $\square$

THEOREM 9.6.5. *Let $V$ be right $A$-module. For any vector $v \in V$, vector generated by the diagram of representations* (9.6.2) *has the following form*

$$(9.6.9) \qquad v(a + n) = va + vn \quad a \in A \quad n \in D$$

9.6.5.1: *The set of maps*

$$(9.6.10) \qquad a + n : v \in V \to v(a + n) \in V$$

*generates*[9.12] *$D$-algebra $A_{(1)}$ where the sum is defined by the equality*

$$(9.6.11) \qquad (a + n) + (b + m) = (a + b) + (n + m)$$

*and the product is defined by the equality*

$$(9.6.12) \qquad (a + n)(b + m) = (ab + ma + nb) + (nm)$$

*The $D$-algebra $A_{(1)}$ is called* **unital extension** *of the $D$-algebra $A$.*

| If $D$-algebra $A$ has unit, then | $D \subseteq A$ | $A_{(1)} = A$ |
| --- | --- | --- |
| If $D$-algebra $A$ is ideal of $D$, then | $A \subseteq D$ | $A_{(1)} = D$ |
| Otherwise | $A_{(1)} = A \oplus D$ | |

9.6.5.2: *The $D$-algebra $A$ is right ideal of $D$-algebra $A_{(1)}$.*

9.6.5.3: *The set of transormations* (9.6.9) *is right-side representation of $D$-algebra $A_{(1)}$ in Abelian group $V$.*

*We use the notation  $A_{(1)}v$  for the set of vectors generated by vector $v$.*

THEOREM 9.6.6. *Following conditions hold for right $A$-module $V$:*

9.6.6.1: **associative law**

$$(9.6.13) \qquad v(pq) = (vp)q$$

9.6.6.2: **distributive law**

$$(9.6.14) \qquad (v + w)p = vp + wp$$

$$(9.6.15) \qquad v(p + q) = vp + vq$$

9.6.6.3: **unitarity law**

$$(9.6.16) \qquad v1 = v$$

*for any   $p, q \in A_{(1)}, v, w \in V$.*

PROOF OF THEOREMS 9.6.5, 9.6.6. Let $v \in V$.

LEMMA 9.6.7. *Let  $d \in D$,    $a \in A$.    The map* (9.6.10) *is endomorphism of Abelian group $V$.*

---

[9.12] See the definition of unital extension also on the pages [6]-52, [7]-64.

Proof. Statements $vd \in V$, $va \in V$ follow from the theorems 6.1.4, 9.6.3. Since $V$ is Abelian group, then

$$vd + va \in V \quad d \in D \quad a \in A$$

Therefore, for any $D$-number $d$ and for any $A$-number $a$, we defined the map (9.6.10). Since transformation $g_{1,4}(d)$ and right-side transformation $g_{3,4}(a)$ are endomorphisms of Abelian group $V$, then the map (9.6.10) is endomorphism of Abelian group $V$.                                                                   ⊙

Let $A_{(1)}$ be the set of maps (9.6.10). The equality (9.6.14) follows from the lemma 9.6.7.

Let $p = a + n \in A_{(1)}$, $q = b + m \in A_{(1)}$. According to the statement 9.3.3.3, we define the sum of $A_{(1)}$-numbers $p$ and $q$ by the equality (9.6.15). The equality

$$(9.6.17) \qquad v((a+n) + (b+m)) = v(a+n) + v(b+m)$$

follows from the equality (9.6.15). Since representation $g_{1,4}$ is homomorphism of the aditive group of ring $D$, we obtain the equality

$$(9.6.18) \qquad v(n+m) = vn + vm$$

Since right-side representation $g_{3,4}$ is homomorphism of the aditive group of $D$-algebra $A$, we obtain the equality

$$(9.6.19) \qquad v(a+b) = va + vb$$

Since $V$ is Abelian group, then the equality

$$(9.6.20) \qquad \begin{aligned} v((a+n) + (b+m)) &= va + vn + vb + vm = va + vb + vn + vm \\ &= v(a+b) + v(n+m) = v((a+b) + (n+m)) \end{aligned}$$

follows from equalities (9.6.17), (9.6.18), (9.6.19). From the equality (9.6.20), it follows that the definition (9.6.11) of sum on the set $A_{(1)}$ does not depend on vector $v$.

Equalities (9.6.13), (9.6.16) follow from the statement 9.6.5.3. Let $p = a + n \in A_{(1)}$, $q = b + m \in A_{(1)}$. Since the product in $D$-algebra $A$ can be non associative, then, based on the theorem 9.6.6, we consider product of $A_{(1)}$-numbers $p$ and $q$ as bilinear map

$$f : A_{(1)} \times A_{(1)} \to A_{(1)}$$

such that following equalities are true

$$(9.6.21) \qquad f(a,b) = ab \quad a,b \in A$$

$$(9.6.22) \qquad f(1,p) = f(p,1) = p \quad p \in A_{(1)} \quad 1 \in D_{(1)}$$

The equality

$$(9.6.23) \qquad \begin{aligned} (a+n)(b+m) &= f(a+n, b+m) \\ &= f(a,b) + f(a,m) + f(n,b) + f(n,m) \\ &= f(a,b) + mf(a,1) + nf(1,b) + nf(1,m) \\ &= ab + ma + nb + nm \end{aligned}$$

follows from equalities (9.6.21), (9.6.22). The equality (9.6.12) follows from the equality (9.6.23).

The statement 9.6.5.2 follows from the equality (9.6.12).                                      □

Bilinear map

$$(v,a) \in V \times A \to va \in V$$

generated by right-side representation $g_{2,3}$ is called **right-side product** of vector over scalar.

THEOREM 9.6.8. *Let $V$ be right $A$-module. The set of vectors generated by the set of vectors $v = (v_i \in V, i \in I)$ has form*[9.13]

$$(9.6.24) \qquad J(v) = \left\{ w : w = \sum_{i \in I} v_i c^i, c^i \in A_{(1)}, |\{i : c^i \neq 0\}| < \infty \right\}$$

PROOF. We prove the theorem by induction based on the theorem 6.1.4, According to the theorem 6.1.4, we need to prove following statements:

9.6.8.1: $v_k \in X_0 \subseteq J(v)$

9.6.8.2: $v_k c^k \in J(v)$, $c^k \in A_{(1)}$, $k \in I$

9.6.8.3: $\sum_{k \in I} v_k c^k \in J(v)$, $c^k \in A_{(1)}$, $|\{i : c^i \neq 0\}| < \infty$

9.6.8.4: $w_1, w_2 \in J(v) \Rightarrow w_1 + w_2 \in J(v)$

9.6.8.5: $a \in A$, $w \in J(v) \Rightarrow aw \in J(v)$

- For any $v_k \in v$, let $c^i = \delta^i_k \in A_{(1)}$. Then

$$(9.6.25) \qquad v_k = \sum_{i \in I} v_i c^i$$

  The statement 9.6.8.1 follows from (9.6.24), (9.6.25).
- The statement 9.6.8.2 follow from the theorems 6.1.4, 9.6.5 and from the statement 9.6.8.1.
- Since $V$ is Abelian group, then the statement 9.6.8.3 follows from the statement 9.6.8.2 and from theorems 6.1.4, 9.2.3.
- Let $w_1, w_2 \in X_k \subseteq J(v)$. Since $V$ is Abelian group, then, according to the statement 6.1.4.3,

$$(9.6.26) \qquad w_1 + w_2 \in X_{k+1}$$

  According to the equality (9.6.24), there exist $A_{(1)}$-numbers $w_1^i, w_2^i, i \in I$, such that

$$(9.6.27) \qquad w_1 = \sum_{i \in I} v_i w_1^i \qquad w_2 = \sum_{i \in I} v_i w_2^i$$

  where sets

$$(9.6.28) \qquad H_1 = \{i \in I : w_1^i \neq 0\} \qquad H_2 = \{i \in I : w_2^i \neq 0\}$$

  are finite. Since $V$ is Abelian group, then from the equality (9.6.27) it follows that

$$(9.6.29) \qquad w_1 + w_2 = \sum_{i \in I} v_i w_1^i + \sum_{i \in I} v_i w_2^i = \sum_{i \in I} (v_i w_1^i + v_i w_2^i)$$

  The equality

$$(9.6.30) \qquad w_1 + w_2 = \sum_{i \in I} v_i (w_1^i + w_2^i)$$

---

[9.13] For a set $A$, we denote by $|A|$ the cardinal number of the set $A$. The notation $|A| < \infty$ means that the set $A$ is finite.

follows from equalities (9.6.15), (9.6.29). From the equality (9.6.28), it follows that the set

$$\{i \in I : w_1^i + w_2^i \neq 0\} \subseteq H_1 \cup H_2$$

is finite.

- Let $w \in X_k \subseteq J(v)$. According to the statement 6.1.4.4, for any $A_{(1)}$-number $a$,

$$(9.6.31) \qquad\qquad wa \in X_{k+1}$$

According to the equality (9.6.24), there exist $A_{(1)}$-numbers $w^i$, $i \in I$, such that

$$(9.6.32) \qquad\qquad w = \sum_{i \in I} v_i w^i$$

where

$$(9.6.33) \qquad\qquad |\{i \in I : w^i \neq 0\}| < \infty$$

From the equality (9.6.32) it follows that

$$(9.6.34) \qquad wa = \left(\sum_{i \in I} v_i w^i\right) a = \sum_{i \in I} (v_i w^i) a = \sum_{i \in I} (v_i w^i a)$$

From the statement (9.6.33), it follows that the set $\{i \in I : w^i a \neq 0\}$ is finite.

From equalities (9.6.26), (9.6.30), (9.6.31), (9.6.34), it follows that $X_{k+1} \subseteq J(v)$.
□

DEFINITION 9.6.9. *Let $v = (v_i \in V, i \in I)$ be set of vectors. The expression $v_i w^i$ is called* **linear combination** *of vectors $v_i$. A vector $\overline{w} = v_i w^i$ is called* **linearly dependent** *on vectors $v_i$.* □

We represent the set of $A_{(1)}$-numbers $w^i$, $i \in I$, as matrix

$$w = \begin{pmatrix} w^I \\ ... \\ w^n \end{pmatrix}$$

We represent the set of vectors $v_i$, $i \in I$, as matrix

$$v = \begin{pmatrix} v_I & ... & v_n \end{pmatrix}$$

Then we can represent linear combination of vectors $\overline{w} = v_i w^i$ as

$$\overline{w} = v_*{}^* w$$

THEOREM 9.6.10. *Let $A$ be associative division $D$-algebra. Since the equation*

$$v_i w^i = 0$$

*implies existence of index $i = j$ such that $w^j \neq 0$, then the vector $v_j$ linearly depends on rest of vectors $v$.*

Proof. The theorem follows from the equality

$$v_j = \sum_{i \in I \setminus \{j\}} v_i w^i (w^j)^{-1}$$

and from the definition 9.6.9. $\qquad \square$

It is evident that for any set of vectors $v_i$

$$w^i = 0 \Rightarrow v_* {}^* w = 0$$

Definition 9.6.11. *The set of vectors*[9.14] $v_i$, $i \in I$, *of right $A$-module $V$ is* **linearly independent** *if $w = 0$ follows from the equation*

$$v_i w^i = 0$$

*Otherwise the set of vectors $v_i$, $i \in I$, is* **linearly dependent**. $\qquad \square$

The following definition follows from the theorems 9.6.8, 6.1.4 and from the definition 6.1.5.

Definition 9.6.12. *$J(v)$ is called* **submodule generated by set** *$v$, and $v$ is a* **generating set** *of submodule $J(v)$. In particular, a* **generating set** *of right $D$-module $V$ is a subset $X \subset V$ such that $J(X) = V$.* $\qquad \square$

The following definition follows from the theorems 9.6.8, 6.1.4 and from the definition 6.2.6.

Definition 9.6.13. *If the set $X \subset V$ is generating set of right $D$-module $V$, then any set $Y$, $X \subset Y \subset V$ also is generating set of right $D$-module $V$. If there exists minimal set $X$ generating the right $D$-module $V$, then the set $X$ is called* **basis** *of right $D$-module $V$.* $\qquad \square$

Theorem 9.6.14. *The set of vectors $\overline{\overline{e}} = (e_i, i \in I)$ is basis of right $A$-module $V$, if following statements are true.*

9.6.14.1: *Arbitrary vector $v \in V$ is linear combination of vectors of the set $\overline{\overline{e}}$.*

9.6.14.2: *Vector $e_i$ cannot be represented as a linear combination of the remaining vectors of the set $\overline{\overline{e}}$.*

Proof. According to the statement 9.6.14.1, the theorem 9.6.8 and the definition 9.6.9, the set $\overline{\overline{e}}$ generates right $A$-module $V$ (the definition 9.6.12). According to the statement 9.6.14.2, the set $\overline{\overline{e}}$ is minimal set generating right $A$-module $V$. According to the definitions 9.6.13, the set $\overline{\overline{e}}$ is a basis of right $A$-module $V$. $\qquad \square$

Theorem 9.6.15. *Let $A$ be associative division $D$-algebra. The set of vectors $\overline{\overline{e}} = (e_i, i \in I)$ is a* **basis of right $A$-vector space** *$V$ if vectors $e_i$ are linearly independent and any vector $v \in V$ linearly depends on vectors $e_i$.*

Proof. Let the set of vectors $e_i$, $i \in I$, be linear dependent. Then the equation

$$e_i w^i = 0$$

implies existence of index $i = j$ such that $w^j \neq 0$. According to the theorem 9.6.10, the vector $e_j$ linearly depends on rest of vectors of the set $\overline{\overline{e}}$. According to the definition 9.6.13, the set of vectors $e_i$, $i \in I$, is not a basis for right $A$-vector space $V$.

---

[9.14] I follow to the definition in [2], page 130.

Therefore, if the set of vectors $e_i,\ i \in I,$ is a basis, then these vectors are linearly independent. Since an arbitrary vector $v \in V$ is linear combination of vectors $e_i, i \in I,$, then the set of vectors $v,\ e_i, i \in I,$ is not linearly independent. $\square$

DEFINITION 9.6.16. *Let $\overline{\overline{e}}$ be the basis of right $A$-module $V$ and vector $\overline{v} \in V$ has expansion*

$$\overline{v} = e_* {}^* v = e_i v^i$$

*with respect to the basis $\overline{\overline{e}}$. $A_{(1)}$-numbers $v^i$ are called **coordinates** of vector $\overline{v}$ with respect to the basis $\overline{\overline{e}}$. Matrix of $A_{(1)}$-numbers $v = (v^i, i \in I)$ is called **coordinate matrix of vector** $\overline{v}$ in basis $\overline{\overline{e}}$.* $\square$

THEOREM 9.6.17. *Let $A$ be associative $D$-algebra. Let $\overline{\overline{e}}$ be basis of right $A$-module $V$. Let*

$$(9.6.35) \qquad\qquad e_i w^i = 0$$

*be linear dependence of vectors of the basis $\overline{\overline{e}}$. Then*

9.6.17.1: *$A_{(1)}$-number $w^i,\ i \in I,$ does not have inverse element in $D$-algebra $A_{(1)}$.*

9.6.17.2: *The set $A'$ of matrices $w = (w^i, i \in I)$ generates right $A$-module.*

PROOF. Let there exist matrix $w = (w^i, i \in I)$ such that the equality (9.6.35) is true and there exist index $i = j$ such that $w^j \neq 0$. If we assume that $A_{(1)}$-number $c^j$ has inverse one, then the equality

$$e_j = \sum_{i \in I \setminus \{j\}} e_i w^i (w^j)^{-1}$$

follows from the equality (9.6.35). Therefore, the vector $e_j$ is linear combination of other vectors of the set $\overline{\overline{e}}$ and the set $\overline{\overline{e}}$ is not basis. Therefore, our assumption is false, and $A_{(1)}$-number $c^j$ does not have inverse.

Let matrices $b = (b^i, i \in I) \in A',\quad c = (c^i, i \in I) \in A'$. From equalities

$$e_i b^i = 0$$
$$e_i c^i = 0$$

it follows that

$$e_i (b^i + c^i) = 0$$

Therefore, the set $A'$ is Abelian group.

Let matrix $c = (c^i, i \in I) \in A'$ and $a \in A$. From the equality

$$e_i w^i = 0$$

it follows that

$$e_i (c^i a) = 0$$

Therefore, Abelian group $A'$ is right $A$-module. $\square$

THEOREM 9.6.18. *Let right $A$-module $V$ have the basis $\overline{\overline{e}}$ such that in the equality*

$$(9.6.36) \qquad\qquad e_i w^i = 0$$

*there exists index $i = j$ such that $w^j \neq 0$. Then*

9.6.18.1: *The matrix $w = (w^i, i \in I)$ determines coordinates of vector $0 \in V$ with respect to basis $\overline{\overline{e}}$.*

9.6.18.2: *Coordinates of vector $\overline{v}$ with respect to basis $\overline{\overline{e}}$ are uniquely determined up to a choice of coordinates of vector $0 \in V$.*

PROOF. The statement 9.6.18.1 follows from the equality (9.6.36) and from the definition 9.6.16.

Let vector $\overline{v}$ have expansion

$$(9.6.37) \qquad \overline{v} = e_*{}^* v = e_i v^i$$

with respect to basis $\overline{\overline{e}}$. The equality

$$(9.6.38) \qquad \overline{v} = \overline{v} + 0 = e_i v^i + e_i c^i = e_i (v^i + c^i)$$

follows from equalities (9.6.36), (9.6.37). The statement 9.6.18.2 follows from equalities (9.6.37), (9.6.38) and from the definition 9.6.16. $\qquad \square$

DEFINITION 9.6.19. *The right $A$-module $V$ is **free right** $A$-**module**,[9.15] if right $A$-module $V$ has basis and vectors of the basis are linearly independent.* $\qquad \square$

THEOREM 9.6.20. *Coordinates of vector $v \in V$ relative to basis $\overline{\overline{e}}$ of free right $A$-module $V$ are uniquely defined.*

PROOF. The theorem follows from the theorem 9.6.18 and from definitions 9.6.11, 9.6.19. $\qquad \square$

## 9.7. Left Module over Nonassociative Algebra

Theorems 9.6.5, 9.6.6 consider the structure of module over associative $D$-algebra $A$. It is easy to see that, considering some corrections, these theorems remain true if $A$ is non associative $D$-algebra. However, because the product in $D$-algebra $A$ is non associative and product of transformations in module over $D$-algebra $A$ is associative, then the map $g_{34}$ cannot be a representation of non associative $D$-algebra $A$.

We have come to that verge where universal algebra representation theory is defined. In order to maintain the ability to use the tool considered in this book, we can agree that the map

$$g_{34} : A \times V \to V$$

is a representation when the map $g_{34}$ is bilinear map. There are new questions that are beyond the scope of this book.

However, we may consider this problem from other point of view. If the map $g_{34}$ does not conserve the operation of the product, then we assume that the map $g_{34}$ is representation of $D$-algebra $A$, in which product is not defined. In other words, the map $g_{34}$ is representation of $D$-module. Therefore, diagram of representations will have the following form

$$(9.7.1)$$

$$g_{12}(d) : a \to d\,a$$
$$g_{34}(a) : v \to a\,v$$
$$g_{14}(d) : v \to d\,v$$

---

[9.15] I follow to the definition in [2], page 135.

However, we lost the structure of $D$-algebra $A$ in diagram of representations (9.7.1). Therefore, proper diagram of representations will have the following form

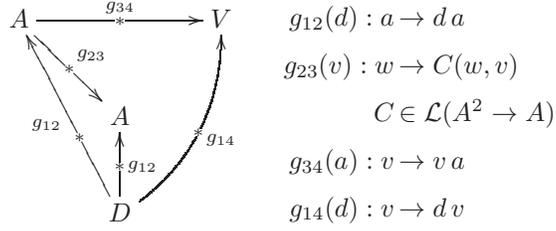

$$g_{12}(d) : a \to d\,a$$

$$g_{23}(v) : w \to C(w, v)$$

$$C \in \mathcal{L}(A^2 \to A)$$

$$g_{34}(a) : v \to v\,a$$

$$g_{14}(d) : v \to d\,v$$

# Examples of Diagram of Representations: Affinne Geometry

## 10.1. About this Chapter

In the chapter 9, we considered examples of diagram of representations associated with module over ring. If representation theory were reduced to studying of modules, it would hardly be an interesting theory.

In this chapter, I considered examples of diagram of representations associated with affinne geometry. This simple algebraic construction turned out to be a rich source of inspiration for me. I met interesting ideas in this area of mathematics twice. At first during study of affine geometry, I discovered that I can describe affine geometry using tower of representations. Afterwards, during similar study of algebra over commutative ring, I began to study a diagram of representations.

However, the second discovery came to me by chance. When I was looking through the calculus textbook, I have met a definition familiar from childhood. This is sum of vectors. The definition is extremely simple. When we define manifold with affine connection, we have the sum of vectors in tangent plane. However, at this time I realized that I can define sum of vectors using parallelogram from geodesic lines. It gave me ability to build affine geometry on affine manifold.

One more step, and I switched from manifold with affine connection to metric -affine manifold. Since parallelogram from geodesic lines is not closed, then sum of vectors in metric affine manifold is not commutative. Without a doubt, this is research which is beyond the scope of this book; and I hope to return to this research in the future. However, I decided to write a sketch of this theory in the section 10.4 to show the reader the limits of theory presented in this book.

The representation theory is natural extension of universal algebra theory. We assume that binary operation on universal algebra $A$ is defined for any two $A$-numbers. However, it is evident that sum of vectors in affine geometry On differentiable manifold is well defined only in enough small neighborhood.

I met similar problem in the paper [11] where I and Alexandre Laugier studied orthogonal transformations in Minkowski space. We discovered that the product of orthogonal transformations not always is orthogonal transformation; therefore, the set of orthogonal transformations is not a group.

## 10.2. Representation of Group on the Set

Let $G$ be Abelian group, and $M$ be a set. Consider effective representation of group $G$ on the set $M$. For given $a \in G$, $A \in M$ we assume $A \to A + a$. We also use notation $a = \vec{AB}$ if

$$(10.2.1) \qquad\qquad B = A + a$$

Then we can represent action of group as

(10.2.2) $$B = A + \overrightarrow{AB}$$

Since the representation is effective, then from equalities (10.2.1), (10.2.2) and the equality

$$D = C + a$$

it folows that

(10.2.3) $$\overrightarrow{AB} = \overrightarrow{CD}$$

$G$-number $a$ and corresponding transformation $\overrightarrow{AB}$ are called vector. We interpret the equality (10.2.3) as the independence of the vector $a$ from the choice of $M$-number $A$.

We can consider the set $M$ as union of orbits of the representation of the group $G$. We can select for basis of the representation the set of points such that one and only one point belongs to each orbit. If $X$ is the basis of representation, $A \in X$, $g \in G$, then $\Omega_2$-word has form $A+g$. Since there is no operations on the set $M$, then there is no $\Omega_2$-word containing different elements of the basis. If representation of group $G$ is single transitive, then basis of representation consists of one point. Any point of the set $M$ can be such point.

THEOREM 10.2.1. *Let the representation $A \to A + a$ of Abelian group $G$ on the set $M$ be single transitive. Then for any $M$-numbers $A$, $B$, $C$, we determine sum of vectors $\overrightarrow{AB}$ and $\overrightarrow{BC}$ and sum of vectors satisfies to the following equality*

(10.2.4) $$\overrightarrow{AB} + \overrightarrow{BC} = \overrightarrow{AC}$$

PROOF. Since the representation is single transitive, then, for any $M$-numbers $A$, $B$, $C$, there exist vectors $\overrightarrow{AB}$, $\overrightarrow{BC}$ such that

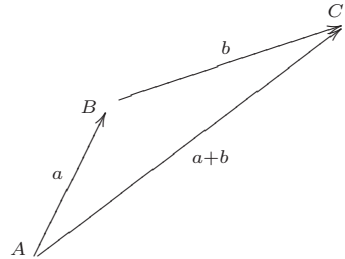

(10.2.5) $$B = A + \overrightarrow{AB}$$

(10.2.6) $$C = B + \overrightarrow{BC}$$

The equality

(10.2.7) $$C = (A + \overrightarrow{AB}) + \overrightarrow{BC} = A + (\overrightarrow{AB} + \overrightarrow{BC})$$

follows from equalities (10.2.5), (10.2.6) and from associativity of sum in Abelian group $G$. Since the representation is single transitive, then the equality (10.2.4) follows from the equality (10.2.7) and from the equality

$$C = A + \overrightarrow{AC}$$

This definition of sum is called the triangle law.                    □

REMARK 10.2.2. *Since $G$ is Abelian group, then statements 10.2.2.1, 10.2.2.2, follow from the theorem 10.2.1.*

10.2.2.1: $\overrightarrow{AA} = 0$

10.2.2.2: $\overrightarrow{AB} = -\overrightarrow{BA}$

10.2.2.3: *Addition is commutative.*

10.2.2.4: *Addition is associative.*

$\square$

THEOREM 10.2.3. *For given $a$, $b \in G$ and $A \in M$, we consider following set of M-numbers.*

- $B = A + a$
- $C = B + b$
- $D = A + b$
- $E = D + a$

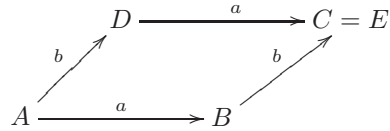

PROOF. The theorem follows from the statement 10.2.2.3.

$\square$

THEOREM 10.2.4. *If $\overrightarrow{AB} = \overrightarrow{CD}$, then $\overrightarrow{AC} = \overrightarrow{BD}$.*

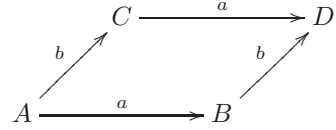

PROOF. Let $\overrightarrow{AB} = \overrightarrow{CD} = a$, $\overrightarrow{AC} = b$. According to the statement 10.2.2.2, $\overrightarrow{BA} = -a$. The theorem follows from the equality

$$D = B + \overrightarrow{BD} = B + \overrightarrow{BA} + \overrightarrow{AD} = B + \overrightarrow{BA} + \overrightarrow{AC} + \overrightarrow{CD}$$
$$= B - a + b + a = B + b$$

$\square$

## 10.3. Affine Space

DEFINITION 10.3.1. *Let $D$ be commutative ring and $V$ be free $D$-module. A set of points $\overset{\circ}{V}$ is called* **affine space** *over $D$-module $V$, if the set of points $\overset{\circ}{V}$ satisfies to following axioms.* [10.1]

10.3.1.1: *There exists at least one point*

10.3.1.2: *One and only one vector is in correspondence to any tuple of points $(A, B)$. We denote this vector as $\overrightarrow{AB}$. The vector $\overrightarrow{AB}$ has tail in the point $A$ and head in the point $B$.*

10.3.1.3: *For any point $A$ and any vector $a$ there exists one and only one point $B$ such, that $\overrightarrow{AB} = a$. We will use notation* [10.2]

(10.3.1) $$B = A + a$$

10.3.1.4: *(Axiom of parallelogram.) If $\overrightarrow{AB} = \overrightarrow{CD}$, then $\overrightarrow{AC} = \overrightarrow{BD}$.*

*A set $V$ is called a set of free vectors. $\overset{\circ}{V}$-number is called point of affine space $\overset{\circ}{V}$.*

$\square$

DEFINITION 10.3.2. *Let $A \in \overset{\circ}{V}$ be arbitrary point.*

---

[10.1]I wrote definitions and theorems in this section according to definition of affine space in [4], pp. 86 - 93.

[10.2][21], p. 9.

*Let $v$ be vector. According to the axiom 10.3.1.3, there exists $B \in \overset{\circ}{V}$, $B = A + v$.*

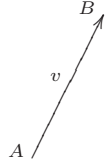

*Let $w$ be vector. According to the axiom 10.3.1.3, there exists $C \in \overset{\circ}{V}$, $C = B + w$.*

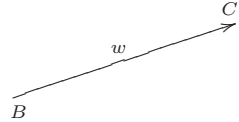

*According to the axiom 10.3.1.2, there exists vector $\overrightarrow{AC}$. Vector $\overrightarrow{AC}$ is called sum of vectors $v$ and $w$*

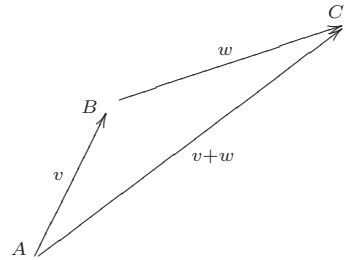

$$(10.3.2) \qquad v + w = \overrightarrow{AC}$$

*This definition of sum is called the triangle law.*

□

THEOREM 10.3.3. *Vector $\overrightarrow{AA}$ is zero with respect to addition and does not depend on point $A$. Vector $\overrightarrow{AA}$ is called zero-vector and we assume $\overrightarrow{AA} = 0$.*

PROOF. We can write rule of addition (10.3.2) in form of the equality

$$(10.3.3) \qquad \overrightarrow{AB} + \overrightarrow{BC} = \overrightarrow{AC}$$

If $B = C$, then from the equality (10.3.3) it follows that

$$(10.3.4) \qquad \overrightarrow{AB} + \overrightarrow{BB} = \overrightarrow{AB}$$

From the equality (10.3.4), it follows that the vector $\overrightarrow{BB}$ is zero with respect to addition. If $C = A$, $B = D$, then from axiom 10.3.1.4, it follows that $\overrightarrow{AA} = \overrightarrow{BB}$. Therefore, a zero-vector $\overrightarrow{AA}$ does not depend on a point $A$. □

THEOREM 10.3.4. *Let $a = \overrightarrow{AB}$. Then*

$$(10.3.5) \qquad \overrightarrow{BA} = -a$$

*and this equality does not depend on point $A$.*

PROOF. From the equality (10.3.3) and the theorem 10.3.3, it follows that

$$(10.3.6) \qquad \overrightarrow{AB} + \overrightarrow{BA} = \overrightarrow{AA} = 0$$

The equality (10.3.5) follows from the equality (10.3.6). Applying axiom 10.3.1.4 to the equality $\overrightarrow{AB} = \overrightarrow{CD}$ we get $\overrightarrow{AC} = \overrightarrow{BD}$, or (this is equivalent)

$$(10.3.7) \qquad \overrightarrow{BD} = \overrightarrow{AC}$$

From the equality (10.3.7) and the axiom 10.3.1.4, it follows that $\overrightarrow{BA} = \overrightarrow{DC}$. Therefore, the equality (10.3.5) does not depend on point $A$. □

Theorem 10.3.5. *Sum of vectors $v$ and $w$ does not depend on point $A$.*

Proof. Let

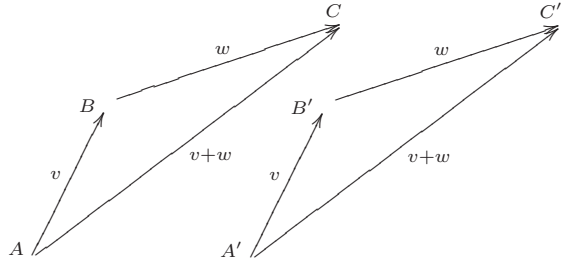

$$(10.3.8) \qquad v = \overrightarrow{AB} = \overrightarrow{A'B'}$$

$$(10.3.9) \qquad w = \overrightarrow{BC} = \overrightarrow{B'C'}$$

We define sum of vectors $v$ and $w$ according to the definition 10.3.2.

$$\overrightarrow{AB} + \overrightarrow{BC} = \overrightarrow{AC}$$
$$\overrightarrow{A'B'} + \overrightarrow{B'C'} = \overrightarrow{A'C'}$$

According to axiom 10.3.1.4, from equalities (10.3.8), (10.3.9), it follows that

$$(10.3.10) \qquad \overrightarrow{A'A} = \overrightarrow{B'B} = \overrightarrow{C'C}$$

Applying axiom 10.3.1.4 to outermost members of equality (10.3.10), we get

$$(10.3.11) \qquad \overrightarrow{A'C'} = \overrightarrow{AC}$$

From the equality (10.3.11) the statement of theorem follows.  □

Theorem 10.3.6. *Sum of vectors is associative.*

Proof. Let $v = \overrightarrow{AB}$, $w = \overrightarrow{BC}$, $u = \overrightarrow{CD}$. From the equality

$$\begin{aligned} v \;+\; w &= \overrightarrow{AC} \\ \overrightarrow{AB} \;+\; \overrightarrow{BC} &= \overrightarrow{AC} \end{aligned}$$

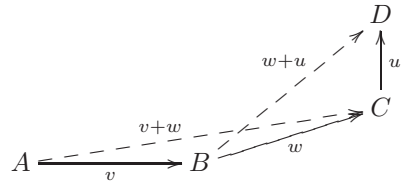

it follows that

$$(10.3.12) \qquad \begin{aligned} (v+w) \;+\; u &= \overrightarrow{AD} \\ \overrightarrow{AC} \;+\; \overrightarrow{CD} &= \overrightarrow{AD} \end{aligned}$$

From the equality

$$\begin{aligned} w \;+\; u &= \overrightarrow{BD} \\ \overrightarrow{BC} \;+\; \overrightarrow{CD} &= \overrightarrow{BD} \end{aligned}$$

it follows that

$$(10.3.13) \qquad \begin{aligned} v \;+\; (w+u) &= \overrightarrow{AD} \\ \overrightarrow{AB} \;+\; \overrightarrow{BD} &= \overrightarrow{AD} \end{aligned}$$

The theorem follows from comparison of equalities (10.3.12) and (10.3.13).  □

Theorem 10.3.7. *The structure of Abelian group is defined on the set $V$.*

PROOF. From theorems 10.3.3, 10.3.4, 10.3.5, 10.3.6 it follows that sum of vectors determines group.

Let  $v = \overrightarrow{AB}, \quad w = \overrightarrow{BC}$.

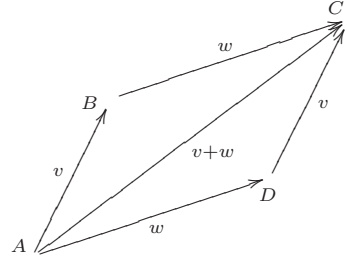

$$(10.3.14) \qquad \begin{array}{cccc} v & + & w & = \overrightarrow{AC} \\ \overrightarrow{AB} & + & \overrightarrow{BC} & = \overrightarrow{AC} \end{array}$$

According to axiom 10.3.1.3, there exists the point $D$ such that

$$w = \overrightarrow{AD} = \overrightarrow{BC}$$

The parallelogram law.

According to axiom 10.3.1.4,  $\overrightarrow{AB} = \overrightarrow{DC} = v$.  According to definition of sum of vectors

$$(10.3.15) \qquad \begin{array}{cccc} \overrightarrow{AD} & + & \overrightarrow{DC} & = \overrightarrow{AC} \\ w & + & v & = \overrightarrow{AC} \end{array}$$

Commutativity of sum follows from equalities (10.3.14) and (10.3.15).  □

THEOREM 10.3.8. *The map*

$$(10.3.16) \qquad\qquad V \to \mathrm{End}(\emptyset, \overset{\circ}{V})$$

*defined by the equality* (10.3.1)*, is a single transitive representation of Abelian group* $V$.

PROOF. The axiom 10.3.1.3 determines the map (10.3.16). From theorem 10.3.5, it follows that the map (10.3.16) is a representation. Efficiency of the representation follows from theorem 10.3.3 and axiom 10.3.1.2. From the axiom 10.3.1.2, it also follows that representation is transitive. Effective and transitive representation is single transitive.  □

If we compare the theorem 10.3.8 and statements of the section 10.2, then we see that a single transitive representation of Abelian group $V$ on the set $\overset{\circ}{V}$ is equivalent to axioms of affine space. However, if we use the theorem 10.3.8 as a definition of affine space, we lose many important constructions in affine space. For instance, vector generates parallel translation in affine space. However, we do not have a tool to define rotation of affine space.

If we look carefully at the definition 10.3.1, then we will see that Abelian group $V$ has additional structure since Abelian group $V$ is $D$-module. Thus, we get the following theorem.

THEOREM 10.3.9. *Let $D$ be commutative ring, $V$ be Abelian group, and $\overset{\circ}{V}$ be any set. If $A \in \overset{\circ}{V}$ and $v \in V$, then we use an expression $A + v$ to denote the action of vector $v$ at the point $A$.* **Affine space** *over $D$-module $V$ is the diagram of representations*

$$\overset{\rightarrow}{V} : D \xrightarrow{\ f_{12}\ } V \xrightarrow{\ f_{23}\ } \overset{\circ}{V} \qquad \begin{array}{l} f_{12}(d) : v \to d\, v \\[4pt] f_{23}(v) : A \to A + v \end{array}$$

*where $f_{12}$ is effective representation of commutative ring $D$ in Abelian group $V$ and $f_{23}$ is single transitive right-side representation of Abelian group $V$ in the set $\overset{\circ}{V}$.*

PROOF. We assume that the set $\overset{\circ}{V}$ is not empty; therefore the set $\overset{\circ}{V}$ satisfies the axiom 10.3.1.1. Since $v \in V$ generates the transformation of the set, then, for any $A \in M$, $B \in M$ is defined uniquely such that

$$B = A + v$$

This statement proves the axiom 10.3.1.3. Since the representation $f_{23}$ is single transitive, then for any $A$, $B \in \overset{\circ}{V}$ there exists unique $v \in V$ such that

$$B = A + v$$

This statement allows us to introduce notation $\overrightarrow{AB} = a$, as well this statement proves the axiom 10.3.1.2. The axiom 10.3.1.4 follows from the statement of the theorem 10.2.4.

The representation $f_{12}$ assures that Abelian group $V$ is $D$-module.          □

The Abelian group $V$ acts single transitive on the set $\overset{\circ}{V}$. From construction in section 10.2, it follows that the basis of the set $\overset{\circ}{V}$ relative to representation of the Abelian group $V$ consists of one point. This point is usually denoted by the letter $O$ and is called **origin of coordinate system of affine space**. Therefore, an arbitrary point $A \in \overset{\circ}{V}$ can be represented using vector $\overrightarrow{OA} \in V$

Let $\overline{\overline{e}}$ be the basis of $D$-module $V$. Then the vector $\overrightarrow{OA}$ has form

$$\overrightarrow{OA} = a^i e_i$$

The set $(a_i, i \in I)$ is called **coordinates of point $A$ of affine space $\overset{\circ}{A}$ relative to basis $(O, \overline{\overline{e}})$**.

## 10.4. Affine Space on Differentiable Manifold

In the section 10.3 we considered the definition of affine geometry. Below we consider a model of affine space in a metric-affine manifold. When we consider connection $\Gamma_{ij}^k$ in Riemann space, we impose a constraint on connection, [10.3] that the torsion

(10.4.1) $$T_{kl}^i = \Gamma_{lk}^i - \Gamma_{kl}^i$$

is 0 (symmetry of connection) and parallel transport does not change scalar product of vectors. If a metric tensor and an arbitrary connection are defined on a differentiable manifold, then this manifold is called **metric-affine manifold**. [10.4] In particular, connection in metric-affine manifold has torsion.

---

[10.3] See the definition of affine connection in Riemann space on the page [4]-443.
[10.4] See also the definition [9]-5.4.1.

In Riemann space, we use geodesics instead of straight lines. So we can represent the vector $v$ using segment $AB$ of geodesic $L_v$ such that vector $v$ is tangent to geodesic $L_v$ at the point $A$ and the length of segment $AB$ equals to the length of the vector $v$.

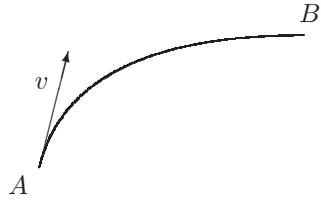

This definition allows us to identify the vector $v$ and the segment $AB$ of geodesic $L_v$.

For given vectors $v$ and $w$ in tangent plane at the point $A$, let $\rho > 0$ be the length of the vector $v$ and $\sigma > 0$ be the length of the vector $w$. Let $V$ be unit vector collinear to the vector $v$

(10.4.2) $\qquad V^k \rho = v^k$

Let $W$ be unit vector collinear to the vector $w$

(10.4.3) $\qquad W^k \sigma = w^k$

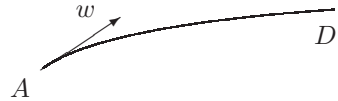

We draw geodesic $L_v$ through the point $A$ using the vector $v$ as a tangent vector to $L_v$ in the point $A$. Let $\tau$ be the canonical parameter on $L_v$ and

$$\frac{dx^k}{d\tau} = V^k$$

We transfer the vector $w$ along the geodesic $L_v$ from the point $A$ into point $B$ that defined by value of the parameter $\tau = \rho$. We mark the result as $w'$.

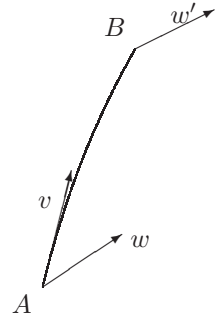

We draw geodesic $L_{w'}$ through the point $B$ using the vector $W'$ as a tangent vector to $L_{w'}$ in the point $B$. Let $\varphi'$ be the canonical parameter on $L_{w'}$ and

$$\frac{dx^k}{d\varphi'} = W'^k$$

We define point $C$ on the geodesic $L_{w'}$ by parameter value $\varphi' = \sigma$

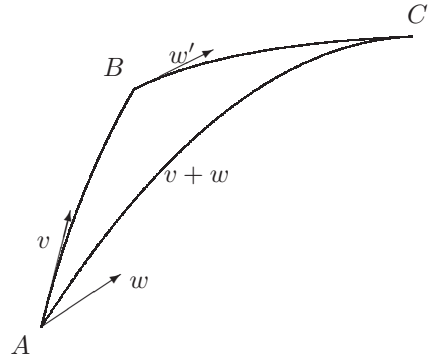

I assume that length of vectors $v$ and $w$ is small. Then there exists unique geodesic $L_u$ from point $A$ to point $C$. I will identify segment $AC$ of geodesic $L_u$ and vector $v + w$.

The same way, I draw triangle $ADE$ to find vector $w + v$.

We draw geodesic $L_w$ through the point $A$ using the vector $w$ as a tangent vector to $L_w$ in the point $A$. Let $\varphi$ be the canonical parameter on $L_w$ and

$$\frac{dx^k}{d\varphi} = W^k$$

We transfer the vector $v$ along the geodesic $L_w$ from the point $A$ into point $D$ that defined by value of the parameter $\varphi = \sigma$. We mark the result as $v'$.

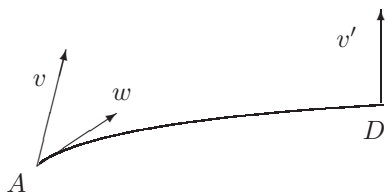

We draw geodesic $L_{v'}$ through the point $D$ using the vector $v'$ as a tangent vector to $L_{v'}$ in the point $D$. Let $\tau'$ be the canonical parameter on $L_{v'}$ and

$$\frac{dx^k}{d\tau'} = V'^k$$

We define point $E$ on the geodesic $L_{v'}$ by parameter value $\tau' = \rho$

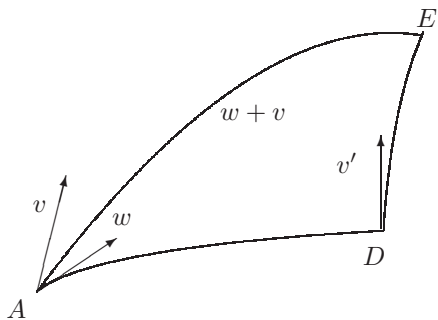

There exists unique geodesic $L_u$ from point $A$ to point $E$. I will identify segment $AE$ of geodesic $L_u$ and vector $w + v$.

Formally the lines $AB$ and $DE$ as well as the lines $AD$ and $BC$ are parallel lines. The lengths of $AB$ and $DE$ are the same, and the lengths of $AD$ and $BC$ are the same as well. We call this figure a **parallelogram** based on vectors $v$ and $w$ with the origin in the point $A$.

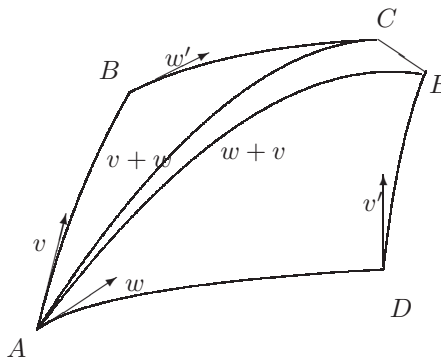

LEMMA 10.4.1. *Let $L_v$ be a geodesic through the point $A$ and the vector $v$ be a tangent vector to $L_v$ in the point $A$. An increase of coordinate $x^k$ along geodesic $L_v$ is*

$$(10.4.4) \qquad \Delta x^k = \frac{dx^k}{d\tau}\tau - \frac{1}{2}\Gamma^k_{mn}\frac{dx^m}{d\tau}\frac{dx^n}{d\tau}\tau^2 + O(\tau^2)$$

*where $\tau$ is canonical parameter and we take values of derivatives and components $\Gamma^k_{mn}$ in the initial point.*

PROOF. The system of differential equations of geodesic $L_v$ has the following form

$$(10.4.5) \qquad \frac{d^2x^i}{d\tau^2} = -\Gamma^i_{kl}\frac{dx^k}{d\tau}\frac{dx^l}{d\tau}$$

We write Taylor expansion of solution of the system of differential equations (10.4.5) in the following form

$$(10.4.6) \qquad \begin{aligned}\Delta x^k &= \frac{dx^k}{d\tau}\tau + \frac{1}{2}\frac{d^2x^k}{d\tau^2}\tau^2 + O(\tau^2) = \\ &= \frac{dx^k}{d\tau}\tau - \frac{1}{2}\Gamma^k_{mn}\frac{dx^m}{d\tau}\frac{dx^n}{d\tau}\tau^2 + O(\tau^2)\end{aligned}$$

The equality (10.4.4) follows from the equality (10.4.6). $\square$

THEOREM 10.4.2. *Suppose CBADE is a parallelogram with an origin in the point A; then the resulting figure will not be closed [1]. The value of the difference of coordinates of points C and E is equal to surface integral of the torsion over this parallelogram*

$$\Delta_{CE}x^k = \iint T^k_{mn}dx^m \wedge dx^n$$

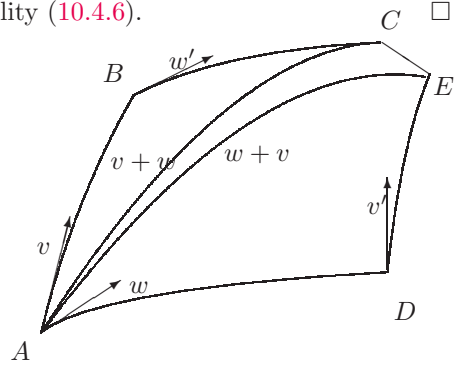

PROOF. According to the lemma 10.4.1, an increase of coordinate $x^k$ along the geodesic $L_v$ has the following form [10.5]

$$\Delta_{CE}x^k = \iint T^k_{mn}dx^m \wedge dx^n$$

$$\Delta_{AB}x^k = V^k\rho - \frac{1}{2}\Gamma^k_{mn}(A)V^mV^n\rho^2 + O(\rho^2)$$

and an increase of coordinate $x^k$ along the geodesic $L_{b'}$ has the following form

$$(10.4.7) \qquad \Delta_{BC}x^k = W'^k\sigma - \frac{1}{2}\Gamma^k_{mn}(B)W'^mW'^n\sigma^2 + O(\sigma^2)$$

Here

$$(10.4.8) \qquad \begin{aligned}W'^k &= W^k - \Gamma^k_{mn}(A)W^m\Delta_{AB}x^n + O(dx) \\ &= W^k - \Gamma^k_{mn}(A)W^mV^n\rho + O(\rho)\end{aligned}$$

is the result of parallel transport of the vector $w$ from $A$ to $B$ and

$$(10.4.9) \qquad \begin{aligned}\Gamma^k_{mn}(B) &= \Gamma^k_{mn}(A) + \partial_p\Gamma^k_{mn}(B)\Delta_{AB}x^p \\ &= \Gamma^k_{mn}(A) + \partial_p\Gamma^k_{mn}(B)V^p\rho\end{aligned}$$

with precision of small value of first level. Putting (10.4.8), (10.4.9) into (10.4.7) we will receive

$$\Delta_{BC}x^k = W^k\sigma - \Gamma^k_{mn}(A)W^mV^n\sigma\rho - \frac{1}{2}\Gamma^k_{mn}(A)W^mW^n\sigma^2 + O(\rho^2)$$

---

[10.5] Proof of this statement I found in [3]

Total increase of coordinate $x^K$ along the way $ABC$ has form

$$
\begin{aligned}
\Delta_{ABC} x^k = {} & \Delta_{AB} x^k + \Delta_{BC} x^k \\
(10.4.10) \qquad = {} & V^k \rho + W^k \sigma - \Gamma^k_{mn}(A) W^m V^n \sigma \rho - \\
& - \frac{1}{2} \Gamma^k_{mn}(A) W^m W^n \sigma^2 - \frac{1}{2} \Gamma^k_{mn}(A) V^m V^n \rho^2 + O(dx^2)
\end{aligned}
$$

In a similar way, total increase of coordinate $x^K$ along the way $ADE$ has form

$$
\begin{aligned}
\Delta_{ADE} x^k = {} & \Delta_{AD} x^k + \Delta_{DE} x^k = \\
(10.4.11) \qquad = {} & W^k \sigma + V^k \rho - \Gamma^k_{mn}(A) V^m W^n \rho \sigma - \\
& - \frac{1}{2} \Gamma^k_{mn}(A) V^m V^n \rho^2 - \frac{1}{2} \Gamma^k_{mn}(A) W^m W^n \sigma^2 + O(dx^2)
\end{aligned}
$$

From (10.4.10) and (10.4.11), it follows that

$$
\begin{aligned}
\Delta_{ADE} x^k - \Delta_{ABC} x^k = {} & - \Gamma^k_{mn}(A) V^m W^n \rho \sigma \\
& - \frac{1}{2} \Gamma^k_{mn}(A) V^m V^n \rho^2 - \frac{1}{2} \Gamma^k_{mn}(A) W^m W^n \sigma^2 \\
& + \Gamma^k_{mn}(A) W^m V^n \sigma \rho \\
& + \frac{1}{2} \Gamma^k_{mn}(A) W^m W^n \sigma^2 + \frac{1}{2} \Gamma^k_{mn}(A) V^m V^n \rho^2
\end{aligned}
$$

and we get integral sum for expression

$$
\Delta_{ADE} x^k - \Delta_{ABC} x^k = \iint_\Sigma (\Gamma^k_{nm} - \Gamma^k_{mn}) dx^m \wedge dx^n
$$

$\square$

THEOREM 10.4.3. *In Riemann space the parallelogram $ABCD$ is closed. At the point $A$, geodesic $AC$ has a tangent vector $u$ which is sum of vectors $v$ and $w$*

$$
(10.4.12) \qquad u^k = v^k + w^k
$$

*Therefore, a sum of vectors in Riemann space is commutative.*

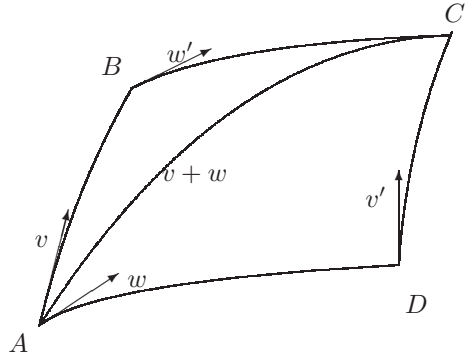

PROOF. Let $\pi$ be the length of the vector $u$. Let $U$ be unit vector collinear to the vector $u$

$$
(10.4.13) \qquad U^k \pi = u^k
$$

According to the lemma 10.4.1, an increase of coordinate $x^k$ along the geodesic $L_u$ has the following form

$$
(10.4.14) \qquad \Delta_{AC} x^k = U^k \pi - \frac{1}{2} \Gamma^k_{mn}(A) U^m U^n \pi^2 + O(\pi^2)
$$

Equalities

(10.4.15)
$$U^k \pi = V^k \rho + W^k \sigma$$

(10.4.16)
$$\Gamma_{mn}^k(A) U^m U^n \pi^2$$
$$= 2\Gamma_{mn}^k(A) W^m V^n \sigma\rho + \Gamma_{mn}^k(A) W^m W^n \sigma^2 + \Gamma_{mn}^k(A) V^m V^n \rho^2$$

follow from equalities (10.4.10), (10.4.14). The equality (10.4.12) follows from equalities (10.4.2), (10.4.3), (10.4.13), (10.4.15). The equality

(10.4.17)
$$\Gamma_{mn}^k(A) u^m u^n = 2\Gamma_{mn}^k(A) w^m v^n + \Gamma_{mn}^k(A) w^m w^n + \Gamma_{mn}^k(A) v^m v^n$$

follows from equalities (10.4.2), (10.4.3), (10.4.13), (10.4.16). The equality

(10.4.18)
$$\Gamma_{mn}^k(A) u^m u^n = \Gamma_{mn}^k(A)(v^m + w^m)(v^n + w^n)$$
$$= \Gamma_{mn}^k(A)(v^m v^n + v^m w^n + w^m v^n + w^m w^n)$$

follows from the equality (10.4.12). The equality (10.4.17) follows from the equality (10.4.18) and from the symmetry of connection. Therefore, the geodesic $AC$ is sum of geodesic $AB$ and $BC$.          □

If connection is not symmetric, then geodesic $L_u$ does not contain points $C$ and $E$. Therefore, sum of vectors in metric-affine manifold is noncommutative.

THEOREM 10.4.4. *There exists vector $t$ such that*

(10.4.19)
$$(v + w)^k = v^k + w^k + t^k$$

(10.4.20)
$$(w + v)^k = v^k + w^k - t^k$$

*Coordinates of the vector $t$ satisfy system of equations*

(10.4.21)  $\Gamma_{mn}^k(A) t^m t^n + (\Gamma_{mn}^k(A) + \Gamma_{nm}^k(A))(v^m + w^m) t^n + 2 T_{mn}^k(A) v^m w^n = 0$

PROOF. We first consider the vector $v + w$. The equality

(10.4.22)
$$v^k + w^k + t^k - \frac{1}{2}\Gamma_{mn}^k(A)(v^m + w^m + t^m)(v^n + w^n + t^n)$$
$$= v^k + w^k - \Gamma_{mn}^k(A) w^m v^n - \frac{1}{2}\Gamma_{mn}^k(A) w^m w^n - \frac{1}{2}\Gamma_{mn}^k(A) v^m v^n$$

follows from the equality (10.4.10) and the lemma 10.4.1. The equality

(10.4.23)
$$v^k + w^k + t^k$$
$$- \frac{1}{2}\Gamma_{mn}^k(A) v^m v^n + \frac{1}{2}\Gamma_{mn}^k(A) v^m w^n + \frac{1}{2}\Gamma_{mn}^k(A) v^m t^n$$
$$+ \frac{1}{2}\Gamma_{mn}^k(A) w^m v^n + \frac{1}{2}\Gamma_{mn}^k(A) w^m w^n + \frac{1}{2}\Gamma_{mn}^k(A) w^m t^n$$
$$+ \frac{1}{2}\Gamma_{mn}^k(A) t^m v^n + \frac{1}{2}\Gamma_{mn}^k(A) t^m w^n + \frac{1}{2}\Gamma_{mn}^k(A) t^m t^n$$
$$= v^k + w^k - \Gamma_{mn}^k(A) w^m v^n - \frac{1}{2}\Gamma_{mn}^k(A) w^m w^n - \frac{1}{2}\Gamma_{mn}^k(A) v^m v^n$$

follows from the equality (10.4.22). The equality

$$
(10.4.24) \quad \begin{aligned}
&t^k - \frac{1}{2}\Gamma^k_{mn}(A)v^m w^n - \frac{1}{2}\Gamma^k_{mn}(A)v^m t^n - \frac{1}{2}\Gamma^k_{mn}(A)w^m t^n \\
&-\frac{1}{2}\Gamma^k_{mn}(A)t^m v^n - \frac{1}{2}\Gamma^k_{mn}(A)t^m w^n - \frac{1}{2}\Gamma^k_{mn}(A)t^m t^n \\
&= -\frac{1}{2}\Gamma^k_{mn}(A)w^m v^n
\end{aligned}
$$

follows from the equality (10.4.23). The equality

$$
(10.4.25) \quad \begin{aligned}
&\Gamma^k_{mn}(A)t^m t^n + (\Gamma^k_{mn}(A)v^m + \Gamma^k_{mn}(A)w^m \\
&+ \Gamma^k_{nm}(A)v^m + \Gamma^k_{nm}(A)w^m - 2\delta^k_n)t^n \\
&+ 2T^k_{mn}(A)v^m w^n = 0
\end{aligned}
$$

follows from the equality (10.4.24). The equality (10.4.21) follows from the equality (10.4.25).

Now we consider the vector $w + v$. The equality

$$
(10.4.26) \quad \begin{aligned}
&v^k + w^k - t^k - \frac{1}{2}\Gamma^k_{mn}(A)(v^m + w^m - t^m)(v^n + w^n - t^n) \\
&= w^k + v^k - \Gamma^k_{mn}(A)v^m w^n - -\frac{1}{2}\Gamma^k_{mn}(A)v^m v^n - \frac{1}{2}\Gamma^k_{mn}(A)W^m W^n
\end{aligned}
$$

follows from the equality (10.4.11) and the lemma 10.4.1. The equality

$$
(10.4.27) \quad \begin{aligned}
&v^k + w^k - t^k \\
&-\frac{1}{2}\Gamma^k_{mn}(A)v^m v^n - \frac{1}{2}\Gamma^k_{mn}(A)v^m w^n + \frac{1}{2}\Gamma^k_{mn}(A)v^m t^n \\
&-\frac{1}{2}\Gamma^k_{mn}(A)w^m v^n - \frac{1}{2}\Gamma^k_{mn}(A)w^m w^n + \frac{1}{2}\Gamma^k_{mn}(A)w^m t^n \\
&+\frac{1}{2}\Gamma^k_{mn}(A)t^m v^n + \frac{1}{2}\Gamma^k_{mn}(A)t^m w^n - \frac{1}{2}\Gamma^k_{mn}(A)t^m t^n \\
&= w^k + v^k - \Gamma^k_{mn}(A)v^m w^n - -\frac{1}{2}\Gamma^k_{mn}(A)v^m v^n - \frac{1}{2}\Gamma^k_{mn}(A)W^m W^n
\end{aligned}
$$

follows from the equality (10.4.26). The equality

$$
(10.4.28) \quad \begin{aligned}
&- t^k \\
&+\frac{1}{2}\Gamma^k_{mn}(A)v^m t^n \\
&-\frac{1}{2}\Gamma^k_{mn}(A)w^m v^n + \frac{1}{2}\Gamma^k_{mn}(A)w^m t^n \\
&+\frac{1}{2}\Gamma^k_{mn}(A)t^m v^n + \frac{1}{2}\Gamma^k_{mn}(A)t^m w^n - \frac{1}{2}\Gamma^k_{mn}(A)t^m t^n \\
&= -\frac{1}{2}\Gamma^k_{mn}(A)v^m w^n
\end{aligned}
$$

follows from the equality (10.4.27). The equality

$$
(10.4.29) \quad \begin{aligned}
&\Gamma^k_{mn}(A)t^m t^n + (\Gamma^k_{mn}(A)v^m + \Gamma^k_{mn}(A)w^m \\
&+ \Gamma^k_{nm}(A)v^m + \Gamma^k_{nm}(A)w^m - 2\delta^k_n)t^n \\
&+ 2T^k_{nm}(A)v^m w^n = 0
\end{aligned}
$$

follows from the equality (10.4.28). The equality (10.4.21) follows from the equality (10.4.29). ☐

It is not a simple question to answer whether the system of equations (10.4.21) has a solution. However there is another way to find coordinates of vector $t$.

We draw geodesic $L_{v+w}$ through the point $A$ using the vector $v + w$ as a tangent vector to $L_{v+w}$ in the point $A$. We draw geodesic $L_{w+v}$ through the point $A$ using the vector $w + v$ as a tangent vector to $L_{w+v}$ in the point $A$. We draw geodesic $L_u$ through the point $A$ using the vector $u$

$$u^k = v^k + w^k$$

as a tangent vector to $L_u$ in the point $A$.

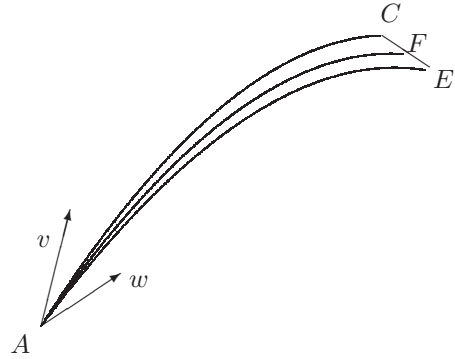

According to theorems 10.4.2, 10.4.4, the point $F$ is the middle of the segment $EC$. Therefore it is possible to consider the segment $AF$ as the median of the triangle $ACE$. According to the theorem 10.4.4, we identify the segment $FC$ and vektor $t$. Therefore, the theorem 10.4.2 gives us the way to find coordinates of vector $t$.

## 10.5. Noncommutative module

In the section 10.4, we considered opportunity to study affine geometry on affine manifold. This geometry has two features. The set of vectors is not closed relative sum and addition operation may be noncommutative.

We are not ready to consider first problem; however we can consider questions related to the noncommutativity of sum of vectors. The representation

$$f : D \longrightarrow{}_{*} G$$

of commutative ring $D$ in arbitrary group $G$ is called non-commutative module. This representation is much like a module, so all theorems about the structure of module are true. However, the question about the structure of basis remains open.

In general

$$av + bw \neq bw + av$$

Therefore, the question arises: what set of group $G$ we want to consider as a basis.

We can construct a basis the same way as we do a basis of module. Then this basis should permit expression

$$av + bw + cv$$

Or we may require items of basis to be in strict order in linear combination. In this case we assume that if $(v, w)$ is a basis of non-commutative module $V$, then for any expression $bw + av$ there exist $c, d \in D$ such that

$$cv + dw = bw + av$$

# Index

# Special Symbols and Notations

# Диаграмма представлений универсальных алгебр

Александр Клейн

Aleks_Kleyn@MailAPS.org
http://AleksKleyn.dyndns-home.com:4080/
http://sites.google.com/site/AleksKleyn/
http://arxiv.org/a/kleyn_a_1
http://AleksKleyn.blogspot.com/

Аннотация. Теория представлений универсальной алгебры является естественным развитием теории универсальной алгебры. В книге рассмотрены представление универсальной алгебры, диаграммы представлений и примеры представления. Морфизм представления - это отображение, сохраняющее структуру представления. Изучение морфизмов представлений ведёт к понятиям множества образующих и базиса представления.

# Оглавление

# Предисловие

## 1.1. Теория представлений

В статьях я часто рассматриваю вопросы, связанные с представлением универсальной алгебры. Вначале это были небольшие наброски, которые я многократно исправлял и переписывал. Но постепенно появлялись новые наблюдения. В результате вспомогательный инструмент превратился в стройную теорию.

Я это понял, когда я работал над книгой [10], и решил посвятить отдельную книгу вопросам, связанным с представлением универсальной алгебры. Изучение теории представлений универсальной алгебры показывает, что эта теория имеет много общего с теорией универсальной алгебры.

Основным толчком к более глубокому изучению представлений универсальной алгебры послужило определение векторного пространства как представление поля в абелевой группе. Я обратил внимание, что это определение меняет роль линейного отображения. По сути, линейное отображение - это отображение, которое сохраняет структуру представления. Эту конструкцию легко обобщить на произвольное представление универсальной алгебры. Таким образом появилось понятие морфизма представлений.

Множество невырожденных автоморфизмов векторного пространства порождает группу. Эта группа действует однотранзитивно на множестве базисов векторного пространства. Это утверждение является фундаментом теории инвариантов векторного пространства.

Возникает естественный вопрос. Можно ли обобщить эту конструкцию на произвольное представление? Базис - это не единственное множество, которое порождает векторное пространство. Если мы к множеству векторов базиса добавим произвольный вектор, то новое множество по прежнему порождает тоже самое векторное пространство, но базисом не является. Это утверждение является исходной точкой, от которой я начал изучение множества образующих представления. Множество образующих представления - это ещё одна интересная параллель теории представлений с теорией универсальной алгебры.

Множество автоморфизмов представления является лупой. Неассоциативность произведения порождает многочисленные вопросы, которые требуют дополнительное исследование. Все эти вопросы ведут к необходимости понимания теории инвариантов заданного представления.

Если мы рассматриваем теорию представлений универсальной алгебры как расширение теории универсальной алгебры, то почему не рассмотреть представление одного представления в другом представлении. Так появилась концепция башни представлений. Самый удивительный факт - это то, что все отображения в башне представлений действуют согласовано.

## 1.2. На грани теории

На протяжении многих лет я считал, что теория представлений является основным инструментом для изучения принципа общековариантности. Однако в процессе подготовки этой книги я неожиданно оказался на грани применимости теории представлений. Я не мог пройти мимо этого крайне важного события.

Точнее это было два разных открытия, связанных между собой темой некоммутативного сложения. Сначала я обнаружил, что я могу моделировать аффинную геометрию на многообразии аффинной связности (Тоже мне открытие. Думаю люди об этом знали со времён Декарта и Гауса). Здесь самым главным для меня было утверждение, что сумма определена не для любой пары векторов. Похожую задачу я видел, когда изучал многообразие базисов пространства Минковского ([11]). Если связность на аффинном многообразии имеет ненулевое кручение, то сумма векторов становится некоммутативной.

Позднее я решил исследовать представление кольца в неабелевой группе. Хотя алгебра замкнута относительно операции, я вижу возможность дальнейшего развития теории представлений. Мы можем пользоваться определением базиса из этой книги, однако некоторые важные детали будут спрятаны. Для меня интересна версия, что элементы базиса могут иметь заданный порядок, но сейчас я недостаточно ясно представляю какие из этого могут быть следствия.

# Предварительные определения

В этой главе собраны определения и теоремы, которые необходимы для понимания текста предлагаемой книги. Поэтому читатель может обращаться к утверждениям из этой главы по мере чтения основного текста книги.

## 2.1. Отношение эквивалентности

Определение 2.1.1. *Соответствие* $\Phi \in A \times A$ *называется* **отношением эквивалентности**, *если* [2.1]

2.1.1.1: *соответствие* $\Phi$ **рефлексивно**

$$(a, a) \in \Phi$$

2.1.1.2: *соответствие* $\Phi$ **симметрично**

$$(a, b) \in \Phi \Rightarrow (b, a) \in \Phi$$

2.1.1.3: *соответствие* $\Phi$ **транзитивно**

$$(a, b), (b, c) \in \Phi \Rightarrow (a, c) \in \Phi$$

$\square$

Теорема 2.1.2. *Для отображения*

$$f : A \to B$$

*множество*

$$(2.1.1) \qquad \ker f = \{(a, b) : a, b \in A, f(a) = f(b)\}$$

*является отношением эквивалентности и называется* **ядром отображения**. [2.2]

Доказательство.

Лемма 2.1.3. *Соответствие* $\ker f$ *рефлексивно*.

Доказательство. Из равенства

$$f(a) = f(a)$$

и определения (2.1.1) следует, что

$$(2.1.2) \qquad (a, a) \in \ker f$$

Лемма является следствием утверждения (2.1.2) и определения 2.1.1.1. $\qquad \odot$

Лемма 2.1.4. *Соответствие* $\ker f$ *симметрично*.

---

Доказательство. Равенство

$$(2.1.3) \qquad\qquad f(a) = f(b)$$

является следствием утверждения

$$(a, b) \in \ker f$$

и определения (2.1.1). Равенство

$$(2.1.4) \qquad\qquad f(b) = f(a)$$

является следствием равенства (2.1.3). Утверждение

$$(b, a) \in \ker f$$

является следствием равенства (2.1.4) и определения (2.1.1). Следовательно, мы доказали утверждение

$$(2.1.5) \qquad\qquad (a, b) \in \ker f \Rightarrow (b, a) \in \ker f$$

Лемма является следствием утверждения (2.1.5) и определения 2.1.1.2.    ⊙

Лемма 2.1.5. *Соответствие* $\ker f$ *транзитивно.*

Доказательство. Равенство

$$(2.1.6) \qquad\qquad f(a) = f(b)$$

является следствием утверждения

$$(a, b) \in \ker f$$

и определения (2.1.1). Равенство

$$(2.1.7) \qquad\qquad f(b) = f(c)$$

является следствием утверждения

$$(b, c) \in \ker f$$

и определения (2.1.1). Равенство

$$(2.1.8) \qquad\qquad f(a) = f(c)$$

является следствием равенств (2.1.6), (2.1.7). Утверждение

$$(a, c) \in \ker f$$

является следствием равенства (2.1.8) и определения (2.1.1). Следовательно, мы доказали утверждение

$$(2.1.9) \qquad\qquad (a, b), (b, c) \in \ker f \Rightarrow (a, c) \in \ker f$$

Лемма является следствием утверждения (2.1.9) и определения 2.1.1.2.    ⊙

Утверждение теоремы является следствием лемм 2.1.3, 2.1.4, 2.1.5 и определения 2.1.1.    □

Теорема 2.1.6. *Пусть $N$ - отношение эквивалентности на множестве $A$. Рассмотрим категорию $\mathcal{A}$ объектами которой являются отображения*[2.3]

$$f_1 : A \rightarrow S_1 \quad \ker f_1 \supseteq N$$

$$f_2 : A \rightarrow S_2 \quad \ker f_2 \supseteq N$$

---

*Мы определим морфизм $f_1 \to f_2$ как отображение $h : S_1 \to S_2$, для которого коммутативна диаграмма*

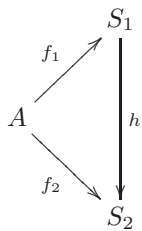

*Отображение*

$$\operatorname{nat} N : A \to A/N$$

*является универсально отталкивающим в категории $\mathcal{A}$.* [2.4]

Доказательство. Рассмотрим диаграмму

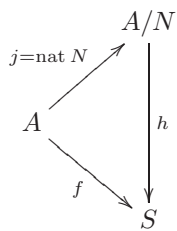

(2.1.10)                          $\ker f \supseteq N$

Из утверждения (2.1.10) и равенства

$$j(a_1) = j(a_2)$$

следует

$$f(a_1) = f(a_2)$$

Следовательно, мы можем однозначно определить отображение $h$ с помощью равенства

$$h(j(b)) = f(b)$$

$\square$

## 2.2. Универсальная алгебра

Определение 2.2.1. *Для любых множеств* [2.5] *$A$, $B$,* **декартова степень $B^A$** *- это множество отображений*

$$f : A \to B$$

$\square$

Определение 2.2.2. *Пусть дано множество $A$ и целое число $n \geq 0$. Отображение* [2.6]

$$\omega : A^n \to A$$

---

[2.4] Определение универсального объекта смотри в определении на с. [2]-47.

[2.5] Я следую определению из примера (iV), [14], страницы 17, 18.

[2.6] Определения 2.2.2, 2.2.7 опираются на определение в примере (vi), страница [14]-26.

*называется n-арной операцией на множестве* $A$ *или просто* **операцией на множестве** $A$. *Для любых* $a_1, ..., a_n \in A$, *мы пользуемся любой из форм записи* $\omega(a_1, ..., a_n)$, $a_1...a_n\omega$ *для обозначения образа отображения* $\omega$. □

Замечание 2.2.3. *Согласно определениям* 2.2.1, 2.2.2, *n-арная операция* $\omega \in A^{A^n}$. □

Определение 2.2.4. **Область операторов** - *это множество операторов*[2.7] $\Omega$ *вместе с отображением*

$$a : \Omega \to N$$

*Если* $\omega \in \Omega$, *то* $a(\omega)$ *называется* **арностью** *оператора* $\omega$. *Если* $a(\omega) = n$, *то оператор* $\omega$ *называется n-арным. Мы пользуемся обозначением*

$$\Omega(n) = \{\omega \in \Omega : a(\omega) = n\}$$

*для множества n-арных операторов.* □

Определение 2.2.5. *Пусть* $A$ - *множество, а* $\Omega$ - *область операторов.*[2.8] *Семейство отображений*

$$\Omega(n) \to A^{A^n} \quad n \in N$$

*называется структурой* $\Omega$-*алгебры на* $A$. *Множество* $A$ *со структурой* $\Omega$-*алгебры называется* $\Omega$-**алгеброй** $A_\Omega$ *или* **универсальной алгеброй**. *Множество* $A$ *называется* **носителем** $\Omega$-*алгебры.* □

Область операторов $\Omega$ описывает множество $\Omega$-алгебр. Элемент множества $\Omega$ называется оператором, так как операция предполагает некоторое множество. Согласно замечанию 2.2.3 и определению 2.2.5, каждому оператору $\omega \in \Omega(n)$ сопоставляется $n$-арная операция $\omega$ на $A$.

Теорема 2.2.6. *Пусть множество* $B$ *является* $\Omega$-*алгеброй. Тогда множество* $B^A$ *отображений*

$$f : A \to B$$

*также является* $\Omega$-*алгеброй.*

Доказательство. Пусть $\omega \in \Omega(n)$. Для отображений $f_1, ..., f_n \in B^A$, мы определим операцию $\omega$ равенством

$$(f_1...f_n\omega)(x) = f_1(x)...f_n(x)\omega$$

□

Определение 2.2.7. *Пусть* $B \subseteq A$. *Если для любых* $b_1, ..., b_n \in B$, $b_1...b_n\omega \in B$, *то мы говорим, что* $B$ **замкнуто относительно** $\omega$ *или что* $B$ **допускает операцию** $\omega$. □

Определение 2.2.8. $\Omega$-*алгебра* $B_\Omega$ *является* **подалгеброй** $\Omega$-*алгебры* $A_\Omega$, *если верны следующие утверждения*[2.9]

2.2.8.1: $B \subseteq A$.

2.2.8.2: *Если оператор* $\omega \in \Omega$ *определяет операции* $\omega_A$ *на* $A$ *и* $\omega_B$ *на* $B$, *то*

$$\omega_A|B = \omega_B$$

□

Определение 2.2.9. *Пусть $A$, $B$ - $\Omega$-алгебры и $\omega \in \Omega(n)$. Отображение*[2.10]

$$f : A \to B$$

**согласовано с операцией** *$\omega$, если, для любых $a_1$, ..., $a_n \in A$,*

(2.2.1)                      $$f(a_1)...f(a_n)\omega = f(a_1...a_n\omega)$$

*Отображение $f$ называется* **гомоморфизмом** *$\Omega$-алгебры $A$ в $\Omega$-алгебру $B$, если $f$ согласовано с каждым $\omega \in \Omega$. Мы обозначим* Hom$(\Omega; A \to B)$ *множество гомоморфизмов $\Omega$-алгебры $A$ в $\Omega$-алгебру $B$.*                    □

Теорема 2.2.10. *Если область операторов пуста, то гомоморфизм $\Omega$-алгебры $A$ в $\Omega$-алгебру $B$ - это отображение*

$$f : A \to B$$

*Следовательно,* Hom$(\emptyset; A \to B) = B^A$.

Доказательство. Теорема является следствием определений 2.2.1, 2.2.9.                    □

Определение 2.2.11. *Гомоморфизм $f$ называется*[2.11] **изоморфизмом** *между $A$ и $B$, если соответствие $f^{-1}$ является гомоморфизмом. Если существует изоморфизм между $A$ и $B$, то говорят, что $A$ и $B$ изоморфны, и пишут $A \cong B$. Инъективный гомоморфизм называется* **мономорфизмом**. *Сюръективный гомоморфизм называется* **эпиморфизмом**.                    □

Определение 2.2.12. *Гомоморфизм, источником и целью которого является одна и таже алгебра, называется* **эндоморфизмом**. *Мы обозначим* End$(\Omega; A)$ *множество эндоморфизмов $\Omega$-алгебры $A$. Эндоморфизм, который является изоморфизмом, называется* **автоморфизмом**.                    □

Теорема 2.2.13. End$(\Omega; A) = $ Hom$(\Omega; A \to A)$

Доказательство. Теорема является следствием определений 2.2.9, 2.2.12.                    □

Теорема 2.2.14. *Если область операторов пуста, то эндоморфизм множества $A$ - это отображение*

$$t : A \to A$$

*Следовательно,* End$(\emptyset; A) = A^A$.

Доказательство. Теорема является следствием теорем 2.2.10, 2.2.13.    □

Определение 2.2.15. *Если существует мономорфизм $\Omega$-алгебры $A$ в $\Omega$-алгебру $B$, то говорят, что $A$* **может быть вложена в** $B$.                    □

Определение 2.2.16. *Если существует эпиморфизм из $A$ в $B$, то $B$ называется* **гомоморфным образом** *алгебры $A$.*                    □

---

[2.10] Я следую определению на странице [14]-63.

[2.11] Я следую определению на странице [14]-63.

## 2.3. Декартово произведение универсальных алгебр

ОПРЕДЕЛЕНИЕ 2.3.1. *Пусть $\mathcal{A}$ - категория. Пусть $\{B_i, i \in I\}$ - множество объектов из $\mathcal{A}$. Объект*

$$P = \prod_{i \in I} B_i$$

*и множество морфизмов*

$$\{f_i : P \to B_i, i \in I\}$$

*называется* **произведением множества объектов** $\{B_i, i \in I\}$ **в категории** $\mathcal{A}$ [2.12], *если для любого объекта $R$ и множество морфизмов*

$$\{g_i : R \to B_i, i \in I\}$$

*существует единственный морфизм*

$$h : R \to P$$

*такой, что диаграмма*

$$f_i \circ h = g_i$$

*коммутативна для всех $i \in I$.*

*Если $|I| = n$, то для произведения множества объектов $\{B_i, i \in I\}$ в $\mathcal{A}$ мы так же будем пользоваться записью*

$$P = \prod_{i=1}^{n} B_i = B_1 \times ... \times B_n$$

$\square$

ПРИМЕР 2.3.2. *Пусть $\mathcal{S}$ - категория множеств.* [2.13] *Согласно определению 2.3.1, декартово произведение*

$$A = \prod_{i \in I} A_i$$

*семейства множеств $(A_i, i \in I)$ и семейство проекций на $i$-й множитель*

$$p_i : A \to A_i$$

*являются произведением в категории $\mathcal{S}$.*

$\square$

ТЕОРЕМА 2.3.3. *Произведение существует в категории $\mathcal{A}$ $\Omega$-алгебр. Пусть $\Omega$-алгебра $A$ и семейство морфизмов*

$$p_i : A \to A_i \quad i \in I$$

*является произведением в категории $\mathcal{A}$. Тогда*

2.3.3.1: *Множество $A$ является декартовым произведением семейства множеств $(A_i, i \in I)$*

2.3.3.2: *Гомоморфизм $\Omega$-алгебры*

$$p_i : A \to A_i$$

*является проекцией на $i$-й множитель.*

2.3.3.3: *Любое $A$-число $a$ может быть однозначно представлено в виде кортежа $(p_i(a), i \in I)$ $A_i$-чисел.*

2.3.3.4: *Пусть $\omega \in \Omega$ - $n$-арная операция. Тогда операция $\omega$ определена покомпонентно*

(2.3.1) $$a_1...a_n\omega = (a_{1i}...a_{ni}\omega, i \in I)$$

*где $a_1 = (a_{1i}, i \in I)$, ..., $a_n = (a_{ni}, i \in I)$ .*

Доказательство. Пусть

$$A = \prod_{i \in I} A_i$$

декартово произведение семейства множеств $(A_i, i \in I)$ и, для каждого $i \in I$, отображение

$$p_i : A \to A_i$$

является проекцией на $i$-й множитель. Рассмотрим диаграмму морфизмов в категории множеств $\mathcal{S}$

(2.3.2)
$$\begin{array}{ccc} A & \xrightarrow{p_i} & A_i \\ \omega \uparrow & \nearrow g_i & \\ A^n & & \end{array} \qquad p_i \circ \omega = g_i$$

где отображение $g_i$ определено равенством

$$g_i(a_1, ..., a_n) = p_i(a_1)...p_i(a_n)\omega$$

Согласно определению 2.3.1, отображение $\omega$ определено однозначно из множества диаграмм (2.3.2)

(2.3.3) $$a_1...a_n\omega = (p_i(a_1)...p_i(a_n)\omega, i \in I)$$

Равенство (2.3.1) является следствием равенства (2.3.3). $\qquad \square$

Определение 2.3.4. *Если $\Omega$-алгебра $A$ и семейство морфизмов*

$$p_i : A \to A_i \quad i \in I$$

*является произведением в категории $\mathcal{A}$, то $\Omega$-алгебра $A$ называется* **прямым** *или* **декартовым произведением** $\Omega$-**алгебр** $(A_i, i \in I)$ . $\qquad \square$

Теорема 2.3.5. *Пусть множество $A$ является декартовым произведением множеств $(A_i, i \in I)$ и множество $B$ является декартовым произведением множеств $(B_i, i \in I)$ . Для каждого $i \in I$, пусть*

$$f_i : A_i \to B_i$$

*является отображением множества $A_i$ в множество $B_i$. Для каждого $i \in I$, рассмотрим коммутативную диаграмму*

(2.3.4)

$$
\begin{array}{ccc}
B & \xrightarrow{\ p'_i\ } & B_i \\
{\scriptstyle f}\big\uparrow & & \big\uparrow{\scriptstyle f_i} \\
A & \xrightarrow[\ p_i\ ]{} & A_i
\end{array}
$$

*где отображения $p_i$, $p'_i$ являются проекцией на $i$-й множитель. Множество коммутативных диаграмм* (2.3.4) *однозначно определяет отображение*

$$f : A \to B$$

$$f(a_i, i \in I) = (f_i(a_i), i \in I)$$

Доказательство. Для каждого $i \in I$, рассмотрим коммутативную диаграмму

(2.3.5)

Пусть $a \in A$. Согласно утверждению 2.3.3.3, $A$-число $a$ может быть представлено в виде кортежа $A_i$-чисел

(2.3.6)          $a = (a_i, i \in I) \quad a_i = p_i(a) \in A_i$

Пусть

(2.3.7)                        $b = f(a) \in B$

Согласно утверждению 2.3.3.3, $B$-число $b$ может быть представлено в виде кортежа $B_i$-чисел

(2.3.8)          $b = (b_i, i \in I) \quad b_i = p'_i(b) \in B_i$

Из коммутативности диаграммы (1) и из равенств (2.3.7), (2.3.8) следует, что

(2.3.9)                         $b_i = g_i(b)$

Из коммутативности диаграммы (2) и из равенства (2.3.6) следует, что

$$b_i = f_i(a_i)$$

$\hfill\square$

Теорема 2.3.6. *Пусть $\Omega$-алгебра $A$ является декартовым произведением $\Omega$-алгебр $(A_i, i \in I)$ и $\Omega$-алгебра $B$ является декартовым произведением $\Omega$-алгебр $(B_i, i \in I)$. Для каждого $i \in I$, пусть отображение*

$$f_i : A_i \to B_i$$

*является гомоморфизмом $\Omega$-алгебры. Тогда отображение*

$$f : A \to B$$

*определённое равенством*

(2.3.10)                    $f(a_i, i \in I) = (f_i(a_i), i \in I)$

*является гомоморфизмом $\Omega$-алгебры.*

Доказательство. Пусть $\omega \in \Omega$ - n-арная операция. Пусть $a_1 = (a_{1i}, i \in I)$, ..., $a_n = (a_{ni}, i \in I)$ и $b_1 = (b_{1i}, i \in I)$, ..., $b_n = (b_{ni}, i \in I)$. Из равенств (2.3.1), (2.3.10) следует, что

$$f(a_1...a_n\omega) = f(a_{1i}...a_{ni}\omega, i \in I)$$
$$= (f_i(a_{1i}...a_{ni}\omega), i \in I)$$
$$= ((f_i(a_{1i}))...(f_i(a_{ni})), i \in I)$$
$$= (b_{1i}...b_{ni}\omega, i \in I)$$
$$f(a_1)...f(a_n)\omega = b_1...b_n\omega = (b_{1i}...b_{ni}\omega, i \in I)$$

□

ОПРЕДЕЛЕНИЕ 2.3.7. *Эквивалентность на $\Omega$-алгебре $A$, которая является подалгеброй $\Omega$-алгебры $A^2$, называется* **конгруэнцией** [2.14] *на $A$.*          □

ТЕОРЕМА 2.3.8 (первая теорема об изоморфизмах). *Пусть*

$$f : A \to B$$

*гомоморфизм $\Omega$-алгебр с ядром $s$. Тогда отображение $f$ имеет разложение*

2.3.8.1: **Ядро гомоморфизма** $\ker f = f \circ f^{-1}$ *является конгруэнцией на $\Omega$-алгебре $A$.*

2.3.8.2: *Множество $A/\ker f$ является $\Omega$-алгеброй.*

2.3.8.3: *Отображение*

$$p : a \in A \to a^{\ker f} \in A/\ker f$$

*является эпиморфизмом и называется* **естественным гомоморфизмом**.

2.3.8.4: *Отображение*

$$q : p(a) \in A/\ker f \to f(a) \in f(A)$$

*является изоморфизмом.*

2.3.8.5: *Отображение*

$$r : f(a) \in f(A) \to f(a) \in B$$

*является мономорфизмом.*

Доказательство. Утверждение 2.3.8.1 является следствием предложения II.3.4 ([14], страница 72). Утверждения 2.3.8.2, 2.3.8.3 являются следствием теоремы II.3.5 ([14], страница 72) и последующего определения. Утверждения 2.3.8.4, 2.3.8.5 являются следствием теоремы II.3.7 ([14], страница 74).          □

---

[2.14] Я следую определению на странице [14]-71.

## 2.4. Полугруппа

Обычно операция $\omega \in \Omega(2)$ называется произведением

$$ab\omega = ab$$

либо суммой

$$ab\omega = a + b$$

ОПРЕДЕЛЕНИЕ 2.4.1. *Пусть A является $\Omega$-алгеброй и $\omega \in \Omega(2)$. A-число e называется* **нейтральным элементом операции** $\omega$, *если для любого A-числа a верны равенства*

(2.4.1) $$ea\omega = a$$

(2.4.2) $$ae\omega = a$$

$\square$

ОПРЕДЕЛЕНИЕ 2.4.2. *Пусть A является $\Omega$-алгеброй. Операция $\omega \in \Omega(2)$ называется* **ассоциативной**, *если верно равенство*

$$a(bc\omega)\omega = (ab\omega)c\omega$$

$\square$

ОПРЕДЕЛЕНИЕ 2.4.3. *Пусть A является $\Omega$-алгеброй. Операция $\omega \in \Omega(2)$ называется* **коммутативной**, *если верно равенство*

$$ab\omega = ba\omega$$

$\square$

ОПРЕДЕЛЕНИЕ 2.4.4. *Пусть $\Omega = \{\omega\}$. Если операция $\omega \in \Omega(2)$ ассоциативна, то $\Omega$-алгебра называется* **полугруппой**. *Если операция в полугруппе коммутативна, то полугруппа называется* **абелевой полугруппой**. $\square$

# Представление универсальной алгебры

## 3.1. Представление универсальной алгебры

Определение 3.1.1. *Пусть множество $A_2$ является $\Omega_2$-алгеброй. Пусть на множестве преобразований $\operatorname{End}(\Omega_2, A_2)$ определена структура $\Omega_1$-алгебры. Гомоморфизм*

$$f : A_1 \to \operatorname{End}(\Omega_2; A_2)$$

*$\Omega_1$-алгебры $A_1$ в $\Omega_1$-алгебру $\operatorname{End}(\Omega_2, A_2)$ называется* **представлением $\Omega_1$-алгебры** *$A_1$ или $A_1$*-**представлением** *в $\Omega_2$-алгебре $A_2$.* □

Диаграмма

$$A_2 \xrightarrow{\ f(a)\ } A_2$$
$$\Big\Uparrow f$$
$$A_1$$

означает, что мы рассматриваем представление $\Omega_1$-алгебры $A_1$. Отображение $f(a)$ является образом $a \in A_1$. Мы будем также пользоваться записью

$$f : A_1 \xrightarrow{\quad * \quad} A_2$$

для обозначения представления $\Omega_1$-алгебры $A_1$ в $\Omega_2$-алгебре $A_2$.

Существует несколько способов описать представление. Мы можем указать отображение $f$, имея в виду что область определения - это $\Omega_1$-алгебра $A_1$ и область значений - это $\Omega_1$-алгебра $\operatorname{End}(\Omega_2, A_2)$. Либо мы можем указать $\Omega_1$-алгебру $A_1$ и $\Omega_2$-алгебру $A_2$, имея в виду что нам известна структура отображения $f$.[3.1]

Определение 3.1.2. *Мы будем называть представление*

$$f : A_1 \xrightarrow{\quad * \quad} A_2$$

*$\Omega_1$-алгебры $A_1$* **эффективным**,[3.2] *если отображение*

$$f : A_1 \to \operatorname{End}(\Omega_2; A_2)$$

*является изоморфизмом $\Omega_1$-алгебры $A_1$ в $\operatorname{End}(\Omega_2, A_2)$.* □

Теорема 3.1.3. *Представление*

$$f : A_1 \xrightarrow{\quad * \quad} A_2$$

---

[3.1] Например, мы рассматриваем векторное пространство $V$ над полем $D$ (раздел 9.3).

[3.2] Аналогичное определение эффективного представления группы смотри в [18], страница 16, [19], страница 111, [15], страница 51 (Кон называет такое представление точным). Смотри также теорему 5.4.2.

*эффективно тогда и только тогда, когда из утверждения $a_1 \neq b_1$, $a_1$, $b_1 \in$*
*$A_1$, следует существование $a_2 \in A_2$ такого, что*[3.3]

$$f(a_1)(a_2) \neq f(b_1)(a_2)$$

Доказательство. Пусть представление $f$ эффективно и $a_1 \neq b_1$. Если для любого $a_2 \in A_2$ верно равенство

$$f(a_1)(a_2) = f(b_1)(a_2)$$

то

$$f(a_1) = f(b_1)$$

Это противоречит утверждению, что представление $f$ эффективно.

Пусть из утверждения $a_1 \neq b_1$, $a_1$, $b_1 \in A_1$, следует существование $a_2 \in A_2$ такого, что

$$f(a_1)(a_2) \neq f(b_1)(a_2)$$

Следовательно, из утверждения $a_1 \neq b_1$, $a_1$, $b_1 \in A_1$, следует, что

$$f(a) \neq f(b)$$

Согласно определению 3.1.2, представление $f$ эффективно.           □

Определение 3.1.4. *Мы будем называть представление*

$$f : A_1 \longrightarrow_* A_2$$

$\Omega_1$-*алгебры* $A_1$ **свободным**,[3.4] *если из утверждения*

$$f(a_1)(a_2) = f(b_1)(a_2)$$

*для любого $a_2 \in A_2$ следует, что $a_1 = b_1$.*           □

Теорема 3.1.5. *Мы будем называть представление*

$$f : A_1 \longrightarrow_* A_2$$

$\Omega_1$-*алгебры* $A_1$ **свободным**, *если из утверждения $f(a_1) = f(b_1)$ следует, что*
$a_1 = b_1$.

Доказательство. Утверждение $f(a_1) = f(b_1)$ верно тогда и только тогда, когда

$$f(a_1)(a_2) = f(b_1)(a_2)$$

для любого $a_2 \in A_2$.           □

Теорема 3.1.6. *Свободное представление эффективно.*

---

[3.3] Для группы теорема 3.1.3 имеет следующий вид. *Представление*

$$f : A_1 \longrightarrow_* A_2$$

*эффективно тогда и только тогда, когда для любого $A_1$-числа $a_1 \neq e$ существует $a_2 \in A_2$*
*такое, что*

$$f(a_1)(a_2) \neq a_2$$

[3.4] Аналогичное определение свободного представления группы смотри в [18], страница 16.
Смотри также теорему 5.5.2.

Доказательство. Пусть отображение

$$f : A_1 \longrightarrow_* A_2$$

является свободным представлением. Пусть $a, b \in A_1$. Согласно определению 3.1.4, из утверждения

$$f(a_1)(a_2) = f(b_1)(a_2)$$

для любого $a_2 \in A_2$ следует, что $a_1 = b_1$. Следовательно, если $a_1 \neq b_1$, то существует $a_2 \in A_2$ такое, что

$$f(a_1)(a_2) \neq f(b_1)(a_2)$$

Согласно теореме 3.1.3, представление $f$ эффективно. $\qquad\square$

Замечание 3.1.7. *Представление группы вращений в аффинном пространстве эффективно, но не свободно, так как начало координат является неподвижной точкой любого преобразования.* $\qquad\square$

Определение 3.1.8. *Мы будем называть представление*

$$f : A_1 \longrightarrow_* A_2$$

$\Omega_1$-*алгебры* **транзитивным**, [3.5] *если для любых* $a, b \in A_2$ *существует такое* $g$, *что*

$$a = f(g)(b)$$

*Мы будем называть представление* $\Omega_1$-*алгебры* **однотранзитивным**, *если оно транзитивно и свободно.* $\qquad\square$

Теорема 3.1.9. *Представление однотранзитивно тогда и только тогда, когда для любых* $a, b \in A_2$ *существует одно и только одно* $g \in A_1$ *такое, что* $a = f(g)(b)$

Доказательство. Следствие определений 3.1.4 и 3.1.8. $\qquad\square$

Теорема 3.1.10. *Пусть*

$$f : A_1 \longrightarrow_* A_2$$

*однотранзитивное представление* $\Omega_1$-*алгебры* $A_1$ *в* $\Omega_2$-*алгебре* $A_2$. *Существует структура* $\Omega_1$-*алгебры на множестве* $A_2$.

Доказательство. Пусть $b \in A_2$, $\omega \in \Omega_1(n)$. Для любых $A_2$-чисел $b_1, ..., b_n$, существуют $A_1$-числа $a_1, ..., a_n$ такие, что

$$b_1 = f(a_1)(b) \quad ... \quad b_n = f(a_n)(b)$$

Мы определим операцию $\omega$ на множестве $A_2$ равенством

$$(3.1.1) \qquad b_1...b_n\omega = f(a_1...a_n\omega)(b)$$

Мы также требуем, что выбор $A_2$-числа $b$ не зависит от операции $\omega$. $\qquad\square$

---

[3.5] Аналогичное определение транзитивного представления группы смотри в [19], страница 110, [15], страница 51.

Теорема 3.1.11. *Пусть*

$$f : A_1 \overset{*}{\longrightarrow} A_2$$

*эффективное представление $\Omega_1$-алгебры $A_1$ в $\Omega_2$-алгебре $A_2$. Пусть $\omega \in \Omega_1(n) \cap \Omega_2(n)$. Тогда*

(3.1.2)                    $$f(a_1...a_n\omega)(b) = f(a_1)(b)...f(a_n)(b)\omega$$

## 3.2. Морфизм представлений универсальной алгебры

Теорема 3.2.1. *Пусть $A_1$ и $B_1$ - $\Omega_1$-алгебры. Представление $\Omega_1$-алгебры $B_1$*

$$g : B_1 \overset{*}{\longrightarrow} A_2$$

*и гомоморфизм $\Omega_1$-алгебры*

$$h : A_1 \to B_1$$

*определяют представление $f$ $\Omega_1$-алгебры $A_1$*

(3.2.1)

Доказательство. Отображение $f$ является гомоморфизмом $\Omega_1$-алгебры $A_1$ в $\Omega_1$-алгебру $\mathrm{End}(\Omega_2, A_2)$, так как отображение $g$ является гомоморфизмом $\Omega_1$-алгебры $B_1$ в $\Omega_1$-алгебру $\mathrm{End}(\Omega_2, A_2)$.                    $\square$

Мы будем также пользоваться диаграммой

вместо диаграммы (3.2.1).

Если мы изучаем представление $\Omega_1$-алгебры в $\Omega_2$-алгебрах $A_2$ и $B_2$, то нас интересуют отображения $A_2 \to B_2$, сохраняющие структуру представления.

Определение 3.2.2. *Пусть*

$$f : A_1 \overset{*}{\longrightarrow} A_2$$

*представление $\Omega_1$-алгебры $A_1$ в $\Omega_2$-алгебре $A_2$ и*

$$g : B_1 \overset{*}{\longrightarrow} B_2$$

*представление $\Omega_1$-алгебры $B_1$ в $\Omega_2$-алгебре $B_2$. Для $i = 1, 2$, пусть отображение*

$$r_i : A_i \to B_i$$

*является гомоморфизмом $\Omega_i$-алгебры. Кортеж отображений $r = (r_1, r_2)$ таких, что*

(3.2.2)                    $$r_2 \circ f(a) = g(r_1(a)) \circ r_2$$

*называется **морфизмом представлений** из $f$ в $g$. Мы также будем говорить, что определён **морфизм представлений $\Omega_1$-алгебры в $\Omega_2$-алгебре**.*                    $\square$

Замечание 3.2.3. *Мы можем рассматривать пару отображений* $r_1$, $r_2$ *как отображение*

$$F : A_1 \cup A_2 \to B_1 \cup B_2$$

*такое, что*

$$F(A_1) = B_1 \qquad F(A_2) = B_2$$

*Поэтому в дальнейшем кортеж отображений* $r = (r_1, r_2)$ *мы будем также называть отображением и пользоваться записью*

$$r : f \to g$$

*Пусть* $a = (a_1, a_2)$ *- кортеж A-чисел. Мы будем пользоваться записью*

$$r(a) = (r_1(a_1), r_2(a_2))$$

*для образа кортежа A-чисел при морфизме представлений* $r$.    □

Определение 3.2.4. *Если представления* $f$ *и* $g$ *совпадают, то морфизм представлений* $r = (r_1, r_2)$ *называется* **морфизмом представления** $f$.    □

Теорема 3.2.5. *Пусть*

$$f : A_1 \longrightarrow_* A_2$$

*представление* $\Omega_1$-*алгебры* $A_1$ *в* $\Omega_2$-*алгебре* $A_2$ *и*

$$g : B_1 \longrightarrow_* B_2$$

*представление* $\Omega_1$-*алгебры* $B_1$ *в* $\Omega_2$-*алгебре* $B_2$. *Отображение*

$$(r_1 : A_1 \to B_1, \ r_2 : A_2 \to B_2)$$

*является морфизмом представлений тогда и только тогда, когда*

$$(3.2.3) \qquad r_2(f(a)(m)) = g(r_1(a))(r_2(m))$$

Доказательство. Для произвольного $m \in A_2$ равенство (3.2.2) имеет вид (3.2.3).    □

Замечание 3.2.6. *Рассмотрим морфизм представлений*

$$(r_1 : A_1 \to B_1, \ r_2 : A_2 \to B_2)$$

*Мы можем обозначать элементы множества* $B_1$, *пользуясь буквой по образцу* $b \in B_1$. *Но если мы хотим показать, что* $b$ *является образом элемента* $a \in A_1$, *мы будем пользоваться обозначением* $r_1(a)$. *Таким образом, равенство*

$$r_1(a) = r_1(a)$$

*означает, что* $r_1(a)$ *(в левой части равенства) является образом* $a \in A_1$ *(в правой части равенства). Пользуясь подобными соображениями, мы будем обозначать элемент множества* $B_2$ *в виде* $r_2(m)$. *Мы будем следовать этому соглашению, изучая соотношения между гомоморфизмами* $\Omega_1$-*алгебр и отображениями между множествами, где определены соответствующие представления.*    □

Замечание 3.2.7. *Мы можем интерпретировать* (3.2.3) *двумя способами*

- *Пусть преобразование* $f(a)$ *отображает* $m \in A_2$ *в* $f(a)(m)$. *Тогда преобразование* $g(r_1(a))$ *отображает* $r_2(m) \in B_2$ *в* $r_2(f(a)(m))$.

- *Мы можем представить морфизм представлений из $f$ в $g$, пользуясь диаграммой*

(3.2.4)

*Из* (3.2.2) *следует, что диаграмма* (1) *коммутативна.*

*Мы будем также пользоваться диаграммой*

(3.2.5)

*вместо диаграммы* (3.2.4).                                                                    □

Теорема 3.2.8. *Рассмотрим представление*

$$f : A_1 \longrightarrow A_2$$

$\Omega_1$-*алгебры* $A_1$ *и представление*

$$g : B_1 \longrightarrow B_2$$

$\Omega_1$-*алгебры* $B_1$. *Морфизм*

$$(r_1 : A_1 \to B_1, \ r_2 : A_2 \to B_2)$$

*представлений из* $f$ *в* $g$ *удовлетворяет соотношению*

(3.2.6)          $r_2 \circ (f(a_1)...f(a_n)\omega) = (g(r_1(a_1))...g(r_1(a_n))\omega) \circ r_2$

*для произвольной операции* $\omega \in \Omega_1(n)$.

Доказательство. Так как $f$ - гомоморфизм, мы имеем

(3.2.7)                $r_2 \circ (f(a_1)...f(a_n)\omega) = r_2 \circ f(a_1...a_n\omega)$

Из (3.2.2) и (3.2.7) следует

(3.2.8)                $r_2 \circ (f(a_1)...f(a_n)\omega) = g(r_1(a_1...a_n\omega)) \circ r_2$

Так как $r_1$ - гомоморфизм, из (3.2.8) следует

(3.2.9)          $r_2 \circ (f(a_1)...f(a_n)\omega) = g(r_1(a_1)...r_1(a_n)\omega) \circ r_2$

Так как $g$ - гомоморфизм, из (3.2.9) следует (3.2.6).                                      □

Теорема 3.2.9. *Пусть отображение*

$$(r_1 : A_1 \to B_1, \ r_2 : A_2 \to B_2)$$

*является морфизмом из представления*

$$f : A_1 \longrightarrow_* \longrightarrow A_2$$

$\Omega_1$-*алгебры* $A_1$ *в представление*

$$g : B_1 \longrightarrow_* \longrightarrow B_2$$

$\Omega_1$-*алгебры* $B_1$. *Если представление* $f$ *эффективно, то отображение*

$$r_2^* : \mathrm{End}(\Omega_2; A_2) \to \mathrm{End}(\Omega_2; B_2)$$

*определённое равенством*

(3.2.10) $$r_2^*(f(a)) = g(r_1(a))$$

*является гомоморфизмом* $\Omega_1$-*алгебры.*

Доказательство. Так как представление $f$ эффективно, то для выбранного преобразования $f(a)$ выбор элемента $a$ определён однозначно. Следовательно, преобразование $g(r_1(a))$ в равенстве (3.2.10) определено корректно.

Так как $f$ - гомоморфизм, мы имеем

(3.2.11) $$r_2^*(f(a_1)...f(a_n)\omega) = r_2^*(f(a_1...a_n\omega))$$

Из (3.2.10) и (3.2.11) следует

(3.2.12) $$r_2^*(f(a_1)...f(a_n)\omega) = g(r_1(a_1...a_n\omega))$$

Так как $h$ - гомоморфизм, из (3.2.12) следует

(3.2.13) $$r_2^*(f(a_1)...f(a_n)\omega) = g(r_1(a_1)...r_1(a_n)\omega)$$

Так как $g$ - гомоморфизм,

$$r_2^*(f(a_1)...f(a_n)\omega) = g(r_1(a_1))...g(r_1(a_n))\omega = r_2^*(f(a_1))...r_2^*(f(a_n))\omega$$

следует из (3.2.13). Следовательно, отображение $r_2^*$ является гомоморфизмом $\Omega_1$-алгебры. $\square$

Теорема 3.2.10. *Пусть*

$$f : A_1 \longrightarrow_* \longrightarrow A_2$$

*однотранзитивное представление* $\Omega_1$-*алгебры* $A_1$ *и*

$$g : B_1 \longrightarrow_* \longrightarrow B_2$$

*однотранзитивное представление* $\Omega_1$-*алгебры* $B_1$. *Если отображение*

$$r_1 : A_1 \to B_1$$

*является гомоморфизмом* $\Omega_1$-*алгебры, то существует морфизм представлений из* $f$ *в* $g$

$$(r_1 : A_1 \to B_1, \ r_2 : A_2 \to B_2)$$

Доказательство. Выберем гомоморфизм $r_1$. Выберем элемент $m \in A_2$ и элемент $n \in B_2$. Чтобы построить отображение $r_2$, рассмотрим следующую диаграмму

Из коммутативности диаграммы (1) следует

$$r_2(f(a)(m)) = g(r_1(a))(r_2(m))$$

Для произвольного $m' \in A_2$ однозначно определён $a \in A_1$ такой, что $m' = f(a)(m)$. Следовательно, мы построили отображении $r_2$, которое удовлетворяет равенству (3.2.2). $\square$

Теорема 3.2.11. *Если представление*

$$f : A_1 \longrightarrow\!\!\ast\!\!\longrightarrow A_2$$

$\Omega_1$-*алгебры* $A_1$ *однотранзитивно и представление*

$$g : B_1 \longrightarrow\!\!\ast\!\!\longrightarrow B_2$$

$\Omega_1$-*алгебры* $B_1$ *однотранзитивно, то для заданного гомоморфизма* $\Omega_1$-*алгебры*

$$r_1 : A_1 \to B_1$$

*гомоморфизм* $\Omega_2$-*алгебры*

$$r_2 : A_2 \to B_2$$

*такой, что* $r = (r_1, r_2)$ *является морфизмом представлений из* $f$ *в* $g$, *определён однозначно с точностью до выбора образа* $n = r_2(m) \in B_2$ *заданного элемента* $m \in A_2$.

Доказательство. Из доказательства теоремы 3.2.10 следует, что выбор гомоморфизма $r_1$ и элементов $m \in A_2$, $n \in B_2$ однозначно определяет отображение $r_2$. $\square$

Теорема 3.2.12. *Если представление*

$$f : A_1 \longrightarrow\!\!\ast\!\!\longrightarrow A_2$$

$\Omega_1$-*алгебры* $A_1$ *однотранзитивно, то для любого эндоморфизма* $r_1 \in \mathrm{End}(\Omega_1; A_1)$ *существует морфизм представления* $f$

$$(r_1 : A_1 \to A_1, \ r_2 : A_2 \to A_2)$$

Доказательство. Рассмотрим следующую диаграмму

Утверждение теоремы является следствием теоремы 3.2.10.       □

### 3.3. Теорема о разложении морфизмов расслоений

Теорема 3.3.1. *Пусть*

$$f : A_1 \xrightarrow{\;\;*\;\;} A_2$$

*представление* $\Omega_1$*-алгебры* $A_1$,

$$g : B_1 \xrightarrow{\;\;*\;\;} B_2$$

*представление* $\Omega_1$*-алгебры* $B_1$,

$$h : C_1 \xrightarrow{\;\;*\;\;} C_2$$

*представление* $\Omega_1$*-алгебры* $C_1$. *Пусть определены морфизмы представлений* $\Omega_1$*-алгебры*

$$(p_1 : A_1 \to B_1, \; p_2 : A_2 \to B_2)$$
$$(q_1 : B_1 \to C_1, \; q_2 : B_2 \to C_2)$$

*Тогда определён морфизм представлений* $\Omega_1$*-алгебры*

$$(r_1 : A_1 \to C_1, \; r_2 : A_2 \to C_2)$$

*где* $r_1 = q_1 \circ p_1$, $r_2 = q_2 \circ p_2$. *Мы будем называть морфизм* $r = (r_1, r_2)$ *представлений из* $f$ *в* $h$ **произведением морфизмов** $p = (p_1, p_2)$ **и** $q = (q_1, q_2)$ **представлений универсальной алгебры**.

Доказательство. Мы можем представить утверждение теоремы, пользуясь диаграммой

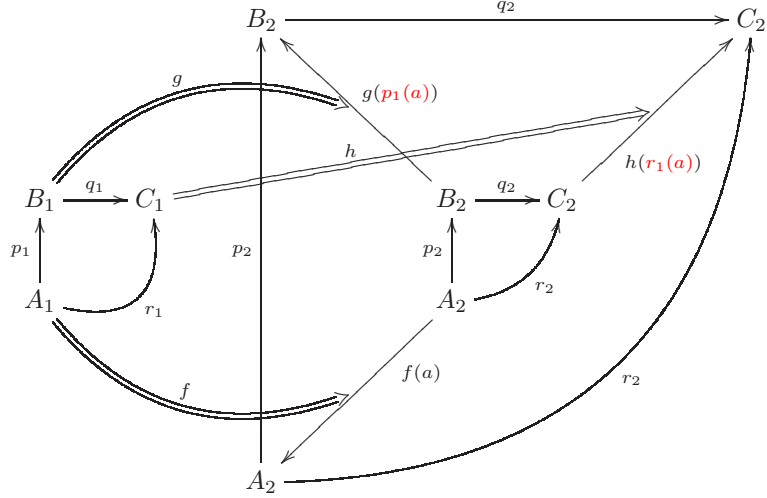

Отображение $r_1$ является гомоморфизмом $\Omega_1$-алгебры $A_1$ в $\Omega_1$-алгебру $C$. Нам надо показать, что отображение $r = (r_1, r_2)$ удовлетворяет (3.2.2):

$$r_2(f(a)(m)) = (q_2 \circ p_2)(f(a)(m))$$
$$= q_2(g(p_1(a))(p_2(m)))$$
$$= h((q_1 \circ p_1)(a))((q_2 \circ p_2)(m)))$$
$$= h(r_1(a))(r_2(m))$$

$\square$

Определение 3.3.2. *Пусть на множестве $A_2$ определена эквивалентность $S$. Преобразование $f$ называется* **согласованным с эквивалентностью** $S$, *если из условия $m_1 \equiv m_2 (\mathrm{mod} S)$ следует $f(m_1) \equiv f(m_2) (\mathrm{mod}\ S)$.* $\square$

Теорема 3.3.3. *Пусть на множестве $A_2$ определена эквивалентность $S$. Пусть на множестве $\mathrm{End}(\Omega_2, A_2)$ определена $\Omega_1$-алгебра. Если любое преобразование $f \in \mathrm{End}(\Omega_2; A_2)$, согласованно с эквивалентностью $S$, то мы можем определить структуру $\Omega_1$-алгебры на множестве $\mathrm{End}(\Omega_2; A_2/S)$.*

Доказательство. Пусть $h = \mathrm{nat}\ S$. Если $m_1 \equiv m_2 (\mathrm{mod} S)$, то $h(m_1) = h(m_2)$. Поскольку $f \in \mathrm{End}(\Omega_2; A_2)$ согласованно с эквивалентностью $S$, то $h(f(m_1)) = h(f(m_2))$. Это позволяет определить преобразование $F$ согласно правилу

$$(3.3.1) \qquad\qquad F([m]) = h(f(m))$$

Пусть $\omega$ - n-арная операция $\Omega_1$-алгебры. Пусть $f_1, ..., f_n \in \mathrm{End}(\Omega_2; A_2)$ и

$$F_1([m]) = h(f_1(m)) \quad ... \quad F_n([m]) = h(f_n(m))$$

Согласно условию теоремы, преобразование

$$f = f_1 ... f_n \omega \in \mathrm{End}(\Omega_2; A_2)$$

согласованно с эквивалентностью $S$. Следовательно, из условия $m_1 \equiv m_2(\mathrm{mod}\, S)$ и определения 3.3.2 следует

$$(3.3.2) \qquad \begin{aligned} f(m_1) &\equiv f(m_2)(\mathrm{mod}\, S) \\ (f_1...f_n\omega)(m_1) &\equiv (f_1...f_n\omega)(m_2)(\mathrm{mod}\, S) \end{aligned}$$

Следовательно, мы можем определить операцию $\omega$ на множестве $\mathrm{End}(\Omega_2; A_2/S)$ по правилу

$$(3.3.3) \qquad (F_1...F_n\omega)([m]) = h((f_1...f_n\omega)(m))$$

Из определения (3.3.1) и равенства (3.3.2) следует, что мы корректно определили операцию $\omega$ на множестве $\mathrm{End}(\Omega_2; A_2/S)$. $\qquad \square$

Определение 3.3.4. *Пусть*

$$f : A_1 \longrightarrow\!\!\!* \longrightarrow A_2$$

*представление $\Omega_1$-алгебры $A_1$,*

$$g : B_1 \longrightarrow\!\!\!* \longrightarrow B_2$$

*представление $\Omega_1$-алгебры $B_1$. Пусть*

$$(r_1 : A_1 \to B_1,\ r_2 : A_2 \to B_2)$$

*морфизм представлений из $f$ в $g$ такой, что $r_1$ - изоморфизм $\Omega_1$-алгебры и $r_2$ - изоморфизм $\Omega_2$-алгебры. Тогда отображение $r = (r_1, r_2)$ называется* **изоморфизмом представлений**. $\qquad \square$

Теорема 3.3.5. *Пусть*

$$f : A_1 \longrightarrow\!\!\!* \longrightarrow A_2$$

*представление $\Omega_1$-алгебры $A_1$,*

$$g : B_1 \longrightarrow\!\!\!* \longrightarrow B_2$$

*представление $\Omega_1$-алгебры $B_1$. Пусть*

$$(t_1 : A_1 \to B_1,\ t_2 : A_2 \to B_2)$$

*морфизм представлений из f в g. Тогда для отображений $t_1$, $t_2$ существуют разложения, которые можно описать диаграммой*

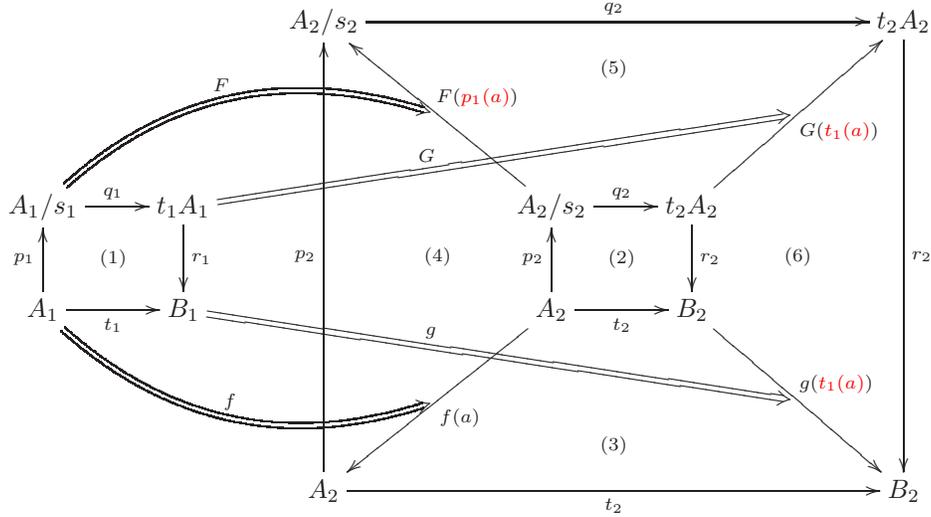

3.3.5.1: *Ядро гомоморфизма* $\ker t_i = t_i \circ t_i^{-1}$ *является конгруэнцией на $\Omega_i$-алгебре $A_i$, $i = 1, 2$.*

3.3.5.2: *Существует разложение гомоморфизма $t_i$, $i = 1, 2$,*

(3.3.4)                          $$t_i = r_i \circ q_i \circ p_i$$

3.3.5.3: *Отображения*

$$p_1(a) = a^{\ker t_1}$$
$$p_2(a) = a^{\ker t_2}$$

*являются естественными гомоморфизмами.*

3.3.5.4: *Отображения*

(3.3.5)                          $$q_1(p_1(a)) = t_1(a)$$

(3.3.6)                          $$q_2(p_2(a)) = t_2(a)$$

*являются изоморфизмами.*

3.3.5.5: *Отображения*

$$r_1 : t_1(a) \in f(A_1) \to t_1(a) \in B_1$$
$$r_2 : t_2(a) \in f(A_2) \to t_2(a) \in B_2$$

*являются мономорфизмами.*

3.3.5.6: *F - представление $\Omega_1$-алгебры $A_1/s$ в $A_2/s_2$*

3.3.5.7: *G - представление $\Omega_1$-алгебры $t_1 A_1$ в $t_2 A_2$*

3.3.5.8: *Отображение $p = (p_1, p_2)$ является морфизмом представлений f и F*

3.3.5.9: *Отображение $q = (q_1, q_2)$ является изоморфизмом представлений F и G*

3.3.5.10: *Отображение $r = (r_1, r_2)$ является морфизмом представлений G и g*

3.3.5.11: *Существует разложение морфизма представлений*

$$(3.3.7) \qquad (t_1, t_2) = (r_1, r_2) \circ (q_1, q_2) \circ (p_1, p_2)$$

Доказательство. Утверждения 3.3.5.1, 3.3.5.2, 3.3.5.3, 3.3.5.4, 3.3.5.5 являются следствием теоремы 2.3.8. Следовательно, диаграммы (1) и (2) коммутативны.

Мы начнём с диаграммы (4).

Пусть $m_1 \equiv m_2 (\operatorname{mod} \ker t_2)$. Следовательно,

$$(3.3.8) \qquad t_2(m_1) = t_2(m_2)$$

Если $a_1 \equiv a_2 (\operatorname{mod} \ker t_1)$, то

$$(3.3.9) \qquad t_1(a_1) = t_1(a_2)$$

Следовательно, $p_1(a_1) = p_1(a_2)$. Так как отображение $(t_1, t_2)$ - морфизм представлений, то

$$(3.3.10) \qquad t_2(f(a_1)(m_1)) = g(t_1(a_1))(t_2(m_1))$$
$$(3.3.11) \qquad t_2(f(a_2)(m_2)) = g(t_1(a_2))(t_2(m_2))$$

Из (3.3.8), (3.3.9), (3.3.10), (3.3.11) следует

$$(3.3.12) \qquad t_2(f(a_1)(m_1)) = t_2(f(a_2)(m_2))$$

Из (3.3.12) следует

$$(3.3.13) \qquad f(a_1)(m_1) \equiv f(a_2)(m_2)(\operatorname{mod} \ker t_2)$$

и, следовательно,

$$(3.3.14) \qquad p_2(f(a_1)(m_1)) = p_2(f(a_2)(m_2))$$

Из (3.3.14) следует, что отображение

$$(3.3.15) \qquad F(p_1(a))(p_2(m)) = p_2(f(a)(m))$$

определено корректно и является преобразованием множества $A_2/\ker t_2$.

Из равенства (3.3.13) (в случае $a_1 = a_2$) следует, что, для любого $a$, преобразование согласованно с эквивалентностью $\ker t_2$. Из теоремы 3.3.3 следует, что на множестве $\operatorname{End}(\Omega_2; A_2/\ker t_2)$. определена структура $\Omega_1$-алгебры. Рассмотрим $n$-арную операцию $\omega$ и $n$ преобразований

$$F(p_1(a_i))(p_2(m)) = p_2(f(a_i)(m)) \quad i = 1, ..., n$$

пространства $A_2/\ker t_2$. Мы положим

$$(F(p_1(a_1))...F(p_1(a_n))\omega)(p_2(m)) = p_2((f(a_1)...f(a_n)\omega)(m))$$

Следовательно, отображение $F$ является представлением $\Omega_1$-алгебры $A_1/\ker t_1$.

Согласно теореме 3.2.5, утверждение 3.3.5.8 является следствием (3.3.15).

Рассмотрим диаграмму (5).

Лемма 3.3.6. *Отображение $q = (q_1, q_2)$ является морфизмом представлений $F$ и $G$.*

Доказательство. Так как $q_2$ - биекция, то мы можем отождествить элементы множества $A_2/\ker t_2$ и множества $t_2(A_2)$, причём это отождествление имеет вид

$$(3.3.16) \qquad q_2(p_2(m)) = t_2(m)$$

Мы можем записать преобразование $F(p_1(a))$ множества $A_2/\ker t_2$ в виде

$$(3.3.17) \qquad F(p_1(a)) : p_2(m) \to F(p_1(a))(p_2(m))$$

Так как $T$ - биекция, то мы можем определить преобразование

$$(3.3.18) \qquad q_2(p_2(m)) \to q_2(F(p_1(a))(p_2(m)))$$

множества $RA_2$. Преобразование (3.3.18) зависит от $p_1(a) \in A_1/\ker t_1$. Так как $q_1$ - биекция, то мы можем отождествить элементы множества $A_1/\ker t_1$ и множества $t_1(A_1)$, причём это отождествление имеет вид

$$(3.3.19) \qquad q_1(p_1(a)) = t_1(a)$$

Следовательно, мы определили отображение

$$G : t_1(A_1) \to \operatorname{End}(\Omega_2; t_2(A_2))$$

согласно равенству

$$(3.3.20) \qquad G(q_1(p_1(a)))(q_2(p_2(m))) = q_2(F(p_1(a))(p_2(m)))$$

Рассмотрим $n$-арную операцию $\omega$ и $n$ преобразований

$$G(t_1(a_i))(t_2(m)) = q_2(F(p_1(a_i))(p_2(m))) \quad i = 1, ..., n$$

множества $t_2(A_2)$. Мы положим

$$(3.3.21) \quad (G(t_1(a_1))...G(t_1(a_n))\omega)(t_2(m)) = q_2((F(p_1(a_1))...F(p_1(a_n))\omega)(p_2(m)))$$

Согласно (3.3.20) операция $\omega$ корректно определена на множестве $\operatorname{End}(\Omega_2; t_2(A_2))$. Следовательно, отображение $G$ является представлением $\Omega_1$-алгебры.

Согласно теореме 3.2.5, лемма является следствием (3.3.20).          ⊙

Лемма 3.3.7. *Отображение* $(q_1^{-1}, q_2^{-1})$ *является морфизмом представлений $G$ и $F$.*

Доказательство. Так как $q_2$ - биекция, то из равенства (3.3.16) следует

$$(3.3.22) \qquad p_2(m) = q_2^{-1}(t_2(m))$$

Мы можем записать преобразование $G(t_1(a))$ множества $t_2(A_2)$ в виде

$$(3.3.23) \qquad G(t_1(a)) : t_2(m) \to G(t_1(a))(t_2(m))$$

Так как $q_2$ - биекция, то мы можем определить преобразование

$$(3.3.24) \qquad q_2^{-1}(t_2(m)) \to q_2^{-1}(G(t_1(a))(t_2(m)))$$

множества $A_2/\ker t_2$. Преобразование (3.3.24) зависит от $t_1(a) \in t_1(A_1)$. Так как $q_1$ - биекция, то из равенства (3.3.19) следует

$$(3.3.25) \qquad p_1(a) = q_1^{-1}(t_1(a))$$

Так как по построению диаграмма (5) коммутативна, то преобразование (3.3.24) совпадает с преобразованием (3.3.17). Равенство (3.3.21) можно записать в виде

$$(3.3.26) \qquad \begin{aligned} &q_2^{-1}((G(t_1(a_1))...G(t_1(a_n))\omega)(t_2(m))) \\ &= (F(p_1(a_1))...F(p_1(a_n))\omega)(p_2(m)) \end{aligned}$$

Согласно теореме 3.2.5, лемма является следствием (3.3.20), (3.3.22), (3.3.25).
⊙

Утверждение 3.3.5.9 является следствием определения 3.3.4 и лемм 3.3.6 и 3.3.7.

Диаграмма (6) является самым простым случаем в нашем доказательстве. Поскольку отображение $r_2$ является вложением и диаграмма (2) коммутативна, мы можем отождествить $n \in B_2$ и $t_2(m)$, если $n \in \mathrm{Im} t_2$. Аналогично, мы можем отождествить соответствующие преобразования.

$$(3.3.27) \qquad g'(r_1(t_1(a)))(r_2(t_2(m))) = r_2(G(t_1(a))(t_2(m)))$$

$$(g'(t_1(a_1))...g'(t_1(a_n))\omega)(t_2(m)) = r_2((G(t_1(a_1))...G(t_1(a_n))\omega)(t_2(m)))$$

Следовательно, $r = (r_1, r_2)$ является морфизмом представлений $G$ и $g$ (утверждение 3.3.5.10).

Для доказательства утверждения 3.3.5.11 осталось показать, что определённое в процессе доказательства представление $g'$ совпадает с представлением $g$, а операции над преобразованиями совпадают с соответствующими операциями на $\mathrm{End}(\Omega_2, B_2)$.

$$
\begin{aligned}
g'(r_1(t_1(a)))(r_2(t_2(m))) &= r_2(G(t_1(a))(t_2(m))) && \text{by (3.3.27)}\\
&= r_2(G(q_1(p_1(a)))(q_2(p_2(m)))) && \text{by (3.3.5), (3.3.6),}\\
&= r_2 \circ q_2(F(p_1(a))(p_2(m))) && \text{by (3.3.20)}\\
&= r_2 \circ q_2 \circ p_2(f(a)(m)) && \text{by (3.3.15)}\\
&= t_2(f(a)(m)) && \text{by (3.3.4), } i = 2\\
&= g(t_1(a))(t_2(m)) && \text{by (3.2.2)}
\end{aligned}
$$

$$
\begin{aligned}
(G(t_1(a_1))...G(t_1(a_n))\omega)(t_2(m)) &= q_2(F(p_1(a_1))...F(p_1(a_n))\omega)(p_2(m)))\\
&= q_2(F(p_1(a_1)...p_1(a_n)\omega)(p_2(m)))\\
&= q_2(F(p_1(a_1...a_n\omega))(p_2(m)))\\
&= q_2(p_2(f(a_1...a_n\omega)(m)))
\end{aligned}
$$

□

## 3.4. Приведенный морфизм представлений

Из теоремы 3.3.5 следует, что мы можем свести задачу изучения морфизма представлений $\Omega_1$-алгебры к случаю, описываемому диаграммой

$$(3.4.1)$$

$$
\begin{array}{ccc}
A_2 & \xrightarrow{p_2} & A_2/\ker t_2 \\
\uparrow{\scriptstyle f} & & \uparrow{\scriptstyle F} \\
A_1 & \xrightarrow{p_1} & A_1/\ker t_1
\end{array}
$$

Теорема 3.4.1. *Диаграмма* (3.4.1) *может быть дополнена представлением* $F_1$ $\Omega_1$-*алгебры* $A_1$ *в* $\Omega_2$-*алгебре* $A_2/\ker t_2$ *так, что диаграмма*

(3.4.2)

*коммутативна. При этом множество преобразований представления* $F$ *и множество преобразований представления* $F_1$ *совпадают.*

Доказательство. Для доказательства теоремы достаточно положить

$$F_1(a) = F(p_1(a))$$

Так как отображение $p_1$ - сюрьекция, то $\mathrm{Im}\,F_1 = \mathrm{Im}\,F$. Так как $p_1$ и $F$ - гомоморфизмы $\Omega_1$-алгебры, то $F_1$ - также гомоморфизм $\Omega_1$-алгебры.                    $\square$

Теорема 3.4.1 завершает цикл теорем, посвящённых структуре морфизма представлений $\Omega_1$-алгебры. Из этих теорем следует, что мы можем упростить задачу изучения морфизма представлений $\Omega_1$-алгебры и ограничиться морфизмом представлений вида

$$(\mathrm{id} : A_1 \to A_1,\ r_2 : A_2 \to B_2)$$

Определение 3.4.2. *Пусть*

$$f : A_1 \xrightarrow{\ \ *\ \ } A_2$$

*представление* $\Omega_1$-*алгебры* $A_1$ *в* $\Omega_2$-*алгебре* $A_2$ *и*

$$g : A_1 \xrightarrow{\ \ *\ \ } B_2$$

*представление* $\Omega_1$-*алгебры* $A_1$ *в* $\Omega_2$-*алгебре* $B_2$. *Пусть*

$$(\mathrm{id} : A_1 \to A_1,\ r_2 : A_2 \to B_2)$$

*морфизм представлений. В этом случае мы можем отождествить морфизм* $(\mathrm{id}, r_2)$ *представлений* $\Omega_1$-*алгебры и соответствующий гомоморфизм* $r_2$ $\Omega_2$-*алгебры и будем называть гомоморфизм* $r_2$ **приведенным морфизмом представлений**. *Мы будем пользоваться диаграммой*

(3.4.3)

*для представления приведенного морфизма $r_2$ представлений $\Omega_1$-алгебры. Из диаграммы следует*

$$(3.4.4) \qquad\qquad r_2 \circ f(a) = g(a) \circ r_2$$

*Мы будем также пользоваться диаграммой*

$$\begin{array}{ccc} A_2 & \xrightarrow{\ r_2\ } & B_2 \\ & & \\ {\scriptstyle f}\nwarrow\!{\scriptstyle *} & & {\scriptstyle *}\!\nearrow{\scriptstyle g} \\ & A_1 & \end{array}$$

*вместо диаграммы* (3.4.3). □

Теорема 3.4.3. *Пусть*

$$f : A_1 \longrightarrow_* A_2$$

*представление $\Omega_1$-алгебры $A_1$ в $\Omega_2$-алгебре $A_2$ и*

$$g : A_1 \longrightarrow_* B_2$$

*представление $\Omega_1$-алгебры $A_1$ в $\Omega_2$-алгебре $B_2$. Отображение*

$$r_2 : A_2 \to B_2$$

*является приведенным морфизмом представлений тогда и только тогда, когда*

$$(3.4.5) \qquad\qquad r_2(f(a)(m)) = g(a)(r_2(m))$$

Доказательство. Равенство (3.4.5) следует из равенства (3.4.4). □

Теорема 3.4.4. *Пусть отображение*

$$r_2 : A_2 \to B_2$$

*является приведенным морфизмом из представления*

$$f : A_1 \longrightarrow_* A_2$$

$\Omega_1$-*алгебры $A_1$ в представление*

$$g : A_1 \longrightarrow_* B_2$$

$\Omega_1$-*алгебры $A_1$. Если представление $f$ эффективно, то отображение*

$$r_2^* : \mathrm{End}(\Omega_2; A_2) \to \mathrm{End}(\Omega_2; B_2)$$

*определённое равенством*

$$(3.4.6) \qquad\qquad r_2^*(f(a)) = g(a)$$

*является гомоморфизмом $\Omega_1$-алгебры.*

Доказательство. Теорема является следствием теоремы 3.2.9, если мы положим $h = id$. □

Теорема 3.4.5. *Пусть представления*

$$f : A_1 \longrightarrow_* A_2$$

$$g : A_1 \longrightarrow_* B_2$$

$\Omega_1$-*алгебры* $A_1$ *однотранзитивны. Существует приведенный морфизм представлений из* $f$ *в* $g$

$$r_2 : A_2 \to B_2$$

Доказательство. Выберем элемент $m \in A_2$ и элемент $n \in B_2$. Чтобы построить отображение $r_2$, рассмотрим следующую диаграмму

Из коммутативности диаграммы (1) следует

$$r_2(f(a)(m)) = g(a)(r_2(m))$$

Для произвольного $m' \in A_2$ однозначно определён $a \in A_1$ такой, что $m' = f(a)(m)$. Следовательно, мы построили отображении $r_2$, которое удовлетворяет равенству (3.4.4). $\square$

Теорема 3.4.6. *Пусть представления*

$$f : A_1 \longrightarrow_* A_2$$

$$g : A_1 \longrightarrow_* B_2$$

$\Omega_1$-*алгебры* $A_1$ *однотранзитивны. Приведенный морфизм представлений из* $f$ *в* $g$

$$r_2 : A_2 \to B_2$$

*определён однозначно с точностью до выбора образа* $n = r_2(m) \in B_2$ *заданного элемента* $m \in A_2$.

Доказательство. Из доказательства теоремы 3.4.5 следует, что выбор элементов $m \in A_2$, $n \in B_2$ однозначно определяет отображение $r_2$. $\square$

Теорема 3.4.7. *Пусть*

$$f : A \longrightarrow_* B$$

*представление* $\Omega_1$-*алгебры* $A$ *в* $\Omega_2$-*алгебре* $B$. *Пусть* $N$ - *такая конгруэнция*[3.6] *на* $\Omega_2$-*алгебре* $B$, *что любое преобразование* $h \in \mathrm{End}(\Omega_2, B)$ *согласованно с конгруэнцией* $N$. *Существует представление*

$$f_1 : A \longrightarrow_* B/N$$

---

[3.6]Смотри определение конгруэнции на с. [14]-71.

$\Omega_1$-алгебры $A$ в $\Omega_2$-алгебре $B/N$ и отображение

$$\operatorname{nat} N : B \to B/N$$

*является приведенным морфизмом представления $f$ в представление $f_1$*

$$
\begin{array}{ccc}
B & \xrightarrow{\ j\ } & B/N \\
\end{array}
\qquad j = \operatorname{nat} N
$$

(с диаграммой: $B \xrightarrow{j} B/N$, $f: A \to B$, $f_1: A \to B/N$, $A$ внизу)

Доказательство. Любой элемент множества $B/N$ мы можем представить в виде $j(a)$, $a \in B$.

Согласно теореме [14]-II.3.5, мы можем определить единственную структуру $\Omega_2$-алгебры на множестве $B/N$. Если $\omega \in \Omega_2(p)$, то мы определим операцию $\omega$ на множестве $B/N$ согласно равенству (3) на странице [14]-73

$$(3.4.7) \qquad j(b_1)...j(b_p)\omega = j(b_1...b_p\omega)$$

Также как в доказательстве теоремы 3.3.5, мы можем определить представление

$$f_1 : A \longrightarrow B/N$$

с помощью равенства

$$(3.4.8) \qquad f_1(a) \circ j(b) = j(f(a) \circ b)$$

Равенство (3.4.8) можно представить с помощью диаграммы

$$(3.4.9)$$

$$
\begin{array}{ccc}
B & \xrightarrow{\ j\ } & B/N \\
\uparrow{\scriptstyle f(a)} & & \uparrow{\scriptstyle f_1(a)} \\
B & \xrightarrow{\ j\ } & B/N
\end{array}
$$

Пусть $\omega \in \Omega_2(p)$. Так как отображения $f(a)$ и $j$ являются гомоморфизмами $\Omega_2$-алгебры, то

$$
\begin{aligned}
f_1(a) \circ (j(b_1)...j(b_p)\omega) &= f_1(a) \circ j(b_1...b_p\omega) \\
&= j(f(a) \circ (b_1...b_p\omega)) \\
&= j((f(a) \circ b_1)...(f(a) \circ b_p)\omega) \\
&= j(f(a) \circ b_1)...j(f(a) \circ b_p)\omega \\
&= (f_1(a) \circ j(b_1))...(f_1(a) \circ j(b_p))\omega
\end{aligned}
$$

$$(3.4.10)$$

Из равенства (3.4.10) следует, что отображение $f_1(a)$ является гомоморфизмом $\Omega_2$-алгебры. Из равенства (3.4.8), согласно определению 3.4.2, следует, что отображение $j$ является приведенным морфизмом представления $f$ в представление $f_1$. $\qquad \square$

Теорема 3.4.8. *Пусть*

$$f : A \longrightarrow B$$

*представление $\Omega_1$-алгебры $A$ в $\Omega_2$-алгебре $B$. Пусть $N$ - такая конгруэнция на $\Omega_2$-алгебре $B$, что любое преобразование $h \in \operatorname{End}(\Omega_2, B)$ согласованно* *с конгруэнцией $N$. Рассмотрим категорию $\mathcal{A}$ объектами которой являются приведенные морфизмы представлений*[3.7]

$$R_1 : B \to S_1 \quad \ker R_1 \supseteq N$$

$$R_2 : B \to S_2 \quad \ker R_2 \supseteq N$$

*где $S_1$, $S_2$ - $\Omega_2$-алгебры и*

$$g_1 : A \dashrightarrow S_1 \qquad g_2 : A \dashrightarrow S_2$$

*представления $\Omega_1$-алгебры $A$. Мы определим морфизм $R_1 \to R_2$ как приведенный морфизм представлений $h : S_1 \to S_2$, для которого коммутативна диаграмма*

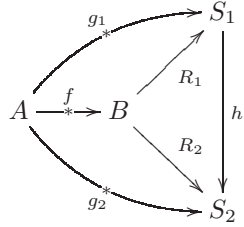

*Приведенный морфизм* $\operatorname{nat} N$ *представления $f$ в представление $f_1$ (теорема 3.4.7) является универсально отталкивающим в категории $\mathcal{A}$.*[3.8]

Доказательство. Существование и единственность отображения $h$, для которого коммутативна диаграмма

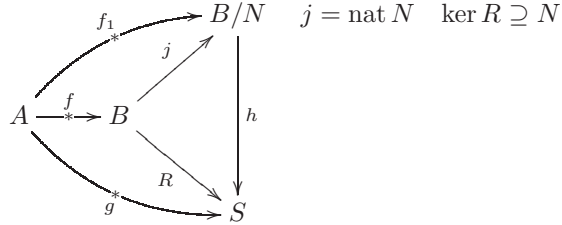

$$j = \operatorname{nat} N \quad \ker R \supseteq N$$

следует из теоремы 2.1.6. Следовательно, мы можем однозначно определить отображение $h$ с помощью равенства

$$(3.4.11) \qquad\qquad h(j(b)) = R(b)$$

Пусть $\omega \in \Omega_2(p)$. Так как отображения $R$ и $j$ являются гомоморфизмами $\Omega_2$-алгебры, то

$$(3.4.12) \qquad h(j(b_1)...j(b_p)\omega) = h(j(b_1...b_p\omega)) = R(b_1...b_p\omega) = R(b_1)...R(b_p)\omega$$
$$= h(j(b_1))...h(j(b_p))\omega$$

Из равенства (3.4.12) следует, что отображение $h$ является гомоморфизмом $\Omega_2$-алгебры.

---

[3.7] Утверждение леммы аналогично утверждению на странице [2]-94.

[3.8] Определение универсального объекта смотри в определении на с. [2]-47.

Так как отображение $R$ является приведенным морфизмом представления $f$ в представление $g$, то верно равенство

$$(3.4.13) \qquad g(a)(R(b)) = R(f(a)(b))$$

Из равенства (3.4.11) следует

$$(3.4.14) \qquad g(a)(h(j(b))) = g(a)(R(b))$$

Из равенств (3.4.13), (3.4.14) следует

$$(3.4.15) \qquad g(a)(h(j(b))) = R(f(a)(b))$$

Из равенств (3.4.11), (3.4.15) следует

$$(3.4.16) \qquad g(a)(h(j(b))) = h(j(f(a)(b)))$$

Из равенств (3.4.8), (3.4.16) следует

$$(3.4.17) \qquad g(a)(h(j(b))) = h(f_1(a)(j(b)))$$

Из равенства (3.4.17) следует, что отображение $h$ является приведенным морфизмом представления $f_1$ в представление $g$. $\qquad \square$

## 3.5. Автоморфизм представления универсальной алгебры

ОПРЕДЕЛЕНИЕ 3.5.1. *Пусть*

$$f : A_1 \longrightarrow_* A_2$$

*представление $\Omega_1$-алгебры $A_1$ в $\Omega_2$-алгебре $A_2$. Приведенный морфизм представлений $\Omega_1$-алгебры*

$$r_2 : A_2 \to A_2$$

*такой, что $r_2$ - эндоморфизм $\Omega_2$-алгебры называется* **эндоморфизмом представления** $f$. $\qquad \square$

ТЕОРЕМА 3.5.2. *Если представление*

$$f : A_1 \longrightarrow_* A_2$$

*$\Omega_1$-алгебры $A_1$ однотранзитивно, то для любых $a_{21}$, $a_{22} \in A_2$ существует единственный эндоморфизм*

$$r_2 : A_2 \to A_2$$

*представления $f$ такой, что $r_2(a_{21}) = a_{22}$.*

ДОКАЗАТЕЛЬСТВО. Рассмотрим следующую диаграмму

Существование эндоморфизма является следствием теоремы 3.2.10. Единственность эндоморфизма для заданных $p$, $q \in A_2$ является следствием теоремы 3.2.11, когда $r_1 = \text{id}$. $\qquad \square$

ТЕОРЕМА 3.5.3. *Эндоморфизмы представления $f$ порождают полугруппу.*

ДОКАЗАТЕЛЬСТВО. Из теоремы 3.3.1 следует, что произведение эндоморфизмов $(id, p_2)$, $(id, r_2)$ представления $f$ является эндоморфизмом $(id, p_2 \circ r_2)$ представления $f$. $\qquad \square$

Определение 3.5.4. *Пусть*

$$f : A_1 \xrightarrow{\;*\;} A_2$$

*представление $\Omega_1$-алгебры $A_1$ в $\Omega_2$-алгебре $A_2$. Морфизм представлений $\Omega_1$-алгебры*

$$r_2 : A_2 \to A_2$$

*такой, что $r_2$ - автоморфизм $\Omega_2$-алгебры называется* **автоморфизмом представления** *$f$.*  □

Теорема 3.5.5. *Пусть*

$$f : A_1 \xrightarrow{\;*\;} A_2$$

*представление $\Omega_1$-алгебры $A_1$ в $\Omega_2$-алгебре $A_2$. Множество автоморфизмов представления $f$ порождает* **группу** $GA(f)$.

Доказательство. Пусть $r_2$, $p_2$ - автоморфизмы представления $f$. Согласно определению 3.5.4, отображения $r_2$, $p_2$ являются автоморфизмами $\Omega_2$-алгебры $A_2$. Согласно теореме II.3.2, ([14], с. 60), отображение $r_2 \circ p_2$ является автоморфизмом $\Omega_2$-алгебры $A_2$. Из теоремы 3.3.1 и определения 3.5.4 следует, что произведение автоморфизмов $r_2 \circ p_2$ представления $f$ является автоморфизмом представления $f$.

Пусть $r_2, p_2, q_2$ - автоморфизмы представления $f$. Из цепочки равенств

$$((r_2 \circ p_2) \circ q_2)(a) = (r_2 \circ p_2)(q_2(a)) = r_2(p_2(q_2(a)))$$
$$= r_2((p_2 \circ q_2)(a)) = (r_2 \circ (p_2 \circ q_2))(a)$$

следует ассоциативность произведения для отображений [3.9] $r_2$, $p_2$, $q_2$.

Пусть $r_2$ - автоморфизм представления $f$. Согласно определению 3.5.4 отображение $r_2$ является автоморфизмом $\Omega_2$-алгебры $A_2$. Следовательно, отображение $r_2^{-1}$ является автоморфизмом $\Omega_2$-алгебры $A_2$. Для автоморфизма $r_2$ представления справедливо равенство (3.2.3). Положим $m' = r_2(m)$. Так как $r_2$ - автоморфизм $\Omega_2$-алгебры, то $m = r_2^{-1}(m')$ и равенство (3.2.3) можно записать в виде

$$(3.5.1) \qquad r_2(f(a')(r_2^{-1}(m'))) = f(a')(m')$$

Так как отображение $r_2$ является автоморфизмом $\Omega_2$-алгебры $A_2$, то из равенства (3.5.1) следует

$$(3.5.2) \qquad f(a')(r_2^{-1}(m')) = r_2^{-1}(f(a')(m'))$$

Равенство (3.5.2) соответствует равенству (3.2.3) для отображения $r_2^{-1}$. Следовательно, отображение $r_2^{-1}$ является автоморфизмом представления $f$.  □

---

[3.9] При доказательстве ассоциативности произведения я следую примеру полугруппы из [5], с. 20, 21.

# $\Omega$-группа

## 4.1. Множество гомоморфизмов $\Omega$-алгебры

Теорема 4.1.1. *Пусть множества $A$, $B$ являются $\Omega$-алгебрами. Множество* $\mathrm{Hom}(\Omega; A \to B)$ *является $\Omega$-алгеброй, если для любых операций* $\omega_1 \in \Omega(m)$, $\omega_2 \in \Omega(n)$, *верно следующее равенство*

$$(4.1.1) \qquad (a_{11}...a_{1n}\omega_2)...(a_{m1}...a_{mn}\omega_2)\omega_1 = (a_{11}...a_{m1}\omega_1)...(a_{1n}...a_{mn}\omega_1)\omega_2$$

Доказательство. Согласно теореме 2.2.6, множество $B^A$ является $\Omega$-алгеброй. Пусть $\omega \in \Omega(n)$. Для отображений $f_1$, ..., $f_n \in B^A$, мы определим операцию $\omega$ равенством

$$(4.1.2) \qquad (f_1...f_n\omega)(x) = f_1(x)...f_n(x)\omega$$

Пусть $\omega_1 \in \Omega(m)$, $\omega_2 \in \Omega(n)$. Пусть отображения $f_1$, ..., $f_m \in \mathrm{Hom}(\Omega; A \to B)$ являются гомоморфизмами $\Omega$-алгебры $A$ в $\Omega$-алгебру $B$. В частности, для любых $a_1$, ..., $a_n \in A$

$$(4.1.3) \qquad \begin{aligned} f_1(a_1...a_n\omega_2) &= f_1(a_1)...f_1(a_n)\omega_2 \\ ... &= ... \\ f_m(a_1...a_n\omega_2) &= f_m(a_1)...f_m(a_n)\omega_2 \end{aligned}$$

Если мы требуем, что отображение $f_1...f_m\omega_1$ является гомоморфизмом $\Omega$-алгебры $A$ в $\Omega$-алгебру $B$, то

$$(4.1.4) \qquad (f_1...f_m\omega_1)(a_1...a_n\omega_2) = ((f_1...f_m\omega_1)(a_1))...((f_1...f_m\omega_1)(a_n))\omega_2$$

Согласно определению (4.1.2), равенство

$$(4.1.5) \qquad \begin{aligned} &f_1(a_1...a_n\omega_2)...f_m(a_1...a_n\omega_2)\omega_1 \\ &= (f_1(a_1)...f_m(a_1)\omega_1)...(f_1(a_n)...f_m(a_n)\omega_1)\omega_2 \end{aligned}$$

является следствием равенства (4.1.4). Равенство

$$(4.1.6) \qquad \begin{aligned} &(f_1(a_1)...f_1(a_n)\omega_2)...(f_m(a_1)...f_m(a_n)\omega_2)\omega_1 \\ &= (f_1(a_1)...f_m(a_1)\omega_1)...(f_1(a_n)...f_m(a_n)\omega_1)\omega_2 \end{aligned}$$

является следствием равенств (4.1.3), (4.1.5). Положим

$$(4.1.7) \qquad a_{ij} = f_i(a_j)$$

Равенство (4.1.1) является следствием равенств (4.1.6), (4.1.7). $\qquad\square$

Не всякая $\Omega$-алгебра удовлетворяет условиям теоремы 4.1.1.

Теорема 4.1.2. *Если $G_1$, $G_2$ - абелевые полугруппы, то множество* $\mathrm{Hom}(\{+\}; G_1 \to G_2)$ *также является абелевой полугруппой.*

Доказательство. Поскольку операция сложения в абелевой полугруппе коммутативна и ассоциативна, то теорема является следствием теоремы 4.1.1. □

Теорема 4.1.3. *Множество* $\mathrm{End}(\{+\}; A)$ *эндоморфизмов абелевой группы* $A$ *является абелевой группой.*

Доказательство. Теорема является следствием теорем 2.2.13, 4.1.2 и утверждения, что уравнение

$$x + a = 0$$

в абелевой группе имеет решение. □

Теорема 4.1.4. *Если* $D_1$, $D_2$ - *кольца, то множество* $\mathrm{Hom}(\{+, *\}; D_1 \to D_2)$, *вообще говоря, кольцом не является.*

Доказательство. В кольце определены две операции: сложение, которое коммутативно и ассоциативно, и произведение, которое дистрибутивно относительно сложения. Согласно теореме 4.1.1, сложение и произведение должны удовлетворять равенству

$$(4.1.8) \qquad a_{11}a_{21} + a_{12}a_{22} = (a_{11} + a_{12})(a_{21} + a_{22})$$

Однако правая часть равенства (4.1.8) имеет вид

$$(a_{11} + a_{12})(a_{21} + a_{22}) = (a_{11} + a_{12})a_{21} + (a_{11} + a_{12})a_{22}$$
$$= a_{11}a_{21} + a_{12}a_{21} + a_{11}a_{22} + a_{12}a_{22}$$

Следовательно, равенство (4.1.8) не верно. □

Анализ теорем 4.1.2, 4.1.4 говорит о том, что множество Ω-алгебр, удовлетворяющих условиям теоремы 4.1.1, невелико.

Вопрос 4.1.5. *Существует ли универсальная алгебра, отличная от абелевой полугруппы и удовлетворяющая условиям теоремы 4.1.1?* □

Из нашего опыта следует, что многие Ω-алгебры содержат операцию, которая соло порождала бы полугруппу. Поэтому мы изменим формулировку теоремы 4.1.1.

Теорема 4.1.6. *Пусть множества* $A$, $B$ *являются* Ω-*алгебрами. Пусть* $\omega \in \Omega(n)$. *Множество* $\mathrm{Hom}(\Omega; A \to B)$ *замкнуто относительно операции* $\omega$, *если верно следующее равенство*

$$(4.1.9) \qquad (a_{11}...a_{1n}\omega)...(a_{n1}...a_{nn}\omega)\omega$$
$$= (a_{11}...a_{n1}\omega)...(a_{1n}...a_{nn}\omega)\omega$$

Доказательство. Вообще говоря, мы рассматриваем множество $\mathrm{Hom}(\{\omega\}; A \to B)$. Теорема является следствием теоремы 4.1.1. □

Теорема 4.1.7. *Пусть операция* $\omega \in \Omega(2)$ *коммутативна и ассоциативна. Множество* $\mathrm{Hom}(\Omega; A \to B)$ *замкнуто относительно операции* $\omega$.

Доказательство. Так как операция $\omega \in \Omega(2)$ коммутативна и ассоциативна, то

$$(4.1.10)\quad \begin{aligned}(a_{11}a_{12}\omega)(a_{21}a_{22}\omega)\omega &= a_{11}(a_{12}(a_{21}a_{22}\omega)\omega)\omega\\ &= a_{11}((a_{12}a_{21}\omega)a_{22}\omega)\omega\\ &= a_{11}((a_{21}a_{12}\omega)a_{22}\omega)\omega\\ &= a_{11}(a_{21}(a_{12}a_{22}\omega)\omega)\omega\\ &= (a_{11}a_{21}\omega)(a_{12}...a_{22}\omega)\omega\end{aligned}$$

Теорема является следствием равенства (4.1.10) и теоремы 4.1.6.  □

Теорема 4.1.8. *Пусть операция* $\omega \in \Omega(2)$ *имеет нейтральный элемент и множество* $\mathrm{Hom}(\Omega; A \to B)$ *замкнуто относительно операции* $\omega$. *Тогда операция* $\omega$ *коммутативна и ассоциативна.*

Доказательство. Равенства

$$(4.1.11)\qquad ab\omega = (ea\omega)(be\omega)\omega = (eb\omega)(ae\omega)\omega = ba\omega$$

$$(4.1.12)\qquad a(bc\omega) = (ae\omega)(bc\omega)\omega = (ab\omega)(ec\omega)\omega = (ab\omega)c\omega$$

являются следствием равенств (2.4.1), (2.4.2), (4.1.9). Ассоциативность операции $\omega$ является следствием равенства (4.1.11). Коммутативность операции $\omega$ является следствием равенства (4.1.12).  □

Вопрос 4.1.9. *Существует ли область операторов* $\Omega$, *для которой верны следующие утверждения?*

- *Множество* $\mathrm{Hom}(\Omega; A \to B)$ *замкнуто относительно операции* $\omega \in \Omega(2)$.
- *Операция* $\omega$ *не является коммутативной или ассоциативной.*

□

## 4.2. Ω-группа

Пусть в $\Omega_2$-алгебре $A_2$ определена операция $\omega \in \Omega_2(2)$, которая коммутативна и ассоциативна. Мы будем отождествлять операцию $\omega$ с суммой. Мы пользуемся символом $+$ для обозначения операции суммы. Положим

$$\Omega = \Omega_2 \setminus \{+\}$$

Определение 4.2.1. *Отображение*

$$f : A_2 \to B_2$$

$\Omega_2$-*алгебры* $A_2$ *в* $\Omega_2$-*алгебру* $B_2$ *называется* **аддитивным отображением**, *если*

$$f(a + b) = f(a) + f(b)$$

*Обозначим* $\mathcal{A}(A_2 \to B_2)$ *множество аддитивных отображений* $\Omega_2$-*алгебры* $A_2$ *в* $\Omega_2$-*алгебру* $B_2$.  □

Теорема 4.2.2. $\mathcal{A}(A_2 \to B_2) = \mathrm{Hom}(\{+\}; A_2 \to B_2)$.

Доказательство. Теорема является следствием определений 2.2.9, 4.2.1.  □

Определение 4.2.3. *Отображение*

$$g : A^n \to A$$

*называется* **полиаддитивным отображением**, *если для любого* $i$, $i = 1, ..., n$,

$$f(a_1, ..., a_i + b_i, ..., a_n) = f(a_1, ..., a_i, ..., a_n) + f(a_1, ..., b_i, ..., a_n)$$

□

Теорема 4.2.4. *Пусть отображение*

$$f : A_1 \xrightarrow{\quad*\quad} A_2$$

*является эффективным представлением Ω-алгебры $A_1$ в абелевой полугруппе $A_2$.*

4.2.4.1: *На множестве $A_1$ можно определить структуру абелевой полугруппы*

(4.2.1)                    $$f(a_1 + b_1)(a_2) = f(a_1)(a_2) + f(b_1)(a_2)$$

4.2.4.2: *Отображение $f$ является аддитивным отображением.*

4.2.4.3: *Отображение $f$ является представлением $\Omega_1$-алгебры $A_1$, где $\Omega_1 = \Omega \cup \{+\}$.*

Доказательство. Согласно теоремам 2.2.13, 4.1.7, множество $\text{End}(\{+\}, A_2)$ является абелевой полугруппой. Поскольку представление $f$ эффективно, то, согласно теоремам 3.1.3, 4.1.1, для любых $A_1$-чисел $a$, $b$ существует единственное $A_1$-число $c$ такое, что

(4.2.2)                    $$f(c)(m) = f(a)(m) + f(b)(m)$$

Опираясь на равенство (4.2.2), мы определяем сумму $A_1$-чисел

(4.2.3)                    $$c = a + b$$

Равенство (4.2.1) является следствием равенств (4.2.2), (4.2.3).

Лемма 4.2.5. *Сумма $A_1$-чисел коммутативна.*

Доказательство. Поскольку сумма $A_2$-чисел коммутативна, то равенство

(4.2.4)      $$\begin{aligned} f(a_1 + b_1)(a_2) &= f(a_1)(a_2) + f(b_1)(a_2) = f(b_1)(a_2) + f(a_1)(a_2) \\ &= f(b_1 + a_1)(a_2) \end{aligned}$$

является следствием равенства (4.2.1). Лемма является следствием равенства (4.2.4).                                                                          ⊙

Лемма 4.2.6. *Сумма $A_1$-чисел ассоциативна.*

Доказательство. Поскольку сумма $A_2$-чисел ассоциативна, то равенство

(4.2.5)   $$\begin{aligned} f((a_1 + b_1) + c_1)(a_2) &= f(a_1 + b_1)(a_2) + f(c_1)(a_2) \\ &= (f(a_1)(a_2) + f(b_1)(a_2)) + f(c_1)(a_2) \\ &= f(a_1)(a_2) + (f(b_1)(a_2) + f(c_1)(a_2)) \\ &= f(a_1)(a_2) + f(b_1 + c_1)(a_2) \\ &= f(a_1 + (b_1 + c_1))(a_2) \end{aligned}$$

является следствием равенства (4.2.1). Лемма является следствием равенства (4.2.5). ⊙

Утверждение 4.2.4.1 является следствием равенства (4.2.3), лемм 4.2.5, 4.2.6 и определения 2.4.4.

Утверждение 4.2.4.2 является следствием равенства (4.2.3). Утверждение 4.2.4.3 является следствием утверждения 4.2.4.2, так как отображение $f$ является гомоморфизмом Ω-алгебры. □

Теорема 4.2.7. *Пусть* $\omega \in \Omega(n)$, $\omega_1 \in \Omega(m)$. *Отображение*

$$(4.2.6) \qquad g : a_i \to a_1...a_n\omega$$

*согласовано с операцией* $\omega_1$, *если верно следующее равенство*

$$(4.2.7) \qquad a_1...(a_{i1}...a_{im}\omega_1)...a_n\omega = (a_1...a_{i1}...a_n\omega)...(a_1...a_{im}...a_n\omega)\omega_1$$

Доказательство. Равенство

$$(4.2.8) \qquad \begin{aligned} g(a_{i1}...a_{im}\omega_1) &= a_1...(a_{i1}...a_{im}\omega_1)...a_n\omega \\ &= (a_1...a_{i1}...a_n\omega)...(a_1...a_{im}...a_n\omega)\omega_1 \\ &= g(a_{i1})...g(a_{im})\omega_1 \end{aligned}$$

является следствием равенств (4.2.7), (4.2.6). Теорема является следствием определения 2.2.9 и равенства (4.2.8). □

Равенство (4.2.7) является менее жёстким, чем равенство (4.1.1). Тем не менее, также как и в случае теоремы 4.1.1, большинство операций универсальной алгебры не удовлетворяет условиям теоремы 4.2.7. Поскольку операция сложения удовлетворяет условиям теоремы 4.1.1, мы ожидаем, что существуют условия, когда операция сложения удовлетворяет условиям теоремы 4.2.7.

Теорема 4.2.8. *Пусть* $\omega \in \Omega(n)$. *Если отображение*

$$(4.2.9) \qquad g : a_i \to a_1...a_n\omega$$

*согласовано со сложением для любого* $i$, *то операция* $\omega$ *является полиаддитивным отображением.*

Доказательство. Согласно теореме 4.2.7, если отображение (4.2.9) согласовано со сложением, то верно следующее равенство

$$(4.2.10) \qquad a_1...(a_{i1} + a_{i2})...a_n\omega = (a_1...a_{i1}...a_n\omega) + (a_1...a_{i2}...a_n\omega)$$

Теорема является следствием равенства (4.2.10) и определения 4.2.3. □

Теорема 4.2.9. *Пусть* $\omega \in \Omega(n)$ *- полиаддитивное отображение. Операция* $\omega$ **дистрибутивна** *относительно сложения*

$$a_1...(a_i + b_i)...a_n\omega = a_1...a_i...a_n\omega + a_1...b_i...a_n\omega \quad i = 1, ..., n$$

Доказательство. Теорема является следствием теоремы 4.2.8. □

Определение 4.2.10. *Пусть в* $\Omega_1$-*алгебре* $A$ *определена операция сложения, которая не обязательно коммутативна. Мы пользуемся символом* + *для обозначения операции суммы. Положим*

$$\Omega = \Omega_1 \setminus \{+\}$$

*Если* $\Omega_1$-*алгебра* $A$ *является группой относительно операции сложения и любая операция* $\omega \in \Omega$ *является полиаддитивным отображением, то* $\Omega_1$-*алгебра*

*A называется* Ω-**группой**. *Если* Ω-*группа A является ассоциативной группой относительно операции сложения, то* $\Omega_1$-*алгебра A называется* **ассоциативной** Ω-**группой**. *Если* Ω-*группа A является абелевой группой относительно операции сложения, то* $\Omega_1$-*алгебра A называется* **абелевой** Ω-**группой**.    □

Пример 4.2.11. *Группа является наиболее очевидным примером* Ω-*группы.*

*Кольцо является* Ω-*группой.*

*Бикольцо матриц над телом ([8]) является* Ω-*группой.*    □

Замечание 4.2.12. *Бурбаки рассматривают похожее определение, а именно группы с операторами (смотри определение 10 в [16] на странице 100).*    □

Теорема 4.2.13. *Пусть A - Ω-группа. Пусть* $\omega \in \Omega(n)$. *Отображение*

$$g : a_i \to a_1...a_n\omega$$

*является эндоморфизмом аддитивной группы A.*

Доказательство. Теорема является следствием теоремы 4.2.9 и определения 4.2.10.    □

Теорема 4.2.14. *Пусть отображение*

$$g : A_1 \longrightarrow{}^* A_2$$

*является представлением* Ω-*группы* $A_1$. *Тогда отображение*

$$\left( a_i \to a_1...a_n\omega \quad f(a_i) \to f(a_1)...f(a_n)\omega \right)$$

*является морфизмом представления f аддитивной группы* $A_1$.

Доказательство. Теорема является следствием теоремы 4.2.13 и определений 3.1.1, 3.2.2.    □

### 4.3. Декартово произведение представлений

Лемма 4.3.1. *Пусть*

$$A = \prod_{i \in I} A_i$$

*декартово произведение семейства* $\Omega_2$-*алгебр* $(A_i, i \in I)$ . *Для каждого* $i \in I$, *пусть множество* $\text{End}(\Omega_2, A_i)$ *является* $\Omega_1$-*алгеброй. Тогда множество*

(4.3.1)        $${}^{\circ}A = \{ f \in \text{End}(\Omega_2; A) : f(a_i, i \in I) = (f_i(a_i), i \in I) \}$$

*является декартовым произведением* $\Omega_1$-*алгебр* $\text{End}(\Omega_2, A_i)$.

Доказательство. Согласно определению (4.3.1), мы можем представить отображение $f \in {}^{\circ}A$ в виде кортежа

$$f = (f_i, i \in I)$$

отображений $f_i \in \text{End}(\Omega_2; A_i)$. Согласно определению (4.3.1),

$$(f_i, i \in I)(a_i, i \in I) = (f_i(a_i), i \in I)$$

Пусть $\omega \in \Omega_2$ - n-арная операция. Мы определим операцию $\omega$ на множестве ${}^{\circ}A$ равенством

$$((f_{1i}, i \in I)...(f_{ni}, i \in I)\omega)(a_i, i \in I) = ((f_{1i}(a_i))...(f_{ni}(a_i))\omega, i \in I)$$

□

Определение 4.3.2. *Пусть $\mathcal{A}_1$ - категория $\Omega_1$-алгебр. Пусть $\mathcal{A}_2$ - категория $\Omega_2$-алгебр. Мы определим* **категорию** $\mathcal{A}_1(\mathcal{A}_2)$ **представлений**. *Объектами этой категории являются представления $\Omega_1$-алгебры в $\Omega_2$-алгебре. Морфизмами этой категории являются морфизмы соответствующих представлений.* □

Теорема 4.3.3. *В категории $\mathcal{A}_1(\mathcal{A}_2)$ существует произведение однотранзитивных представлений $\Omega_1$-алгебры в $\Omega_2$-алгебре.*

Доказательство. Для $j = 1, 2$, пусть

$$P_j = \prod_{i \in I} B_{ji}$$

произведение семейства $\Omega_j$-алгебр $\{B_{ji}, i \in I\}$ и для любого $i \in I$ отображение

$$t_{ji} : P_j \longrightarrow B_{ji}$$

является проекцией на множитель $i$. Для каждого $i \in I$, пусть

$$h_i : B_{1i} \dashrightarrow B_{2i}$$

однотранзитивное $B_{1i}$-представление в $\Omega_2$-алгебре $B_{2i}$.

Пусть $b_1 \in P_1$. Согласно утверждению 2.3.3.3, $P_1$-число $b_1$ может быть представлено в виде кортежа $B_{1i}$-чисел

(4.3.2) $\qquad b_1 = (b_{1i}, i \in I) \quad b_{1i} = t_{1i}(b_1) \in B_{1i}$

Пусть $b_2 \in P_2$. Согласно утверждению 2.3.3.3, $P_2$-число $b_2$ может быть представлено в виде кортежа $B_{2i}$-чисел

(4.3.3) $\qquad b_2 = (b_{2i}, i \in I) \quad b_{2i} = t_{2i}(b_2) \in B_{2i}$

Лемма 4.3.4. *Для каждого $i \in I$, рассмотрим диаграмму отображений*

(4.3.4)

*Пусть отображение*

$$g : P_1 \to \mathrm{End}(\Omega_2; P_2)$$

*определено равенством*

(4.3.5) $\qquad g(b_1)(b_2) = (h_i(b_{1i})(b_{2i}), i \in I)$

*Тогда отображение $g$ является однотранзитивным $P_1$-представлением в $\Omega_2$-алгебре $P_2$*

$$g : P_1 \dashrightarrow P_2$$

*Отображение $(t_{1i}, t_{2i})$ является морфизмом представления $g$ в представление $h_i$.*

Доказательство.

4.3.4.1: Согласно определениям 3.1.1, отображение $h_i(b_{1i})$ является гомоморфизмом $\Omega_2$-алгебры $B_{2i}$. Согласно теореме 2.3.6, из коммутативности диаграммы (1) для каждого $i \in I$, следует, что отображение

$$g(b_1) : P_2 \to P_2$$

определённое равенством (4.3.5) является гомоморфизмом $\Omega_2$-алгебры $P_2$.

4.3.4.2: Согласно определению 3.1.1, множество $\operatorname{End}(\Omega_2; B_{2i})$ является $\Omega_1$-алгеброй. Согласно лемме 4.3.1, множество $\,^\circ P_2 \subseteq \operatorname{End}(\Omega_2; P_2)$ является $\Omega_1$-алгеброй.

4.3.4.3: Согласно определению 3.1.1, отображение

$$h_i : B_{1i} \to \operatorname{End}(\Omega_2; B_{2i})$$

является гомоморфизмом $\Omega_1$-алгебры. Согласно теореме 2.3.6, отображение

$$g : P_1 \to \operatorname{End}(\Omega_2; P_2)$$

определённое равенством

$$g(b_1) = (h_i(b_{1i}), i \in I)$$

является гомоморфизмом $\Omega_1$-алгебры.

Согласно утверждениям 4.3.4.1, 4.3.4.3 и определению 3.1.1, отображение $g$ является $P_1$-представлением в $\Omega_2$-алгебре $P_2$.

Пусть $b_{21}$, $b_{22} \in P_2$. Согласно утверждению 2.3.3.3, $P_2$-числа $b_{21}$, $b_{22}$ могут быть представлены в виде кортежей $B_{2i}$-чисел

$$(4.3.6) \qquad \begin{aligned} b_{21} &= (b_{21i}, i \in I) \quad b_{21i} = t_{2i}(b_{21}) \in B_{2i} \\ b_{22} &= (b_{22i}, i \in I) \quad b_{22i} = t_{2i}(b_{22}) \in B_{2i} \end{aligned}$$

Согласно теореме 3.1.9, поскольку представление $h_i$ однотранзитивно, то существует единственное $B_{1i}$-число $b_{1i}$ такое, что

$$b_{22i} = h_i(b_{1i})(b_{21i})$$

Согласно определениям (4.3.2), (4.3.5), (4.3.6), существует единственное $P_1$-число $b_1$ такое, что

$$b_{22} = g(b_1)(b_{21})$$

Согласно теореме 3.1.9, представление $g$ однотранзитивно.

Из коммутативности диаграммы (1) и определения 3.2.2, следует, что отображение $(t_{1i}, t_{2i})$ является морфизмом представления $g$ в представление $h_i$.                                                                                      ⊙

Пусть

$$(4.3.7) \qquad\qquad d_2 = g(b_1)(b_2) \quad d_2 = (d_{2i}, i \in I)$$

Из равенств (4.3.5), (4.3.7) следует, что

$$(4.3.8) \qquad\qquad d_{2i} = h_i(b_{1i})(b_{2i})$$

Для $j = 1, 2$, пусть $R_j$ - другой объект категории $\mathcal{A}_j$. Для любого $i \in I$, пусть отображение

$$r_{1i} : R_1 \longrightarrow B_{1i}$$

является морфизмом из $\Omega_1$-алгебра $R_1$ в $\Omega_1$-алгебру $B_{1i}$. Согласно определению 2.3.1, существует единственный морфизм $\Omega_1$-алгебры

$$s_1 : R_1 \longrightarrow P_1$$

такой, что коммутативна диаграмма

(4.3.9)

$$t_{1i}(s_1) = r_{1i}$$

Пусть $a_1 \in R_1$. Пусть

(4.3.10)                     $b_1 = s_1(a_1) \in P_1$

Из коммутативности диаграммы (4.3.9) и утверждений (4.3.10), (4.3.2) следует, что

(4.3.11)                     $b_{1i} = r_{1i}(a_1)$

Пусть

$$f : R_1 \dashrightarrow R_2$$

однотранзитивное $R_1$-представление в $\Omega_2$-алгебре $R_2$. Согласно теореме 3.2.11, морфизм $\Omega_2$-алгебры

$$r_{2i} : R_2 \longrightarrow B_{2i}$$

такой, что отображение $(r_{1i}, r_{2i})$ является морфизмом представлений из $f$ в $h_i$, определён однозначно с точностью до выбора образа $R_2$-числа $a_2$. Согласно замечанию 3.2.7, в диаграмме отображений

(4.3.12)

диаграмма (2) коммутативна. Согласно определению 2.3.1, существует единственный морфизм $\Omega_2$-алгебры

$$s_2 : R_2 \longrightarrow P_2$$

такой, что коммутативна диаграмма

(4.3.13)

$$P_2 \xrightarrow{t_{2i}} B_{2i} \qquad t_{2i}(s_2) = r_{2i}$$

with $s_2$ from $R_2$ and $r_{2i}$

Пусть $a_2 \in R_2$. Пусть

(4.3.14)                         $b_2 = s_2(a_2) \in P_2$

Из коммутативности диаграммы (4.3.13) и утверждений (4.3.14), (4.3.3) следует, что

(4.3.15)                         $b_{2i} = r_{2i}(a_2)$

Пусть

(4.3.16)                         $c_2 = f(a_1)(a_2)$

Из коммутативности диаграммы (2) и равенств (4.3.8), (4.3.15), (4.3.16) следует, что

(4.3.17)                         $d_{2i} = r_{2i}(c_2)$

Из равенств (4.3.8), (4.3.17) следует, что

(4.3.18)                         $d_2 = s_2(c_2)$

что согласуется с коммутативносью диаграммы (4.3.13).

Для каждого $i \in I$, мы объединим диаграммы отображений (4.3.4), (4.3.9), (4.3.13), (4.3.12)

Из равенств (4.3.7) (4.3.14) и из равенств (4.3.16), (4.3.18), следует коммутативность диаграммы (3). Следовательно, отображение $(s_1, s_2)$ является морфизмом представлений из $f$ в $g$. Согласно теореме 3.2.11, морфизм $(s_1, s_2)$ определён однозначно, так как мы требуем (4.3.18).

Согласно определению 2.3.1, представление $g$ и семейство морфизмов представления $((t_{1i}, t_{2i}), i \in I)$ является произведением в категории $\mathcal{A}_1(\mathcal{A}_2)$.  $\square$

Определение 4.3.5. *Пусть $A_1$, ..., $A_n$, $A$ - $\Omega_1$-алгебры. Пусть $B_1$, ..., $B_n$, $B$ - $\Omega_2$-алгебры. Пусть, для любого $k$, $k = 1$, ..., $n$,*

$$f_k : A_k \longrightarrow_* B_k$$

*представление $\Omega_1$-алгебры $A_k$ в $\Omega_2$-алгебре $B_k$. Пусть*

$$f : A \longrightarrow_* B$$

*представление $\Omega_1$-алгебры $A$ в $\Omega_2$-алгебре $B$. Отображение*

$$\Big( r_{1k} : A_k \to A \quad k = 1, ..., n \quad r_2 : B_1 \times ... \times B_n \to B \Big)$$

*называется* **полиморфизмом представлений** *$f_1$, ..., $f_n$ в представление $f$, если, для любого $k$, $k = 1$, ..., $n$, при условии, что все переменные кроме переменных $a_k \in A_k$, $b_k \in B_k$ имеют заданное значение, отображение $(r_{1k}, r_2)$ является морфизмом представления $f_k$ в представление $f$.*

*Если $f_1 = ... = f_n$, то мы будем говорить, что отображение $((r_{1,k}, k = 1, ..., n)$ $r_2)$ является полиморфизмом представления $f_1$ в представление $f$.*

*Если $f_1 = ... = f_n = f$, то мы будем говорить, что отображение $((r_{1,k}, k = 1, ..., n)$ $r_2)$ является полиморфизмом представления $f$.* □

Мы также будем говорить, что отображение $r = (r_1, r_2)$ является полиморфизмом представлений в $\Omega_2$-алгебрах $B_1$, ..., $B_n$ в представление в $\Omega_2$-алгебре $B$.

Теорема 4.3.6. *Пусть отображение $((r_{1,k}, k = 1, ..., n)$ $r_2)$ является полиморфизмом представлений $f_1$, ..., $f_n$ в представление $f$. Для любого $k$, $k = 1$, ..., $n$, отображение $(r_{1k}, r_2)$ удовлетворяет равенству*

$$(4.3.19) \qquad r_2(m_1, ..., f_k(a_k)(m_k), ..., m_n) = f(r_{1k}(a_k))(r_2(m_1, ..., m_n))$$

*Пусть $\omega_1 \in \Omega_1(p)$. Для любого $k$, $k = 1$, ..., $n$, отображение $r_{1k}$ удовлетворяет равенству*

$$(4.3.20) \qquad r_{1k}(a_{k \cdot 1}...a_{k \cdot p}\omega_1) = r_{1k}(a_{k \cdot 1})...r_{1k}(a_{k \cdot p})\omega_1$$

*Пусть $\omega_2 \in \Omega_2(p)$. Для любого $k$, $k = 1$, ..., $n$, отображение $r_2$ удовлетворяет равенству*

$$(4.3.21) \qquad r_2(m_1, ..., m_{k \cdot 1}...m_{k \cdot p}\omega_2, ..., m_n)$$
$$= r_2(m_1, ..., m_{k \cdot 1}, ..., m_n)...r_2(m_1, ..., m_{k \cdot p}, ..., m_n)\omega_2$$

Доказательство. Равенство (4.3.19) следует из определения 4.3.5 и равенства (2.3.3). Равенство (4.3.20) следует из утверждения, что, для любого $k$, $k = 1$, ..., $n$, при условии, что все переменные кроме переменной $x_k \in A_k$ имеют заданное значение, отображение $r_1$ является гомоморфизмом $\Omega_1$-алгебры $A_k$ в $\Omega_1$-алгебру $A$. Равенство (4.3.21) следует из утверждения, что, для любого $k$, $k = 1$, ..., $n$, при условии, что все переменные кроме переменной $m_k \in B_k$ имеют заданное значение, отображение $r_2$ является гомоморфизмом $\Omega_2$-алгебры $B_k$ в $\Omega_2$-алгебру $B$. □

### 4.4. Приведенное декартово произведение представлений

ОПРЕДЕЛЕНИЕ 4.4.1. *Пусть $A_1$ - $\Omega_1$-алгебра. Пусть $\mathcal{A}_2$ - категория $\Omega_2$-алгебр. Мы определим* **категорию** $A_1(\mathcal{A}_2)$ **представлений** *$\Omega_1$-алгебры $A_1$ в $\Omega_2$-алгебре. Обьектами этой категории являются представления $\Omega_1$-алгебры $A_1$ в $\Omega_2$-алгебре. Морфизмами этой категории являются приведенные морфизмы соответствующих представлений.*                                                                          □

ТЕОРЕМА 4.4.2. *В категории $A_1(\mathcal{A}_2)$ существует произведение эффективных представлений $\Omega_1$-алгебры $A_1$ в $\Omega_2$-алгебре и это произведение является эффективным представлением $\Omega_1$-алгебры $A_1$.*

ДОКАЗАТЕЛЬСТВО. Пусть

$$A_2 = \prod_{i \in I} A_{2i}$$

произведение семейства $\Omega_2$-алгебр $\{A_{2i}, i \in I\}$ и для любого $i \in I$ отображение

$$t_i : A_2 \longrightarrow A_{2i}$$

является проекцией на множитель $i$. Для каждого $i \in I$, пусть

$$h_i : A_1 \relbar\joinrel\relbar\joinrel\twoheadrightarrow A_{2i}$$

эффективное $A_1$-представление в $\Omega_2$-алгебре $A_{2i}$.

Пусть $b_1 \in A_1$. Пусть $b_2 \in A_2$. Согласно утверждению 2.3.3.3, $A_2$-число $b_2$ может быть представлено в виде кортежа $A_{2i}$-чисел

(4.4.1)              $b_2 = (b_{2i}, i \in I) \quad b_{2i} = t_i(b_2) \in A_{2i}$

ЛЕММА 4.4.3. *Для каждого $i \in I$, рассмотрим диаграмму отображений*

(4.4.2)

*Пусть отображение*

$$g : A_1 \to \operatorname{End}(\Omega_2; A_2)$$

*определено равенством*

(4.4.3)              $g(b_1)(b_2) = (h_i(b_1)(b_{2i}), i \in I)$

*Тогда отображение $g$ является эффективным $A_1$-представлением в $\Omega_2$-алгебре $A_2$*

$$g : A_1 \relbar\joinrel\relbar\joinrel\twoheadrightarrow A_2$$

*Отображение $t_i$ является приведенным морфизмом представления $g$ в представление $h_i$.*

ДОКАЗАТЕЛЬСТВО.

4.4.3.1: Согласно определениям 3.1.1, отображение $h_i(b_1)$ является гомоморфизмом $\Omega_2$-алгебры $A_{2i}$. Согласно теореме 2.3.6, из коммутативности диаграммы (1) для каждого $i \in I$, следует, что отображение

$$g(b_1) : A_2 \to A_2$$

определённое равенством (4.4.3) является гомоморфизмом $\Omega_2$-алгебры $A_2$.

4.4.3.2: Согласно определению 3.1.1, множество $\mathrm{End}(\Omega_2; A_{2i})$ является $\Omega_1$-алгеброй. Согласно лемме 4.3.1, множество $^\circ A_2 \subseteq \mathrm{End}(\Omega_2; A_2)$ является $\Omega_1$-алгеброй.

4.4.3.3: Согласно определению 3.1.1, отображение

$$h_i : A_1 \to \mathrm{End}(\Omega_2, A_{2i})$$

является гомоморфизмом $\Omega_1$-алгебры. Согласно теореме 2.3.6, отображение

$$g : A_1 \to \mathrm{End}(\Omega_2; A_2)$$

определённое равенством

$$g(b_1) = (h_i(b_1), i \in I)$$

является гомоморфизмом $\Omega_1$-алгебры.

Согласно утверждениям 4.4.3.1, 4.4.3.3 и определению 3.1.1, отображение $g$ является $A_1$-представлением в $\Omega_2$-алгебре $A_2$.

Для любого $i \in I$, согласно определению 3.1.2, $A_1$-число $a_1$ порождает единственное преобразование

$$(4.4.4) \qquad\qquad b_{22i} = h_i(b_1)(b_{21i})$$

Пусть $b_{21}$, $b_{22} \in A_2$. Согласно утверждению 2.3.3.3, $A_2$-числа $b_{21}$, $b_{22}$ могут быть представлены в виде кортежей $A_{2i}$-чисел

$$(4.4.5) \qquad \begin{aligned} b_{21} &= (b_{21i}, i \in I) & b_{21i} &= t_i(b_{21}) \in A_{2i} \\ b_{22} &= (b_{22i}, i \in I) & b_{22i} &= t_i(b_{22}) \in A_{2i} \end{aligned}$$

Согласно определению (4.4.3) представления $g$, из равенств (4.4.4), (4.4.5) следует, что $A_1$-число $a_1$ порождает единственное преобразование

$$(4.4.6) \qquad\qquad b_{22} = (h_i(b_1)(b_{21i}), i \in I) = g(b_1)(b_{21})$$

Согласно определению 3.1.2, представление $g$ эффективно.

Из коммутативности диаграммы (1) и определения 3.2.2, следует, что отображение $t_i$ является приведенным морфизмом представления $g$ в представление $h_i$. $\odot$

Пусть

$$(4.4.7) \qquad\qquad d_2 = g(b_1)(b_2) \quad d_2 = (d_{2i}, i \in I)$$

Из равенств (4.4.3), (4.4.7) следует, что

$$(4.4.8) \qquad\qquad d_{2i} = h_i(b_1)(b_{2i})$$

Пусть $R_2$ - другой объект категории $\mathcal{A}_2$. Пусть

$$f : A_1 \relbar\joinrel\rightarrow R_2$$

эффективное $A_1$-представление в $\Omega_2$-алгебре $R_2$. Для любого $i \in I$, пусть существует морфизм

$$r_i : R_2 \longrightarrow A_{2i}$$

представлений из $f$ в $h_i$. Согласно замечанию 3.2.7, в диаграмме отображений

(4.4.9)

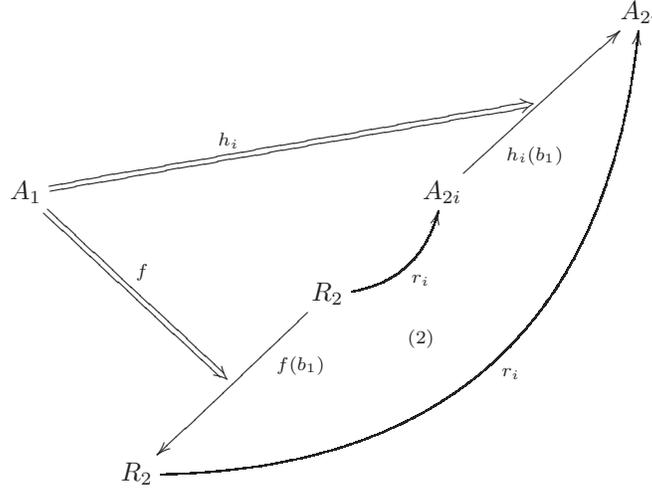

диаграмма (2) коммутативна. Согласно определению 2.3.1, существует единственный морфизм $\Omega_2$-алгебры

$$s : R_2 \longrightarrow A_2$$

такой, что коммутативна диаграмма

(4.4.10)

$$A_2 \xrightarrow{\ t_i\ } A_{2i} \qquad t_i(s) = r_i$$

Пусть $a_2 \in R_2$. Пусть

(4.4.11)          $$b_2 = s(a_2) \in A_2$$

Из коммутативности диаграммы (4.4.10) и утверждений (4.4.11), (4.4.1) следует, что

(4.4.12)          $$b_{2i} = r_i(a_2)$$

Пусть

(4.4.13)          $$c_2 = f(a_1)(a_2)$$

Из коммутативности диаграммы (2) и равенств (4.4.8), (4.4.12), (4.4.13) следует, что

(4.4.14)          $$d_{2i} = r_i(c_2)$$

Из равенств (4.4.8), (4.4.14) следует, что

(4.4.15)          $$d_2 = s(c_2)$$

что согласуется с коммутативносью диаграммы (4.4.10).

Для каждого $i \in I$, мы объединим диаграммы отображений (4.4.2), (4.4.10), (4.4.9)

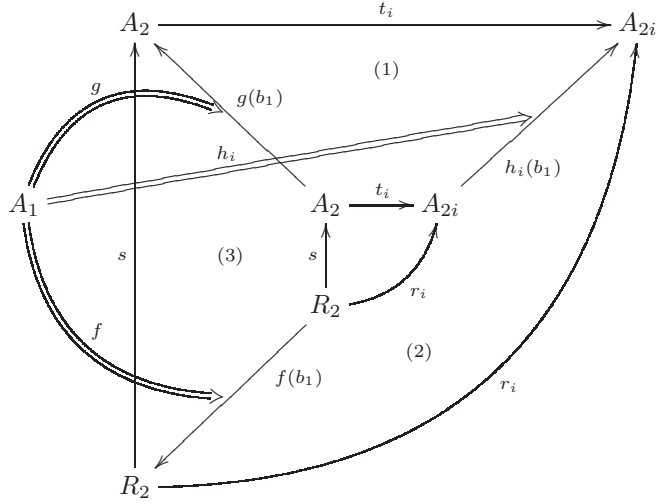

Из равенств (4.4.7), (4.4.11) и из равенств (4.4.13), (4.4.15), следует коммутативность диаграммы (3). Следовательно, отображение $s$ является приведенным морфизмом представлений из $f$ в $g$. Согласно определению 3.4.2, отображение $s$ является гомоморфизмом $\Omega_2$ алгебры. Согласно теореме 2.3.3 и определению 2.3.1, приведенный морфизм $s$ определён однозначно.

Согласно определению 2.3.1, представление $g$ и семейство морфизмов представления $(t_i, i \in I)$ является произведением в категории $A_1(\mathcal{A}_2)$.          $\square$

Определение 4.4.4. *Пусть $A$, $B_1$, ..., $B_n$, $B$ - универсальные алгебры. Пусть, для любого $k$, $k = 1, ..., n$,*

$$f_k : A \longrightarrow_* B_k$$

*эффективное представление $\Omega_1$-алгебры $A$ в $\Omega_2$-алгебре $B_k$. Пусть*

$$f : A \longrightarrow_* B$$

*эффективное представление $\Omega_1$-алгебры $A$ в $\Omega_2$-алгебре $B$. Отображение*

$$r_2 : B_1 \times ... \times B_n \to B$$

*называется* **приведенным полиморфизмом представлений** $f_1$, ..., $f_n$ *в представление $f$, если для любого $k$, $k = 1, ..., n$, при условии, что все переменные кроме переменной $x_k \in B_k$ имеют заданное значение, отображение $r_2$ является приведенным морфизмом представления $f_k$ в представление $f$.*

*Если $f_1 = ... = f_n$, то мы будем говорить, что отображение $r_2$ является приведенным полиморфизмом представления $f_1$ в представление $f$.*

*Если $f_1 = ... = f_n = f$, то мы будем говорить, что отображение $r_2$ является приведенным полиморфизмом представления $f$.*          $\square$

Теорема 4.4.5. *Пусть отображение $r_2$ - приведенный полиморфизм эффективных представлений $f_1$, ..., $f_n$ в эффективное представление $f$.*

- *Для любого $k$, $k = 1, ..., n$, отображение $r_2$ удовлетворяет равенству*

(4.4.16)          $r_2(m_1, ..., f_k(a)(m_k), ..., m_n) = f(a)(r_2(m_1, ..., m_n))$

● *Для любого* $k$, $l$, $k = 1$, ..., $n$, $l = 1$, ..., $n$, *отображение* $r_2$ *удовлетворяет равенству*

(4.4.17)
$$r_2(m_1, ..., f_k(a)(m_k), ..., m_l, ..., m_n)$$
$$= r_2(m_1, ..., m_k, ..., f_l(a)(m_l), ..., m_n)$$

● *Пусть* $\omega_2 \in \Omega_2(p)$. *Для любого* $k$, $k = 1$, ..., $n$, *отображение* $r_2$ *удовлетворяет равенству*

(4.4.18)
$$r_2(m_1, ..., m_{k \cdot 1} ... m_{k \cdot p} \omega_2, ..., m_n)$$
$$= r_2(m_1, ..., m_{k \cdot 1}, ..., m_n) ... r_2(m_1, ..., m_{k \cdot p}, ..., m_n)\omega_2$$

Доказательство. Равенство (4.4.16) следует из определения 4.4.4 и равенства (3.4.4). Равенство (4.4.17) следует из равенства (4.4.16). Равенство (4.4.18) следует из утверждения, что, для любого $k$, $k = 1$, ..., $n$, при условии, что все переменные кроме переменной $m_k \in B_k$ имеют заданное значение, отображение $r_2$ является гомоморфизмом $\Omega_2$-алгебры $B_k$ в $\Omega_2$-алгебру $B$.    □

Мы также будем говорить, что отображение $r_2$ является приведенным полиморфизмом представлений в $\Omega_2$-алгебрах $B_1$, ..., $B_n$ в представление в $\Omega_2$-алгебре $B$.

## 4.5. Мультипликативная Ω-группа

Пусть отображение
$$f : A \longrightarrow\!\!\!* \ B$$

является представлением $\Omega_1$-алгебры $A$ в $\Omega_2$ алгебре $B$. Согласно теореме 3.5.3, множество $\operatorname{End}(A(\Omega_2); B)$ является полугруппой. В тоже время [4.1]

(4.5.1)
$$\operatorname{End}(A(\Omega_2); B) \subseteq \operatorname{End}(\Omega_2; B)$$

Согласно определению 3.1.1, множество $\operatorname{End}(\Omega_2, B)$ является $\Omega_2$-алгеброй. Однако из утверждения (4.5.1) не следует, что множество $\operatorname{End}(A(\Omega_2); B)$ является $\Omega_2$-алгеброй.

Чтобы понять, при каких условиях множество $\operatorname{End}(A(\Omega_2); B)$ является $\Omega_2$-алгеброй, мы рассмотрим связь между множеством представлений $\Omega_1$-алгебры $A$ в $\Omega_2$-алгебре $B$ и множеством приведенных морфизмов этих представлений.

Теорема 4.5.1. *Пусть отображение*
$$r : B \to B$$

*является приведенным эндоморфизмом представления*
$$f : A \longrightarrow\!\!\!* \ B$$

$\Omega_1$-*алгебры* $A$ *в* $\Omega_2$ *алгебре* $B$. *Отображение*

(4.5.2)
$$rf : a \in A \to r \circ f(a) \in \operatorname{End}(\Omega_2; B)$$

*является представлением* $\Omega_1$-*алгебры* $A$ *в* $\Omega_2$ *алгебре* $B$ *тогда и только тогда, когда на множестве* $f(A) \subseteq \operatorname{End}(\Omega_2, B)$ *произведение* $\circ$ *отображений дистрибутивно слева относительно произвольной операции* $\omega \in \Omega_1$

(4.5.3)
$$r \circ (f(a_1) ... f(a_p)\omega) = (r \circ f(a_1)) ... (r \circ f(a_p))\omega$$

---

[4.1] В утверждении (4.5.1), я обозначил $\Omega_2$ категорию $\Omega_2$-алгебр и $A(\Omega_2)$ категорию представлений $\Omega_1$-алгебры $A$ в $\Omega_2$-алгебре.

Доказательство. Согласно определению 3.1.1, отображение $f(a)$ является эндоморфизмом $\Omega_2$-алгебры $B$. Согласно определениям 3.2.2, 3.4.2, отображение $r$ является эндоморфизмом $\Omega_2$-алгебры $B$. Следовательно, отображение $r \circ f(a)$ является эндоморфизмом $\Omega_2$-алгебры $B$.

4.5.1.1: Согласно определению 3.1.1, отображение $rf$ является представлением $\Omega_1$-алгебры $A$ в $\Omega_2$ алгебре $B$ тогда и только тогда, когда отображение $rf$ является гомоморфизмом $\Omega_1$-алгебры.

4.5.1.2: Утверждение 4.5.1.1 означает, что для любой операции $\omega \in \Omega_1$ верно равенство

$$(4.5.4) \qquad \begin{aligned} r \circ f(a_1...a_p\omega) &= (rf)(a_1...a_p\omega) = ((rf)(a_1))...((rf)(a_p))\omega \\ &= (r \circ f(a_1))...(r \circ f(a_p))\omega \end{aligned}$$

Поскольку отображение $f$ является представлением $\Omega_1$-алгебры $A$ в $\Omega_2$ алгебре $B$, то, согласно определению 3.1.1, отображение $f$ является гомоморфизмом $\Omega_1$-алгебры

$$(4.5.5) \qquad r \circ f(a_1...a_p\omega) = r \circ (f(a_1)...f(a_p)\omega)$$

Равенство (4.5.3) является следствием равенств (4.5.4), (4.5.5).

Теорема является следствием утверждения 4.5.1.2. $\qquad\square$

Теорема 4.5.2. *Пусть отображение*

$$f : A \longrightarrow B$$

*является представлением $\Omega_1$-алгебры $A$ в $\Omega_2$ алгебре $B$. Пусть*

$$(4.5.6) \qquad f(A) = \operatorname{End}(A(\Omega_2); B)$$

4.5.2.1: *Произведение $\circ$ в полугруппе $\operatorname{End}(A(\Omega_2); B)$ коммутативно.*

4.5.2.2: *Произведение $\circ$ в полугруппе $\operatorname{End}(A(\Omega_2); B)$ порождает произведение $*$ в $\Omega_1$-алгебре $A$ таким образом, что*

$$(4.5.7) \qquad f(a * b) = f(a) \circ f(b)$$

4.5.2.3: *Полугруппа $\operatorname{End}(A(\Omega_2); B)$ является $\Omega_1$-алгеброй.*

Доказательство. Пусть отображение $h$ является эндоморфизмом представления $f$. Согласно утверждению (4.5.6), существует $b \in A$ такое, что $h = f(b)$. Следовательно, равенство

$$(4.5.8) \qquad f(a) \circ f(b) = f(b) \circ f(a)$$

является следствием равенства (3.4.4). Согласно утверждению (4.5.6), отображения $f(a)$, $f(b)$ являются эндоморфизмами представления $f$. Следовательно, произведение $\circ$ в полугруппе $\operatorname{End}(A(\Omega_2); B)$ коммутативно.

Согласно теореме 3.5.3, произведение эндоморфизмов $f(a)$, $f(b)$ представления $f$ является эндоморфизмом $h$ представления $f$. Согласно утверждению (4.5.6), существует $c \in A$ такое, что $h = f(c)$. Бинарная операция $*$ на множестве $A$ определена равенством

$$c = a * b$$

Следовательно, утверждение 4.5.2.2 верно.

Пусть отображения $h_1, ..., h_p$ являются эндоморфизмами представления $f$. Согласно утверждению (4.5.6), существуют $A$-числа $a_1, ..., a_p$ такие, что

$$h_1 = f(a_1) \quad ... \quad h_n = f(a_n)$$

Поскольку отображение $f$ является представлением $\Omega_1$-алгебры $A$ в $\Omega_2$ алгебре $B$, то, согласно определению 3.1.1, отображение $f$ является гомоморфизмом $\Omega_1$-алгебры $A$

$$(4.5.9) \qquad h_1...h_p\omega = f(a_1)...f(a_p)\omega = f(a_1...a_p\omega)$$

Согласно утверждению (4.5.6), $h_1...h_p\omega \in \text{End}(A(\Omega_2); B)$. Следовательно, утверждение 4.5.2.3 верно. $\qquad\square$

Согласно теореме 4.5.2, если утверждение (4.5.6) выполнено, то на множестве $\text{End}(A(\Omega_2); B)$ определены две алгебраические структуры. А именно, множество $\text{End}(A(\Omega_2); B)$ является полугруппой и в тоже время это множество является $\Omega_1$-алгеброй. Аналогичное утверждение верно для $\Omega_1$-алгебры $A$. Однако мы не можем утверждать, что операция произведения в $\Omega_1$-алгебре $A$ дистрибутивна по отношению к произвольной операции $\omega \in \Omega_1$ (смотри теорему 4.5.1).

Теорема 4.5.3. *Пусть отображение*

$$f : A \longrightarrow B$$

*является представлением $\Omega_1$-алгебры $A$ в $\Omega_2$ алгебре $B$ и*

$$f(A) = \text{End}(A(\Omega_2); B)$$

*Произведение $*$, определённое в $\Omega_1$-алгебре $A$, дистрибутивно относительно произвольной операции $\omega \in \Omega_1$ тогда и только тогда, когда отображение*

$$(4.5.10) \qquad f(b*a) : a \in A \to f(b*a) \in \text{End}(\Omega_2; B)$$

*является представлением $\Omega_1$-алгебры $A$ в $\Omega_2$ алгебре $B$*

Доказательство. Согласно утверждению 4.5.2.2, для нас не имеет значение рассматриваем ли мы $\Omega_1$-алгебру $A$ или мы рассматриваем $\Omega_1$-алгебру $\text{End}(A(\Omega_2); B)$. Теорема является следствием определения (4.5.7) произведения $*$ в $\Omega_1$-алгебре $A$, а также теоремы 4.5.1 и утверждений 4.5.2.1, 4.5.2.3. $\qquad\square$

В теореме 4.5.3, мы встречаем универсальную алгебру, похожую на $\Omega$-группу, однако эта алгебра отличается от $\Omega$-группы. Поскольку эта универсальная алгебра играет важную роль в теории представлений, мы рассмотрим определения 4.5.4, 4.5.5.

Определение 4.5.4. *Пусть произведение*

$$c_1 = a_1 * b_1$$

*является операцией $\Omega_1$-алгебры $A$. Положим $\Omega = \Omega_1 \setminus \{*\}$. Если $\Omega_1$-алгебра $A$ является группой относительно произведения и для любой операции $\omega \in \Omega(n)$ умножение дистрибутивно относительно операция $\omega$*

$$a*(b_1...b_n\omega) = (a*b_1)...(a*b_n)\omega$$

$$(b_1...b_n\omega)*a = (b_1*a)...(b_n*a)\omega$$

*то $\Omega_1$-алгебра $A$ называется* **мультипликативной $\Omega$-группой**. $\qquad\square$

Определение 4.5.5. *Если*

$$(4.5.11) \qquad a*b = b*a$$

*то мультипликативная $\Omega$-группа называется* **абелевой**. $\qquad\square$

Определение 4.5.6. *Если*

$$(4.5.12) \qquad a * (b * c) = (a * b) * c$$

*то мультипликативная Ω-группа называется* **асоциативной**.    □

Теорема 4.5.7. *Пусть* $A$, $B_1$, ..., $B_n$, $B$ *- универсальные алгебры. Пусть, для любого* $k$, $k = 1$, ..., $n$,

$$f_k : A \longrightarrow\!\!\!* \longrightarrow B_k$$

*представление* $\Omega_1$-*алгебры* $A$ *в* $\Omega_2$-*алгебре* $B_k$. *Пусть*

$$f : A \longrightarrow\!\!\!* \longrightarrow B$$

*представление* $\Omega_1$-*алгебры* $A$ *в* $\Omega_2$-*алгебре* $B$. *Пусть отображение*

$$r_2 : B_1 \times ... \times B_n \to B$$

*является приведенным полиморфизмом представлений* $f_1$, ..., $f_n$ *в представление* $f$. *Произведение* $\circ$, *определённое в* $\Omega_1$-*алгебре* $f(A)$, *коммутативно.*

*Представление*

$$f : A \longrightarrow\!\!\!* \longrightarrow B$$

*допускает приведенный полиморфизм представлений тогда и только тогда, когда следующие условия выполнены*

4.5.7.1: *Произведение* $\circ$, *определённое в* $\Omega_1$-*алгебре* $\text{End}(A(\Omega_2); B)$, *дистрибутивно относительно произвольной операции* $\omega \in \Omega_1$

4.5.7.2: $f(a * b) = f(a) \circ f(b)$

Доказательство. Пользуясь равенством (4.4.16), мы можем записать выражение

$$(4.5.13) \qquad r_2(m_1, ..., f_k(a_k)(m_k), ..., f_l(a_l)(m_l), ..., m_n)$$

либо в виде

$$(4.5.14) \quad \begin{aligned} &r_2(m_1, ..., f_k(a_k)(m_k), ..., f_l(a_l)(m_l), ..., m_n) \\ &= f(a_k)(r_2(m_1, ..., m_k, ..., f_l(a_l)(m_l), ..., m_n)) \\ &= f(a_k)(f(a_l)(r_2(m_1, ..., m_k, ..., m_l, ..., m_n))) \\ &= (f(a_k) \circ f(a_l))(r_2(m_1, ..., m_k, ..., m_l, ..., m_n)) \end{aligned}$$

либо в виде

$$(4.5.15) \quad \begin{aligned} &r_2(m_1, ..., f_k(a_k)(m_k), ..., f_l(a_l)(m_l), ..., m_n) \\ &= f(a_l)(r_2(m_1, ..., f_k(a_k)(m_k), ..., m_l, ..., m_n)) \\ &= f(a_l)(f(a_k)(r_2(m_1, ..., m_k, ..., m_l, ..., m_n))) \\ &= (f(a_l) \circ f(a_k))(r_2(m_1, ..., m_k, ..., m_l, ..., m_n)) \end{aligned}$$

Коммутативность произведения $\circ$ следует из равенств (4.5.14), (4.5.15).    □

Теорема 4.5.8. *Пусть*

$$f : A \longrightarrow\!\!\!* \longrightarrow B$$

*представление* $\Omega_1$-*алгебры* $A$ *в* $\Omega_2$-*алгебре* $B$ *и*

$$(4.5.16) \qquad f(A) = \text{End}(A(\Omega_2); B)$$

*Тогда представление* $f$ *допускает приведенный полиморфизм представлений.*

*Пусть* $\Omega = \Omega_1 \setminus \{*\}$. *Представление*

$$h : A_1 \to \mathrm{End}(\Omega; A_1) \quad h(a) : b \in A_1 \to a * b \in A_1$$

*полугруппы $A_1$ в Ω-алгебре $A_1$ существует тогда и только тогда, когда для любой операции $\omega \in \Omega(n)$ умножение* **дистрибутивно** *относительно операции $\omega$*

(4.5.17)                    $a * (b_1 ... b_n \omega) = (a * b_1)...(a * b_n)\omega$

(4.5.18)                    $(b_1 ... b_n \omega) * a = (b_1 * a)...(b_n * a)\omega$

Доказательство. Согласно определению 3.1.1, равенства (4.5.17), (4.5.18) верны тогда и только тогда, когда отображение $h$ является представлением полугруппы $A_1$ в Ω-алгебре $A_1$. Одновременно равенства (4.5.17), (4.5.18) выражают закон дистрибутивности умножения относительно операции $\omega$.          □

В $\Omega_1$-алгебре $A_1$, мы определили произведение, согласованное с однотранзитивным представлением в $\Omega_2$-алгебре $A_2$. Эту конструкцию можно построить в случае произвольного представления при условии, что произведение в $\Omega_1$-алгебре $A_1$ определено однозначно. Однако в общем случае произведение может быть некоммутативным.

Теорема 4.5.9. *Пусть*

$$A \xrightarrow{\;\;\ast\;\;} B_1 \qquad A \xrightarrow{\;\;\ast\;\;} B_2 \qquad A \xrightarrow{\;\;\ast\;\;} B$$

*эффективные представления абелевой мультиплиткативной $\Omega_1$-группы $A$ в $\Omega_2$-алгебрах $B_1$, $B_2$, $B$. Допустим $\Omega_2$-алгебра имеет 2 операции, а именно $\omega_1 \in \Omega(m)$, $\omega_2 \in \Omega(n)$. Необходимым условием существования приведенного полиморфизма*

$$R : B_1 \times B_2 \to B$$

*является равенство*

(4.5.19)      $(a_{11}...a_{1n}\omega_2)...(a_{m1}...a_{mn}\omega_2)\omega_1 = (a_{11}...a_{m1}\omega_1)...(a_{1n}...a_{mn}\omega_1)\omega_2$

Доказательство. Пусть $a_1, ..., a_p \in B_1$, $b_1, ..., b_q \in B_2$. Согласно равенству (4.4.18), выражение

(4.5.20)                    $r_2(a_1...a_p\omega_1, b_1...b_q\omega_2)$

может иметь 2 значения

$$r_2(a_1...a_m\omega_1, b_1...b_n\omega_2)$$

(4.5.21)   $= r_2(a_1, b_1...b_n\omega_2)...r_2(a_m, b_1...b_n\omega_2)\omega_1$

$$= (r_2(a_1, b_1)...r_2(a_1, b_n)\omega_2)...(r_2(a_m, b_1)...r_2(a_m, b_n)\omega_2)\omega_1$$

$$r_2(a_1...a_m\omega_1, b_1...b_n\omega_2)$$

(4.5.22)   $= r_2(a_1...a_m\omega_1, b_1)...r_2(a_1...a_m\omega_1, b_n)\omega_2$

$$= (r_2(a_1, b_1)...r_2(a_m, b_1)\omega_1)...(r_2(a_1, b_n)...r_2(a_m, b_n)\omega_1)\omega_2$$

Из равенств (4.5.21), (4.5.22) следует, что

(4.5.23)   $(r_2(a_1, b_1)...r_2(a_1, b_n)\omega_2)...(r_2(a_m, b_1)...r_2(a_m, b_n)\omega_2)\omega_1$

$$= (r_2(a_1, b_1)...r_2(a_m, b_1)\omega_1)...(r_2(a_1, b_n)...r_2(a_m, b_n)\omega_1)\omega_2$$

Следовательно, выражение (4.5.20) определенно корректно тогда и только тогда, когда равенство (4.5.23) верно. Положим

(4.5.24)                     $a_{i \cdot j} = r_2(a_i, b_j) \in A$

Равенство (4.5.19) является следствием равенств (4.5.23), (4.5.24).            □

ТЕОРЕМА 4.5.10. *Существует приведенный полиморфизм эффективного представления абелевой мультиплиткативной Ω-группы в абелевой группе.*

ДОКАЗАТЕЛЬСТВО. Поскольку операция сложения в абелевой группе коммутативна и ассоциативна, то теорема является следствием теоремы 4.5.9.            □

ТЕОРЕМА 4.5.11. *Не существует приведенный полиморфизм эффективного представления абелевой мультиплиткативной Ω-группы в кольце.*

ДОКАЗАТЕЛЬСТВО. В кольце определены две операции: сложение, которое коммутативно и ассоциативно, и произведение, которое дистрибутивно относительно сложения. Согласно теореме 4.5.9, если существует полиморфизм эффективного представления в кольцо, то сложение и произведение должны удовлетворять равенству

(4.5.25)               $a_{11}a_{21} + a_{12}a_{22} = (a_{11} + a_{12})(a_{21} + a_{22})$

Однако правая часть равенства (4.5.25) имеет вид

$$(a_{11} + a_{12})(a_{21} + a_{22}) = (a_{11} + a_{12})a_{21} + (a_{11} + a_{12})a_{22}$$
$$= a_{11}a_{21} + a_{12}a_{21} + a_{11}a_{22} + a_{12}a_{22}$$

Следовательно, равенство (4.5.25) не верно.            □

ВОПРОС 4.5.12. *Возможно, что полиморфизм представлений существует только для эффективного представления в Абелевая группе. Однако это утверждение пока не доказано.*            □

## 4.6. Ω-кольцо

ОПРЕДЕЛЕНИЕ 4.6.1. *Пусть сложение*

$$c_1 = a_1 + b_1$$

*которое не обязательно коммутативно, и произведение*

$$c_1 = a_1 * b_1$$

*являются операциями $\Omega_1$-алгебры $A$. Положим  $\Omega = \Omega_1 \setminus \{+, *\}$. Если $\Omega_1$-алгебра $A$ является $\Omega \cup \{*\}$-группой и мультипликативной $\Omega \cup \{+\}$-группой, то $\Omega_1$-алгебра $A$ называется* Ω-**кольцом**.            □

ТЕОРЕМА 4.6.2. *Произведение в Ω-кольце* **дистрибутивно** *относительно сложения*

$$a * (b_1 + b_2) = a * b_1 + a * b_2$$
$$(b_1 + b_2) * a = b_1 * a + b_2 * a$$

ДОКАЗАТЕЛЬСТВО. Теорема является следствием определений 4.2.10, 4.5.4, 4.6.1.            □

Определение 4.6.3. *Пусть $A$ - Ω-кольцо.* **Матрица** *над Ω-кольцом $A$ - это таблица $A$-чисел $a_j^i$, где индекс $i$ - это номер строки и индекс $j$ - это номер столбца.* □

Соглашение 4.6.4. *Мы будем пользоваться соглашением Эйнштейна о сумме. Это означает, что, когда индекс присутствует в выражении дважды (один вверху и один внизу) и множество индексов известно, это выражение подразумевает сумму по повторяющемуся индексу. В этом случае предполагается известным множество индекса суммирования и знак суммы опускается*

$$a^i v_i = \sum_{i \in I} a^i v_i$$

*Я буду явно указывать множество индексов, если это необходимо.* □

Произведение матриц связано с произведением гомоморфизмов векторных пространств над полем. Согласно традиции произведение матриц $a$ и $b$ определено как произведение строк матрицы $a$ и столбцов матрицы $b$.

Пример 4.6.5. *Пусть $\overline{\overline{e}}$ - базис правого векторного пространства $V$ над $D$-алгеброй $A$ (смотри определение 9.6.2 и теорему 9.6.15). Мы представим базис $\overline{\overline{e}}$ как строку матрицы*

$$e = \begin{pmatrix} e_1 & ... & e_n \end{pmatrix}$$

*Мы можем представить координаты вектора $v$ как вектор столбец*

$$v = \begin{pmatrix} v^1 \\ ... \\ v^n \end{pmatrix}$$

*Поэтому мы можем представить вектор $v$ как традиционное произведение матриц*

$$v = \begin{pmatrix} e_1 & ... & e_n \end{pmatrix} \begin{pmatrix} v^1 \\ ... \\ v^n \end{pmatrix} = e_i v^i$$

*Линейный гомоморфизм правого векторного пространства $V$ может быть представлен с помощью матрицы*

$$(4.6.1) \qquad v'^i = f_j^i v^j$$

*Равенство (4.6.1) выражает традиционное произведение матриц $f$ и $v$.* □

Пример 4.6.6. *Пусть $\overline{\overline{e}}$ - базис левого векторного пространства $V$ над $D$-алгеброй $A$ (смотри определение 9.5.2 и теорему 9.5.15). Мы представим базис $\overline{\overline{e}}$ как строку матрицы*

$$e = \begin{pmatrix} e_1 & ... & e_n \end{pmatrix}$$

*Мы можем представить координаты вектора $v$ как вектор столбец*

$$v = \begin{pmatrix} v^1 \\ ... \\ v^n \end{pmatrix}$$

*Однако мы не можем представить вектор*

$$v = v^i e_i$$

*как традиционное произведение матриц*

$$v = \begin{pmatrix} v^1 \\ ... \\ v^n \end{pmatrix} \qquad e = \begin{pmatrix} e_1 & ... & e_n \end{pmatrix}$$

*так как это произведение не определено. Линейный гомоморфизм левого векторного пространства $V$ может быть представлен с помощью матрицы*

(4.6.2)
$$v'^i = v^j f^i_j$$

*Равенство (4.6.2) не может быть выражено как традиционное произведение матриц $v$ и $f$.*                                                      □

Из примеров 4.6.5, 4.6.6 следует, что мы не можем ограничиться традиционным произведением матриц и нам нужно определить два вида произведения матриц. Чтобы различать эти произведения, мы вводим новые обозначения. Для совместимости обозначений с существующими мы будем иметь в виду $_*{}^*$-произведение, когда нет явных обозначений.

ОПРЕДЕЛЕНИЕ 4.6.7. *Пусть число столбцов матрицы $a$ равно числу строк матрицы $b$. $_*{}^*$-произведение матриц $a$ и $b$ имеет вид*

(4.6.3)
$$\begin{cases} a_*{}^*b = \left( a^i_k b^k_j \right) \\ (a_*{}^*b)^i_j = a^i_k b^k_j \end{cases}$$

*и может быть выражено как произведение строк матрицы $a$ и столбцов матрицы $b$.*[4.2]                                                      □

ОПРЕДЕЛЕНИЕ 4.6.8. *Пусть число строк матрицы $a$ равно числу столбцов матрицы $b$. $^*{}_*$-произведение матриц $a$ и $b$ имеет вид*

(4.6.4)
$$\begin{cases} a^*{}_*b = \left( a^k_i b^j_k \right) \\ (a^*{}_*b)^i_j = a^k_i b^j_k \end{cases}$$

---

[4.2] Мы будем пользоваться символом $_*{}^*$- в последующей терминологии и обозначениях. Мы будем читать символ $_*{}^*$ как $rc$-произведение или произведение строки на столбец. Символ произведения строки на столбец сформирован из двух символов операции произведения, которые записываются на месте индекса суммирования. Например, если произведение $A$-чисел имеет вид $a \circ b$, то $_*{}^*$-произведение матриц $a$ и $b$ имеет вид $a_\circ{}^\circ b$.

*и может быть выражено как произведение столбцов матрицы a и строк матрицы b.* [4.3]                                                                    □

Мы так же определим следующие операции на множестве матриц.

ОПРЕДЕЛЕНИЕ 4.6.9. *Транспонирование $a^T$ матрицы a меняет местами строки и столбцы*

$$(4.6.5) \qquad\qquad (a^T)^i_j = a^j_i$$

□

ОПРЕДЕЛЕНИЕ 4.6.10. *Сумма матриц a и b определена равенством*

$$(a+b)^i_j = a^i_j + b^i_j$$

□

ЗАМЕЧАНИЕ 4.6.11. *Мы будем пользоваться символом $_*^*$- или $_*_*$- в имени свойств каждого произведения и в обозначениях. Мы можем читать символы $_*^*$ и $_*_*$ как rc-произведение и cr-произведение. Это правило мы распространим на последующую терминологию.*                                        □

ТЕОРЕМА 4.6.12.

$$(4.6.6) \qquad\qquad (a_*{}^*b)^T = a^T{}_*{}^*b^T$$

ДОКАЗАТЕЛЬСТВО. Цепочка равенств

$$(4.6.7) \qquad ((a_*{}^*b)^T)^j_i = (a_*{}^*b)^j_i = a^i_k b^k_j = (a^T)^k_i (b^T)^j_k = ((a^T)^*{}_*(b^T))^j_i$$

следует из (4.6.5), (4.6.3) и (4.6.4). Равенство (4.6.6) следует из (4.6.7).       □

ОПРЕДЕЛЕНИЕ 4.6.13. **Бикольцо** $\mathcal{A}$ *- это множество, на котором мы определили унарную операцию, называемую транспозицией, и три бинарных операции, называемые $_*^*$-произведение, $_*_*$-произведение и сумма, такие что*

- *$_*^*$-произведение и сумма определяют структуру кольца на $\mathcal{A}$*
- *$_*_*$-произведение и сумма определяют структуру кольца на $\mathcal{A}$*
- *оба произведения имеют общую единицу $\delta$*
- *произведения удовлетворяют равенству*

$$(a_*{}^*b)^T = a^T{}_*{}^*b^T$$

- *транспозиция единицы есть единица*

$$(4.6.8) \qquad\qquad \delta^T = \delta$$

- *двойная транспозиция есть исходный элемент*

$$(4.6.9) \qquad\qquad (a^T)^T = a$$

□

ТЕОРЕМА 4.6.14 (принцип двойственности для бикольца). *Пусть $\mathcal{A}$ - истинное утверждение о бикольце A. Если мы заменим одновременно*

- *$a \in A$ и $a^T$*

---

[4.3] Мы будем пользоваться символом $_*_*$- в последующей терминологии и обозначениях. Мы будем читать символ $_*_*$ как cr-произведение или произведение столбца на строку. Символ произведение столбца на строку сформирован из двух символов операции произведения, которые записываются на месте индекса суммирования. Например, если произведение A-чисел имеет вид $a \circ b$, то $_*_*$-произведение матриц a и b имеет вид $a^\circ_\circ b$.

- $*^*$-*произведение и* $^*_*$-*произведение*

*то мы снова получим истинное утверждение.*

**Теорема 4.6.15** (принцип двойственности для бикольца матриц). *Пусть $A$ является бикольцом матриц. Пусть $\mathcal{A}$ - истинное утверждение о матрицах. Если мы заменим одновременно*

- *строки и столбцы всех матриц*
- $*^*$-*произведение и* $^*_*$-*произведение*

*то мы снова получим истинное утверждение.*

Доказательство. Непосредственное следствие теоремы 4.6.14. $\square$

Замечание 4.6.16. *Если произведение в $\Omega$-кольце коммутативно, то*

(4.6.10) $$a_*{}^*b = (a_i^k b_k^j) = (b_k^j a_i^k) = b^*{}_*a$$

**Приводимое бикольцо** - *это бикольцо, в котором выполняется* **условие приводимости произведений** (4.6.10). *Поэтому в приводимом бикольце достаточно рассматривать только* $*^*$-*произведение. Однако в тех случаях, когда порядок сомножителей существенен, мы будем пользоваться также* $^*_*$-*произведением.* $\square$

## 4.7. Тензорное произведение представлений

Определение 4.7.1. *Пусть $A$ является абелевой мультипликативной $\Omega_1$-группой. Пусть $A_1$, ..., $A_n$ - $\Omega_2$-алгебры.[4.4] Пусть для любого $k$, $k = 1$, ..., $n$,*

$$f_k : A \longrightarrow{*} A_k$$

*эффективное представление мультипликативной $\Omega_1$-группы $A$ в $\Omega_2$-алгебре $A_k$. Рассмотрим категорию $\mathcal{A}$ объектами которой являются приведенные полиморфизмы представлений $f_1$, ..., $f_n$*

$$r_1 : B_1 \times ... \times B_n \longrightarrow S_1 \qquad r_2 : B_1 \times ... \times B_n \longrightarrow S_2$$

*где $S_1$, $S_2$ - $\Omega_2$-алгебры и*

$$g_1 : A \longrightarrow{*} S_1 \qquad g_2 : A \longrightarrow{*} S_2$$

*эффективные представления мультипликативной $\Omega_1$-группы $A$. Мы определим морфизм $r_1 \to r_2$ как приведенный морфизм представлений $h : S_1 \to S_2$, для которого коммутативна диаграмма*

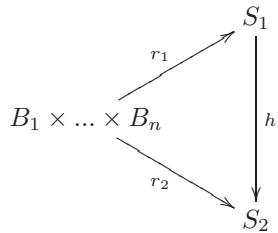

*Универсальный объект $B_1 \otimes ... \otimes B_n$ категории $\mathcal{A}$ называется* **тензорным произведением** *представлений $A_1$, ..., $A_n$.* $\square$

---

[4.4] Я определяю тензорное произведение представлений универсальной алгебры по аналогии с определением в [2], с. 456 - 458.

Теорема 4.7.2. *Если тензорное произведение эффективных представлений существует, то тензорное произведение определено однозначно с точностью до изоморфизма представлений.*

Доказательство. Пусть $A$ является абелевой мультипликативной $\Omega_1$-группой. Пусть $A_1$, ..., $A_n$ - $\Omega_2$-алгебры. Пусть для любого $k$, $k = 1$, ..., $n$,

$$f_k : A \longrightarrow_* B_k$$

эффективное представление мультипликативной $\Omega_1$-группы $A$ в $\Omega_2$-алгебре $B_k$. Пусть эффективные представления

$$g_1 : A \longrightarrow_* S_1 \qquad g_2 : A \longrightarrow_* S_2$$

являются тензорным произведением представлений $B_1$, ..., $B_n$. Из коммутативности диаграммы

(4.7.1)

следует, что

(4.7.2)
$$R_1 = h_2 \circ h_1 \circ R_1$$
$$R_2 = h_1 \circ h_2 \circ R_2$$

Из равенств (4.7.2) следует, что морфизмы представления $h_1 \circ h_2$, $h_2 \circ h_1$ являются тождественными отображениями. Следовательно, морфизмы представления $h_1$, $h_2$ являются изоморфизмами. □

Соглашение 4.7.3. *Алгебры $S_1$, $S_2$ могут быть различными множествами. Однако они неразличимы для нас, если мы рассматриваем их как изоморфные представления. В этом случае мы будем писать $S_1 = S_2$.* □

Определение 4.7.4. *Тензорное произведение*

$$B^{\otimes n} = B_1 \otimes ... \otimes B_n \quad B_1 = ... = B_n = B$$

*называется* **тензорной степенью** *представления $B$.* □

Теорема 4.7.5. *Если существует полиморфизм представлений, то тензорное произведение представлений существует.*

Доказательство. Пусть

$$f : A \longrightarrow_* M$$

представление $\Omega_1$-алгебры $A$, порождённое декартовым произведением $B_1 \times \ldots \times B_n$ множеств $B_1$, ..., $B_n$.[4.5] Инъекция

$$i : B_1 \times \ldots \times B_n \longrightarrow M$$

определена по правилу[4.6]

(4.7.3)                              $i \circ (b_1, ..., b_n) = (b_1, ..., b_n)$

Пусть $N$ - отношение эквивалентности, порождённое равенствами[4.7]

(4.7.4)        $(b_1, ..., b_{i \cdot 1}...b_{i \cdot p}\omega, ..., b_n) = (b_1, ..., b_{i \cdot 1}, ..., b_n)...(b_1, ..., b_{i \cdot p}, ..., b_n)\omega$

(4.7.5)        $(b_1, ..., f_i(a) \circ b_i, ..., b_n) = f(a) \circ (b_1, ..., b_i, ..., b_n)$

$$b_k \in B_k \quad k = 1, ..., n \quad b_{i \cdot 1}, ..., b_{i \cdot p} \in B_i \quad \omega \in \Omega_2(p) \quad a \in A$$

Лемма 4.7.6. *Пусть* $\omega \in \Omega_2(p)$. *Тогда*

(4.7.6)
$$f(c) \circ (b_1, ..., b_{i \cdot 1}...b_{i \cdot p}\omega, ..., b_n)$$
$$= f(c) \circ ((b_1, ..., b_{i \cdot 1}, ..., b_n)...(b_1, ..., b_{i \cdot p}, ..., b_n)\omega)$$

Доказательство. Из равенства (4.7.5) следует

(4.7.7)        $f(c) \circ (b_1, ..., b_{i \cdot 1}...b_{i \cdot p}\omega, ..., b_n) = (b_1, ..., f_i(c \circ (b_{i \cdot 1}...b_{i \cdot p}\omega)), ..., b_n)$

Так как $f_i(c)$ - эндоморфизм $\Omega_2$-алгебры $B_i$, то из равенства (4.7.7) следует

(4.7.8) $f(c) \circ (b_1, ..., b_{i \cdot 1}...b_{i \cdot p}\omega, ..., b_n) = (b_1, ..., (f_i(c) \circ b_{i \cdot 1})...(f_i(c) \circ b_{i \cdot p})\omega, ..., b_n)$

Из равенств (4.7.8), (4.7.4) следует

(4.7.9)
$$f(c) \circ (b_1, ..., b_{i \cdot 1}...b_{i \cdot p}\omega, ..., b_n)$$
$$= (b_1, ..., f_i(c) \circ b_{i \cdot 1}, ..., b_n)...(b_1, ..., f_i(c) \circ b_{i \cdot p}, ..., b_n)\omega$$

Из равенств (4.7.9), (4.7.5) следует

(4.7.10)
$$f(c) \circ (b_1, ..., b_{i \cdot 1}...b_{i \cdot p}\omega, ..., b_n)$$
$$= (f(c) \circ (b_1, ..., b_{i \cdot 1}, ..., b_n))...(f(c) \circ (b_1, ..., b_{i \cdot p}, ..., b_n))\omega$$

Так как $f(c)$ - эндоморфизм $\Omega_2$-алгебры $B$, то равенство (4.7.6) следует из равенства (4.7.10).                                                                ⊙

Лемма 4.7.7.

(4.7.11)        $f(c) \circ (b_1, ..., f_i(a) \circ b_i, ..., b_n) = f(c) \circ (f(a) \circ (b_1, ..., b_i, ..., b_n))$

---

[4.5]Согласно теоремам 2.3.3, 4.4.2, множество, порождённое приведенным декартовым произведением представлений $B_1$, ..., $B_n$ совпадает с декартовым произведением $B_1 \times \ldots \times B_n$ множеств $B_1$, ..., $B_n$. В этом месте доказательства нас не интересует алгебраическая структура на множестве $B_1 \times \ldots \times B_n$.

[4.6]Равенство (4.7.3) утверждает, что мы отождествляем базис представления $M$ с множеством $B_1 \times \ldots \times B_n$.

[4.7] Я рассматриваю формирование элементов представления из элементов базиса согласно теореме 6.1.4. Теорема 4.7.11 требует выполнения условий (4.7.4), (4.7.5).

Доказательство. Из равенства (4.7.5) следует, что

$$f(c) \circ (b_1, ..., f_i(a) \circ b_i, ..., b_n) = (b_1, ..., f_i(c) \circ (f_i(a) \circ b_i), ..., b_n)$$
$$= (b_1, ..., (f_i(c) \circ f_i(a)) \circ b_i, ..., b_n)$$
(4.7.12)
$$= (f(c) \circ f(a)) \circ (b_1, ..., b_i, ..., b_n)$$
$$= f(c) \circ (f(a) \circ (b_1, ..., b_i, ..., b_n))$$

Равенство (4.7.11) следует из равенства (4.7.12).                              ⊙

Лемма 4.7.8. *Для любого* $c \in A$ *эндоморфизм* $f(c)$ $\Omega_2$-*алгебры* $M$ *согласовано с эквивалентностью* $N$.

Доказательство. Утверждение леммы следует из лемм 4.7.6, 4.7.7 и определения 3.3.2.                              ⊙

Из леммы 4.7.8 и теоремы 3.3.3 следует, что на множестве $^*M/N$ определена $\Omega_1$-алгебра. Рассмотрим диаграмму

Согласно лемме 4.7.8, из условия

$$j \circ b_1 = j \circ b_2$$

следует

$$j \circ (f(a) \circ b_1) = j \circ (f(a) \circ b_2)$$

Следовательно, преобразование $F(a)$ определено корректно и

(4.7.13)                              $$F(a) \circ j = j \circ f(a)$$

Если $\omega \in \Omega_1(p)$, то мы положим

$$(F(a_1)...F(a_p)\omega) \circ (J \circ b) = J \circ ((f(a_1)...f(a_p)\omega) \circ b)$$

Следовательно, отображение $F$ является представлением $\Omega_1$-алгебры $A$. Из (4.7.13) следует, что $j$ является приведенным морфизмом представлений $f$ и $F$.

Рассмотрим коммутативную диаграмму

(4.7.14)

Из коммутативности диаграммы (4.7.14) и равенства (4.7.3) следует, что

$$(4.7.15) \qquad g_1 \circ (b_1, ..., b_n) = j \circ (b_1, ..., b_n)$$

Из равенств (4.7.3), (4.7.4), (4.7.5) следует

$$(4.7.16) \qquad g_1 \circ (b_1, ..., b_{i\cdot 1}...b_{i\cdot p}\omega, ..., b_n)$$
$$= (g_1 \circ (b_1, ..., b_{i\cdot 1}, ..., b_n))...(g_1 \circ (b_1, ..., b_{i\cdot p}, ..., b_n))\omega$$

$$(4.7.17) \qquad g_1 \circ (b_1, ..., f_i(a) \circ b_i, ..., b_n) = f(a) \circ (g_1 \circ (b_1, ..., b_i, ..., b_n))$$

Из равенств (4.7.16) и (4.7.17) следует, что отображение $g_1$ является приведённым полиморфизмом представлений $f_1, ..., f_n$.

Поскольку $B_1 \times ... \times B_n$ - базис представления $M$ $\Omega_1$-алгебры $A$, то, согласно теореме 6.2.10, для любого представления

$$A \xrightarrow{\ *\ } V$$

и любого приведённого полиморфизма

$$g_2 : B_1 \times ... \times B_n \longrightarrow V$$

существует единственный морфизм представлений $k : M \to V$, для которого коммутативна следующая диаграмма

$$(4.7.18)$$

Так как $g_2$ - приведённый полиморфизм, то $\ker k \supseteq N$.

Согласно теореме 3.4.8 отображение $j$ универсально в категории морфизмов представления $f$, ядро которых содержит $N$. Следовательно, определён морфизм представлений

$$h : M/N \to V$$

для которого коммутативна диаграмма

$$(4.7.19)$$

Объединяя диаграммы (4.7.14), (4.7.18), (4.7.19), получим коммутативную диаграмму

Так как $\operatorname{Im} g_1$ порождает $M/N$, то отображение $h$ однозначно определено. $\qquad \square$

Согласно доказательству теоремы 4.7.5

$$B_1 \otimes ... \otimes B_n = M/N$$

Для $d_i \in A_i$ будем записывать

(4.7.20)                    $j \circ (d_1, ..., d_n) = d_1 \otimes ... \otimes d_n$

Из равенств (4.7.15), (4.7.20) следует, что

(4.7.21)                    $g_1 \circ (d_1, ..., d_n) = d_1 \otimes ... \otimes d_n$

Теорема 4.7.9. *Отображение*

$$(x_1, ..., x_n) \in B_1 \times ... \times B_n \to x_1 \otimes ... \otimes x_n \in B_1 \otimes ... \otimes B_n$$

*является полиморфизмом.*

Доказательство. Теорема является следствием определений 4.4.4, 4.7.1. □

Теорема 4.7.10. *Пусть $B_1$, ..., $B_n$ - $\Omega_2$-алгебры. Пусть*

$$f : B_1 \times ... \times B_n \to B_1 \otimes ... \otimes B_n$$

*приведенный полиморфизм, определённый равенством*

(4.7.22)                    $f \circ (b_1, ..., b_n) = b_1 \otimes ... \otimes b_n$

*Пусть*

$$g : B_1 \times ... \times B_n \to V$$

*приведенный полиморфизм в $\Omega_2$-алгебру $V$. Существует морфизм представлений*

$$h : B_1 \otimes ... \otimes B_n \to V$$

*такой, что диаграмма*

*коммутативна.*

Доказательство. Равенство (4.7.22) следует из равенств (4.7.3) и (4.7.20). Существование отображения $h$ следует из определения 4.7.1 и построений, выполненных при доказательстве теоремы 4.7.5. □

Теорема 4.7.11. *Пусть*

$$b_k \in B_k \quad k = 1, ..., n \quad b_{i \cdot 1}, ..., b_{i \cdot p} \in B_i \quad \omega \in \Omega_2(p) \quad a \in A$$

*Тензорное произведение дистрибутивно относительно операции $\omega$*

(4.7.23)
$$b_1 \otimes ... \otimes (b_{i \cdot 1}...b_{i \cdot p}\omega) \otimes ... \otimes b_n$$
$$= (b_1 \otimes ... \otimes b_{i \cdot 1} \otimes ... \otimes b_n)...(b_1 \otimes ... \otimes b_{i \cdot p} \otimes ... \otimes b_n)\omega$$

*Представление мультипликативной $\Omega_1$-группы $A$ в тензорном произведении определено равенством*

$$(4.7.24) \qquad b_1 \otimes \ldots \otimes (f_i(a) \circ b_i) \otimes \ldots \otimes b_n = f(a) \circ (b_1 \otimes \ldots \otimes b_i \otimes \ldots \otimes b_n)$$

Доказательство. Равенство (4.7.23) является следствием равенства (4.7.16) и определения (4.7.21). Равенство (4.7.24) является следствием равенства (4.7.17) и определения (4.7.21). $\qquad \square$

## 4.8. Ассоциативность тензорного произведения

Пусть $A$ является мультипликативной $\Omega_1$-группой. Пусть $B_1$, $B_2$, $B_3$ - $\Omega_2$-алгебры. Пусть для $k = 1, 2, 3$

$$f_k : A \longrightarrow B_k$$

эффективное представление мультипликативной $\Omega_1$-группы $A$ в $\Omega_2$-алгебре $B_k$.

Лемма 4.8.1. *Для заданного значения $x_3 \in B_3$, отображение*

$$(4.8.1) \qquad h_{12} : (B_1 \otimes B_2) \times B_3 \to B_1 \otimes B_2 \otimes B_3$$

*определённое равенством*

$$(4.8.2) \qquad h_{12}(x_1 \otimes x_2, x_3) = x_1 \otimes x_2 \otimes x_3$$

*является приведенным морфизмом представления $B_1 \otimes B_2$ в представление $B_1 \otimes B_2 \otimes B_3$ .*

Доказательство. Согласно теореме 4.7.9, для заданного значения $x_3 \in B_3$, отображение

$$(4.8.3) \qquad (x_1, x_2, x_3) \in B_1 \times B_2 \times B_3 \to x_1 \otimes x_2 \otimes x_3 \in B_1 \otimes B_2 \otimes B_3$$

является полиморфизмом по переменным $x_1 \in B_1$, $x_2 \in B_2$. Следовательно, для заданного значения $x_3 \in B_3$, лемма является следствием теоремы 4.7.10.
$\qquad \square$

Лемма 4.8.2. *Для заданного значения $x_{12} \in B_1 \otimes B_2$ отображение $h_{12}$ является приведенным морфизмом представления $B_3$ в представление $B_1 \otimes B_2 \otimes B_3$ .*

Доказательство. Согласно теореме 4.7.9 и равенству (4.7.21), для заданного значения $x_1 \in B_1$, $x_2 \in B_2$, отображение

$$(4.8.4) \qquad (x_1 \otimes x_2, x_3) \in B_1 \times B_2 \times B_3 \to x_1 \otimes x_2 \otimes x_3 \in B_1 \otimes B_2 \otimes B_3$$

является морфизмом по переменной $x_3 \in B_3$. Следовательно, теорема является следсвием равенства (4.4.16) и теоремы 4.5.9. $\qquad \square$

Лемма 4.8.3. *Существует приведенный морфизм представлений*

$$h : (B_1 \otimes B_2) \otimes B_3 \to B_1 \otimes B_2 \otimes B_3$$

Доказательство. Согласно леммам 4.8.1, 4.8.2 и определению 4.4.4, отображение $h_{12}$ является приведенным полиморфизмом представлений. Утверждение леммы является следствием теоремы 4.7.10. $\qquad \square$

Лемма 4.8.4. *Существует приведенный морфизм представлений*

$$g : B_1 \otimes B_2 \otimes B_3 \to (B_1 \otimes B_2) \otimes B_3$$

Доказательство. Отображение

$$(x_1, x_2, x_3) \in B_1 \times B_2 \times B_3 \to (x_1 \otimes x_2) \otimes x_3 \in (B_1 \otimes B_2) \otimes B_3$$

является полиморфизмом по переменным $x_1 \in B_1$, $x_2 \in B_2$, $x_3 \in B_3$. Следовательно, лемма является следствием теоремы 4.7.10.                    □

Теорема 4.8.5.

(4.8.5)        $(A_1 \otimes A_2) \otimes A_3 = A_1 \otimes (A_2 \otimes A_3) = A_1 \otimes A_2 \otimes A_3$

Доказательство. Согласно лемме 4.8.3, существует приведённый морфизм представлений

$$h : (B_1 \otimes B_2) \otimes B_3 \to B_1 \otimes B_2 \otimes B_3$$

Согласно лемме 4.8.4, существует приведённый морфизм представлений

$$g : B_1 \otimes B_2 \otimes B_3 \to (B_1 \otimes B_2) \otimes B_3$$

Следовательно, приведённые морфизмы представлений $h$, $g$ являются изоморфизмами, откуда следует равенство

(4.8.6)        $(B_1 \otimes B_2) \otimes B_3 = B_1 \otimes B_2 \otimes B_3$

Аналогично мы можем доказать равенство

$$B_1 \otimes (B_2 \otimes B_3) = B_1 \otimes B_2 \otimes B_3$$

□

Замечание 4.8.6. *Очевидно, что структура $\Omega_2$-алгебр $(B_1 \otimes B_2) \otimes B_3$, $B_1 \otimes B_2 \otimes B_3$ слегка различна. Мы записываем равенство* (4.8.6), *опираясь на соглашение 4.7.3 и это позволяет нам говорить об ассоциативности тензорного произведения представлений.*                    □

# Представление мультипликативной $\Omega$-группы

## 5.1. Представление мультипликативной $\Omega$-группы

Согласованность произведения в мультипликативной $\Omega$-группе $G$ и соответствующих преобразований представления $f$ позволяет нам рассмотреть больше деталей представления $f$. Однако конструкция, рассмотренная в теореме 4.5.7, не полна в случае некоммутативного произведения.

Если для заданного представления

$$g : A_1 \longrightarrow A_2$$

для любых $A_1$-чисел $a_1$, $b_1$, однозначно определено $A_1$-число $c_1$ такое, что

$$f(c_1) = f(a_1) \circ f(b_1)$$

то какой формат произведения мы должны выбрать:

$$(5.1.1) \qquad\qquad c_1 = a_1 * b_1$$

или

$$(5.1.2) \qquad\qquad c_1 = b_1 * a_1$$

ПРИМЕР 5.1.1. *Пусть*

$$\overline{\overline{e}} = \begin{pmatrix} e_1 & ... & e_n \end{pmatrix}$$

*базис левого векторного пространства $V$ над ассоциативной алгеброй с делением $A$. Мы можем представить произвольный вектор $\overline{v} \in V$ как $*_*$-произведение матриц*

$$(5.1.3) \qquad \overline{v} = v {}^*_* \overline{\overline{e}} = \begin{pmatrix} v^1 \\ ... \\ v^n \end{pmatrix} {}^*_* \begin{pmatrix} e_1 & ... & e_n \end{pmatrix}$$

*где*

$$v = \begin{pmatrix} v^1 \\ ... \\ v^n \end{pmatrix}$$

*матрица координат вектора $\overline{v}$ относительно базиса $\overline{\overline{e}}$.*

Рассмотрим однотранзитивное действие группы $G$, определённое равенством

$$(5.1.4) \qquad g^*{}_*\overline{\overline{e}} = \begin{pmatrix} g_1^1 & ... & g_n^1 \\ ... & ... & ... \\ g_1^n & ... & g_n^n \end{pmatrix} {}^*{}_* \begin{pmatrix} e_1 & ... & e_n \end{pmatrix}$$

где мы отождествляем $G$-число $g$ и невырожденную матрицу

$$\begin{pmatrix} g_1^1 & ... & g_n^1 \\ ... & ... & ... \\ g_1^n & ... & g_n^n \end{pmatrix}$$

Действие группы $G$ на многообразии базисов является представлением, так как верно равенство

$$(5.1.5) \qquad g_1{}^*{}_*(g_2{}^*{}_*\overline{\overline{e}}) = (g_1{}^*{}_*g_2){}^*{}_*\overline{\overline{e}}$$

Пусть

$$(5.1.6) \qquad v_i = \begin{pmatrix} v_i^1 \\ ... \\ v_i^n \end{pmatrix}$$

матрица координат вектора $\overline{v}$ относительно базиса $\overline{\overline{e}}_i$, $i = 1,\ 2,\ 3$. Тогда

$$(5.1.7) \qquad \overline{v} = v_1{}^*{}_*\overline{\overline{e}}_1 = v_2{}^*{}_*\overline{\overline{e}}_2 = v_3{}^*{}_*\overline{\overline{e}}_3$$

Пусть $G$-число $g_1$ отображает базис $\overline{\overline{e}}_1$ в базис $\overline{\overline{e}}_2$

$$(5.1.8) \qquad \overline{\overline{e}}_2 = g_1{}^*{}_*\overline{\overline{e}}_1$$

Пусть $G$-число $g_2$ отображает базис $\overline{\overline{e}}_2$ в базис $\overline{\overline{e}}_3$

$$(5.1.9) \qquad \overline{\overline{e}}_3 = g_2{}^*{}_*\overline{\overline{e}}_2$$

Равенство

$$(5.1.10) \qquad \overline{\overline{e}}_3 = (g_2{}^*{}_*g_1){}^*{}_*\overline{\overline{e}}_1$$

является следствием равенств (5.1.8), (5.1.9). Равенство

$$(5.1.11) \qquad v_1{}^*{}_*\overline{\overline{e}}_1 = v_2{}^*{}_*g_1{}^*{}_*\overline{\overline{e}}_1$$

является следствием равенств (5.1.7), (5.1.8). Равенство

$$(5.1.12) \qquad v_1 = v_2{}^*{}_*g_1$$

является следствием равенства (5.1.11), так как координаты вектора $\overline{v}$ определены однозначно относительно базиса $\overline{\overline{e}}_1$. Равенство

$$(5.1.13) \qquad v_2 = v_1{}^*{}_*g_1^{-1}$$

является следствием равенства (5.1.12). Аналогично, равенство

$$(5.1.14) \qquad v_3 = v_2{}^*{}_*g_2^{-1}$$

является следствием равенств (5.1.7), (5.1.9) и равенство

$$(5.1.15) \qquad v_3 = v_1{}^*{}_*(g_2{}^*{}_*g_1)^{-1}$$

*является следствием равенств* (5.1.7), (5.1.10). *Равенство*

$$(5.1.16) \qquad v_3 = v_1{}^*{}_* g_1^{-1}{}^*{}_* g_2^{-1}$$

*является следствием равенств* (5.1.13), (5.1.14).                                        □

Пример 5.1.2. *Пусть $V$ - левый модуль над кольцом $D$. Это значит, что определено представление*

$$f : D \longrightarrow{\!\!\!*\!\!\!}\longrightarrow V \qquad f(d) : v \to d\,v$$

*такое, что*

$$(d_1 + d_2)v = d_1 v + d_2 v$$
$$d(v_1 + v_2) = dv_1 + dv_2$$
$$d_1(d_2 v) = (d_1 d_2)v$$

*Отображение*

$$w : V \to D$$

*называется аддитивным, если*

$$w(v_1 + v_2) = w(v_1) + w(v_2)$$

*Мы пользуемся записью*

$$(w, v) = w(v)$$

*для образа аддитивного отображения. Мы определим сумму аддитивных отображений равенством*

$$(w_1 + w_2, v) = (w_1, v) + (w_2, v)$$

*Нетрудно показать, что множество $W$ аддитивных отображений является абелевой группой.*

*Мы определим отображение*

$$h : D \longrightarrow{\!\!\!*\!\!\!}\longrightarrow W \qquad h(d) : w \to w\,d$$

*равенством*

$$(wd, v) = (w, dv)$$

*Из равенств*

$$((w_1 + w_2)d, v) = (w_1 + w_2, dv) = (w_1, dv) + (w_2, dv)$$
$$= (w_1 d, v) + (w_2 d, v)$$
$$= (w_1 d + w_2 d, v)$$
$$(w(d_1 + d_2), v) = (w, (d_1 + d_2)v) = (w, d_1 v + d_2 v)$$
$$= (w, d_1 v) + (w, d_2 v) = (wd_1, v) + (wd_2, v)$$
$$= (wd_1 + wd_2, v)$$
$$((wd_1)d_2, v) = (wd_1, d_2 v) = (w, d_1(d_2 v)) = (w, (d_1 d_2)v)$$
$$(5.1.17) \qquad = (w(d_1 d_2), v)$$

*следует, что отображение $h$ является представлением группы $G$. Однако мы можем записать равенство* (5.1.17) *в виде*

$$((h(d_2) \circ h(d_1)(w)), v) = ((h(d_2)h(d_1)(w)), v) = (h(d_1 d_2)(w), v)$$

*откуда следует, что отображение $h$ не является гомоморфизмом группы $G$.*

                                                                                            □

Мы предполагаем, что преобразования представления мультипликативной Ω-группы $A_1$ могут действовать на $A_2$-числа либо слева, либо справа. В этом случае нам достаточно ограничиться произведением (5.1.1) в мультипликативной Ω-группе $A_1$. Таким образом, концепция представления мультипликативной Ω-группы состоит в том, что в каком порядке мы перемножаем элементы мультипликативной Ω-группы, в том же порядке перемножаются соответствующие преобразования представления. Эта точка зрения отражена в примере 5.1.2. Мы также видим, что мы должны изменить формат записи, прежде чем мы можем перейти к этой точке зрения. Вместо того, чтобы рассматривать $f \in \mathrm{End}(\Omega_2; A_2)$, как отображение

$$f : a_2 \in A_2 \to f(a_2) \in A_2$$

мы должны рассматривать эндоморфизм $f$ как оператор.

ОПРЕДЕЛЕНИЕ 5.1.3. *Пусть* $\mathrm{End}(\Omega_2, A_2)$ *- мультипликативная Ω-группа с произведением*[5.1]

$$(f, g) \to f \bullet g$$

*Пусть эндоморфизм $f$ действует на $A_2$-число а слева. Мы будем пользоваться записью*

$$(5.1.18) \qquad\qquad f(a_2) = f \bullet a_2$$

*Пусть $A_1$ - мультипликативная Ω-группа с произведением*

$$(a, b) \to a * b$$

*Мы будем называть гомоморфизм мультипликативной Ω-группы*

$$(5.1.19) \qquad\qquad f : A_1 \to \mathrm{End}(\Omega_2, A_2)$$

**левосторонним представлением** *мультипликативной Ω-группы $A_1$ или* **левосторонним** *$A_1$-**представлением** в $\Omega_2$-алгебре $A_2$, если отображение $f$ удовлетворяет условиям*

$$(5.1.20) \qquad\qquad f(a_1 * b_1) \bullet a_2 = (f(a_1) \bullet f(b_1)) \bullet a_2$$

*Мы будем отождествлять $A_1$-число $a_1$ с его образом $f(a_1)$ и записывать левостороннее преобразование, порождённое $A_1$-числом $a_1$, в форме*

$$a_2' = f(a_1) \bullet a_2 = a_1 * a_2$$

*В этом случае равенство* (5.1.20) *принимает вид*

$$(5.1.21) \qquad\qquad f(a_1 * b_1) \bullet a_2 = (a_1 * b_1) * a_2$$

*Отображение*

$$(a_1, a_2) \in A_1 \times A_2 \to a_1 * a_2 \in A_2$$

*порождённое левосторонним представлением $f$, называется* **левосторонним произведением** *$A_2$-числа $a_2$ на $A_1$-число $a_1$.* □

---

[5.1] Очень часто произведение в мультипликативной Ω-группе $\mathrm{End}(\Omega_2, A_2)$ совпадает с суперпозицией эндоморфизмов

$$f \bullet g = f \circ g$$

Однако, как мы увидим в примере 5.2.5, произведение в мультипликативной Ω-группе $\mathrm{End}(\Omega_2, A_2)$ может отличаться от суперпозиции эндоморфизмов. Согласно определению 4.6.13, мы можем рассматривать две операции произведения в универсальной алгебре $A$.

Пусть

$$f : A_2 \to A_2$$
$$g : A_2 \to A_2$$

эндоморфизмы $\Omega_2$-алгебры $A_2$. Пусть произведение в мультипликативной $\Omega$-группе $\mathrm{End}(\Omega_2, A_2)$ является композицией эндоморфизмов. Так как произведение отображений $f$ и $g$ определено в том же порядке, как эти отображения действуют на $A_2$-число, то мы можем рассматривать равенство

(5.1.22) $$(f \circ g) \circ a = f \circ (g \circ a)$$

как **закон ассоциативности**, который позволяет записывать равенство (5.1.22) без использования скобок

$$f \circ g \circ a = f \circ (g \circ a) = (f \circ g) \circ a$$

а также записать равенство (5.1.20) в виде

(5.1.23) $$f(a_1 * b_1) \circ a_2 = f(a_1) \circ f(b_1) \circ a_2$$

Из равенства (5.1.21) следует, что

(5.1.24) $$(a_1 * b_1) * a_2 = a_1 * (b_1 * a_2)$$

Мы можем рассматривать равенство (5.1.24) как **закон ассоциативности**,

ЗАМЕЧАНИЕ 5.1.4. *Пусть отображение*

$$f : A_1 \longrightarrow\!\!\!* A_2$$

*является левосторонним представлением мультипликативной $\Omega$-группы $A_1$ в $\Omega_2$-алгебре $A_2$. Пусть отображение*

$$g : B_1 \longrightarrow\!\!\!* B_2$$

*является левосторонним представлением мультипликативной $\Omega$-группы $B_1$ в $\Omega_2$-алгебре $B_2$. Пусть отображение*

$$(r_1 : A_1 \to B_1, \; r_2 : A_2 \to B_2)$$

*является морфизмом представлений. Мы будем пользоваться записью*

$$r_2(a_2) = r_2 \circ a_2$$

*для образа $A_2$-числа $a_2$ при отображении $r_2$. Тогда мы можем записать равенство (3.2.3) следующим образом*

$$r_2 \circ (a_1 * a_2) = r_1(a_1) * (r_2 \circ a_2)$$

□

ОПРЕДЕЛЕНИЕ 5.1.5. *Пусть $\mathrm{End}(\Omega_2, A_2)$ - мультипликативная $\Omega$-группа с произведением*[5.2]

$$(f, g) \to f \bullet g$$

---

[5.2] Очень часто произведение в мультипликативной $\Omega$-группе $\mathrm{End}(\Omega_2, A_2)$ совпадает с суперпозицией эндоморфизмов

$$f \bullet g = f \circ g$$

Однако, как мы увидим в примере 5.2.5, произведение в мультипликативной $\Omega$-группе $\mathrm{End}(\Omega_2, A_2)$ может отличаться от суперпозиции эндоморфизмов. Согласно определению 4.6.13, мы можем рассматривать две операции произведения в универсальной алгебре $A$.

*Пусть эндоморфизм $f$ действует на $A_2$-число $a$ справа. Мы будем пользоваться записью*

$$(5.1.25) \qquad\qquad f(a_2) = a_2 \bullet f$$

*Пусть $A_1$ - мультипликативная Ω-группа с произведением*

$$(a, b) \to a * b$$

*Мы будем называть гомоморфизм мультипликативной Ω-группы*

$$(5.1.26) \qquad\qquad f : A_1 \to \text{End}(\Omega_2, A_2)$$

**правосторонним представлением** *мультипликативной Ω-группы $A_1$ или* **правосторонним** $A_1$-**представлением** *в Ω₂-алгебре $A_2$, если отображение $f$ удовлетворяет условиям*

$$(5.1.27) \qquad\qquad a_2 \bullet f(a_1 * b_1) = a_2 \bullet (f(a_1) \bullet f(b_1))$$

*Мы будем отождествлять $A_1$-число $a_1$ с его образом $f(a_1)$ и записывать правостороннее преобразование, порождённое $A_1$-числом $a_1$, в форме*

$$a_2' = a_2 \bullet f(a_1) = a_2 * a_1$$

*В этом случае равенство* (5.1.27) *принимает вид*

$$(5.1.28) \qquad\qquad a_2 \bullet f(a_1 * b_1) = a_2 * (a_1 * b_1)$$

*Отображение*

$$(a_1, a_2) \in A_1 \times A_2 \to a_2 * a_1 \in A_2$$

*порождённое правосторонним представлением $f$, называется* **правосторонним произведением** *$A_2$-числа $a_2$ на $A_1$-число $a_1$.*                    □

Пусть

$$f : A_2 \to A_2$$
$$g : A_2 \to A_2$$

эндоморфизмы Ω₂-алгебры $A_2$. Пусть произведение в мультипликативной Ω-группе $\text{End}(\Omega_2, A_2)$ является композицией эндоморфизмов. Так как произведение отображений $f$ и $g$ определено в том же порядке, как эти отображения действуют на $A_2$-число, то мы можем рассматривать равенство

$$(5.1.29) \qquad\qquad a \circ (g \circ f) = (a \circ g) \circ f$$

как **закон ассоциативности**, который позволяет записывать равенство (5.1.29) без использования скобок

$$a \circ g \circ f = (a \circ g) \circ f = a \circ (g \circ f)$$

а также записать равенство (5.1.27) в виде

$$(5.1.30) \qquad\qquad a_2 \circ f(a_1 * b_1) = a_2 \circ f(a_1) \circ f(b_1)$$

Из равенства (5.1.28) следует, что

$$(5.1.31) \qquad\qquad a_2 * (a_1 * b_1) = (a_2 * a_1) * b_1$$

Мы можем рассматривать равенство (5.1.31) как **закон ассоциативности**,

ЗАМЕЧАНИЕ 5.1.6. *Пусть отображение*

$$f : A_1 \longrightarrow_* A_2$$

*является правосторонним представлением мультипликативной Ω-группы $A_1$ в Ω$_2$-алгебре $A_2$. Пусть отображение*

$$g : B_1 \longrightarrow_* B_2$$

*является правосторонним представлением мультипликативной Ω-группы $B_1$ в Ω$_2$-алгебре $B_2$. Пусть отображение*

$$(r_1 : A_1 \to B_1, \ r_2 : A_2 \to B_2)$$

*является морфизмом представлений. Мы будем пользоваться записью*

$$r_2(a_2) = r_2 \circ a_2$$

*для образа $A_2$-числа $a_2$ при отображении $r_2$. Тогда мы можем записать равенство* (3.2.3) *следующим образом*

$$r_2 \circ (a_2 * a_1) = (r_2 \circ a_2) * r_1(a_1)$$

□

Если мультипликативная Ω-группа $A_1$ - абелевая, то нет разницы между левосторонним и правосторонним представлениями.

ОПРЕДЕЛЕНИЕ 5.1.7. *Пусть $A_1$ - абелевая мультипликативная Ω-группа. Мы будем называть гомоморфизм мультипликативной Ω-группы*

$$f : A_1 \to \text{End}(\Omega_2, A_2) \tag{5.1.32}$$

**представлением** *мультипликативной Ω-группы $A_1$ или $A_1$-**представлением** в Ω$_2$-алгебре $A_2$, если отображение $f$ удовлетворяет условиям*

$$f(a_1 * b_1) \bullet a_2 = (f(a_1) \bullet f(b_1)) \bullet a_2 \tag{5.1.33}$$

□

Обычно мы отождествляем представление абелевой мультипликативной Ω-группы $A_1$ и левостороннее представление мультипликативной Ω-группы $A_1$. Однако, если это необходимо нам, мы можем отождествить представление абелевой мультипликативной Ω-группы $A_1$ и правостороннее представление мультипликативной Ω-группы $A_1$.

Из анализа примера 5.1.2 следует, что выбор между левосторонним и правосторонним представлением зависит от рассматриваемой модели. Так как левостороннее представление и правостороннее представление опирается на гомоморфизм Ω-группы, то верно следующее утверждение.

ТЕОРЕМА 5.1.8 (принцип двойственности для представления мультипликативной Ω-группы). *Любое утверждение, справедливое для левостороннего представления мультипликативной Ω-группы $A_1$, будет справедливо для правостороннего представления мультипликативной Ω-группы $A_1$, если мы будем пользоваться правосторонним произведением на $A_1$-число $a_1$ вместо левостороннего произведения на $A_1$-число $a_1$.* □

Замечание 5.1.9. *Если* $\Omega_1$*-алгебра не является мультипликативной* $\Omega$*-группой, то мы не можем сказать, действует ли представление слева или справа. В этом случае мы сохраним функциональную запись* $f(a_1)(a_2)$ *для представления* $\Omega_1$*-алгебры.*        □

Из анализа равенств (5.1.15), (5.1.16) следует, что действие группы $G$ на множестве координат вектора $\overline{v}$ (пример 5.1.1) не соответствует ни левостороннему, ни правостороннему представлению. следует, что у нас есть два выбора. Мы согласны, что в мультипликативной Ω-группе $A_1$ мы можем определить оба варианта произведения: (5.1.1) и (5.1.2) - с целью согласовать произведение в мультипликативной Ω-группе $A_1$ и произведение преобразований представления мультипликативной Ω-группы $A_1$. Эта точка зрения отражена в определениях 5.1.10, 5.1.11.

Определение 5.1.10. *Левостороннее представление*

$$f : A_1 \longrightarrow\!\!\ast\!\longrightarrow A_2$$

*называется* **ковариантным**, *если равенство*

$$a_1 * (b_1 * a_2) = (a_1 * b_1) * a_2$$

*верно.*        □

Определение 5.1.11. *Левостороннее представление*

$$f : A_1 \longrightarrow\!\!\ast\!\longrightarrow A_2$$

*называется* **контравариантным**, *если равенство*

$$(5.1.34) \qquad a_1^{-1} * (b_1^{-1} * a_2) = (b_1 * a_1)^{-1} * a_2$$

*верно.*        □

Если тип представления не указан, мы будем предполагать, что представление ковариантно. Из равенств (5.1.15), (5.1.16) следует, что действие группы $G$ на множестве координат вектора $\overline{v}$ (пример 5.1.1) является контравариантным правосторонним представлением.

Насколько велика разница между ковариантным и контравариантным представлениями. Поскольку

$$(b_1 * a_1)^{-1} = a_1^{-1} * b_1^{-1}$$

то равенство

$$(5.1.35) \qquad a_1^{-1} * (b_1^{-1} * a_2) = (a_1^{-1} * b_1^{-1}) * a_2$$

является следствием равенства (5.1.34). Из равенства (5.1.35) следует, что мы можем рассматривать контравариантное представление группы $G$ как ковариантное представление группы $G$, порождённое $G$-числами вида $a^{-1}$. Так же как в примере 5.1.1, мы рассматриваем два согласованных представления группы $G$

$$f : G \longrightarrow\!\!\ast\!\longrightarrow A_2$$

$$h : G \longrightarrow\!\!\ast\!\longrightarrow B_2$$

причём $G$-число $g$ порождает преобразование

$$a_1 \in G : a_2 \in A_2 \to a_1 * a_2 \in A_2$$

в универсальной алгебре $A_2$ и преобразование

$$a_1 \in G : b_2 \in B_2 \rightarrow a_1^{-1} * b_2 \in B_2$$

в универсальной алгебре $B_2$.

## 5.2. Левый и правый сдвиги

Теорема 5.2.1. *Произведение*

$$(a, b) \rightarrow a * b$$

*в мультипликативной $\Omega$-группе $A$ определяет два различных представления.*

- **Левый сдвиг**

$$a' = L(b) \circ a = b * a$$

*является левосторонним представлением мультипликативной $\Omega$-группы $A$ в $\Omega$-алгебре $A$*

(5.2.1) $$L(c * b) = L(c) \circ L(b)$$

- **Правый сдвиг**

$$a' = a \circ R(b) = a * b$$

*является правосторонним представлением мультипликативной $\Omega$-группы $A$ в $\Omega$-алгебре $A$*

(5.2.2) $$R(b * c) = R(b) \circ R(c)$$

Доказательство. Согласно определению 4.5.4, левый и правый сдвиги являются эндоморфизмами $\Omega$-алгебры $A$. Согласно определению 4.5.4, мы можем определить $\Omega$-алгебру на множестве левых сдвигов. Согласно определению мультипликативной группы, [5.3] равенство $a_1 = a_2$ является следствием равенства

$$L(a_1) \circ x = a_1 * x = a_2 * x = L(a_2) \circ x$$

для любого $x$. В частности, равенство (5.2.1) является следствием равенства

$$L(c * b) \circ a = (c * b) * a = c * (b * a) = L(c) \circ (L(b) \circ a) = L(c) \circ L(b) \circ a$$

Следовательно, отображение

$$a \in A \rightarrow L(a)$$

является левовосторонним представлением мультипликативной $\Omega$-группы $A$ в $\Omega$-алгебре $A$. Аналогичное рассуждение верно для правого сдвига. $\qquad\square$

Ассоциативная $D$-алгебра является мультипликативной $\Omega$-группой. Неассоциативная $D$-алгебра $A$ не является $\Omega$-группой, так как относительно произведения $A$ является группоидом. Однако нас также будет интересовать представление неассоциативной $D$-алгебры.

Определение 5.2.2. *Пусть произведение*

$$c_1 = a_1 * b_1$$

_________________________

[5.3] Смотри, например, определение на страницах [2]-17, [2]-21.

*является операцией $\Omega_1$-алгебры $A$. Положим $\Omega = \Omega_1 \setminus \{*\}$. Если $\Omega_1$-алгебра $A$ является группоидом относительно произведения и для любой операции $\omega \in \Omega(n)$ умножение дистрибутивно относительно операции $\omega$*

$$a * (b_1 ... b_n \omega) = (a * b_1) ... (a * b_n) \omega$$

$$(b_1 ... b_n \omega) * a = (b_1 * a) ... (b_n * a) \omega$$

*то $\Omega_1$-алгебра $A$ называется $\Omega$-**группоидом**.*                            $\square$

Мы будем пользоваться тем же форматом записи для представления $\Omega$-группоида, что мы пользуемся для представления мультипликативной $\Omega$-группы.

Теорема 5.2.3. *Произведение в неассоциативном $\Omega$-группоиде $A$ определяет два различных представления.*

- *Левый сдвиг*
$$a' = L(b) \circ a = b * a$$
  *является представлением $\Omega$-алгебры $A$ в $\Omega$-алгебре $A$.*
- *Правый сдвиг*
$$a' = a \circ R(b) = a * b$$
  *является представлением $\Omega$-алгебры $A$ в $\Omega$-алгебре $A$.*

Доказательство. Согласно определению 4.5.4, левый и правый сдвиги являются эндоморфизмами $\Omega$-алгебры $A$. Согласно определению 5.2.2, мы можем определить $\Omega$-алгебру на множестве левых сдвигов. Следовательно, отображение

$$a \in A \to L(a)$$

является представлением $\Omega$-алгебры $A$ в $\Omega$-алгебре $A$.              $\square$

Теорема 5.2.4. *Пусть*

$$L : A \xrightarrow{\quad * \quad} A$$

*представление неассоциативного $\Omega$-группоида $A$ в $\Omega$-алгебре $A$. Тогда на множестве $\mathrm{End}(\Omega, A)$ определена операция произведения, отличная от суперпозиции эндоморфизмов.*

Доказательство. Рассмотрим отображение

$$L : A \to \mathrm{End}(\Omega, A) \quad L(a) : b \to ab$$

Поскольку произведение в $\Omega$-группоиде $A$ не ассоциативно, то, вообще говоря,

$$L(a) \circ (L(b) \circ c) = a * (b * c) \neq (a * b) * c = L(a * b) \circ c$$

Следовательно, $L(ab) \neq L(a) \circ L(b)$.                            $\square$

Согласно теореме 5.2.1, если $A$ - мультипликативная $\Omega$-группа, то равенство (5.2.1) гарантирует, что левый сдвиг порождает левостороннее представление мультипликативной $\Omega$-группы $A$ в $\Omega$-алгебре $A$. Согласно теореме 5.2.4 это равенство не верно в неассоциативном $\Omega$-группоиде $A$. Однако теоремы 5.2.3, 5.2.4 не отвечают на вопрос о возможности рассмотрения левостороннего представления неассоциативного $\Omega$-группоида $A$ в $\Omega$-алгебре $A$. Согласно примеру 5.2.5, существует возможность подобного представления, даже если произведение в $\Omega$-группоиде неассоциативно.

Пример 5.2.5. *Пусть A - алгебре Ли. Произведение*[5.4] $[a, b]$ *A-чисел a, b удовлетворяет равенству*

$$(5.2.3) \qquad\qquad [a, b] = -[b, a]$$

*а также тождеству Ли*

$$(5.2.4) \qquad\qquad [c, [b, a]] + [b, [a, c]] + [a, [c, b]] = 0$$

*Левый сдвиг на алгебре Ли A определён равенством*

$$(5.2.5) \qquad\qquad L(b) \circ a = [b, a]$$

*Из равенства* (5.2.5) *следует, что*

$$(5.2.6) \qquad L(c) \circ L(b) \circ a = L(c) \circ (L(b) \circ a) = [c, [b, a]]$$

*Равенство*

$$(5.2.7) \qquad \begin{aligned} L(c) \circ L(b) \circ a - L(b) \circ L(c) \circ a &= [c, [b, a]] - [b, [c, a]] \\ &= [c, [b, a]] + [b, [a, c]] \end{aligned}$$

*является следствием равенств* (5.2.3), (5.2.6). *Равенство*

$$(5.2.8) \qquad\qquad [c, [b, a]] + [b, [a, c]] = -[a, [c, b]] = [[c, b], a]$$

*является следствием равенств* (5.2.3), (5.2.4). *Равенство*

$$(5.2.9) \qquad L(c) \circ L(b) \circ a - L(b) \circ L(c) \circ a = L([c, b]) \circ a$$

*является следствием равенств* (5.2.5), (5.2.7), (5.2.8).

*Если я определю произведение Ли*

$$[L(c), L(b)] \circ a = L(c) \circ L(b) \circ a - L(b) \circ L(c) \circ a$$

*на множестве левых сдвигов, то равенство* (5.2.9) *принимает вид*

$$(5.2.10) \qquad\qquad [L(c), L(b)] \circ a = L([c, b]) \circ a$$

*Следовательно, алгебра Ли A с произведением* $[a, b]$ *порождает представление в векторном пространстве A.* □

## 5.3. Орбита представления мультипликативной Ω-группы

Теорема 5.3.1. *Пусть отображение*

$$f : A_1 \longrightarrow\!\!\!* \longrightarrow A_2$$

*является левосторонним представлением мультипликативной $\Omega_1$-группы $A_1$ и e - единица мультипликативной $\Omega_1$-группы $A_1$. Тогда*

$$f(e) = \delta$$

*где $\delta$ - тождественное преобразование $\Omega_2$-алгебры $A_2$.*

Доказательство. Теорема является следствием равенства

$$f(a) = f(a * e) = f(a) \circ f(e)$$

для любого $A_1$-числа. □

---

Теорема 5.3.2. *Пусть отображение*

$$g : A_1 \xrightarrow{\quad * \quad} A_2$$

*является левосторонним представлением мультипликативной Ω-группы $A_1$. Для любого $g \in A_1$ преобразование $f(g)$ имет обратное отображение и удовлетворяет равенству*

(5.3.1)                          $$f(g^{-1}) = f(g)^{-1}$$

Доказательство. Пусть $e$ - единица мультипликативной Ω-группы $A_1$ и $\delta$ - тождественное преобразование множества $A_2$. На основании (5.1.20) и теоремы 5.3.1 мы можем записать

$$u = \delta \circ u = f(gg^{-1}) \circ u = f(g) \circ f(g^{-1}) \circ u$$

Это завершает доказательство.                                        □

Определение 5.3.3. *Пусть $A_1$ является Ω-группоидом с произведением*

$$(a, b) \to a * b$$

*Пусть отображение*

$$f : A_1 \xrightarrow{\quad * \quad} A_2$$

*является левосторонним представлением Ω-группоида $A_1$ в $\Omega_2$-алгебре $A_2$. Для любого $a_2 \in A_2$, мы определим* **орбиту представления** *Ω-группоида $A_1$ как множество*

$$A_1 * a_2 = \{b_2 = a_1 * a_2 : a_1 \in A_1\}$$

                                                                    □

Определение 5.3.4. *Пусть $A_1$ является Ω-группоидом с произведением*

$$(a, b) \to a * b$$

*Пусть отображение*

$$f : A_1 \xrightarrow{\quad * \quad} A_2$$

*является правосторонним представлением Ω-группоида $A_1$ в $\Omega_2$-алгебре $A_2$. Для любого $a_2 \in A_2$, мы определим* **орбиту представления** *Ω-группоида $A_1$ как множество*

$$a_2 * A_1 = \{b_2 = a_2 * a_1 : a_1 \in A_1\}$$

                                                                    □

Теорема 5.3.5. *Пусть отображение*

$$f : A_1 \xrightarrow{\quad * \quad} A_2$$

*является левосторонним представлением мультипликативной Ω-группы $A_1$. Then $a_2 \in A_1 * a_2$.*

Доказательство. Согласно теореме 5.3.1

$$a_2 = e * a_2 = f(e) \circ a_2$$

                                                                    □

Теорема 5.3.6. *Пусть*

$$L : A \longrightarrow\!\!\!*\!\!\!\longrightarrow A$$

*представление алгебры Ли, порождённое множеством левых сдвигов. Тогда* $a \notin [A, a]$.

Доказательство. Теорема является следствием отсутствия единицы в алгебре Ли. Кроме того, самый простой пример алгебры Ли - это множество векторов трёх мерного пространства на котором определена операция векторного произведения. Очевидно, что не существует вектора $b$ такого, что

$$a = b \times a$$

$\square$

Теорема 5.3.7. *Пусть отображение*

$$f : A_1 \longrightarrow\!\!\!*\!\!\!\longrightarrow A_2$$

*является левосторонним представлением мультипликативной Ω-группы* $A_1$.
*Если*

(5.3.2) $$b_2 \in A_1 * a_2$$

*то*

(5.3.3) $$A_1 * a_2 = A_1 * b_2$$

Доказательство. Из (5.3.2) следует существование $a_1 \in A_1$ такого, что

(5.3.4) $$b_2 = a_1 * a_2$$

Если $c_2 \in A_1 * b_2$, то существует $b_1 \in A_1$ такой, что

(5.3.5) $$c_2 = b_1 * b_2$$

Подставив (5.3.4) в (5.3.5), мы получим

(5.3.6) $$c_2 = b_1 * a_1 * a_2$$

На основании (5.1.20) из (5.3.6) следует, что $c_2 \in A_1 * a_2$. Таким образом,

(5.3.7) $$A_1 * b_2 \subseteq A_1 * a_2$$

На основании (5.3.1) из (5.3.4) следует, что

(5.3.8) $$a_2 = a_1^{-1} * b_2$$

Равенство (5.3.8) означает, что $a_2 \in A_1 * b_2$ и, следовательно,

(5.3.9) $$A_1 * a_2 \subseteq A_1 * b_2$$

Равенство (5.3.3) является следствием утверждений (5.3.7), (5.3.9). $\square$

Таким образом, левостороннее представление мультипликативной Ω-группы $A_1$ в $\Omega_2$-алгебре $A_2$ порождает отношение эквивалентности $S$ и орбита $A_1 * a_2$ является классом эквивалентности. Мы будем пользоваться обозначением $A_2/A_1$ для фактор множества $A_2/S$ и мы будем называть это множество **пространством орбит левостороннего представления** $f$.

## 5.4. Представление в $\Omega$-группе

Теорема 5.4.1. *Мы будем называть ядром неэффективности левосторон-него представления мультипликативной $\Omega$-группы $A_1$ в $\Omega_2$-алгебре $A_2$ мно-жество*

$$K_f = \{a_1 \in A_1 : f(a_1) = \delta\}$$

*Ядро неэффективности левостороннего представления - это подгруппа муль-типликативной группы $A_1$.*

Доказательство. Допустим $f(a_1) = \delta$ и $f(a_2) = \delta$. Тогда

$$f(a_1 * a_2) = f(a_1) \bullet f(a_2) = \delta$$
$$f(a_1^{-1}) = (f(a_1))^{-1} = \delta$$

$\square$

Теорема 5.4.2. *Левостороннее представление мультипликативной $\Omega$-группы $A_1$ в $\Omega_2$-алгебре $A_2$ **эффективно** тогда и только тогда, когда ядро неэффективности $K_f = \{e\}$.*

Доказательство. Утверждение является следствием определения 3.1.2 и теоремы 5.4.1. $\square$

Теорема 5.4.3. *Если представление*

$$f : A_1 \longrightarrow\!\!\!* A_2$$

*мультипликативной $\Omega$-группы $A_1$ в $\Omega_2$-алгебре $A_2$ не эффективно, мы можем перейти к эффективному заменив мультипликативную $\Omega$-группу $A_1$ муль-типликативной $\Omega$-группой $A_1' = A_1/K_f$.*

Доказательство. Пусть операция $\omega \in \Omega(n)$. Чтобы доказать теорему, мы должны показать, что равенство

(5.4.1)     $f(a_1...a_n\omega) = f(b_1...b_n\omega)$

является следствием утверждения $f(a_1) = f(b_1)$, ..., $f(a_n) = f(b_n)$. Действи-тельно, равенство (5.4.1) является следствием равенства

$$f(a_1...a_n\omega) = f(a_1)...f(a_n)\omega = f(b_1)...f(b_n)\omega = f(b_1...b_n\omega)$$

$\square$

Теорема 5.4.3 означает, что мы можем изучать только эффективное дей-ствие.

## 5.5. Однотранзитивное правостороннее представление группы

Теорема 5.5.1. *Пусть отображение*

$$g : A_1 \longrightarrow\!\!\!* A_2$$

*является левосторонним представлением мультипликативной $\Omega$-группы $A_1$ в $\Omega_2$-алгебре $A_2$. **Малая группа** или **группа стабилизации** элемента $a_2 \in A_2$ - это множество*

$$A_{1a_2} = \{a_1 \in A_1 : a_1 * a_2 = a_2\}$$

*Представление $f$ **свободно** тогда и только тогда, когда для любого $a_2 \in A_2$ группа стабилизации $A_{1a_2} = \{e\}$.*

Доказательство. Согласно определению 3.1.4, представление $f$ свободно тогда и только тогда, когда равенство $a_1 = b_1$ является следсвием утверждения

(5.5.1) $$f(a_1) = f(b_1)$$

Равенство (5.5.1) эквивалентно равенству

(5.5.2) $$f(b_1^{-1} * a_1) = \delta$$

Равенство $a_1 = b_1$ является следствием утверждения (5.5.2) тогда и только тогда, когда для любого $a_2 \in A_2$ группа стабилизации $A_{1a_2} = \{e\}$. $\quad\square$

Теорема 5.5.2. *Пусть отображение*

$$f : A_1 \xrightarrow{\quad * \quad} A_2$$

*является свободным левосторонним представлением мультипликативной $\Omega$-группы $A_1$ в $\Omega_2$-алгебре $A_2$. Существует взаимно однозначное соответствие между любыми двумя орбитами представления, а также между орбитой представления и мультипликативной $\Omega$-группой $A_1$.*

Доказательство. Допустим для точки $a_2 \in A_2$ существуют $a_1, b_1 \in A_1$

(5.5.3) $$a_1 * a_2 = b_1 * a_2$$

Умножим обе части равенства (5.5.3) на $a_1^{-1}$

$$a_2 = a_1^{-1} * b_1 * a_2$$

Поскольку представление свободное, $a_1 = b_1$. Теорема доказана, так как мы установили взаимно однозначное соответствие между орбитой и мультипликативной $\Omega$-группой $A_1$. $\quad\square$

Теорема 5.5.3. *Левостороннее представление*

$$g : A_1 \xrightarrow{\quad * \quad} A_2$$

*мультипликативной $\Omega$-группы $A_1$ в $\Omega_2$-алгебре $A_2$ **однотранзитивно** тогда и только тогда, когда для любых $a_2, b_2 \in A_2$ существует одно и только одно $a_1 \in A_1$ такое, что $a_2 = a_1 * b_2$.*

Доказательство. Следствие определений 3.1.2 и 3.1.8. $\quad\square$

Теорема 5.5.4. *Если существует однотранзитивное представление*

$$f : A_1 \xrightarrow{\quad * \quad} A_2$$

*мультипликативной $\Omega$-группы $A_1$ в $\Omega_2$-алгебре $A_2$, то мы можем однозначно определить координаты на $A_2$, пользуясь $A_1$-числами.*

*Если $f$ - левостороннее однотранзитивное однотранзитивное представление, то $f(a)$ эквивалентно левому сдвигу $L(a)$ на группе $A_1$. Если $f$ - правостороннее однотранзитивное представление, то $f(a)$ эквивалентно правому сдвигу $R(a)$ на группе $A_1$.*

Доказательство. Пусть $f$ - левостороннее представление. Мы выберем $A_2$-число $a_2$ и определим координаты $A_2$-числа $b_2$ как $A_1$-число $a_1$ такое, что

$$b_2 = a_1 * a_2 = (a_1 * e) * a_2 = (L(a_1) \circ e) * a_2$$

Координаты, определённые таким образом, однозначны с точностью до выбора $A_2$-числа $a_2$, так как действие эффективно. Для левостороннего однотранзитивного представления мы будем также пользоваться записью

$$f(a_1) \bullet a_2 = L(a_1) \circ a_2 = (L(a_1) \circ e) * a_2$$

Мы пользуемся записью $L(a_1) \circ a_2$ для левостороннего однотранзитивного представления $f$ так как, согласно теореме 5.2.1, произведение левых сдвигов совпадает с их композицией.

Пусть $f$ - правостороннее однотранзитивное представление. Мы выберем $A_2$-число $a_2$ и определим координаты $A_2$-числа $b_2$ как $A_1$-число $a_1$ такое, что

$$b_2 = a_2 * a_1 = a_2 * (e * a_1) = a_2 * (e \circ R(a_1))$$

Координаты, определённые таким образом, однозначны с точностью до выбора $A_2$-числа $a_2$, так как действие эффективно. Для правостороннего однотранзитивного представления мы будем также пользоваться записью

$$a_2 \bullet f(a_1) = a_2 \circ R(a_1) = a_2 * (e \circ R(a_1))$$

Мы пользуемся записью $a_2 \circ R(a_1)$ для правостороннего однотранзитивного представления $f$ так как, согласно теореме 5.2.1, произведение правых сдвигов совпадает с их композицией.                                                                   □

Определение 5.5.5. *Мы будем называть $\Omega_2$-алгебру $A_2$* **однородным пространством** *мультипликативной $\Omega$-группы $A_1$, если существует однотранзитивное левостороннее представление*

$$f : A_1 \xrightarrow{\quad * \quad} A_2$$

□

Теорема 5.5.6. *Свободное левостороннее представление мультипликативной $\Omega$-группы $A_1$ в $\Omega_2$-алгебре $A_2$ однотранзитивно на орбите.*

Доказательство. Следствие теоремы 5.5.2.                                   □

Теорема 5.5.7. *Правый и левый сдвиги на мультипликативной $\Omega$-группе $A_1$ перестановочны.*

Доказательство. Теорема является следствием ассоциативности произведения в мультипликативной $\Omega$-группе $A_1$

$$(L(a) \circ c) \circ R(b) = (a * c) * b = a * (c * b) = L(a) \circ (c \circ R(b))$$

□

Теорема 5.5.7 может быть сформулирована следующим образом.

Теорема 5.5.8. *Пусть $A_1$ - мультипликативная $\Omega$-группа. Для любого $a_1 \in A_1$ отображение $L(a_1)$ является автоморфизмом представления $R$.*

Доказательство. Согласно теореме 5.5.7

(5.5.4)                          $$L(a_1) \circ R(b_1) = R(b_1) \circ L(a_1)$$

Равенство (5.5.4) совпадает с равенством (3.2.2) из определения 3.2.2 при условии $r_1 = id$, $r_2 = L(a_1)$.                                             □

Теорема 5.5.9. *Пусть левостороннее $A_1$-представление $f$ на $\Omega_2$-алгебре $A_2$ однотранзитивно. Тогда мы можем однозначно определить однотранзитивное правостороннее $A_1$-представление $h$ на $\Omega_2$-алгебре $A_2$ такое, что диаграмма*

$$\begin{array}{ccc} A_2 & \xrightarrow{\ h(a_1)\ } & A_2 \\ \downarrow{\scriptstyle f(b_1)} & & \downarrow{\scriptstyle f(b_1)} \\ A_2 & \xrightarrow[\ h(a_1)\ ]{} & A_2 \end{array}$$

*коммутативна для любых $a_1$, $b_1 \in A_1$.* [5.5]

Доказательство. Мы будем пользоваться групповыми координатами для $A_2$-чисел $a_2$. Тогда согласно теореме 5.5.4 мы можем записать левый сдвиг $L(a_1)$ вместо преобразования $f(a_1)$.

Пусть $a_2$, $b_2 \in A_2$. Тогда мы можем найти одно и только одно $a_1 \in A_1$ такое, что

$$b_2 = a_2 * a_1 = a_2 \circ R(a_1)$$

Мы предположим

$$h(a) = R(a)$$

Существует $b_1 \in A_1$ такое, что

$$c_2 = f(b_1) \bullet a_2 = L(b_1) \circ a_2 \qquad d_2 = f(b_1) \bullet b_2 = L(b_1) \circ b_2$$

Согласно теореме 5.5.7 диаграмма

(5.5.5)
$$\begin{array}{ccc} a_2 & \xrightarrow{\ h(a_1)=R(a_1)\ } & b_2 \\ \downarrow{\scriptstyle f(b_1)=L(b_1)} & & \downarrow{\scriptstyle f(b_1)=L(b_1)} \\ c_2 & \xrightarrow[\ h(a_1)=R(a_1)\ ]{} & d_2 \end{array}$$

коммутативна.

Изменяя $b_1$ мы получим, что $c_2$ - это произвольное $A_2$-число.

Мы видим из диаграммы, что, если $a_2 = b_2$, то $c_2 = d_2$ и следовательно $h(e) = \delta$. С другой стороны, если $a_2 \neq b_2$, то $c_2 \neq d_2$ потому, что левостороннее $A_1$-представление $f$ однотранзитивно. Следовательно правостороннее $A_1$-представление $h$ эффективно.

Таким же образам мы можем показать, что для данного $c_2$ мы можем найти $a_1$ такое, что $d_2 = c_2 \bullet h(a_1)$. Следовательно правостороннее $A_1$-представление $h$ однотранзитивно.

В общем случае, произведение преобразований левостороннего $A_1$-представления $f$ не коммутативно и следовательно правосторонним $A_1$-представление $h$ отлично от левостороннего $A_1$-представления $f$. Таким же образом мы можем создать левостороннее $A_1$-представление $f$, пользуясь правосторонним $A_1$-представлением $h$. □

Мы будем называть представления $f$ и $h$ **парными представлениями** мультипликативной $\Omega$-группы $A_1$.

---

[5.5]Это утверждение можно также найти в [4].

Замечание 5.5.10. *Очевидно, что преобразования $L(a)$ и $R(a)$ отличаются, если мультипликативная Ω-группа $A_1$ неабелева. Тем не менее, они являются отображениями на. Теорема 5.5.9 утверждает, что, если оба представления правого и левого сдвига существуют на множестве $A_2$, то мы можем определить два перестановочных представления на множестве $A_2$. Только правый или левый сдвиг не может представлять оба типа представления. Чтобы понять почему это так, мы можем изменить диаграмму (5.5.5) и предположить*

$$h(a_1) \bullet a_2 = L(a_1) \circ a_2 = b_2$$

*вместо*

$$a_2 \bullet h(a_1) = a_2 \circ R(a_1) = b_2$$

*и проанализировать, какое выражение $h(a_1)$ имеет в точке $c_2$. Диаграмма*

*эквивалентна диаграмме*

*и мы имеем $d_2 = b_1 b_2 = b_1 a_1 a_2 = b_1 a_1 b_1^{-1} c_2$. Следовательно*

$$h(a_1) \bullet c_2 = (b_1 a_1 b_1^{-1}) c_2$$

*Мы видим, что представление $h$ зависит от его аргумента.* □

Теорема 5.5.11. *Пусть $f$ и $h$ - парные преставления мультипликативной Ω-группы $A_1$. Для любого $a_1 \in A_1$ отображение $h(a_1)$ является автоморфизмом представления $f$.*

Доказательство. Следствие теорем 5.5.8 и 5.5.9. □

Вопрос 5.5.12. *Существует ли морфизм представлений из $L$ в $L$, отличный от автоморфизма $R(a_1)$? Если мы положим*

$$r_1(a_1) = c_1 a_1 c_1^{-1}$$

$$r_2(a_1) \circ a_2 = c_1 a_2 a_1 c_1^{-1}$$

*то нетрудно убедиться, что отображение $(r_1 \ r_2(a_1))$ является морфизмом представлений из $L$ в $L$. Но это отображение не является автоморфизмом представления $L$, так как $r_1 \neq \mathrm{id}$.* □

# Базис представления универсальной алгебры

## 6.1. Множество образующих представления

ОПРЕДЕЛЕНИЕ 6.1.1. *Пусть*

$$f : A_1 \xrightarrow{\quad *\quad} A_2$$

*представление $\Omega_1$-алгебры $A_1$ в $\Omega_2$-алгебре $A_2$. Множество $B_2 \subset A_2$ называется **стабильным множеством представления** $f$, если $f(a)(m) \in B_2$ для любых $a \in A_1$, $m \in B_2$.* □

Мы также будем говорить, что множество $A_2$ стабильно относительно представления $f$.

ТЕОРЕМА 6.1.2. *Пусть*

$$f : A_1 \xrightarrow{\quad *\quad} A_2$$

*представление $\Omega_1$-алгебры $A_1$ в $\Omega_2$-алгебре $A_2$. Пусть множество $B_2 \subset A_2$ является подалгеброй $\Omega_2$-алгебры $A_2$ и стабильным множеством представления $f$. Тогда существует представление*

$$f_{B_2} : A_1 \xrightarrow{\quad *\quad} B_2$$

*такое, что $f_{B_2}(a) = f(a)|_{B_2}$. Представление $f_{B_2}$ называется **подпредставлением** представления $f$.*

ДОКАЗАТЕЛЬСТВО. Пусть $\omega_1$ - $n$-арная операция $\Omega_1$-алгебры $A_1$. Тогда для любых $a_1, ..., a_n \in A_1$ и любого $b \in B_2$

$$(f_{B_2}(a_1)...f_{B_2}(a_n)\omega_1)(b) = (f(a_1)...f(a_n)\omega_1)(b) = f(a_1...a_n\omega_1)(b)$$
$$= f_{B_2}(a_1...a_n\omega_1)(b)$$

Пусть $\omega_2$ - $n$-арная операция $\Omega_2$-алгебры $A_2$. Тогда для любых $b_1, ..., b_n \in B_2$ и любого $a \in A_1$

$$f_{B_2}(a)(b_1)...f_{B_2}(a)(b_n)\omega_2 = f(a)(b_1)...f(a)(b_n)\omega_2 = f(a)(b_1...b_n\omega_2)$$
$$= f_{B_2}(a)(b_1...b_n\omega_2)$$

Утверждение теоремы доказано. □

Из теоремы 6.1.2 следует, что если $f_{B_2}$ - подпредставление представления $f$, то отображение

$$(\mathrm{id} : A \to A, \mathrm{id}_{B_2} : B_2 \to A_2)$$

является морфизмом представлений.

ТЕОРЕМА 6.1.3. *Множество*[6.1] $\mathcal{B}_f$ *всех подпредставлений представления $f$ порождает систему замыканий на $\Omega_2$-алгебре $A_2$ и, следовательно, является полной структурой.*

ДОКАЗАТЕЛЬСТВО. Пусть $(K_\lambda)_{\lambda \in \Lambda}$ - семейство подалгебр $\Omega_2$-алгебры $A_2$, стабильных относительно представления $f$. Операцию пересечения на множестве $\mathcal{B}_f$ мы определим согласно правилу

$$\bigcap f_{K_\lambda} = f_{\cap K_\lambda}$$

Операция пересечения подпредставлений определена корректно. $\cap K_\lambda$ - подалгебра $\Omega_2$-алгебры $A_2$. Пусть $m \in \cap K_\lambda$. Для любого $\lambda \in \Lambda$ и для любого $a \in A_1$, $f(a)(m) \in K_\lambda$. Следовательно, $f(a)(m) \in \cap K_\lambda$. Следовательно, $\cap K_\lambda$ - стабильное множество представления $f$. $\qquad \square$

Обозначим соответствующий оператор замыкания через $J[f]$. Таким образом, $J[f, X]$ является пересечением всех подалгебр $\Omega_2$-алгебры $A_2$, содержащих $X$ и стабильных относительно представления $f$.

ТЕОРЕМА 6.1.4. *Пусть*[6.2]

$$g : A_1 \longrightarrow A_2$$

*представление $\Omega_1$-алгебры $A_1$ в $\Omega_2$-алгебре $A_2$. Пусть $X \subset A_2$. Определим подмножество $X_k \subset A_2$ индукцией по $k$.*

6.1.4.1: $X_0 = X$
6.1.4.2: $x \in X_k => x \in X_{k+1}$
6.1.4.3: $x_1 \in X_k$, ..., $x_n \in X_k$, $\omega \in \Omega_2(n) => x_1...x_n\omega \in X_{k+1}$
6.1.4.4: $x \in X_k$, $a \in A => f(a)(x) \in X_{k+1}$

*Тогда*

$$(6.1.1) \qquad \bigcup_{k=0}^{\infty} X_k = J[f, X]$$

ДОКАЗАТЕЛЬСТВО. Если положим $U = \cup X_k$, то по определению $X_k$ имеем $X_0 \subset J[f, X]$, и если $X_k \subset J[f, X]$, то $X_{k+1} \subset J[f, X]$. По индукции следует, что $X_k \subset J[f, X]$ для всех $k$. Следовательно,

$$(6.1.2) \qquad U \subset J[f, X]$$

Если $a \in U^n$, $a = (a_1, ..., a_n)$, где $a_i \in X_{k_i}$, и если $k = \max\{k_1, ..., k_n\}$, то $a_1...a_n\omega \in X_{k+1} \subset U$. Следовательно, $U$ является подалгеброй $\Omega_2$-алгебры $A_2$.

Если $m \in U$, то $m \in X_k$ для некоторого $k$. Следовательно, $f(a)(m) \in X_{k+1} \subset U$ для любого $a \in A_1$. Следовательно, $U$ - стабильное множество представления $f$.

Так как $U$ - подалгебра $\Omega_2$-алгебры $A_2$ и стабильное множество представления $f$, то определено подпредставление $f_U$. Следовательно,

$$(6.1.3) \qquad J[f, X] \subset U$$

---

[6.1] Это определение аналогично определению структуры подалгебр ([14], стр. 93, 94). Вообще говоря, в этой и последующих теоремах этой главы необходимо рассмотреть структуру универсальных алгебр $A_1$ и $A_2$. Так как основная задача этой главы - это изучение структуры представления, я сознательно упростил теоремы, чтобы детали не заслоняли основные утверждения. Более подробно эта тема будет раскрыта в главе 8, где теоремы будут сформулированы в общем виде.

[6.2] Утверждение теоремы аналогично утверждению теоремы 5.1, [14], страница 94.

Из (6.1.2), (6.1.3), следует $J[f, X] = U$. $\qquad\square$

ОПРЕДЕЛЕНИЕ 6.1.5. $J[f, X]$ *называется* **подпредставлением**, *порождённым множеством* $X$, *а* $X$ *- множеством образующих подпредставления* $J[f, X]$. *В частности,* **множеством образующих** *представления* $f$ *будет такое подмножество* $X \subset A_2$, *что* $J[f, X] = A_2$. $\qquad\square$

Следующее определение является следствием теоремы 6.1.4.

ОПРЕДЕЛЕНИЕ 6.1.6. *Пусть* $X \subset A_2$. *Для любого* $m \in J[f, X]$ *существует* $\Omega_2$**-слово** $w[f, X, m]$, *определённое согласно следующему правилам.*

6.1.6.1: *Если* $m \in X$, *то* $m$ - $\Omega_2$-слово.

6.1.6.2: *Если* $m_1, ..., m_n$ - $\Omega_2$-слова и $\omega \in \Omega_2(n)$, *то* $m_1...m_n\omega$ - $\Omega_2$-слово.

6.1.6.3: *Если* $m$ - $\Omega_2$-слово и $a \in A_1$, *то* $f(a)(m)$ - $\Omega_2$-слово.

*Мы будем отождествлять элемент* $m \in J[f, X]$ *и соответствующее ему* $\Omega_2$-слово, *выражая это равенством*

$$m = w[f, X, m]$$

*Аналогично, для произвольного множества* $B \subset J[f, X]$ *рассмотрим множество* $\Omega_2$-слов [6.3]

$$w[f, X, B] = \{w[f, X, m] : m \in B\}$$

*Мы будем также пользоваться записью*

$$w[f, X, B] = (w[f, X, m], m \in B)$$

*Обозначим* $w[f, X]$ **множество** $\Omega_2$**-слов представления** $J[f, X]$. $\qquad\square$

ТЕОРЕМА 6.1.7. $w[f, X, X] = X$.

ДОКАЗАТЕЛЬСТВО. Теорема является следствием утверждения 6.1.6.1. $\square$

ТЕОРЕМА 6.1.8. *Пусть* $X$, $Y$ - *множества образующих представления*

$$f : A_1 \xrightarrow{\quad * \quad} A_2$$

*Пусть* $w[f, X, m]$ - $\Omega_2$-слово $A_2$-числа $m$ *относительно множества образующих* $X$. *Пусть* $w[f, Y, X]$ - *множество* $\Omega_2$-слов множества $X$ *относительно множества образующих* $Y$. *Если в слове* $w[f, X, m]$ *вместо каждого* $x \in X$ *подставить его образ* $w[f, Y, x]$, *то мы получим* $\Omega_2$-слово $w[f, Y, m]$ $A_2$-числа $m$ *относительно множества образующих* $Y$.

*Преобразование* $\Omega_2$-слов

$$w[f, X, m] \to w[f, Y, m]$$

$$w[f, Y, m] = w[f, Y, X] \circ w[f, X, m]$$

*называется суперпозицией координат.*

---

[6.3]Выражение $w[f, X, m]$ является частным случаем выражения $w[f, X, B]$, а именно

$$w[f, X, \{m\}] = \{w[f, X, m]\}$$

Доказательство. Мы будем доказывать теорему индукцией по сложности $\Omega_2$-слова.

Если $m \in X$, то $w[f, X, m] = m$. Если вместо $m$ подставить его образ $w[f, Y, m]$, то мы получим $\Omega_2$-слово $w[f, Y, m]$ $A_2$-числа $m$ относительно множества образующих $Y$.

Пусть $\Omega_2$-слово $w[f, X, m]$ $A_2$-числа $m$ имеет вид

$$(6.1.4) \qquad w[f, X, m] = w[f, X, m_1]...w[f, X, m_n]\omega$$

где $\omega \in \Omega_2(n)$ и для каждого $A_2$-числа $m_i$ мы определили отображение

$$w[f, X, m_i] \to w[f, Y, m_i]$$

Согласно утверждению 6.1.6.2 выражение

$$w[f, Y, m_1]...w[f, Y, m_n]\omega$$

является $\Omega_2$-словом $w[f, Y, m]$ $A_2$-числа $m$ относительно множества образующих $Y$. Следовательно, мы определили отображение

$$w[f, X, m] \to w[f, Y, m]$$

для $A_2$-числа $m$.

Пусть $\Omega_2$-слово $w[f, X, m]$ $A_2$-числа $m$ имеет вид

$$(6.1.5) \qquad w[f, X, m] = f(a)(w[f, X, m_1])$$

где для $A_2$-числа $m_1$ мы определили отображение

$$w[f, X, m_1] \to w[f, Y, m_1]$$

Согласно утверждению 6.1.6.3 выражение

$$f(a)(w[f, Y, m_1])$$

является $\Omega_2$-словом $w[f, Y, m]$ $A_2$-числа $m$ относительно множества образующих $Y$. Следовательно, мы определили отображение

$$w[f, X, m] \to w[f, Y, m]$$

для $A_2$-числа $m$. $\qquad\qquad\qquad\qquad\qquad\qquad\qquad\qquad\qquad\qquad\square$

Выбор $\Omega_2$-слова относительно множества образующих $X$ неоднозначен. Поэтому, если $\Omega_2$-число имеет различные $\Omega_2$-слова, то мы, чтобы их отличать, будем пользоваться индексами: $w[f, X, m]$, $w_1[f, X, m]$, $w_2[f, X, m]$.

Определение 6.1.9. *Множество образующих $X$ представления $f$ порождает отношение эквивалентности*

$$\rho[f, X] = \{(w[f, X, m], w_1[f, X, m]) : m \in A_2\}$$

*на множестве $\Omega_2$-слов.* $\qquad\qquad\qquad\qquad\qquad\qquad\qquad\qquad\qquad\square$

Согласно определению 6.1.9, два $\Omega_2$-слова относительно множества образующих $X$ представления $f$ эквивалентны тогда и только тогда, когда они соответствуют одному и тому же $A_2$-числу. Когда мы будем записывать равенство двух $\Omega_2$-слов относительно множества образующих $X$ представления $f$, мы будем иметь в виду, что это равенство верно с точностью до отношения эквивалентности $\rho[f, X]$.

Теорема 6.1.10. *Пусть*

$$f : A_1 \longrightarrow A_2$$

*представление* $\Omega_1$-*алгебры* $A_1$ *в* $\Omega_2$-*алгебре* $A_2$. *Пусть*

$$g : A_1 \longrightarrow B_2$$

*представление* $\Omega_1$-*алгебры* $A_1$ *в* $\Omega_2$-*алгебре* $B_2$. *Пусть* $X$ - *множество образующих представления* $f$. *Пусть*

$$R : A_2 \to B_2$$

*приведенный морфизм представления*[6.4] *и* $X' = R(X)$. *Приведенный морфизм* $R$ *представления порождает отображение* $\Omega_2$-*слов*

$$w[f \to g, X, R] : w[f, X] \to w[g, X']$$

*такое, что*

6.1.10.1: *Если* $m \in X$, $m' = R(m)$, *то*

$$w[f \to g, X, R](m) = m'$$

6.1.10.2: *Если*

$$m_1, ..., m_n \in w[f, X]$$

$$m'_1 = w[f \to g, X, R](m_1) \quad ... \quad m'_n = w[f \to g, X, R](m_n)$$

*то для операции* $\omega \in \Omega_2(n)$ *справедливо*

$$w[f \to g, X, R](m_1...m_n\omega) = m'_1...m'_n\omega$$

6.1.10.3: *Если*

$$m \in w[f, X] \quad m' = w[f \to g, X, R](m) \quad a \in A_1$$

*то*

$$w[f \to g, X, R](f(a)(m)) = g(a)(m')$$

Доказательство. Утверждения 6.1.10.1, 6.1.10.2 справедливы в силу определения приведенного морфизма $R$. Утверждение 6.1.10.3 является следствием равенства (3.4.5). $\square$

Замечание 6.1.11. *Пусть*

$$R : A_2 \to B_2$$

*приведенный морфизм представления. Пусть*

$$m \in J[f, X] \quad m' = R(m) \quad X' = R(X)$$

*Теорема 6.1.10 утверждает, что* $m' \in J[g, X']$. *Теорема 6.1.10 также утверждает, что* $\Omega_2$-*слово, представляющее* $m$, *относительно* $X$ *и* $\Omega_2$-*слово, представляющее* $m'$, *относительно* $X'$ *формируются согласно одному и тому же*

---

*алгоритму. Это позволяет рассматривать множество $\Omega_2$-слов $w[g, X', m']$ как отображение*

(6.1.6) $\qquad W[f, X, m] : (g, X') \to (g, X') \circ W[f, X, m] = w[g, X', m']$

*где*

$$X' = R(X) \quad m' = R(m)$$

*для некоторого приведенного морфизма $R$.*

*Если $f = g$, то вместо отображения (6.1.6) мы будем рассматривать отображение*

$$W[f, X, m] : X' \to X' \circ W[f, X, m] = w[f, X', m']$$

*такое, что, если для некоторого эндоморфизма $R$*

$$X' = R(X) \quad m' = R(m)$$

*то*

$$W[f, X, m](X') = X' \circ W[f, X, m] = w[f, X', m'] = m'$$

*Отображение $W[f, X, m]$ называется **координатами** $A_2$-**числа** $m$ относительно множества $X$. Аналогично, мы можем рассмотреть координаты множества $B \subset J[f, X]$ относительно множества $X$*

$$W[f, X, B] = \{W[f, X, m] : m \in B\} = (W[f, X, m], m \in B)$$

*Обозначим*

$$W[f, X] = \{W[f, X, m] : m \in J[f, X]\} = (W[f, X, m], m \in J[f, X])$$

**множество координат представления $J[f, X]$.** $\qquad\square$

Теорема 6.1.12. *На множестве координат $W[f, X]$ определена структура $\Omega_2$-алгебры.*

Доказательство. Пусть $\omega \in \Omega_2(n)$. Тогда для любых $m_1, ..., m_n \in J[f, X]$ положим

(6.1.7) $\qquad W[f, X, m_1]...W[f, X, m_n]\omega = W[f, X, m_1...m_n\omega]$

Согласно замечанию 6.1.11, из равенства (6.1.7) следует

(6.1.8) $\qquad \begin{aligned} X \circ (W[f, X, m_1]...W[f, X, m_n]\omega) &= X \circ W[f, X, m_1...m_n\omega] \\ &= w[f, X, m_1...m_n\omega] \end{aligned}$

Согласно правилу 6.1.6.2, из равенства (6.1.8) следует

(6.1.9) $\qquad \begin{aligned} X \circ (W[f, X, m_1]...&W[f, X, m_n]\omega) \\ = w[f, X, m_1]...&w[f, X, m_n]\omega \\ = (X \circ W[f, X, m_1])...&(X \circ W[f, X, m_n])\omega \end{aligned}$

Из равенства (6.1.9) следует корректность определения (6.1.7) операции $\omega$ на множестве координат. $\qquad\square$

Теорема 6.1.13. *Определено представление $\Omega_1$-алгебры $A_1$ в $\Omega_2$-алгебре $W[f, X]$.*

Доказательство. Пусть $a \in A_1$. Тогда для любого $m \in J[f, X]$ положим

(6.1.10) $$f(a)(W[f, X, m]) = W[f, X, f(a)(m)]$$

Согласно замечанию 6.1.11, из равенства (6.1.10) следует

(6.1.11) $$X \circ (f(a)(W[f, X, m])) = X \circ W[f, X, f(a)(m)] = w[f, X, f(a)(m)]$$

Согласно правилу 6.1.6.3, из равенства (6.1.11) следует

(6.1.12) $$X \circ (f(a)(W[f, X, m])) = f(a)(w[f, X, m]) = f(a)(X \circ W[f, X, m])$$

Из равенства (6.1.12) следует корректность определения (6.1.10) представления $\Omega_1$-алгебры $A_1$ в $\Omega_2$-алгебре $W[f, X]$.   □

Теорема 6.1.14. *Пусть*

$$f : A_1 \longrightarrow_* A_2$$

*представление $\Omega_1$-алгебры $A_1$ в $\Omega_2$-алгебре $A_2$. Пусть*

$$g : A_1 \longrightarrow_* B_2$$

*представление $\Omega_1$-алгебры $A_1$ в $\Omega_2$-алгебре $B_2$. Для заданных множеств $X \subset A_2$, $X' \subset B_2$, пусть отображение*

$$R_1 : X \to X'$$

*согласовано со структурой представления $f$, т. е.*

$$\omega \in \Omega_2(n) \ x_1, \ ..., \ x_n, \ x_1...x_n\omega \in X, \ R_1(x_1...x_n\omega) \in X'$$
$$=> R_1(x_1...x_n\omega) = R_1(x_1)...R_1(x_n)\omega$$
$$x \in X, \ a \in A, \ R_1(f(a)(x)) \in X'$$
$$=> R_1(f(a)(x)) = g(a)(R_1(x))$$

*Рассмотрим отображение $\Omega_2$-слов*

$$w[f \to g, X, X', R_1] : w[f, X] \to w[g, X']$$

*удовлетворяющее условиям 6.1.10.1, 6.1.10.2, 6.1.10.3, и такое, что*

$$x \in X => w[f \to g, X, X', R_1](x) = R_1(x)$$

*Существует единственное отображение*

$$R : A_2 \to B_2$$

*определённое правилом*

$$R(m) = w[f \to g, X, X', R_1](w[f, X, m])$$

*которое является приведенным морфизмом представлений $J[f, X]$ и $J[g, X']$.*

Доказательство. Мы будем доказывать теорему индукцией по сложности $\Omega_2$-слова.

Если $w[f, X, m] = m$, то $m \in X$. Согласно условию 6.1.10.1,

$$R(m) = w[f \to g, X, X', R_1](w[f, X, m]) = w[f, X, R_1](m) = R_1(m)$$

Следовательно, на множестве $X$ отображения $R$ и $R_1$ совпадают, и отображение $R$ согласовано со структурой представления $f$.

Пусть $\omega \in \Omega_2(n)$. Пусть отображение $R$ определено для $m_1, ..., m_n \in J[f, X]$. Пусть

$$w_1 = w[f, X, m_1] \quad ... \quad w_n = w[f, X, m_n]$$

сли $m = m_1...m_n\omega$, то согласно правилу 6.1.6.2,

$$w[f, X, m] = w_1...w_n\omega$$

Согласно условию 6.1.10.2,

$$\begin{aligned}
R(m) &= w[f \to g, X, X', R_1](w[f, X, m]) = w[f \to g, X, X', R_1](w_1...w_n\omega) \\
&= w[f \to g, X, X', R_1](w_1)...w[f \to g, X, X', R_1](w_n)\omega \\
&= R(m_1)...R(m_n)\omega
\end{aligned}$$

Следовательно, отображение $R$ является эндоморфизмом $\Omega_2$-алгебры $A_2$.

Пусть отображение $R$ определено для $m_1 \in J[f, X]$, $w_1 = w[f, X, m_1]$. Пусть $a \in A_1$. Если $m = f(a)(m_1)$, то согласно правилу 6.1.6.3,

$$w[f, X, f(a)(m_1)] = f(a)(w_1)$$

Согласно условию 6.1.10.3,

$$\begin{aligned}
R(m) &= w[f \to g, X, X', R_1](w[f, X, m]) = w[f \to g, X, X', R_1](f(a)(w_1)) \\
&= f(a)(w[f \to g, X, X', R_1](w_1)) = f(a)(R(m_1))
\end{aligned}$$

Из равенства (3.2.3) следует, что отображение $R$ является морфизмом представления $f$.

Единственность эндоморфизма $R$, а следовательно, корректность его определения, следует из следующего рассуждения. Допустим, $m \in A_2$ имеет различные $\Omega_2$-слова относительно множества $X$, например

(6.1.13) $$m = x_1...x_n\omega = f(a)(x)$$

Так как $R$ - эндоморфизм представления, то из равенства (6.1.13) следует

(6.1.14) $$R(m) = R(x_1...x_n\omega) = R(x_1)...R(x_n)\omega = R(f(a)(x)) = f(a)(R(x))$$

Из равенства (6.1.14) следует

(6.1.15) $$R(m) = R(x_1)...R(x_n)\omega = f(a)(R(x))$$

Из равенств (6.1.13), (6.1.15) следует, что равенство (6.1.13) сохраняется при отображении. Следовательно, образ $A_2$ не зависит от выбора координат.    $\square$

Замечание 6.1.15. *Теорема 6.1.14 - это теорема о продолжении отображения. Единственное, что нам известно о множестве $X$ - это то, что $X$ - множество образующих представления $f$. Однако, между элементами множества $X$ могут существовать соотношения, порождённые либо операциями $\Omega_2$-алгебры $A_2$, либо преобразованиями представления $f$. Поэтому произвольное отображение множества $X$, вообще говоря, не может быть продолжено до приведенного морфизма представления $f$.*[6.5] *Однако, если отображение $R_1$ согласованно со структурой представления на множестве $X$, то мы можем построить продолжение этого отображения, которое является приведенным морфизмом представления $f$.*    $\square$

Определение 6.1.16. *Пусть $X$ - множество образующих представления*

$$f : A_1 \longrightarrow_* A_2$$

*$\Omega_1$-алгебры $A_1$ в $\Omega_2$-алгебре $A_2$. Пусть $Y$ - множество образующих представления*

$$g : A_1 \longrightarrow_* B_2$$

*$\Omega_1$-алгебры $A_1$ в $\Omega_2$-алгебре $B_2$. Пусть*

$$R : A_2 \to B_2$$

*приведенный морфизм представления $f$. Множество координат $W[g, Y, R(X)]$ называется **координатами приведенного морфизма представления**.* □

Из определений 6.1.6, 6.1.16 следует, что

$$W[g, Y, R(X)] = (W[g, Y, R(x)], x \in X)$$

Пусть $m \in A_2$. Если в слове $w[f, X, m]$ вместо каждого $x \in X$ подставить его образ $w[g, Y, R(x)]$, то, согласно теореме 6.1.14, мы получим $\Omega_2$-слово $w[g, Y, R(m)]$. Из этого утверждения следует определение 6.1.17.

Определение 6.1.17. *Пусть $X$ - множество образующих представления*

$$f : A_1 \longrightarrow_* A_2$$

*$\Omega_1$-алгебры $A_1$ в $\Omega_2$-алгебре $A_2$. Пусть $Y$ - множество образующих представления*

$$g : A_1 \longrightarrow_* B_2$$

*$\Omega_1$-алгебры $A_1$ в $\Omega_2$-алгебре $B_2$. Пусть $R$*

$$R : A_2 \to B_2$$

*приведенный морфизм представления $f$. Пусть $m \in A_2$. Мы определим **суперпозицию координат** приведенного морфизма $R$ представления $f$ и $A_2$-числа $m$ как координаты, определённые согласно правилу*

(6.1.16) $\qquad W[g, Y, R(X)] \circ W[f, X, m] = W[g, Y, R(m)]$

*Мы определим суперпозицию координат приведенного морфизма $R$ представления $f$ и множества $B \subseteq A_2$ согласно правилу*

(6.1.17) $\qquad W[g, Y, R(X)] \circ W[f, X, B] = (W[g, Y, R(X)] \circ W[f, X, m], m \in B)$

$$W[g, Y, R(X)] \circ w[f, X, B] = w[g, Y, R(X)] \circ W[f, X, B] = w[g, Y, R(B)]$$

□

Теорема 6.1.18. *Пусть $X$ - множество образующих представления*

$$f : A_1 \longrightarrow_* A_2$$

*$\Omega_1$-алгебры $A_1$ в $\Omega_2$-алгебре $A_2$. Пусть $Y$ - множество образующих представления*

$$g : A_1 \longrightarrow_* B_2$$

*$\Omega_1$-алгебры $A_1$ в $\Omega_2$-алгебре $B_2$. Приведенный морфизм представления*

$$R : A_2 \to B_2$$

*порождает отображение координат представления*

$$(6.1.18) \qquad W[f \to g, X, Y, R] : W[f, X] \to W[g, Y]$$

*такое, что*

$$(6.1.19) \qquad W[f, X, m] \to W[f \to g, X, Y, R] \circ W[f, X, m] = W[g, Y, R(m)]$$

Доказательство. Согласно замечанию 6.1.11, мы можем рассматривать равенства (6.1.16), (6.1.18) относительно заданных множеств образующих $X$, $Y$. При этом координатам $W[f, X, m]$ соответствует слово

$$(6.1.20) \qquad X \circ W[f, X, m] = w[f, X, m]$$

а координатам $W[g, Y, R(m)]$ соответствует слово

$$(6.1.21) \qquad Y \circ W[g, Y, R(m)] = w[g, Y, R(m)]$$

Поэтому для того, чтобы доказать теорему, нам достаточно показать, что отображению $W[f, X, R]$ соответствует отображение $w[f, X, R]$. Мы будем доказывать это утверждение индукцией по сложности $\Omega_2$-слова.

Если $m \in X$, $m' = R(m)$, то, согласно равенствам (6.1.20), (6.1.21), отображения $W[f, X, R]$ и $w[f, X, R]$ согласованы.

Пусть для $m_1, ..., m_n \in X$ отображения $W[f, X, R]$ и $w[f, X, R]$ согласованы. Пусть $\omega \in \Omega_2(n)$. Согласно теореме 6.1.12

$$(6.1.22) \qquad W[f, X, m_1...m_n\omega] = W[f, X, m_1]...W[f, X, m_n]\omega$$

Так как $R$ - эндоморфизм $\Omega_2$-алгебры $A_2$, то из равенства (6.1.22) следует

$$(6.1.23) \qquad \begin{aligned} W[f, X, R \circ (m_1...m_n\omega)] &= W[f, X, (R \circ m_1)...(R \circ m_n)\omega] \\ &= W[f, X, R \circ m_1]...W[f, X, R \circ m_n]\omega \end{aligned}$$

Из равенств (6.1.22), (6.1.23) и предположения индукции следует, что отображения $W[f, X, R]$ и $w[f, X, R]$ согласованы для $m = m_1...m_n\omega$.

Пусть для $m_1 \in A_2$ отображения $W[f, X, R]$ и $w[f, X, R]$ согласованы. Пусть $a \in A_1$. Согласно теореме 6.1.13

$$(6.1.24) \qquad W[f, X, f(a)(m_1)] = f(a)(W[f, X, m_1])$$

Так как $R$ - эндоморфизм представления $f$, то из равенства (6.1.24) следует

$$(6.1.25) \quad W[f, X, R \circ f(a)(m_1)] = W[f, X, f(a)(R \circ m_1)] = f(a)(W[f, X, R \circ m_1])$$

Из равенств (6.1.24), (6.1.25) и предположения индукции следует, что отображения $W[f, X, R]$ и $w[f, X, R]$ согласованы для $m = f(a)(m_1)$.    $\square$

Следствие 6.1.19. *Пусть $X$ - множество образующих представления $f$. Пусть $R$ - эндоморфизм представления $f$. Отображение $W[f, X, R]$ является эндоморфизмом представления $\Omega_1$-алгебры $A_1$ в $\Omega_2$-алгебре $W[f, X]$.*    $\square$

В дальнейшем мы будем отождествлять отображение $W[f, X, R]$ и множество координат $W[f, X, R \circ X]$.

Теорема 6.1.20. *Пусть $X$ - множество образующих представления $f$. Пусть $R$ - эндоморфизм представления $f$. Пусть $Y \subset A_2$. Тогда*

$$(6.1.26) \qquad W[f, X, R(X)] \circ W[f, X, Y] = W[f, X, R(Y)]$$

Доказательство. Равенство (6.1.26) является следствием равенства

$$R \circ Y = (R \circ m, m \in Y)$$

а также равенств (6.1.16), (6.1.17).          □

Теорема 6.1.21. *Пусть $X$ - множество образующих представления $f$. Пусть $R$, $S$ - эндоморфизмы представления $f$. Тогда*

$$(6.1.27) \qquad W[f, X, R] \circ W[f, X, S] = W[f, X, R \circ S]$$

Доказательство. Равенство (6.1.27) следует из равенства (6.1.26), если положить $Y = S \circ X$.          □

Концепция суперпозиции координат очень проста и напоминает своеобразную машину Тьюринга. Если элемент $m \in A_2$ имеет вид

$$m = m_1 ... m_n \omega$$

или

$$m = f(a)(m_1)$$

то мы ищем координаты элементов $m_i$, для того чтобы подставить их в соответствующее выражение. Как только элемент $m \in A_2$ принадлежит множеству образующих $\Omega_2$-алгебры $A_2$, мы выбираем координаты соответствующего элемента из второго множителя. Поэтому мы требуем, чтобы второй множитель в суперпозиции был множеством координат образа множества образующих $X$.

Следующие формы записи образа множества $Y$ при эндоморфизме $R$ эквивалентны.

$$(6.1.28) \qquad R \circ Y = (R(X)) \circ W[f, X, Y] = (X \circ W[f, X, R]) \circ W[f, X, Y]$$

Из равенств (6.1.26), (6.1.28) следует, что

$$(6.1.29) \qquad X \circ (W[f, X, R] \circ W[f, X, Y]) = (X \circ W[f, X, R]) \circ W[f, X, Y]$$

Равенство (6.1.29) является законом ассоциативности для операции композиции и позволяет записать выражение

$$X \circ W[f, X, R] \circ W[f, X, Y]$$

без использования скобок.

Определение 6.1.22. *Пусть $X \subset A_2$ - множество образующих представления*

$$f : A_1 \longrightarrow A_2$$

*Пусть отображение*

$$H : A_2 \rightarrow A_2$$

*является эндоморфизмом представления $f$. Пусть множество $X' = H \circ X$ является образом множества $X$ при отображении $H$. Эндоморфизм $H$ представления $f$ называется невырожденным на множестве образующих $X$, если множество $X'$ является множеством образующих представления $f$. В противном случае, эндоморфизм $H$ представления $f$ называется вырожденным на множестве образующих $X$.*          □

Определение 6.1.23. *Эндоморфизм представления $f$ называется **невырожденным**, если он невырожден на любом множестве образующих. В противном случае, эндоморфизм $H$ представления $f$ называется **вырожденным**.*          □

Теорема 6.1.24. *Автоморфизм $R$ представления*

$$f : A_1 \xrightarrow{\quad * \quad} A_2$$

*является невырожденным эндоморфизмом.*

Доказательство. Пусть $X$ - множество образующих представления $f$. Пусть $X' = R(X)$.

Согласно теореме 6.1.10 эндоморфизм $R$ порождает отображение $\Omega_2$-слов $w[f \to g, X, R]$.

Пусть $m' \in A_2$. Так как $R$ - автоморфизм, то существует $m \in A_2$, $R \circ m = m'$. Согласно определению 6.1.6, $w[f, X, m]$ - $\Omega_2$-слово, представляющее $A_2$ относительно множества образующих $X$. Согласно теореме 6.1.10, $w[f, X', m']$ - $\Omega_2$-слово, представляющее $m'$ относительно множества образующих $X'$

$$w[f, X', m'] = w[f \to g, X, R](w[f, X, m])$$

Следовательно, $X'$ - множество образующих представления $f$. Согласно определению 6.1.23, автоморфизм $R$ - невырожден.  □

## 6.2. Базис представления

Определение 6.2.1. *Пусть*

$$f : A_1 \xrightarrow{\quad * \quad} A_2$$

*представление $\Omega_1$-алгебры $A_1$ в $\Omega_2$-алгебре $A_2$ и*

$$Gen[f] = \{X \subseteq A_2 : J[f, X] = A_2\}$$

*Если для множества $X \subset A_2$ верно $X \in Gen[f]$, то для любого множества $Y$, $X \subset Y \subset A_2$, также верно $Y \in Gen[f]$. Если существует минимальное множество $X \in Gen[f]$, то множество $X$ называется **квазибазисом представления** $f$.*  □

Теорема 6.2.2. *Если множество $X$ является квазибазисом представления $f$, то для любого $m \in X$ множество $X \setminus \{m\}$ не является множеством образующих представления $f$.*

Доказательство. Пусть $X$ - квазибазис представления $f$. Допустим для некоторого $m \in X$ существует $\Omega_2$-слово

$$w = w[f, X \setminus \{m\}, m]$$

Рассмотрим $A_2$-число $m'$, для которого $\Omega_2$-слово $w' = w[f, X, m']$ зависит от $m$. Согласно определению 6.1.6, любое вхождение $A_2$-числа $m$ в $\Omega_2$-слово $w'$ может быть заменено $\Omega_2$-словом $w$. Следовательно, $\Omega_2$-слово $w'$ не зависит от $m$, а множество $X \setminus \{m\}$ является множеством образующих представления $f$. Следовательно, $X$ не является квазибазисом расслоения $f$.  □

Замечание 6.2.3. *Доказательство теоремы 6.2.2 даёт нам эффективный метод построения квазибазиса представления $f$. Выбрав произвольное множество образующих, мы шаг за шагом исключаем те элементы множества, которые имеют координаты относительно остальных элементов множества. Если множество образующих представления бесконечно, то рассмотренная операция может не иметь последнего шага. Если представление имеет конечное множество образующих, то за конечное число шагов мы можем построить квазибазис этого представления.*  □

Мы ввели $\Omega_2$-слово элемента $x \in A_2$ относительно множества образующих $X$ в определении 6.1.6. Из теоремы 6.2.2 следует, что если множество образующих $X$ не является квазибазисом, то выбор $\Omega_2$-слова относительно множества образующих $X$ неоднозначен. Но даже если множество образующих $X$ является квазибазисом, то представление $m \in A_2$ в виде $\Omega_2$-слова неоднозначно.

Замечание 6.2.4. *Существует три источника неоднозначности в записи $\Omega_2$-слова.*

6.2.4.1: *В $\Omega_i$-алгебре $A_i$, $i = 1, 2$, могут быть определены равенства. Например, если $e$ - единица мультипликативной группы $A_i$, то верно равенство*

$$ae = a$$

*для любого $a \in A_i$.*

6.2.4.2: *Неоднозначность выбора $\Omega_2$-слова может быть связана со свойствами представления. Например, если $m_1, ..., m_n$ - $\Omega_2$-слова, $\omega \in \Omega_2(n)$ и $a \in A_1$, то[6.6]*

$$(6.2.1) \qquad f(a)(m_1...m_n\omega) = (f(a)(m_1))...(f(a)(m_n))\omega$$

*В тоже время, если $\omega$ является операцией $\Omega_1$-алгебры $A_1$ и операцией $\Omega_2$-алгебры $A_2$, то мы можем потребовать, что $\Omega_2$-слова $f(a_1...a_n\omega)(x)$ и $(f(a_1)(x))...(f(a_n)(x))\omega$ описывают один и тот же элемент $\Omega_2$-алгебры $A_2$.[6.7]*

$$(6.2.2) \qquad f(a_1...a_n\omega)(x) = (f(a_1)(x))...(f(a_n)(x))\omega$$

6.2.4.3: *Равенства вида (6.2.1), (6.2.2) сохраняются при морфизме представления. Поэтому мы можем игнорировать эту форму неоднозначности записи $\Omega_2$-слова. Однако возможна принципиально другая форма неоднозначности, пример которой можно найти в теоремах 9.3.15, 9.3.16.*

*Таким образом, мы видим, что на множестве $\Omega_2$-слов можно определить различные отношения эквивалентности.[6.8] Наша задача - найти максимальное отношение эквивалентности на множестве $\Omega_2$-слов, которое сохраняется при морфизме представления.*

*Аналогичное замечание касается отображения $W[f, X, m]$, определённого в замечании 6.1.11.[6.9]* □

---

[6.6] Например, пусть $\{e_1, e_2\}$ - базис векторного пространства над полем $k$. Равенство (6.2.1) принимает форму закона дистрибутивности

$$a(b^1 e_1 + b^2 e_2) = (ab^1)e_1 + (ab^2)e_2$$

[6.7] Для векторного пространства это требование принимает форму закона дистрибутивности

$$(a + b)e_1 = ae_1 + be_1$$

[6.8] Очевидно, что каждое из равенств (6.2.1), (6.2.2) порождает некоторое отношение эквивалентности.

[6.9] Если базис векторного пространства - конечен, то мы можем представить базис в виде матрицы строки

$$\overline{\overline{e}} = \begin{pmatrix} e_1 & ... & e_2 \end{pmatrix}$$

Теорема 6.2.5. *Пусть $X$ - квазибазис представления*

$$f : A_1 \overset{*}{\longrightarrow} A_2$$

*Рассмотрим отношение эквивалентности*

$$\lambda[f, X] \subseteq w[f, X] \times w[f, X]$$

*которое порождено исключительно следующими утверждениями.*

6.2.5.1: *Если в $\Omega_2$-алгебре $A_2$ существует равенство*

$$w_1[f, X, m] = w_2[f, X, m]$$

*определяющее структуру $\Omega_2$-алгебры, то*

$$(w_1[f, X, m], w_2[f, X, m]) \in \lambda[f, X]$$

6.2.5.2: *Если в $\Omega_1$-алгебре $A_1$ существует равенство*

$$w_1[f, X, m] = w_2[f, X, m]$$

*определяющее структуру $\Omega_1$-алгебры, то*

$$(f(w_1)(w[f, X, m]), f(w_2)(w[f, X, m])) \in \lambda[f, X]$$

6.2.5.3: *Для любой операции $\omega \in \Omega_1(n)$,*

$$(f(a_{11}...a_{1n}\omega)(a_2), (f(a_{11})...f(a_{1n})\omega)(a_2)) \in \lambda[f, X]$$

6.2.5.4: *Для любой операции $\omega \in \Omega_2(n)$,*

$$(f(a_1)(a_{21}...a_{2n}\omega), f(a_1)(a_{21})...f(a_1)(a_{2n})\omega) \in \lambda[f, X]$$

6.2.5.5: *Пусть $\omega \in \Omega_1(n) \cap \Omega_2(n)$. Если представление $f$ удовлетворяет равенству*[6.10]

$$f(a_{11}...a_{1n}\omega)(a_2) = (f(a_{11})(a_2))...(f(a_{1n})(a_2))\omega$$

*то мы можем предположить, что верно равенство*

$$(f(a_{11}...a_{1n}\omega)(a_2), (f(a_{11})(a_2))...(f(a_{1n})(a_2))\omega) \in \lambda[f, X]$$

Мы можем представить отображение $W[f, \overline{\overline{e}}](v)$ в виде матрицы столбца

$$W[f, \overline{\overline{e}}, v] = \begin{pmatrix} v^1 \\ ... \\ v^n \end{pmatrix}$$

Тогда

$$W[f, \overline{\overline{e}}, v](\overline{\overline{e}}') = W[f, \overline{\overline{e}}, v] \begin{pmatrix} e_1' & ... & e_n' \end{pmatrix} = \begin{pmatrix} v^1 \\ ... \\ v^n \end{pmatrix} \begin{pmatrix} e_1' & ... & e_n' \end{pmatrix}$$

имеет вид произведения матриц.

[6.10] Рассмотрим представление коммутативного кольца $D$ в $D$-алгебре $A$. Мы будем пользоваться записью

$$f(a)(v) = av$$

В обеих алгебрах определены операции сложения и умножения. Однако равенство

$$f(a + b)(v) = f(a)(v) + f(b)(v)$$

верно, а равенство

$$f(ab)(v) = f(a)(v)f(b)(v)$$

является ошибочным.

Доказательство. Теорема верна, так как рассмотренные равенства сохраняются при гомоморфизмах универсальных алгебр $A_1$ и $A_2$. □

Определение 6.2.6. *Квазибазис $\overline{\overline{e}}$ представления $f$ такой, что*

$$\rho[f, \overline{\overline{e}}] = \lambda[f, \overline{\overline{e}}]$$

*называется* **базисом представления** *$f$.* □

Замечание 6.2.7. *Как отметил Кон в [14], стр. 96, 97, представление может иметь неэквивалентные базисы. Например, циклическая группа шестого порядка имеет базисы $\{a\}$ и $\{a^2, a^3\}$, которые нельзя отобразить один в другой эндоморфизмом представления.* □

Замечание 6.2.8. *Мы будем записывать базис также в виде*

$$\overline{\overline{e}} = (e, e \in \overline{\overline{e}})$$

*Если базис - конечный, то мы будем также пользоваться записью*

$$\overline{\overline{e}} = (e_i, i \in I) = (e_1, ..., e_n)$$

□

Теорема 6.2.9. *Автоморфизм представления $f$ отображает базис представления $f$ в базис.*

Доказательство. Пусть отображение $R$ - автоморфизм представления $f$. Пусть множество $\overline{\overline{e}}$ - базис представления $f$. Пусть [6.11] $\overline{\overline{e}}' = R \circ \overline{\overline{e}}$. Допустим множество $\overline{\overline{e}}'$ не является базисом. Согласно теореме 6.2.2 существует $e' \in \overline{\overline{e}}'$ такое, что $\overline{\overline{e}}' \backslash \{e'\}$ является множеством образующих представления $f$. Согласно теореме 3.5.5 отображение $R^{-1}$ является автоморфизмом представления $f$. Согласно теореме 6.1.24 и определению 6.1.23, множество $\overline{\overline{e}} \backslash \{e\}$ является множеством образующих представления $f$. Полученное противоречие доказывает теорему. □

Теорема 6.2.10. *Пусть $\overline{\overline{e}}$ - базис представления $f$. Пусть*

$$R_1 : \overline{\overline{e}} \to \overline{\overline{e}}'$$

*произвольное отображение множества $X$. Рассмотрим отображение $\Omega_2$-слов*

$$w[f \to g, \overline{\overline{e}}, \overline{\overline{e}}', R_1] : w[f, \overline{\overline{e}}] \to w[g, \overline{\overline{e}}']$$

*удовлетворяющее условиям 6.1.10.1, 6.1.10.2, 6.1.10.3, и такое, что*

$$e \in \overline{\overline{e}} => w[f \to g, \overline{\overline{e}}, \overline{\overline{e}}', R_1](e) = R_1(e)$$

*Существует единственный эндоморфизм представления $f$ [6.12]*

$$r_2 : A_2 \to A_2$$

*определённый правилом*

$$R(m) = w[f \to g, \overline{\overline{e}}, \overline{\overline{e}}', R_1](w[f, \overline{\overline{e}}, m])$$

Доказательство. Утверждение теоремы является следствием теорем 6.1.10, 6.1.14. □

---

[6.11] Согласно определениям 5.1.3, 6.4.1, мы будем пользоваться записью $R(\overline{\overline{e}}) = R \circ \overline{\overline{e}}$.
[6.12] Это утверждение похоже на теорему [2]-1, с. 104.

Следствие 6.2.11. *Пусть $\overline{\overline{e}}$, $\overline{\overline{e}}'$ - базисы представления $f$. Пусть $R$ - автоморфизм представления $f$ такой, что $\overline{\overline{e}}' = R \circ \overline{\overline{e}}$. Автоморфизм $R$ определён однозначно.*

Замечание 6.2.12. *Теорема 6.2.10, так же как и теорема 6.1.14, является теоремой о продолжении отображения. Однако здесь $\overline{\overline{e}}$ - не произвольное множество образующих представления, а базис. Согласно замечанию 6.2.3, мы не можем определить координаты любого элемента базиса через остальные элементы этого же базиса. Поэтому отпадает необходимость в согласованности отображения базиса с представлением.* $\square$

Теорема 6.2.13. *Набор координат $W[f, \overline{\overline{e}}, \overline{\overline{e}}]$ соответствует тождественному преобразованию*

$$W[f, \overline{\overline{e}}, E] = W[f, \overline{\overline{e}}, \overline{\overline{e}}]$$

Доказательство. Утверждение теоремы следует из равенства

$$m = \overline{\overline{e}} \circ W[f, \overline{\overline{e}}, m] = \overline{\overline{e}} \circ W[f, \overline{\overline{e}}, \overline{\overline{e}}] \circ W[f, \overline{\overline{e}}, m]$$

$\square$

Теорема 6.2.14. *Пусть $W[f, \overline{\overline{e}}, R \circ \overline{\overline{e}}]$ - множество координат автоморфизма $R$. Определено множество координат $W[f, R \circ \overline{\overline{e}}, \overline{\overline{e}}]$, соответствующее автоморфизму $R^{-1}$. Множество координат $W[f, R \circ \overline{\overline{e}}, \overline{\overline{e}}]$ удовлетворяют равенству*

$$(6.2.3) \qquad W[f, \overline{\overline{e}}, R \circ \overline{\overline{e}}] \circ W[f, R \circ \overline{\overline{e}}, \overline{\overline{e}}] = W[f, \overline{\overline{e}}, \overline{\overline{e}}]$$

$$W[f \to f, \overline{\overline{e}}, \overline{\overline{e}}, R^{-1}] = W[f \to f, \overline{\overline{e}}, \overline{\overline{e}}, R]^{-1} = W[f, R \circ \overline{\overline{e}}, \overline{\overline{e}}]$$

Доказательство. Поскольку $R$ - автоморфизм представления $f$, то, согласно теореме 6.2.9, множество $R \circ \overline{\overline{e}}$ - базис представления $f$. Следовательно, определено множество координат $W[f, R \circ \overline{\overline{e}}, \overline{\overline{e}}]$. Равенство (6.2.3) следует из цепочки равенств

$$W[f, \overline{\overline{e}}, R \circ \overline{\overline{e}}] \circ W[f, R \circ \overline{\overline{e}}, \overline{\overline{e}}] = W[f, \overline{\overline{e}}, R \circ \overline{\overline{e}}] \circ W[f, \overline{\overline{e}}, R^{-1} \circ \overline{\overline{e}}]$$
$$= W[f, \overline{\overline{e}}, R \circ R^{-1} \circ \overline{\overline{e}}] = W[f, \overline{\overline{e}}, \overline{\overline{e}}]$$

$\square$

Замечание 6.2.15. *В $\Omega_2$-алгебре $A_2$ не существует универсального алгоритма определения множества координат $W[f, R \circ \overline{\overline{e}}, \overline{\overline{e}}]$ для заданного множества $W[f, \overline{\overline{e}}, R \circ \overline{\overline{e}}]$.* [6.13] *Мы полагаем, что в теореме 6.2.14 этот алгоритм задан неявно. Очевидно также, что множество $\Omega_2$-слов*

$$(6.2.4) \qquad \overline{\overline{e}} \circ W[f, R \circ \overline{\overline{e}}, \overline{\overline{e}}] \circ W[f, \overline{\overline{e}}, R \circ \overline{\overline{e}}]$$

*вообще говоря, не совпадает с множеством $\Omega_2$-слов*

$$(6.2.5) \qquad \overline{\overline{e}} \circ W[f, \overline{\overline{e}}, \overline{\overline{e}}]$$

*Теорема 6.2.14 утверждает, что множества $\Omega_2$-слов (6.2.4) и (6.2.5) совпадают с точностью до отношения эквивалентности, порождённой представлением $f$.* $\square$

---

[6.13]В векторном пространстве линейному преобразованию соответствует матрица чисел. Соответственно, обратному преобразованию соответствует обратная матрица.

Теорема 6.2.16. *Пусть* $W[f, \overline{e}, R \circ \overline{e}]$ - *множество координат автоморфизма* $R$. *Пусть* $W[f, \overline{e}, S \circ \overline{e}]$- *множество координат автоморфизма* $S$. *Множество координат автоморфизма* $(R \circ S)^{-1}$ *удовлетворяет равенству*

(6.2.6)  $W[f, (R \circ S) \circ \overline{e}, \overline{e}] = W[f, S \circ (R \circ \overline{e}), \overline{e}] = W[f, S \circ \overline{e}, \overline{e}] \circ W[f, R \circ \overline{e}, \overline{e}]$

Доказательство. Равенство

$$W[f, (R \circ S) \circ \overline{e}, \overline{e}] = W[f, \overline{e}, (R \circ S)^{-1} \circ \overline{e}] = W[f, \overline{e}, S^{-1} \circ R^{-1} \circ \overline{e}]$$

(6.2.7)
$$= W[f, \overline{e}, S^{-1} \circ e] \circ W[f, \overline{e}, R^{-1} \circ \overline{e}]$$
$$= W[f, S \circ \overline{e}, \overline{e}] \circ W[f, R \circ \overline{e}, \overline{e}]$$
$$= W[f, S \circ (R \circ \overline{e}), \overline{e}]$$

является следствием теорем 6.1.21, 6.2.14. Равенство (6.2.6) является следствием равенства (6.2.7). □

Теорема 6.2.17. *Группа автоморфизмов* $GA(f)$ *эффективного представления* $f$ *в* $\Omega_2$-*алгебре* $A_2$ *порождает эффективное левостороннее представление в* $\Omega_2$-*алгебре* $A_2$.

Доказательство. Из следствия 6.2.11 следует, что если автоморфизм $R$ отображает базис $\overline{e}$ в базис $\overline{e}'$, то множество координат $W[f, \overline{e}, \overline{e}']$ однозначно определяет автоморфизм $R$. Из теоремы 6.1.18 следует, что множество координат $W[f, \overline{e}, \overline{e}']$ определяет правило отображения координат относительно базиса $\overline{e}$ при автоморфизме представления $f$. Из равенства (6.1.28) следует, что автоморфизм $R$ действует слева на элементы $\Omega_2$-алгебры $A_2$. Из равенства (6.1.27) следует, что представление группы является левосторонним представлением. Согласно теореме 6.2.13 набор координат $W[f, \overline{e}, \overline{e}]$ соответствует тождественному преобразованию. Из теоремы 6.2.14 следует, что набор координат $W[f, R \circ \overline{e}, \overline{e}]$ соответствует преобразованию, обратному преобразованию $W[f, \overline{e}, R \circ \overline{e}]$. □

## 6.3. Свободное представление

В разделе 3.1 мы рассмотрели определение 3.1.4 свободного представления. Однако мы можем рассмотреть другое определение, которое аналогично определению свободного модуля.

Определение 6.3.1. *Представление*

$$f : A_1 \longrightarrow_* A_2$$

*называется свободным, если это представление имеет базис.* □

Теорема 6.3.2. *Пусть*

$$f : A_1 \longrightarrow_* A_2$$

*свободное представление согласно определению* 6.3.1. *Тогда представление* $f$ *свободно согласно определению* 3.1.4.

Доказательство. Пусть $\overline{e}$ - базис представления $f$ и $m \in \overline{e}$. Пусть существуют $A_1$-числа $a_1$, $b_1$ такие, что $f(a_1) = f(b_1)$. Согласно предположению $f(a_1)(m) = f(b_1)(m)$. Однако, если $a_1 \neq b_1$, то $f(a_1)(m)$ и $f(b_1)(m)$ различные $\Omega_2$-слова. Следовательно, $\overline{\overline{e}}$ не является базисом. Из полученного противоречия следует, что $a_1 = b_1$. Следовательно представление $f$ свободно согласно определению 3.1.4.                                                              □

Теорема 6.3.3. *Пусть*

$$f : A_1 \xrightarrow{\quad*\quad} A_2$$

*свободное представление согласно определению 3.1.4. Тогда представление $f$ свободно согласно определению 6.3.1.*

Вопрос 6.3.4. *Очень важно найти доказательство теоремы 6.3.3 или найти пример, когда эта теорема не верна. Мы увидим в главе 7 какую роль играет свободное представление согласно определению 6.3.1. Так как в дальнейшем я буду предполагать, что представление всегда имеет базис, то в рамках этой книги я могу ограничиться теоремой 6.3.2.*                    □

## 6.4. Многообразие базисов представления

Множество $\mathcal{B}[f]$ базисов представления $f$ называется **многообразием базисов** представления $f$.

Определение 6.4.1. *Согласно теоремам 6.1.20, 6.2.9, автоморфизм $R$ представления $f$ порождает преобразование*

$$(6.4.1) \qquad \begin{aligned} R &: \overline{\overline{h}} \to R \circ \overline{\overline{h}} \\ R \circ \overline{\overline{h}} &= W[f, \overline{\overline{e}}, R \circ \overline{\overline{e}}] \circ \overline{\overline{h}} \end{aligned}$$

*многообразия базисов представления. Это преобразование называется **активным**. Согласно теореме 3.5.5, определено левостороннее представление*

$$A(f) : GA(f) \xrightarrow{\quad*\quad} \mathcal{B}[f]$$

*группы $GA(f)$ в многообразии базисов $\mathcal{B}[f]$. Представление $A(f)$ называется **активным представлением**. Согласно следствию 6.2.11, это представление однотранзитивно.*                                                              □

Замечание 6.4.2. *Согласно замечанию 6.2.3, могут существовать базисы представления $f$, не связанные активным преобразованием. В этом случае мы в качестве многообразия базисов будем рассматривать орбиту выбранного базиса. Следовательно, представление $f$ может иметь различные многообразия базисов. Мы будем предполагать, что мы выбрали многообразие базисов.*

Теорема 6.4.3. *Существует однотранзитивное правостороннее представление*

$$P(f) : GA(f) \xrightarrow{\quad*\quad} \mathcal{B}[f]$$

*группы $GA(f)$ в многообразии базисов $\mathcal{B}[f]$. Представление $P(f)$ называется **пассивным представлением**.*

Доказательство. Поскольку $A(f)$ - однотранзитивное левостороннее представление группы $GA(f)$, то однотранзитивное правостороннее представление $P(f)$ определено однозначно согласно теореме 5.5.9.                    □

Теорема 6.4.4. *Преобразование представления $P(f)$ называется* **пассивным преобразованием многообразия базисов** *представления. Мы будем пользоваться записью*

$$S(\overline{\overline{e}}) = \overline{\overline{e}} \circ S$$

*для обозначения образа базиса $\overline{\overline{e}}$ при пассивном преобразовании $S$. Пассивное преобразование базиса имеет вид*

(6.4.2)
$$S : \overline{\overline{h}} \rightarrow \overline{\overline{h}} \circ S$$
$$\overline{\overline{h}} \circ S = \overline{\overline{h}} \circ W[f, \overline{\overline{e}}, \overline{\overline{e}} \circ S]$$

Доказательство. Согласно равенству (6.4.1), активное преобразование действует на координаты базиса слева. Равенство (8.3.2) следует из теорем 5.5.8, 5.5.9, 5.5.11, согласно которым пассивное преобразование действует на координаты базиса справа. $\square$

Теорема 6.4.5. *Пассивное преобразование многообразия базисов является автоморфизмом представления $A(f)$.*

Доказательство. Теорема является следствием теоремы 5.5.11. $\square$

Теорема 6.4.6. *Пусть $S$ - пассивное преобразование многообразия базисов представления $f$. Пусть $\overline{\overline{e}}_1$ - базис представления $f$, $\overline{\overline{e}}_2 = \overline{\overline{e}}_1 \circ S$. Пусть для базиса $\overline{\overline{e}}_3$ существует активное преобразование $R$ такое, что $\overline{\overline{e}}_3 = R \circ \overline{\overline{e}}_1$. Положим $\overline{\overline{e}}_4 = R \circ \overline{\overline{e}}_2$. Тогда $\overline{\overline{e}}_4 = \overline{\overline{e}}_3 \circ S$.*

Доказательство. Согласно равенству (6.4.1), активное преобразование координат базиса $\overline{\overline{e}}_3$ имеет вид

(6.4.3)
$$\overline{\overline{e}}_4 = W[f, \overline{\overline{e}}_1, \overline{\overline{e}}_3] \circ \overline{\overline{e}}_2 = W[f, \overline{\overline{e}}_1, \overline{\overline{e}}_3] \circ \overline{\overline{e}}_1 \circ W[f, \overline{\overline{e}}_1, \overline{\overline{e}}_2]$$

Пусть $\overline{\overline{e}}_5 = \overline{\overline{e}}_3 \circ S$. Из равенства (6.4.2) следует, что

(6.4.4)
$$\overline{\overline{e}}_5 = \overline{\overline{e}}_3 \circ W[f, \overline{\overline{e}}_1, \overline{\overline{e}}_2] = W[f, \overline{\overline{e}}_1, \overline{\overline{e}}_3] \circ \overline{\overline{e}}_1 \circ W[f, \overline{\overline{e}}_1, \overline{\overline{e}}_2]$$

Из совпадения выражений в равенствах (6.4.3), (6.4.4) следует, что $\overline{\overline{e}}_4 = \overline{\overline{e}}_5$. Следовательно, коммутативна диаграмма

$$
\begin{array}{ccc}
\overline{\overline{e}}_1 \in \mathcal{B}[f] & \xrightarrow{\ \ R\ \ } & \overline{\overline{e}}_3 \in \mathcal{B}[f] \\
\Big\downarrow{\scriptstyle S} & & \Big\downarrow{\scriptstyle S} \\
\overline{\overline{e}}_2 \in \mathcal{B}[f] & \xrightarrow{\ \ R\ \ } & \overline{\overline{e}}_4 \in \mathcal{B}[f]
\end{array}
$$

$\square$

## 6.5. Геометрический объект представления универсальной алгебры

Активное преобразование изменяет базис представления и $\Omega_2$-число согласовано и координаты $\Omega_2$-числа относительно базиса не меняются. Пассивное преобразование меняет только базис, и это ведёт к изменению координат $\Omega_2$-числа относительно базиса.

Теорема 6.5.1. *Допустим пассивное преобразование $S \in GA(f)$ отображает базис $\overline{\overline{e}}_1 \in \mathcal{B}[f]$ в базис $\overline{\overline{e}}_2 \in \mathcal{B}[f]$*

(6.5.1)
$$\overline{\overline{e}}_2 = \overline{\overline{e}}_1 \circ S = \overline{\overline{e}}_1 \circ W[f, \overline{\overline{e}}_1, \overline{\overline{e}}_1 \circ S]$$

*Допустим $A_2$-число $m$ имеет $\Omega_2$-слово*

$$(6.5.2) \qquad m = \overline{\overline{e}}_1 \circ W[f, \overline{\overline{e}}_1, m]$$

*относительно базиса $\overline{\overline{e}}_1$ и имеет $\Omega_2$-слово*

$$(6.5.3) \qquad m = \overline{\overline{e}}_2 \circ W[f, \overline{\overline{e}}_2, m]$$

*относительно базиса $\overline{\overline{e}}_2$. Преобразование координат*

$$(6.5.4) \qquad W[f, \overline{\overline{e}}_2, m] = W[f, \overline{\overline{e}}_1 \circ S, \overline{\overline{e}}_1] \circ W[f, \overline{\overline{e}}_1, m]$$

*не зависит от $A_2$-числа $m$ или базиса $\overline{\overline{e}}_1$, а определенно исключительно координатами $A_2$-числа $m$ относительно базиса $\overline{\overline{e}}_1$.*

Доказательство. Из (6.5.1) и (6.5.3) следует, что

$$(6.5.5) \qquad \begin{aligned} \overline{\overline{e}}_1 \circ W[f, \overline{\overline{e}}_1, m] &= \overline{\overline{e}}_2 \circ W[f, \overline{\overline{e}}_2, m] = \overline{\overline{e}}_1 \circ W[f, \overline{\overline{e}}_1, \overline{\overline{e}}_2] \circ W[f, \overline{\overline{e}}_2, m] \\ &= \overline{\overline{e}}_1 \circ W[f, \overline{\overline{e}}_1, \overline{\overline{e}}_1 \circ S] \circ W[f, \overline{\overline{e}}_2, m] \end{aligned}$$

Сравнивая (6.5.2) и (6.5.5) получаем, что

$$(6.5.6) \qquad W[f, \overline{\overline{e}}_1, m] = W[f, \overline{\overline{e}}_1, \overline{\overline{e}}_1 \circ S] \circ W[f, \overline{\overline{e}}_2, m]$$

Так как $S$ - автоморфизм представления, то равенство (6.5.4) следует из (6.5.6) и теоремы 6.2.14. $\qquad \square$

Теорема 6.5.2. *Преобразования координат (6.5.4) порождают эффективное контравариантное правостороннее представление группы $GA(f)$, называемое* **координатным представлением** *в $\Omega_2$-алгебре.*

Доказательство. Согласно следствию 6.1.19, преобразование (6.5.4) является эндоморфизмом представления [6.14]

$$f : A_1 \longrightarrow\!\!\!* W[f, \overline{\overline{e}}_1]$$

Допустим мы имеем два последовательных пассивных преобразования $S$ и $T$. Преобразование координат

$$(6.5.7) \qquad W[f, \overline{\overline{e}}_2, m] = W[f, \overline{\overline{e}}_1 \circ S, \overline{\overline{e}}_1] \circ W[f, \overline{\overline{e}}_1, m]$$

соответствует пассивному преобразованию $S$. Преобразование координат

$$(6.5.8) \qquad W[f, \overline{\overline{e}}_2, m] = W[f, \overline{\overline{e}}_1 \circ T, \overline{\overline{e}}_1] \circ W[f, \overline{\overline{e}}_1, m]$$

соответствует пассивному преобразованию $T$. Согласно теореме 6.4.3, произведение преобразований координат (6.5.7) и (6.5.8) имеет вид

$$\begin{aligned} W[f, \overline{\overline{e}}_3, m] &= W[f, \overline{\overline{e}}_1 \circ T, \overline{\overline{e}}_1] \circ W[f, \overline{\overline{e}}_1 \circ S, \overline{\overline{e}}_1] \circ W[f, \overline{\overline{e}}_1, m] \\ &= W[f, \overline{\overline{e}}_1 \circ T \circ S, \overline{\overline{e}}_1] \circ W[f, \overline{\overline{e}}_1, m] \end{aligned}$$

и является координатным преобразованием, соответствующим пассивному преобразованию $S \circ T$. Согласно теоремам 6.2.14, 6.2.16 и определению 5.1.11 преобразования координат порождают правостороннее контравариантное представление группы $GA(f)$.

---

[6.14] Это преобразование не порождает эндоморфизма представления $f$. Координаты меняются, поскольку меняется базис, относительно которого мы определяем координаты. Однако $A_2$-число, координаты которого мы рассматриваем, не меняется.

Если координатное преобразование не изменяет координаты выбранного базиса, то ему соответствует единица группы $GA(f)$, так как пассивное представление однотранзитивно. Следовательно, координатное представление эффективно.                                                                                                                □

Пусть $f$ - представление $\Omega_1$-алгебры $A_1$ в $\Omega_2$-алгебре $A_2$. Пусть $g$ - представление $\Omega_1$-алгебры $A_1$ в $\Omega_3$-алгебре $A_3$. Пассивное представление $P(g)$ согласовано с пассивным представлением $P(f)$, если существует гомоморфизм $h$ группы $GA(f)$ в группу $GA(g)$. Рассмотрим диаграмму

$$
\begin{array}{ccc}
\mathrm{End}(B[f]) & \xrightarrow{\;H\;} & \mathrm{End}(B[g]) \\[2pt]
\big\uparrow{\scriptstyle P(f)} & {\scriptstyle f'}\nearrow & \big\uparrow{\scriptstyle P(g)} \\[2pt]
GA(f) & \xrightarrow[\;h\;]{} & GA(g)
\end{array}
$$

Так как отображения $P(f)$, $P(g)$ являются изоморфизмами группы, то отображение $H$ является гомоморфизмом групп. Следовательно, отображение $f'$ является представлением группы $GA(f)$ в многообразии базисов $\mathcal{B}(g)$. Согласно построению, пассивному преобразованию $S$ многообразии базисов $\mathcal{B}(f)$ соответствует пассивное преобразование $H(S)$ многообразия базисов $\mathcal{B}(g)$

$$(6.5.9) \qquad\qquad \overline{\overline{e}}_{g1} = \overline{\overline{e}}_g \circ H(S)$$

Тогда координатное преобразование в представлении $g$ принимает вид

$$(6.5.10) \qquad W[g, \overline{\overline{e}}_{g1}, m] = W[g, \overline{\overline{e}}_g \circ H(S), \overline{\overline{e}}_g] \circ W[g, \overline{\overline{e}}_g, m]$$

Определение 6.5.3. *Мы будем называть орбиту*

$$
\mathcal{O}(f, g, \overline{\overline{e}}_g, m) = H(GA(f)) \circ W[g, \overline{\overline{e}}_g, m]
$$
$$
= (W[g, \overline{\overline{e}}_g \circ H(S), \overline{\overline{e}}_g] \circ W[g, \overline{\overline{e}}_g, m], \overline{\overline{e}}_f \circ S, S \in GA(f))
$$

**геометрическим объектом в координатном представлении**, *определённом в представлении* $f$. *Для любого базиса* $\overline{\overline{e}}_{f1} = \overline{\overline{e}}_f \circ S$  *соответствующая точка* (6.5.10) *орбиты определяет* **координаты геометрического объекта** *относительно базиса* $\overline{\overline{e}}_{f1}$.                      □

Определение 6.5.4. *Мы будем называть орбиту*

$$
\mathcal{O}(f, g, m) = (W[g, \overline{\overline{e}}_g \circ H(S), \overline{\overline{e}}_g] \circ W[g, \overline{\overline{e}}_g, m], \overline{\overline{e}}_g \circ H(S), \overline{\overline{e}}_f \circ S, S \in GA(f))
$$

**геометрическим объектом**, *определённым в представлении* $f$. *Мы будем также говорить, что* $m$ - *это* **геометрический объект типа** $H$. *Для любого базиса* $\overline{\overline{e}}_{f1} = \overline{\overline{e}}_f \circ S$  *соответствующая точка* (6.5.10) *орбиты определяет* $A_2$-*число*

$$
m = \overline{\overline{e}}_g \circ W[g, \overline{\overline{e}}_g, m]
$$

*которое мы называем* **представителем геометрического объекта** *в представлении* $f$.                                                                                                □

Так как геометрический объект - это орбита представления, то согласно теореме 5.3.7 определение геометрического объекта корректно.

Определение 6.5.3 строит геометрический объект в координатном пространстве. Определение 6.5.4 предполагает, что мы выбрали базис представления $g$. Это позволяет использовать представитель геометрического объекта вместо его координат.

Теорема 6.5.5 (принцип инвариантности). *Представитель геометрического объекта не зависит от выбора базиса* $\overline{\overline{e}}_f$.

Доказательство. Чтобы определить представителя геометрического объекта, мы должны выбрать базис $\overline{\overline{e}}_f$ представления $f$, базис $\overline{\overline{e}}_g$ представления $g$ и координаты геометрического объекта $W[g, \overline{\overline{e}}_g, n]$. Соответствующий представитель геометрического объекта имеет вид

$$n = \overline{\overline{e}}_g \circ W[g, \overline{\overline{e}}_g, n]$$

Базис $\overline{\overline{e}}_{f1}$ связан с базисом $\overline{\overline{e}}_f$ пассивным преобразованием

$$\overline{\overline{e}}_{f1} = \overline{\overline{e}}_f \circ S$$

Согласно построению это порождает пассивное преобразование (6.5.9) и координатное преобразование (6.5.10). Соответствующий представитель геометрического объекта имеет вид

$$\begin{aligned}
n' &= \overline{\overline{e}}_{g1} \circ W[g, \overline{\overline{e}}_{g1}, n'] \\
&= \overline{\overline{e}}_g \circ W[g, \overline{\overline{e}}_g \circ H(S)] \circ W[g, \overline{\overline{e}}_g \circ H(S), \overline{\overline{e}}_g] \circ W[g, \overline{\overline{e}}_g, n] \\
&= \overline{\overline{e}}_g \circ W[g, \overline{\overline{e}}_g, n] = n
\end{aligned}$$

Следовательно, представитель геометрического объекта инвариантен относительно выбора базиса.    □

Теорема 6.5.6. *Множество геометрических объектов типа $H$ является* $\Omega_3$-*алгеброй.*

Доказательство. Пусть

$$m_i = \overline{\overline{e}}_g \circ W[g, \overline{\overline{e}}_g, m_i] \quad i = 1, ..., n$$

Для операции $\omega \in \Omega_3(n)$ мы положим

$$(6.5.11) \qquad m_1...m_n\omega = \overline{\overline{e}}_g \circ (W[g, \overline{\overline{e}}_g, m_1]...W[g, \overline{\overline{e}}_g, m_n]\omega)$$

Так как отображение $W[g, \overline{\overline{e}}_g, \overline{\overline{e}}_g \circ H(S)]$ для произвольного эндоморфизма $S$ $\Omega_2$-алгебры $A_2$ является эндоморфизмом $\Omega_3$-алгебры $A_3$, то определение (6.5.11) корректно.    □

Теорема 6.5.7. *Определено представление* $\Omega_1$-*алгебры* $A_1$ *в* $\Omega_3$-*алгебре* $N$ *геометрических объектов типа* $H$.

Доказательство. Пусть

$$m = \overline{\overline{e}}_g \circ W[g, \overline{\overline{e}}_g, m]$$

Для $a \in A_1$, мы положим

$$(6.5.12) \qquad f(a)(m) = \overline{\overline{e}}_g \circ f(a)(W[g, \overline{\overline{e}}_g, m])$$

Так как отображение $W[g, \overline{\overline{e}}_g, \overline{\overline{e}}_g \circ H(S)]$ для произвольного эндоморфизма $S$ $\Omega_2$-алгебры $A_2$ является эндоморфизмом представления $g$, то определение (6.5.12) корректно.    □

# Диаграмма представлений универсальных алгебр

## 7.1. Диаграмма представлений универсальных алгебр

Из сравнения теорем 6.1.4 и [14]-5.1 следует, что нет жёсткой границы между универсальной алгеброй и представлением универсальной алгебры. Отсюда следует возможность обобщения теории представлений универсальной алгебры.

Самая простая конструкция возникает следующим образом. Пусть

$$f_{12} : A_1 \longrightarrow A_2$$

представление $\Omega_1$-алгебры $A_1$ в $\Omega_2$-алгебре $A_2$. Если мы вместо $\Omega_2$-алгебры $A_2$ рассмотрим представление

$$f_{23} : A_2 \longrightarrow A_3$$

$\Omega_2$-алгебры $A_2$ в $\Omega_3$-алгебре $A_3$, то мы получим диаграмму вида

$$(7.1.1) \qquad A_1 \xrightarrow{f_{12}} A_2 \xrightarrow{f_{23}} A_3$$

Очевидно, что в диаграмме (7.1.1) мы можем положить, что $A_3$ является представлением

$$f_{34} : A_3 \longrightarrow A_4$$

Цепочку представлений универсальных алгебр можно сделать сколь угодно длинной. Таким образом мы получаем следующее определение.

Определение 7.1.1. *Рассмотрим множество $\Omega_k$-алгебр $A_k$, $k = 1, ..., n$. Положим $A = (A_1, ..., A_n)$. Положим $f = (f_{1\,2}, ..., f_{n-1\,n})$. Множество представлений $f_{k\,k+1}$, $k = 1, ..., n$, $\Omega_k$-алгебры $A_k$ в $\Omega_{k+1}$-алгебре $A_{k+1}$ называется* **башней** $(f, A)$ **представлений** $\Omega$-**алгебр**. $\qquad\square$

Башню представлений $(f, A)$ можно описать с помощью диаграммы

$$A_1 \xrightarrow{f_{1\,2}} A_2 \xrightarrow{f_{2\,3}} ... \xrightarrow{f_{n-1\,n}} A_n$$

Рассматривая башню представлений мы можем снова предположить, что $A_2$ или $A_3$ являются представлениями универсальных алгебр или башнями представлений. В этом случае диаграмма (7.1.1) примет вид

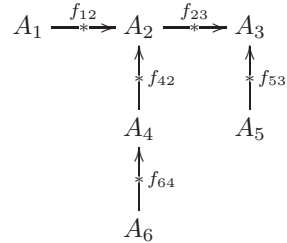

либо

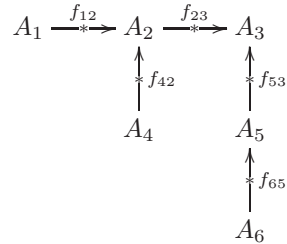

Мы также допускаем, что некоторые алгебры или отображения на диаграмме совпадают. Таким образом, мы будем считать, что диаграммы

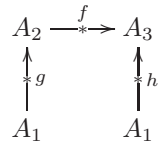

и

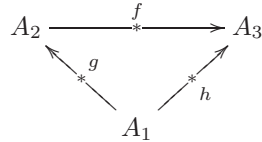

эквивалентны.

ОПРЕДЕЛЕНИЕ 7.1.2. **Диаграмма** $(f, A)$ **представлений универсальных алгебр** - *это такой ориентированный граф, что*

7.1.2.1: *вершина* $A_k$, $k = 1, ..., n$, *является* $\Omega_k$-*алгеброй;*

7.1.2.2: *ребро* $f_{kl}$ *является представлением* $\Omega_k$-*алгебры* $A_k$ *в* $\Omega_l$-*алгебре* $A_l$;

*Мы будем требовать, чтобы этот граф был связным и не содержал циклов. Мы будем полагать, что* $A_{[0]}$ *- это множество начальных вершин графа. Мы будем полагать, что* $A_{[k]}$ *- это множество вершин графа, для которых максимальный путь от начальных вершин равен* $k$. □

ЗАМЕЧАНИЕ 7.1.3. *Так как в разных вершинах графа может быть одна и таже алгебра, то мы обозначим* $A = (A_{(1)} \ ... \ A_{(n)})$ *множество универсальных алгебр, которые попарно различны. Из равенства*

$$A = (A_{(1)} \ ... \ A_{(n)}) = (A_1 \ ... \ A_n)$$

*следует, что для любого индекса $(i)$ существует по крайней мере один индекс $i$ такой, что $A_{(i)} = A_i$. Если даны два набора множеств $A = (A_{(1)} \ \dots \ A_{(n)})$, $B = (B_{(1)} \ \dots \ B_{(n)})$ и определено отображение*

$$h_{(i)} : A_{(i)} \to B_{(i)}$$

*для некоторого индекса $(i)$, то также определено отображение*

$$h_i : A_i \to B_i$$

*для любого индекса $i$ такого, что $A_{(i)} = A_i$ и в этом случае $h_i = h_{(i)}$.* $\qquad \square$

Теорема 7.1.4 (индукция по диаграмме представлений). *Пусть теорема $\mathcal{T}$ верна для множества универсальных алгебр $A_{[0]}$ диаграммы $(f, A)$ представлений универсальных алгебр. Пусть из утверждения, что теорема $\mathcal{T}$ верна для множества универсальных алгебр $A_{[k]}$ диаграммы $(f, A)$ представлений, следует утверждение, что теорема $\mathcal{T}$ верна для множества универсальных алгебр $A_{[k+1]}$ диаграммы $(f, A)$ представлений. Тогда теорема $\mathcal{T}$ верна для множества универсальных алгебр диаграммы $(f, A)$ представлений.*

Доказательство. Теорема является следствием принципа математической индукции. $\qquad \square$

Определение 7.1.5. *Диаграмма $(f, A)$ представлений универсальных алгебр называется* **коммутативной**, *если выполнено следующее условие. для каждой пары представлений*

$$f_{ik} : A_i \longrightarrow A_k$$

$$f_{jk} : A_j \longrightarrow A_k$$

*следующее равенство верно*[7.1]

$$(7.1.2) \qquad\qquad f_{ik}(a_i)(f_{jk}(a_j)(a_k)) = f_{jk}(a_j)(f_{ik}(a_i)(a_k))$$

$\qquad \square$

Теорема 7.1.6. *Пусть*

$$f_{ij} : A_i \longrightarrow A_j$$

*представление $\Omega_i$-алгебры $A_i$ в $\Omega_j$-алгебре $A_j$. Пусть*

$$f_{jk} : A_j \longrightarrow A_k$$

*представление $\Omega_j$-алгебры $A_j$ в $\Omega_k$-алгебре $A_k$. Мы можем описать фрагмент*[7.2]

$$A_i \xrightarrow{f_{ij}} A_j \xrightarrow{f_{jk}} A_k$$

---

[7.1] Образно говоря, представления $f_{ik}$ и $f_{jk}$ прозрачны друг для друга.

[7.2] Теорема 7.1.6 утверждает, что преобразования в башне представлений согласованы.

*диаграммы представлений с помощью диаграммы*

(7.1.3)

*Отображение*

$$f_{ijk} : A_i \to \mathrm{End}(\Omega_j, \mathrm{End}(\Omega_k, A_k))$$

*определено равенством*

(7.1.4) $$f_{ijk}(a_i)(f_{jk}(a_j)) = f_{jk}(f_{ij}(a_i)(a_j))$$

*где $a_i \in A_i$, $a_j \in A_j$. Если представление $f_{jk}$ эффективно и представление $f_{ij}$ свободно, то отображение $f_{ijk}$ является свободным представлением*

$$f_{ijk} : A_i \longrightarrow \mathrm{End}(\Omega_k, A_k)$$

$\Omega_i$-*алгебры $A_i$ в $\Omega_j$-алгебре $\mathrm{End}(\Omega_k, A_k)$.*

Доказательство.

**Лемма 7.1.7.** *Отображение $f_{ijk}$ является инъекцией.*

Доказательство. Пусть $(a_i, b_i) \in \ker f_{ijk}$. Тогда

(7.1.5) $$f_{jk}(f_{ijk}(a_i)(a_j)) = f_{ijk}(a_i)(f_{jk}(a_j)) = f_{ijk}(b_i)(f_{jk}(a_j))$$
$$= f_{jk}(f_{ij}(b_i)(a_j))$$

Если представление $f_{jk}$ эффективно, то равенство

(7.1.6) $$f_{ij}(a_i)(a_j) = f_{ij}(b_i)(a_j)$$

является следствием определения 3.1.2 и равенства (7.1.5) для любого $a_j \in A_j$. Утверждение $a_i = b_i$ следует из определения 3.1.4. ⊙

**Лемма 7.1.8.** *На множестве $\mathrm{End}(\Omega_j, \mathrm{End}(\Omega_k, A_k))$ определена структура $\Omega_i$-алгебры.*

Доказательство. Пусть $\omega \in \Omega_i$. Пусть $a_1, ..., a_m \in A_i$. Мы определим операцию $\omega$ на множестве $\mathrm{End}(\Omega_j, \mathrm{End}(\Omega_k, A_k))$ с помощью равенства

(7.1.7) $$f_{ijk}(a_1)...f_{ijk}(a_m)\omega = f_{ijk}(a_1...a_m\omega)$$

Согласно лемме 7.1.7, операция $\omega$ корректно определена равенством (7.1.7). ⊙

Следствие 7.1.9. *Отображение $f_{ijk}$ является гомоморфизмом $\Omega_i$-алгебры.*          ⊙

Лемма 7.1.10. *Отображение $f_{ijk}(a)$ является гомоморфизмом $\Omega_j$-алгебры.*

Доказательство. Пусть $b_1, ..., b_m \in A_j$. Тогда равенство

$$(7.1.8) \quad f_{ijk}(a)(f_{jk}(b_1))...f_{ijk}(a)(f_{jk}(b_m))\omega = f_{jk}(f_{ij}(a)(b_1))...f_{jk}(f_{ij}(a)(b_m))\omega$$

является следствием равенства (7.1.4). Так как отображения $f_{ij}(a)$, $f_{jk}$ являются гомоморфизмами $\Omega_j$-алгебры, то равенство

$$(7.1.9) \quad \begin{aligned} &f_{ijk}(a)(f_{jk}(b_1))...f_{ijk}(a)(f_{jk}(b_m))\omega \\ &= f_{jk}(f_{ij}(a)(b_1)...f_{ij}(a)(b_m)\omega) \\ &= f_{jk}(f_{ij}(a)(b_1...b_m\omega)) \end{aligned}$$

является следствием равенства (7.1.8). Равенство

$$(7.1.10) \quad f_{ijk}(a)(f_{jk}(b_1))...f_{ijk}(a)(f_{jk}(b_m))\omega = f_{ijk}(a)(f_{jk}(b_1...b_m\omega))$$

является следствием равенств (7.1.4), (7.1.9). Так как отображение $f_{jk}$ является гомоморфизмом $\Omega_2$-алгебры, то равенство

$$(7.1.11) \quad f_{ijk}(a)(f_{jk}(b_1))...f_{ijk}(a)(f_{jk}(b_m))\omega = f_{ijk}(a)(f_{jk}(b_1)...f_{jk}(b_m))\omega$$

является следствием равенства (7.1.10).          ⊙

Теорема является следствием следствия 7.1.9 и леммы 7.1.10.          □

Теорема 7.1.11. *Отображение $f_{jk}$ является приведенным морфизмом представлений из $f_{ij}$ в $f_{ijk}$.*

Доказательство. Рассмотрим более детально диаграмму (7.1.3).

$$(7.1.12)$$

Утверждение теоремы следует из равенства (7.1.4) и определения 3.4.2.          □

Теорема 7.1.12. *Пусть*

$$f_{ij} : A_i \longrightarrow_* A_j$$

*представление $\Omega_i$-алгебры $A_i$ в $\Omega_j$-алгебре $A_j$. Пусть*

$$f_{jk} : A_j \longrightarrow_* A_k$$

*представление $\Omega_j$-алгебры $A_j$ в $\Omega_k$-алгебре $A_k$. Тогда существует представление*

$$f_{ij,k} : A_i \times A_j \longrightarrow_* A_k$$

*множества[7.3] $A_i \times A_j$ в $\Omega_k$-алгебре $A_k$.*

---

[7.3] Так как $\Omega_i$-алгебра $A_i$ и $\Omega_j$-алгебра $A_j$ имеют различный набор операций, мы не можем определить структуру универсальной алгебры на множестве $A_i \times A_j$.

Доказательство. Мы можем описать фрагмент

$$A_i \xrightarrow{f_{ij}} A_j \xrightarrow{f_{jk}} A_k$$

диаграммы представлений с помощью диаграммы

(7.1.13)

Из диаграммы (7.1.13) следует, что отображение $f_{ij,k}$ определено равенством

$$f_{ij,k}(a_i, a_j) = {\color{teal} f_{jk}(f_{ij}(a_i)(a_j))}$$

$\square$

## 7.2. Морфизм диаграммы представлений

Определение 7.2.1. *Пусть* $(f, A)$ *- диаграмма представлений, где* $A = (A_{(1)} \ ... \ A_{(n)})$ *- множество универсальных алгебр. Пусть* $(B, g)$ *- диаграмма представлений, где* $B = (B_{(1)} \ ... \ B_{(n)})$ *- множество универсальных алгебр. Множество отображений* $h = (h_{(1)} \ ... \ h_{(n)})$

$$h_{(i)} : A_{(i)} \to B_{(i)}$$

*называется* **морфизмом из диаграммы представлений** $(f, A)$ **в диаграмму представлений** $(B, g)$, *если для любых индексов* $(i)$, $(j)$, $i$, $j$ *таких, что* $A_{(i)} = A_i$, $A_{(j)} = A_j$, *и для каждого представления*

$$f_{ji} : A_j \xrightarrow{\quad\quad} A_i$$

*пара отображений* $(h_j \ h_i)$ *является морфизмом представлений из* $f_{ji}$ *в* $g_{ji}$. $\square$

Мы будем пользоваться записью

$$h : A \to B$$

если кортеж отображений $h$ является морфизмом из диаграммы представлений $(f, A)$ в диаграмму представлений $(B, g)$.

Очень часто при изучении морфизма представления универсальной алгебры, мы предполагаем, что алгебра, порождающая представление, задана. Поэтому нас не интересует отображение этой алгебры, и это соглашение упрощает структуру морфизма. Такой морфизм представления мы называем приведенным морфизмом представления.

Похожая задача встречается при изучении морфизма диаграммы представлений. Для каждой универсальной алгебры из диаграммы представлений существует множество алгебр, предшествующих этой алгебре в соответствующем графе. Мы можем предположить, что некоторые из этих алгебр заданы и не рассматривать соответствующие гомоморфизмы. Соответствующий морфизм диаграммы представлений также называется приведенным. Однако в виду сложности структуры диаграммы представлений, мы не будем рассматривать приведенные морфизмы диаграммы представлений.

Для любого представления $f_{ij}$, $i = 1, ..., n$, $j = 1, ..., n$, мы имеем диаграмму

$$(7.2.1)$$

Равенства

$$(7.2.2) \qquad h_j \circ f_{ij}(a_i) = g_{ij}(h_i(a_i)) \circ h_j$$

$$(7.2.3) \qquad h_j(f_{ij}(a_i)(a_j)) = g_{ij}(h_i(a_i))(h_j(a_j))$$

выражают коммутативность диаграммы (1).

Пусть определены представления $f_{ij}$ и $f_{jk}$ универсальных алгебр. Учитывая диаграмму (7.2.1) для представлений $f_{ij}$ и $f_{jk}$, мы получим следующую диаграмму

$$(7.2.4)$$

Очевидно, что существует морфизм из $\operatorname{End}(\Omega_k, A_k)$ в $\operatorname{End}(\Omega_k, B_k)$, отображающий $f_{ijk}(a_i)$ в $g_{ijk}(h_k(a_i))$.

Теорема 7.2.2. *Если представление $f_{jk}$ эффективно и представление $f_{ij}$ свободно, то*[7.4] *($h_i, h_k^*$) является морфизмом представлений из представления $f_{ijk}$ в представление $g_{ijk}$ $\Omega_i$-алгебры.*

Доказательство. Рассмотрим диаграмму

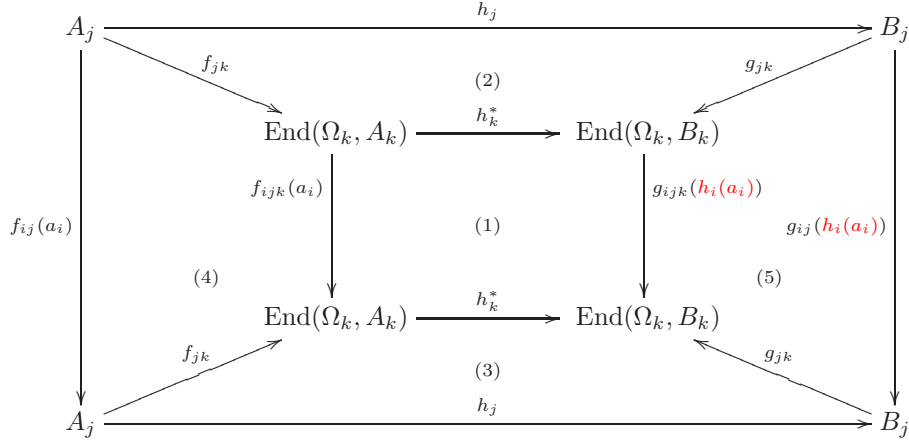

Существование отображения $h_k^*$ и коммутативность диаграмм (2) и (3) следует из эффективности отображения $f_{jk}$ и теоремы 3.2.9. Коммутативность диаграмм (4) и (5) следует из теоремы 7.1.11.

Из коммутативности диаграммы (4) следует

$$(7.2.5) \qquad f_{jk} \circ f_{ij}(a_i) = f_{ijk}(a_i) \circ f_{jk}$$

Из равенства (7.2.5) следует

$$(7.2.6) \qquad h_k^* \circ f_{jk} \circ f_{ij}(a_i) = h_k^* \circ f_{ijk}(a_i) \circ f_{jk}$$

Из коммутативности диаграммы (3) следует

$$(7.2.7) \qquad h_k^* \circ f_{jk} = g_{jk} \circ h_j$$

Из равенства (7.2.7) следует

$$(7.2.8) \qquad h_k^* \circ f_{jk} \circ f_{ij}(a_i) = g_{jk} \circ h_j \circ f_{ij}(a_i)$$

Из равенств (7.2.6) и (7.2.8) следует

$$(7.2.9) \qquad h_k^* \circ f_{ijk}(a_i) \circ f_{jk} = g_{jk} \circ h_j \circ f_{ij}(a_i)$$

Из коммутативности диаграммы (5) следует

$$(7.2.10) \qquad g_{jk} \circ g_{ij}(h_i(a_i)) = g_{ijk}(h_i(a_i)) \circ g_{jk}$$

Из равенства (7.2.10) следует

$$(7.2.11) \qquad g_{jk} \circ g_{ij}(h_i(a_i)) \circ h_j = g_{ijk}(h_i(a_i)) \circ g_{jk} \circ h_j$$

Из коммутативности диаграммы (2) следует

$$(7.2.12) \qquad h_k^* \circ f_{jk} = g_{jk} \circ h_j$$

Из равенства (7.2.12) следует

$$(7.2.13) \qquad g_{ijk}(h_i(a_i)) \circ h_k^* \circ f_{jk} = g_{ijk}(h_i(a_i)) \circ g_{jk} \circ h_j$$

---

[7.4] Смотри определение отображения $h^*$ в теореме 3.2.9.

Из равенств (7.2.11) и (7.2.13) следует

$$(7.2.14) \qquad g_{jk} \circ g_{ij}(h_i(a_i)) \circ h_j = g_{ij}(h_i(a_i)) \circ h_k^* \circ f_{jk}$$

Внешняя диаграмма является диаграммой (7.2.1) при $i = 1$. Следовательно, внешняя диаграмма коммутативна

$$(7.2.15) \qquad h_j \circ f_{ij}(a_i) = g_{ij}(h_i(a_i)) \circ h_j$$

Из равенства (7.2.15) следует

$$(7.2.16) \qquad g_{jk} \circ h_J \circ f_{ij}(a_i) = g_{jk} \circ g_{ij}(h_i(a_i)) \circ h_j(a_j)$$

Из равенств (7.2.9), (7.2.14) и (7.2.16) следует

$$(7.2.17) \qquad h_k^* \circ f_{ijk}(a_i) \circ f_{jk} = g_{ijk}(h_i(a_i)) \circ h_k^* \circ f_{jk}$$

Так как отображение $f_{i+1,i+2}$ - инъекция, то из равенства (7.2.17) следует

$$(7.2.18) \qquad h_k^* \circ f_{ijk}(a_i) = g_{ijk}(h_i(a_i)) \circ h_k^*$$

Из равенства (7.2.18) следует коммутативность диаграммы (1), откуда следует утверждение теоремы. $\qquad\qquad\square$

Теорема 7.2.2 утверждает, что неизвестное отображение на диаграмме (7.2.4) является отображением $h_k^*$. Смысл теорем 7.1.11 и 7.2.2 состоит в том, что, если все представления свободны, то все отображения в диаграмме представлений действуют согласованно.

Теорема 7.2.3. *Рассмотрим множество $\Omega_i$-алгебр $A_{(i)}$, $B_{(i)}$, $C_{(i)}$, $(i) = (1), ..., (n)$. Пусть определены морфизмы диаграмм представлений*

$$p :(f, A) \to (g, B)$$
$$q :(g, B) \to (h, C)$$

*Тогда определён морфизм представлений $\Omega$-алгебры*

$$r : (f, A) \to (h, C)$$

*где $r_{(k)} = q_{(k)} \circ p_{(k)}$, $(k) = (1), ..., (n)$. Мы будем называть морфизм $r$ диаграммы представлений из $f$ в $h$* **произведением морфизмов $p$ и $q$ диаграммы представлений**.

Доказательство. Для любых $i$, $j$ таких, что $A_{(j)} = A_j$, если существует представление $f_{ij}$, мы можем представить утверждение теоремы, пользуясь

диаграммой

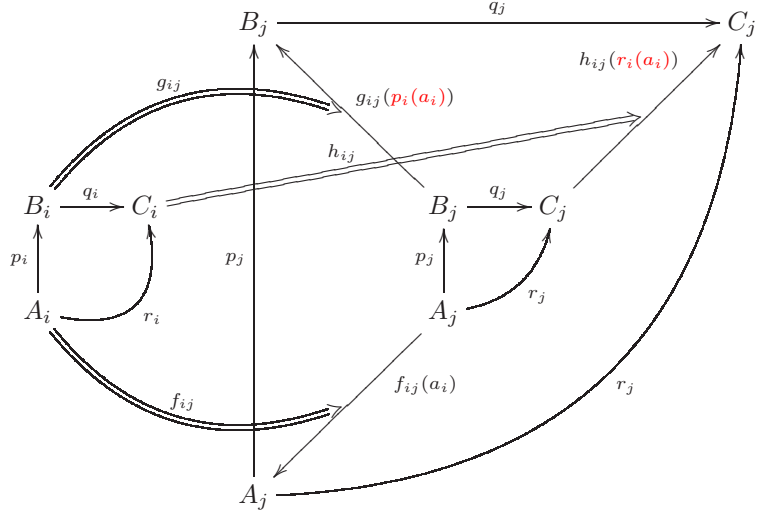

Отображение $r_i$ является гомоморфизмом $\Omega_i$-алгебры $A_i$ в $\Omega_i$-алгебру $C_i$. Нам надо показать, что пара отображений $(r_i, r_j)$ удовлетворяет (7.2.2):

$$r_k(f_{ij}(a_i)(a_j)) = q_j \circ p_j(f_{ij}(a_i)(a_j))$$
$$= q_k(g_{ij}(p_i(a_i))(p_j(a_j)))$$
$$= h_{ij}(q_i \circ p_i(a_i))(q_j \circ p_j(a_j))$$
$$= h_{ij}(r(a_i))(r_j(a_j))$$

□

## 7.3. Автоморфизм диаграммы представлений

Определение 7.3.1. *Пусть* $(f, A)$ *- диаграмма представлений универсальных алгебр. Морфизм диаграммы представлений* $(h_{(1)}, ..., h_{(n)})$ *такой, что для любого* $(k)$, $(k) = (1), ..., (n)$, $h_{(k)}$ *является эндоморфизмом* $\Omega_{(k)}$*-алгебры* $A_{(k)}$, *называется* **эндоморфизмом диаграммы представлений**. □

Определение 7.3.2. *Пусть* $(f, A)$ *- диаграмма представлений универсальных алгебр. Морфизм диаграммы представлений* $(h_{(1)}, ..., h_{(n)})$ *такой, что для любого* $(k)$, $(k) = (1), ..., (n)$, $h_{(k)}$ *является автоморфизмом* $\Omega_{(k)}$*-алгебры* $A_{(k)}$, *называется* **автоморфизмом диаграммы представлений**. □

Теорема 7.3.3. *Пусть* $(f, A)$ *- диаграмма представлений универсальных алгебр. Множество автоморфизмов диаграммы* $(f, A)$ *представлений порождает* **группу** $GA(f)$.

Доказательство. Пусть $r = (r_{(1)}, ..., r_{(n)})$, $p = (p_{(1)}, ..., p_{(n)})$ - автоморфизмы диаграммы представлений $(f, A)$. Согласно определению 7.3.2 для любого $(k)$, $(k) = (1), ..., (n)$, отображения $r_{(k)}$, $p_{(k)}$ являются автоморфизмами $\Omega_{(k)}$-алгебры $A_{(k)}$. Согласно теореме II.3.2 ([14], с. 60) для любого $(k)$, $(k) = (1), ..., (n)$, отображение $r_{(k)} \circ p_{(k)}$ является автоморфизмом $\Omega_{(k)}$-алгебры $A_{(k)}$. Из теоремы 7.2.3 и определения 7.3.2 следует, что произведение автоморфизмов $r \circ p$ диаграммы представлений $(f, A)$ является автоморфизмом диаграммы представлений $(f, A)$.

Согласно доказательству теоремы 3.5.5, для любого $(k)$, $(k) = (1)$, ..., $(n)$, произведение автоморфизмов $\Omega_{(k)}$-алгебры ассоциативно. Следовательно, ассоциативно произведение автоморфизмов диаграммы представлений.

Пусть $r = (r_{(1)}, ..., r_{(n)})$ - автоморфизм диаграммы представлений $(f, A)$. Согласно определению 7.3.2 для любого $(k)$, $(k) = (1)$, ..., $(n)$, отображение $r_{(k)}$ является автоморфизмом $\Omega_{(k)}$-алгебры $A_{(k)}$. Следовательно, для любого $(k)$, $(k) = (1)$, ..., $(n)$, отображение $r_{(k)}^{-1}$ является автоморфизмом $\Omega_{(k)}$-алгебры $A_{(k)}$. Для автоморфизма $r = (r_{(1)}, ..., r_{(n)})$ справедливо равенство (7.2.3). Пусть $j$ - индекс такой, что $j = (k)$. Положим $a'_j = r_j(a_j)$. Так как $r_j$ - автоморфизм, то $a_j = r_j^{-1}(a'_j)$ и равенство (7.2.3) можно записать в виде

$$(7.3.1) \qquad h_j(f_{ij}(h_i^{-1}(a'_i))(h_j^{-1}(a'_j))) = g_{ij}(a'_i)(a'_j)$$

для любых $i$, $j$ при условии, что представление $f_{ij}$ существует. Аналогично, из равенства (7.3.1) следует

$$(7.3.2) \qquad f_{ij}(h_i^{-1}(a'_i))(h_j^{-1}(a'_j)) = h_j^{-1}(g_{ij}(a'_i)(a'_j))$$

Равенство (7.3.2) соответствует равенству (7.2.3) для отображения $r^{-1}$. Следовательно, отображение $r^{-1}$ является автоморфизмом диаграммы представлений $(f, A)$. $\qquad \square$

# Базис диаграммы представлений универсальной алгебры

## 8.1. Множество образующих диаграммы представлений

Мы строим базис диаграммы представлений по той же схеме, что мы построили базис представления в секции 6.2.

ОПРЕДЕЛЕНИЕ 8.1.1. *Пусть* $(f, A)$ *- диаграмма представлений. Кортеж множеств*

$$N = (N_{(1)} \subset A_{(1)}, ..., N_{(n)} \subset A_{(n)}) = (N_1 \subset A_1, ..., N_n \subset A_n)$$

*называется* **кортежем стабильных множеств диаграммы представлений** $(f, A)$, *если*

$$f_{ij}(a_i)(a_j) \in N_j \quad i, j = 1, ..., n$$

*для любых* $a_1 \in N_1$, ..., $a_n \in N_n$, *при условии, что существует представление* $f_{ij}$. *Мы также будем говорить, что кортеж множеств*

$$N = (N_{(1)} \subset A_{(1)}, ..., N_{(n)} \subset A_{(n)}) = (N_1 \subset A_1, ..., N_n \subset A_n)$$

*стабилен относительно диаграммы представлений* $(f, A)$. □

ТЕОРЕМА 8.1.2. *Пусть* $(f, A)$ *- диаграмма представлений. Пусть множество* $N_{(i)} \subset A_{(i)}$ *является подалгеброй* $\Omega_{(i)}$-*алгебры* $A_{(i)}$, $(i) = (1)$, ..., $(n)$. *Пусть кортеж множеств*

$$N = (N_{(1)} \subset A_{(1)}, ..., N_{(n)} \subset A_{(n)}) = (N_1 \subset A_1, ..., N_n \subset A_n)$$

*стабилен относительно диаграммы представлений* $(f, A)$. *Тогда существует диаграмма представлений*

$$(8.1.1) \qquad (N, f_N = (f_{Nij}))$$

*такая, что*

$$f_{Nij}(a_i) = f_{ij}(a_i)|_{N_j} \quad i = 1, ..., n \quad j = 1, ..., n$$

*Диаграмма представлений* (8.1.1) *называется* **подпредставлением** *диаграммы представлений* $(f, A)$.

ДОКАЗАТЕЛЬСТВО. Пусть $\omega$ - $m$-арная операция $\Omega_i$-алгебры $A_i$, $i = 1$, ..., $n$. Тогда для любых $a_{i,1}$, ..., $a_{i,m} \in N_i$ и любого $a_j \in N_j$

$$(f_{Nij}(a_{i,1})...f_{Nij}(a_{i,m})\omega)(a_j) = (f_{ij}(a_{i,1})...f_{ij}(a_{i,m})\omega)(a_j)$$
$$= f_{ij}(a_{i,1}...a_{i,m}\omega)(a_j)$$
$$= f_{Nij}(a_{i,1}...a_{i,m}\omega)(a_i)$$

Пусть $\omega$ - $m$-арная операция $\Omega_j$-алгебры $A_j$, $j = 1, ..., n$. Тогда для любых $a_{j,1}$, ..., $a_{j,m} \in N_j$ и любого $a_i \in N_i$

$$f_{Nij}(a_i)(a_{j,1})...f_{Nij}(a_i)(a_{i,m})\omega = f_{ij}(a_i)(a_{j,1})...f_{ij}(a_i)(a_{j,m})\omega$$
$$= f_{ij}(a_i)(a_{j,1}...a_{j,m}\omega)$$
$$= f_{Nij}(a_i)(a_{j,1}...a_{j,m}\omega)$$

Утверждение теоремы доказано. $\qquad\square$

Из теоремы 8.1.2 следует, что если диаграммы представлений (8.1.1) являются диаграммой подпредставлений диаграммы представлений $(f, A)$, то отображение

$$(id_{(1)} : N_{(1)} \to A_{(1)}, ..., id_{(n)} : N_{(n)} \to A_{(n)})$$

является морфизмом диаграмм представлений.

ТЕОРЕМА 8.1.3. *Множество*[8.1] $\mathcal{B}[f, A]$ *всех диаграмм подпредставлений диаграммы представлений* $(f, A)$ *порождает систему замыканий на диаграмме представлений* $(f, A)$ *и, следовательно, является полной структурой.*

ДОКАЗАТЕЛЬСТВО. Пусть для данного $\lambda \in \Lambda$,

$$K_\lambda = (K_{\lambda,(1)} \subset A_{(1)}, ..., K_{\lambda,(n)} \subset A_{(n)})$$

кортеж множеств, стабильных относительно диаграммы представлений $(f, A)$. Операцию пересечения на множестве $\mathcal{B}[f, A]$ мы определим согласно правилу

$$\bigcap f_{K_\lambda ij} = f_{\cap K_\lambda ij} \quad i, j = 1, ..., n$$
$$\bigcap K_\lambda = \left(K_{(1)} = \bigcap K_{\lambda,(1)}, ..., K_{(n)} = \bigcap K_{\lambda,(n)}\right)$$

$\cap K_{\lambda,(i)}$ - подалгебра $\Omega_{(i)}$-алгебры $A_{(i)}$. Пусть

$$a_j \in \bigcap K_{\lambda,j} = \bigcap K_{\lambda,(j)}$$

Для любого $\lambda \in \Lambda$ и для любого $a_i \in K_i = K_{(i)}$

$$f_{ij}(a_i)(a_j) \in K_{\lambda,j} = K_{\lambda,(j)}$$

Следовательно,

$$f_{ij}(a_i)(a_j) \in K_j = K_{(j)}$$

Следовательно, операция пересечения диаграмм подпредставлений определена корректно. $\qquad\square$

Обозначим соответствующий оператор замыкания через $J[f]$. Если

$$X = (X_{(1)} \subset A_{(1)}, ..., X_{(n)} \subset A_{(n)}) = (X_1 \subset A_1, ..., X_n \subset A_n)$$

кортеж множеств, то

$$J[f, X] = (J_{(1)}[f, X], ..., J_{(n)}[f, X]) = (J_1[f, X], ..., J_n[f, X])$$

является пересечением всех кортежей

$$K = (K_{(1)} \subset A_{(1)}, ..., K_{(n)} \subset A_{(n)}) = (K_1 \subset A_1, ..., K_n \subset A_n)$$

стабильных относительно диаграммы представлений $(f, A)$ и таких, что для $(i) = (1), ..., (n)$, $K_{(i)}$ - подалгебра $\Omega_{(i)}$-алгебры $A_{(i)}$, содержащая $X_{(i)}$.

---

[8.1]Эта теорема аналогична определению структуры подалгебр ([14], стр. 93, 94)

Теорема 8.1.4. *Пусть*[8.2] *$(f, A)$ - диаграмма представлений. Пусть*

$$X = (X_{(1)} \subset A_{(1)}, ..., X_{(n)} \subset A_{(n)}) = (X_1 \subset A_1, ..., X_n \subset A_n)$$

*Для каждого значения $(i)$, $(i) = (1), ..., (n)$, определим подмножества $X_{(i)k} \subset A_{(i)}$ индукцией по $k$.*

8.1.4.1: $X_{(i)0} = X_{(i)}$

8.1.4.2: $x \in X_{(i)k} => x \in X_{(i)k+1}$

8.1.4.3: $x_1 \in X_{(i)k}, ..., x_p \in X_{(i)k}, \omega \in \Omega_{(i)}(p) => x_1...x_p\omega \in X_{(i)k+1}$

8.1.4.4: $x_i \in X_{ik} = X_{(i)k}, x_j \in X_{jk} = X_{(j)k} => f_{ji}(x_j)(x_i) \in X_{(i)k+1}$

*Для каждого значения $(i)$ положим*

$$Y_{(i)} = \bigcup_{m=0}^{\infty} X_{(i)m}$$

*Тогда*

$$J_{(i)}[f, X] = Y_{(i)} \qquad (i) = (1), ..., (n)$$

Доказательство. Для каждого значения $(i)$ доказательство теоремы совпадает с доказательством теоремы 6.1.4. □

$J[f, X]$ называется **подпредставлением** диаграммы представлений $(f, A)$, порождённым кортежем множеств $X$, а $X$ - **множеством образующих** диаграммы представлений J[f,X]. В частности, множество образующих диаграммы представлений $(f, A)$ будет такой кортеж

$$X = (X_{(1)} \subset A_{(1)}, ..., X_{(n)} \subset A_{(n)}) = (X_1 \subset A_1, ..., X_n \subset A_n)$$

что $J[f, X] = A$.

Из теоремы 8.1.4 следует следующее определение.

Определение 8.1.5. *Пусть*

$$X = (X_{(1)} \subset A_{(1)}, ..., X_{(n)} \subset A_{(n)}) = (X_1 \subset A_1, ..., X_n \subset A_n)$$

*кортеж множеств. Для любого кортежа $A$-чисел $a \in J[f, X]$*

$$a = (a_{(1)} \quad ... \quad a_{(n)}) = (a_1 \quad ... \quad a_n)$$

*существует* **кортеж $\Omega$-слов**

$$w[f, X, a] = (w_{(1)}[f, X, a_{(1)}], ..., w_{(n)}[f, X, a_{(n)}])$$
$$= (w_1[f, X, a_1], ..., w_n[f, X, a_n])$$

*определённых согласно следующему правилу.*

8.1.5.1: *Если* $a_{(i)} \in X_{(i)}$, $(i) = (1), ..., (n)$, *то $a_{(i)}$ - $\Omega_{(i)}$-слово*

$$w_{(i)}[f, X, a_{(i)}] = a_{(i)}$$

8.1.5.2: *Если* $a_{(i)1}, ..., a_{(i)p}$ *- $\Omega_{(i)}$-слова, $(i) = (1), ..., (n)$, и $\omega \in \Omega_{(i)}(p)$, то $a_{(i)1}...a_{(i)p}\omega$ - $\Omega_{(i)}$-слово.*

8.1.5.3: *Пусть $a_i = a_{(i)}$ - $\Omega_{(i)}$-слово, $a_j = a_{(j)}$ - $\Omega_{(j)}$-слово. Пусть существует представление $f_{ij}$. Тогда $f_{ij}(a_i)(a_j)$ - $\Omega_{(j)}$-слово.*

*Обозначим* $w[f, X]$ **множество кортежей $\Omega$-слов** *диаграммы представлений J[f,X].* □

---

[8.2] Утверждение теоремы аналогично утверждению теоремы 5.1, [14], стр. 94.

Мы рассматриваем кортеж $A$-чисел в определении 8.1.5, так как нам нужен алгоритм формирования кортежа $\Omega$-слов. Однако при решении конкретных задач нам может понадобиться только некоторое подмножество кортежа $A$-чисел. Например, в аффинном пространстве мы можем рассматривать либо множество точек, либо множество векторов.

Выбор $\Omega_{(i)}$-слова относительно множества образующих $X$ неоднозначен. Поэтому, если $\Omega_{(i)}$-число имеет различные $\Omega_{(i)}$-слова, то мы, чтобы их отличать, будем пользоваться индексами: $w_{(i)}[f, X, m], w_{(i)1}[f, X, m], w_{(i)2}[f, X, m]$.

ОПРЕДЕЛЕНИЕ 8.1.6. *Множество образующих $X$ диаграммы представлений $(f, A)$ порождает кортеж отношений эквивалентности*

$$\rho[f, X] = (\rho_{(1)}[f, X], ..., \rho_{(n)}[f, X])$$
$$\rho_{(i)}[f, X] = \{(w_{(i)}[f, X, m_{(i)}], w_{(i)1}[f, X, m_{(i)}]) : m_{(i)} \in A_{(i)}\}$$

*на множестве кортежей $\Omega$-слов.* ☐

Согласно определению 8.1.6, два $\Omega_{(i)}$-слова относительно множества образующих $X$ диаграммы представлений $(f, A)$ эквивалентны тогда и только тогда, когда они соответствуют одному и тому же $A_{(i)}$-числу. Когда мы будем записывать равенство двух $\Omega_{(i)}$-слов относительно множества образующих $X$ диаграммы представлений $(f, A)$, мы будем иметь в виду, что это равенство верно с точностью до отношения эквивалентности $\rho_{(i)}[f, X]$.

Мы будем пользоваться записью

$$r(a) = (r_{(1)}(a_{(1)}), ..., r_{(n)}(a_{(n)}))$$

для образа кортежа элементов $a = (a_{(1)}, ..., a_{(n)})$ при морфизме диаграммы представлений.

ТЕОРЕМА 8.1.7. *Пусть $X$ - множество образующих диаграммы представлений $(f, A)$. Пусть $Y$ - множество образующих диаграммы представлений $(g, B)$. Морфизм $r$ диаграммы представлений $(f, A)$ порождает отображение $\Omega$-слов*

$$w[f \to g, X, Y, r] : w[f, X] \to w[g, Y]$$

$$X_{(i)} \subset A_{(i)} \quad Y_{(i)} = r_{(i)}(X_{(i)}) \quad (i) = (1), ..., (n)$$

*такое, что для любого $(i)$, $(i) = (1), ..., (n)$,*

8.1.7.1: *Если $a_{(i)} \in X_{(i)}, a'_{(i)} = r_{(i)}(a_{(i)})$, то*

$$w_{(i)}[f \to g, X, Y, r](a_{(i)}) = a'_{(i)}$$

8.1.7.2: *Если*

$$a_{(i)1}, ..., a_{(i)p} \in w_{(i)}[f, X]$$

$$a'_{(i)1} = w_{(i)}[f \to g, X, Y, r](a_{(i)1}) \quad ... \quad a'_{(i)p} = w_{(i)}[f \to g, X, Y, r](a_{(i)p})$$

*то для операции $\omega \in \Omega_{(i)}(p)$ справедливо*

$$w_{(i)}[f \to g, X, Y, r](a_{(i),1}...a_{(i),p}\omega) = a'_{(i),1}...a'_{(i),p}\omega$$

8.1.7.3: *Если*

$$a_i = a_{(i)} \in w_{(i)}[f, X] \quad a'_{(i)} = w_{(i)}[f \to g, X, Y, r](a_{(i)})$$

$$a_j = a_{(j)} \in w_{(j)}[f, X] \quad a'_j = a'_{(j)} = w_{(j)}[f, X, r](a_{(j)})$$

*то*

$$w_{(i)}[f \to g, X, Y, r](f_{ji}(a_j)(a_i)) = g_{ji}(a'_j)(a'_i)$$

Доказательство. Утверждения 8.1.7.1, 8.1.7.2 справедливы в силу определения морфизма $r$. Утверждение 8.1.7.3 следует из равенства (7.2.3).    □

Замечание 8.1.8. *Пусть $r$ - морфизм диаграммы представлений $(f, A)$ в диаграмму представлений $(g, B)$. Пусть*

$$a \in J[f, X] \quad a' = r(a) \quad Y = r(X)$$

*Теорема 8.1.7 утверждает, что $a' \in J[g, Y]$. Теорема 8.1.7 также утверждает, что кортеж $\Omega$-слов, представляющий $a$ относительно $X$, и кортеж $\Omega$-слов, представляющий $a'$ относительно $Y$, формируются согласно одному и тому же алгоритму. Это позволяет рассматривать кортеж $\Omega$-слов $w[g, Y, a']$ как кортеж отображений*

$$W[f, X, a] = (W_{(1)}[f, X, a], ..., W_{(n)}[f, X, a]) = (W_1[f, X, a], ..., W_n[f, X, a])$$

(8.1.2)    $W_{(k)}[f, X, a] : (g, X') \to (g, X') \circ W_{(k)}[f, X, a] = w_{(k)}[g, X', a']$

*Если $f = g$, то вместо отображения (8.1.2) мы будем рассматривать отображение*

$$W_{(k)}[f, X, a] : Y \to Y \circ W_{(k)}[f, X, a] = w_{(k)}[f, Y, a']$$
$$W_{(k)}[f, X, a](Y) = Y \circ W_{(k)}[f, X, a]$$

*такое, что, если для некоторого морфизма $r$*

$$Y = r(X) \quad a' = r(a)$$

*то*

$$W_{(k)}[f, X, a](Y) = Y \circ W[f, X, a] = w[f, Y, a'] = a'$$

*Отображение $W_{(k)}[f, X, a]$ называется* **координатами** $A_{(k)}$-**числа** $a_{(k)}$ *относительно кортежа множеств $X$. Аналогично, мы можем рассмотреть координаты множества $B \subset J_{(k)}[f, X]$ относительно множества $X$*

$$W_{(k)}[f, X, B] = \{W_{(k)}[f, X, a] : a \in B\} = (W_{(k)}[f, X, a], a \in B)$$

*Обозначим*

$$W[f, X] = (W_{(1)}[f, X], ..., W_{(n)}[f, X]) = (W_1[f, X], ..., W_n[f, X])$$

$$W_{(k)}[f, X] = \{W_{(k)}[f, X, a] : a \in J_{(k)}[f, X]\} = (W_{(k)}[f, X, a], a \in J_{(k)}[f, X])$$

**множество координат представления** $J[f, X]$.    □

Теорема 8.1.9. *На множестве координат $W_{(k)}[f, X]$ определена структура $\Omega_{(k)}$-алгебры.*

Доказательство. Пусть $\omega \in \Omega_{(k)}(n)$. Тогда для любых $m_1, ..., m_n \in J_{(k)}[f, X]$ положим

$$(8.1.3) \qquad W_{(k)}(f, X, m_1)...W_{(k)}(f, X, m_n)\omega = W_{(k)}(f, X, m_1...m_n\omega)$$

Согласно замечанию 8.1.8, из равенства (8.1.3) следует

$$(8.1.4) \qquad \begin{aligned} X \circ (W_{(k)}[f, X, m_1]...W_{(k)}[f, X, m_n]\omega) &= X \circ W_{(k)}[f, X, m_1...m_n\omega] \\ &= w_{(k)}[f, X, m_1...m_n\omega] \end{aligned}$$

Согласно правилу 8.1.5.2, из равенства (8.1.4) следует

$$(8.1.5) \qquad \begin{aligned} & X \circ (W_{(k)}[f, X, m_1]...W_{(k)}[f, X, m_n]\omega) \\ &= w_{(k)}[f, X, m_1]...w_{(k)}[f, X, m_n]\omega \\ &= (X \circ W_{(k)}[f, X, m_1])...(X \circ W_{(k)}[f, X, m_n])\omega \end{aligned}$$

Из равенства (8.1.5) следует корректность определения (8.1.3) операции $\omega$ на множестве координат $W_{(k)}[f, X]$. $\qquad \square$

Теорема 8.1.10. *Если определено представление $f_{jk}$ $\Omega_j$-алгебры $A_j$ в $\Omega_k$-алгебре $A_k$, то определено представление $F_{jk}$ $\Omega_j$-алгебры $W_j[f, X]$ в $\Omega_k$-алгебре $W_k[f, X]$.*

Доказательство. Пусть $a_j \in J_j[f, X]$. Тогда для любого $a_k \in J_k[f, X]$, положим

$$(8.1.6) \qquad F_{jk}(W_j[f, X, a_j])(W_k[f, X, a_k]) = W_k[f, X, f_{jk}(a_j)(a_k)]$$

Согласно замечанию 8.1.8, из равенства (8.1.6) следует

$$(8.1.7) \qquad \begin{aligned} X \circ (F_{jk}(W_j[f, X, a_j])(W_k[f, X, a_k])) &= X \circ W_k[f, X, f_{jk}(a_j)(a_k)] \\ &= w_k[f, X, f_{jk}(a_j)(a_k)] \end{aligned}$$

Согласно правилу 8.1.5.3, из равенства (8.1.7) следует

$$(8.1.8) \qquad \begin{aligned} & X \circ (F_{jk}(W_j[f, X, a_j])(W_k[f, X, a_k])) \\ &= f_{jk}(w_j[f, X, a_j])(w_k[f, X, a_k]) \\ &= f_{jk}(X \circ W_j(f, X, a_j))(X \circ W_k(f, X, a_k)) \end{aligned}$$

Из равенства (8.1.8) следует корректность определения (8.1.6) представления $\Omega_j$-алгебры $W_j[f, X]$ в $\Omega_k$-алгебре $W_k[f, X]$. $\qquad \square$

Следствие 8.1.11. *Кортеж $\Omega$-алгебр*

$$W[f, X] = (W_{(1)}[f, X], ..., W_{(n)}[f, X])$$

*и множество представлений $F$ порождает диаграмму представлений $(F, W[f, X])$.* $\qquad \square$

Теорема 8.1.12. *Пусть $(f, A)$, $(g, B)$ - диаграммы представлений. Для заданных множеств $X_{(k)} \subset A_{(k)}$, $Y_{(k)} \subset B_{(k)}$, $(k) = (1), ..., (n)$, рассмотрим кортеж отображений*

$$R = (R_{(1)}, ..., R_{(n)})$$

*таких, что для любого $(k)$, $(k) = (1), ..., (n)$, отображение*

$$R_{(k)} : X_{(k)} \to Y_{(k)}$$

*согласовано со структурой диаграммы представлений, т. е.*

(8.1.9)
$$\begin{cases} \omega \in \Omega_{(k)}(p), \ x_{(k)1}, \ ..., \ x_{(k)p}, \ x_{(k)1}...x_{(k)p}\omega \in X_{(k)}, \\ R_{(k)}(x_{(k)1}...x_{(k)p}\omega) \in Y_{(k)} \\ => R_{(k)}(x_{(k)1}...x_{(k)p}\omega) = R_{(k)}(x_{(k)1})...R_{(k)}(x_{(k)p})\omega \end{cases}$$

(8.1.10)
$$\begin{cases} a_j \in X_j, \ a_k \in X_k, \ R_k(f_{jk}(a_j)(a_k)) \in Y_k \\ => R_k(f_{jk}(a_j)(a_k)) = g_{jk}(R_j(a_j))(R_k(a_k)) \end{cases}$$

*Рассмотрим кортеж отображений $\Omega$-слов*

$$w_{(k)}[f \to g, \overline{\overline{e}}, Y, R] : w_{(k)}[f, \overline{\overline{e}}] \to w_{(k)}[g, Y]$$

*удовлетворяющее условиям 8.1.7.1, 8.1.7.2, 8.1.7.3, и такое, что*

$$e_{(k)i} \in \overline{\overline{e}}_{(k)} => w_{(k)}[f \to g, \overline{\overline{e}}, Y, R](e_{(k)i}) = R_{(k)}(e_{(k)i})$$

*Для каждого $(k)$, $(k)$, $(k) = (1), ..., (n)$, существует гомоморфизм $\Omega_{(k)}$-алгебры*

$$r_{(k)} : A_{(k)} \to B_{(k)}$$

*определённый правилом*

(8.1.11)        $$r_{(k)}(a_{(k)}) = w_{(k)}[f \to g, X, Y, R](w_{(k)}[f, X, a_{(k)}])$$

*Кортеж гомоморфизмов*

$$r = (r_{(1)} \ ... \ r_{(n)}) = (r_1 \ ... \ r_n)$$

*является морфизмом диаграмм представлений $J[f, X]$ и $J[g, Y]$.*

Доказательство. Для любого $(k)$, $(k) = (1), ..., (n)$, рассмотрим отображение

$$r_{(k)} : A_{(k)} \to B_{(k)}$$

Лемма 8.1.13. *Для любого $(k)$, $(k) = (1), ..., (n)$, на множестве $X_{(k)}$ отображения $r_{(k)}$ и $R_{(k)}$ совпадают, и отображение $r_{(k)}$ согласовано со структурой $\Omega_{(k)}$-алгебры.*

Доказательство. Если

(8.1.12)        $$w_{(k)}[f, X, a_{(k)}] = a_{(k)}$$

то $a_{(k)} \in X_{(k)}$. Согласно условию 8.1.7.1, равенство

(8.1.13)
$$r_{(k)}(a_{(k)}) = w_{(k)}[f \to g, X, Y, R](w_{(k)}[f, X, a_{(k)}]) = w_{(k)}[f \to g, X, Y, R](a_{(k)})$$
$$= R_{(k)}(a_{(k)})$$

является следствием равенств (8.1.11), (8.1.12). Лемма является следствием равенства (8.1.13).        $\odot$

Лемма 8.1.14. *Пусть $\omega \in \Omega_{(k)}(p)$.*

(8.1.14)        $$r_{(k)}(x_{(k)1}...x_{(k)p}\omega) = r_{(k)}(x_{(k)1})...r_{(k)}(x_{(k)p})\omega$$

Доказательство. Мы будем доказывать лемму индукцией по сложности $\Omega_{(k)}$-слова.

Если

$$x_{(k)1},\ ...,\ x_{(k)p},\ x_{(k)1}...x_{(k)p}\omega \in X_{(k)}$$

то равенство (8.1.14) является следствием утверждения (8.1.9).

Пусть предположение индукции верно для

$$a_{(k)1}, ..., a_{(k)p} \in J_{(k)}[f, X]$$

Пусть

(8.1.15)     $$w_{(k)1} = w_{(k)}[f, X, a_{(k)1}] \quad ... \quad w_{(k)p} = w_{(k)}[f, X, a_{(k)p}]$$

Согласно предположению индукции, равенство

$$r_{(k)}(a_{(k)1}) = w_{(k)}[f \to g, X, Y, R](w_{(k)1})$$

(8.1.16)     $$... = ...$$

$$r_{(k)}(a_{(k)p}) = w_{(k)}[f \to g, X, Y, R](w_{(k)p})$$

является следствием равенств (8.1.11), (8.1.15). Если

(8.1.17)     $$a_{(k)} = a_{(k)1}...a_{(k)p}\omega$$

то согласно условию 8.1.5.2,

$$w_{(k)}[f, X, a_{(k)}] = w_{(k)1}...w_{(k)p}\omega$$

Согласно условию 8.1.7.2, равенство

$$r_{(k)}(a_{(k)}) = w_{(k)}[f \to g, X, Y, R](w_{(k)}[f, X, a_{(k)}])$$

(8.1.18)     $$= w_{(k)}[f \to g, X, Y, R](w_{(k)1}...w_{(k)p}\omega)$$

$$= w_{(k)}[f \to g, X, Y, R](w_{(k)1})...w_{(k)}[f \to g, X, Y, R](w_{(k)p})\omega$$

$$= (r_{(k)}(a_{(k)1}))...(r_{(k)}(a_{(k)p)})\omega$$

является следствием равенств (8.1.11), (8.1.17), (8.1.16). Равенство (8.1.14) является следствием равенства (8.1.18).                    ⊙

Согласно лемме 8.1.13, отображения $r_{(k)}$ и $R_{(k)}$ совпадают на множестве $X_{(k)}$. Согласно лемме 8.1.14, отображение $r_{(k)}$ является гомоморфизмом $\Omega_{(k)}$-алгебры $A_{(k)}$ в $\Omega_{(k)}$-алгебру $B_{(k)}$. Для доказательства теоремы достаточно показать, что если существует представление

$$f_{ji} : A_j \longrightarrow^* A_i$$

то пара отображений $(r_j \quad r_i)$ является морфизмом представлений из $f_{ji}$ в $g_{ji}$ (определение 7.2.1).

Мы будем доказывать теорему индукцией по сложности $\Omega_i$-слова.

Если $a_i \in X_i, a_j \in X_j$, то предположение индукции является следствием утверждения (8.1.10).

Пусть предположение индукции верно для

$$a_j \in J_j[f, X] \quad w_j[f, X, a_j] = m_j$$

$$a_i \in J_i[f, X] \quad w_i[f, X, a_i] = m_i$$

Согласно условию 8.1.5.3,

(8.1.19)     $$w_i(f, X, f_{ji}(a_j)(a_i)) = f_{ji}(m_j)(m_i)$$

Согласно условию [8.1.7.3], равенство

$$
\begin{aligned}
r_i(f_{ji}(a_j)(a_i)) &= w_i[f \to g, X, Y, R](w_i[f, X, f_{ji}(a_j)(a_i)]) \\
&= w_i[f \to g, X, Y, R](f_{ji}(m_j)(m_i)) \\
&= g_{ji}(w_j[g, Y, r_j(a_j)])(w_i[g, Y, r_i(a_i)]) \\
&= g_{ji}(r_j(a_j))(r_i(a_i))
\end{aligned}
$$

(8.1.20)

является следствием равенств [(8.1.11)], [(8.1.19)], Из равенств [(7.2.3)], [(8.1.20)] следует, что отображение $r$ является морфизмом диаграммы представлений $(f, A)$.          $\square$

ЗАМЕЧАНИЕ 8.1.15. *Теорема [8.1.12] - это теорема о продолжении отображения. Единственное, что нам известно о кортеже множеств $X$ - это то, что $X$ - кортеж множеств образующих диаграммы представлений $(f, A)$. Однако, между элементами множества $X_{(k)}$, $(k) = (1), ..., (n)$, могут существовать соотношения, порождённые либо операциями $\Omega_{(k)}$-алгебры $A_{(k)}$, либо преобразованиями представления $f_{jk}$. Поэтому произвольное отображение кортежа множеств $X$, вообще говоря, не может быть продолжено до эндоморфизма диаграммы представлений $(f, A)$.[8.3] Однако, если для каждого $(k)$, $(k) = (1), ..., (n)$, отображение $R_{(k)}$ согласованно со структурой диаграммы представлений, то мы можем построить продолжение этого отображения, которое является морфизмом диаграммы представлений $(f, A)$.*          $\square$

ОПРЕДЕЛЕНИЕ 8.1.16. *Пусть $X$ - кортеж множеств образующих диаграммы представлений $(f, A)$. Пусть $Y$ - кортеж множеств образующих диаграммы представлений $(g, B)$. Пусть $r$ - морфизм диаграммы представлений $(f, A)$ в диаграмму представлений $(g, B)$. Множество координат $W[g, Y, r(X)]$ называется **координатами морфизма диаграммы представлений**.*          $\square$

ОПРЕДЕЛЕНИЕ 8.1.17. *Пусть $X$ - кортеж множеств образующих диаграммы представлений $(f, A)$. Пусть $Y$ - кортеж множеств образующих диаграммы представлений $(g, B)$. Пусть $r$ - морфизм диаграммы представлений $(f, A)$ в диаграмму представлений $(g, B)$. Пусть для $(k) = (1), ..., (n)$, $a_{(k)} \in A_{(k)}$. Мы определим **суперпозицию координат** морфизма $r$ диаграммы представлений и $A_{(k)}$-числа $a_{(k)}$ как координаты, определённые согласно правилу*

(8.1.21)     $W_{(k)}[g, Y, r_{(k)}(X_{(k)})] \circ W_{(k)}[f, X, a_{(k)}] = W_{(k)}[g, Y, r_{(k)}(a_{(k)})]$

*Пусть $Y_{(k)} \subset A_{(k)}$. Мы определим суперпозицию координат морфизма $r$ диаграммы представлений и множества $Y_{(k)}$ согласно правилу*

(8.1.22)
$$
\begin{aligned}
&W_{(k)}[g, Y, r_{(k)}(X_{(k)})] \circ W_{(k)}[f, X, Y_{(k)}] \\
&= (W_{(k)}[g, Y, r_{(k)}(X_{(k)})] \circ W_{(k)}[f, X, a_{(k)}], a_{(k)} \in Y_{(k)})
\end{aligned}
$$

          $\square$

Теорема 8.1.18. *Морфизм $r$ диаграммы представлений $(f, A)$ в диаграмму представлений $(g, B)$ порождает отображение координат диаграммы представлений*

$$(8.1.23) \qquad W_{(k)}[f \to g, X, Y, r] : W_{(k)}[f, X] \to W_{(k)}[g, Y]$$

$(k) = (1), ..., (n),$ *такое, что*

$$(8.1.24) \qquad \begin{aligned} W_{(k)}[f, X, a] &\to W_{(k)}[f \to g, X, Y, r] \circ W_{(k)}[f, X, a_{(k)}] \\ &= W_{(k)}[g, Y, r_{(k)}(a_{(k)})] \end{aligned}$$

$$(8.1.25) \qquad \begin{aligned} W_{(k)}[f &\to g, X, Y, r] \circ W_{(k)}[f, X, a_{(k)}] \\ &= W_{(k)}[g, Y, r_{(k)}(X_{(k)})] \circ W_{(k)}[f, X, a_{(k)}] \end{aligned}$$

Доказательство. Согласно замечанию 8.1.8, мы можем рассматривать равенства (8.1.21), (8.1.23) относительно заданного кортежа множеств образующих $X$. При этом координатам $W_{(k)}[f, X, a_{(k)}]$ соответствует кортеж слов

$$(8.1.26) \qquad X \circ W_{(k)}[f, X, a_{(k)}] = w_{(k)}[f, X, a_{(k)}]$$

а координатам $W_{(k)}[g, Y, r_{(k)}(a_{(k)})]$ соответствует кортеж слов

$$(8.1.27) \qquad Y \circ W_{(k)}[g, Y, r_{(k)}(a_{(k)})] = w_{(k)}[g, Y, r_{(k)}(a_{(k)})]$$

Поэтому для того, чтобы доказать теорему, нам достаточно показать, что отображению $W_{(k)}[f \to g, X, Y, r]$ соответствует отображение $w_{(k)}[f \to g, X, Y, r]$.

Мы будем доказывать теорему индукцией по сложности $\Omega_{(k)}$-слова.

Если $a_{(k)} \in X_{(k)}$, $a'_{(k)} = r_{(k)}(a_{(k)})$, то, согласно равенствам (8.1.26), (8.1.27), отображения $W_{(k)}[f \to g, X, Y, r]$ и $w_{(k)}[f \to g, X, Y, r]$ согласованы.

Пусть для $a_{(k)1}, ..., a_{(k)p} \in X_{(k)}$ отображения $W_{(k)}[f \to g, X, Y, r]$ и $w_{(k)}[f \to g, X, Y, r]$ согласованы. Пусть $\omega \in \Omega_{(k)}(p)$. Согласно теореме 6.1.12

$$(8.1.28) \qquad W_{(k)}[f, X, a_{(k)1}...a_{(k)p}\omega] = W_{(k)}[f, X, a_{(k)1}]...W_{(k)}[f, X, a_{(k)p}]\omega$$

Так как отображение

$$r_{(k)} : A_{(k)} \to B_{(k)}$$

является гомоморфизмом $\Omega_{(k)}$-алгебры, то из равенства (8.1.28) следует

$$(8.1.29) \qquad \begin{aligned} W_{(k)}[g, &Y, r_{(k)}(a_{(k)1}...a_{(k)p}\omega)] \\ &= W_{(k)}[g, Y, (r(a_{(k)1})...(r(a_{(k)p}))\omega] \\ &= W_{(k)}[g, Y, r_{(k)}(a_{(k)1})]...W_{(k)}[g, Y, r_{(k)}(a_{(k)p})]\omega \end{aligned}$$

Из равенств (8.1.28), (8.1.29) и предположения индукции следует, что отображения $W_{(k)}[f \to g, X, Y, r]$ и $w_{(k)}[f \to g, X, Y, r]$ согласованы для $a_{(k)} = a_{(k)1}...a_{(k)p}\omega$.

Пусть для $a_{j1} \in A_j$ отображения $W_j[f \to g, X, Y, r]$ и $w_j[f \to g, X, Y, r]$ согласованы. Пусть для $a_{i1} \in A_i$ отображения $W_i[f \to g, X, Y, r]$ и $w_i[f \to g, X, Y, r]$ согласованы. Согласно теореме 8.1.10

$$(8.1.30) \qquad W_i(f, X, f_{ji}(a_j)(a_i)) = F_{ji}(W_j(f, X, a_j))(W_i(f, X, a_i))$$

Так как отображение $(r_j, r_i)$ является морфизмом представления $f_{ji}$ в представление $F_{ji}$, то из равенства (8.1.30) следует

$$(8.1.31) \qquad \begin{aligned} W_i[g, &Y, r_i(f_{ji}(a_j)(a_i))] = W_i[g, Y, g_{ji}(r_j(a_j))(r_i(a_i))] \\ &= G_{ji}(W_j[g, Y, r_j(a_j)])(W_i[g, Y, r_i(a_{n,1})]) \end{aligned}$$

Из равенств (8.1.30), (8.1.31) и предположения индукции следует, что отображения $W_i[f \rightarrow g, X, Y, r]$ и $w_i[f \rightarrow g, X, Y, r]$ согласованы для $b_i = f_{ji}(a_j)(a_i)$. □

Следствие 8.1.19. *Пусть $X$ - кортеж множеств образующих диаграммы представлений $(f, A)$. Пусть $Y$ - кортеж множеств образующих диаграммы представлений $(g, B)$. Пусть $r$ - морфизм диаграммы представлений $(f, A)$ в диаграмму представлений $(g, B)$. Отображение*

$$W[f \rightarrow g, X, Y, r] = (W_{(1)}[f \rightarrow g, X, Y, r], ..., W_{(n)}[f \rightarrow g, X, Y, r])$$

*является морфизмом диаграммы представлений $(F, W[f, X])$ в диаграмму представлений $(G, W[g, Y])$.* □

В дальнейшем мы будем отождествлять отображение $W[f \rightarrow g, X, Y, r]$ и множество координат $W[g, Y, r(X)]$ .

Теорема 8.1.20. *Пусть $X$ - кортеж множеств образующих диаграммы представлений $(f, A)$. Пусть $Y$ - кортеж множеств образующих диаграммы представлений $(g, B)$. Пусть $r$ - морфизм диаграммы представлений $(f, A)$ в диаграмму представлений $(g, B)$. Пусть $Y \subset A$. Тогда*

$$(8.1.32) \qquad W[g, Y, r(X)] \circ W[f, X, X'] = W[g, Y, r(X')]$$

$$(8.1.33) \qquad W[f \rightarrow g, X, Y, r] \circ W[f, X, X'] = W[g, Y, r(X')]$$

Доказательство. Равенство (8.1.32) является следствием равенства

$$r(X') = (r(a), a \in X')$$

а также равенств (8.1.21), (8.1.22). Равенство (8.1.33) является следствием равенств (8.1.32), (8.1.24). □

Теорема 8.1.21. *Пусть $X$ - кортеж множеств образующих диаграммы представлений $(f, A)$. Пусть $Y$ - кортеж множеств образующих диаграммы представлений $(g, B)$. Пусть $Z$ - кортеж множеств образующих диаграммы представлений $(h, C)$. Пусть $r$ - морфизм диаграммы представлений $(f, A)$ в диаграмму представлений $(g, B)$. Пусть $s$ - морфизм диаграммы представлений $(g, B)$ в диаграмму представлений $(h, C)$. Тогда*

$$(8.1.34) \qquad W[h, Z, s(Y)] \circ W[g, Y, r(X)] = W[h, Z, (s \circ r)(X)]$$

$$(8.1.35) \qquad W[g \rightarrow h, Y, Z, s] \circ W[f \rightarrow g, X, Y, r] = W[f \rightarrow h, X, Z, s \circ r]$$

Доказательство. Равенство

$$(8.1.36) \qquad W[h, Z, s(Y')] \circ W[g, Y, Y'] = W[h, z, s(Y')]$$

следует из равенства (8.1.32). Равенство (8.1.34) следует из равенства (8.1.36), если положить $Y' = r(X)$. Равенство (8.1.35) следует из равенства (8.1.34). □

Определение 8.1.22. *Мы можем обобщить определение суперпозиции координат и предположить, что один из множителей является кортежем множеств $\Omega$-слов. Соответственно, определение суперпозиции координат имеет вид*

$$W[g, Y, r(X)] \circ w[f, X, X'] = w[g, Y, r(X)] \circ W[f, X, X'] = w[g, Y, r(X')]$$

□

Следующие формы записи образа кортежа множеств $X'$ при морфизме $r$ диаграммы представлений эквивалентны

$$(8.1.37) \qquad r(X') = r(X) \circ W[f, X, X']$$
$$= (Y \circ W[g, Y, r(X)]) \circ W[f, X, X']$$

Из равенств (8.1.32), (8.1.37) следует, что

$$(8.1.38) \qquad Y \circ (W[g, Y, r(X) \circ W[f, X, X']])$$
$$= (Y \circ W[g, Y, r(X)]) \circ W[f, X, X']$$

Равенство (8.1.38) является законом ассоциативности для операции композиции и позволяет записать выражение

$$Y \circ W[g, Y, r(X) \circ W[f, X, X']$$

без использования скобок.

ОПРЕДЕЛЕНИЕ 8.1.23. *Пусть*

$$X = (X_{(1)} \subset A_{(1)}, ..., X_{(n)} \subset A_{(n)}) = (X_1 \subset A_1, ..., X_n \subset A_n)$$

*множество образующих диаграммы представлений* $(f, A)$. *Пусть отображение* $r$ *является эндоморфизмом диаграммы представлений* $(f, A)$. *Пусть кортеж множеств* $Y = r(X)$ *является образом кортежа множеств* $X$ *при отображении* $r$. *Эндоморфизм* $r$ *диаграммы представлений* $(f, A)$ *называется невырожденным на кортеже множеств образующих* $X$, *если кортеж множеств* $Y$ *является кортежем множеств образующих диаграммы представлений* $(f, A)$. *В противном случае, эндоморфизм* $r$ *называется вырожденным на кортеже множеств образующих* $X$. $\square$

ОПРЕДЕЛЕНИЕ 8.1.24. *Эндоморфизм* $r$ *диаграммы представлений* $(f, A)$ *называется* **невырожденным**, *если он невырожден на любом кортеже множеств образующих. В противном случае, эндоморфизм* $r$ *называется* **вырожденным**. $\square$

ТЕОРЕМА 8.1.25. *Автоморфизм* $r$ *диаграммы представлений* $(f, A)$ *является невырожденным эндоморфизмом.*

ДОКАЗАТЕЛЬСТВО. Пусть $X$ - кортеж множеств образующих диаграммы представлений $(f, A)$. Пусть $Y = r(X)$. Согласно теореме 8.1.18 эндоморфизм $r$ порождает отображение $\Omega$-слов $w[f \to g, X, Y, r]$. Пусть $a' \in A$. Так как $r$ - автоморфизм, то существует $a \in A$, $r(a) = a'$. Согласно определению 8.1.5, $w[f, X, a]$ - кортеж $\Omega$-слов, представляющих $a$ относительно кортежа множеств образующих $X$. Согласно теореме 8.1.18, $w[f, X', a']$ - кортеж $\Omega$-слов, представляющих $a'$ относительно кортежа множеств $Y$

$$w[f, Y, a'] = w[f \to g, X, Y, r](w[f, X, a])$$

Следовательно, $Y$ - множество образующих диаграммы представлений $(f, A)$. Согласно определению 8.1.24, автоморфизм $r$ - невырожден. $\square$

## 8.2. Базис диаграммы представлений

Определение 8.2.1. *Пусть* $(f, A)$ *- диаграмма представлений и*

$$Gen[f, A] = \{X = (X_{(1)}, ..., X_{(n)}) : X_{(k)} \subseteq A_{(k)}, J_{(k)}[f, X] = A_{(k)}\}$$

*Если для кортежа множеств* $X \subset A_2$ *верно* $X \in Gen[f, A]$, *то для любого кортежа множеств* $Y$, $X_{(k)} \subset Y_{(k)} \subset A_{(k)}$, $(k) = (1), ..., (n)$ *m акже верно* $Y \in Gen[f, A]$. *Если существует минимальный кортеж множеств* $X \in Gen[f, A]$, *то такой кортеж множеств* $X$ *называется* **квазибазисом диаграммы представлений** $(f, A)$. $\qquad\square$

Теорема 8.2.2. *Если кортеж множеств* $X$ *является квазибазисом диаграммы представлений* $(f, A)$, *то, для любого* $(k)$, $(k) = (1), ..., (n)$, *и любого* $m \in X_{(k)}$, *кортеж множеств*

$$X' = (X_{(1)}, ..., X'_{(k)} = X_{(k)} \setminus \{m\}, ..., X_{(n)})$$

*не является множеством образующих диаграммы представлений* $(f, A)$.

Доказательство. Пусть $X$ - квазибазис диаграммы представлений $(f, A)$. Допустим для некоторого $m \in X_{(k)}$ существует $\Omega_{(k)}$-слово

$$w = w_{(k)}[f, X', m]$$

Рассмотрим $A_{(k)}$-число $m'$, для которого $\Omega_{(k)}$-слово $w' = w_{(k)}[f, X, m']$ зависит от $m$. Согласно определению 8.1.5, любое вхождение $A_{(k)}$-числа $m$ в $\Omega_{(k)}$-слово $w'$ может быть заменено $\Omega_{(k)}$-словом $w$. Следовательно, $\Omega_{(k)}$-слово $w'$ не зависит от $m$, а кортеж множеств $X'$ является множеством образующих диаграммы представлений $(f, A)$. Следовательно, $X$ не является квазибазисом диаграммы представлений $(f, A)$. $\qquad\square$

Замечание 8.2.3. *Доказательство теоремы 8.2.2 даёт нам эффективный метод построения квазибазиса диаграммы представлений* $(f, A)$. *Квазибазис диаграммы представлений определён индукцией по диаграмме представлений. Мы начинаем строить квазибазис в* $\Omega$*-алгебрах из множества* $A_{[0]}$. *Когда квазибазис построен в* $\Omega$*-алгебрах из множества* $A_{[i]}$, *мы можем перейти к построению квазибазиса в* $\Omega$*-алгебрах из множества* $A_{[i+1]}$. $\qquad\square$

Для каждого $(k)$, $(k) = (1), ..., (n)$, мы ввели $\Omega_{(k)}$-слово $A_{(k)}$-числа $x$ относительно множества образующих $X$ в определении [8.4] 8.1.5. Из теоремы 8.2.2 следует, что если множество образующих $X$ не является квазибазисом, то выбор $\Omega_{(k)}$-слова относительно множества образующих $X$ неоднозначен. Но даже если множество образующих $X$ является квазибазисом, то представление $m \in A_{(k)}$ в виде $\Omega_{(k)}$-слова неоднозначно.

Замечание 8.2.4. *Существует три источника неоднозначности в записи* $\Omega_{(k)}$*-слова.*

8.2.4.1: *В* $\Omega_{(k)}$*-алгебре* $A_{(k)}$, $(k) = (1), ..., (n)$, *могут быть определены равенства. Например, если е - единица мультипликативной группы* $A_{(k)}$, *то верно равенство*

$$ae = a$$

*для любого* $a \in A_{(k)}$.

---

[8.4] Рассуждения в начале этого раздела естественно повторяют рассуждения в начале раздела 6.2 и я сохранил эти рассуждения для полноты текста.

8.2.4.2: *Неоднозначность выбора $\Omega_{(k)}$-слова может быть связана со свойствами представления. Например, допустим существует представление $f_{ik}$ $\Omega_i$-алгебры $A_i$ в $\Omega_k$-алгебре $A_k$. Если $m_1$, ..., $m_n$ - $\Omega_k$-слова, $\omega \in \Omega_k(n)$ и $a$ - $\Omega_i$-слово, то[8.5]*

$$(8.2.1) \qquad f_{ik}(a)(m_1...m_n\omega) = (f_{ik}(a)(m_1))...(f_{ik}(a)(m_n))\omega$$

*В тоже время, если $\omega$ является операцией $\Omega_i$-алгебры $A_i$ и операцией $\Omega_k$-алгебры $A_k$, то мы можем потребовать, что $\Omega_k$-слова $f(a_1...a_n\omega)(x)$ и $(f(a_1)(x))...(f(a_n)(x))\omega$ описывают один и тот же элемент $\Omega_k$-алгебры $A_k$.[8.6]*

$$(8.2.2) \qquad f(a_1...a_n\omega)(x) = (f(a_1)(x))...(f(a_n)(x))\omega$$

8.2.4.3: *Равенства вида (8.2.1), (8.2.2) сохраняются при морфизме диаграммы представлений. Поэтому мы можем игнорировать эту форму неоднозначности записи $\Omega_{(k)}$-слова. Однако возможна принципиально другая форма неоднозначности, пример которой можно найти в теоремах 9.3.15, 9.3.16.*

*Таким образом, мы видим, что на множестве $\Omega_{(k)}$-слов можно определить различные отношения эквивалентности.[8.7] Наша задача - найти максимальное отношение эквивалентности на множестве $\Omega_{(k)}$-слов, которое сохраняется при морфизме представления.*

*Аналогичное замечание касается отображения $W[f, X, m]$, определённого в замечании 8.1.8.[8.8]* □

---

[8.5] Например, пусть $\{e_1, e_2\}$ - базис векторного пространства над полем $k$. Равенство (8.2.1) принимает форму закона дистрибутивности

$$a(b^1 e_1 + b^2 e_2) = (ab^1)e_1 + (ab^2)e_2$$

[8.6] Для векторного пространства это требование принимает форму закона дистрибутивности

$$(a + b)e_1 = ae_1 + be_1$$

[8.7] Очевидно, что каждое из равенств (8.2.1), (8.2.2) порождает некоторое отношение эквивалентности.

[8.8] Если базис векторного пространства - конечен, то мы можем представить базис в виде матрицы строки

$$\overline{\overline{e}} = \begin{pmatrix} e_1 & ... & e_2 \end{pmatrix}$$

Мы можем представить отображение $W[f, \overline{\overline{e}}](v)$ в виде матрицы столбца

$$W[f, \overline{\overline{e}}, v] = \begin{pmatrix} v^1 \\ ... \\ v^n \end{pmatrix}$$

Тогда

$$W[f, \overline{\overline{e}}, v](\overline{\overline{e}}') = W[f, \overline{\overline{e}}, v] \begin{pmatrix} e'_1 & ... & e'_n \end{pmatrix} = \begin{pmatrix} v^1 \\ ... \\ v^n \end{pmatrix} \begin{pmatrix} e'_1 & ... & e'_n \end{pmatrix}$$

имеет вид произведения матриц.

ТЕОРЕМА 8.2.5. *Пусть $X$ - квазибазис диаграммы представлений $(f, A)$. Рассмотрим кортеж отношений эквивалентности*

$$\lambda[f, X] = (\lambda_{(1)}[f, X], ..., \lambda_{(n)}[f, X])$$

$$\lambda_{(k)}[f, X] \subseteq w_{(k)}[f, X] \times w_{(k)}[f, X]$$

*которое порождено исключительно следующими утверждениями.*

8.2.5.1: *Если в $\Omega_{(k)}$-алгебре $A_{(k)}$ существует равенство*

$$w_{(k)1}[f, X, m] = w_{(k)2}[f, X, m]$$

*определяющее структуру $\Omega_{(k)}$-алгебры, то*

$$(w_{(k)1}[f, X, m], w_{(k)2}[f, X, m]) \in \lambda_{(k)}[f, X]$$

8.2.5.2: *Если существует представление $f_{ik}$ и в $\Omega_i$-алгебре $A_i$ существует равенство*

$$w_{i1}[f, X, m] = w_{i2}[f, X, m]$$

*определяющее структуру $\Omega_i$-алгебры, то*

$$(f_{ik}(w_{i1})(w_k[f, X, m]), f_{ik}(w_{i2})(w_k[f, X, m])) \in \lambda_k[f, X]$$

8.2.5.3: *Если существует представление $f_{ik}$, то для любой операции $\omega \in \Omega_i(n)$,*

$$(f_{ik}(a_{i1}...a_{in}\omega)(a_2), (f_{ik}(a_{i1})...f_{ik}(a_{in})\omega)(a_2)) \in \lambda_k[f, X]$$

8.2.5.4: *Если существует представление $f_{ik}$, то для любой операции $\omega \in \Omega_k(n)$,*

$$(f_{ik}(a_i)(a_{k1}...a_{kn}\omega), f_{ik}(a_i)(a_{k1})...f_{ik}(a_i)(a_{kn})\omega) \in \lambda_k[f, X]$$

8.2.5.5: *Если существует представление $f_{ik}$, $\omega \in \Omega_i(n) \cap \Omega_k(n)$ и представление $f_{ik}$ удовлетворяет равенству*[8.9]

$$f(a_{i1}...a_{in}\omega)(a_k) = (f(a_{i1})(a_k))...(f(a_{in})(a_k))\omega$$

*то мы можем предположить, что верно равенство*

$$(f(a_{i1}...a_{in}\omega)(a_k), (f(a_{i1})(a_k))...(f(a_{in})(a_k))\omega) \in \lambda_k[f, X]$$

ДОКАЗАТЕЛЬСТВО. Теорема верна, так как рассмотренные равенства сохраняются при гомоморфизмах универсальных алгебр $A_{(k)}$. ☐

ОПРЕДЕЛЕНИЕ 8.2.6. *Квазибазис $\overline{\overline{e}}$ диаграммы представлений $(f, A)$ такой, что*

$$\rho[f, \overline{\overline{e}}] = \lambda[f, \overline{\overline{e}}]$$

*называется* **базисом диаграммы представлений** $(f, A)$. ☐

---

[8.9] Рассмотрим представление коммутативного кольца $D$ в $D$-алгебре $A$. Мы будем пользоваться записью

$$f(a)(v) = av$$

В обеих алгебрах определены операции сложения и умножения. Однако равенство

$$f(a + b)(v) = f(a)(v) + f(b)(v)$$

верно, а равенство

$$f(ab)(v) = f(a)(v)f(b)(v)$$

является ошибочным.

Замечание 8.2.7. *Мы будем записывать базис также в виде*

$$\overline{\overline{e}} = (\overline{\overline{e}}_{(1)}, ..., \overline{\overline{e}}_{(n)})$$

$$\overline{\overline{e}}_{(k)} = (e_{(k)l}, e_{(k)l} \in \overline{\overline{e}}_{(k)}) \quad (k) = (1), ..., (n)$$

*Если базис - конечный, то мы будем также пользоваться записью*

$$\overline{\overline{e}}_{(k)} = (e_{(k)i}, i \in I_{(k)}) = (e_{(k)1}, ..., e_{(k)p_{(k)}}) \quad (k) = (1), ..., (n)$$

<div align="right">□</div>

Теорема 8.2.8. *Автоморфизм диаграммы представлений* $(f, A)$ *отображает базис диаграммы представлений* $(f, A)$ *в базис.*

Доказательство. Пусть отображение $r$ - автоморфизм диаграммы представлений $(f, A)$. Пусть кортеж множеств $\overline{\overline{e}}$ - базис диаграммы представлений $(f, A)$. Пусть [8.10] $\overline{\overline{e}}' = r \circ \overline{\overline{e}}$. Допустим кортеж множеств $\overline{\overline{e}}'$ не является базисом. Согласно теореме 8.2.2 существуют $(k)$, $(k) = (1), ..., (n)$, и $e'_{(k)i} \in \overline{\overline{e}}'_{(k)}$ такие, что кортеж множества

$$Z = (\overline{\overline{e}}'_{(1)}, ..., Z_{(k)} = \overline{\overline{e}}'_{(k)} \setminus \{e'_{(k)i}\}, ..., \overline{\overline{e}}'_{(n)})$$

является множеством образующих диаграммы представлений $(f, A)$. Согласно теореме 7.3.3 отображение $r^{-1}$ является автоморфизмом диаграммы представлений $(f, A)$. Согласно теореме 8.1.25 и определению 8.1.24, кортеж множеств

$$X = (\overline{\overline{e}}_{(1)}, ..., X_{(k)} = \overline{\overline{e}}_{(k)} \setminus \{r_{(k)}^{-1}(e'_{(k)i})\}, ..., \overline{\overline{e}}_{(n)})$$

является множеством образующих диаграммы представлений $(f, A)$. Полученное противоречие доказывает теорему. □

Теорема 8.2.9. *Пусть* $\overline{\overline{e}}$ - *базис диаграммы представлений* $(f, A)$. *Пусть* $(g, B)$ - *диаграмма представлений. Пусть*

$$R : \overline{\overline{e}} \to Y$$

*произвольное отображение кортежа множеств* $\overline{\overline{e}}$, $Y_{(k)} \subseteq B_{(k)}$, $(k) = (1)$, ..., $(n)$. *Рассмотрим кортеж отображений*

$$w_{(k)}[f \to g, \overline{\overline{e}}, Y, R] : w_{(k)}[f, \overline{\overline{e}}] \to w_{(k)}[g, Y]$$

*удовлетворяющих условиям 8.1.7.1, 8.1.7.2, 8.1.7.3, и такое, что*

$$e_{(k)i} \in \overline{\overline{e}}_{(k)} => w_{(k)}[f \to g, \overline{\overline{e}}, Y, R](e_{(k)i}) = R_{(k)}(e_{(k)i})$$

*Существует единственный морфизм диаграммы представлений* [8.11]

$$r : A \to B$$

*определённый правилом*

$$r(a) = w[f \to g, \overline{\overline{e}}, Y, R](w[f, \overline{\overline{e}}, a])$$

Доказательство. Утверждение теоремы является следствием теорем 6.1.10, 6.1.14. □

---

[8.10] Согласно определениям 5.1.3, 8.3.1, мы будем пользоваться записью $r(\overline{\overline{e}}) = r \circ \overline{\overline{e}}$.

[8.11] Это утверждение похоже на теорему [2]-1, с. 104.

Следствие 8.2.10. *Пусть $\overline{\overline{e}}$, $\overline{\overline{e}}'$ - базисы представления $(f, A)$. Пусть $r$ - автоморфизм представления $(f, A)$ такой, что $\overline{\overline{e}}' = r \circ \overline{\overline{e}}$. Автоморфизм $r$ определён однозначно.* □

Замечание 8.2.11. *Теорема 8.2.9, так же как и теорема 8.1.12, является теоремой о продолжении отображения. Одако здесь $\overline{\overline{e}}$ - не произвольное множество образующих диаграммы представлений, а базис. Согласно замечанию 8.2.3, мы не можем определить координаты любого элемента базиса через остальные элементы этого же базиса. Поэтому отпадает необходимость в согласованности отображения базиса с представлением.* □

Теорема 8.2.12. *Набор координат $W[f, \overline{\overline{e}}, \overline{\overline{e}}]$ соответствует тождественному преобразованию*

$$W[f, \overline{\overline{e}}, E] = W[f, \overline{\overline{e}}, \overline{\overline{e}}]$$

Доказательство. Утверждение теоремы следует из равенства

$$a = \overline{\overline{e}} \circ W[f, \overline{\overline{e}}, a] = \overline{\overline{e}} \circ W[f, \overline{\overline{e}}, \overline{\overline{e}}] \circ W[f, \overline{\overline{e}}, a]$$

□

Теорема 8.2.13. *Пусть $W[f, \overline{\overline{e}}, r \circ \overline{\overline{e}}]$ - множество координат автоморфизма $r$. Определено множество координат $W[f, r \circ \overline{\overline{e}}, \overline{\overline{e}}]$, соответствующее автоморфизму $r^{-1}$. Множество координат $W[f, r \circ \overline{\overline{e}}, \overline{\overline{e}}]$ удовлетворяет равенству*[8.12]

$$(8.2.3) \qquad W[f, \overline{\overline{e}}, r \circ \overline{\overline{e}}] \circ W[f, r \circ \overline{\overline{e}}, \overline{\overline{e}}] = W[f, \overline{\overline{e}}, \overline{\overline{e}}]$$

$$W[f \to f, \overline{\overline{e}}, \overline{\overline{e}}, r^{-1}] = W[f \to f, \overline{\overline{e}}, \overline{\overline{e}}, r]^{-1} = W[f, r \circ \overline{\overline{e}}, \overline{\overline{e}}]$$

Доказательство. Поскольку $r$ - автоморфизм диаграммы представлений $(f, A)$, то, согласно теореме 8.2.8, множество $r \circ \overline{\overline{e}}$ - базис диаграммы представлений $(f, A)$. Следовательно, определено множество координат $W[f, r \circ \overline{\overline{e}}, \overline{\overline{e}}]$. Равенство (8.2.3) следует из цепочки равенств

$$W[f, \overline{\overline{e}}, r \circ \overline{\overline{e}}] \circ W[f, r \circ \overline{\overline{e}}, \overline{\overline{e}}] = W[f, \overline{\overline{e}}, r \circ \overline{\overline{e}}] \circ W[f, \overline{\overline{e}}, r^{-1} \circ \overline{\overline{e}}]$$
$$= W[f, \overline{\overline{e}}, r \circ r^{-1} \circ \overline{\overline{e}}] = W[f, \overline{\overline{e}}, \overline{\overline{e}}]$$

□

Теорема 8.2.14. *Пусть $W[f, \overline{\overline{e}}, r \circ \overline{\overline{e}}]$ - множество координат автоморфизма $r$. Пусть $W[f, \overline{\overline{e}}, s \circ \overline{\overline{e}}]$- множество координат автоморфизма $s$. Множество координат автоморфизма $(r \circ s)^{-1}$ удовлетворяет равенству*

$$(8.2.4) \quad W[f, (r \circ s) \circ \overline{\overline{e}}, \overline{\overline{e}}] = W[f, s \circ (r \circ \overline{\overline{e}}), \overline{\overline{e}}] = W[f, s \circ \overline{\overline{e}}, \overline{\overline{e}}] \circ W[f, r \circ \overline{\overline{e}}, \overline{\overline{e}}]$$

Доказательство. Равенство

$$W[f, (r \circ s) \circ \overline{\overline{e}}, \overline{\overline{e}}] = W[f, \overline{\overline{e}}, (r \circ s)^{-1} \circ \overline{\overline{e}}] = W[f, \overline{\overline{e}}, s^{-1} \circ r^{-1} \circ \overline{\overline{e}}]$$
$$= W[f, \overline{\overline{e}}, s^{-1} \circ \ e] \circ W[f, \overline{\overline{e}}, r^{-1} \circ \overline{\overline{e}}]$$
$$= W[f, s \circ \overline{\overline{e}}, \overline{\overline{e}}] \circ W[f, r \circ \overline{\overline{e}}, \overline{\overline{e}}]$$
$$= W[f, s \circ (r \circ \overline{\overline{e}}), \overline{\overline{e}}]$$

$$(8.2.5)$$

является следствием теорем 8.1.21, 8.2.13. Равенство (8.2.4) является следствием равенства (8.2.5). □

---

[8.12] Смотри также замечание 6.2.15.

Теорема 8.2.15. *Группа автоморфизмов $GA(f)$ диаграммы эффективных представлений $(f, A)$ порождает эффективное левостороннее представление в диаграмме представлений $(f, A)$.*

Доказательство. Из следствия 8.2.10 следует, что если автоморфизм $r$ отображает базис $\overline{e}$ в базис $\overline{e}'$, то множество координат $W[f, \overline{e}, \overline{e}']$ однозначно определяет автоморфизм $r$. Из теоремы 8.1.18 следует, что множество координат $W[f, \overline{e}, \overline{e}']$ определяет правило отображения координат относительно базиса $\overline{e}$ при автоморфизме диаграммы представлений $(f, A)$. Из равенства (8.1.37) следует, что автоморфизм $r$ действует слева на элементы $\Omega_{(k)}$-алгебры $A_{(k)}$, $(k) = (1), ..., (n)$. Из равенства (8.1.34) следует, что представление группы является левосторонним представлением. Согласно теореме 8.2.12 набор координат $W[f, \overline{e}, \overline{e}]$ соответствует тождественному преобразованию. Из теоремы 8.2.13 следует, что набор координат $W[f, r \circ \overline{e}, \overline{e}]$ соответствует преобразованию, обратному преобразованию $W[f, \overline{e}, r \circ \overline{e}]$. $\square$

## 8.3. Многообразие базисов диаграммы представлений

Множество $\mathcal{B}[f]$ базисов диаграммы представлений $(f, A)$ называется **многообразием базисов** диаграммы представлений $(f, A)$.

Определение 8.3.1. *Согласно теореме 8.2.8 и определению 8.1.22, автоморфизм $r$ диаграммы представлений $(f, A)$ порождает преобразование*

$$
\begin{aligned}
(8.3.1) \qquad & r : \overline{\overline{h}} \to r \circ \overline{\overline{h}} \\
& r \circ \overline{\overline{h}} = W[f, \overline{e}, r \circ \overline{e}] \circ \overline{\overline{h}}
\end{aligned}
$$

*многообразия базисов диаграммы представлений. Это преобразование называется* **активным**. *Согласно теореме 7.3.3, определено левостороннее представление*

$$
A(f) : GA(f) \xrightarrow{\quad\quad} \mathcal{B}[f]
$$

*группы $GA(f)$ в многообразии базисов $\mathcal{B}[f]$. Представление $A(f)$ называется* **активным представлением**. *Согласно следствию 8.2.10, это представление однотранзитивно.* $\square$

Замечание 8.3.2. *Согласно замечанию 8.2.3, могут существовать базисы диаграммы представлений $(f, A)$, не связанные активным преобразованием. В этом случае мы в качестве многообразия базисов будем рассматривать орбиту выбранного базиса. Следовательно, диаграмма представлений $(f, A)$ может иметь различные многообразия базисов. Мы будем предполагать, что мы выбрали многообразие базисов.* $\square$

Теорема 8.3.3. *Существует однотранзитивное правостороннее представление*

$$
P(f) : GA(f) \xrightarrow{\quad*\quad} \mathcal{B}[f]
$$

*группы $GA(f)$ в многообразии базисов $\mathcal{B}[f]$. Представление $P(f)$ называется* **пассивным представлением**.

Доказательство. Поскольку $A(f)$ - однотранзитивное левостороннее представление группы $GA(f)$, то однотранзитивное правостороннее представление $P(f)$ определено однозначно согласно теореме 5.5.9. $\square$

Теорема 8.3.4. *Преобразование представления* $P(f)$ *называется* **пассивным преобразованием многообразия базисов** *диаграммы представлений. Мы будем пользоваться записью*

$$s(\overline{\overline{e}}) = \overline{\overline{e}} \circ s$$

*для обозначения образа базиса* $\overline{\overline{e}}$ *при пассивном преобразовании* $s$. *Пассивное преобразование базиса имеет вид*

$$(8.3.2) \qquad \begin{aligned} s : \overline{\overline{h}} \to \overline{\overline{h}} \circ s \\ \overline{\overline{h}} \circ s = \overline{\overline{h}} \circ W[f, \overline{\overline{e}}, \overline{\overline{e}} \circ s] \end{aligned}$$

Доказательство. Согласно равенству (8.3.1), активное преобразование действует на координаты базиса слева. Равенство (8.3.2) следует из теорем 5.5.8, 5.5.9, 5.5.11, согласно которым пассивное преобразование действует на координаты базиса справа. $\qquad \square$

Теорема 8.3.5. *Пассивное преобразование многообразия базисов является автоморфизмом представления* $A(f)$.

Доказательство. Теорема является следствием теоремы 5.5.11. $\qquad \square$

Теорема 8.3.6. *Пусть* $s$ - *пассивное преобразование многообразия базисов диаграммы представлений* $(f, A)$. *Пусть* $\overline{\overline{e}}_1$ - *базис диаграммы представлений* $(f, A)$, $\overline{\overline{e}}_2 = \overline{\overline{e}}_1 \circ s$. *Пусть для базиса* $\overline{\overline{e}}_3$ *существует активное преобразование* $r$ *такое, что* $\overline{\overline{e}}_3 = r \circ \overline{\overline{e}}_1$. *Положим* $\overline{\overline{e}}_4 = r \circ \overline{\overline{e}}_2$. *Тогда* $\overline{\overline{e}}_4 = \overline{\overline{e}}_3 \circ s$.

Доказательство. Согласно равенству (8.3.1), активное преобразование координат базиса $\overline{\overline{e}}_3$ имеет вид

$$(8.3.3) \qquad \overline{\overline{e}}_4 = W[f, \overline{\overline{e}}_1, \overline{\overline{e}}_3] \circ \overline{\overline{e}}_2 = W[f, \overline{\overline{e}}_1, \overline{\overline{e}}_3] \circ \overline{\overline{e}}_1 \circ W[f, \overline{\overline{e}}_1, \overline{\overline{e}}_2]$$

Пусть $\overline{\overline{e}}_5 = \overline{\overline{e}}_3 \circ s$. Из равенства (8.3.2) следует, что

$$(8.3.4) \qquad \overline{\overline{e}}_5 = \overline{\overline{e}}_3 \circ W[f, \overline{\overline{e}}_1, \overline{\overline{e}}_2] = W[f, \overline{\overline{e}}_1, \overline{\overline{e}}_3] \circ \overline{\overline{e}}_1 \circ W[f, \overline{\overline{e}}_1, \overline{\overline{e}}_2]$$

Из совпадения выражений в равенствах (8.3.3), (8.3.4) следует, что $\overline{\overline{e}}_4 = \overline{\overline{e}}_5$. Следовательно, коммутативна диаграмма

$$
\begin{array}{ccc}
\overline{\overline{e}}_1 \in \mathcal{B}[f] & \xrightarrow{\ r\ } & \overline{\overline{e}}_3 \in \mathcal{B}[f] \\
{\scriptstyle s}\downarrow & & \downarrow{\scriptstyle s} \\
\overline{\overline{e}}_2 \in \mathcal{B}[f] & \xrightarrow{\ r\ } & \overline{\overline{e}}_4 \in \mathcal{B}[f]
\end{array}
$$

$\qquad \square$

## 8.4. Геометрический объект диаграммы представлений

Активное преобразование изменяет базис диаграммы представлений и кортеж $\Omega$-чисел согласовано и координаты кортежа $\Omega$-чисел относительно базиса не меняются. Пассивное преобразование меняет только базис, и это ведёт к изменению координат кортежа $\Omega$-чисел относительно базиса.

Теорема 8.4.1. *Допустим пассивное преобразование* $s \in GA(f)$ *отображает базис* $\overline{\overline{e}}_1 \in \mathcal{B}[f]$ *в базис* $\overline{\overline{e}}_2 \in \mathcal{B}[f]$

$$(8.4.1) \qquad \overline{\overline{e}}_2 = \overline{\overline{e}}_1 \circ s = \overline{\overline{e}}_1 \circ W[f, \overline{\overline{e}}_1, \overline{\overline{e}}_1 \circ s]$$

*Допустим кортеж A-чисел a имеет кортеж Ω-слов*

$$(8.4.2) \qquad a = \overline{\overline{e}}_1 \circ W[f, \overline{\overline{e}}_1, a]$$

*относительно базиса* $\overline{\overline{e}}_1$ *и имеет кортеж Ω-слов*

$$(8.4.3) \qquad a = \overline{\overline{e}}_2 \circ W[f, \overline{\overline{e}}_2, a]$$

*относительно базиса* $\overline{\overline{e}}_2$. *Преобразование координат*

$$(8.4.4) \qquad W[f, \overline{\overline{e}}_2, a] = W[f, \overline{\overline{e}}_1 \circ s, \overline{\overline{e}}_1] \circ W[f, \overline{\overline{e}}_1, a]$$

*не зависит от кортежа A-чисел a или базиса* $\overline{\overline{e}}_1$, *а определенно исключительно координатами кортежа A-чисел a относительно базиса* $\overline{\overline{e}}_1$.

Доказательство. Из (8.4.1) и (8.4.3) следует, что

$$(8.4.5) \quad \begin{aligned} \overline{\overline{e}}_1 \circ W[f, \overline{\overline{e}}_1, a] &= \overline{\overline{e}}_2 \circ W[f, \overline{\overline{e}}_2, a] = \overline{\overline{e}}_1 \circ W[f, \overline{\overline{e}}_1, \overline{\overline{e}}_2] \circ W[f, \overline{\overline{e}}_2, a] \\ &= \overline{\overline{e}}_1 \circ W[f, \overline{\overline{e}}_1, \overline{\overline{e}}_1 \circ s] \circ W[f, \overline{\overline{e}}_2, a] \end{aligned}$$

Сравнивая (8.4.2) и (8.4.5) получаем, что

$$(8.4.6) \qquad W[f, \overline{\overline{e}}_1, a] = W[f, \overline{\overline{e}}_1, \overline{\overline{e}}_1 \circ s] \circ W[f, \overline{\overline{e}}_2, a]$$

Так как $s$ - автоморфизм представления, то равенство (8.4.4) следует из (8.4.6) и теоремы 8.2.13. $\qquad\square$

Теорема 8.4.2. *Преобразования координат* (8.4.4) *порождают эффективное контравариантное правостороннее представление группы* $GA(f)$, *называемое* **координатным представлением** *в кортеже Ω-алгебр.*

Доказательство. Согласно следствию 8.1.19, преобразование (8.4.4) является эндоморфизмом диаграммы представлений [8.13] $(F, W[f, \overline{\overline{e}}_1])$.

Допустим мы имеем два последовательных пассивных преобразования $s$ и $t$. Преобразование координат

$$(8.4.7) \qquad W[f, \overline{\overline{e}}_2, a] = W[f, \overline{\overline{e}}_1 \circ s, \overline{\overline{e}}_1] \circ W[f, \overline{\overline{e}}_1, a]$$

соответствует пассивному преобразованию $s$. Преобразование координат

$$(8.4.8) \qquad W[f, \overline{\overline{e}}_2, a] = W[f, \overline{\overline{e}}_1 \circ t, \overline{\overline{e}}_1] \circ W[f, \overline{\overline{e}}_1, a]$$

соответствует пассивному преобразованию $t$. Согласно теореме 8.3.3, произведение преобразований координат (8.4.7) и (8.4.8) имеет вид

$$\begin{aligned} W[f, \overline{\overline{e}}_3, a] &= W[f, \overline{\overline{e}}_1 \circ t, \overline{\overline{e}}_1] \circ W[f, \overline{\overline{e}}_1 \circ s, \overline{\overline{e}}_1] \circ W[f, \overline{\overline{e}}_1, a] \\ &= W[f, \overline{\overline{e}}_1 \circ t \circ s, \overline{\overline{e}}_1] \circ W[f, \overline{\overline{e}}_1, a] \end{aligned}$$

и является координатным преобразованием, соответствующим пассивному преобразованию $s \circ t$. Согласно теоремам 8.2.13, 8.2.14 и определению 5.1.11 преобразования координат порождают правостороннее контравариантное представление группы $GA(f)$.

Если координатное преобразование не изменяет координаты выбранного базиса, то ему соответствует единица группы $GA(f)$, так как пассивное представление однотранзитивно. Следовательно, координатное представление эффективно. $\qquad\square$

---

[8.13] Это преобразование не порождает эндоморфизма диаграммы представлений $(f, A)$. Координаты меняются, поскольку меняется базис, относительно которого мы определяем координаты. Однако кортеж A-чисел, координаты которого мы рассматриваем, не меняется.

Рассмотрим диаграммы представлени $(f, A)$, $(B, g)$. Пассивное представление $P(g)$ согласовано с пассивным представлением $P(f)$, если существует гомоморфизм $h$ группы $GA(f)$ в группу $GA(g)$. Рассмотрим диаграмму

$$
\begin{array}{ccc}
\mathrm{End}(B[f]) & \xrightarrow{\ H\ } & \mathrm{End}(B[g]) \\[4pt]
\Big\uparrow{\scriptstyle P(f)} & \nearrow{\scriptstyle f'} & \Big\uparrow{\scriptstyle P(g)} \\[4pt]
GA(f) & \xrightarrow{\ h\ } & GA(g)
\end{array}
$$

Так как отображения $P(f)$, $P(g)$ являются изоморфизмами группы, то отображение $H$ является гомоморфизмом групп. Следовательно, отображение $f'$ является представлением группы $GA(f)$ в многообразии базисов $\mathcal{B}(g)$. Согласно построению, пассивному преобразованию $s$ многообразии базисов $\mathcal{B}(f)$ соответствует пассивное преобразование $H(s)$ многообразия базисов $\mathcal{B}(g)$

$$(8.4.9) \qquad \overline{\overline{e}}_{g1} = \overline{\overline{e}}_g \circ H(s)$$

Тогда координатное преобразование в диаграмме представлений $(B, g)$ принимает вид

$$(8.4.10) \qquad W[g, \overline{\overline{e}}_{g1}, a] = W[g, \overline{\overline{e}}_g \circ H(s), \overline{\overline{e}}_g] \circ W[g, \overline{\overline{e}}_g, a]$$

Определение 8.4.3. *Мы будем называть орбиту*

$$
\begin{aligned}
\mathcal{O}(f, g, \overline{\overline{e}}_g, a) &= H(GA(f)) \circ W[g, \overline{\overline{e}}_g, a] \\
&= (W[g, \overline{\overline{e}}_g \circ H(s), \overline{\overline{e}}_g] \circ W[g, \overline{\overline{e}}_g, a], \overline{\overline{e}}_f \circ s, s \in GA(f))
\end{aligned}
$$

**геометрическим объектом в координатном представлении**, *определённом в диаграмме представлений* $(f, A)$. *Для любого базиса* $\overline{\overline{e}}_{f1} = \overline{\overline{e}}_f \circ s$ *соответствующая точка* (8.4.10) *орбиты определяет* **координаты геометрического объекта** *относительно базиса* $\overline{\overline{e}}_{f1}$. □

Определение 8.4.4. *Мы будем называть орбиту*

$$\mathcal{O}(f, g, a) = (W[g, \overline{\overline{e}}_g \circ H(s), \overline{\overline{e}}_g] \circ W[g, \overline{\overline{e}}_g, a], \overline{\overline{e}}_g \circ H(s), \overline{\overline{e}}_f \circ s, s \in GA(f))$$

**геометрическим объектом**, *определённым в диаграмме представлений* $(f, A)$. *Мы будем также говорить, что $a$ - это* **геометрический объект типа** $H$. *Для любого базиса* $\overline{\overline{e}}_{f1} = \overline{\overline{e}}_f \circ s$ *соответствующая точка* (8.4.10) *орбиты определяет кортеж $A$-чисел*

$$a = \overline{\overline{e}}_g \circ W[g, \overline{\overline{e}}_g, a]$$

*который мы называем* **представителем геометрического объекта** *в диаграмме представлений* $(f, A)$. □

Так как геометрический объект - это орбита представления, то согласно теореме 5.3.7 определение геометрического объекта корректно.

Определение 8.4.3 строит геометрический объект в координатном пространстве. Определение 8.4.4 предполагает, что мы выбрали базис представления $g$. Это позволяет использовать представитель геометрического объекта вместо его координат.

Теорема 8.4.5 (принцип инвариантности). *Представитель геометрического объекта не зависит от выбора базиса* $\overline{\overline{e}}_f$.

Доказательство. Чтобы определить представителя геометрического объекта, мы должны выбрать базис $\overline{\overline{e}}_f$ диаграммы представлений $(f, A)$, базис $\overline{\overline{e}}_g$ диаграммы представлений $(B, g)$ и координаты геометрического объекта $W[g, \overline{\overline{e}}_g, b]$. Соответствующий представитель геометрического объекта имеет вид

$$b = \overline{\overline{e}}_g \circ W[g, \overline{\overline{e}}_g, b]$$

Базис $\overline{\overline{e}}_{f1}$ связан с базисом $\overline{\overline{e}}_f$ пассивным преобразованием

$$\overline{\overline{e}}_{f1} = \overline{\overline{e}}_f \circ s$$

Согласно построению это порождает пассивное преобразование (8.4.9) и координатное преобразование (8.4.10). Соответствующий представитель геометрического объекта имеет вид

$$b' = \overline{\overline{e}}_{g1} \circ W[g, \overline{\overline{e}}_{g1}, b']$$
$$= \overline{\overline{e}}_g \circ W[g, \overline{\overline{e}}_g, \overline{\overline{e}}_g \circ H(s)] \circ W[g, \overline{\overline{e}}_g \circ H(s), \overline{\overline{e}}_g] \circ W[g, \overline{\overline{e}}_g, b]$$
$$= \overline{\overline{e}}_g \circ W[g, \overline{\overline{e}}_g, b] = b$$

Следовательно, представитель геометрического объекта инвариантен относительно выбора базиса. □

# Примеры диаграммы представлений: модуль

## 9.1. Об этой главе

Теория представлений универсальной алгебры - это важный инструмент, которым я пользуюсь на протяжении многих лет для изучения алгебры, геометрии, математического анализа. Основная задача этой и следующей глав - показать как работает теория представлений универсальной алгебры в различных разделах математики.

Примеры в этой главе имеют отношение к различным конструкциям, связанным с модулем над кольцом.

Первый пример - это абелевая группа. Модуль - это эффективное представление кольца в абелевой группе. Поэтому существует параллель между абелевой группой и модулем. Я рассматриваю это сходство в разделе 9.2.

Модуль над коммутативным кольцом - это относительно простая конструкция. С другой стороны, многие определения теории представлений (базис представления, морфизм представлений, свободное представление) опираются на аналогичные определения в модуле. Поэтому раздел 9.3 посвящён детальному рассмотрению модуля над коммутативным кольцом.

Я рассматриваю алгебру над коммутативным кольцом в разделе 9.4 и левый модуль над $D$-алгеброй в разделе 9.5. Мы можем рассматривать модуль над некоммутативным кольцом также как мы рассматривали модуль над коммутативным кольцом. Однако мы встречаем серьёзные проблемы при изучении линейного отображения.

Рассмотрение некоммутативного кольца как алгебры над центром кольца существенно меняет картину. Анализ диаграммы представлений, описывающей модуль $V$ над $D$-алгеброй $A$, позволяет рассмотреть различные группы отображений, сохраняющих структуру алгебры. Среди этих отображений мы выделяем линейные отображения $A$-модуля $V$ (приведенный морфизм $D$-модуля $V$) и гомоморфизм $A$-модуля $V$ (приведенный морфизм диаграммы представлений). Такое определение линейного отображения позволяет рассмотреть полилинейное отображение модуля над $D$-алгеброй $A$.

Если $D$-алгебра $A$ является банаховой, то мы получаем инструмент для изучения математического анализа функций нескольких переменных. К сожалению, структура линейного отображения некоммутативной алгебры лежит вне рамок этой главы. Подробнее эту тему читатель может изучить в книге [12].

## 9.2. Абелева группа

Определение 9.2.1. *Мы определим действие кольца целых чисел $Z$ в абелевой группе $G$ согласно правилу*

$$(9.2.1) \qquad 0g = 0$$

$$(9.2.2) \qquad (n+1)g = ng + g$$

$$(9.2.3) \qquad (n-1)g = ng - g$$

□

Теорема 9.2.2. *Действие кольца целых чисел $Z$ в абелевой группе $G$, рассмотренное в определении 9.2.1, является представлением. Верны следующие равенства*

$$(9.2.4) \qquad 1a = a$$

$$(9.2.5) \qquad (nm)a = n(ma)$$

$$(9.2.6) \qquad (m+n)a = ma + na$$

$$(9.2.7) \qquad (m-n)a = ma - na$$

$$(9.2.8) \qquad n(a+b) = na + nb$$

Доказательство. Равенство (9.2.4) является следствием равенства (9.2.1) и равенства (9.2.2), когда $n = 0$.

Из равенства (9.2.1) следует, что равенство (9.2.6) верно, когда $n = 0$.

- Пусть равенство (9.2.6) верно, когда $n = k \geq 0$. Тогда

$$(m+k)a = ma + ka$$

Равенство

$$(m+(k+1))a = ((m+k)+1)a = (m+k)a + a = ma + ka + a$$
$$= ma + (k+1)a$$

является следствием равенства (9.2.2). Следовательно, равенство (9.2.6) верно, когда $n = k+1$. Согласно принципу математической индукции, равенство (9.2.6) верно для любого $n \geq 0$.

- Пусть равенство (9.2.6) верно, когда $n = k \leq 0$. Тогда

$$(m+k)a = ma + ka$$

Равенство

$$(m+(k-1))a = ((m+k)-1)a = (m+k)a - a = ma + ka - a$$
$$= ma + (k-1)a$$

является следствием равенства (9.2.3). Следовательно, равенство (9.2.6) верно, когда $n = k-1$. Согласно принципу математической индукции, равенство (9.2.6) верно для любого $n \leq 0$.

- Следовательно, равенство (9.2.6) верно для любого $n \in Z$.

Равенство

$$(9.2.9) \qquad (k+n)a - na = ka$$

является следствием равенства (9.2.6). Равенство (9.2.7) является следствием равенства (9.2.9), если мы положим $m = k+n$, $k = m-n$.

Из равенства (9.2.1) следует, что равенство (9.2.5) верно, когда $n = 0$.

- Пусть равенство (9.2.5) верно, когда $n = k \geq 0$. Тогда

$$(km)a = k(ma)$$

Равенство

$$((k+1)m)a = (km+m)a = (km)a + ma = k(ma) + ma$$
$$= (k+1)(ma)$$

является следствием равенств (9.2.2), (9.2.6). Следовательно, равенство (9.2.5) верно, когда $n = k+1$. Согласно принципу математической индукции, равенство (9.2.5) верно для любого $n \geq 0$.

- Пусть равенство (9.2.6) верно, когда $n = k \leq 0$. Тогда

$$(km)a = k(ma)$$

Равенство

$$((k-1)m)a = (km-m)a = (km)a - ma = k(ma) - ma$$
$$= (k-1)(ma)$$

является следствием равенств (9.2.3), (9.2.7). Следовательно, равенство (9.2.5) верно, когда $n = k-1$. Согласно принципу математической индукции, равенство (9.2.5) верно для любого $n \leq 0$.

- Следовательно, равенство (9.2.5) верно для любого $n \in Z$.

Из равенства (9.2.1) следует, что равенство (9.2.8) верно, когда $n = 0$.

- Пусть равенство (9.2.8) верно, когда $n = k \geq 0$. Тогда

$$k(a+b) = ka + kb$$

Равенство

$$(k+1)(a+b) = k(a+b) + a + b = ka + kb + a + b$$
$$= ka + a + kb + b$$
$$= (k+1)a + (k+1)b$$

является следствием равенства (9.2.2). Следовательно, равенство (9.2.8) верно, когда $n = k+1$. Согласно принципу математической индукции, равенство (9.2.8) верно для любого $n \geq 0$.

- Пусть равенство (9.2.6) верно, когда $n = k \leq 0$. Тогда

$$k(a+b) = ka + kb$$

Равенство

$$(k-1)(a+b) = k(a+b) - (a+b) = ka + kb - a - b$$
$$= ka - a + kb - b$$
$$= (k-1)a + (k-1)b$$

является следствием равенства (9.2.3). Следовательно, равенство (9.2.8) верно, когда $n = k-1$. Согласно принципу математической индукции, равенство (9.2.8) верно для любого $n \leq 0$.

- Следовательно, равенство (9.2.8) верно для любого $n \in Z$.

Из равенства (9.2.8) следует, что отображение

$$\varphi(n) : a \in G \to na \in G$$

является эндоморфизмом абелевой группы $G$. Из равенств (9.2.6), (9.2.5) следует, что отображение

$$\varphi : Z \to \mathrm{End}(Ab, G)$$

является гомоморфизмом кольца $Z$. Согласно определению 3.1.1, отображение $\varphi$ является представлением кольца целых чисел $Z$ в абелевой группе $G$.     □

Теорема 9.2.3. *Пусть $G$ - абелева группа. Множество $G$-чисел, порождённое множеством $S = \{s_i : i \in I\}$, имеет вид*

$$(9.2.10) \qquad J(S) = \left\{ g : g = \sum_{i \in I} g^i s_i, g^i \in Z \right\}$$

*где множество $\{i \in I : g^i \neq 0\}$ конечно.*

Доказательство. Мы докажем теорему по индукции, опираясь на теоремы [14]-5.1, страница 94, и 6.1.4.

Для произвольного $s_k \in S$, положим $g^i = \delta^i_k$. Тогда

$$(9.2.11) \qquad s_k = \sum_{i \in I} g^i s_i$$

$s_k \in J(S)$ следует из (9.2.10), (9.2.11).

Пусть $g_1, g_2 \in X_k \subseteq J(S)$. Так как $G$ является абелевой группой, то, согласно утверждению 6.1.4.3, $g_1 + g_2 \in J(S)$. Согласно равенству (9.2.10), существуют $Z$-числа $g^i_1, g^i_2, i \in I$, такие, что

$$(9.2.12) \qquad g_1 = \sum_{i \in I} g^i_1 v_i \qquad g_2 = \sum_{i \in I} g^i_1 v_i$$

где множества

$$(9.2.13) \qquad H_1 = \{i \in I : g^i_1 \neq 0\} \qquad H_2 = \{i \in I : g^i_2 \neq 0\}$$

конечны. Из равенства (9.2.12) следует, что

$$(9.2.14) \qquad g_1 + g_2 = \sum_{i \in I} g^i_1 v_i + \sum_{i \in I} g^i_2 v_i = \sum_{i \in I} (g^i_1 v_i + g^i_2 v_i)$$

Равенство

$$(9.2.15) \qquad g_1 + g_2 = \sum_{i \in I} (g^i_1 + g^i_2) v_i$$

является следствием равенств (9.2.6), (9.2.14). Из равенства (9.2.13) следует, что множество

$$\{i \in I : g^i_1 + g^i_2 \neq 0\} \subseteq H_1 \cup H_2$$

конечно. Из равенства (9.2.15) следует, что $g_1 + g_2 \in J(S)$.     □

### 9.3. Векторное пространство

#### 9.3.1. Модуль над коммутативным кольцом.

Определение 9.3.1. *Эффективное представление коммутативного кольца $D$ в абелевой группе $V$*

$$(9.3.1) \qquad f : D \longrightarrow V \qquad f(d) : v \to dv$$

*называется* **модулем над кольцом** $D$ *или* $D$**-модулем.** $V$*-число называется* **вектором.** $\quad\square$

Теорема 9.3.2. *Следующая диаграмма представлений описывает $D$-модуль $V$*

$$(9.3.2)$$

$$
\begin{array}{c}
D \xrightarrow{\;g_2\;} V \\
\big\uparrow {\scriptstyle g_1} \\
Z
\end{array}
$$

*В диаграмме представлений* (9.3.2) *верна* **коммутативность представлений** *кольца целых чисел $Z$ и коммутативного кольца $D$ в абелевой группе $V$*

$$(9.3.3) \qquad a(nv) = n(av)$$

Доказательство. Диаграмма представлений (9.3.2) является следствием определения 9.3.1 и теоремы 9.2.2. Равенство (9.3.3) является следствием утверждения, что преобразование $g_2(a)$ является эндоморфизмом $Z$-модуля $V$. $\quad\square$

Теорема 9.3.3. *Пусть $V$ является $D$-модулем. Для любого вектора $v \in V$, вектор, порождённый диаграммой представлений* (9.3.2), *имеет следующий вид*

$$(9.3.4) \qquad (a + n)v = av + nv \quad a \in D \quad n \in Z$$

9.3.3.1: *Множество отображений*

$$(9.3.5) \qquad a + n : v \in V \to (a + n)v \in V$$

*порождает* [9.1] *кольцо $D_{(1)}$ где сложение определено равенством*

$$(9.3.6) \qquad (a + n) + (b + m) = (a + b) + (n + m)$$

*и произведение определено равенством*

$$(9.3.7) \qquad (a + n)(b + m) = (ab + ma + nb) + (nm)$$

*Кольцо $D_{(1)}$ называется* **унитальным расширением** *кольца $D$.*

| | | |
|---|---|---|
| *Если кольцо $D$ имеет единицу, то* | $Z \subseteq D$ | $D_{(1)} = D$ |
| *Если кольцо $D$ является идеалом $Z$, то* | $D \subseteq Z$ | $D_{(1)} = Z$ |
| *В противном случае* | | $D_{(1)} = D \oplus Z$ |

9.3.3.2: *Кольцо $D$ является идеалом кольца $D_{(1)}$.*

---

[9.1] Смотри определение унитального расширения также на страницах [6]-52, [7]-64.

9.3.3.3: *Множество преобразований* (9.3.4) *порождает представление кольца* $D_{(1)}$ *в абелевой группе* $V$.

Мы будем пользоваться обозначением $D_{(1)}v$ *для множества векторов, порождённых вектором* $v$.

Теорема 9.3.4. *Элементы $D$-модуля $V$ удовлетворяют соотношениям*

9.3.4.1: **закон ассоциативности**

$$(9.3.8) \qquad (pq)v = p(qv)$$

9.3.4.2: **закон дистрибутивности**

$$(9.3.9) \qquad p(v + w) = pv + pw$$

$$(9.3.10) \qquad (p + q)v = pv + qv$$

9.3.4.3: **закон унитарности**

$$(9.3.11) \qquad 1v = v$$

*для любых* $p, q \in D_{(1)}$, $v, w \in V$.

Доказательство теорем 9.3.3, 9.3.4. Пусть $v \in V$.

Лемма 9.3.5. *Пусть* $n \in Z$, $a \in D$. *Отображение* (9.3.5) *является эндоморфизмом абелевой группы* $V$.

Доказательство. Утверждения $nv \in V$, $av \in V$ являются следствием теорем 6.1.4, 9.3.2. Так как $V$ является абелевой группой, то

$$nv + av \in V \quad n \in Z \quad a \in D$$

Следовательно, для любого $Z$-числа $n$ и любого $D$-числа $a$, мы определили отображение (9.3.5). Поскольку преобразование $g_1(n)$ и преобразование $g_2(a)$ являются эндоморфизмами абелевой группы $V$, то отображение (9.3.5) является эндоморфизмом абелевой группы $V$.                ⊙

Пусть $D_{(1)}$ - множество отображений (9.3.5). Равенство (9.3.9) является следствием леммы 9.3.5.

Пусть $p = a + n \in D_{(1)}$, $q = b + m \in D_{(1)}$. Согласно утверждению 9.3.3.3, мы определим сумму $D_{(1)}$-чисел $p$ и $q$ равенством (9.3.10). Равенство

$$(9.3.12) \qquad ((a + n) + (b + m))v = (a + n)v + (b + m)v$$

является следствием равенства (9.3.10). Равенство

$$(9.3.13) \qquad (n + m)v = nv + mv$$

является следствием утверждения, что представление $g_1$ является гомоморфизмом аддитивной группы кольца $Z$. Равенство

$$(9.3.14) \qquad (a + b)v = av + bv$$

является следствием утверждения, что представление $g_2$ является гомоморфизмом аддитивной группы кольца $D$. Так как $V$ является абелевой группой, то равенство

$$(9.3.15) \qquad \begin{aligned} ((a + n) + (b + m))v &= av + nv + bv + mv = av + bv + nv + mv \\ &= (a + b)v + (n + m)v = ((a + b) + (n + m))v \end{aligned}$$

является следствием равенств (9.3.12), (9.3.13), (9.3.14). Из равенства (9.3.15) следует, что определение (9.3.6) суммы на множестве $D_{(1)}$ не зависит от вектора $v$.

Равенства (9.3.8), (9.3.11) являются следствием утверждения 9.3.3.3. Пусть $p = a + n \in D_{(1)}$, $q = b + m \in D_{(1)}$. Равенство

$$(9.3.16) \qquad (mn)v = m(nv)$$

является следствием утверждения, что представление $g_1$ является представлением мультипликативной группы кольца $Z$. Равенство

$$(9.3.17) \qquad (ab)v = a(bv)$$

является следствием утверждения, что представление $g_2$ является представлением мультипликативной группы кольца $D$. Равенство

$$(9.3.18) \qquad (md)v = m(dv)$$

является следствием утверждения, что кольцо $D$ является абелевой группой. Равенство

$$(9.3.19) \qquad \begin{aligned} ((a+n)(b+m))v &= (a+n)((b+m)v) = (a+n)(bv+mv) \\ &= a(bv+mv) + n(bv+mv) \\ &= a(bv) + a(mv) + n(bv) + n(mv) \\ &= (ab)v + m(av) + +n(bv) + (nm)v \\ &= (ab)v + (ma)v + +(nb)v + (nm)v \\ &= ((ab+ma+nb)+nm)v \end{aligned}$$

является следствием равенств (9.3.3), (9.3.4), (9.3.8), (9.3.16), (9.3.17), (9.3.18). Равенство (9.3.7) является следствием равенства (9.3.19).

Утверждение 9.3.3.2 является следствием равенства (9.3.7).                    $\square$

Теорема 9.3.6. *Пусть $V$ - $D$-модуль. Множество векторов, порождённое множеством векторов $v = (v_i \in V, i \in I)$, имеет вид*[9.2]

$$(9.3.20) \qquad J(v) = \left\{ w : w = \sum_{i \in I} c^i v_i, c^i \in D_{(1)}, |\{i : c^i \neq 0\}| < \infty \right\}$$

Доказательство. Мы докажем теорему по индукции, опираясь на теорему 6.1.4, Согласно теореме 6.1.4, мы должны доказать следующие утверждения:

9.3.6.1: $v_k \in X_0 \subseteq J(v)$

9.3.6.2: $c^k v_k \in J(v)$, $c^k \in D_{(1)}$, $k \in I$

9.3.6.3: $\displaystyle\sum_{k \in I} c^k v_k \in J(v)$, $c^k \in D_{(1)}$, $|\{i : c^i \neq 0\}| < \infty$

9.3.6.4: $w_1, w_2 \in J(v) \Rightarrow w_1 + w_2 \in J(v)$

9.3.6.5: $a \in D$, $w \in J(v) \Rightarrow aw \in J(v)$

• Для произвольного $v_k \in v$, положим $c^i = \delta^i_k \in D_{(1)}$. Тогда

$$(9.3.21) \qquad v_k = \sum_{i \in I} c^i v_i$$

___

[9.2] Для множества $A$, мы обозначим $|A|$ мощность множества $A$. Запись $|A| < \infty$ означает, что множество $A$ конечно.

Утверждение 9.3.6.1 следует из (9.3.20), (9.3.21).

- Утверждение 9.3.6.2 являются следствием теорем 6.1.4, 9.3.3 и утверждения 9.3.6.1.
- Так как $V$ является абелевой группой, то утверждение 9.3.6.3 следует из утверждения 9.3.6.2 и теорем 6.1.4, 9.2.3.
- Пусть $w_1$, $w_2 \in X_k \subseteq J(v)$. Так как $V$ является абелевой группой, то, согласно утверждению 6.1.4.3,

$$(9.3.22) \qquad w_1 + w_2 \in X_{k+1}$$

Согласно равенству (9.3.20), существуют $D_{(1)}$-числа $w_1^i$, $w_2^i$, $i \in I$, такие, что

$$(9.3.23) \qquad w_1 = \sum_{i \in I} w_1^i v_i \qquad w_2 = \sum_{i \in I} w_1^i v_i$$

где множества

$$(9.3.24) \qquad H_1 = \{ i \in I : w_1^i \neq 0 \} \qquad H_2 = \{ i \in I : w_2^i \neq 0 \}$$

конечны. Так как $V$ является абелевой группой, то из равенства (9.3.23) следует, что

$$(9.3.25) \qquad w_1 + w_2 = \sum_{i \in I} w_1^i v_i + \sum_{i \in I} w_2^i v_i = \sum_{i \in I} (w_1^i v_i + w_2^i v_i)$$

Равенство

$$(9.3.26) \qquad w_1 + w_2 = \sum_{i \in I} (w_1^i + w_2^i) v_i$$

является следствием равенств (9.3.10), (9.3.25). Из равенства (9.3.24) следует, что множество

$$\{ i \in I : w_1^i + w_2^i \neq 0 \} \subseteq H_1 \cup H_2$$

конечно.

- Пусть $w \in X_k \subseteq J(v)$. Согласно утверждению 6.1.4.4, для любого $D_{(1)}$-числа $a$,

$$(9.3.27) \qquad aw \in X_{k+1}$$

Согласно равенству (9.3.20), существуют $D_{(1)}$-числа $w^i$, $i \in I$, такие, что

$$(9.3.28) \qquad w = \sum_{i \in I} w^i v_i$$

где

$$(9.3.29) \qquad |\{ i \in I : w^i \neq 0 \}| < \infty$$

Из равенства (9.3.28) следует, что

$$(9.3.30) \qquad aw = a \sum_{i \in I} w^i v_i = \sum_{i \in I} a(w^i v_i) = \sum_{i \in I} (aw^i) v_i$$

Из утверждения (9.3.29) следует, что множество $\{ i \in I : aw^i \neq 0 \}$ конечно.

Из равенств (9.3.22), (9.3.26), (9.3.27), (9.3.30) следует, что $X_{k+1} \subseteq J(v)$. $\qquad \square$

Определение 9.3.7. *Пусть* $v = (v_i \in V, i \in I)$ *- множество векторов. Выражение* $w^i v_i$ *называется* **линейной комбинацией** *векторов* $v_i$. *Вектор* $\overline{w} = w^i v_i$ *называется* **линейно зависимым** *от векторов* $v_i$.    □

Представим множество $D_{(1)}$-чисел $w^i, i \in I$, в виде матрицы

$$w = \begin{pmatrix} w^1 \\ ... \\ w^n \end{pmatrix}$$

Представим множество векторов $v_i, i \in I$, в виде матрицы

$$v = \begin{pmatrix} v_1 & ... & v_n \end{pmatrix}$$

Тогда мы можем записать линейную комбинацию векторов $\overline{w} = w^i v_i$ в виде

$$\overline{w} = w^* {}_* v$$

Теорема 9.3.8. *Пусть $D$ - поле. Если уравнение*

$$w^i v_i = 0$$

*предполагает существования индекса $i = j$ такого, что $w^j \neq 0$, то вектор $v_j$ линейно зависит от остальных векторов $v$.*

Доказательство. Теорема является следствием равенства

$$v_j = \sum_{i \in I \setminus \{j\}} \frac{w^i}{w^j} v_i$$

и определения 9.3.7.                                          □

Очевидно, что для любого множества векторов $v_i$,

$$w^i = 0 \Rightarrow w^* {}_* v = 0$$

Определение 9.3.9. *Множество векторов* [9.3] $v_i, i \in I$, *$D$-модуля $V$* **линейно независимо***, если $w = 0$ следует из уравнения*

$$w^i v_i = 0$$

*В противном случае, множество векторов $v_i, i \in I$,* **линейно зависимо***.*
                                                             □

Следующее определение является следствием теорем 9.3.6, 6.1.4 и определения 6.1.5.

Определение 9.3.10. *$J(v)$ называется* **подмодулем, порождённым множеством** *$v$, а $v$ -* **множеством образующих** *подмодуля $J(v)$. В частности,* **множеством образующих** *$D$-модуля $V$ будет такое подмножество $X \subset V$, что $J(X) = V$.*    □

Следующее определение является следствием теорем 9.3.6, 6.1.4 и определения 6.2.6.

---

[9.3] Я следую определению в [2], страница 100.

Определение 9.3.11. *Если множество $X \subset V$ является множеством образующих $D$-модуля $V$, то любое множество $Y$, $X \subset Y \subset V$ также является множеством образующих $D$-модуля $V$. Если существует минимальное множество $X$, порождающее $D$-модуль $V$, то такое множество $X$ называется* **базисом** *$D$-модуля $V$.* □

Теорема 9.3.12. *Множество векторов $\overline{\overline{e}} = (e_i, i \in I)$ является базисом $D$-модуля $V$, если верны следующие утверждения.*

9.3.12.1: *Произвольный вектор $v \in V$ является линейной комбинацией векторов множества $\overline{\overline{e}}$.*

9.3.12.2: *Вектор $e_i$ нельзя представить в виде линейной комбинации остальных векторов множества $\overline{\overline{e}}$.*

Доказательство. Согласно утверждению 9.3.12.1, теореме 9.3.6 и определению 9.3.7, множество $\overline{\overline{e}}$ порождает $D$-модуль $V$ (определение 9.3.10). Согласно утверждению 9.3.12.2, множество $\overline{\overline{e}}$ является минимальным множеством, порождающим $D$-модуль $V$. Согласно определению 9.3.11, множество $\overline{\overline{e}}$ является базисом $D$-модуля $V$. □

Теорема 9.3.13. *Пусть $D$ - поле. Множество векторов $\overline{\overline{e}} = (e_i, i \in I)$ является* **базисом $D$-векторного пространства** *$V$, если векторы $e_i$ линейно независимы и любой вектор $v \in V$ линейно зависит от векторов $e_i$.*

Доказательство. Пусть множество векторов $e_i$, $i \in I$, линейно зависимо. Тогда в равенстве

$$w^i e_i = 0$$

существует индекс $i = j$ такой, что $w^j \neq 0$. Согласно теореме 9.3.8, вектор $e_j$ линейно зависит от остальных векторов множества $\overline{\overline{e}}$. Согласно определению 9.3.11, множество векторов $e_i$, $i \in I$, не является базисом $D$-векторного пространства $V$.

Следовательно, если множество векторов $e_i$, $i \in I$, является базисом, то эти векторы линейно независимы. Так как произвольный вектор $v \in V$ является линейной комбинацией векторов $e_i$, $i \in I$, , то множество векторов $v$, $e_i$, $i \in I$, не является линейно независимым. □

Определение 9.3.14. *Пусть $\overline{\overline{e}}$ - базис $D$-модуля $V$, и вектор $\overline{v} \in V$ имеет разложение*

$$\overline{v} = v^*_* e = v^i e_i$$

*относительно базиса $\overline{\overline{e}}$. $D_{(1)}$-числа $v^i$ называются* **координатами** *вектора $\overline{v}$ относительно базиса $\overline{\overline{e}}$. Матрица $D_{(1)}$-чисел $v = (v^i, i \in I)$ называется* **координатной матрицей вектора** *$\overline{v}$ в базисе $\overline{\overline{e}}$.* □

Теорема 9.3.15. *Пусть $D$ - кольцо. Пусть $\overline{\overline{e}}$ - базис $D$-модуля $V$. Пусть*

$$(9.3.31) \qquad w^i e_i = 0$$

*линейная зависимость векторов базиса $\overline{\overline{e}}$. Тогда*

9.3.15.1: *$D_{(1)}$-число $w^i$, $i \in I$, не имеет обратного элемента в кольце $D_{(1)}$.*

9.3.15.2: *Множество $D'$ матриц $w = (w^i, i \in I)$ порождает $D$-модуль.*

Доказательство. Допустим существует матрица $w = (w^i, i \in I)$ такая, что равенство (9.3.31) верно и существует индекс $i = j$ такой, что $w^j \neq 0$. Если мы положим, что $D_{(1)}$-число $c^j$ имеет обратный, то равенство

$$e_j = \sum_{i \in I \setminus \{j\}} \frac{w^i}{w^j} e_i$$

является следствием равенства (9.3.31). Следовательно вектор $e_j$ является линейной комбинацией остальных векторов множества $\overline{\overline{e}}$ и множество $\overline{\overline{e}}$ не является базисом. Следовательно, наше предположение неверно, и $D_{(1)}$-число $c^j$ не имеет обратного.

Пусть матрицы $b = (b^i, i \in I) \in D', \quad c = (c^i, i \in I) \in D'$. Из равенств

$$b^i e_i = 0$$
$$c^i e_i = 0$$

следует

$$(b^i + c^i) e_i = 0$$

Следовательно, множество $D'$ является абелевой группой.

Пусть матрица $c = (c^i, i \in I) \in D'$ и $a \in D$. Из равенства

$$c^i e_i = 0$$

следует

$$(a c^i) e_i = 0$$

Следовательно, абелева группа $D'$ является $D$-модулем.    □

Теорема 9.3.16. *Пусть $D$-модуль $V$ имеет базис $\overline{\overline{e}}$ такой, что в равенстве*

$$(9.3.32) \qquad\qquad\qquad w^i e_i = 0$$

*существует индекс $i = j$ такой, что $w^j \neq 0$. Тогда*

9.3.16.1: *Матрица $w = (w^i, i \in I)$ определяет координаты вектора $0 \in V$ относительно базиса $\overline{\overline{e}}$.*

9.3.16.2: *Координаты вектора $\overline{v}$ относительно базиса $\overline{\overline{e}}$ определены однозначно с точностью до выбора координат вектора $0 \in V$.*

Доказательство. Утверждение 9.3.16.1 является следствием равенства (9.3.32) и определения 9.3.14.

Пусть вектор $\overline{v}$ имеет разложение

$$(9.3.33) \qquad\qquad\qquad \overline{v} = v^* {}_* e = v^i e_i$$

относительно базиса $\overline{\overline{e}}$. Равенство

$$(9.3.34) \qquad\qquad \overline{v} = \overline{v} + 0 = v^i e_i + c^i e_i = (v^i + c^i) e_i$$

является следствием равенств (9.3.32), (9.3.33). Утверждение 9.3.16.2 является следствием равенств (9.3.33), (9.3.34) и определения 9.3.14.    □

Определение 9.3.17. $D$-модуль $V$ - **свободный $D$-модуль**, [9.4] *если $D$-модуль $V$ имеет базис и векторы базиса линейно независимы.*    □

---

[9.4] Я следую определению в [2], страница 103.

Теорема 9.3.18. *Координаты вектора $v \in V$ относительно базиса $\overline{\overline{e}}$ свободного $D$-модуля $V$ определены однозначно.*

Доказательство. Теорема является следствием теоремы 9.3.16 и определений 9.3.9, 9.3.17. □

Пример 9.3.19. *Из теоремы 9.2.2 и определения 9.3.1 следует, что абелева группа $G$ является модулем над кольцом целых чисел $Z$.* □

### 9.3.2. Линейное отображение.

Определение 9.3.20. *Морфизм представлений*

$$\left( h : D_1 \to D_2 \quad f : V_1 \to V_2 \right)$$

*$D_1$-модуля $A_1$ в $D_2$-модуль $A_2$ называется* **линейным отображением** *$D_1$-модуля $A_1$ в $D_2$-модуль $A_2$. Обозначим* $\mathcal{L}(D_1 \to D_2; A_1 \to A_2)$ *множество линейных отображений $D_1$-модуля $A_1$ в $D_2$-модуль $A_2$.* □

Если отображение

$$f : A_1 \to A_2$$

является линейным отображением $D_1$-алгебры $A_1$ в $D_2$-алгебру $A_2$, то я пользуюсь обозначением

$$f \circ a = f(a)$$

для образа отображения $f$.

Теорема 9.3.21. *Линейное отображение*

$$\left( h : D_1 \to D_2 \quad f : A_1 \to A_2 \right)$$

*$D_1$-модуля $A_1$ в $D_2$-модуль $A_2$ удовлетворяет равенствам*[9.5]

$$(9.3.35) \qquad h(d_1 + d_2) = h(d_1) + h(d_2)$$

$$(9.3.36) \qquad h(d_1 d_2) = h(d_1) h(d_2)$$

$$(9.3.37) \qquad f \circ (a + b) = f \circ a + f \circ b$$

$$(9.3.38) \qquad f \circ (da) = h(d)(f \circ a)$$

$$a, b \in A_1 \quad d, d_1, d_2 \in D_1$$

Доказательство. Из определений 3.2.2, 9.3.20 следует, что отображение $h$ является гомоморфизмом кольца $D_1$ в кольцо $D_2$ (равенства (9.3.35), (9.3.36)) и отображение $f$ является гомоморфизмом абелевой группы $A_1$ в абелеву группу $A_2$ (равенство (9.3.37)). Равенство (9.3.38) является следствием равенства (3.2.3). □

---

[9.5] В некоторых книгах (например, на странице [2]-94) теорема 9.3.21 рассматривается как определение.

Теорема 9.3.22. *Пусть*

$$\overline{\overline{e}}_1 = (e_{1 \cdot i}, i \in I)$$

*базис в* $D_1$-*модуле* $A_1$. *Пусть*

$$\overline{\overline{e}}_2 = (e_{2 \cdot j}, j \in J)$$

*базис в* $D_2$-*модуле* $A_2$. *Тогда линейное отображение*

$$\left( h : D_1 \to D_2 \quad \overline{f} : A_1 \to A_2 \right)$$

*имеет представление*

$$(9.3.39) \qquad\qquad b = h(a)^*{}_* f$$

*относительно заданных базисов. Здесь*

- $a$ - *координатная матрица* $A_1$-*числа* $\overline{a}$ *относительно базиса* $\overline{\overline{e}}_1$

$$(9.3.40) \qquad\qquad \overline{a} = a^*{}_* e_1$$

- $h(a) = (h(a_i), i \in I)$ - *матрица* $D_2$-*чисел.*
- $b$ - *координатная матрица вектора*

$$(9.3.41) \qquad\qquad \overline{b} = \overline{f} \circ \overline{a}$$

*относительно базиса* $\overline{\overline{e}}_2$

$$(9.3.42) \qquad\qquad \overline{b} = b^*{}_* e_2$$

- $f$ - *координатная матрица множества векторов* $(\overline{f} \circ e_{1 \cdot i}, i \in I)$ *относительно базиса* $\overline{\overline{e}}_2$. *Мы будем называть матрицу* $f$ **матрицей линейного отображения** $\overline{f}$ *относительно базисов* $\overline{\overline{e}}_1$ *и* $\overline{\overline{e}}_2$.

Доказательство. Так как

$$\left( h : D_1 \to D_2 \quad \overline{f} : A_1 \to A_2 \right)$$

линейное отображение, то равенство

$$(9.3.43) \qquad \overline{b} = \overline{f} \circ \overline{a} = \overline{f} \circ (a^*{}_* e_1) = h(a)^*{}_* (\overline{f} \circ e_1)$$

является следствием равенств (9.3.38), (9.3.40), (9.3.41). $A_2$-число $\overline{f} \circ e_{1 \cdot i}$ имеет разложение

$$(9.3.44) \qquad\qquad \overline{f} \circ e_{1 \cdot i} = f_i{}^*{}_* e_2 = f_i^j e_{2 \cdot j}$$

относительно базиса $\overline{\overline{e}}_2$. Комбинируя (9.3.43) и (9.3.44), мы получаем

$$(9.3.45) \qquad\qquad \overline{b} = h(a)^*{}_* f^*{}_* e_2$$

(9.3.39) следует из сравнения (9.3.42) и (9.3.45) и теоремы 9.3.18.     □

Определение 9.3.23. *Приведенный морфизм представлений*

$$f : A_1 \to A_2$$

$D$-*модуля* $A_1$ *в* $D$-*модуль* $A_2$ *называется* **линейным отображением** $D$-*модуля* $A_1$ *в* $D$-*модуль* $A_2$. *Обозначим* $\mathcal{L}(D; A_1 \to A_2)$ *множество линейных отображений* $D$-*модуля* $A_1$ *в* $D$-*модуль* $A_2$.     □

Теорема 9.3.24. *Линейное отображение*

$$f : A_1 \to A_2$$

*$D$-модуля $A_1$ в $D$-модуль $A_2$ удовлетворяет равенствам[9.6]*

$$(9.3.46) \qquad f \circ (a + b) = f \circ a + f \circ b$$

$$(9.3.47) \qquad f \circ (da) = d(f \circ a)$$

$$a, b \in A_1 \quad d \in D$$

Доказательство. Из определений 3.4.2, 9.3.23 следует, что отображение $f$ является гомоморфизмом абелевой группы $A_1$ в абелеву группу $A_2$ (равенство (9.3.46)). Равенство (9.3.47) является следствием равенства (3.4.4). □

Теорема 9.3.25. *Пусть*

$$\overline{\overline{e}}_1 = (e_{1 \cdot i}, i \in I)$$

*базис в $D$-модуле $A_1$. Пусть*

$$\overline{\overline{e}}_2 = (e_{2 \cdot j}, j \in J)$$

*базис в $D$-модуле $A_2$. Тогда линейное отображение*

$$\overline{f} : A_1 \to A_2$$

*имеет представление*

$$(9.3.48) \qquad b = a^*{}_* f$$

*относительно заданных базисов. Здесь*

- *$a$ - координатная матрица $A_1$-числа $\overline{a}$ относительно базиса $\overline{\overline{e}}_1$*

$$(9.3.49) \qquad \overline{a} = a^*{}_* e_1$$

- *$b$ - координатная матрица вектора*

$$(9.3.50) \qquad \overline{b} = \overline{f} \circ \overline{a}$$

  *относительно базиса $\overline{\overline{e}}_2$*

$$(9.3.51) \qquad \overline{b} = b^*{}_* e_2$$

- *$f$ - координатная матрица множества векторов $(\overline{f} \circ e_{1 \cdot i}, i \in I)$ относительно базиса $\overline{\overline{e}}_2$. Мы будем называть матрицу $f$ **матрицей линейного отображения** $\overline{f}$ относительно базисов $\overline{\overline{e}}_1$ и $\overline{\overline{e}}_2$.*

Доказательство. Так как

$$\overline{f} : A_1 \to A_2$$

линейное отображение, то равенство

$$(9.3.52) \qquad \overline{b} = \overline{f} \circ \overline{a} = \overline{f} \circ (a^*{}_* e_1) = a^*{}_* (\overline{f} \circ e_1)$$

является следствием равенств (9.3.47), (9.3.49), (9.3.50). $A_2$-число $\overline{f} \circ e_{1 \cdot i}$ имеет разложение

$$(9.3.53) \qquad \overline{f} \circ e_{1 \cdot i} = f_i{}^*{}_* e_2 = f_i^j e_{2 \cdot j}$$

---

[9.6] В некоторых книгах (например, на странице [2]-94) теорема 9.3.24 рассматривается как определение.

относительно базиса $\overline{\overline{e}}_2$. Комбинируя (9.3.52) и (9.3.53), мы получаем

$$(9.3.54) \qquad \overline{b} = a^*_{\phantom{*}*} f^*_{\phantom{*}*} e_2$$

(9.3.48) следует из сравнения (9.3.51) и (9.3.54) и теоремы 9.3.18.          $\square$

### 9.3.3. Полилинейное отображение.

ОПРЕДЕЛЕНИЕ 9.3.26. *Пусть $D$ - коммутативное кольцо. Приведенный полиморфизм $D$-модулей $A_1$, ..., $A_n$ в $D$-модуль $S$*

$$f : A_1 \times ... \times A_n \to S$$

*называется* **полилинейным отображением** *$D$-модулей $A_1$, ..., $A_n$ в $D$-модуль $S$. Обозначим $\mathcal{L}(D; A_1 \times ... \times A_n \to S)$ множество полилинейных отображений $D$-модулей $A_1$, ..., $A_n$ в $D$-модуль $S$. Обозначим $\mathcal{L}(D; A^n \to S)$ множество $n$-линейных отображений $D$-модуля $A$ ($A_1 = ... = A_n = A$) в $D$-модуль $S$.*          $\square$

ТЕОРЕМА 9.3.27. *Пусть $D$ - коммутативное кольцо. Полилинейное отображение $D$-модулей $A_1$, ..., $A_n$ в $D$-модуль $S$*

$$f : A_1 \times ... \times A_n \to S$$

*удовлетворяет равенствам*

$$f \circ (a_1, ..., a_i + b_i, ..., a_n) = f \circ (a_1, ..., a_i, ..., a_n) + f \circ (a_1, ..., b_i, ..., a_n)$$

$$f \circ (a_1, ..., pa_i, ..., a_n) = pf \circ (a_1, ..., a_i, ..., a_n)$$

$$f \circ (a_1, ..., a_i + b_i, ..., a_n) = f \circ (a_1, ..., a_i, ..., a_n) + f \circ (a_1, ..., b_i, ..., a_n)$$

$$f \circ (a_1, ..., pa_i, ..., a_n) = pf \circ (a_1, ..., a_i, ..., a_n)$$

$$1 \le i \le n \quad a_i, b_i \in A_i \quad p \in D$$

ДОКАЗАТЕЛЬСТВО. Теорема является следствием определений 4.4.4, 9.3.23, 9.3.26 и теоремы 9.3.24.          $\square$

ТЕОРЕМА 9.3.28. *Пусть $D$ - коммутативное кольцо. Пусть $A_1$, ..., $A_n$, $S$ - $D$-модули. Отображение*

$$(9.3.55) \qquad f + g : A_1 \times ... \times A_n \to S \quad f, g \in \mathcal{L}(D; A_1 \times ... \times A_n \to S)$$

*определённое равенством*

$$(9.3.56) \qquad (f + g) \circ (a_1, ..., a_n) = f \circ (a_1, ..., a_n) + g \circ (a_1, ..., a_n)$$

*называется* **суммой полилинейных отображений** *$f$ и $g$ и является полилинейным отображением. Множество $\mathcal{L}(D; A_1 \times ... \times A_n \to S)$ является абелевой группой относительно суммы отображений.*

ДОКАЗАТЕЛЬСТВО. Согласно теореме 9.3.27

$$(9.3.57) \quad f \circ (a_1, ..., a_i + b_i, ..., a_n) = f \circ (a_1, ..., a_i, ..., a_n) + f \circ (a_1, ..., b_i, ..., a_n)$$

$$(9.3.58) \qquad f \circ (a_1, ..., pa_i, ..., a_n) = pf \circ (a_1, ..., a_i, ..., a_n)$$

$$(9.3.59) \quad g \circ (a_1, ..., a_i + b_i, ..., a_n) = g \circ (a_1, ..., a_i, ..., a_n) + g \circ (a_1, ..., b_i, ..., a_n)$$

$$(9.3.60) \qquad g \circ (a_1, ..., pa_i, ..., a_n) = pg \circ (a_1, ..., a_i, ..., a_n)$$

Равенство

$$(f + g) \circ (x_1, ..., x_i + y_i, ..., x_n)$$
$$= f \circ (x_1, ..., x_i + y_i, ..., x_n) + g \circ (x_1, ..., x_i + y_i, ..., x_n)$$
(9.3.61) $$= f \circ (x_1, ..., x_i, ..., x_n) + f \circ (x_1, ..., y_i, ..., x_n)$$
$$+ g \circ (x_1, ..., x_i, ..., x_n) + g \circ (x_1, ..., y_i, ..., x_n)$$
$$= (f + g) \circ (x_1, ..., x_i, ..., x_n) + (f + g) \circ (x_1, ..., y_i, ..., x_n)$$

является следствием равенств (9.3.56), (9.3.57), (9.3.59). Равенство

$$(f + g) \circ (x_1, ..., px_i, ..., x_n)$$
$$= f \circ (x_1, ..., px_i, ..., x_n) + g \circ (x_1, ..., px_i, ..., x_n)$$
(9.3.62) $$= pf \circ (x_1, ..., x_i, ..., x_n) + pg \circ (x_1, ..., x_i, ..., x_n)$$
$$= p(f \circ (x_1, ..., x_i, ..., x_n) + g \circ (x_1, ..., x_i, ..., x_n))$$
$$= p(f + g) \circ (x_1, ..., x_i, ..., x_n)$$

является следствием равенств (9.3.56), (9.3.58), (9.3.60). Из равенств (9.3.61), (9.3.62) и теоремы 9.3.27 следует, что отображение (9.3.55) является полилинейным отображением $D$-модулей.

Пусть $f$, $g$, $h \in \mathcal{L}(D; A_1 \times ... \times A_2 \to S)$. Для любого $a = (a_1, ..., a_n)$, $a_1 \in A_1, ..., a_n \in A_n$,

$$(f + g) \circ a = f \circ a + g \circ a = g \circ a + f \circ a$$
$$= (g + f) \circ a$$
$$((f + g) + h) \circ a = (f + g) \circ a + h \circ a = (f \circ a + g \circ a) + h \circ a$$
$$= f \circ a + (g \circ a + h \circ a) = f \circ a + (g + h) \circ a$$
$$= (f + (g + h)) \circ a$$

Следовательно, сумма полилинейных отображений коммутативна и ассоциативна.

Из равенства (9.3.56) следует, что отображение

$$0 : v \in A_1 \times ... \times A_n \to 0 \in S$$

является нулём операции сложения

$$(0 + f) \circ (a_1, ..., a_n) = 0 \circ (a_1, ..., a_n) + f \circ (a_1, ..., a_n) = f \circ (a_1, ..., a_n)$$

Из равенства (9.3.56) следует, что отображение

$$-f : (a_1, ..., a_n) \in A_1 \times ... \times A_n \to -(f \circ (a_1, ..., a_n)) \in S$$

является отображением, обратным отображению $f$

$$f + (-f) = 0$$

так как

$$(f + (-f)) \circ (a_1, ..., a_n) = f \circ (a_1, ..., a_n) + (-f) \circ (a_1, ..., a_n)$$
$$= f \circ (a_1, ..., a_n) - f \circ (a_1, ..., a_n)$$
$$= 0 = 0 \circ (a_1, ..., a_n)$$

Из равенства

$$(f + g) \circ (a_1, ..., a_n) = f \circ (a_1, ..., a_n) + g \circ (a_1, ..., a_n)$$
$$= g \circ (a_1, ..., a_n) + f \circ (a_1, ..., a_n)$$
$$= (g + f) \circ (a_1, ..., a_n)$$

следует, что сумма отображений коммутативно. Следовательно, множество $\mathcal{L}(D; A_1 \times ... \times A_n \to S)$ является абелевой группой. $\qquad\square$

СЛЕДСТВИЕ 9.3.29. *Пусть $A_1$, $A_2$ - $D$-модули. Отображение*

$$(9.3.63) \qquad f + g : A_1 \to A_2 \quad f, g \in \mathcal{L}(D; A_1 \to A_2)$$

*определённое равенством*

$$(9.3.64) \qquad (f + g) \circ x = f \circ x + g \circ x$$

*называется* **суммой отображений** *$f$ и $g$ и является линейным отображением. Множество $\mathcal{L}(D; A_1; A_2)$ является абелевой группой относительно суммы отображений.* $\qquad\square$

ТЕОРЕМА 9.3.30. *Пусть $D$ - коммутативное кольцо. Пусть $A_1$, ..., $A_n$, $S$ - $D$-модули. Отображение*

$$(9.3.65) \qquad d\, f : A_1 \times ... \times A_n \to S \quad d \in D \quad f \in \mathcal{L}(D; A_1 \times ... \times A_n \to S)$$

*определённое равенством*

$$(9.3.66) \qquad (d\, f) \circ (a_1, ..., a_n) = d(f \circ (a_1, ..., a_n))$$

*называется* **произведением отображения** *$f$* **на скаляр** *$d$ и является полилинейным отображением. Представление*

$$(9.3.67) \qquad a : f \in \mathcal{L}(D; A_1 \times ... \times A_n \to S) \to af \in \mathcal{L}(D; A_1 \times ... \times A_n \to S)$$

*кольца $D$ в абелевой группе $\mathcal{L}(D; A_1 \times ... \times A_n \to S)$ порождает структуру $D$-модуля.*

ДОКАЗАТЕЛЬСТВО. Согласно теореме 9.3.27

$$(9.3.68) \quad f \circ (a_1, ..., a_i + b_i, ..., a_n) = f \circ (a_1, ..., a_i, ..., a_n) + f \circ (a_1, ..., b_i, ..., a_n)$$

$$(9.3.69) \qquad f \circ (a_1, ..., pa_i, ..., a_n) = pf \circ (a_1, ..., a_i, ..., a_n)$$

Равенство

$$(pf) \circ (x_1, ..., x_i + y_i, ..., x_n)$$
$$= p\, f \circ (x_1, ..., x_i + y_i, ..., x_n)$$
$$(9.3.70) \qquad = p\,(f \circ (x_1, ..., x_i, ..., x_n) + f \circ (x_1, ..., y_i, ..., x_n))$$
$$= p(f \circ (x_1, ..., x_i, ..., x_n)) + p(f \circ (x_1, ..., y_i, ..., x_n))$$
$$= (pf) \circ (x_1, ..., x_i, ..., x_n) + (pf) \circ (x_1, ..., y_i, ..., x_n)$$

является следствием равенств (9.3.66), (9.3.68). Равенство

$$(pf) \circ (x_1, ..., qx_i, ..., x_n)$$
$$(9.3.71) \qquad = p(f \circ (x_1, ..., qx_i, ..., x_n)) = pq(f \circ (x_1, ..., x_i, ..., x_n))$$
$$= qp(f \circ (x_1, ..., x_n)) = q(pf) \circ (x_1, ..., x_n)$$

является следствием равенств (9.3.66), (9.3.69). Из равенств (9.3.70), (9.3.71) и теоремы 9.3.27 следует, что отображение (9.3.65) является полилинейным отображением $D$-модулей.

Равенство

$$(9.3.72) \qquad (p+q)f = pf + qf$$

является следствием равенства

$$\begin{aligned}((p+q)f) \circ (x_1, ..., x_n) &= (p+q)(f \circ (x_1, ..., x_n)) \\ &= p(f \circ (x_1, ..., x_n)) + q(f \circ (x_1, ..., x_n)) \\ &= (pf) \circ (x_1, ..., x_n) + (qf) \circ (x_1, ..., x_n)\end{aligned}$$

Равенство

$$(9.3.73) \qquad p(qf) = (pq)f$$

является следствием равенства

$$\begin{aligned}(p(qf)) \circ (x_1, ..., x_n) &= p \ (qf) \circ (x_1, ..., x_n) = p \ (q \ f \circ (x_1, ..., x_n)) \\ &= (pq) \ f \circ (x_1, ..., x_n) = ((pq)f) \circ (x_1, ..., x_n)\end{aligned}$$

Из равенств (9.3.72), (9.3.73), следует, что отображение (9.3.67) является представлением кольца $D$ в абелевой группе $\mathcal{L}(D; A_1 \times ... \times A_n \to S)$. Так как указанное представление эффективно, то, согласно определению 9.3.1 и теореме 9.3.28, абелева группа $\mathcal{L}(D; A_1 \to A_2)$ является $D$-модулем. □

СЛЕДСТВИЕ 9.3.31. *Пусть $A_1$, $A_2$ - $D$-модули. Отображение*

$$(9.3.74) \qquad d\,f : A_1 \to A_2 \quad d \in D \quad f \in \mathcal{L}(D; A_1 \to A_2)$$

*определённое равенством*

$$(9.3.75) \qquad (d\,f) \circ x = d(f \circ x)$$

*называется* **произведением отображения $f$ на скаляр $d$** *и является линейным отображением. Представление*

$$(9.3.76) \qquad a : f \in \mathcal{L}(D; A_1 \to A_2) \to af \in \mathcal{L}(D; A_1 \to A_2)$$

*кольца $D$ в абелевой группе $\mathcal{L}(D; A_1 \to A_2)$ порождает структуру $D$-модуля.* □

## 9.4. Алгебра над коммутативным кольцом

ОПРЕДЕЛЕНИЕ 9.4.1. *Пусть $D$ - коммутативное кольцо. $D$-модуль $A$ называется* **алгеброй над кольцом $D$** *или $D$-алгеброй, если определена операция произведения[9.7] в $A$*

$$(9.4.1) \qquad v\,w = C \circ (v, w)$$

*где $C$ - билинейное отображение*

$$C : A \times A \to A$$

*Если $A$ является свободным $D$-модулем, то $A$ называется* **свободной алгеброй над кольцом $D$**. □

---

[9.7] Я следую определению, приведенному в [20], страница 1, [13], страница 4. Утверждение, верное для произвольного $D$-модуля, верно также для $D$-алгебры.

Теорема 9.4.2. *Пусть $D$ - коммутативное кольцо и $A$ - абелева группа. Диаграмма представлений*

$$D \xrightarrow{g_{12}} A \xrightarrow{g_{23}} A \qquad g_{12}(d) : v \to d\,v$$
$$\uparrow g_{12} \qquad g_{23}(v) : w \to C \circ (v, w)$$
$$D \qquad C \in \mathcal{L}(D; A^2 \to A)$$

*порождает структуру $D$-алгебры $A$.*

Доказательство. Структура $D$-модуля $A$ порождена эффективным представлением

$$g_{12} : D \longrightarrow A$$

кольца $D$ в абелевой группе $A$.

Лемма 9.4.3. *Пусть в $D$-модуле $A$ определена структура $D$-алгебры $A$, порождённая произведением*

$$v\,w = C \circ (v, w)$$

**Левый сдвиг** $D$-модуля $A$, *определённый равенством*

$$(9.4.2) \qquad l \circ v : w \in A \to v\,w \in A$$

*порождает представление*

$$A \xrightarrow{g_{23}} A \qquad g_{23} : v \to l \circ v$$
$$g_{23} \circ v : w \to (l \circ v) \circ w$$

$D$-*модуля $A$ в $D$-модуле $A$*

Доказательство. Согласно определениям 9.4.1 и 9.3.26, левый сдвиг $D$-модуля $A$ является линейным отображением. Согласно определению 9.3.23, отображение $l \circ v$ является эндоморфизмом $D$-модуля $A$. Равенство

$$(9.4.3) \quad (l \circ (v_1 + v_2)) \circ w = (v_1 + v_2)w = v_1 w + v_2 w = (l \circ v_1) \circ w + (l \circ v_2) \circ w$$

является следствием определения 9.3.26 и равенства (9.4.2). Согласно следствию 9.3.29, равенство

$$(9.4.4) \qquad l \circ (v_1 + v_2) = l \circ v_1 + l \circ v_2$$

является следствием равенства (9.4.3). Равенство

$$(9.4.5) \qquad (l \circ (dv)) \circ w = (dv)w = d(vw) = d((l \circ v) \circ w)$$

является следствием определения 9.3.26 и равенства (9.4.2). 9.3.29, равенство

$$(9.4.6) \qquad l \circ (dv) = d(l \circ v)$$

является следствием равенства (9.4.5). Лемма является следствием равенств (9.4.4), (9.4.6).                                                                 ⊙

Лемма 9.4.4. *Представление*

$$A \xrightarrow{g_{23}} A \qquad g_{23} : v \to l \circ v$$
$$g_{23} \circ v : w \to (l \circ v) \circ w$$

$D$-модуля $A$ в $D$-модуле $A$ определяет произведение в $D$-модуле $A$ согласно правилу

$$ab = (g_{23} \circ a) \circ b$$

Доказательство. Поскольку отображение $g_{23} \circ v$ является эндоморфизмом $D$-модуля $A$, то

(9.4.7)
$$(g_{23} \circ v)(w_1 + w_2) = (g_{23} \circ v) \circ w_1 + (g_{23} \circ v) \circ w_2$$
$$(g_{23} \circ v) \circ (dw) = d((g_{23} \circ v) \circ w)$$

Поскольку отображение $g_{23}$ является линейным отображением

$$g_{23} : A \to \mathcal{L}(D; A \to A)$$

то, согласно следствиям 9.3.29, 9.3.31,

(9.4.8) $(g_{23} \circ (v_1 + v_2)) \circ w = (g_{23} \circ v_1 + g_{23} \circ v_2)(w) = (g_{23} \circ v_1) \circ w + (g_{23} \circ v_2) \circ w$

(9.4.9) $\quad (g_{23} \circ (d\,v)) \circ w = (d\,(g_{23} \circ v)) \circ w = d\,((g_{23} \circ v) \circ w)$

Из равенств (9.4.7), (9.4.8), (9.4.9) и определения 9.3.26, следует, что отображение $g_{23}$ является билинейным отображением. Следовательно, отображение $g_{23}$ определяет произведение в $D$-модуле $A$ согласно правилу

$$ab = (g_{23} \circ a) \circ b$$

$\odot$

Теорема является следствием лем 9.4.3, 9.4.4. $\qquad\qquad\square$

Обычно, когда мы рассматриваем $D$-алгебру $A$, мы выбираем базис $\overline{\overline{e}}$ соответствующего $D$-модуля $A$. Этот выбор удобен, так как если $D$-модулЬ $A$ является свободным $D$-модулем, то разложение вектора однозначно относительно базиса $D$-модуля $A$. Это, в частности, позволяет описать операции произведения, указав структурные константы алгебры относительно заданного базиса.

В общем случае, базис $R$-модуля $A$ может оказаться множеством образующих. Например, если в векторном пространстве $H$, в котором задана алгебра кватернионов над полем действительных чисел, рассмотреть базис

(9.4.10)
$$e_0 = 1 \quad e_1 = i \quad e_2 = j \quad e_3 = k$$

то в алгебре $H$ верно равенство

(9.4.11)
$$e_0 = -e_1 e_1 = -e_2 e_2$$
$$e_3 = e_1 e_2$$

Следовательно, множество $(e_1, e_2)$ является базисом алгебры $H$. Следствием равенства (9.4.11) является неоднозначность представления кватерниона относительно заданного базиса. А именно, кватернион $a \in H$ можно записать в виде

$$a = (a^0 - a^4) e_1 e_1 + a^4 e_2 e_2 + a^1 e_1 + a^2 e_2 + a^3 e_1 e_2$$

где $a^4$ - произвольно.

## 9.5. Левый модуль над алгеброй

Определение 9.5.1. *Эффективное левостороннее представление*

$$(9.5.1) \qquad f : A \dashrightarrow V \quad f(a) : v \in V \to av \in V \quad a \in A$$

*ассоциативной D-алгебры A в D-модуле V называется* **левым модулем** *над D-алгеброй A. Мы также будем говорить, что D-модуль V является* **левым A-модулем** *или A∗-модулем. V-число называется* **вектором**. □

Определение 9.5.2. *Пусть A - алгебра с делением. Эффективное левостороннее представление*

$$f : A \dashrightarrow V \quad f(a) : v \in V \to av \in V \quad a \in A$$

*абелевой группы A в D-модуле V называется* **левым векторным пространством** *над D-алгеброй A. Мы также будем говорить, что D-модуль V является* **левым A-векторным пространством** *или A∗-векторным пространством. V-число называется* **вектором**. □

Теорема 9.5.3. *Следующая диаграмма представлений описывает левый A-модуль V*

$$(9.5.2)$$

$$
\begin{aligned}
g_{12}(d) &: a \to d\,a \\
g_{23}(v) &: w \to C(w, v) \\
C &\in \mathcal{L}(A^2 \to A) \\
g_{3,4}(a) &: v \to v\,a \\
g_{1,4}(d) &: v \to d\,v
\end{aligned}
$$

*В диаграмме представлений* (9.5.2) *верна* **коммутативность представлений** *коммутативного кольца D и D-алгебры A в абелевой группе V*

$$(9.5.3) \qquad a(dv) = d(av)$$

Доказательство. Диаграмма представлений (9.5.2) является следствием определения 9.5.1 и теоремы 9.4.2. Равенство (9.5.3) является следствием утверждения, что левостороннее преобразование $g_{3,4}(a)$ является эндоморфизмом D-модуля V. □

Теорема 9.5.4. *Пусть g - эффективное левостороннее представление D-алгебры A в D-модуле V. Тогда D-алгебра A ассоциативна.*

Доказательство. Пусть $a$, $b$, $c \in A$, $v \in V$. Равенство

$$(9.5.4) \qquad (ab)v = a(bv)$$

является следствием утверждения, что левостороннее представление g является левосторонним представлением мультипликативной группы D-алгебры A. Равенство

$$(9.5.5) \qquad a(b(cv)) = a((bc)v) = (a(bc))v$$

является следствием равенства (9.5.4). Так как $cv \in A$, равенство

$$(9.5.6) \qquad a(b(cv)) = (ab)(cv) = ((ab)c)v$$

является следствием равенства (9.5.4). Равенство

$$(9.5.7) \qquad (a(bc))v = ((ab)c)v$$

является следствием равенств (9.5.5), (9.5.7). Поскольку $v$ - произвольный вектор $A$-модуля $V$, равенство

$$(9.5.8) \qquad a(bc) = (ab)c$$

является следствием равенства (9.5.7). Следовательно, $D$-алгебра $A$ ассоциативна. $\qquad\square$

Теорема 9.5.5. *Пусть $V$ является левым $A$-модулем. Для любого вектора $v \in V$, вектор, порождённый диаграммой представлений (9.5.2), имеет следующий вид*

$$(9.5.9) \qquad (a + n)v = av + nv \quad a \in A \quad n \in D$$

9.5.5.1: *Множество отображений*

$$(9.5.10) \qquad a + n : v \in V \to (a + n)v \in V$$

*порождает[9.8] $D$-алгебру $A_{(1)}$ где сложение определено равенством*

$$(9.5.11) \qquad (a + n) + (b + m) = (a + b) + (n + m)$$

*и произведение определено равенством*

$$(9.5.12) \qquad (a + n)(b + m) = (ab + ma + nb) + (nm)$$

*$D$-алгебра $A_{(1)}$ называется **унитальным расширением** $D$-алгебры $A$.*

| | | |
|---|---|---|
| *Если $D$-алгебра $A$ имеет единицу, то* | $D \subseteq A$ | $A_{(1)} = A$ |
| *Если $D$-алгебра $A$ является идеалом $D$, то* | $A \subseteq D$ | $A_{(1)} = D$ |
| *В противном случае* | | $A_{(1)} = A \oplus D$ |

9.5.5.2: *$D$-алгебра $A$ является левым идеалом $D$-алгебры $A_{(1)}$.*

9.5.5.3: *Множество преобразований (9.5.9) порождает левостороннее представление $D$-алгебры $A_{(1)}$ в абелевой группе $V$.*

Мы будем пользоваться обозначением $A_{(1)}v$ для множества векторов, порождённых вектором $v$.

Теорема 9.5.6. *Элементы левого $A$-модуля $V$ удовлетворяют соотношениям*

9.5.6.1: **закон ассоциативности**

$$(9.5.13) \qquad (pq)v = p(qv)$$

9.5.6.2: **закон дистрибутивности**

$$(9.5.14) \qquad p(v + w) = pv + pw$$

$$(9.5.15) \qquad (p + q)v = pv + qv$$

9.5.6.3: **закон унитарности**

$$(9.5.16) \qquad 1v = v$$

---

[9.8] Смотри определение унитального расширения также на страницах [6]-52, [7]-64.

*для любых   p, q ∈ A₍₁₎, v, w ∈ V.*

Доказательство теорем 9.5.5, 9.5.6. Пусть $v \in V$.

Лемма 9.5.7. *Пусть   $d \in D$,   $a \in A$.   Отображение* (9.5.10) *является эндоморфизмом абелевой группы* $V$.

Доказательство. Утверждения $dv \in V$,   $av \in V$   являются следствием теорем 6.1.4, 9.5.3. Так как $V$ является абелевой группой, то

$$dv + av \in V \quad d \in D \quad a \in A$$

Следовательно, для любого $D$-числа $d$ и любого $A$-числа $a$, мы определили отображение (9.5.10). Поскольку преобразование $g_{1,4}(d)$ и левостороннее преобразование $g_{3,4}(a)$ являются эндоморфизмами абелевой группы $V$, то отображение (9.5.10) является эндоморфизмом абелевой группы $V$.                ⊙

Пусть $A₍₁₎$ - множество отображений (9.5.10). Равенство (9.5.14) является следствием леммы 9.5.7.

Пусть  $p = a + n \in A₍₁₎$, $q = b + m \in A₍₁₎$.  Согласно утверждению 9.3.3.3, мы определим сумму $A₍₁₎$-чисел $p$ и $q$ равенством (9.5.15). Равенство

$$(9.5.17) \qquad ((a + n) + (b + m))v = (a + n)v + (b + m)v$$

является следствием равенства (9.5.15). Равенство

$$(9.5.18) \qquad (n + m)v = cn + dm$$

является следствием утверждения, что представление $g_{1,4}$ является гомоморфизмом аддитивной группы кольца $D$. Равенство

$$(9.5.19) \qquad (a + b)v = av + bv$$

является следствием утверждения, что левостороннее представление $g_{3,4}$ является гомоморфизмом аддитивной группы $D$-алгебры $A$. Так как $V$ является абелевой группой, то равенство

$$(9.5.20) \quad \begin{aligned} ((a + n) + (b + m))v &= av + nv + bv + mv = av + bv + nv + mv \\ &= (a + b)v + (n + m)v = ((a + b) + (n + m))v \end{aligned}$$

является следствием равенств (9.5.17), (9.5.18), (9.5.19). Из равенства (9.5.20) следует, что определение (9.5.11) суммы на множестве $A₍₁₎$ не зависит от вектора $v$.

Равенства (9.5.13), (9.5.16) являются следствием утверждения 9.5.5.3. Пусть $p = a + n \in A₍₁₎$, $q = b + m \in A₍₁₎$.  Так как произведение в $D$-алгебре $A$ может быть неассоциативным, то опираясь на теорему 9.5.6, мы рассмотрим произведение $A₍₁₎$-чисел $p$ и $q$ как билинейное отображение

$$f : A₍₁₎ \times A₍₁₎ \to A₍₁₎$$

такое, что верны равенства

$$(9.5.21) \qquad f(a, b) = ab \quad a, b \in A$$

$$(9.5.22) \qquad f(1, p) = f(p, 1) = p \quad p \in A₍₁₎ \quad 1 \in D₍₁₎$$

Равенство

$$(9.5.23) \quad \begin{aligned} (a+n)(b+m) &= f(a+n, b+m) \\ &= f(a,b) + f(a,m) + f(n,b) + f(n,m) \\ &= f(a,b) + mf(a,1) + nf(1,b) + nf(1,m) \\ &= ab + ma + nb + nm \end{aligned}$$

является следствием равенств (9.5.21), (9.5.22). Равенство (9.5.12) является следствием равенства (9.5.23).

Утверждение 9.5.5.2 является следствием равенства (9.5.12). □

Билинейное отображение

$$(a,v) \in A \times V \to av \in V$$

порождённое левосторонним представлением $g_{2,3}$ называется **левосторонним произведением** вектора на скаляр.

Теорема 9.5.8. *Пусть $V$ - левый $A$-модуль. Множество векторов, порождённое множеством векторов $v = (v_i \in V, i \in I)$, имеет вид*[9.9]

$$(9.5.24) \quad J(v) = \left\{ w : w = \sum_{i \in I} c^i v_i, c^i \in A_{(1)}, |\{i : c^i \neq 0\}| < \infty \right\}$$

Доказательство. Мы докажем теорему по индукции, опираясь на теорему 6.1.4, Согласно теореме 6.1.4, мы должны доказать следующие утверждения:

9.5.8.1: $v_k \in X_0 \subseteq J(v)$
9.5.8.2: $c^k v_k \in J(v)$, $c^k \in A_{(1)}$, $k \in I$
9.5.8.3: $\sum_{k \in I} c^k v_k \in J(v)$, $c^k \in A_{(1)}$, $|\{i : c^i \neq 0\}| < \infty$
9.5.8.4: $w_1, w_2 \in J(v) \Rightarrow w_1 + w_2 \in J(v)$
9.5.8.5: $a \in A$, $w \in J(v) \Rightarrow aw \in J(v)$

- Для произвольного $v_k \in v$, положим $c^i = \delta_k^i \in A_{(1)}$. Тогда

$$(9.5.25) \quad v_k = \sum_{i \in I} c^i v_i$$

  Утверждение 9.5.8.1 следует из (9.5.24), (9.5.25).
- Утверждение 9.5.8.2 являются следствием теорем 6.1.4, 9.5.5 и утверждения 9.5.8.1.
- Так как $V$ является абелевой группой, то утверждение 9.5.8.3 следует из утверждения 9.5.8.2 и теорем 6.1.4, 9.2.3.
- Пусть $w_1, w_2 \in X_k \subseteq J(v)$. Так как $V$ является абелевой группой, то, согласно утверждению 6.1.4.3,

$$(9.5.26) \quad w_1 + w_2 \in X_{k+1}$$

  Согласно равенству (9.5.24), существуют $A_{(1)}$-числа $w_1^i$, $w_2^i$, $i \in I$, такие, что

$$(9.5.27) \quad w_1 = \sum_{i \in I} w_1^i v_i \quad w_2 = \sum_{i \in I} w_2^i v_i$$

---

[9.9] Для множества $A$, мы обозначим $|A|$ мощность множества $A$. Запись $|A| < \infty$ означает, что множество $A$ конечно.

где множества

$$(9.5.28) \qquad H_1 = \{i \in I : w_1^i \neq 0\} \quad H_2 = \{i \in I : w_2^i \neq 0\}$$

конечны. Так как $V$ является абелевой группой, то из равенства (9.5.27) следует, что

$$(9.5.29) \qquad w_1 + w_2 = \sum_{i \in I} w_1^i v_i + \sum_{i \in I} w_2^i v_i = \sum_{i \in I}(w_1^i v_i + w_2^i v_i)$$

Равенство

$$(9.5.30) \qquad w_1 + w_2 = \sum_{i \in I}(w_1^i + w_2^i)v_i$$

является следствием равенств (9.5.15), (9.5.29). Из равенства (9.5.28) следует, что множество

$$\{i \in I : w_1^i + w_2^i \neq 0\} \subseteq H_1 \cup H_2$$

конечно.

- Пусть $w \in X_k \subseteq J(v)$. Согласно утверждению 6.1.4.4, для любого $A_{(1)}$-числа $a$,

$$(9.5.31) \qquad aw \in X_{k+1}$$

Согласно равенству (9.5.24), существуют $A_{(1)}$-числа $w^i$, $i \in I$, такие, что

$$(9.5.32) \qquad w = \sum_{i \in I} w^i v_i$$

где

$$(9.5.33) \qquad |\{i \in I : w^i \neq 0\}| < \infty$$

Из равенства (9.5.32) следует, что

$$(9.5.34) \qquad aw = a \sum_{i \in I} w^i v_i = \sum_{i \in I} a(w^i v_i) = \sum_{i \in I}(aw^i)v_i$$

Из утверждения (9.5.33) следует, что множество $\{i \in I : aw^i \neq 0\}$ конечно.

Из равенств (9.5.26), (9.5.30), (9.5.31), (9.5.34) следует, что $X_{k+1} \subseteq J(v)$. $\square$

Определение 9.5.9. *Пусть $v = (v_i \in V, i \in I)$ - множество векторов. Выражение $w^i v_i$ называется* **линейной комбинацией** *векторов $v_i$. Вектор $\overline{w} = w^i v_i$ называется* **линейно зависимым** *от векторов $v_i$.* $\square$

Представим множество $A_{(1)}$-чисел $w^i$, $i \in I$, в виде матрицы

$$w = \begin{pmatrix} w^1 \\ ... \\ w^n \end{pmatrix}$$

Представим множество векторов $v_i$, $i \in I$, в виде матрицы

$$v = \begin{pmatrix} v_1 & ... & v_n \end{pmatrix}$$

Тогда мы можем записать линейную комбинацию векторов $\overline{w} = w^i v_i$ в виде

$$\overline{w} = w^* {}_* v$$

Теорема 9.5.10. *Пусть $A$ - ассоциативная $D$-алгебра с делением. Если уравнение*

$$w^i v_i = 0$$

*предполагает существования индекса $i = j$ такого, что $w^j \neq 0$, то вектор $v_j$ линейно зависит от остальных векторов $v$.*

Доказательство. Теорема является следствием равенства

$$v_j = \sum_{i \in I \setminus \{j\}} (w^j)^{-1} w^i v_i$$

и определения 9.5.9. $\qquad\square$

Очевидно, что для любого множества векторов $v_i$,

$$w^i = 0 \Rightarrow w^*_* v = 0$$

Определение 9.5.11. *Множество векторов [9.10] $v_i$, $i \in I$, левого $A$-модуля $V$ линейно независимо, если $w = 0$ следует из уравнения*

$$w^i v_i = 0$$

*В противном случае, множество векторов $v_i$, $i \in I$, линейно зависимо.* $\qquad\square$

Следующее определение является следствием теорем 9.5.8, 6.1.4 и определения 6.1.5.

Определение 9.5.12. *$J(v)$ называется подмодулем, порождённым множеством $v$, а $v$ - множеством образующих подмодуля $J(v)$. В частности, множеством образующих левого $D$-модуля $V$ будет такое подмножество $X \subset V$, что $J(X) = V$.* $\qquad\square$

Следующее определение является следствием теорем 9.5.8, 6.1.4 и определения 6.2.6.

Определение 9.5.13. *Если множество $X \subset V$ является множеством образующих левого $D$-модуля $V$, то любое множество $Y$, $X \subset Y \subset V$ также является множеством образующих левого $D$-модуля $V$. Если существует минимальное множество $X$, порождающее левый $D$-модуль $V$, то такое множество $X$ называется базисом левого $D$-модуля $V$.* $\qquad\square$

Теорема 9.5.14. *Множество векторов $\overline{\overline{e}} = (e_i, i \in I)$ является базисом левого $A$-модуля $V$, если верны следующие утверждения.*

9.5.14.1: *Произвольный вектор $v \in V$ является линейной комбинацией векторов множества $\overline{\overline{e}}$.*

9.5.14.2: *Вектор $e_i$ нельзя представить в виде линейной комбинации остальных векторов множества $\overline{\overline{e}}$.*

Доказательство. Согласно утверждению 9.5.14.1, теореме 9.5.8 и определению 9.5.9, множество $\overline{\overline{e}}$ порождает левый $A$-модуль $V$ (определение 9.5.12). Согласно утверждению 9.5.14.2, множество $\overline{\overline{e}}$ является минимальным множеством, порождающим левый $A$-модуль $V$. Согласно определению 9.5.13, множество $\overline{\overline{e}}$ является базисом левого $A$-модуля $V$. $\qquad\square$

---

[9.10] Я следую определению в [2], страница 100.

**Теорема 9.5.15.** *Пусть $A$ - ассоциативная $D$-алгебра с делением. Множество векторов $\overline{\overline{e}} = (e_i, i \in I)$ является **базисом левого $A$-векторного пространства** $V$, если векторы $e_i$ линейно независимы и любой вектор $v \in V$ линейно зависит от векторов $e_i$.*

Доказательство. Пусть множество векторов $e_i, i \in I$, линейно зависимо. Тогда в равенстве

$$w^i e_i = 0$$

существует индекс $i = j$ такой, что $w^j \neq 0$. Согласно теореме 9.5.10, вектор $e_j$ линейно зависит от остальных векторов множества $\overline{\overline{e}}$. Согласно определению 9.5.13, множество векторов $e_i, i \in I$, не является базисом левого $A$-векторного пространства $V$.

Следовательно, если множество векторов $e_i, i \in I$, является базисом, то эти векторы линейно независимы. Так как произвольный вектор $v \in V$ является линейной комбинацией векторов $e_i, i \in I$, , то множество векторов $v, e_i, i \in I$, не является линейно независимым. □

**Определение 9.5.16.** *Пусть $\overline{\overline{e}}$ - базис левого $A$-модуля $V$, и вектор $\overline{v} \in V$ имеет разложение*

$$\overline{v} = v^*_* e = v^i e_i$$

*относительно базиса $\overline{\overline{e}}$. $A_{(1)}$-числа $v^i$ называются **координатами** вектора $\overline{v}$ относительно базиса $\overline{\overline{e}}$. Матрица $A_{(1)}$-чисел $v = (v^i, i \in I)$ называется **координатной матрицей вектора** $\overline{v}$ в базисе $\overline{\overline{e}}$.* □

**Теорема 9.5.17.** *Пусть $A$ - ассоциативная $D$-алгебра. Пусть $\overline{\overline{e}}$ - базис левого $A$-модуля $V$. Пусть*

$$(9.5.35) \qquad\qquad w^i e_i = 0$$

*линейная зависимость векторов базиса $\overline{\overline{e}}$. Тогда*

9.5.17.1: *$A_{(1)}$-число $w^i$, $i \in I$, не имеет обратного элемента в $D$-алгебре $A_{(1)}$.*

9.5.17.2: *Множество $A'$ матриц $w = (w^i, i \in I)$ порождает левый $A$-модуль.*

Доказательство. Допустим существует матрица $w = (w^i, i \in I)$ такая, что равенство (9.5.35) верно и существует индекс $i = j$ такой, что $w^j \neq 0$. Если мы положим, что $A_{(1)}$-число $c^j$ имеет обратный, то равенство

$$e_j = \sum_{i \in I \setminus \{j\}} (w^j)^{-1} w^i e_i$$

является следствием равенства (9.5.35). Следовательно вектор $e_j$ является линейной комбинацией остальных векторов множества $\overline{\overline{e}}$ и множество $\overline{\overline{e}}$ не является базисом. Следовательно, наше предположение неверно, и $A_{(1)}$-число $c^j$ не имеет обратного.

Пусть матрицы $b = (b^i, i \in I) \in A'$, $c = (c^i, i \in I) \in A'$. Из равенств

$$b^i e_i = 0$$
$$c^i e_i = 0$$

следует

$$(b^i + c^i) e_i = 0$$

Следовательно, множество $A'$ является абелевой группой.

Пусть матрица $c = (c^i, i \in I) \in A'$ и $a \in A$. Из равенства

$$c^i e_i = 0$$

следует

$$(ac^i)e_i = 0$$

Следовательно, абелева группа $A'$ является левым $A$-модулем. $\qquad\square$

**Теорема 9.5.18.** *Пусть левый $A$-модуль $V$ имеет базис $\overline{\overline{e}}$ такой, что в равенстве*

$$(9.5.36) \qquad\qquad w^i e_i = 0$$

*существует индекс $i = j$ такой, что $w^j \neq 0$. Тогда*

9.5.18.1: *Матрица $w = (w^i, i \in I)$ определяет координаты вектора $0 \in V$ относительно базиса $\overline{\overline{e}}$.*

9.5.18.2: *Координаты вектора $\overline{v}$ относительно базиса $\overline{\overline{e}}$ определены однозначно с точностью до выбора координат вектора $0 \in V$.*

Доказательство. Утверждение 9.5.18.1 является следствием равенства (9.5.36) и определения 9.5.16.

Пусть вектор $\overline{v}$ имеет разложение

$$(9.5.37) \qquad\qquad \overline{v} = v^*{}_* e = v^i e_i$$

относительно базиса $\overline{\overline{e}}$. Равенство

$$(9.5.38) \qquad\qquad \overline{v} = \overline{v} + 0 = v^i e_i + c^i e_i = (v^i + c^i) e_i$$

является следствием равенств (9.5.36), (9.5.37). Утверждение 9.5.18.2 является следствием равенств (9.5.37), (9.5.38) и определения 9.5.16. $\qquad\square$

**Определение 9.5.19.** *Левый $A$-модуль $V$ - **свободный левый $A$-модуль**,* [9.11] *если левый $A$-модуль $V$ имеет базис и векторы базиса линейно независимы.* $\qquad\square$

**Теорема 9.5.20.** *Координаты вектора $v \in V$ относительно базиса $\overline{\overline{e}}$ свободного левого $A$-модуля $V$ определены однозначно.*

Доказательство. Теорема является следствием теоремы 9.5.18 и определений 9.5.11, 9.5.19. $\qquad\square$

## 9.6. Правый модуль над алгеброй

**Определение 9.6.1.** *Эффективное правостороннее представление*

$$(9.6.1) \qquad f : A \dashrightarrow V \quad f(a) : v \in V \to va \in V \quad a \in A$$

*ассоциативной $D$-алгебры $A$ в $D$-модуле $V$ называется **правым модулем над $D$-алгеброй $A$**. Мы также будем говорить, что $D$-модуль $V$ является **правым $A$-модулем** или $*A$-модулем. $V$-число называется **вектором**.* $\qquad\square$

---

[9.11] Я следую определению в [2], страница 103.

Определение 9.6.2. *Пусть $A$ - алгебра с делением. Эффективное право-*
*стороннее представление*

$$f : A \xrightarrow{\ \ *\ \ } V \quad f(a) : v \in V \to va \in V \quad a \in A$$

*абелевой группы $A$ в $D$-модуле $V$ называется* **правым векторным про-**
**странством** *над $D$-алгеброй $A$. Мы также будем говорить, что $D$-модуль $V$*
*является* **правым $A$-векторным пространством** *или $*A$-векторным про-*
**странством**. *$V$-число называется* **вектором**.                    □

Теорема 9.6.3. *Следующая диаграмма представлений описывает правый*
*$A$-модуль $V$*

$$g_{12}(d) : a \to d\,a$$

(9.6.2)   

$$g_{23}(v) : w \to C(w, v)$$

$$C \in \mathcal{L}(A^2 \to A)$$

$$g_{3,4}(a) : v \to v\,a$$

$$g_{1,4}(d) : v \to v\,d$$

*В диаграмме представлений* (9.6.2) *верна* **коммутативность представле-**
**ний** *коммутативного кольца $D$ и $D$-алгебры $A$ в абелевой группе $V$*

$$(9.6.3) \qquad\qquad (vd)a = (va)d$$

Доказательство. Диаграмма представлений (9.6.2) является следстви-
ем определения 9.6.1 и теоремы 9.4.2. Равенство (9.6.3) является следствием
утверждения, что правостороннее преобразование $g_{3,4}(a)$ является эндомор-
физмом $D$-модуля $V$.                    □

Теорема 9.6.4. *Пусть $g$ - эффективное левостороннее представление $D$-*
*алгебры $A$ в $D$-модуле $V$. Тогда $D$-алгебра $A$ ассоциативна.*

Доказательство. Пусть $a$, $b$, $c \in A$, $v \in V$. Равенство

$$(9.6.4) \qquad\qquad v(ab) = (va)b$$

является следствием утверждения, что правостороннее представление $g$ явля-
ется правосторонним представлением мультипликативной группы $D$-алгебры
$A$. Равенство

$$(9.6.5) \qquad\qquad ((vc)b)a = (v(cb))a = v((cb)a)$$

является следствием равенства (9.6.4). Так как $vc \in A$, равенство

$$(9.6.6) \qquad\qquad ((vc)b)a = (vc)(ba) = v(c(ba))$$

является следствием равенства (9.6.4). Равенство

$$(9.6.7) \qquad\qquad v((cb)a) = v(c(ba))$$

является следствием равенств (9.6.5), (9.6.7). Поскольку $v$ - произвольный век-
тор $A$-модуля $V$, равенство

$$(9.6.8) \qquad\qquad (cb)a = c(ba)$$

является следствием равенства (9.6.7). Следовательно, $D$-алгебра $A$ ассоциа-
тивна.                    □

Теорема 9.6.5. *Пусть $V$ является правым $A$-модулем. Для любого вектора $v \in V$, вектор, порождённый диаграммой представлений* (9.6.2), *имеет следующий вид*

$$(9.6.9) \qquad v(a + n) = va + vn \quad a \in A \quad n \in D$$

9.6.5.1: *Множество отображений*

$$(9.6.10) \qquad a + n : v \in V \to v(a + n) \in V$$

*порождает*[9.12] *$D$-алгебру $A_{(1)}$ где сложение определено равенством*

$$(9.6.11) \qquad (a + n) + (b + m) = (a + b) + (n + m)$$

*и произведение определено равенством*

$$(9.6.12) \qquad (a + n)(b + m) = (ab + ma + nb) + (nm)$$

*$D$-алгебра $A_{(1)}$ называется* **унитальным расширением** *$D$-алгебры $A$.*

| | | |
|---|---|---|
| *Если $D$-алгебра $A$ имеет единицу, то* | $D \subseteq A$ | $A_{(1)} = A$ |
| *Если $D$-алгебра $A$ является идеалом $D$, то* | $A \subseteq D$ | $A_{(1)} = D$ |
| *В противном случае* | $A_{(1)} = A \oplus D$ | |

9.6.5.2: *$D$-алгебра $A$ является правым идеалом $D$-алгебры $A_{(1)}$.*

9.6.5.3: *Множество преобразований* (9.6.9) *порождает правостороннее представление $D$-алгебры $A_{(1)}$ в абелевой группе $V$.*

Мы будем пользоваться обозначением $A_{(1)}v$ для множества векторов, порождённых вектором $v$.

Теорема 9.6.6. *Элементы правого $A$-модуля $V$ удовлетворяют соотношениям*

9.6.6.1: **закон ассоциативности**

$$(9.6.13) \qquad v(pq) = (vp)q$$

9.6.6.2: **закон дистрибутивности**

$$(9.6.14) \qquad (v + w)p = vp + wp$$

$$(9.6.15) \qquad v(p + q) = vp + vq$$

9.6.6.3: **закон унитарности**

$$(9.6.16) \qquad v1 = v$$

*для любых $p, q \in A_{(1)}$, $v, w \in V$.*

Доказательство теорем 9.6.5, 9.6.6. Пусть $v \in V$.

Лемма 9.6.7. *Пусть $d \in D$, $a \in A$. Отображение* (9.6.10) *является эндоморфизмом абелевой группы $V$.*

---

[9.12] Смотри определение унитального расширения также на страницах [6]-52, [7]-64.

Доказательство. Утверждения $vd \in V$, $va \in V$ являются следствием теорем 6.1.4, 9.6.3. Так как $V$ является абелевой группой, то

$$vd + va \in V \quad d \in D \quad a \in A$$

Следовательно, для любого $D$-числа $d$ и любого $A$-числа $a$, мы определили отображение (9.6.10). Поскольку преобразование $g_{1,4}(d)$ и правостороннее преобразование $g_{3,4}(a)$ являются эндоморфизмами абелевой группы $V$, то отображение (9.6.10) является эндоморфизмом абелевой группы $V$.                    ⊙

Пусть $A_{(1)}$ - множество отображений (9.6.10). Равенство (9.6.14) является следствием леммы 9.6.7.

Пусть $p = a + n \in A_{(1)}$, $q = b + m \in A_{(1)}$. Согласно утверждению 9.3.3.3, мы определим сумму $A_{(1)}$-чисел $p$ и $q$ равенством (9.6.15). Равенство

$$(9.6.17) \qquad v((a+n) + (b+m)) = v(a+n) + v(b+m)$$

является следствием равенства (9.6.15). Равенство

$$(9.6.18) \qquad v(n+m) = vn + vm$$

является следствием утверждения, что представление $g_{1,4}$ является гомоморфизмом аддитивной группы кольца $D$. Равенство

$$(9.6.19) \qquad v(a+b) = va + vb$$

является следствием утверждения, что правостороннее представление $g_{3,4}$ является гомоморфизмом аддитивной группы $D$-алгебры $A$. Так как $V$ является абелевой группой, то равенство

$$(9.6.20) \qquad \begin{aligned} v((a+n)+(b+m)) &= va + vn + vb + vm = va + vb + vn + vm \\ &= v(a+b) + v(n+m) = v((a+b)+(n+m)) \end{aligned}$$

является следствием равенств (9.6.17), (9.6.18), (9.6.19). Из равенства (9.6.20) следует, что определение (9.6.11) суммы на множестве $A_{(1)}$ не зависит от вектора $v$.

Равенства (9.6.13), (9.6.16) являются следствием утверждения 9.6.5.3. Пусть $p = a + n \in A_{(1)}$, $q = b + m \in A_{(1)}$. Так как произведение в $D$-алгебре $A$ может быть неассоциативным, то опираясь на теорему 9.6.6, мы рассмотрим произведение $A_{(1)}$-чисел $p$ и $q$ как билинейное отображение

$$f : A_{(1)} \times A_{(1)} \to A_{(1)}$$

такое, что верны равенства

$$(9.6.21) \qquad f(a,b) = ab \quad a,b \in A$$

$$(9.6.22) \qquad f(1,p) = f(p,1) = p \quad p \in A_{(1)} \quad 1 \in D_{(1)}$$

Равенство

$$(9.6.23) \qquad \begin{aligned} (a+n)(b+m) &= f(a+n, b+m) \\ &= f(a,b) + f(a,m) + f(n,b) + f(n,m) \\ &= f(a,b) + mf(a,1) + nf(1,b) + nf(1,m) \\ &= ab + ma + nb + nm \end{aligned}$$

является следствием равенств (9.6.21), (9.6.22). Равенство (9.6.12) является следствием равенства (9.6.23).

Утверждение 9.6.5.2 является следствием равенства (9.6.12).                    □

Билинейное отображение

$$(v, a) \in V \times A \to va \in V$$

порождённое правосторонним представлением $g_{2,3}$ называется **правосторонним произведением** вектора на скаляр.

**Теорема 9.6.8.** *Пусть $V$ - правый $A$-модуль. Множество векторов, порождённое множеством векторов $v = (v_i \in V, i \in I)$, имеет вид*[9.13]

$$(9.6.24) \qquad J(v) = \left\{ w : w = \sum_{i \in I} v_i c^i, c^i \in A_{(1)}, |\{i : c^i \neq 0\}| < \infty \right\}$$

**Доказательство.** Мы докажем теорему по индукции, опираясь на теорему 6.1.4. Согласно теореме 6.1.4, мы должны доказать следующие утверждения:

9.6.8.1: $v_k \in X_0 \subseteq J(v)$
9.6.8.2: $v_k c^k \in J(v),\ c^k \in A_{(1)},\ k \in I$
9.6.8.3: $\sum_{k \in I} v_k c^k \in J(v),\ c^k \in A_{(1)},\ |\{i : c^i \neq 0\}| < \infty$
9.6.8.4: $w_1,\ w_2 \in J(v) \Rightarrow w_1 + w_2 \in J(v)$
9.6.8.5: $a \in A,\ w \in J(v) \Rightarrow aw \in J(v)$

- Для произвольного $v_k \in v$, положим $c^i = \delta_k^i \in A_{(1)}$. Тогда

$$(9.6.25) \qquad v_k = \sum_{i \in I} v_i c^i$$

Утверждение 9.6.8.1 следует из (9.6.24), (9.6.25).
- Утверждение 9.6.8.2 являются следствием теорем 6.1.4, 9.6.5 и утверждения 9.6.8.1.
- Так как $V$ является абелевой группой, то утверждение 9.6.8.3 следует из утверждения 9.6.8.2 и теорем 6.1.4, 9.2.3.
- Пусть $w_1,\ w_2 \in X_k \subseteq J(v)$. Так как $V$ является абелевой группой, то, согласно утверждению 6.1.4.3,

$$(9.6.26) \qquad w_1 + w_2 \in X_{k+1}$$

Согласно равенству (9.6.24), существуют $A_{(1)}$-числа $w_1^i,\ w_2^i,\ i \in I$, такие, что

$$(9.6.27) \qquad w_1 = \sum_{i \in I} v_i w_1^i \qquad w_2 = \sum_{i \in I} v_i w_2^i$$

где множества

$$(9.6.28) \qquad H_1 = \{i \in I : w_1^i \neq 0\} \qquad H_2 = \{i \in I : w_2^i \neq 0\}$$

конечны. Так как $V$ является абелевой группой, то из равенства (9.6.27) следует, что

$$(9.6.29) \qquad w_1 + w_2 = \sum_{i \in I} v_i w_1^i + \sum_{i \in I} v_i w_2^i = \sum_{i \in I} (v_i w_1^i + v_i w_2^i)$$

---

[9.13] Для множества $A$, мы обозначим $|A|$ мощность множества $A$. Запись $|A| < \infty$ означает, что множество $A$ конечно.

Равенство

(9.6.30)
$$w_1 + w_2 = \sum_{i \in I} v_i(w_1^i + w_2^i)$$

является следствием равенств (9.6.15), (9.6.29). Из равенства (9.6.28) следует, что множество

$$\{i \in I : w_1^i + w_2^i \neq 0\} \subseteq H_1 \cup H_2$$

конечно.

- Пусть $w \in X_k \subseteq J(v)$. Согласно утверждению 6.1.4.4, для любого $A_{(1)}$-числа $a$,

(9.6.31)
$$wa \in X_{k+1}$$

Согласно равенству (9.6.24), существуют $A_{(1)}$-числа $w^i$, $i \in I$, такие, что

(9.6.32)
$$w = \sum_{i \in I} v_i w^i$$

где

(9.6.33)
$$|\{i \in I : w^i \neq 0\}| < \infty$$

Из равенства (9.6.32) следует, что

(9.6.34)
$$wa = \left(\sum_{i \in I} v_i w^i\right) a = \sum_{i \in I} (v_i w^i) a = \sum_{i \in I} (v_i w^i a)$$

Из утверждения (9.6.33) следует, что множество $\{i \in I : w^i a \neq 0\}$ конечно.

Из равенств (9.6.26), (9.6.30), (9.6.31), (9.6.34) следует, что $X_{k+1} \subseteq J(v)$.     $\square$

ОПРЕДЕЛЕНИЕ 9.6.9. *Пусть* $v = (v_i \in V, i \in I)$ *- множество векторов. Выражение* $v_i w^i$ *называется* **линейной комбинацией** *векторов* $v_i$. *Вектор* $\overline{w} = v_i w^i$ *называется* **линейно зависимым** *от векторов* $v_i$.     $\square$

Представим множество $A_{(1)}$-чисел $w^i$, $i \in I$, в виде матрицы

$$w = \begin{pmatrix} w^1 \\ ... \\ w^n \end{pmatrix}$$

Представим множество векторов $v_i$, $i \in I$, в виде матрицы

$$v = \begin{pmatrix} v_1 & ... & v_n \end{pmatrix}$$

Тогда мы можем записать линейную комбинацию векторов $\overline{w} = v_i w^i$ в виде

$$\overline{w} = v_*{}^* w$$

ТЕОРЕМА 9.6.10. *Пусть* $A$ *- ассоциативная* $D$-*алгебра с делением. Если уравнение*

$$v_i w^i = 0$$

*предполагает существования индекса* $i = j$ *такого, что* $w^j \neq 0$, *то вектор* $v_j$ *линейно зависит от остальных векторов* $v$.

Доказательство. Теорема является следствием равенства

$$v_j = \sum_{i \in I \setminus \{j\}} v_i w^i (w^j)^{-1}$$

и определения 9.6.9.                                                      □

Очевидно, что для любого множества векторов $v_i$,

$$w^i = 0 \Rightarrow v_* {}^* w = 0$$

Определение 9.6.11. *Множество векторов*[9.14] $v_i$, $i \in I$, *правого A-модуля* $V$ **линейно независимо**, *если* $w = 0$ *следует из уравнения*

$$v_i w^i = 0$$

*В противном случае, множество векторов* $v_i$, $i \in I$, **линейно зависимо**.
                                                                         □

Следующее определение является следствием теорем 9.6.8, 6.1.4 и определения 6.1.5.

Определение 9.6.12. $J(v)$ *называется* **подмодулем, порождённым множеством** $v$, *а* $v$ - **множеством образующих** *подмодуля* $J(v)$. *В частности,* **множество образующих** *правого D-модуля* $V$ *будет такое подмножество* $X \subset V$, *что* $J(X) = V$.                                      □

Следующее определение является следствием теорем 9.6.8, 6.1.4 и определения 6.2.6.

Определение 9.6.13. *Если множество* $X \subset V$ *является множеством образующих правого D-модуля* $V$, *то любое множество* $Y$, $X \subset Y \subset V$ *также является множеством образующих правого D-модуля* $V$. *Если существует минимальное множество* $X$, *порождающее правый D-модуль* $V$, *то такое множество* $X$ *называется* **базисом** *правого D-модуля* $V$.           □

Теорема 9.6.14. *Множество векторов* $\overline{\overline{e}} = (e_i, i \in I)$ *является базисом правого A-модуля* $V$, *если верны следующие утверждения.*

9.6.14.1: *Произвольный вектор* $v \in V$ *является линейной комбинацией векторов множества* $\overline{\overline{e}}$.

9.6.14.2: *Вектор* $e_i$ *нельзя представить в виде линейной комбинации остальных векторов множества* $\overline{\overline{e}}$.

Доказательство. Согласно утверждению 9.6.14.1, теореме 9.6.8 и определению 9.6.9, множество $\overline{\overline{e}}$ порождает правый A-модуль $V$ (определение 9.6.12). Согласно утверждению 9.6.14.2, множество $\overline{\overline{e}}$ является минимальным множеством, порождающим правый A-модуль $V$. Согласно определению 9.6.13, множество $\overline{\overline{e}}$ является базисом правого A-модуля $V$.             □

Теорема 9.6.15. *Пусть* $A$ - *ассоциативная D-алгебра с делением. Множество векторов* $\overline{\overline{e}} = (e_i, i \in I)$ *является* **базисом правого A-векторного пространства** $V$, *если векторы* $e_i$ *линейно независимы и любой вектор* $v \in V$ *линейно зависит от векторов* $e_i$.

---

[9.14] Я следую определению в [2], страница 100.

Доказательство. Пусть множество векторов $e_i,\ i \in I,$ линейно зависимо. Тогда в равенстве

$$e_i w^i = 0$$

существует индекс $i = j$ такой, что $w^j \neq 0$. Согласно теореме 9.6.10, вектор $e_j$ линейно зависит от остальных векторов множества $\overline{\overline{e}}$. Согласно определению 9.6.13, множество векторов $e_i,\ i \in I,$ не является базисом правого $A$-векторного пространства $V$.

Следовательно, если множество векторов $e_i,\ i \in I,$ является базисом, то эти векторы линейно независимы. Так как произвольный вектор $v \in V$ является линейной комбинацией векторов $e_i,\ i \in I,$, то множество векторов $v,\ e_i,\ i \in I,$ не является линейно независимым.    $\square$

Определение 9.6.16. *Пусть $\overline{\overline{e}}$ - базис правого $A$-модуля $V$, и вектор $\overline{v} \in V$ имеет разложение*

$$\overline{v} = e_* {}^* v = e_i v^i$$

*относительно базиса $\overline{\overline{e}}$. $A_{(1)}$-числа $v^i$ называются* **координатами** *вектора $\overline{v}$ относительно базиса $\overline{\overline{e}}$. Матрица $A_{(1)}$-чисел $v = (v^i, i \in I)$ называется* **координатной матрицей вектора** $\overline{v}$ *в базисе $\overline{\overline{e}}$.*    $\square$

Теорема 9.6.17. *Пусть $A$ - ассоциативная $D$-алгебра. Пусть $\overline{\overline{e}}$ - базис правого $A$-модуля $V$. Пусть*

$$(9.6.35) \qquad\qquad e_i w^i = 0$$

*линейная зависимость векторов базиса $\overline{\overline{e}}$. Тогда*

9.6.17.1: *$A_{(1)}$-число $w^i,\ i \in I,$ не имеет обратного элемента в $D$-алгебре $A_{(1)}$.*

9.6.17.2: *Множество $A'$ матриц $w = (w^i, i \in I)$ порождает правый $A$-модуль.*

Доказательство. Допустим существует матрица $w = (w^i, i \in I)$ такая, что равенство (9.6.35) верно и существует индекс $i = j$ такой, что $w^j \neq 0$. Если мы положим, что $A_{(1)}$-число $c^j$ имеет обратный, то равенство

$$e_j = \sum_{i \in I \setminus \{j\}} e_i w^i (w^j)^{-1}$$

является следствием равенства (9.6.35). Следовательно вектор $e_j$ является линейной комбинацией остальных векторов множества $\overline{\overline{e}}$ и множество $\overline{\overline{e}}$ не является базисом. Следовательно, наше предположение неверно, и $A_{(1)}$-число $c^j$ не имеет обратного.

Пусть матрицы $b = (b^i, i \in I) \in A',\quad c = (c^i, i \in I) \in A'$. Из равенств

$$e_i b^i = 0$$
$$e_i c^i = 0$$

следует

$$e_i (b^i + c^i) = 0$$

Следовательно, множество $A'$ является абелевой группой.

Пусть матрица $c = (c^i, i \in I) \in A'$ и $a \in A$. Из равенства

$$e_i c^i = 0$$

следует

$$e_i(c^i a) = 0$$

Следовательно, абелевая группа $A'$ является правым $A$-модулем. $\qquad\square$

Теорема 9.6.18. *Пусть правый $A$-модуль $V$ имеет базис $\overline{\overline{e}}$ такой, что в равенстве*

$$(9.6.36) \qquad\qquad e_i w^i = 0$$

*существует индекс $i = j$ такой, что $w^j \ne 0$. Тогда*

9.6.18.1: *Матрица $w = (w^i, i \in I)$ определяет координаты вектора $0 \in V$ относительно базиса $\overline{\overline{e}}$.*

9.6.18.2: *Координаты вектора $\overline{v}$ относительно базиса $\overline{\overline{e}}$ определены однозначно с точностью до выбора координат вектора $0 \in V$.*

Доказательство. Утверждение 9.6.18.1 является следствием равенства (9.6.36) и определения 9.6.16.

Пусть вектор $\overline{v}$ имеет разложение

$$(9.6.37) \qquad\qquad \overline{v} = e_*{}^* v = e_i v^i$$

относительно базиса $\overline{\overline{e}}$. Равенство

$$(9.6.38) \qquad\qquad \overline{v} = \overline{v} + 0 = e_i v^i + e_i c^i = e_i(v^i + c^i)$$

является следствием равенств (9.6.36), (9.6.37). Утверждение 9.6.18.2 является следствием равенств (9.6.37), (9.6.38) и определения 9.6.16. $\qquad\square$

Определение 9.6.19. *Правый $A$-модуль $V$ - **свободный правый $A$-модуль**,[9.15] если правый $A$-модуль $V$ имеет базис и векторы базиса линейно независимы.* $\qquad\square$

Теорема 9.6.20. *Координаты вектора $v \in V$ относительно базиса $\overline{\overline{e}}$ свободного правого $A$-модуля $V$ определены однозначно.*

Доказательство. Теорема является следствием теоремы 9.6.18 и определений 9.6.11, 9.6.19. $\qquad\square$

## 9.7. Левый модуль над неассоциативной алгеброй

Теоремы 9.6.5, 9.6.6 рассматривают структуру модуля над ассоциативной $D$-алгеброй $A$. Нетрудно заметить, что с учётом некоторых поправок, эти теоремы остаются верны если $A$ - неассоциативная $D$-алгебра. Однако, так как произведение в $D$-алгебре $A$ неассоциативно, а произведение преобразований в модуле над $D$-алгеброй $A$ ассоциативно, то отображение $g_{34}$ не может быть представлением неассоциативной $D$-алгебры $A$.

Мы подошли к той границе, где определена теория представлений универсальной алгебры. Для того, чтобы сохранить возможность применения рассмотренного в этой книге аппарата, мы можем согласиться, что отображение

$$g_{34} : A \times V \to V$$

является представлением, если отображение $g_{34}$ является билинейным отображением. Появляются новые вопросы, рассмотрение которых выходит за рамки этой книги.

---

[9.15] Я следую определению в [2], страница 103.

Однако мы можем рассмотреть эту задачу с другой стороны. Если отображение $g_{34}$ не сохраняет операцию произведения, то мы полагаем, что отображение $g_{34}$ - это представление $D$-алгебры $A$, в которой не определено произведение. Другими словами, отображение $g_{34}$ - это представление $D$-модуля. Следовательно, диаграмма представлений будет иметь вид

$$(9.7.1)$$

$$
\begin{array}{c}
A \xrightarrow{\ g_{34}\ } V \\
{}_{g_{12}}\diagdown \quad \diagup_{g_{14}} \\
D
\end{array}
\qquad
\begin{array}{l}
g_{12}(d): a \to d\,a \\
g_{34}(a): v \to a\,v \\
g_{14}(d): v \to d\,v
\end{array}
$$

Однако мы потеряли структуру $D$-алгебры $A$ в диаграмме представлений (9.7.1). Следовательно, правильная диаграмма представлений будет иметь вид

$$
\begin{array}{c}
A \xrightarrow{\ g_{34}\ } V \\
{}_{g_{23}}\nearrow \quad \\
A \quad \Big\downarrow {}_{g_{14}} \\
{}_{g_{12}}\Big\uparrow {}_{g_{12}} \\
D
\end{array}
\qquad
\begin{array}{l}
g_{12}(d): a \to d\,a \\
g_{23}(v): w \to C(w,v) \\
\qquad C \in \mathcal{L}(A^2 \to A) \\
g_{34}(a): v \to v\,a \\
g_{14}(d): v \to d\,v
\end{array}
$$

# Примеры диаграммы представлений: аффинная геометрия

## 10.1. Об этой главе

В главе 9 мы рассмотрели примеры диаграммы представлений, связанные с модулем над кольцом. Если бы теория представлений сводилась бы к изучению модулей, вряд ли это была бы интересная теория.

В этой главе я рассмотрел примеры диаграммы представлений, связанные с аффинной геометрией. Внешне простая алгебраическая конструкция оказалась для меня богатейшим источником вдохновения. Я дважды встретил интересные идеи в этой области математики. Сперва изучая аффинную геометрию, я обнаружил, что я могу описать аффинную геометрию с помощью башни представлений. Впоследствии, изучая таким же образом алгебру над коммутативным кольцом, я стал изучать диаграмму представлений.

Однако второе открытие пришло ко мне случайно. Когда я просматривал учебник по математическому анализу, я встретил определение знакомое с детства. Сумма векторов. Определение крайне простое. При определении многообразия аффинной связности мы имем сумму векторов в касательной плоскости. Но в этот раз я понял, что я могу определить сумму векторов пользуясь параллелограммом из геодезических. Это позволило построить аффинную геометрию на многообразии аффинной связности.

Ещё один шаг, и от многообразия аффинной связности я перешёл к метрико аффинному многообразию. Так как параллелограм из геодезических не замкнут, то сумма векторов в метрико аффинном многообразии не коммутативна. Без сомнения, это исследование, которое выходит за рамки этой книги и к которому я надеюсь вернуться в будущем. Однако я решил написать набросок этой конструкции в разделе 10.4, чтобы показать читателю границы теории, изложенной в этой книге.

Теория представлений является естественным продолжением теории универсальных алгебр. Предполагается, что бинарная операция на универсальной алгебре $A$ определена для любой пары $A$-чисел. Однако очевидно, что сумма векторов в аффинной геометрии на дифференцируемом многообразии хорошо определена только в достаточно малой окрестности.

С похожей задачей я столкнулся в статье [11], где я и Александр Ложье изучали ортогональные преобразования в пространстве Минковского. Мы обнаружили, что произведение ортогональных преобразований не всегда является ортогональным преобразованием, и следовательно множество ортогональных преобразований не является группой.

## 10.2. Представление группы на множестве

Пусть $G$ - абелева группа, и $M$ - множество. Рассмотрим эфективное представление группы $G$ на множестве $M$. Для заданных $a \in G$, $A \in M$ положим $A \to A + a$. Мы будем также пользоваться записью $a = \overrightarrow{AB}$, если

$$(10.2.1) \qquad B = A + a$$

Тогда действие группы можно представить в виде

$$(10.2.2) \qquad B = A + \overrightarrow{AB}$$

Поскольку представление эффективно, то из равенств ($10.2.1$), ($10.2.2$) и равенства

$$D = C + a$$

следует, что

$$(10.2.3) \qquad \overrightarrow{AB} = \overrightarrow{CD}$$

Мы будем называть вектором $G$-число $a$ и соответствующее преобразование $\overrightarrow{AB}$. Мы интерпретируем равенство ($10.2.3$) как независимость вектора $a$ от выбора $M$-числа $A$.

Мы можем рассматривать множество $M$ как объединение орбит представления группы $G$. В качестве базиса представления можно выбрать множество точек таким образом, что одна и только одна точка принадлежит каждой орбите представления. Если $X$ - базис представления, $A \in X$, $g \in G$, то $\Omega_2$-слово имеет вид $A + g$. Поскольку на множестве $M$ не определены операции, то не существует $\Omega_2$-слово, содержащее различные элементы базиса. Если представление группы $G$ однотранзитивно, то базис представления состоит из одной точки. Этой точкой может быть любая точка множества $M$.

ТЕОРЕМА 10.2.1. *Пусть представление $A \to A + a$ абелевой группы $G$ на множестве $M$ однотранзитивно. Тогда, для любых $M$-чисел $A$, $B$, $C$, определена сумма векторов $\overrightarrow{AB}$ и $\overrightarrow{BC}$ и сумма векторов удовлетворяет равенству*

$$(10.2.4) \qquad \overrightarrow{AB} + \overrightarrow{BC} = \overrightarrow{AC}$$

ДОКАЗАТЕЛЬСТВО. Поскольку представление однотранзитивно, то, для любых $M$-чисел $A$, $B$, $C$, существуют векторы $\overrightarrow{AB}$, $\overrightarrow{BC}$ такие, что

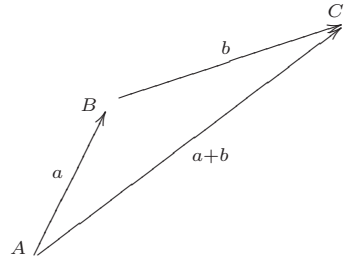

$$(10.2.5) \qquad B = A + \overrightarrow{AB}$$

$$(10.2.6) \qquad C = B + \overrightarrow{BC}$$

Равенство

$$(10.2.7) \qquad C = (A + \overrightarrow{AB}) + \overrightarrow{BC} = A + (\overrightarrow{AB} + \overrightarrow{BC})$$

является следствием равенств ($10.2.5$), ($10.2.6$) и ассоциативности сложения в абелевой группе $G$. Поскольку представление однотранзитивно, то равенство

(10.2.4) является следствием равенства (10.2.7) и равенства

$$C = A + \overrightarrow{AC}$$

Это определение суммы называется правилом треугольника. □

ЗАМЕЧАНИЕ 10.2.2. *Так как $G$ - абелева группа, то утверждения 10.2.2.1, 10.2.2.2, следуют из теоремы 10.2.1*

10.2.2.1: $\overrightarrow{AA} = 0$

10.2.2.2: $\overrightarrow{AB} = -\overrightarrow{BA}$

10.2.2.3: *Сложение коммутативно.*

10.2.2.4: *Сложение ассоциативно.*

□

ТЕОРЕМА 10.2.3. *Для заданных $a$, $b \in G$ и $A \in M$ рассмотрим следующее множество $M$-чисел.*

- $B = A + a$
- $C = B + b$
- $D = A + b$
- $E = D + a$

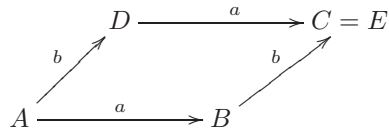

Доказательство. Теорема является следствием утверждения 10.2.2.3. □

ТЕОРЕМА 10.2.4. *Если $\overrightarrow{AB} = \overrightarrow{CD}$, то $\overrightarrow{AC} = \overrightarrow{BD}$.*

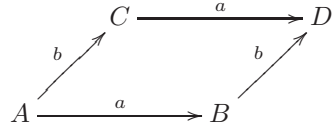

Доказательство. Положим $\overrightarrow{AB} = \overrightarrow{CD} = a$, $\overrightarrow{AC} = b$. Согласно утверждению 10.2.2.2, $\overrightarrow{BA} = -a$. Теорема является следствием равенства

$$D = B + \overrightarrow{BD} = B + \overrightarrow{BA} + \overrightarrow{AD} = B + \overrightarrow{BA} + \overrightarrow{AC} + \overrightarrow{CD}$$
$$= B - a + b + a = B + b$$

□

## 10.3. Аффинное пространство

ОПРЕДЕЛЕНИЕ 10.3.1. *Пусть $D$ - коммутативное кольцо и $V$ - свободный $D$-модуль. Множество точек $\overset{\circ}{V}$ называется **аффинным пространством** над $D$-модулем $V$, если множество точек $\overset{\circ}{V}$ удовлетворяет следующим аксиомам.[10.1]*

10.3.1.1: *Существует по крайней мере одна точка*

10.3.1.2: *Каждой паре точек $(A, B)$ поставлен в соответствие один и только один вектор. Этот вектор мы будем обозначать $\overrightarrow{AB}$. Вектор $\overrightarrow{AB}$ имеет начало в точке $A$ и конец в точке $B$.*

---

[10.1] Я написал определения и теоремы в этом разделе согласно определению аффинного пространства в [4], с. 86 - 93.

10.3.1.3: *Для каждой точки $A$ и любого вектора $a$ существует одна и только одна точка $B$ такая, что $\overrightarrow{AB} = a$. Мы будем также пользоваться записью* [10.2]

$$(10.3.1) \qquad\qquad B = A + a$$

10.3.1.4: *(Аксиома параллелограмма.) Если $\overrightarrow{AB} = \overrightarrow{CD}$, то $\overrightarrow{AC} = \overrightarrow{BD}$.*

*Множество $V$ называется множеством свободных векторов. $\overset{\circ}{V}$-числа называются точками аффинного пространства $\overset{\circ}{V}$.* □

Определение 10.3.2. *Пусть $A \in \overset{\circ}{V}$ - произвольная точка.*

*Пусть $v$ - вектор. Согласно аксиоме 10.3.1.3, существует $B \in \overset{\circ}{V}$, $B = A + v$.*

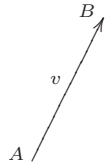

*Пусть $w$ - вектор. Согласно аксиоме 10.3.1.3, существует $C \in \overset{\circ}{V}$, $C = B + w$.*

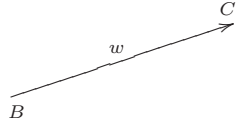

*Согласно аксиоме 10.3.1.2, существует вектор $\overrightarrow{AC}$. Вектор $\overrightarrow{AC}$ называется суммой векторов $v$ и $w$*

$$(10.3.2) \qquad v + w = \overrightarrow{AC}$$

*Это определение суммы называется правилом треугольника.*

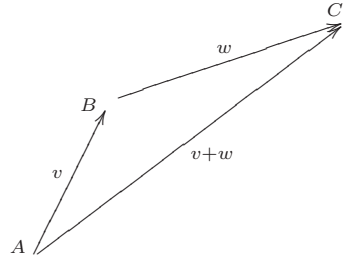

□

Теорема 10.3.3. *Вектор $\overrightarrow{AA}$ является нулём по отношению к операции сложения и не зависит от выбора точки $A$. Вектор $\overrightarrow{AA}$ называется нуль-вектором и мы полагаем $\overrightarrow{AA} = 0$.*

Доказательство. Мы можем записать правило сложения (10.3.2) в виде

$$(10.3.3) \qquad\qquad \overrightarrow{AB} + \overrightarrow{BC} = \overrightarrow{AC}$$

Если $B = C$, то из равенства (10.3.3) следует

$$(10.3.4) \qquad\qquad \overrightarrow{AB} + \overrightarrow{BB} = \overrightarrow{AB}$$

Из равенства (10.3.4) следует, что вектор $\overrightarrow{BB}$ является нулём по отношению к операции сложения. Если $C = A$, $B = D$, то из аксиомы 10.3.1.4 следует $\overrightarrow{AA} = \overrightarrow{BB}$. Следовательно, нуль-вектор $\overrightarrow{AA}$ не зависит от выбора точки $A$. □

Теорема 10.3.4. *Пусть* $a = \overrightarrow{AB}$. *Тогда*

(10.3.5)
$$\overrightarrow{BA} = -a$$

*и это равенство не зависит от выбора точки $A$.*

Доказательство. Из равенства (10.3.3) и теоремы 10.3.3 следует

(10.3.6)
$$\overrightarrow{AB} + \overrightarrow{BA} = \overrightarrow{AA} = 0$$

Равенство (10.3.5) следует из равенства (10.3.6). Применяя аксиому 10.3.1.4 к равенству $\overrightarrow{AB} = \overrightarrow{CD}$ получим $\overrightarrow{AC} = \overrightarrow{BD}$, или, что то же,

(10.3.7)
$$\overrightarrow{BD} = \overrightarrow{AC}$$

Из равенства (10.3.7) и аксиомы 10.3.1.4 следует, что $\overrightarrow{BA} = \overrightarrow{DC}$. Следовательно, равенство (10.3.5) не зависит от выбора точки $A$. $\square$

Теорема 10.3.5. *Сумма векторов $v$ и $w$ не зависит от выбора точки $A$.*

Доказательство. Пусть

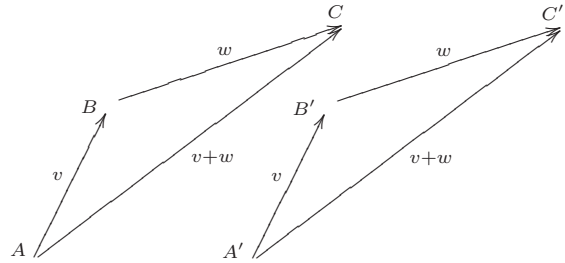

(10.3.8)     $v = \overrightarrow{AB} = \overrightarrow{A'B'}$

(10.3.9)     $w = \overrightarrow{BC} = \overrightarrow{B'C'}$

Сумма векторов $v$ и $w$ определена согласно определению 10.3.2.

$$\overrightarrow{AB} + \overrightarrow{BC} = \overrightarrow{AC}$$
$$\overrightarrow{A'B'} + \overrightarrow{B'C'} = \overrightarrow{A'C'}$$

Согласно аксиоме 10.3.1.4, из равенств (10.3.8), (10.3.9) следует, что

(10.3.10)
$$\overrightarrow{A'A} = \overrightarrow{B'B} = \overrightarrow{C'C}$$

Применяя аксиому 10.3.1.4 к крайним членам равенства (10.3.10), получаем

(10.3.11)
$$\overrightarrow{A'C'} = \overrightarrow{AC}$$

Из равенства (10.3.11) следует утверждение теоремы. $\square$

Теорема 10.3.6. *Сложение векторов ассоциативно.*

Доказательство. Пусть $v = \overrightarrow{AB}$, $w = \overrightarrow{BC}$, $u = \overrightarrow{CD}$. Из равенства

$$v \ + \ w \ = \overrightarrow{AC}$$
$$\overrightarrow{AB} \ + \ \overrightarrow{BC} \ = \overrightarrow{AC}$$

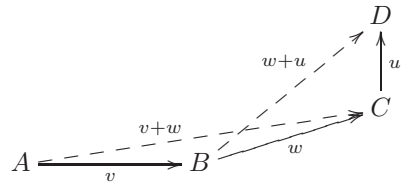

следует

$$(10.3.12) \quad \begin{aligned} (v+w) \quad + \quad u \quad &= \overrightarrow{AD} \\ \overrightarrow{AC} \quad + \quad \overrightarrow{CD} \quad &= \overrightarrow{AD} \end{aligned}$$

Из равенства

$$\begin{aligned} w \quad + \quad u \quad &= \overrightarrow{BD} \\ \overrightarrow{BC} \quad + \quad \overrightarrow{CD} \quad &= \overrightarrow{BD} \end{aligned}$$

следует

$$(10.3.13) \quad \begin{aligned} v \quad + \quad (w+u) \quad &= \overrightarrow{AD} \\ \overrightarrow{AB} \quad + \quad \overrightarrow{BD} \quad &= \overrightarrow{AD} \end{aligned}$$

Теорема следует из сравнения равенств (10.3.12) и (10.3.13). □

**Теорема 10.3.7.** *На множестве $V$ определена структура абелевой группы.*

Доказательство. Из теорем 10.3.3, 10.3.4, 10.3.5, 10.3.6 следует, что операция сложение векторов определяет группу.

Пусть $v = \overrightarrow{AB}$, $w = \overrightarrow{BC}$.

$$(10.3.14) \quad \begin{aligned} v \quad + \quad w \quad &= \overrightarrow{AC} \\ \overrightarrow{AB} \quad + \quad \overrightarrow{BC} \quad &= \overrightarrow{AC} \end{aligned}$$

Согласно аксиоме 10.3.1.3 существует точка $D$ такая, что

$$w = \overrightarrow{AD} = \overrightarrow{BC}$$

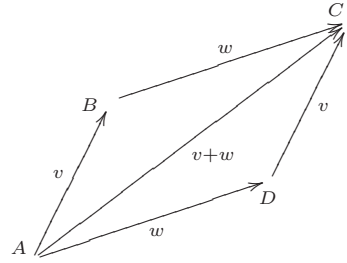

Правило параллелограмма.

Согласно аксиоме 10.3.1.4, $\overrightarrow{AB} = \overrightarrow{DC} = v$. Согласно определению суммы векторов

$$(10.3.15) \quad \begin{aligned} \overrightarrow{AD} \quad + \quad \overrightarrow{DC} \quad &= \overrightarrow{AC} \\ w \quad + \quad v \quad &= \overrightarrow{AC} \end{aligned}$$

Из равенств (10.3.14) и (10.3.15) следует коммутативность сложения. □

**Теорема 10.3.8.** *Отображение*

$$(10.3.16) \quad V \to \operatorname{End}(\emptyset, \overset{\circ}{V})$$

*определённое равенством* (10.3.1), *является однотранзитивным представлением абелевой группы $V$.*

Доказательство. Аксиома 10.3.1.3 определяет отображение (10.3.16). Из теоремы 10.3.5 следует, что отображение (10.3.16) является представлением. Эффективность представления следует из теоремы 10.3.3 и аксиомы 10.3.1.2. Из аксиомы 10.3.1.2 следует также, что представление транзитивно. Эффективное и транзитивное представление однотранзитивно. □

Сравнение теоремы 10.3.8 и утверждений раздела 10.2 видно, что одно-транзитивное представление абелевой группы $V$ на множестве $\overset{\circ}{V}$ эквивалентно аксиомам аффинного пространства. Однако, пользуясь теоремой 10.3.8 как определением аффинного пространства, мы теряем многие важные конструкции в аффинном пространстве. Например, вектор определяет преобразование параллельного переноса в аффинном пространстве. Но у нас нет инструмента, чтобы определить преобразование поворота аффинного пространства.

Если мы внимательно посмотрим на определение 10.3.1, то мы увидим, что абелевая группа $V$ имеет дополнительную структуру, поскольку абелевая группа $V$ является $D$-модулем. Таким образом, мы получаем следующую теорему.

Теорема 10.3.9. *Пусть $D$ - коммутативное кольцо, $V$ - абелева группа и $\overset{\circ}{V}$ - множество. Если $A \in \overset{\circ}{V}$ и $v \in V$, то мы будем обозначать действие вектора $v$ на точку $A$ выражением $A + v$.* **Аффинное пространство** *над $D$-модулем $V$ - это диаграмма представлений*

$$\overrightarrow{V} : D \xrightarrow{f_{12}} V \xrightarrow{f_{23}} \overset{\circ}{V} \qquad \begin{aligned} f_{12}(d) &: v \to d\,v \\ f_{23}(v) &: A \to A + v \end{aligned}$$

*где $f_{12}$ - эффективное представление коммутативного кольца $D$ в абелевой группе $V$ и $f_{23}$ - однотранзитивное правостороннее представление абелевой группы $V$ в множестве $\overset{\circ}{V}$.*

Доказательство. Мы полагаем, что множество $\overset{\circ}{V}$ не пусто; следовательно множество $\overset{\circ}{V}$ удовлетворяет аксиоме 10.3.1.1. Поскольку $v \in V$ порождает преобразование множества, то для любого $A \in M$ однозначно определён $B \in M$ такое, что

$$B = A + v$$

Это утверждение доказывает аксиому 10.3.1.3. Поскольку представление $f_{23}$ однотранзитивно, то для любых $A$, $B \in \overset{\circ}{V}$ существует единственное $v \in V$ такое, что

$$B = A + v$$

Это утверждение позволяет ввести обозначение $\overrightarrow{AB} = a$, а также доказывает аксиому 10.3.1.2. Аксиома 10.3.1.4 следует из утверждения теоремы 10.2.4.

Представление $f_{12}$ гарантирует, что абелевая группа $V$ является $D$-модулем. $\square$

Абелевая группа $V$ действует однотранзитивно на множестве $\overset{\circ}{V}$. Из построений в разделе 10.2 следует, что базис множества $\overset{\circ}{V}$ относительно представления абелевой группы $V$ состоит из одной точки. Эту точку обычно обозначают буквой $O$ и называют **началом системы координат аффинного пространства**. Следовательно, произвольную точку $A \in \overset{\circ}{V}$ можно представить с помощью вектора $\overrightarrow{OA} \in V$

Пусть $\overline{\overline{e}}$ - базис $D$-модуля $V$. Тогда вектор $\overrightarrow{OA}$ имеет вид

$$\overrightarrow{OA} = a^i e_i$$

Множество $(a_i, i \in I)$ называется **координатами точки** $A$ **аффинного пространства** $\overset{\circ}{A}$ **относительно базиса** $(O, \overline{\overline{e}})$.

## 10.4. Аффинное пространство на дифференцируемом многообразии

В разделе 10.3 мы рассмотрели определение аффинной геометрии. Ниже мы рассмотрим модель аффинного пространства в метрико-аффинном многообразии. Когда мы рассматриваем связность $\Gamma_{ij}^k$ в римановом пространстве, мы накладываем на связность ограничение,[10.3] что тензор кручения

$$(10.4.1) \qquad T_{kl}^i = \Gamma_{lk}^i - \Gamma_{kl}^i$$

обращается в 0 (симметрия связности) и скалярное произведение векторов при параллельном переносе не меняется. Если на дифференцируемом многообразии определены метрический тензор и произвольная связность, то это многообразие называется **метрико-аффинным многообразием**.[10.4] В частности, связность в метрико-аффинном многообразии имеет кручение.

В римановом пространстве, мы пользуемся геодезическими вместо прямых. Поэтому вектор $v$ мы можем представить с помощью отрезка $AB$ геодезической $L_v$ при условии, что вектор $v$ касается геодезической $L_v$ в точке $A$ и длина отрезка $AB$ равна длине вектора $v$.

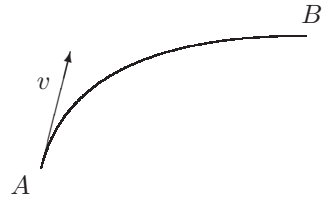

Это определение позволяет отождествить вектор $v$ и отрезок $AB$ геодезической $L_v$.

Для заданных векторов $v$ и $w$ в касательной плоскости к точке $A$, мы будем полагать $\rho > 0$ - длина вектора $v$ и $\sigma > 0$ - длина вектора $w$. Пусть $V$ - единичный вектор, колинеарный вектору $v$

$$(10.4.2) \qquad V^k \rho = v^k$$

Пусть $W$ - единичный вектор, колинеарный вектору $w$

$$(10.4.3) \qquad W^k \sigma = w^k$$

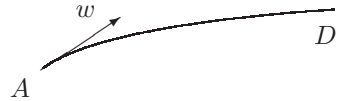

---

[10.3] Смотри определение аффинной связности в римановом пространстве на странице [4]-443.

[10.4] Смотри также определение [9]-5.4.1.

Мы проведём геодезическую $L_v$ через точку $A$, используя вектор $v$ как касательный вектор к $L_v$ в точке $A$. Пусть $\tau$ - канонический параметр на $L_v$ и

$$\frac{dx^k}{d\tau} = V^k$$

Мы перенесём вектор $w$ вдоль геодезической $L_v$ из точки $A$ в точку $B$, определённую значением параметра $\tau = \rho$. Мы обозначим результат $w'$.

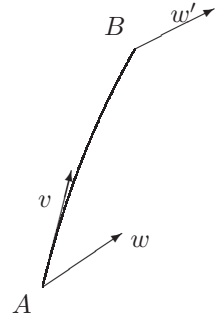

Мы проведём геодезическую $L_{w'}$ через точку $B$, используя вектор $W'$ как касательный вектор к $L_{w'}$ в точке $B$. Пусть $\varphi'$ - канонический параметр на $L_{w'}$ и

$$\frac{dx^k}{d\varphi'} = W'^k$$

Мы определим точку $C$ на геодезической $L_{w'}$ значением параметра $\varphi' = \sigma$

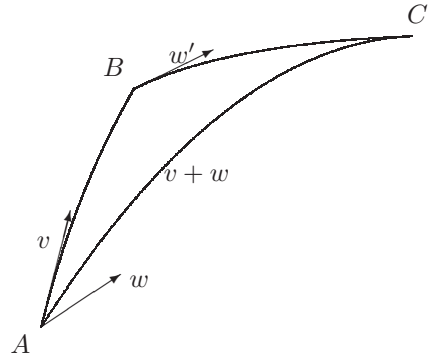

Я полагаю, что длина векторов $v$ и $w$ мала. Тогда существует единственная геодезическая $L_u$ из точки $A$ в точку $C$. Я буду отождествлять отрезок $AC$ геодезической $L_u$ и вектор $v + w$.

Аналогичным образом я строю треугольник $ADE$, чтобы определить вектор $w + v$.

Мы проведём геодезическую $L_w$ через точку $A$, используя вектор $w$ как касательный вектор к $L_w$ в точке $A$. Пусть $\varphi$ - канонический параметр на $L_w$ и

$$\frac{dx^k}{d\varphi} = W^k$$

Мы перенесём вектор $v$ вдоль геодезической $L_w$ из точки $A$ в точку $D$, определённую значением параметра $\varphi = \sigma$. Мы обозначим результат $v'$.

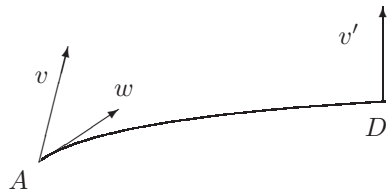

Мы проведём геодезическую $L_{v'}$ через точку $D$, используя вектор $v'$ как касательный вектор к $L_{v'}$ в точке $D$. Пусть $\tau'$ - канонический параметр на $L_{v'}$ и

$$\frac{dx^k}{d\tau'} = V'^k$$

Мы определим точку $E$ на геодезической $L_{v'}$ значением параметра $\tau' = \rho$

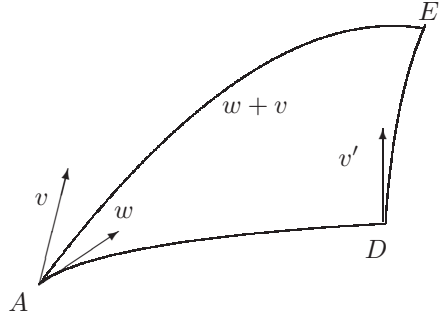

Существует единственная геодезическая $L_u$ из точки $A$ в точку $E$. Я буду отождествлять отрезок $AE$ геодезической $L_u$ и вектор $w + v$.

Формально линии $AB$ и $DE$ так же, как линии $AD$ и $BC$, параллельны. Длины отрезков $AB$ и $DE$ равны так же, как длины отрезков $AD$ и $BC$ равны. Мы называем такую фигуру **параллелограммом**, построенным на векторах $v$ и $w$ с вершиной в точке $A$.

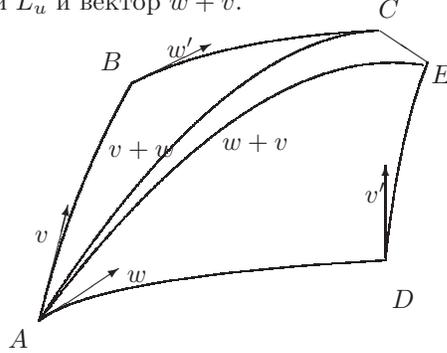

**Лемма 10.4.1.** *Пусть $L_v$ - геодезическая через точку $A$ и вектор $v$ - касательный вектор к $L_v$ в точке $A$. Приращение координаты $x^k$ вдоль геодезической $L_v$ имеет вид*

$$(10.4.4) \qquad \Delta x^k = \frac{dx^k}{d\tau}\tau - \frac{1}{2}\Gamma^k_{mn}\frac{dx^m}{d\tau}\frac{dx^n}{d\tau}\tau^2 + O(\tau^2)$$

*где $\tau$ - канонический параметр и мы вычисляем производные и компоненты $\Gamma^k_{mn}$ в начальной точке.*

Доказательство. Система дифференциальных уравнений геодезической $L_v$ имеет вид

$$(10.4.5) \qquad \frac{d^2x^i}{d\tau^2} = -\Gamma^i_{kl}\frac{dx^k}{d\tau}\frac{dx^l}{d\tau}$$

Мы можем записать решение системы дифференциальных уравнений (10.4.5) в виде ряда Тейлора

$$(10.4.6) \qquad \begin{aligned} \Delta x^k &= \frac{dx^k}{d\tau}\tau + \frac{1}{2}\frac{d^2x^k}{d\tau^2}\tau^2 + O(\tau^2) = \\ &= \frac{dx^k}{d\tau}\tau - \frac{1}{2}\Gamma^k_{mn}\frac{dx^m}{d\tau}\frac{dx^n}{d\tau}\tau^2 + O(\tau^2) \end{aligned}$$

Равенство (10.4.4) является следствием равенства (10.4.6).  $\square$

ТЕОРЕМА 10.4.2. *Предположим CBADE - параллелограм с вершиной в точке A; тогда построенная фигура не будет замкнута [1]. Величина различия координат точек C и E равна поверхностному интегралу кручения над этим параллелограммом*

$$\Delta_{CE} x^k = \iint T^k_{mn} dx^m \wedge dx^n$$

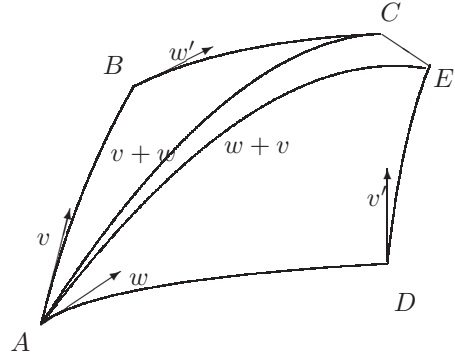

Доказательство. Согласно лемме 10.4.1, приращение координаты $x^k$ вдоль геодезической $L_v$ имеет следующий вид[10.5]

$$\Delta_{CE} x^k = \iint T^k_{mn} dx^m \wedge dx^n$$

$$\Delta_{AB} x^k = V^k \rho - \frac{1}{2} \Gamma^k_{mn}(A) V^m V^n \rho^2 + O(\rho^2)$$

и приращение координаты $x^k$ вдоль геодезической $L_{b'}$ имеет следующий вид

(10.4.7) $$\Delta_{BC} x^k = W'^k \sigma - \frac{1}{2} \Gamma^k_{mn}(B) W'^m W'^n \sigma^2 + O(\sigma^2)$$

Здесь

(10.4.8) $$\begin{aligned} W'^k &= W^k - \Gamma^k_{mn}(A) W^m \Delta_{AB} x^n + O(dx) \\ &= W^k - \Gamma^k_{mn}(A) W^m V^n \rho + O(\rho) \end{aligned}$$

результат параллельного переноса вектора $w$ из $A$ в $B$ и

(10.4.9) $$\begin{aligned} \Gamma^k_{mn}(B) &= \Gamma^k_{mn}(A) + \partial_p \Gamma^k_{mn}(B) \Delta_{AB} x^p \\ &= \Gamma^k_{mn}(A) + \partial_p \Gamma^k_{mn}(B) V^p \rho \end{aligned}$$

с точностью до малой первого порядка. Подставляя (10.4.8), (10.4.9) в (10.4.7), мы получим

$$\Delta_{BC} x^k = W^k \sigma - \Gamma^k_{mn}(A) W^m V^n \sigma \rho - \frac{1}{2} \Gamma^k_{mn}(A) W^m W^n \sigma^2 + O(\rho^2)$$

Общее приращение координаты $x^K$ вдоль пути $ABC$ имеет вид

(10.4.10) $$\begin{aligned} \Delta_{ABC} x^k &= \Delta_{AB} x^k + \Delta_{BC} x^k \\ &= V^k \rho + W^k \sigma - \Gamma^k_{mn}(A) W^m V^n \sigma \rho - \\ &- \frac{1}{2} \Gamma^k_{mn}(A) W^m W^n \sigma^2 - \frac{1}{2} \Gamma^k_{mn}(A) V^m V^n \rho^2 + O(dx^2) \end{aligned}$$

Аналогично общее приращение координаты $x^K$ вдоль пути $ADE$ имеет вид

(10.4.11) $$\begin{aligned} \Delta_{ADE} x^k &= \Delta_{AD} x^k + \Delta_{DE} x^k = \\ &= W^k \sigma + V^k \rho - \Gamma^k_{mn}(A) V^m W^n \rho \sigma - \\ &- \frac{1}{2} \Gamma^k_{mn}(A) V^m V^n \rho^2 - \frac{1}{2} \Gamma^k_{mn}(A) W^m W^n \sigma^2 + O(dx^2) \end{aligned}$$

---

10.5 Доказательство этого утверждения я нашёл в [3]

Из ([10.4.10](#)) и ([10.4.11](#)) следует, что

$$\Delta_{ADE}x^k - \Delta_{ABC}x^k = -\Gamma^k_{mn}(A)V^mW^n\rho\sigma$$
$$-\frac{1}{2}\Gamma^k_{mn}(A)V^mV^n\rho^2 - \frac{1}{2}\Gamma^k_{mn}(A)W^mW^n\sigma^2$$
$$+\Gamma^k_{mn}(A)W^mV^n\sigma\rho$$
$$+\frac{1}{2}\Gamma^k_{mn}(A)W^mW^n\sigma^2 + \frac{1}{2}\Gamma^k_{mn}(A)V^mV^n\rho^2$$

и мы получаем интегральную сумму для выражения

$$\Delta_{ADE}x^k - \Delta_{ABC}x^k = \iint_\Sigma (\Gamma^k_{nm} - \Gamma^k_{mn})dx^m \wedge dx^n$$

□

Теорема 10.4.3. *В римано­вом пространстве параллелограм $ABCD$ замкнут. В точке $A$ геоде­зическая $AC$ имеет касательный вектор $u$, который является сум­мой векторов $v$ и $w$*

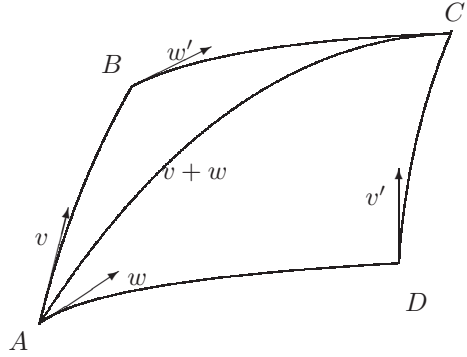

(10.4.12)        $u^k = v^k + w^k$

*Следовательно, сумма векторов в римановом пространстве комму­тативна.*

Доказательство. Мы будем полагать $\pi$ - длина вектора $u$. Пусть $U$ - единичный вектор, колинеарный вектору $u$

(10.4.13)                    $U^k\pi = u^k$

Согласно лемме [10.4.1](#), приращение координаты $x^k$ вдоль геодезической $L_u$ имеет следующий вид

(10.4.14)        $\Delta_{AC}x^k = U^k\pi - \frac{1}{2}\Gamma^k_{mn}(A)U^mU^n\pi^2 + O(\pi^2)$

Равенства

(10.4.15)                    $U^k\pi = V^k\rho + W^k\sigma$

(10.4.16)

$$\Gamma^k_{mn}(A)U^mU^n\pi^2$$
$$= 2\Gamma^k_{mn}(A)W^mV^n\sigma\rho + \Gamma^k_{mn}(A)W^mW^n\sigma^2 + \Gamma^k_{mn}(A)V^mV^n\rho^2$$

являются следствием равенств ([10.4.10](#)), ([10.4.14](#)). Равенство ([10.4.12](#)) является следствием равенств ([10.4.2](#)), ([10.4.3](#)), ([10.4.13](#)), ([10.4.15](#)). Равенство

(10.4.17)        $\Gamma^k_{mn}(A)u^mu^n = 2\Gamma^k_{mn}(A)w^mv^n + \Gamma^k_{mn}(A)w^mw^n + \Gamma^k_{nm}(A)v^mv^n$

является следствием равенств ([10.4.2](#)), ([10.4.3](#)), ([10.4.13](#)), ([10.4.16](#)). Равенство

(10.4.18)        $\Gamma^k_{mn}(A)u^mu^n = \Gamma^k_{mn}(A)(v^m + w^m)(v^n + w^n)$
$$= \Gamma^k_{mn}(A)(v^mv^n + v^mw^n + w^mv^n + w^mw^n)$$

является следствием равенства (10.4.12). Равенство (10.4.17) является следствием равенства (10.4.18) и симметрии связности. Следовательно, геодезическая $AC$ является суммой геодезических $AB$ и $BC$. □

Если связность не симметрична, то геодезическая $L_u$ не содержит точек $C$ и $E$. Следовательно, сумма векторов в метрико-аффинном многообразии некоммутативна.

Теорема 10.4.4. *Существует вектор $t$ такой, что*

$$(10.4.19) \qquad (v+w)^k = v^k + w^k + t^k$$

$$(10.4.20) \qquad (w+v)^k = v^k + w^k - t^k$$

*Координаты вектора $t$ удовлетворяют системе уравнений*

$$(10.4.21) \quad \Gamma^k_{mn}(A)t^m t^n + (\Gamma^k_{mn}(A) + \Gamma^k_{nm}(A))(v^m + w^m)t^n + 2T^k_{mn}(A)v^m w^n = 0$$

Доказательство. Рассмотрим сперва вектор $v + w$. Равенство

$$(10.4.22) \quad \begin{aligned} &v^k + w^k + t^k - \tfrac{1}{2}\Gamma^k_{mn}(A)(v^m + w^m + t^m)(v^n + w^n + t^n) \\ &= v^k + w^k - \Gamma^k_{mn}(A)w^m v^n - \tfrac{1}{2}\Gamma^k_{mn}(A)w^m w^n - \tfrac{1}{2}\Gamma^k_{mn}(A)v^m v^n \end{aligned}$$

является следствием равенства (10.4.10) и леммы 10.4.1. Равенство

$$(10.4.23) \quad \begin{aligned} &v^k + w^k + t^k \\ &- \tfrac{1}{2}\Gamma^k_{mn}(A)v^m v^n + \tfrac{1}{2}\Gamma^k_{mn}(A)v^m w^n + \tfrac{1}{2}\Gamma^k_{mn}(A)v^m t^n \\ &+ \tfrac{1}{2}\Gamma^k_{mn}(A)w^m v^n + \tfrac{1}{2}\Gamma^k_{mn}(A)w^m w^n + \tfrac{1}{2}\Gamma^k_{mn}(A)w^m t^n \\ &+ \tfrac{1}{2}\Gamma^k_{mn}(A)t^m v^n + \tfrac{1}{2}\Gamma^k_{mn}(A)t^m w^n + \tfrac{1}{2}\Gamma^k_{mn}(A)t^m t^n \\ &= v^k + w^k - \Gamma^k_{mn}(A)w^m v^n - \tfrac{1}{2}\Gamma^k_{mn}(A)w^m w^n - \tfrac{1}{2}\Gamma^k_{mn}(A)v^m v^n \end{aligned}$$

является следствием равенства (10.4.22) Равенство

$$(10.4.24) \quad \begin{aligned} &t^k - \tfrac{1}{2}\Gamma^k_{mn}(A)v^m w^n - \tfrac{1}{2}\Gamma^k_{mn}(A)v^m t^n - \tfrac{1}{2}\Gamma^k_{mn}(A)w^m t^n \\ &- \tfrac{1}{2}\Gamma^k_{mn}(A)t^m v^n - \tfrac{1}{2}\Gamma^k_{mn}(A)t^m w^n - \tfrac{1}{2}\Gamma^k_{mn}(A)t^m t^n \\ &= -\tfrac{1}{2}\Gamma^k_{mn}(A)w^m v^n \end{aligned}$$

является следствием равенства (10.4.23) Равенство

$$(10.4.25) \quad \begin{aligned} &\Gamma^k_{mn}(A)t^m t^n + (\Gamma^k_{mn}(A)v^m + \Gamma^k_{mn}(A)w^m \\ &+ \Gamma^k_{nm}(A)v^m + \Gamma^k_{nm}(A)w^m - 2\delta^k_n)t^n \\ &+ 2T^k_{mn}(A)v^m w^n = 0 \end{aligned}$$

является следствием равенства (10.4.24) Равенство (10.4.21) является следствием равенства (10.4.25)

Рассмотрим теперь вектор $w + v$. Равенство

$$(10.4.26) \quad \begin{aligned} &v^k + w^k - t^k - \frac{1}{2}\Gamma^k_{mn}(A)(v^m + w^m - t^m)(v^n + w^n - t^n) \\ &= w^k + v^k - \Gamma^k_{mn}(A)v^m w^n - -\frac{1}{2}\Gamma^k_{mn}(A)v^m v^n - \frac{1}{2}\Gamma^k_{mn}(A)W^m W^n \end{aligned}$$

является следствием равенства (10.4.11) и леммы 10.4.1. Равенство

$$(10.4.27) \quad \begin{aligned} &v^k + w^k - t^k \\ &-\frac{1}{2}\Gamma^k_{mn}(A)v^m v^n - \frac{1}{2}\Gamma^k_{mn}(A)v^m w^n + \frac{1}{2}\Gamma^k_{mn}(A)v^m t^n \\ &-\frac{1}{2}\Gamma^k_{mn}(A)w^m v^n - \frac{1}{2}\Gamma^k_{mn}(A)w^m w^n + \frac{1}{2}\Gamma^k_{mn}(A)w^m t^n \\ &+\frac{1}{2}\Gamma^k_{mn}(A)t^m v^n + \frac{1}{2}\Gamma^k_{mn}(A)t^m w^n - \frac{1}{2}\Gamma^k_{mn}(A)t^m t^n \\ &= w^k + v^k - \Gamma^k_{mn}(A)v^m w^n - -\frac{1}{2}\Gamma^k_{mn}(A)v^m v^n - \frac{1}{2}\Gamma^k_{mn}(A)W^m W^n \end{aligned}$$

является следствием равенства (10.4.26) Равенство

$$(10.4.28) \quad \begin{aligned} &- t^k \\ &+\frac{1}{2}\Gamma^k_{mn}(A)v^m t^n \\ &-\frac{1}{2}\Gamma^k_{mn}(A)w^m v^n + \frac{1}{2}\Gamma^k_{mn}(A)w^m t^n \\ &+\frac{1}{2}\Gamma^k_{mn}(A)t^m v^n + \frac{1}{2}\Gamma^k_{mn}(A)t^m w^n - \frac{1}{2}\Gamma^k_{mn}(A)t^m t^n \\ &= -\frac{1}{2}\Gamma^k_{mn}(A)v^m w^n \end{aligned}$$

является следствием равенства (10.4.27) Равенство

$$(10.4.29) \quad \begin{aligned} &\Gamma^k_{mn}(A)t^m t^n + (\Gamma^k_{mn}(A)v^m + \Gamma^k_{mn}(A)w^m \\ &+\Gamma^k_{nm}(A)v^m + \Gamma^k_{nm}(A)w^m - 2\delta^k_n)t^n \\ &+2T^k_{nm}(A)v^m w^n = 0 \end{aligned}$$

является следствием равенства (10.4.28) Равенство (10.4.21) является следствием равенства (10.4.29) □

Ответить на вопрос, имеет ли система уравнений (10.4.21) решение - задача непростая. Однако есть другой способ найти координаты вектора $t$.

Мы проведём геодезическую
$L_{v+w}$ через точку $A$, используя
вектор $v + w$ как касательный век-
тор к $L_{v+w}$ в точке $A$. Мы прове-
дём геодезическую $L_{w+v}$ через точ-
ку $A$, используя вектор $w + v$ как
касательный вектор к $L_{w+v}$ в точке
$A$. Мы проведём геодезическую $L_u$
через точку $A$, используя вектор $u$

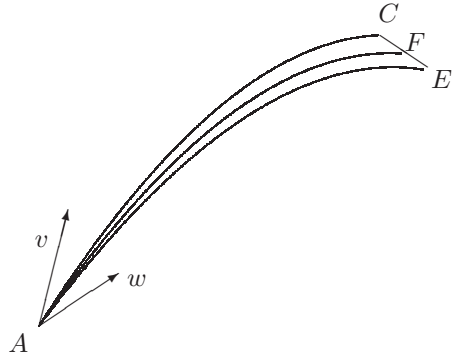

$$u^k = v^k + w^k$$

как касательный вектор к $L_u$ в точ-
ке $A$.

Согласно теоремам 10.4.2, 10.4.4 точка $F$ является серединой отрезка $EC$.
Следовательно, мы можем рассматривать отрезок $AF$ как медиану треуголь-
ника $ACE$. Согласно теореме 10.4.4 мы можем отождествить отрезок $FC$ и
вектор $t$. Следовательно, теорема 10.4.2 даёт нам способ найти координаты
вектора $t$.

## 10.5. Некоммутативный модуль

В разделе 10.4, мы рассмотрели возможность изучения аффинной геомет-
рии на многообразии аффинной связности. Эта геометрия имеет две особен-
ности. Множество векторов не замкнуто относительно сложения и операция
сложения может быть некоммутативной.

Мы ещё не готовы рассмотреть первую проблему, однако мы можем рас-
смотреть вопросы, связанные с некоммутативностью суммы векторов. Пред-
ставление

$$f : D \longrightarrow_* G$$

коммутативного кольца $D$ в произвольной группе $G$ называется некоммута-
тивным модулем. Это представление во многом похоже на модуль, поэтому все
теоремы о структуре модуля верны. Однако вопрос о структуре базиса остаётся
открытым.

Вообще говоря,

$$av + bw \neq bw + av$$

Поэтому возникает вопрос, какое множество группы $G$ мы хотим рассмотреть
в качестве базиса.

Мы можем построить базис аналогично тому, как мы строим базис модуля.
Тогда этот базис должен допускать выражение вида

$$av + bw + cv$$

Либо мы можем потребовать, чтобы элементы базиса в линейной комбина-
ции имели строгий порядок. При этом предполагается, что если $(v, w)$ - базис
некоммутативного модуля $V$, то для любого выражения $bw + av$ существуют $c$,
$d \in D$ такие, что

$$cv + dw = bw + av$$

# Список литературы

# Предметный указатель

# Специальные символы и обозначения